# Considerations about Andrews-Curtis invariants based on sliced 2-complexes

Holger Kaden

## CONTENTS

list of figures

# 1 Abstract

We consider a 2-complex in a particular form, called the Quinn model of a 2-complex. It can be sliced in graphs, where a change from one graph to another can be organized by a sequence of local transitions, which are described in a list of F. Quinn **[Q1]**.
The decomposition of that 2-complex into graphs has to be translated into an algebraic context (for example Topological Quantum field theory (TQFT)) to construct suitable potential invariants under 3-deformations. These invariants are accessible for computation by using a supercomputer and the results may yield a counterexample to the Andrews-Curtis conjecture.
To achieve invariance under 3-deformations, there are obvious topological relations among the local transitions, for example to deform a bubble out of a rectangle.
In this paper our main result is that we contribute a complete list of such topological relations in a totally geometric fashion. One outcome of our considerations is that the corresponding list of F. Quinn **[Q1]** is extended by an additional relation which takes care of locally changing a slicing. We do not know so far whether this relation is a consequence of the remaining ones. But it may be crucial for further work to focus on such subtleties, as algebraic "simplifications", where this question is bypassed, so far have been unable to distinguish between simple homotopy and 3-deformations at all.
In our introduction we summarize some known results on the situation when passing to Algebra; and in § 8 we calculate an example of an algebraic TQFT in order to demonstrate that our additional relation holds.
All considerations are carried out for 2-complexes with two generators and two defining relators. But the results also hold in the general case.

# 2 Introduction

In **[Q1]** F. Quinn proposes a new approach, to use Topological Quantum Field Theory (TQFT) for the problem of detecting counterexamples to the Andrews-Curtis Conjecture. He defines a TQFT on sliced 2-complexes, i.e:
Given a 2-complex in general position, this is sliced into graphs. Each graph is assigned to a state module and each local change from a slice to its neighbouring slice associates a homomorphism between the corresponding state modules. By composing the local changes through the whole 2-complex and by considering that as a bordism from the empty set to the empty set, the resulting homomorphism is a multiplication by an element of R. R is a ring that defines the state module of the empty set.
F. Quinn presents a list of local transitions and relations and sketches a proof, that these are sufficient to get a well defined A-C invariant for 2-complexes.
These invariants survive a stabilization phenomenon:
Let $K^2$ and $L^2$ be simple-homotopy equivalent 2-complexes. By attaching sufficiently many 2-spheres to $K^2$ and $L^2$, they are Andrews-Curtis equivalent, but the associated homomorphism to the 2-sphere is a finite sum of squares, hence it can be set to 0 in a finite field $Z_p$:

$$Z_{(K^2 \vee nS^2)} = Z_{K^2} \otimes nZ_{S^2} = Z_{L^2} \otimes nZ_{S^2} = Z_{(L^2 \vee nS^2)}$$

But since $Z_{S^2} = 0$ we can **not** conclude $Z_{K^2} = Z_{L^2}$.

We present results about TQFT:
I. Bobtcheva gives in **[Bo]** an algebraic proof, that Quinn's invariant is an Andrews-Curtis invariant. She defines a new invariant, expressed as a composition of morphisms in a semisimple tensor category. The morphism can be translated to diagrams. Hence the proof of the A-C invariance can be transformed to prove, that the corresponding diagrams are equivalent by using general identities for a semisimple tensor category. These identities are also presented as equivalent diagrams.
It can be shown by the same methods as above, that the appearance of a generator in a relation induces the same homomorphism in Bobtcheva's invariant as the circulator in the Quinn invariant, which is the associated homomorphism to the sliced topological part in the Quinn model, where an attaching map of a 2-cell wraps around that generator **[Q2]** .
In **[Bo/Q]** the authors show, that TQFT fail to detect A-C counterexamples for 2-complexes $K^2$ with Eulercharacteristic $\chi(K^2) > 0$. They consider reductions of invariants for 4-dimensional thickenings of 2-complexes $K^2$. They can exclude those where the invariant only depends on the 3-dimensional boundary.
If $\chi(K^2) > 0$, then it follows from **[Hu]** , that this boundary appears in a dual 2-complex to $K^2$ in $S^4$ , which is determined by its 1. homology.

In **[BLM]** the authors show, that in the original case of a contractible 2-complex for the A-C Conjecture, the projection of the free generators into a finite testgroup must fail. Since in all computable cases the circulator has finite order, these tests do not detect A-C counterexamples.

K. Müller gives in his theses **[Mül]** excluding criteria for finite testgroups to be TQFT counterexamples to the Andrews-Curtis conjecture. He connects results of a Browning paper with algebraic criteria about higher commutators for simple-homotopy equivalence with the fact, that the circulator has a prime order in these cases.

These notes are divided into several chapters:
In chapter 3 we recall basics about simple-homotopy equivalence. It includes the Theorem of P. Wright to get the connection between the topological and algebraic view for 2-complexes under 3-deformations.
In chapter 4 we present the concept of a sliced 2-complex. We describe the Quinn list of local transitions to change from a slice to its neighbouring slice and introduce the definition of a topological relation among these transitions. We explain the Quinn model, point out and solve a problem for the $T_3$ move. That solution results in a new topological relation, in the sense that its associated algebraic relation is not in the Quinn list.
In chapter 5 we present the sequence for a special case of slices which we use to construct and annihilate a pair of leftside loops. We show how to transfer a rightside loop into a leftside loop and construct the sequence of slices during that process. We provide the pullback of the general case to that special case.
In chapter 6 we consider two different possibilities of organizing the sequence of slices in the Quinn model, if the attaching curve is changing under twists around $\pi$ or a composition (2-times) of these. One possibility is that we realize the twist on the generator cylinder, and one is (by considering the muliplication of two relators) to realize it at the other 2-cell, by sliding one cell through the other.
In chapter 7 we consider a sliced 2-complex in the Quinn model under Q - transformations. These are composed as a sequence of Matveev moves. For each such move, we provide the sequence of slices. We use the result of the former chapter (the twist is realized on the generator cylinder) to omit considerations about the twisted attaching curve. We also include the extended prolongation in our consideration.
In chapter 8 we transfer the $T_3$ move problem (in chapter 4) into the Quinn model, and see that the new topological relation is required for the list of topological relations. It is an open question, if it is also an algebraic relation. We split the sequence of slices into one sequence, that solves the $T_3$ move problem and one, that describes the passing of a vertex in the local vertex model. We set the sequences in the context of roottrees in TQFT and compare their associated homomorphisms for a chosen tensor category. In this example they are the same.
In chapter 9 we confirm in our example, that the associated homomorphism is the identity for the given algebraic relation deforming a bubble out of a rectangle.
We finish with chapter 10 and present stages of further research. We suggest especially to consider sliced 3-cells by connecting an algebraic criterion for simple homotopy equivalence in **[HoMeSier]** with a topological criterion in **[Q3]**.

I want to thank Prof. Metzler for his engagement, support and our discussions. Special thanks to Prof. F. Quinn for many discussions about TQFT and for his encouraging feedback. I would also like to thank Prof. S. Matveev, Dr. C. Hog-Angeloni and Dr. K. Müller of helpful explanations about TQFT and the s-move. Further thanks to the students of the research group of Prof. Metzler for discussions about this notes.

Special thanks  to Jan Hofmann and Ingrid Metzler for their assistance with the English language. Jan Hofmann also contributed a feedback from a viewpoint of a non-expert.

# 3 Basics

In this chapter we present the main topics of simple-homotopy theory. We refer to **[HoMeSier]** for details also to **[Mat]** for Matveev-moves.

## *3.1 Simple-homotopy equivalence and Q**-transformations*

### 3.1.1 CW-complex

A CW-complex K of a space I K I is a decomposition in cells, inductively constructed by increasing dimension.
Start with a discrete set of points $K^0$, these are the 0-cells of K.
Take a disjoint family of line segments and attach them to $K^0$ by identifying their boundary points with points in $K^0$, but project the interior of each line segment homeomorphic. These are the 1-cells and the set of 1-cells is denoted as $K^1$. The union of cells $K^0$ and $K^1$ is called the 1-skeleton.
Let $K^{n-1}$ be constructed.
Take $K^{n-1}$ together with a family of closed n-discs $D_i^n$ and construct the quotient space $K^n$ via continuous maps $\varphi_i: \partial D_i^n \rightarrow K^{n-1}$, such each $x \in \partial D_i^n$ is identified with $\varphi_i(x)$, then the interior of $D_i^n$ projects homeomorphic to an n-cell $e_i^n$ . $\varphi_i$ is called an attaching map for $e_i^n$ .
The n-skeleton is the union of all cells with maximal dimension n.
$I K I = \bigcup I K^n I$ is assigned the weak topology according to the closure of $\bar{e}_i^n$, i.e. a subset of I K I is closed, iff it is intersection with each $\bar{e}_i^n$ is closed.

We formalize the homeomorphic projection on the interior of a cell to the

**Definition of a characteristic map**
A map $\Phi: D^n \rightarrow K$ is characteristic for an n-cell $e^n$ of K, if $\Phi$ maps the interior of $D^n$ homeomorphically onto $e^n$ and $\Phi\, I_{\partial D^n}$ is an attaching map.

As a simple example for a CW-complex we present the standard 2-complex:
Let $P = \langle a_1,\ldots,a_n / R_1,\ldots,R_m \rangle$ be the presentation of a standard 2-complex $K^2$. To construct this complex, start with a point as 0-cell. Attach n disjoint line segments and identify their boundary points with the 0-cell, then we get a bouquet of n loops, which are the oriented 1-cells denote as $e_1^1, e_2^1, \ldots ,e_n^1$. For each relation $R_k$ we choose a characteristic map $\Phi_k$, attach the 2-disc by identfying the edges $a_i$ with the boundary of the disc according their appearance and orientation. P presents the fundamental group $\pi_1(K^2)$.

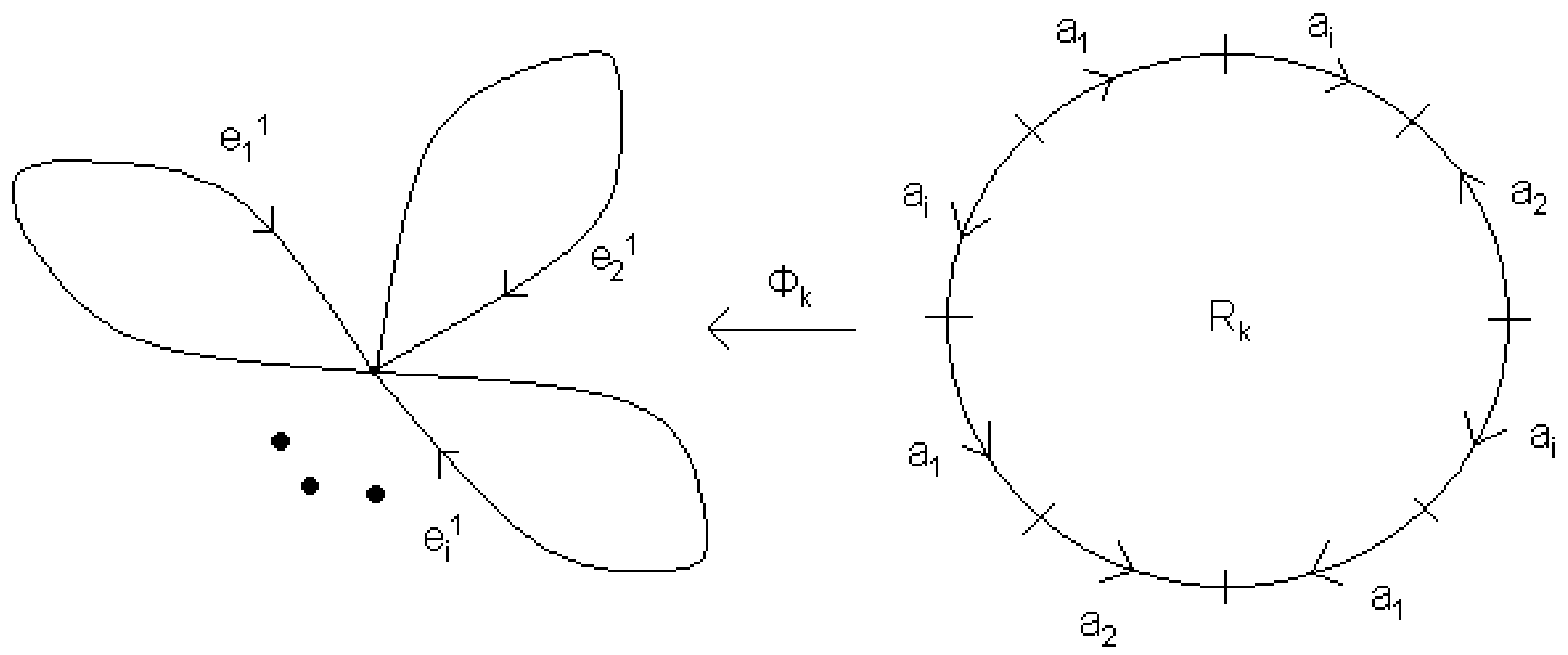

**Figure 1- sh-equivalence- presentation of a standard 2-complex**

### 3.1.2 Elementary expansion and collaps

We want to describe the attachment of an n-ball $D^n$ along an (n – 1)-ball $D^{n-1}$ on its boundary, but we want to point out rather the idea than the exact definition.
Let $\varphi: D^{n-1} \rightarrow L$ an attaching map for an n-cell $e^n$, such that $K = L \vee e^n \vee e^{n-1}$ where $e^n, e^{n-1} \notin L$ with $\partial e^n = e^{n-1} \vee \delta e^n$. That is we consider a characteristic map $\Phi: D^n \rightarrow K$ for $e^n$, such that an (n-1)-disc $D^{n-1}$ in the boundary of $D^n$ maps via $\Phi$ to a union of n -1-cells (denoted as $\delta e^n$) $\in L$, and $\Phi$ is a characteristic map from the remaining boundary $\partial D^n - D^{n-1}$ to $e^{n-1} \notin L$. The cell $e^{n-1}$ is free from n-cells of L and is called a free cell.

- K is an elementary n-dimensional expansion of L via the n-cell $e^n$.
- L is the result of an elementary n-dimensional collaps from K , i.e. start with the free face $e^{n-1}$, push the n-dimensional material of $e^n$ and $e^{n-1}$ to the remaining boundary $\delta e^n$ of $e^n$.

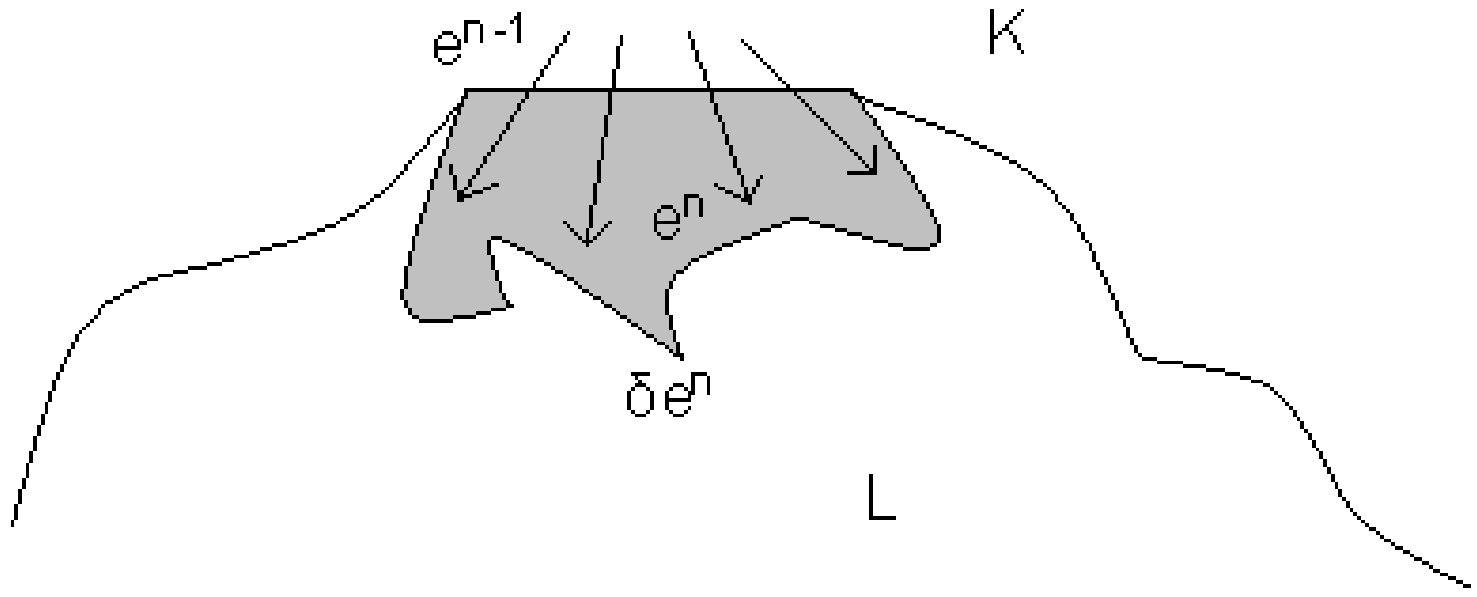

**Figure 2- sh-equivalence- elementary expansion /collaps**

An arbitrary sequence of collapses and expansions is called a deformation, a simple-homotopy equivalence from K to L is a map homotopic to a deformation, we write: $K \nearrow\!\!\searrow L$ and say K and L are simple-homotopy equivalent (sh-equivalent). If the maximal dimension of the cells used in the simple-homotopy equivalence is n, then we add this in our notion as

$$K \overset{n}{\nearrow\!\!\searrow} L$$

If K and L are sh-equivalent, then the next theorem gives an answer about the required dimension.

**Theorem 1 (Wall)**
Let f: K → L a simple-homotopy equivalence of connected, finite CW-complexes and n = max (dim K, dim L), then f is homotopic to a deformation $K \overset{n+1}{\nearrow\!\!\searrow} L$ for $n \geq 3$.

For n = 2 the answer is still unknown and is the topological formulation of the

**Andrews-Curtis Conjecture**
Suppose $K^2$, $L^2$ are finite, 2-dimensional CW-complexes such that $K^2$ and $L^2$ are sh-equivalent , then they are sh-equivalent by a 3-deformation, i.e.

$$K^2 \nearrow\!\!\searrow L^2 \Rightarrow K^2 \overset{3}{\nearrow\!\!\searrow} L^2$$

The next lemma shows the useful fact, that simple-homotopy equivalence is induced by homotopic attaching maps, as formulate in the next lemma:

**Lemma 1**
Let K a n-dimensional complex, and φ, ψ homotopic maps from $S^{n-1}$ into the n – 1 skeleton $K^{n-1}$, attaching the n-cell $e^{\varphi}$ respectively $e^{\Psi}$ to K, then

$K_{\varphi} \overset{n+1}{\nearrow\searrow} K_{\psi}$, where $K_{\varphi} = K \vee e^{\varphi}$ and $K_{\psi} = K \vee e^{\Psi}$

Use the homotopy of the attaching map (indicated from the arrows) to attach the n+1-cell $e^{n+1}$ between $e^{\varphi}$, $e^{\Psi}$ and the region between $\varphi(S^{n-1})$and $\psi(S^{n-1})$, both $e^{\varphi}$, $e^{\Psi}$ are free faces of $e^{n+1}$. Then
$K \vee e^{n+1}$ collaps to $K_{\varphi}$ using the free face $e^{\Psi}$ and
$K \vee e^{n+1}$collaps to $K_{\psi}$ using the free face $e^{\varphi}$:

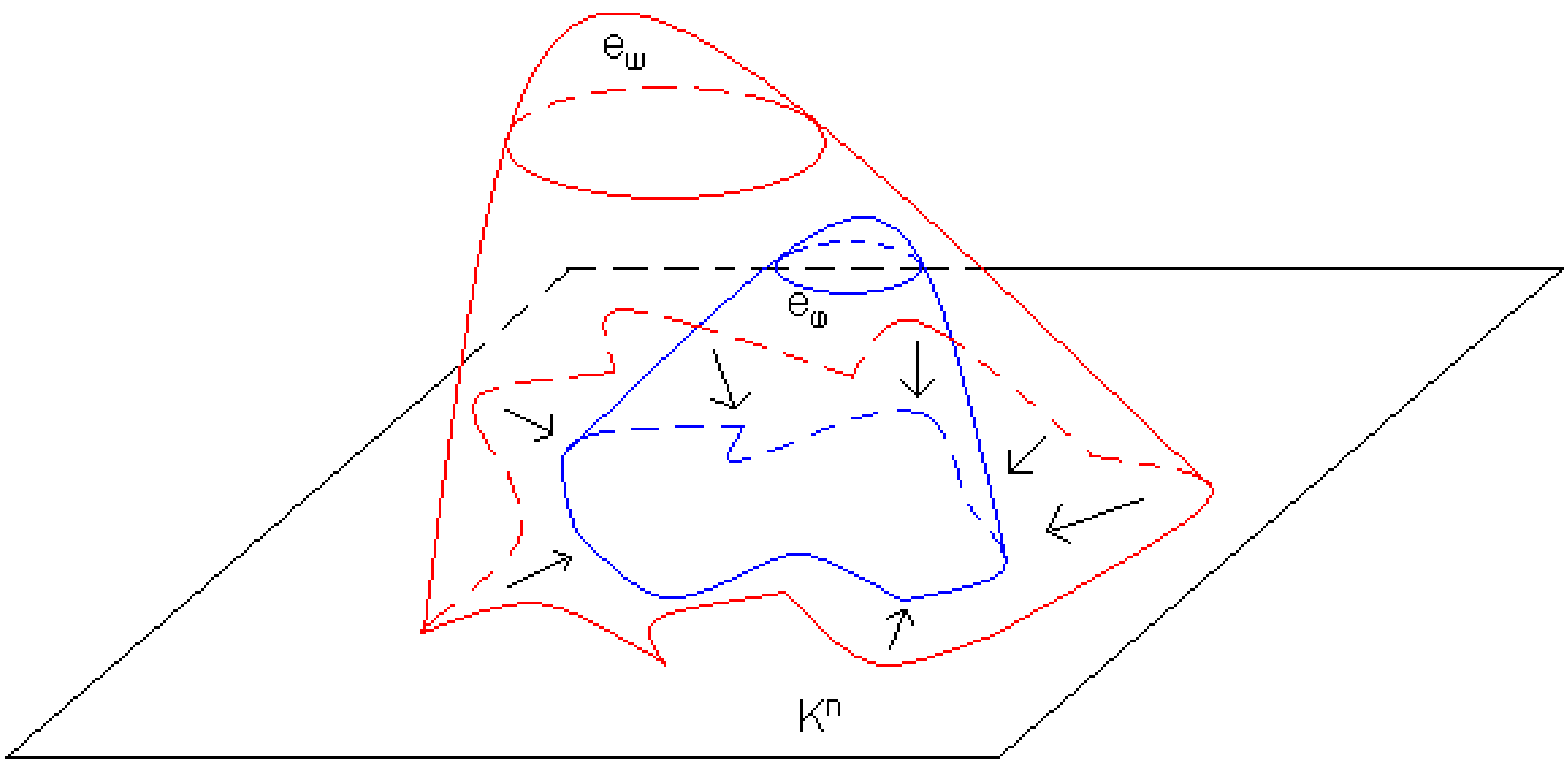

**Figure 3- sh-equivalence- homotopy of attaching maps**

### 3.1.3 Andrews-Curtis conjecture and Q**-transformations

Getting back to the Andrews-Curtis conjecture, there is a translation of 3-deformations of 2-complexes which expresses these as certain transformations of presentations.

Let $P = \langle a_1,\dots,a_n / R_1,\dots,R_m \rangle$ be a finite presentation of the fundamental group of a 2-complex. A Q-transformation of P is a finite sequence of 3 types:

a) $R_i \to wR_iw^{-1}$ for some i and $w \in F(a_k)$ (conjugation)
b) $R_i \to R_i^{-1}$ for some i (inversion)

c) $R_i \rightarrow R_i R_k$ or $R_k R_i$ ($i \neq k$) (multiplication on right or left)

we additionally replace a generator $a_i$ in the relators by a finite sequence of:

d) $a_i \rightarrow a_i^{-1}$

e) $a_i \rightarrow a_i a_k$ or $a_k a_i$ ($i \neq k$)

with the remaining $a_m$ ($i \neq m$) stay unchanged then we call these transformations Q*-transformations.

We can extend the presentation by a new generator a and relator R = a:

f) $P = \langle a_1,\dots,a_n / R_1,\dots,R_m \rangle \rightarrow Q = \langle a_1,\dots,a_n, a / R_1,\dots,R_m, R = a \rangle$

This transition is called a prolongation.

We call the transformations a) – f) and f) inverse Q**-transformations.
The Q**-transformation leads to an equivalence relation on the set of presentations and the equivalence classes are called Q**-classes.

**Theorem 2 (P. Wright)**
There is a bijective map between Q**-classes of finite group presentations and 3-deformation-classes of compact, connected CW-complexes.

Note, that the Q**-transformations d) – f) can be replaced by generalized prolongation and its inverse:

$P = \langle a_1,\dots,a_n / R_1,\dots,R_m \rangle \rightarrow Q = \langle a_1,\dots,a_n, a / R_1,\dots,R_m, R = w^{-1}a \rangle$
where $w \in F(a_k)$ and a is a new generator.

We will use this criterion and the Q-transformations as Q**-transformations.

## *3.2 Matveev-moves*

Using theorem 2 we can study the change of a 2-complex under 3-deformations by observing the change of a presentation of a standard-2-complex under Q-transformations. In this chapter we also study local 3-deformations on a certain type of 2-complexes, called special polyhedron.

**Theorem 3 (S. Matveev)**
A special polyhedron $K^2$ 3-deforms to another special polyhedron $L^2$, if and only if it 3-deforms by a sequence of single local moves and their inverse. These moves called T-moves or Matveev-moves, are:

a) $T_1$ move:

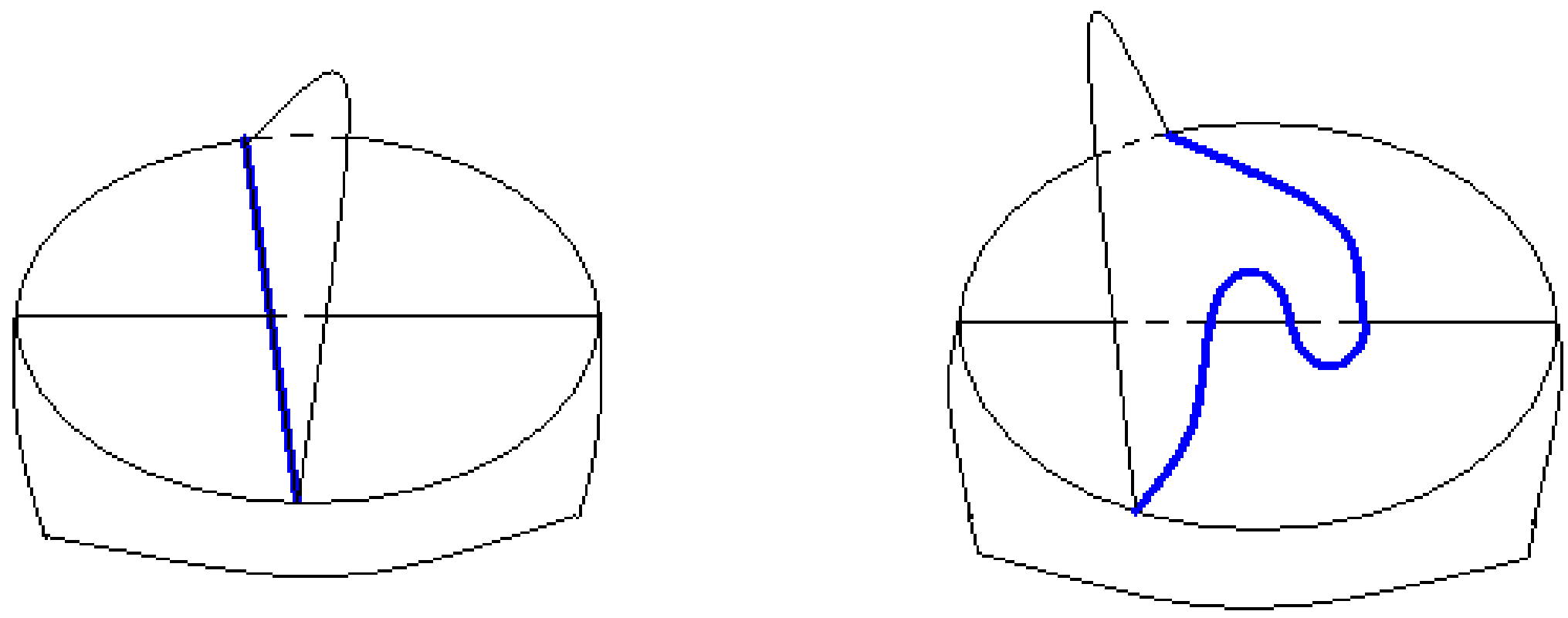

**Figure 4- Matveev moves- $T_1$ move**

b) $T_2$ move:

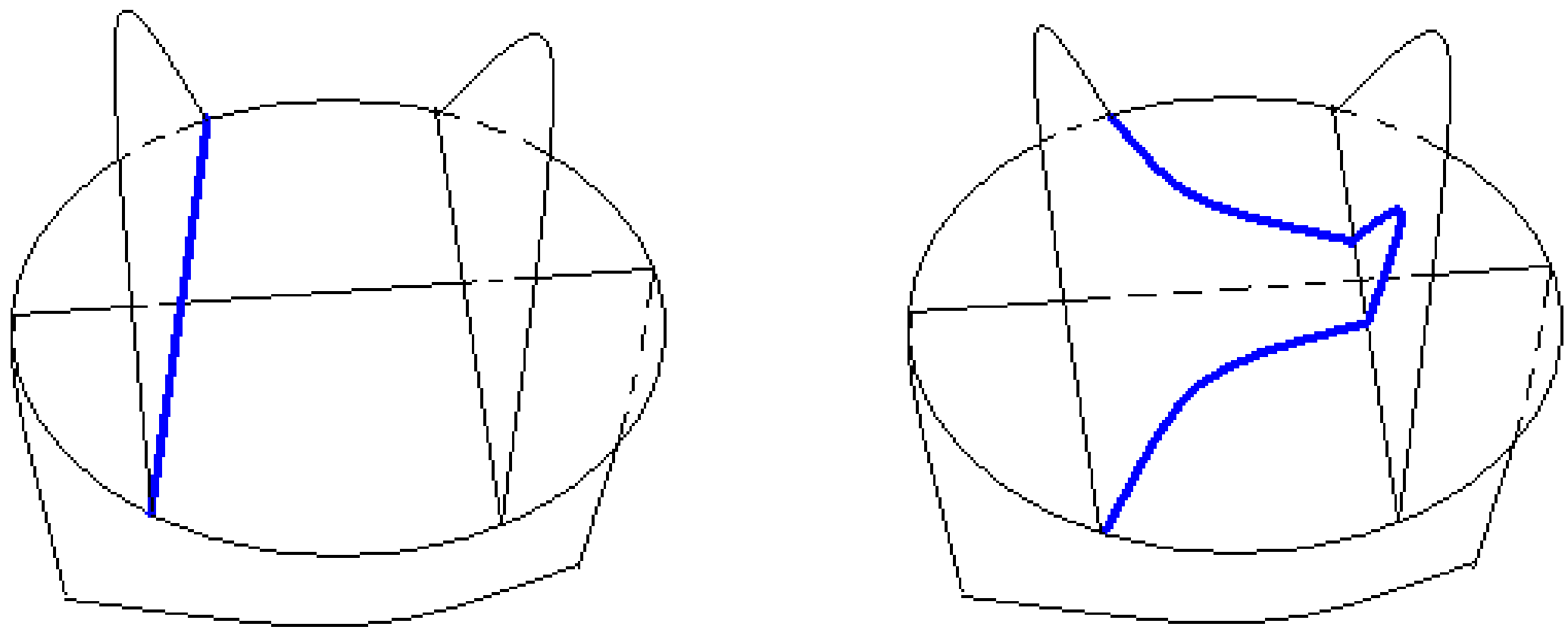

**Figure 5- Matveev moves- $T_2$ move**

c) $T_3$ move (there is also a turn to left):

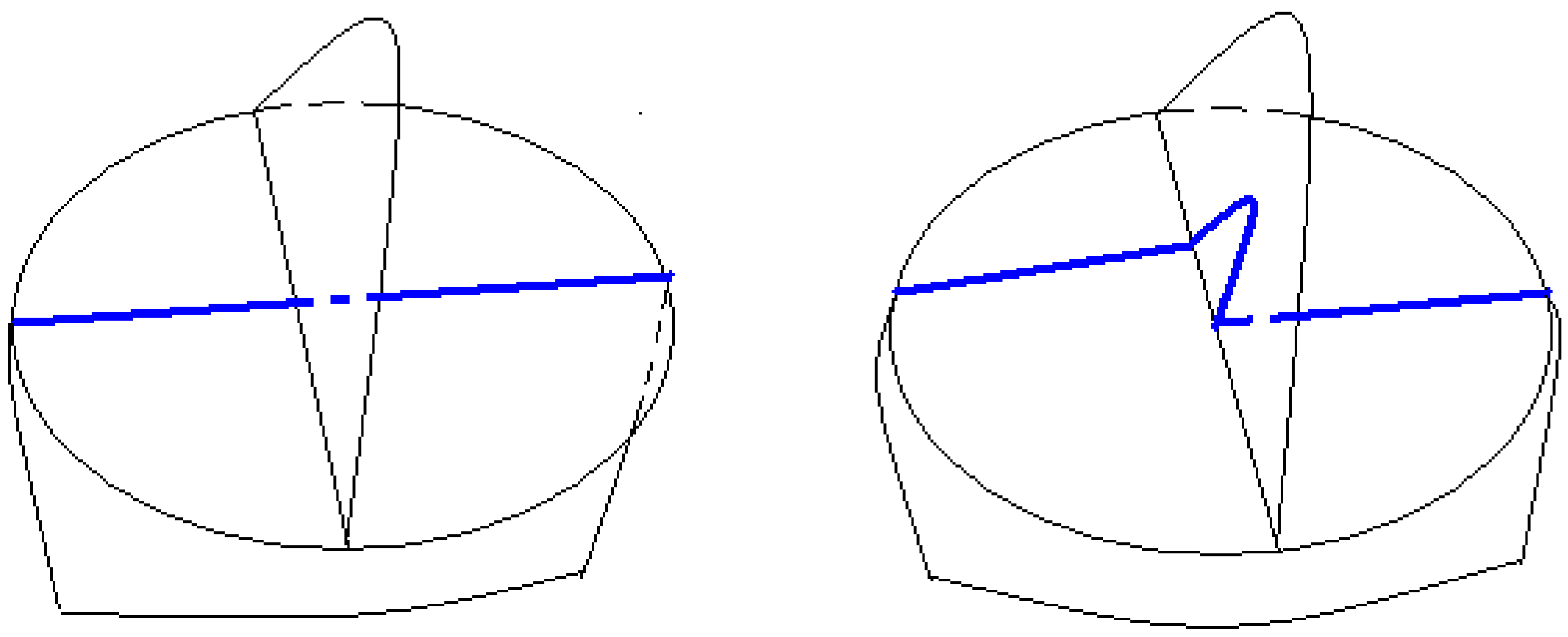

**Figure 6- Matveev moves- $T_3$ move (turn to right)**

We often use the next move, which is a composition of $T_1$ and $T_2^{-1}$, and set $T^* = T_1 T_2^{-1}$:

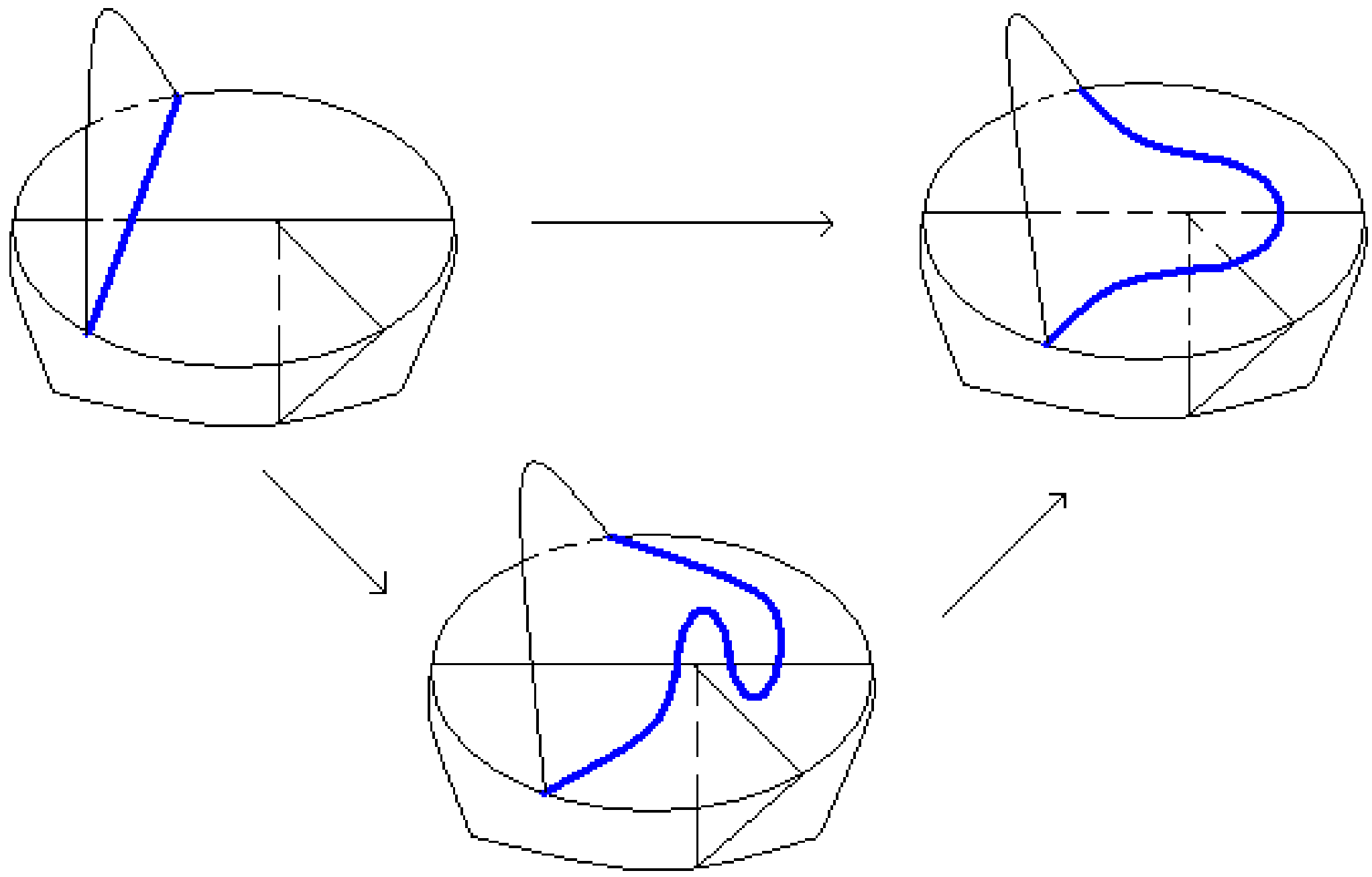

**Figure 7- Matveev moves- T* move is a composition of $T_1$ and $T_2^{-1}$**

## *3.3 Twists and loops*

To study the possible variations of an attaching curve, we first demonstrate the appearance of twists and loops from the Q-transformation “multiplication of 2 relations”.
For simplicity we assume that a relator R is of the form …aa…:
Choose an arc w in the 2-cell from the startpoint of the first a and the endpoint of the second a. Identify a with a , considering only the restricted part of the 2-cell, which has the boundary aaw. This leads to a Klein bottle with a hole in it and hence includes a twisted strip:

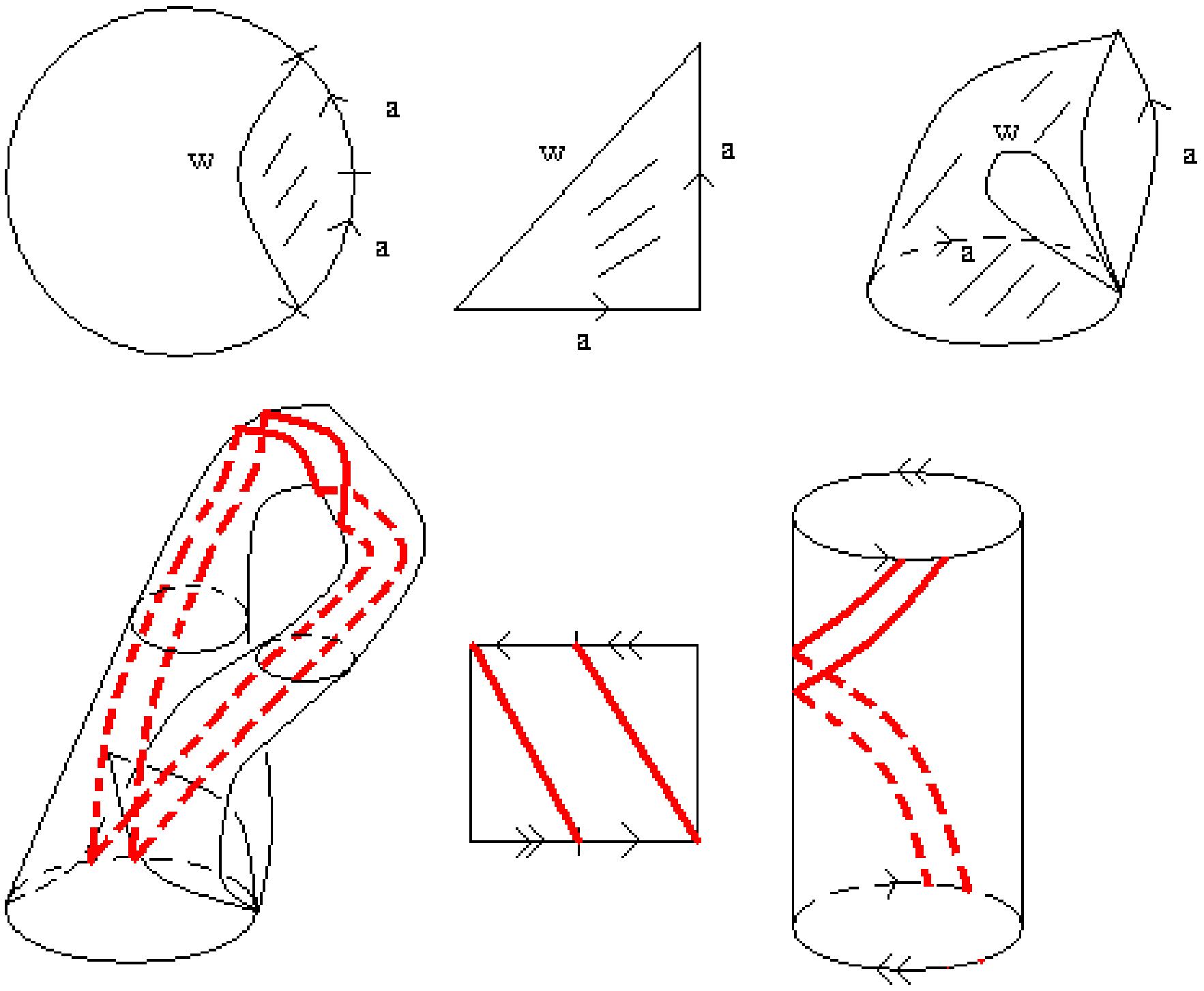

**Figure 8- Twists and loops- a twisted strip in the perforated Klein bottle**

Assume, there is an arc belonging to the attaching curve of a relator S, which is pushed on R to get R→ RS . This arc has 2 ways to pass the twisted strip in R:

In the first case the arc inherits the twist:

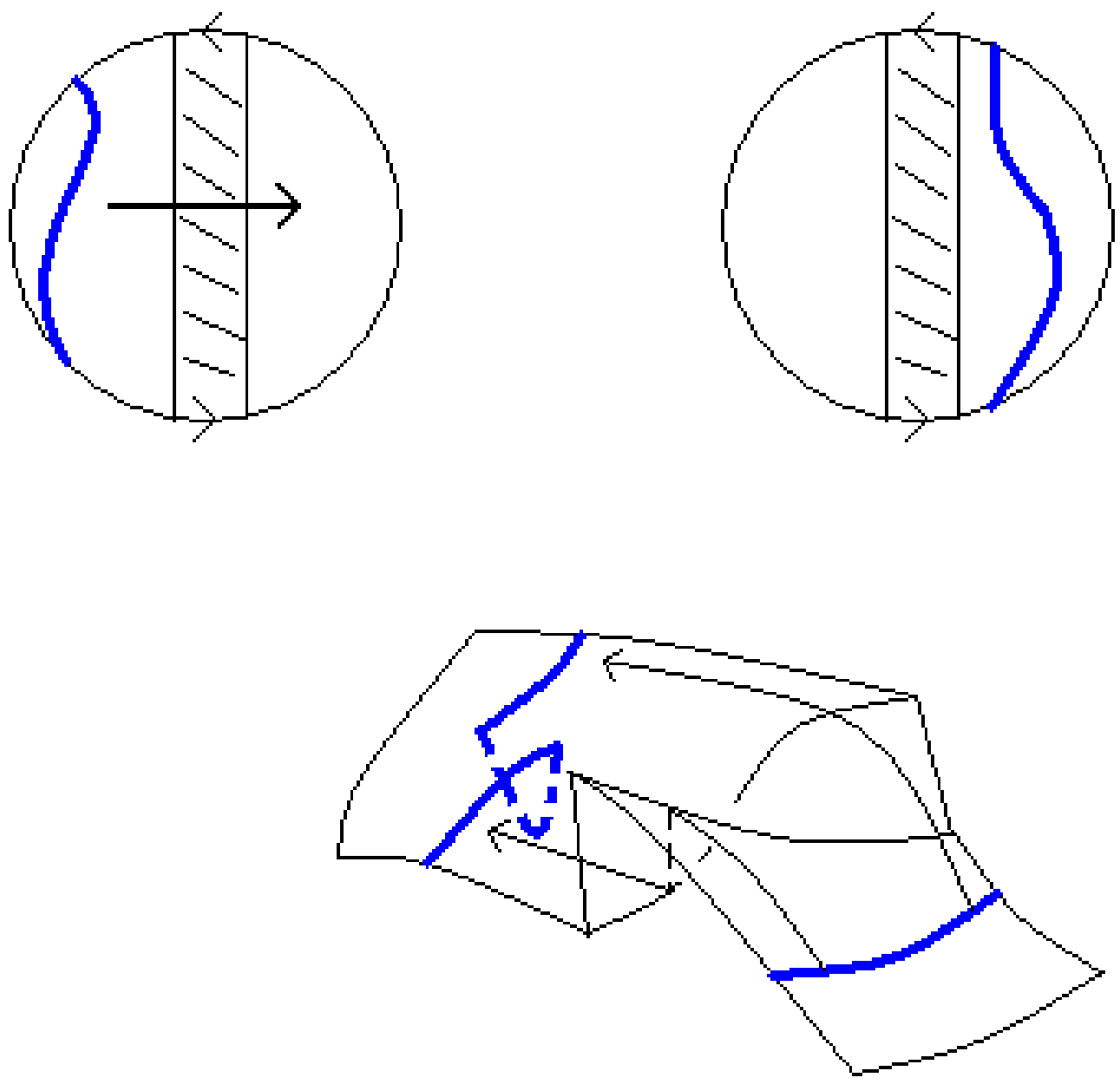

**Figure 9- Twists and loops- pass a twisted strip inherit the twist**

In the second case the arc gets a double twist which can deform to a loop:

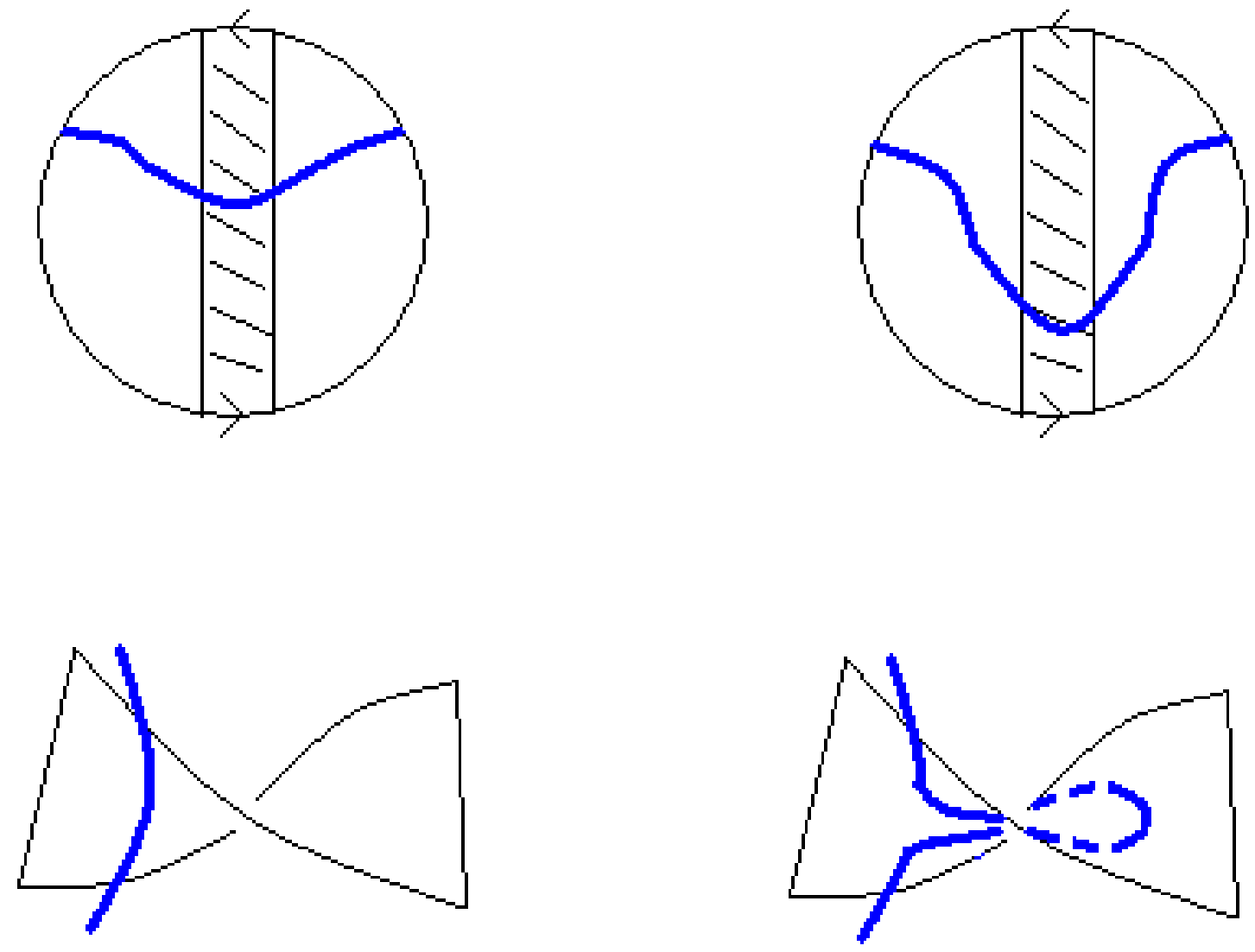

**Figure 10- Twists and loops- result in a loop**

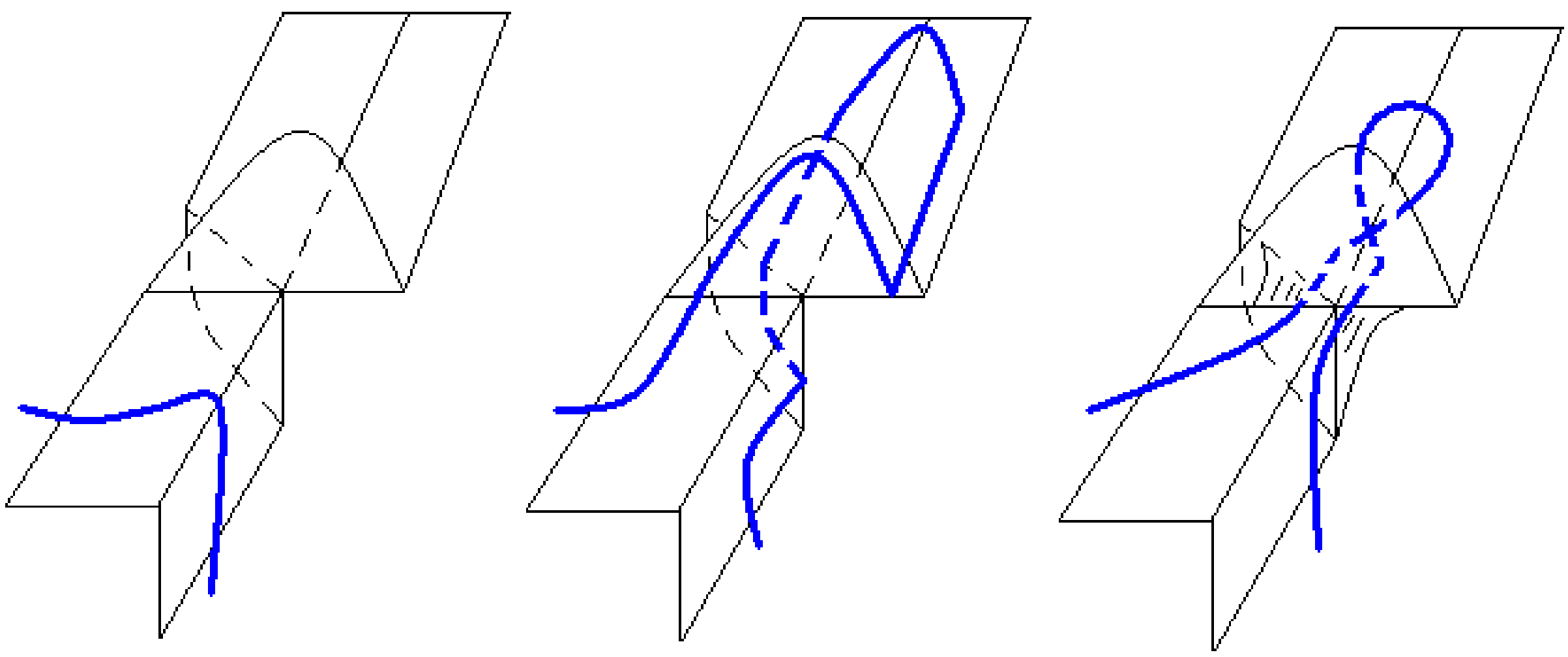

**Figure 11- Twists and loops- details to the loop case**

Our second motivation is based on the relation between “bases-up-to-conjugation“ and prismatic 1-collapsible 2-complexes. This is a 2-complex K in K× I , consisting of a wedge of circles× I with sheets embeded in K× I.
K×I collapses to the sheets and then from the sheets to the wedgepoint, hence K 3-deforms to a point. A part of the proof for the relation above requires to realize embeddings of bands $S^1\times I$ for elementary Nielsen transformations, for example
a → a, b → ba.

To insure the embedding at wedgepoint x I, detours are necessary. For details we refer to **[CoMeSa]** .

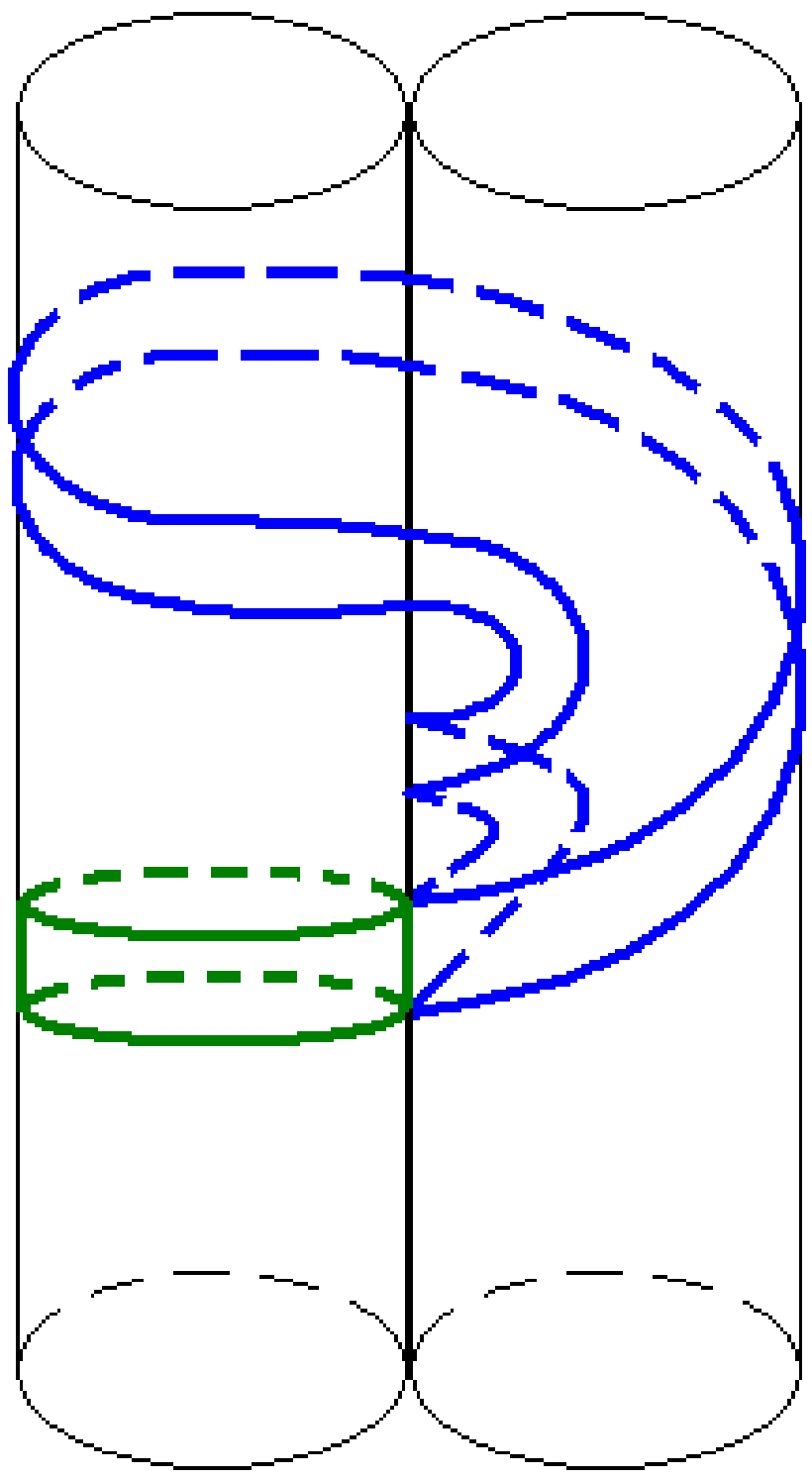

**Figure 12- Twists and loops- Nielsen transformation of $S^1\times I$ in K×I requires a detour**

# 4 Sliced 2-complexes

## 4.1 Local transitions and relations among graphs

A sliced 2-complex is a 2-complex that is sliced into graphs. The graphs have to change by local transitions. A standard example is the torus:
Starting at one end the graph is the empty set, then it becomes a circle. At the first saddlepoint it changes to a wedge of two circles and then it splits in two circles. At the second saddlepoint it becomes a wedge of two circles again, after that it changes back to a circle and finally becoming the empty set again.
In this sense the torus can be seen as a bordism from the empty set to the empty set and as a compostion of elementary bordisms, defined as bordisms from one graph to another by a single local transition.

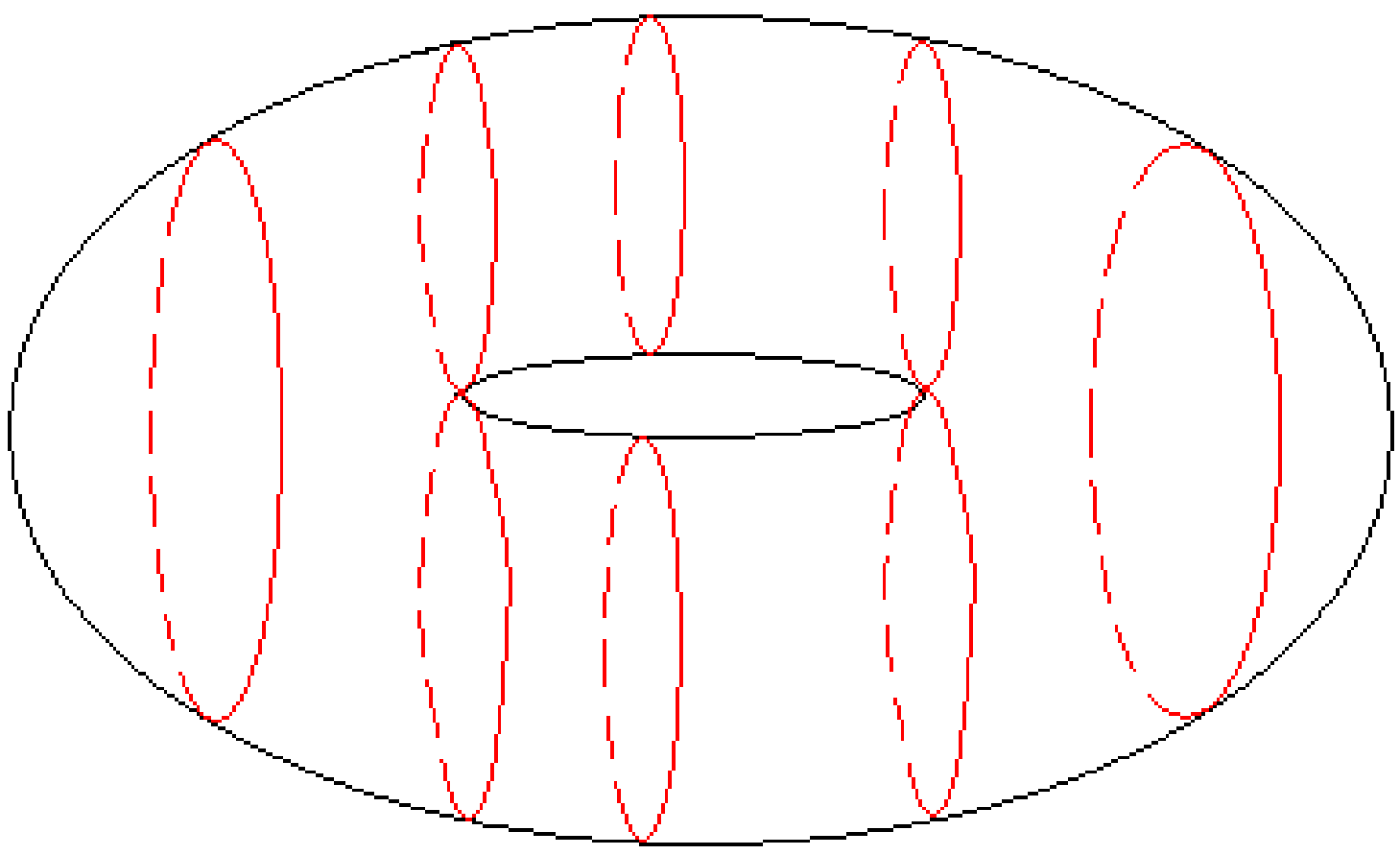

**Figure 13- local transitions- slices of a torus**

Note that in this simplified example we have to add flanges at the saddlepoints to use the local vertex model:
It consists of 4 halfdiscs, where 3 of them attached to a common line segment. The remaining halfdisc is attached on another line segment. The line segments intersect transversally in one point, which is called the vertex of the local model (see the startfigure for $T_1$ move).

Local transitions of the graphs are only permited if they are in the Quinn list. Note that the points at the graphs indicates an arbitrary continuation of the graph:

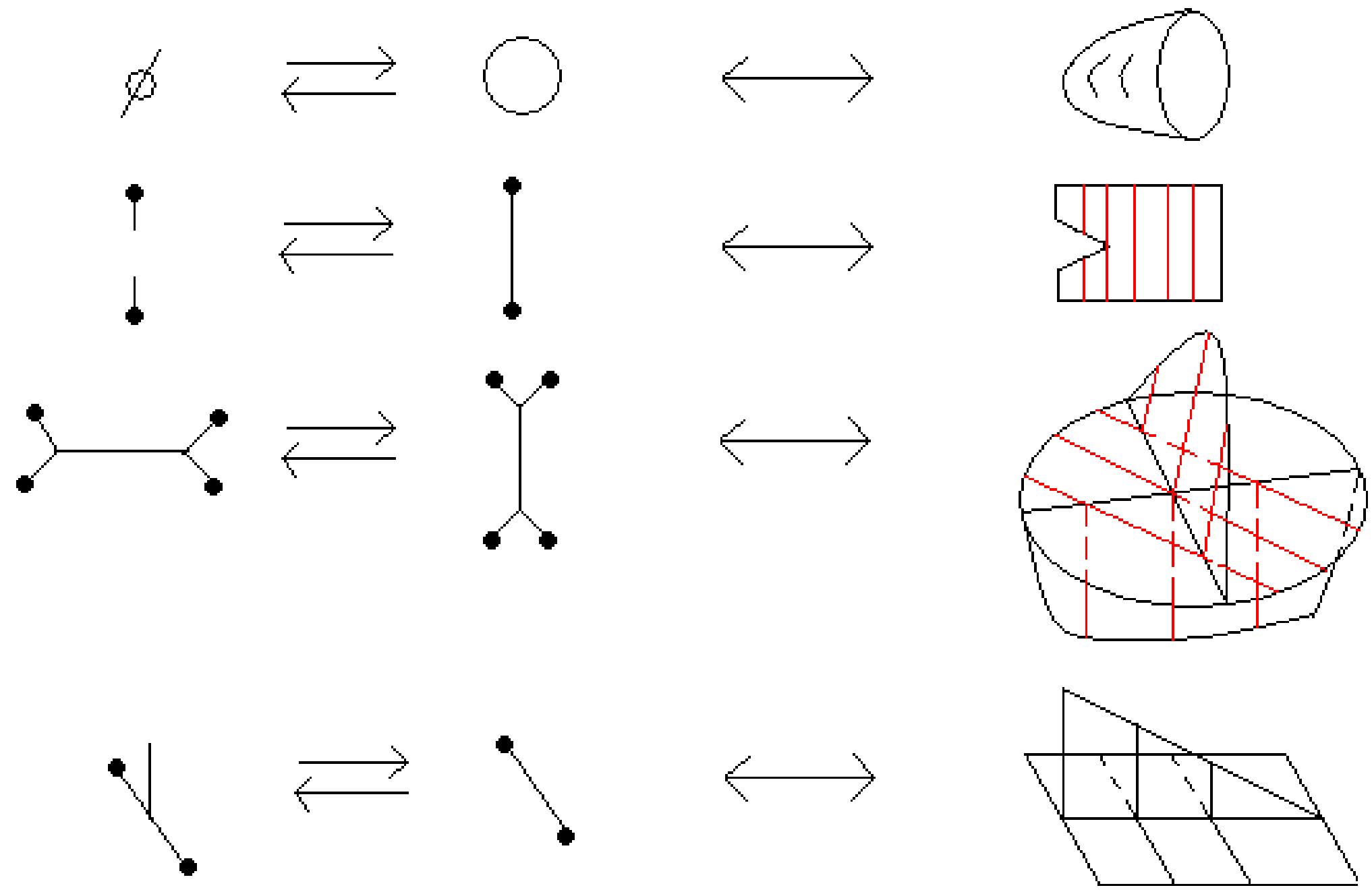

**Figure 14- local transitions- Quinn list**

To construct an A-C-invariant for sliced 2-complexes, the local transitions are settings in an algebraic context. Each 3-deformation decomposes into finitely many local 3-deformations, which are performed on local models of a 2-complex. Hence it is sufficient to study the sequence of slices. We evaluate them in the potential algebraic context.
Our setting defines an invariant, if for every sequence of slices, this sequence and the deformed sequence differ only by an identity in the chosen context. Given a sequence and its deformed sequence. Reverse one of these and compose it with the other. We get a single sequence (with the same start and endgraph), which represents an identity among the transitions in our context. It is called an algebraic relation.
We are considering only the topological sequence of slices without any algebraic context (except for chapter 8,9), which is called a topological relation. It is one of the main tasks to prepare these topological relations for all local 3-deformations.
One should consider as first examples two 3-deformations of a rectangle:

- building a bubble
- expansion with a disk.

Since the rectangle itself has only line segments as slices, we restrict our considerations on the rectangle under deformation, presented below:

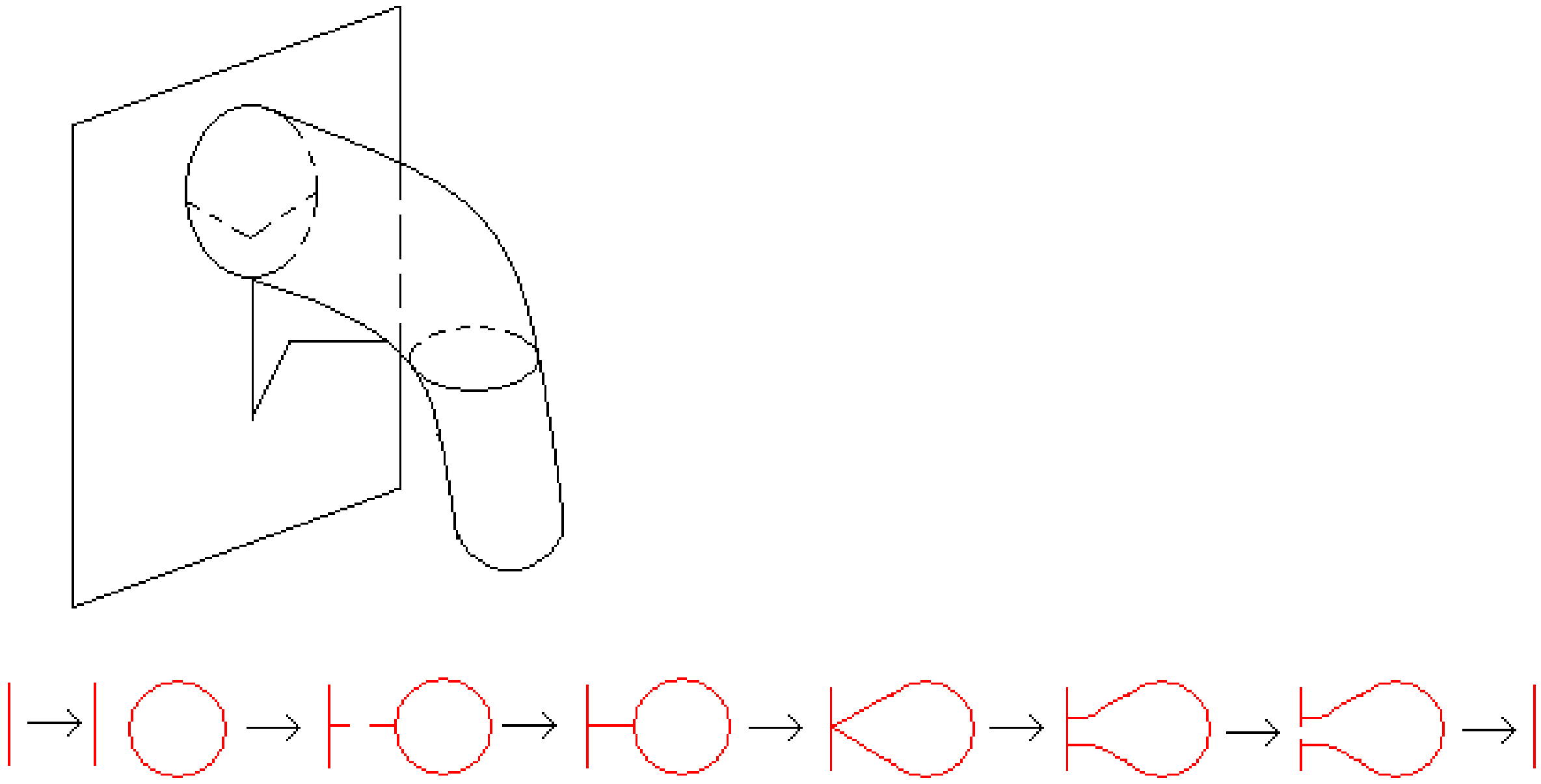

**Figure 15- topological relations- Quinn list- building a bubble**

The sequence of graphs from left to right corresponds to slicing of the figure from bottom to top.

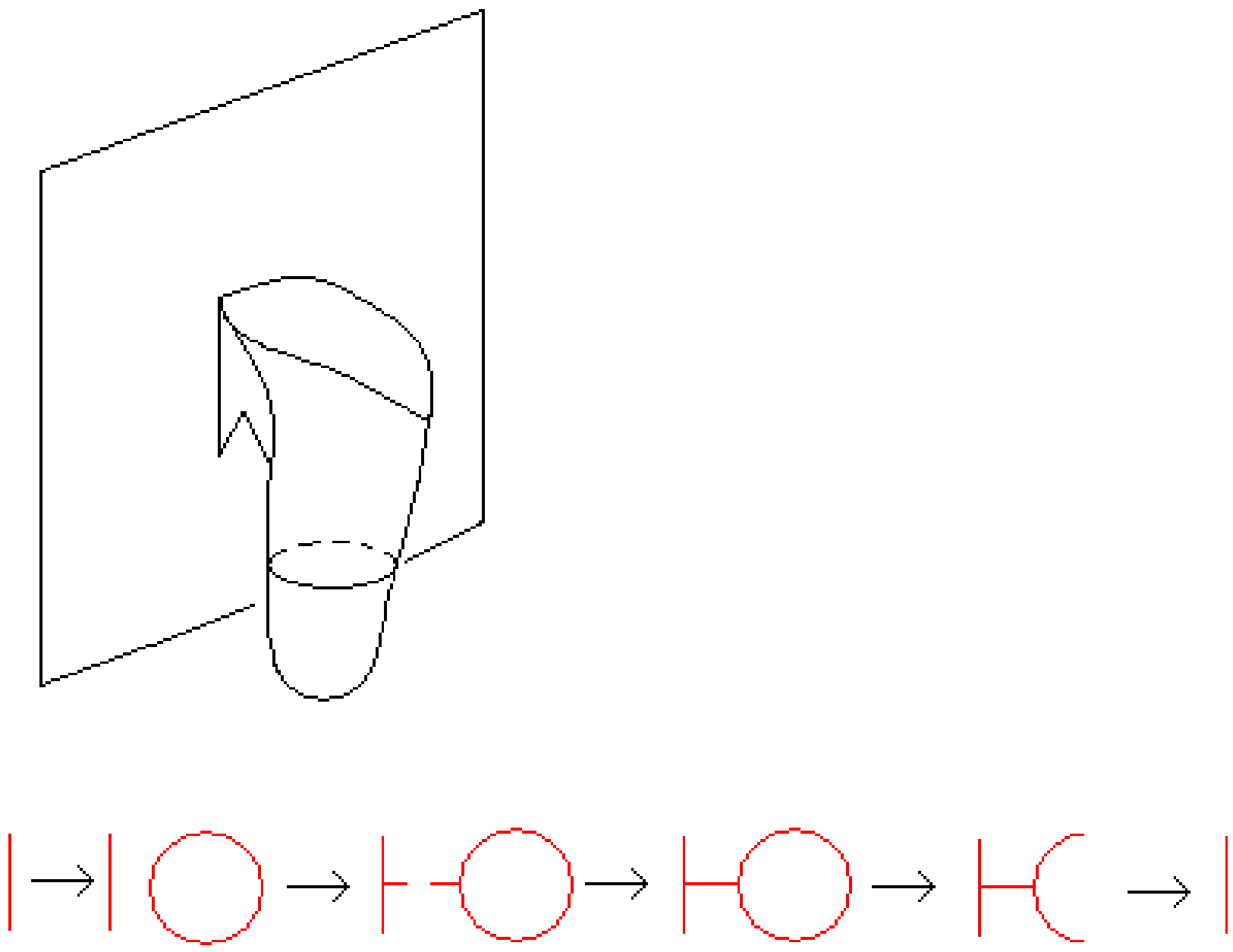

**Figure 16- topological relations- Quinn list- expansion with a disk**

Note, that the 4. local transition in the Quinn list also describes an elementary collaps or extension, hence it has to be an algebraic relation. There are further relations coming from contractible graphs, which will not be discussed in this paper.

## 4.2 The Quinn model of a sliced 2-complex

To construct the Quinn model of a 2-complex, take a presentation P = ⟨a, b / R, S ⟩ for it and glue the relators R, S into the generator cylinders a×I, b×I by their attaching curves. These wrap around the cylinders from bottom to top corresponding the appearance of the generators in the relators.

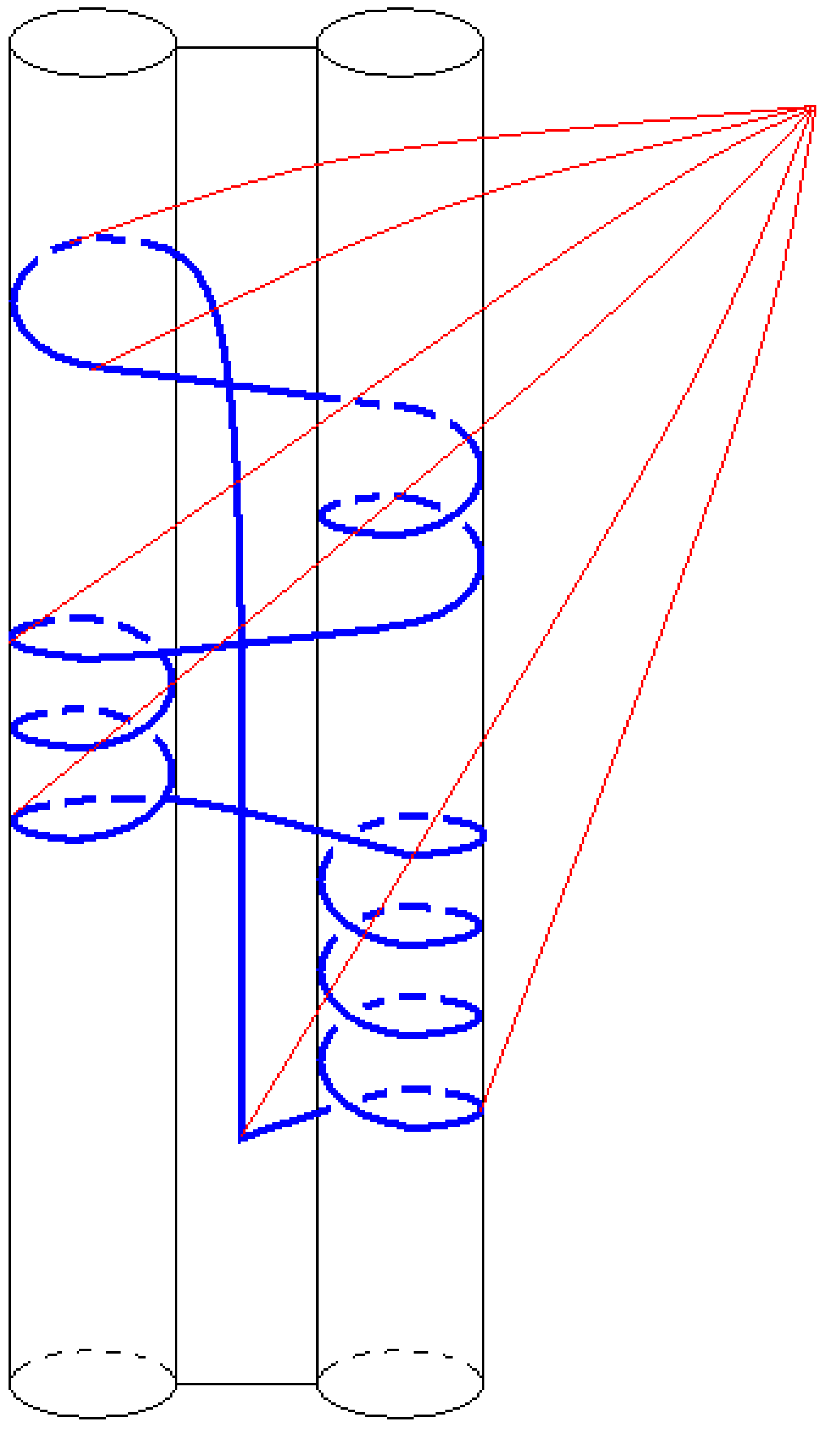

**Figure 17- Quinn model of a sliced 2-complex- attached relation**

Note that with the connecting rectangle between the cylinders the 2-complex is in general position, hence the Quinn-model of a 2-complex is the underlying space of a special polyeder.

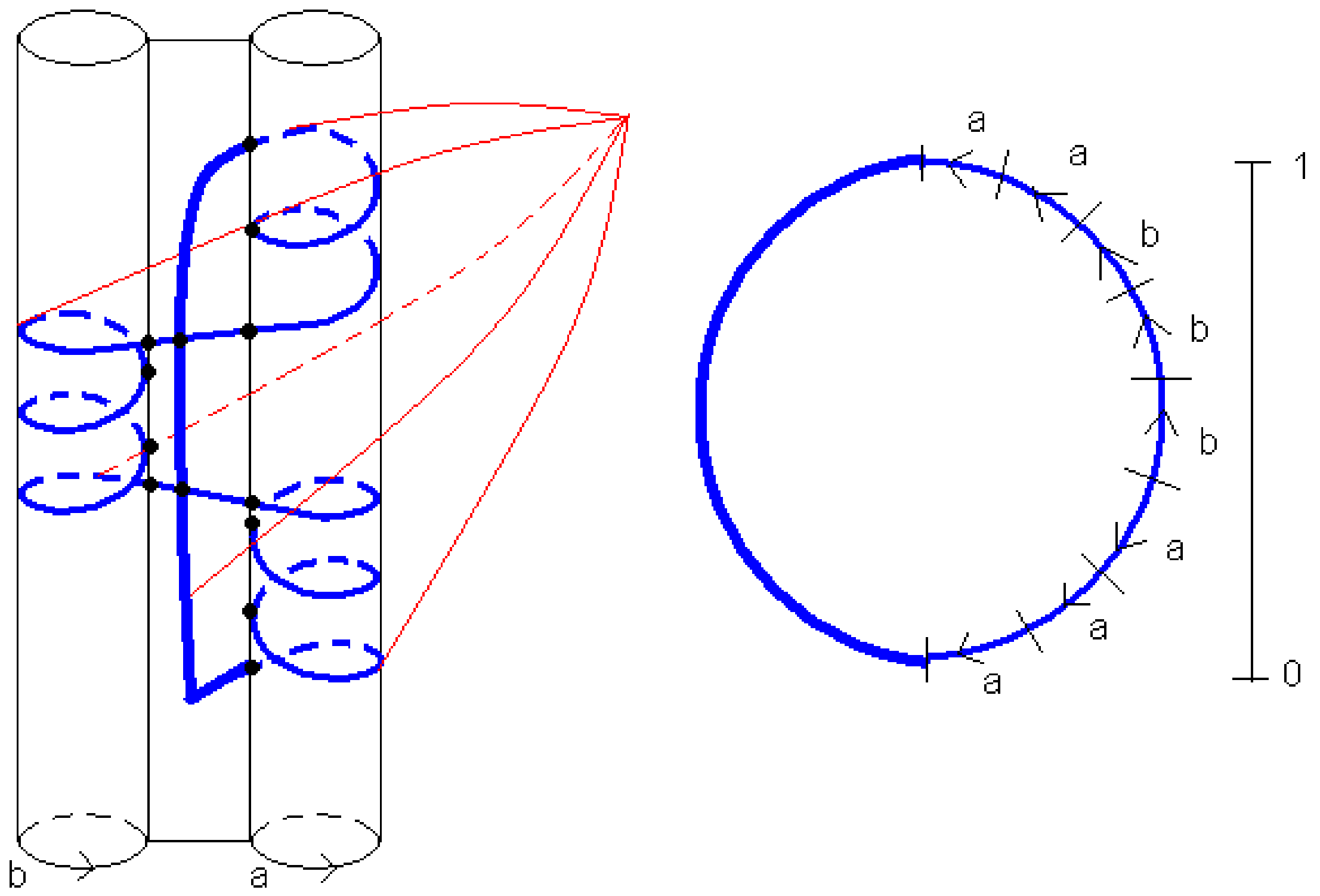

**Figure 18- Quinn model of a sliced 2-complex- attached relation induces vertices**

The intersection points of the attaching curve with the generator cylinder or with itself define the center of a local vertex model. The attaching curve of the relator defines a height function. Now we will describe the height function of the characterisic map corresponding to the 2-cell in more detail.

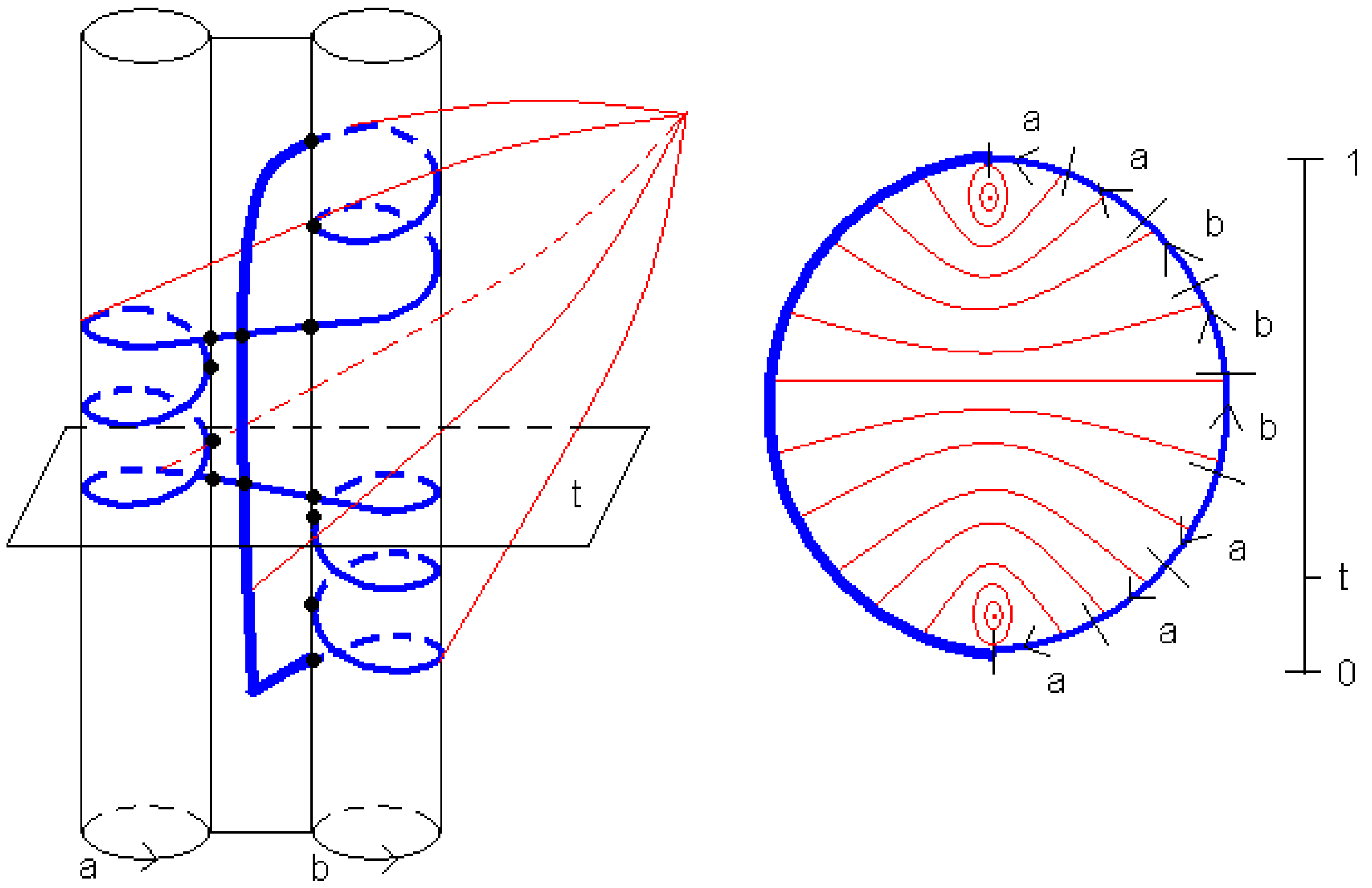

**Figure 19- Quinn model of a sliced 2-complex- height function of the characteristic map**

There is an arc in the 2-cell, called relation arc (start and end as a circle), which slides from minimum along the generators to maximum. This can be seen as a height function from 0 to 1, which determines the sequence of slices for the attached 2-cell. Combine that with the slices of the generator cylinders and the rectangle between them, we get a height function or equivalent to that the sequence of slices, one slice for each level t of the height function, for the whole 2-complex.

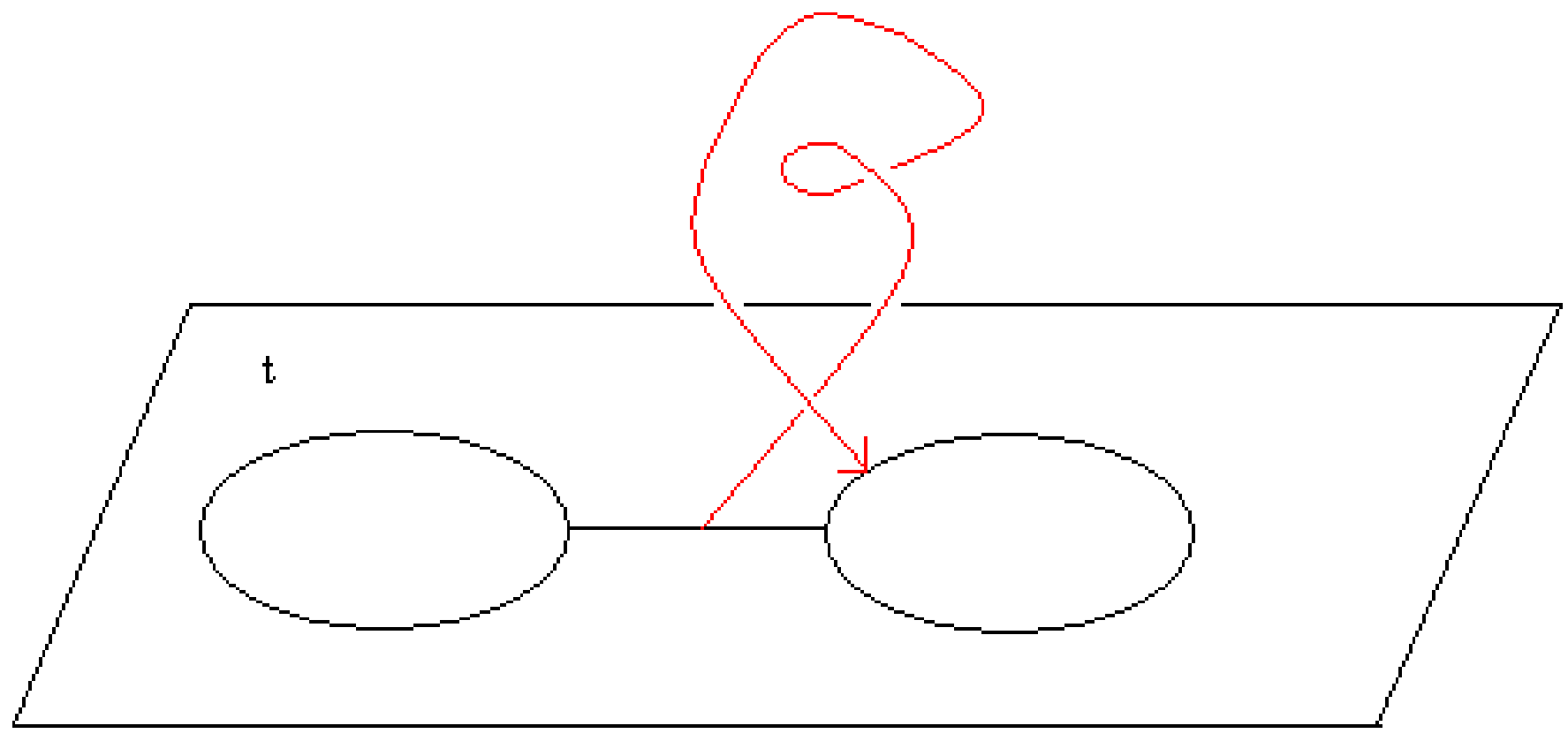

**Figure 20- Quinn model of a sliced 2-complex- a layer with level t- embedded version**

Note that the relation arc is an unknotted curve without any selfintersections. Since we are only looking for the topological difference of graphs, we can ignore the winding and simplify the relation arc:

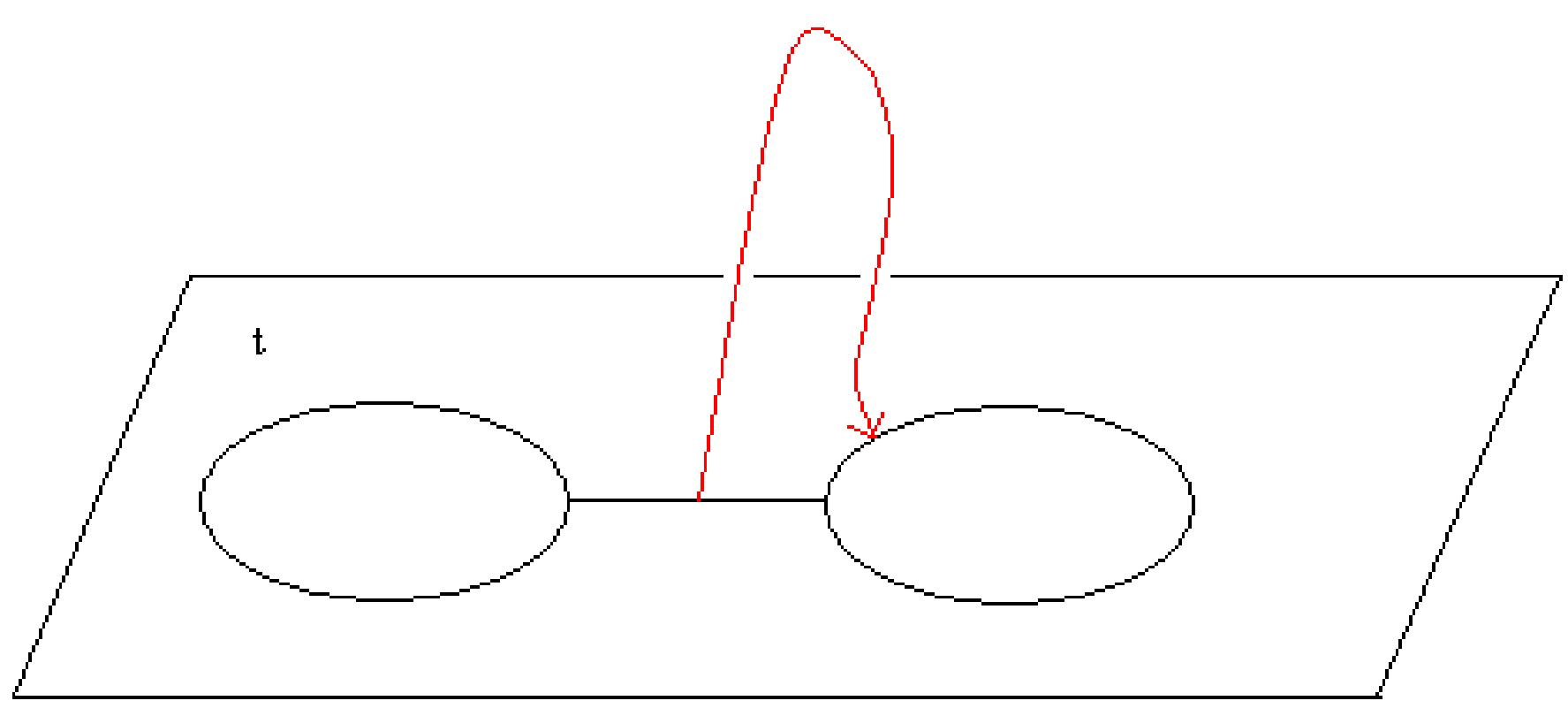

**Figure 21- Quinn model of a sliced 2-complex- a layer with level t- reduced version**

To describe and in the end compute the invariants of the 2-complex in the Quinn model, there are 3 essential parts:

a) entry the 2-cell at the minimum of the height function:

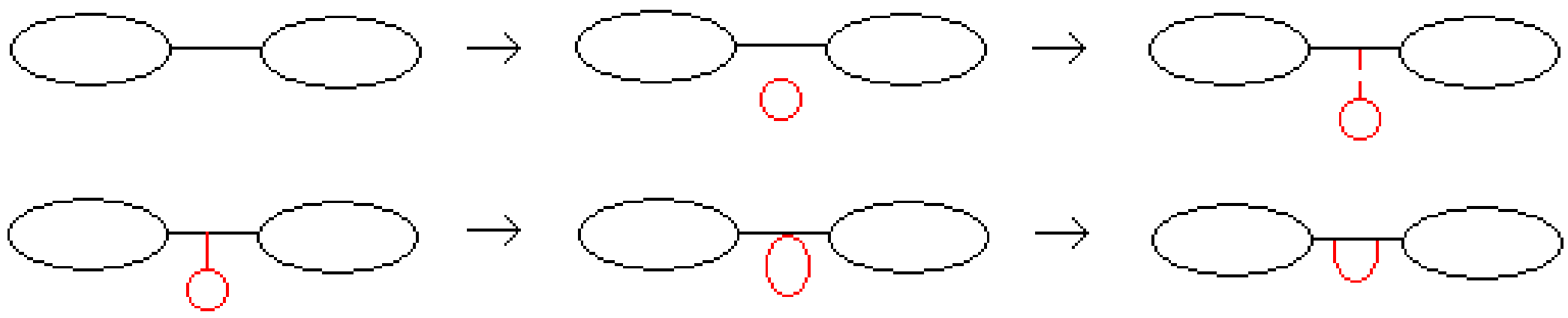

**Figure 22- Quinn model of a sliced 2-complex- entry at minimum**

b) slide around a generator:

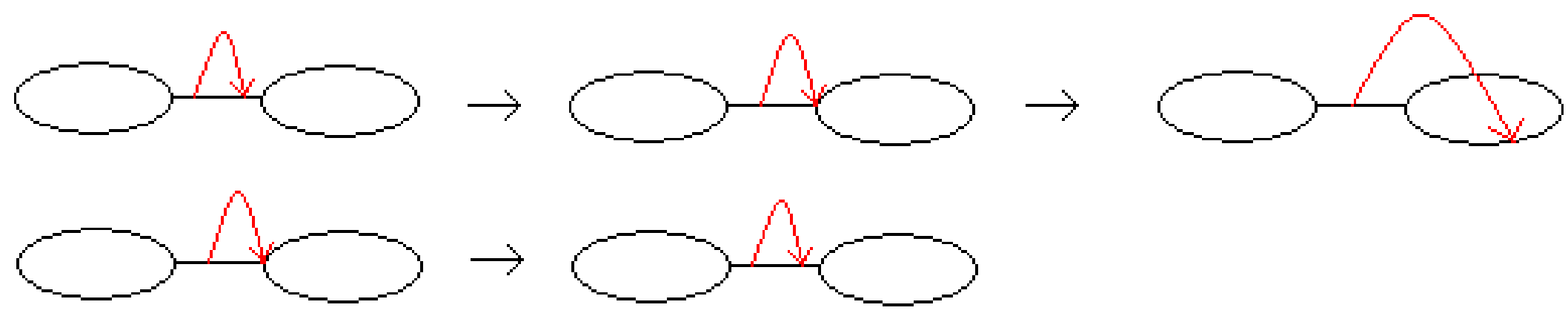

**Figure 23- Quinn model of a sliced 2-complex- slide around generator**

c) exit the 2-cell at the maximum of the height function:

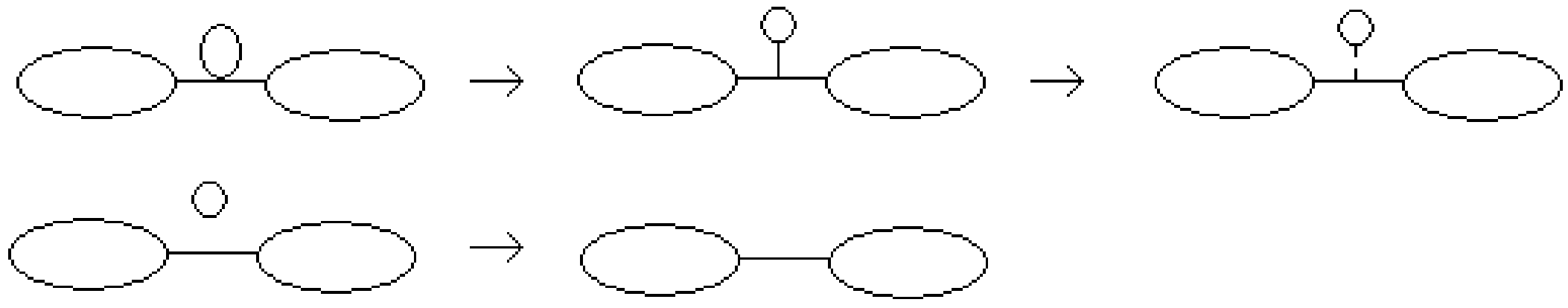

**Figure 24- Quinn model of a sliced 2-complex- exit at maximum**

This sequences of slices can be interpreted in an algebraic context. For example in the Topological Quantum Field Theory (TQFT), the main part is to develop the slide around a generator, which is called the circulator (compare **[Q2]** or **[Mül]**).

### *4.3 Sliced $T_3$ move - good and bad $T_3$ turn*

This Matveev move indicates a general problem that occurs by constructing A-C-invariants based on sliced 2-complexes:
Theorem 2 of P. Wright shows a bijective assignment between Q**-classes of finite group presentations and 3-deformation classes of compact, connected CW-complexes. A central idea in his proof is the association of Q-transformations with the homotopy of attaching maps for 2-cells. The homotopy can be performed on a 2-complex (not sliced) without any restrictions.
However the same homotopy considered on a sliced 2-complex has additional data: the slices. Hence it has to be verified, that for every step of the homotopy the slices of the attached 2-cell are compatibel with the slices on the remaining 2-skeleton of the 2-complex.
The $T_3$ move is an example, that we could run in trouble:
It is clear, that in general the sequence of slices inside and outside a deformed local model have to fit together. Consider the $T_3$ move, the turn of the attaching curve can be deformed to right or left.

turn to right:

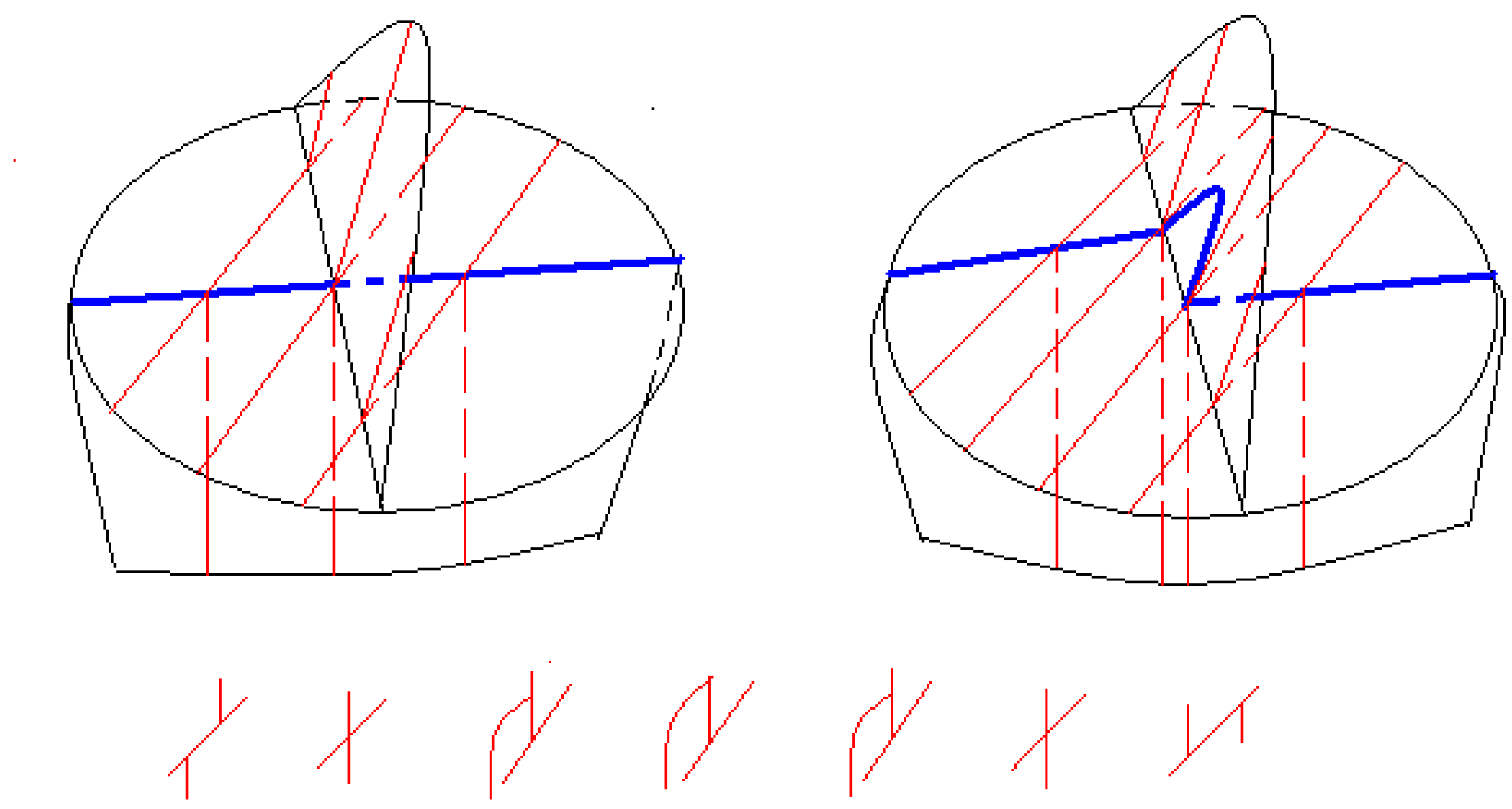

**Figure 25- sliced $T_3$ move- a good $T_3$ turn**

The figure shows that the sequence of slices before and after performing the move fit together, so we call this version of $T_3$ move in association with the given sequence a “good $T_3$ turn”.

turn to left:

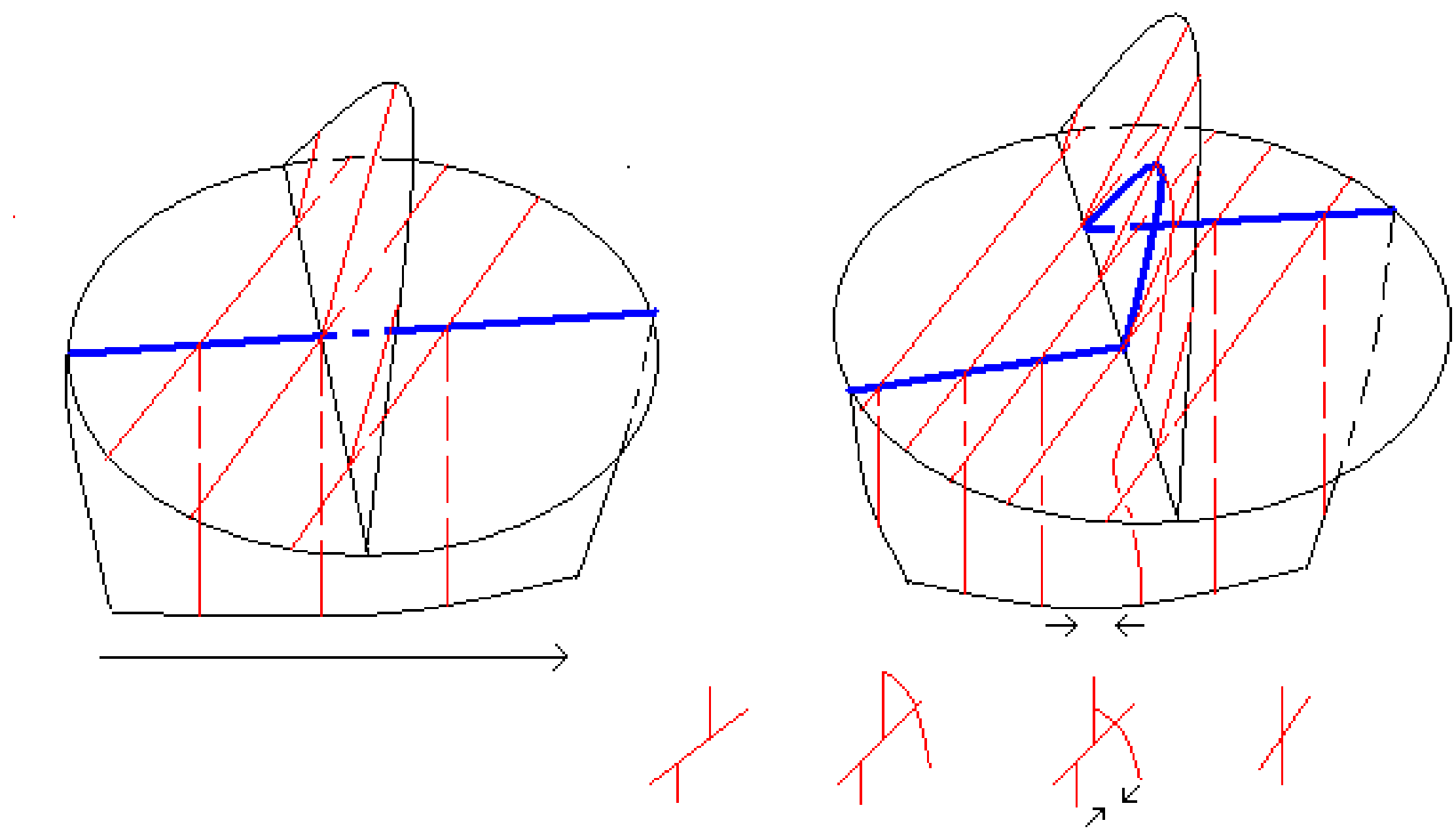

**Figure 26- sliced $T_3$ move- a bad $T_3$ turn**

The arrow on the left figure indicates an increasing height function on the boundary of the bottom component of the local vertex model. Assume, we slice the local vertex model after performing the $T_3$ move (see the right figure) similar to the left one. Then at the vertex the arcs which belong to the bottom component move together, as the two little arrows indicate there:

This would imply a height function on the boundary of the bottom component, restricted to a local neighbourhood of that point, which is connected by an line segment to the vertex, with increasing values from left and increasing values from right.

Hence the height function has been totally changed on the boundary of the bottom component: This would extend on all attached component and so on.

Furthermore it is not pemitted from the Quinn list of local transitions, that 2 disjoint arcs join together. That’s why the choice of slices in association with the direction of the turn the curve is performed in, is called a “bad $T_3$ turn”.

We find an appropriate sequence, if we use the solution for the turn to right:

We change the sequence of slices near the vertex in the local model **before** performing the $T_3$ move:

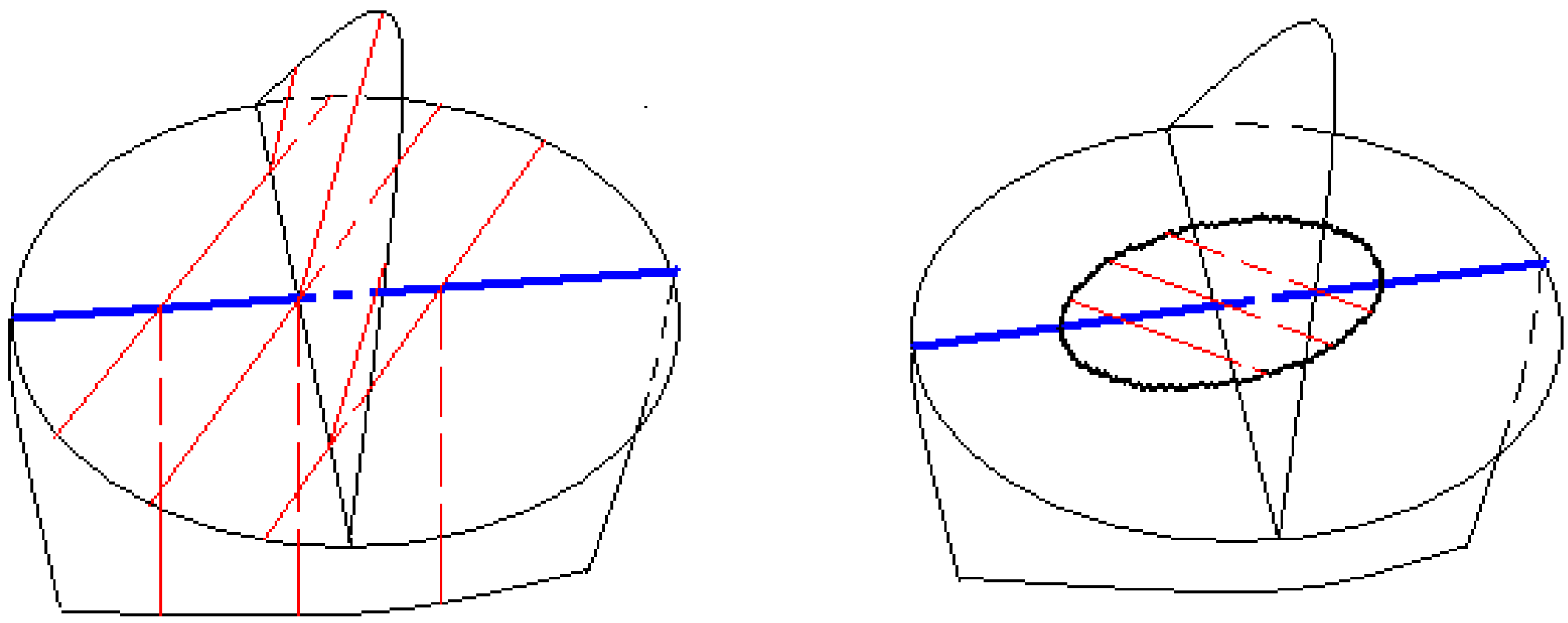

**Figure 27- sliced $T_3$ move- solve bad $T_3$ turn - step 1**

If we consider the sequence of slices in the neighbourhood (indicated by a circle) after performing the $T_3$ move, we see that this corresponds to the case for the “good $T_3$ turn” (the combination of the slices and the direction of the turn is similar to the first case).
Extend the sequence on the whole model as pointed out in the next figure. The slice has to be connected to the endpoints on the boundary and the arrows show that the height function increases as required.

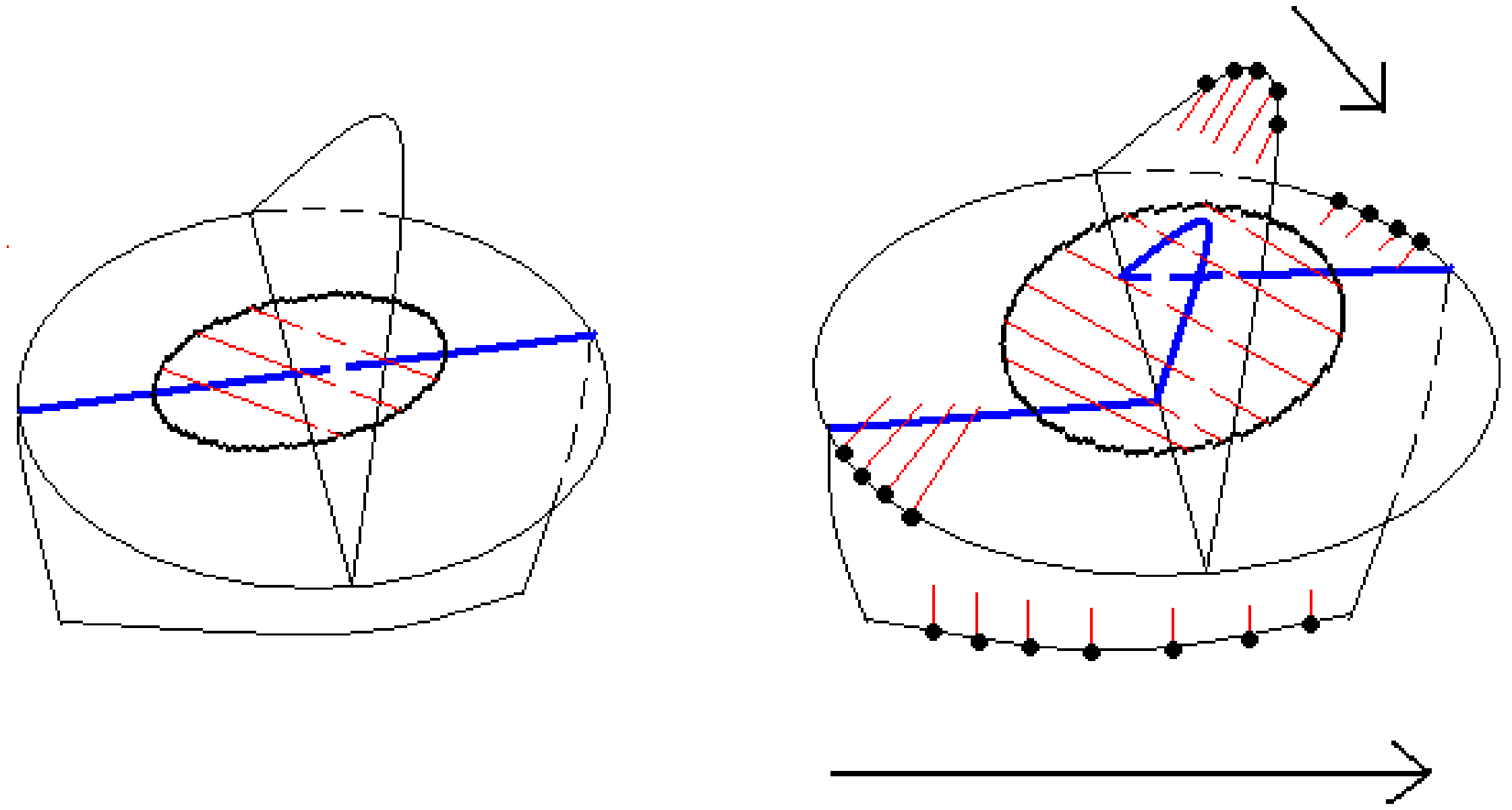

**Figure 28- sliced $T_3$ move- solve bad $T_3$ turn - step 2**

The extension can be realized as follows:

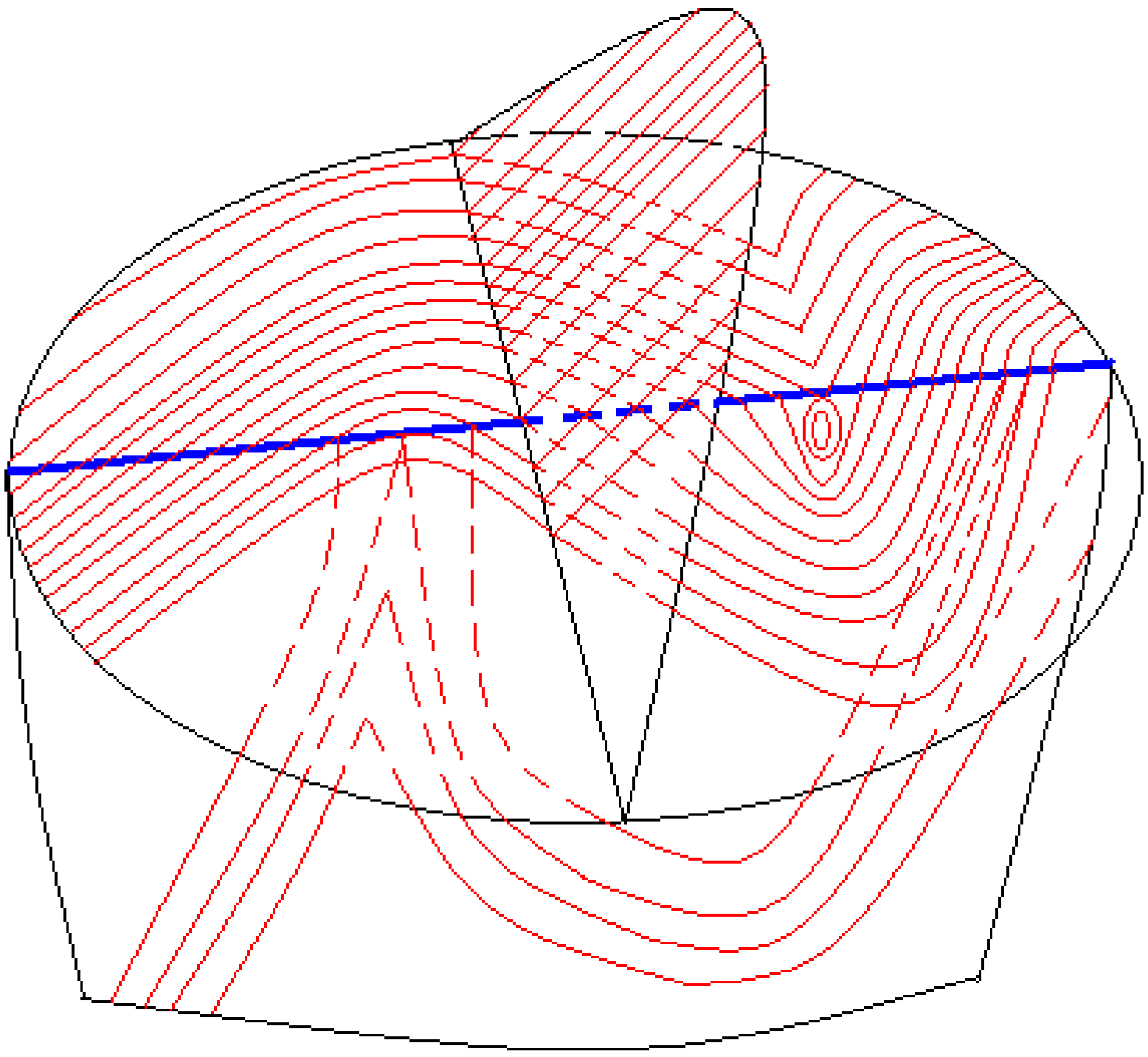

**Figure 29- sliced $T_3$ move - solve bad $T_3$ turn- modified slices of step 1**

We have to find a convenient sequence of slices, however the height function on the slices on the boundary is not exactly the same as before, but clearly the “rhythm” is kept, for example an increasing or decreasing height function. It is possible to change the height function on the slices such that it fits together with the height function on the attached components.
Furhermore it is fundamental to construct a connecting homotopy from the sequence of line segment slices to that in the figure above. It corresponds to a birth and death of a cancellation pair of saddlepoints, more precisely these are births and deaths of maxima and minima in the components or in a union of these. We use this construction as a standard procedure to introduce and cancel pairs of saddlepoints. The homotopy beeing constructed is presented in the next pictures, where we also conserve the “rhythm” of the height function.

We illustrate the homotopy not for the former, but for a similar case:

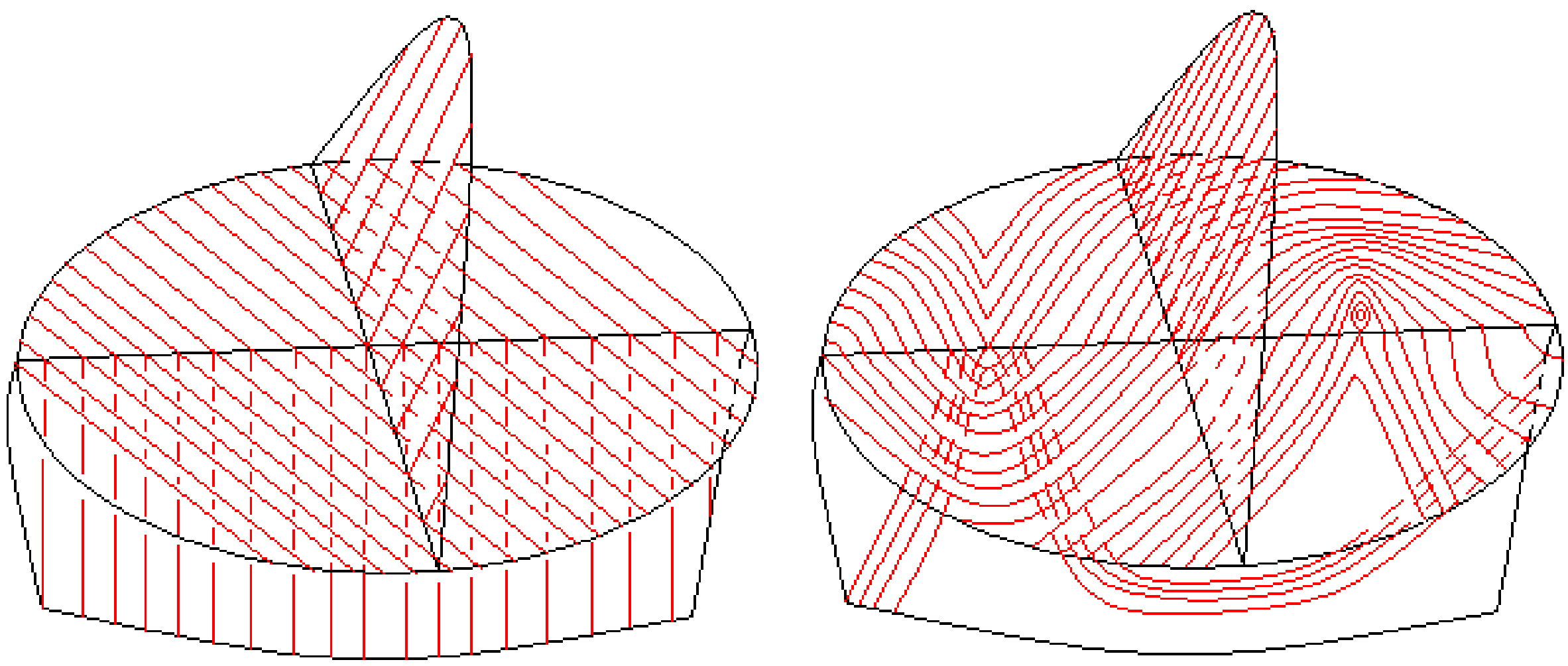

**Figure 30- sliced $T_3$ move- homotopy from straight to modified slices- the case**

We abbreviate:
The attaching curve of the bottom component is called the bottom line.
The attaching curve of the top component is called the top line.

Our homotopy is a composition of the following steps:

1. homotopy (change the straight slice to a waved slice on the bottom line):

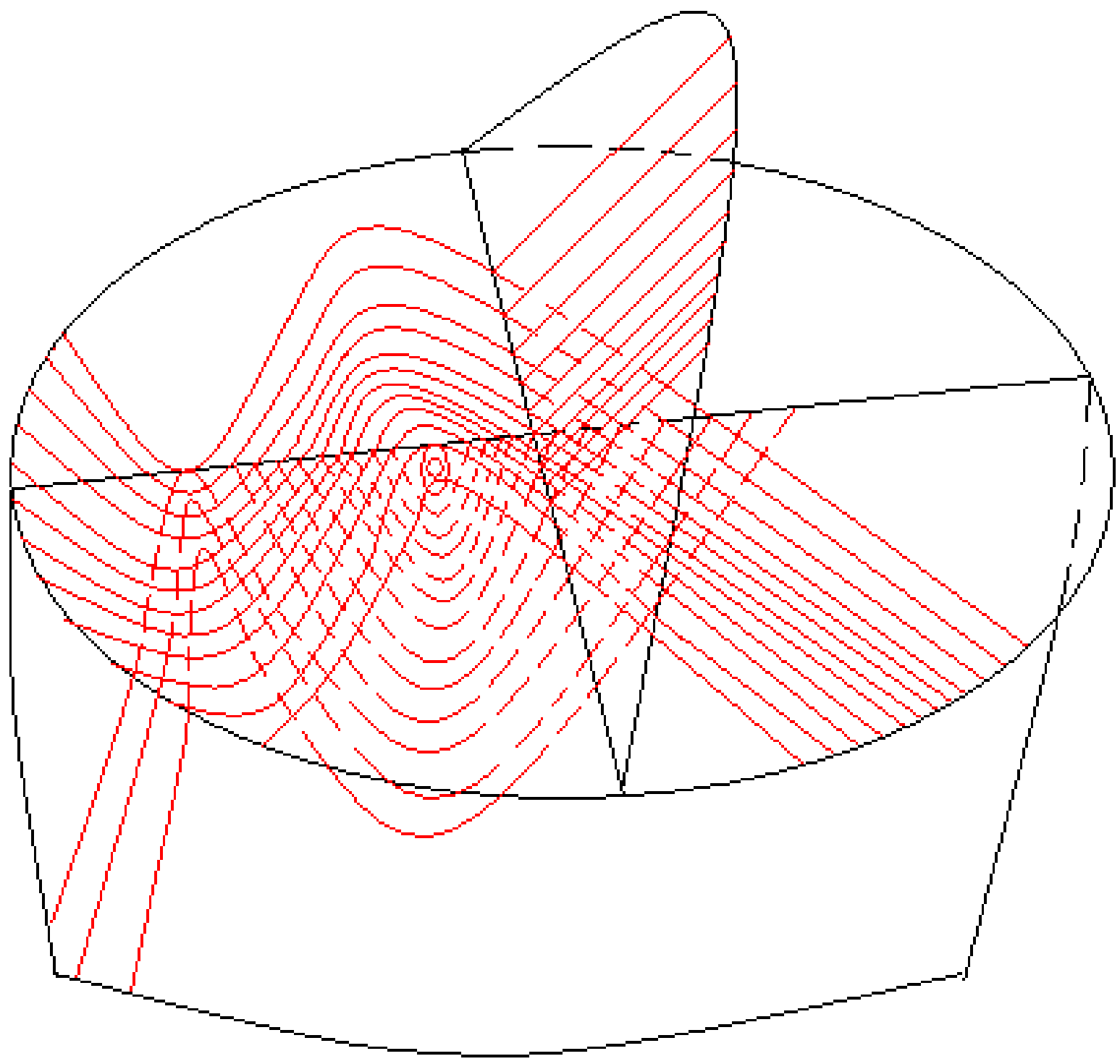

**Figure 31- sliced $T_3$ move- homotopy from straight to modified slices- 1**

2. homotopy (change the straight slice to a waved slice on the top line):

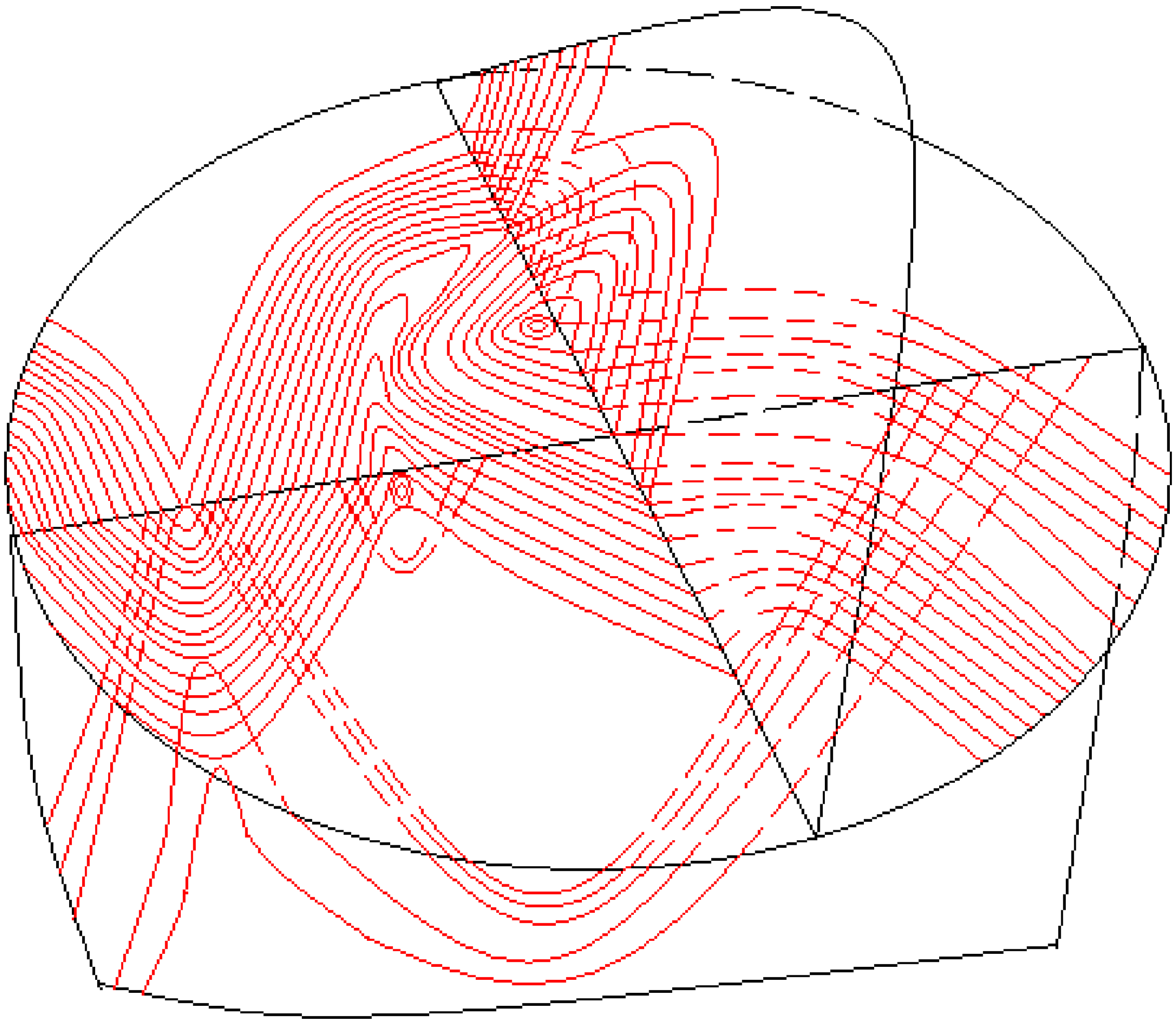

**Figure 32- sliced $T_3$ move- homotopy from straight to modified slices- 2**

3. homotopy (switch at the second saddlepoint on the top line, it corresponds to a local change of slices at the saddlepoint):

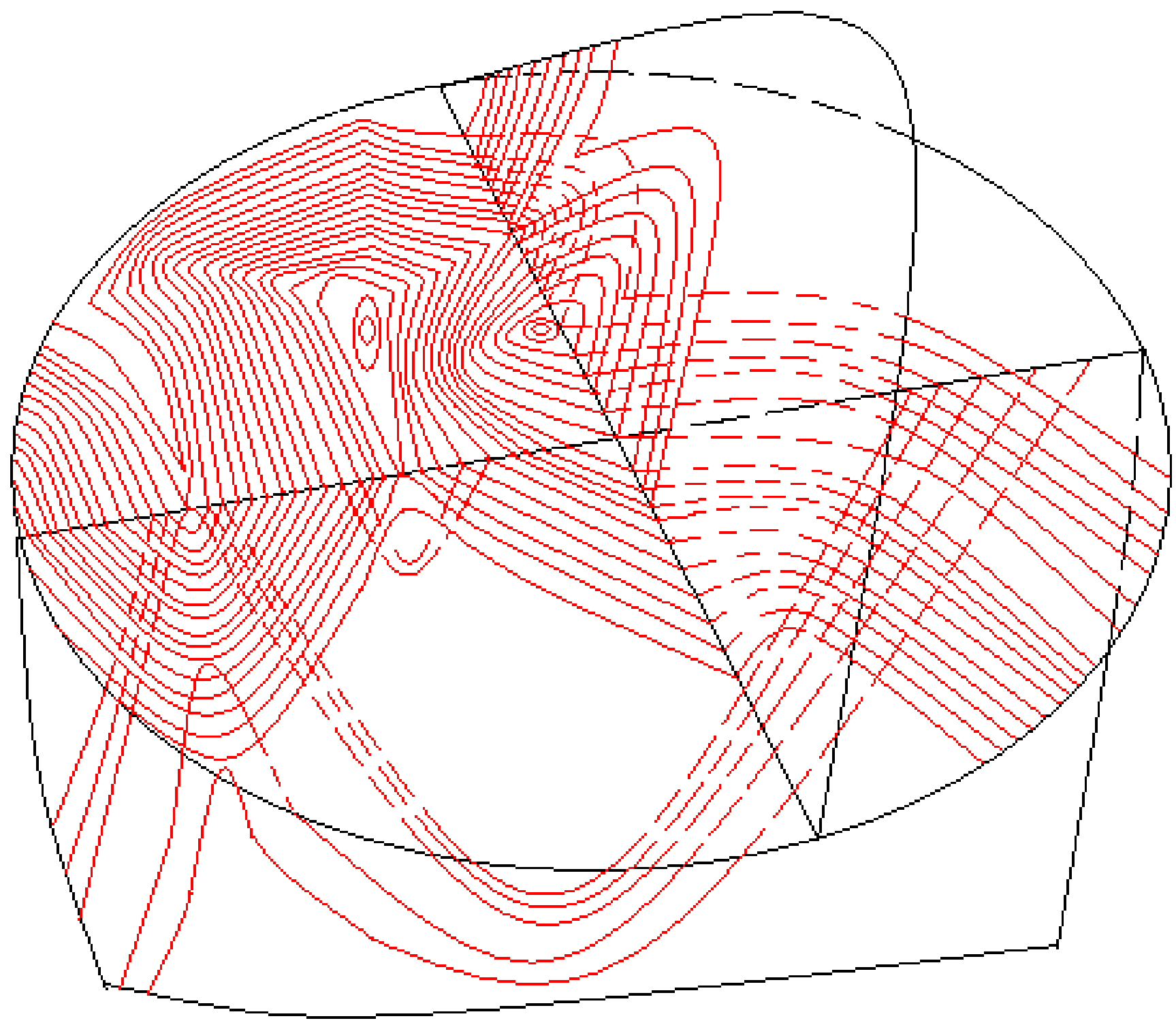

**Figure 33- sliced $T_3$ move- homotopy from straight to modified slices- 3**

4. homotopy:

We describe the next step, it sounds a little bit mystic:
Move the second saddlepoint on the bottom line and also the second saddlepoint on the top line to the central vertex. When they meet there, the saddlepoint of the bottom line should be dominant, i.e. the circle which belongs to the saddlepoint on the top line shrinks to a point (the second figure indicates the transformation near the vertex). However the slice of this saddlepoint , restricted on the top component is still alive at the vertex and fits together with the moved saddlepoint from the bottom line:

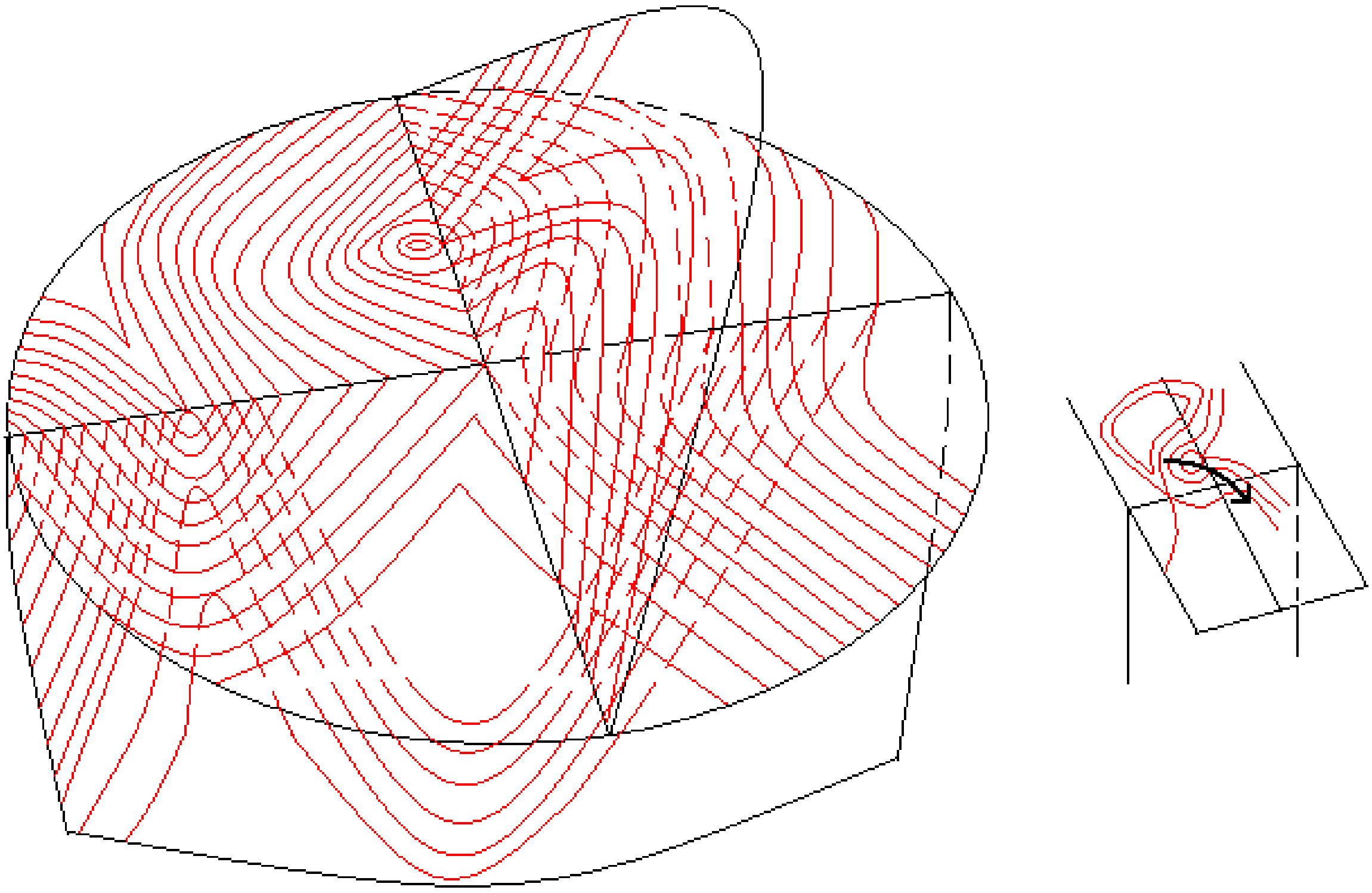

**Figure 34- sliced $T_3$ move- homotopy from straight to modified slices- 4**

5. homotopy:

We move the saddlepoint on the vertex along the bottom line into the backside components. This is the reversed action of the former two steps in the backside components. Hence the result is then similar to the figure before step 3:

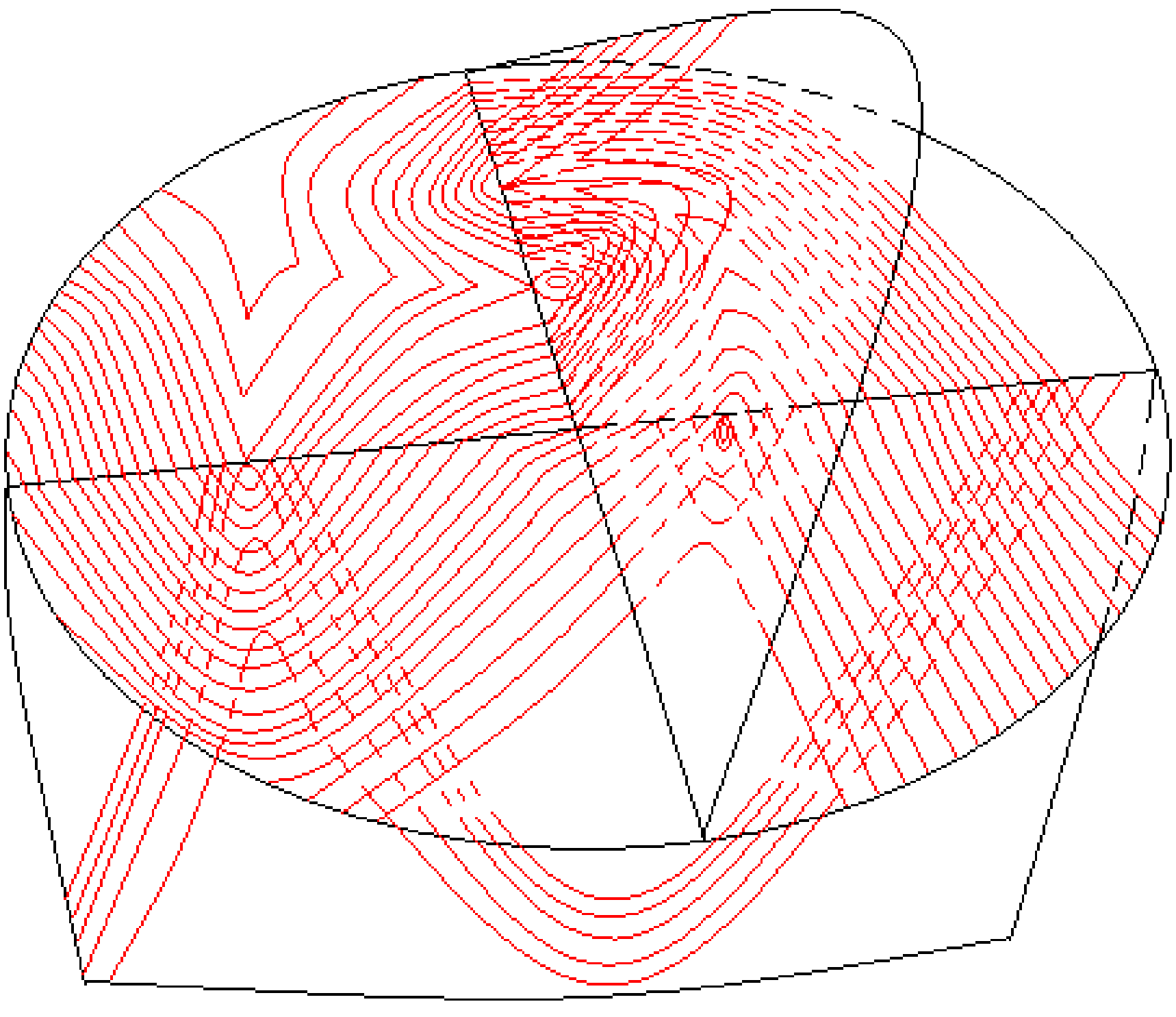

**Figure 35- sliced $T_3$ move- homotopy from straight to modified slices- 5**

6. homotopy:

We change the waved slice in the backside components to a straight slice:

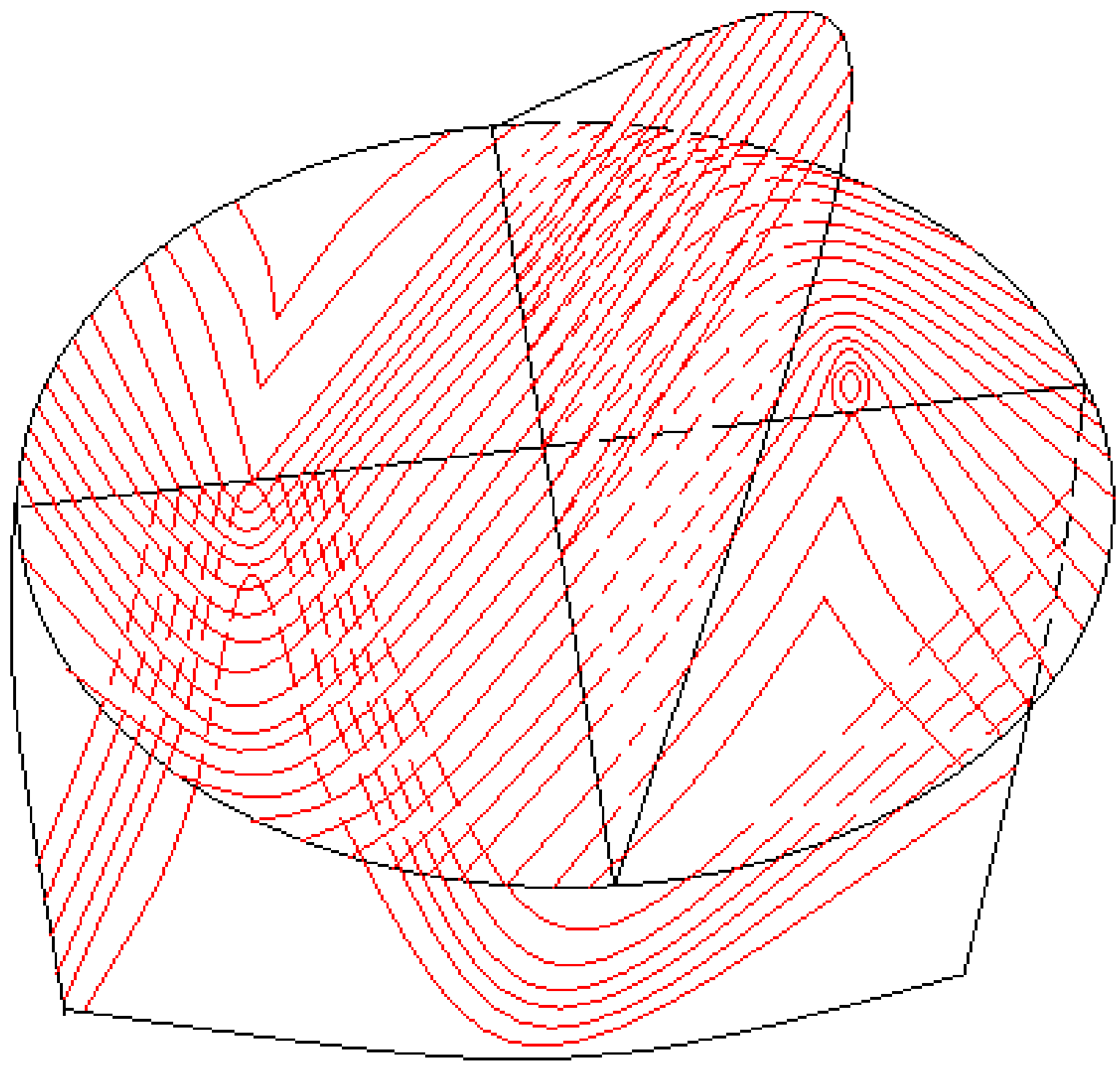

**Figure 36- sliced $T_3$ move- homotopy from straight to modified slices- 6 (end)**

All these figures include a basic sequence, which leads to a new sequence among the local transitions in the Quinn list. The problem of the $T_3$ move shows, that this sequence (shown by the next figure) is required. It can be modified to a new topological relation, if we attach the transition for passing the local vertex model i.e. we extend to the sequence by a transition to get the same start- and endgraph.
We see, that this topological relation could not be a consequence of the known topological relations “building a bubble” and “expansion with a disc”.
Nevertheless, if we view the local transitions in a suitable algebraic context (for example TQFT), the corresponding algebraic relation could be a consequence of the known algebraic relations “building a bubble” and “expansion with a disc” or other algebraic relations.
Therefore it is still an open question, that for example in TQFT the new topological relation represents an (new) algebraic relation in the Quinn list of relations.
In chapter 8 we compute in a special TQFT example, that this new topological relation leads to an algebraic relation.

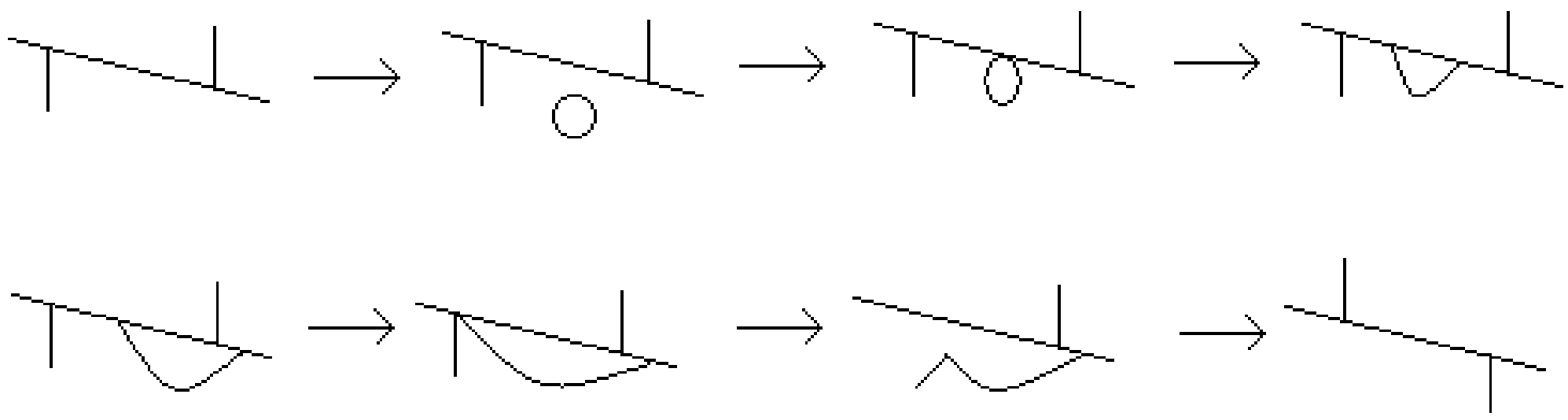

**Figure 37- sliced $T_3$ move- the basic sequence of slices**

For every step of the homotopy we present the sequence of slices, where the basic sequence comes from the first homotopy:

2. homotopy:

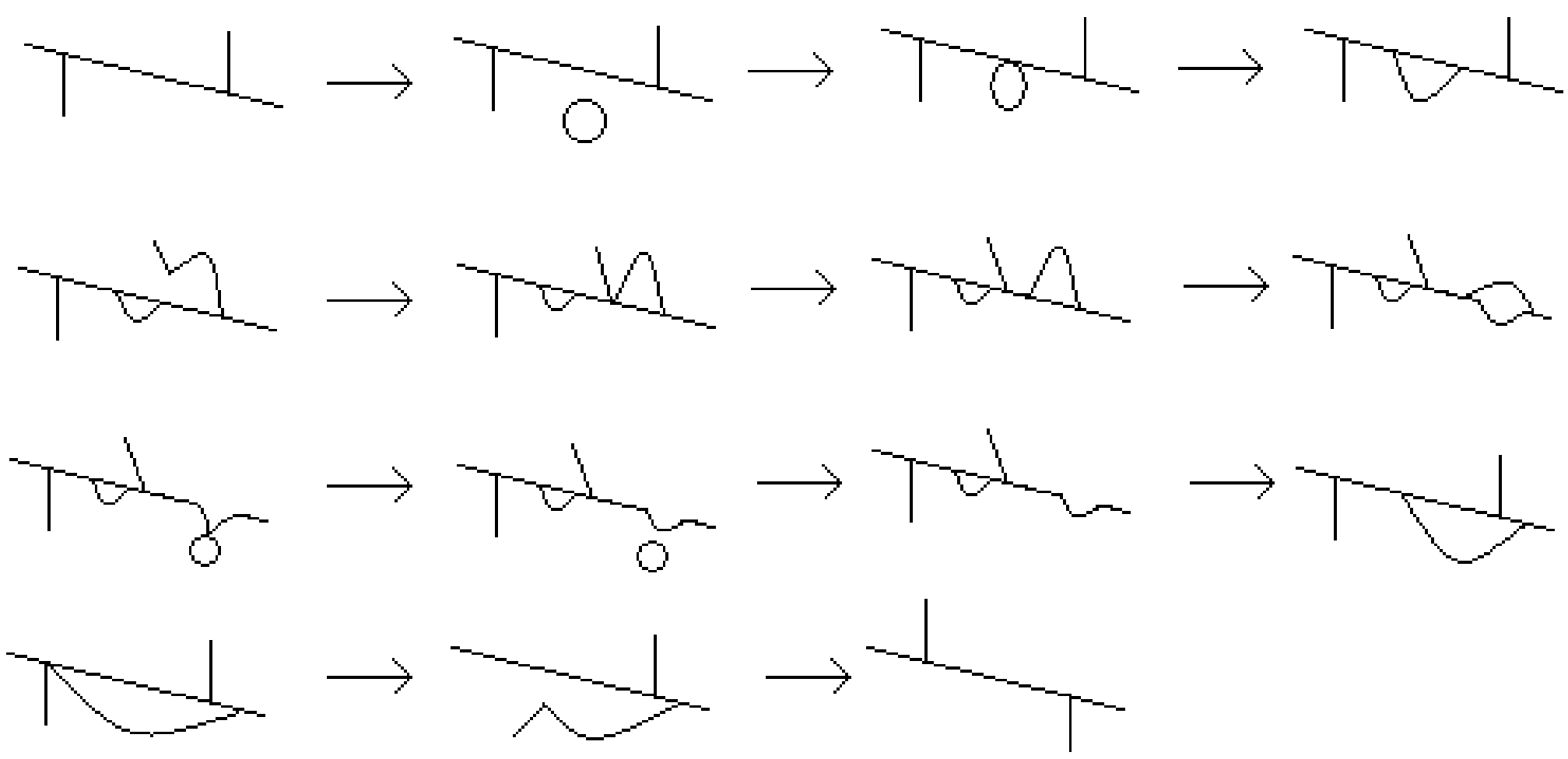

**Figure 38- sliced $T_3$ move- homotopy step 2- the slices**

3. homotopy:

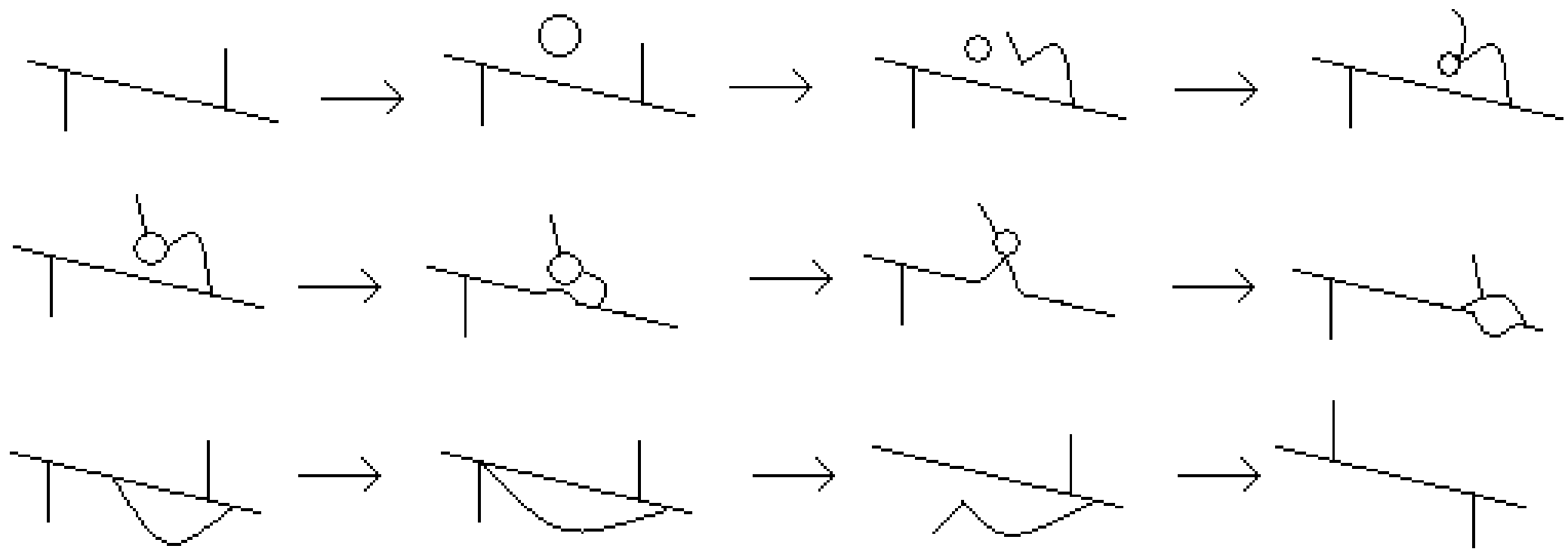

**Figure 39- sliced $T_3$ move- homotopy step 3- the slices**

4. homotopy:

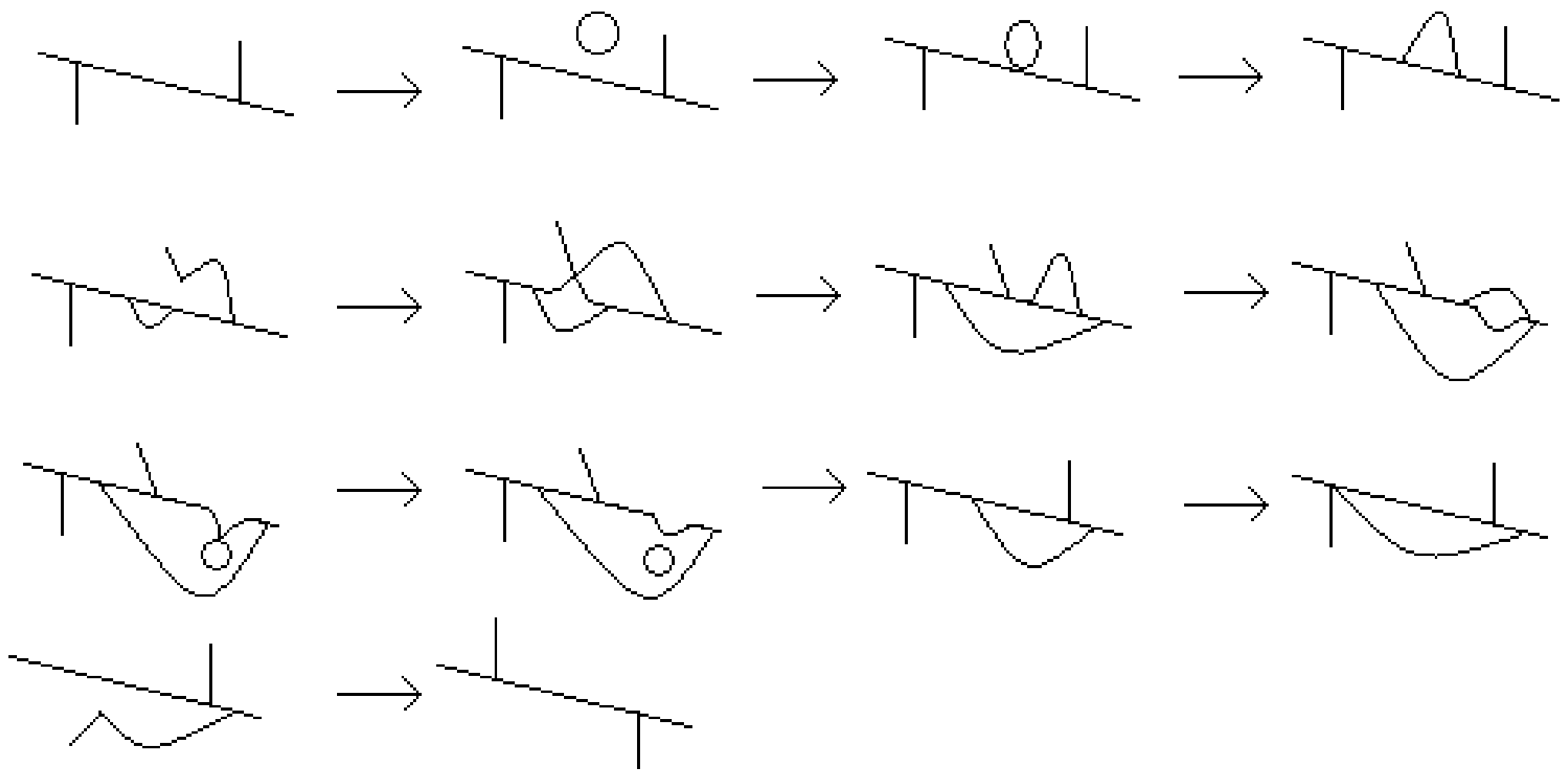

**Figure 40- sliced $T_3$ move- homotopy step 4- the slices**

5. homotopy:

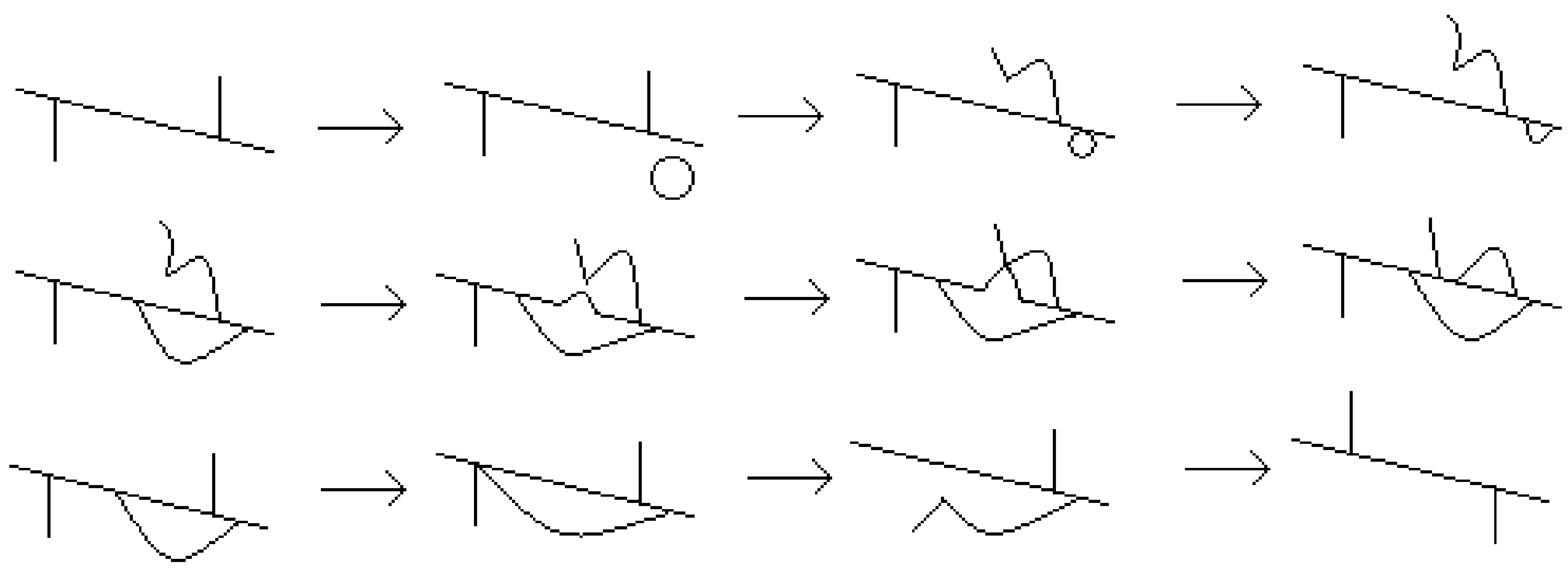

**Figure 41- sliced $T_3$ move- homotopy step 5- the slices**

6. homotopy:

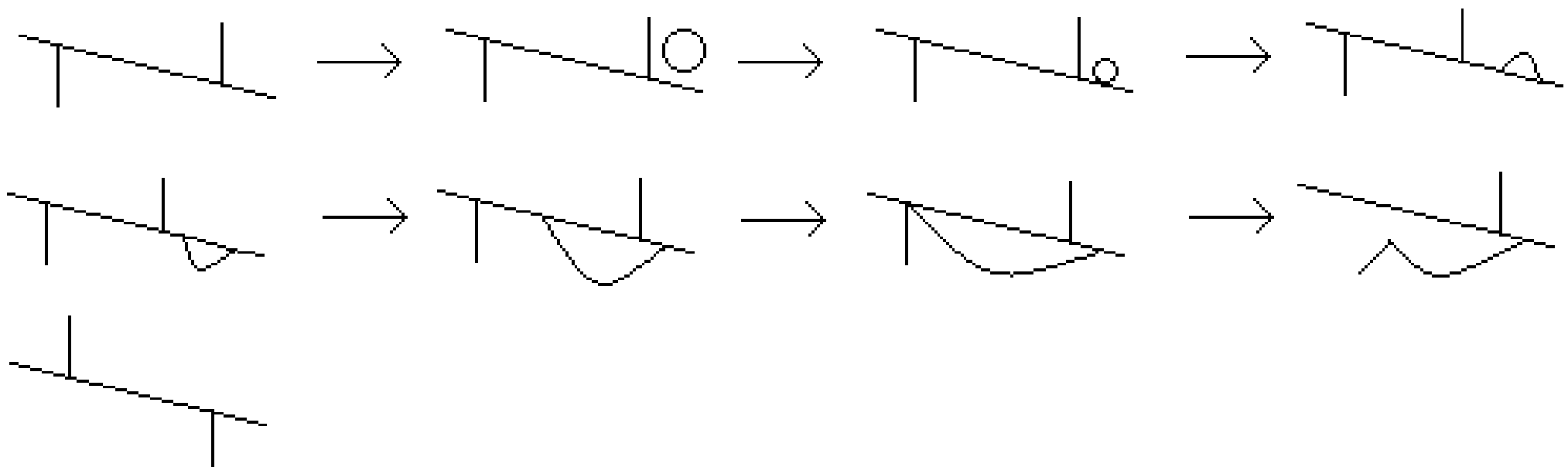

**Figure 42- sliced $T_3$ move- homotopy step 6- the slices**

## *4.4 The sliced leftside loop*

This loop is another example (in addition to the “bad $T_3$ turn”) which shows, that there could not exist a general fixed slice of the 2-complex in the Quinn model, which holds for all changes of the attaching curve. We will see, that depending on the change of that curve, the sequence of slices has been modified. A central part of this work is the subsolution, where we cancel a pair of loops. We refer to chapter 5 for details. Assume we have a sequence of slices for the vertex model (add little flanges to get it), such that the sequence corresponds to the entry of the 2-cell. Further assume, there is an attaching curve that passes the vertex, for example this happens by multiplication of 2 relations.

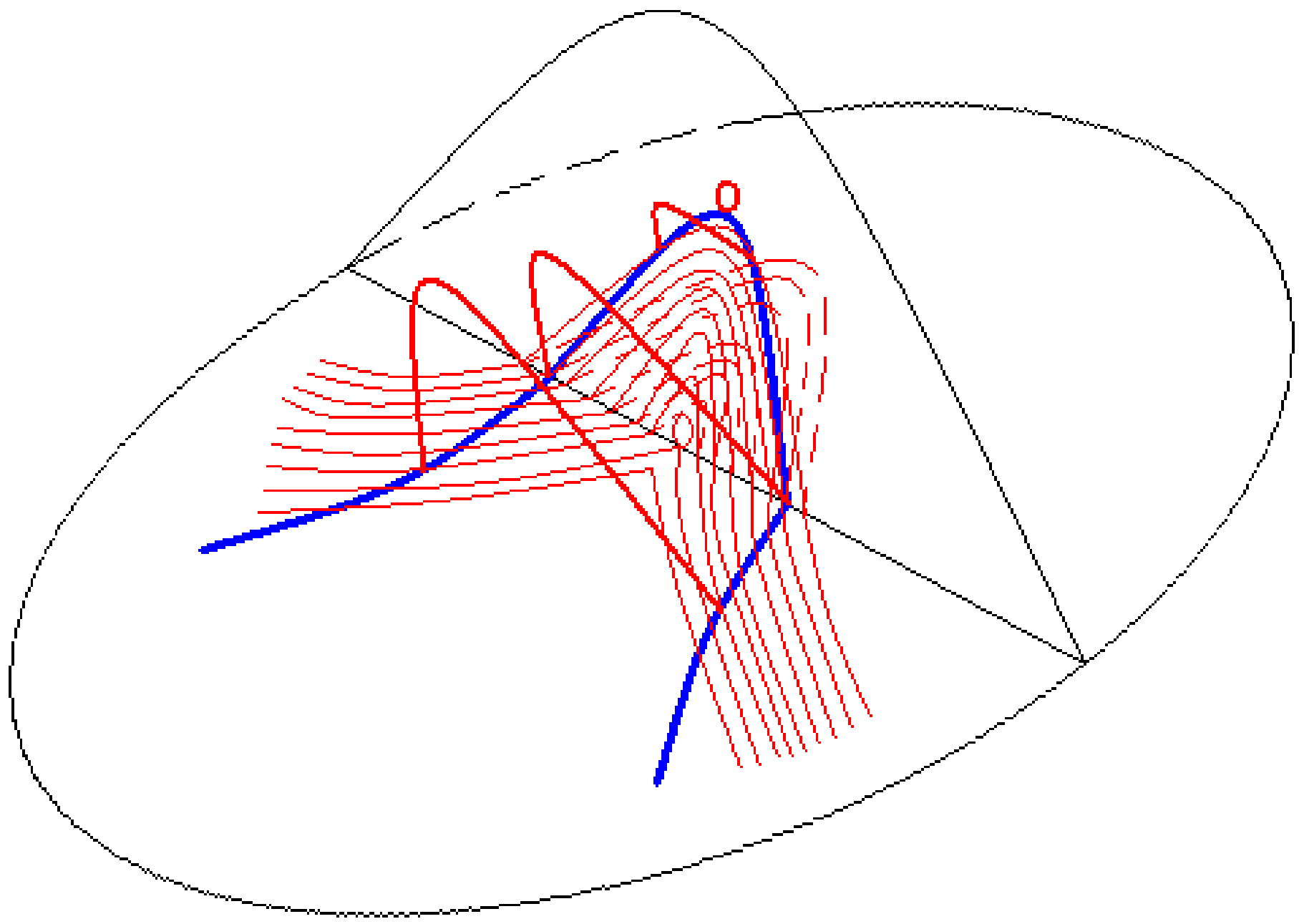

**Figure 43- the sliced leftside loop- start figure**

Here the thicker red lines represent the slices of the 2-cell of the attaching curve. We reduce the picture to the essential parts, where we can easily develop the extension to the rest:

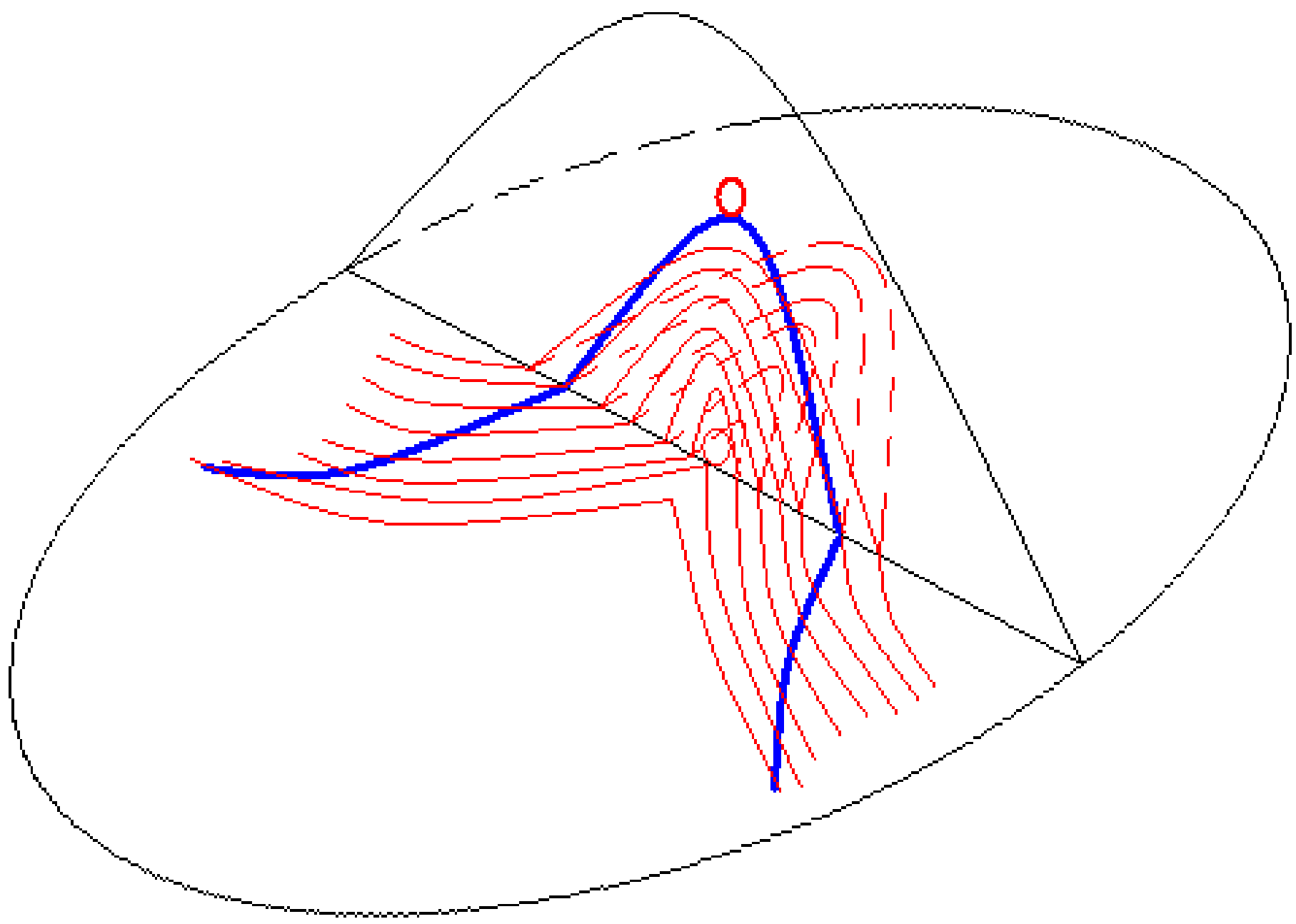

**Figure 44- the sliced leftside loop- reduced startfigure**

First, without including the sequence of slices, the process to build a leftside loop is given by:

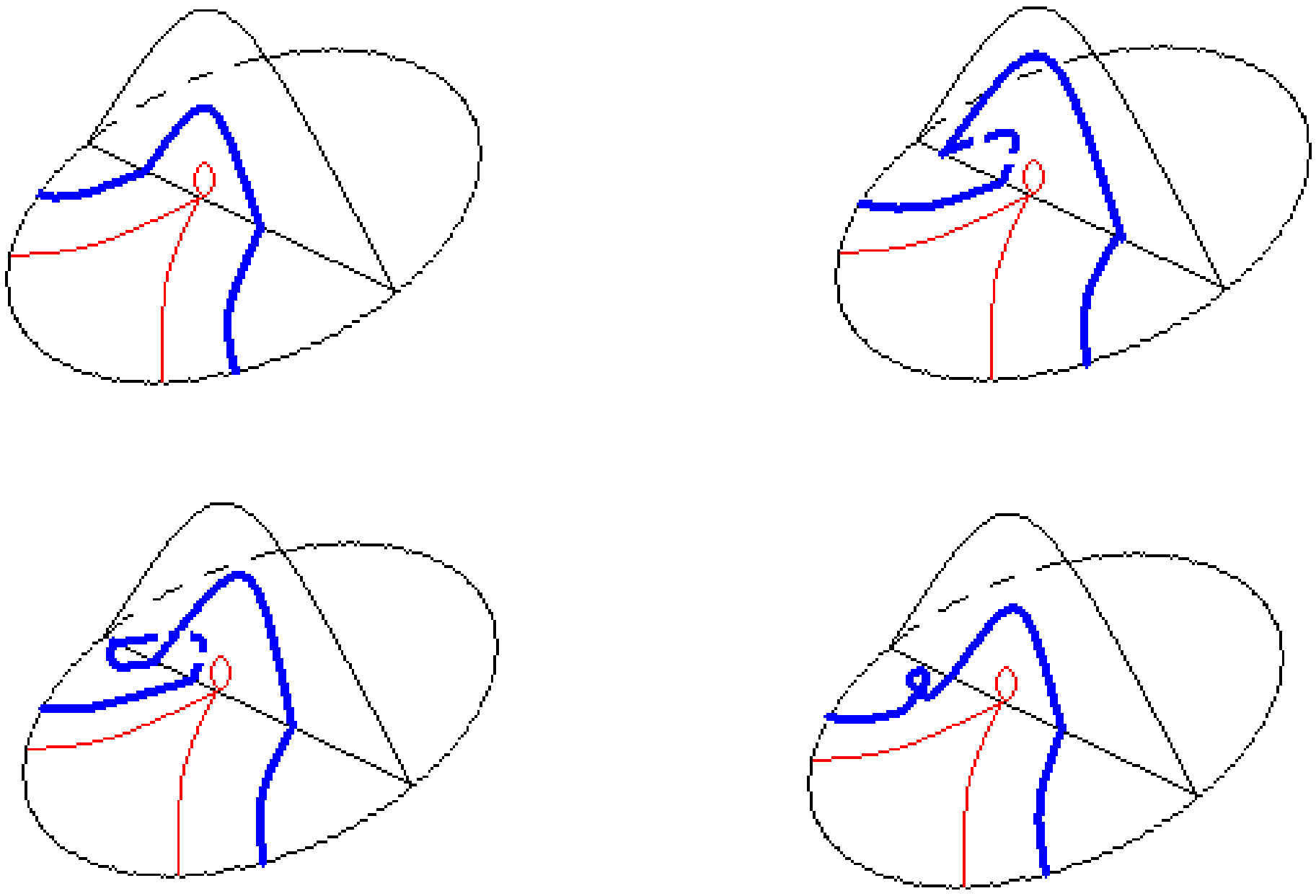

**Figure 45- the sliced leftside loop- preview of leftside loop construction**

It is a compostition of $T_3$ and $T_2$ moves, hence we have to develop the sequence of slices for each step.

Step 1:
We peform the first step and see, that there are no modifcations required. The second local figure shows that we get to a "good $T_3$ turn".

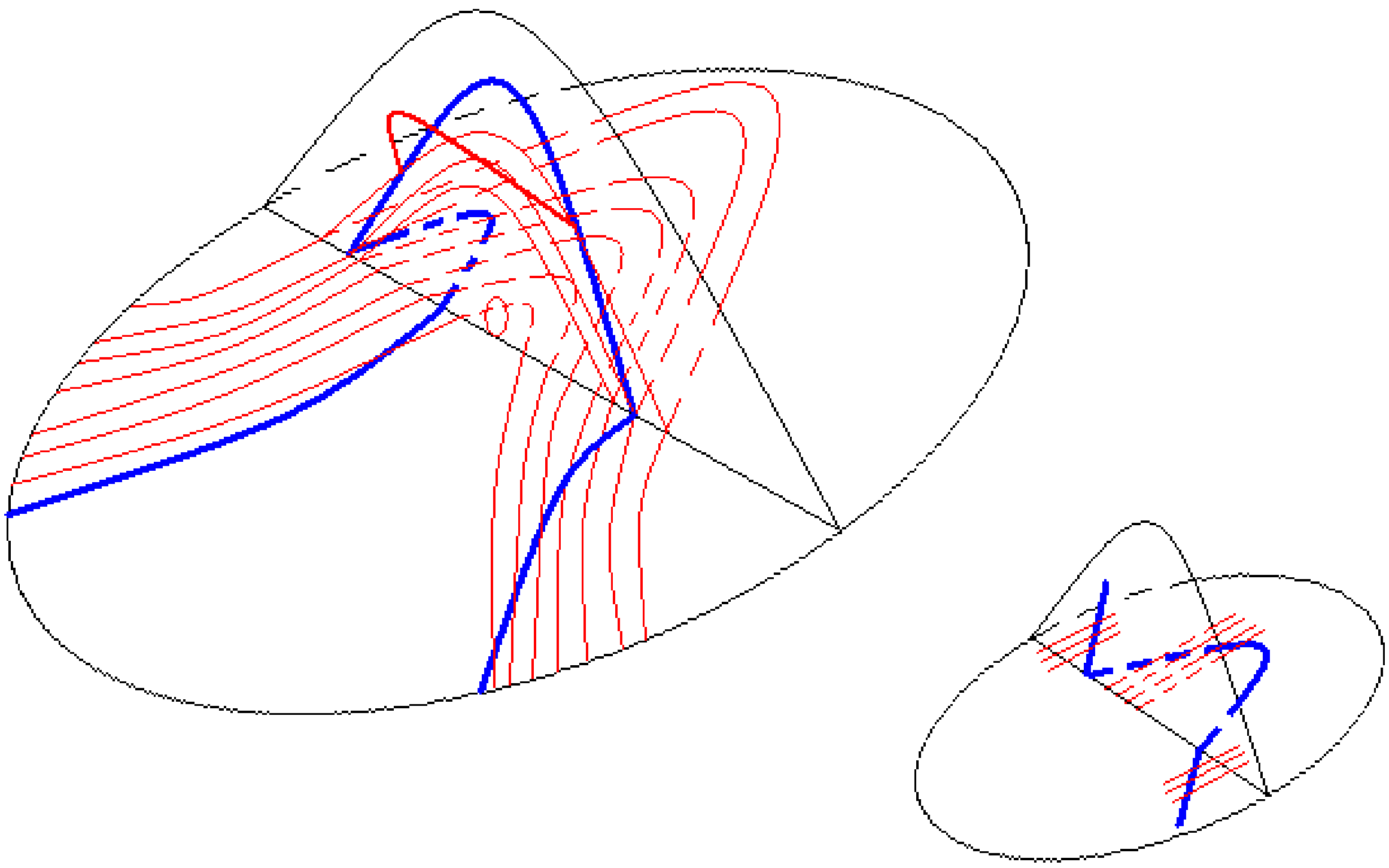

**Figure 46- the sliced leftside loop- 1**

Step 2:
If we perform the next step without modifying the sequence of slices, then we would get a local picture which represents the "bad $T_3$ turn", therefore we have to modify the slices before executing step 2.

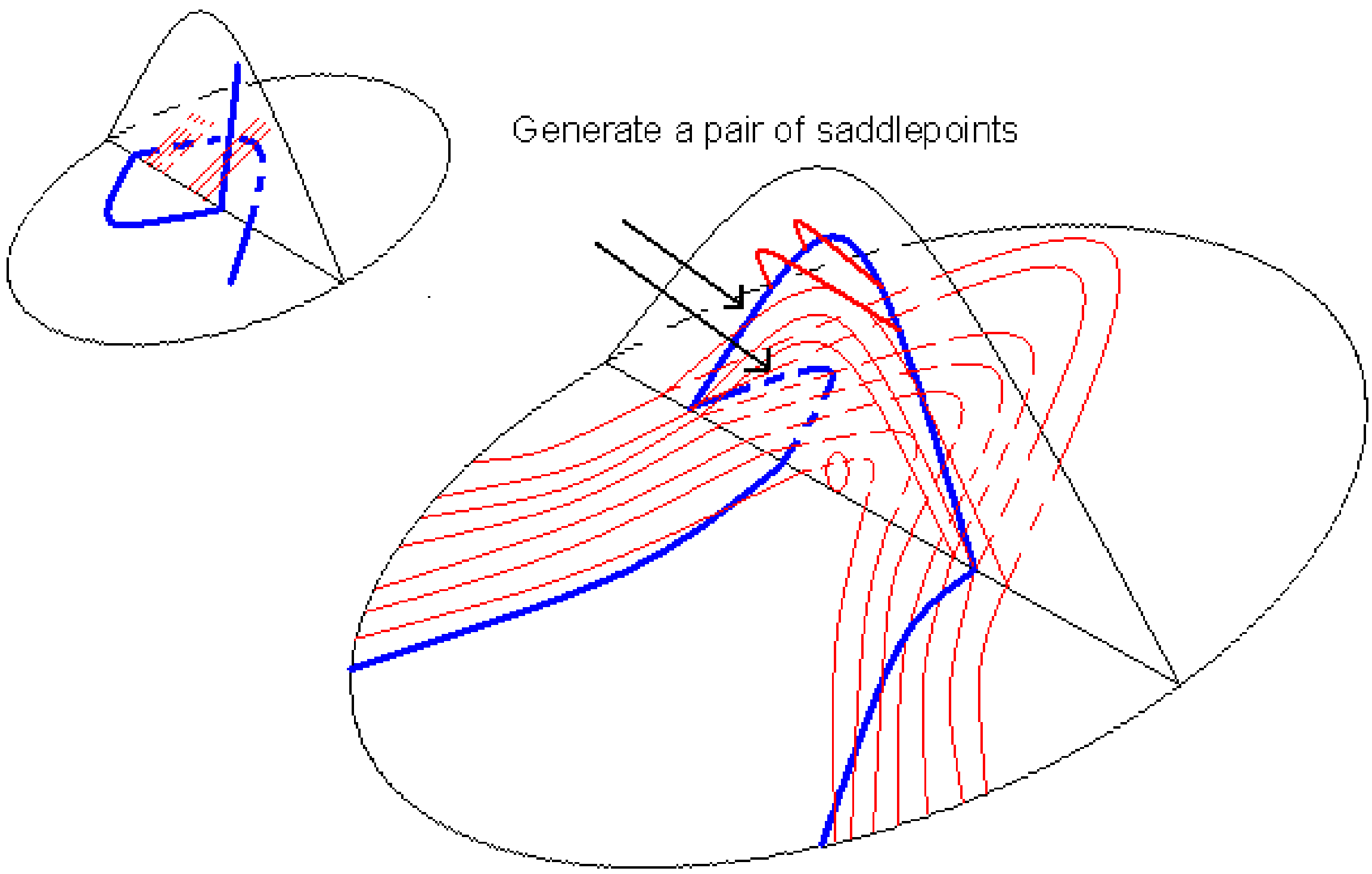

**Figure 47- the sliced leftside loop- 2**

We introduce a pair of saddlepoints as in our standard solution for the “bad $T_3$ turn”. The second figure indicates how to embed them at the component of the attached 2-cell.

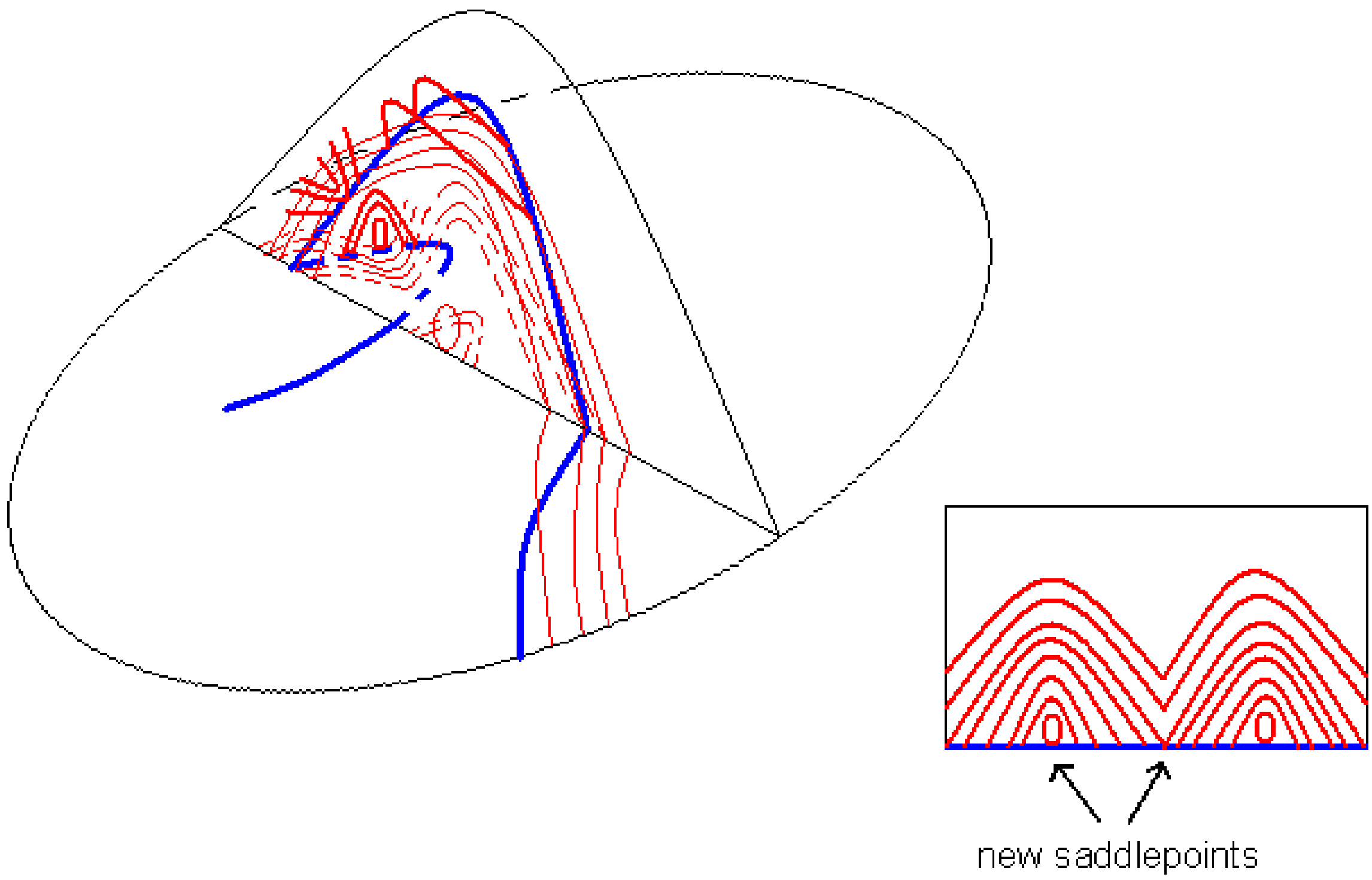

**Figure 48- the sliced leftside loop- 3**

After this preparation we can execute step 2.

However the new saddlepoints help us more than expected. When we perform the next step as indicated from the arrow, we get a loop.
For a loop, we need at least one saddlepoint and for pushing the arc across the vertex, we need a saddlepoint too:

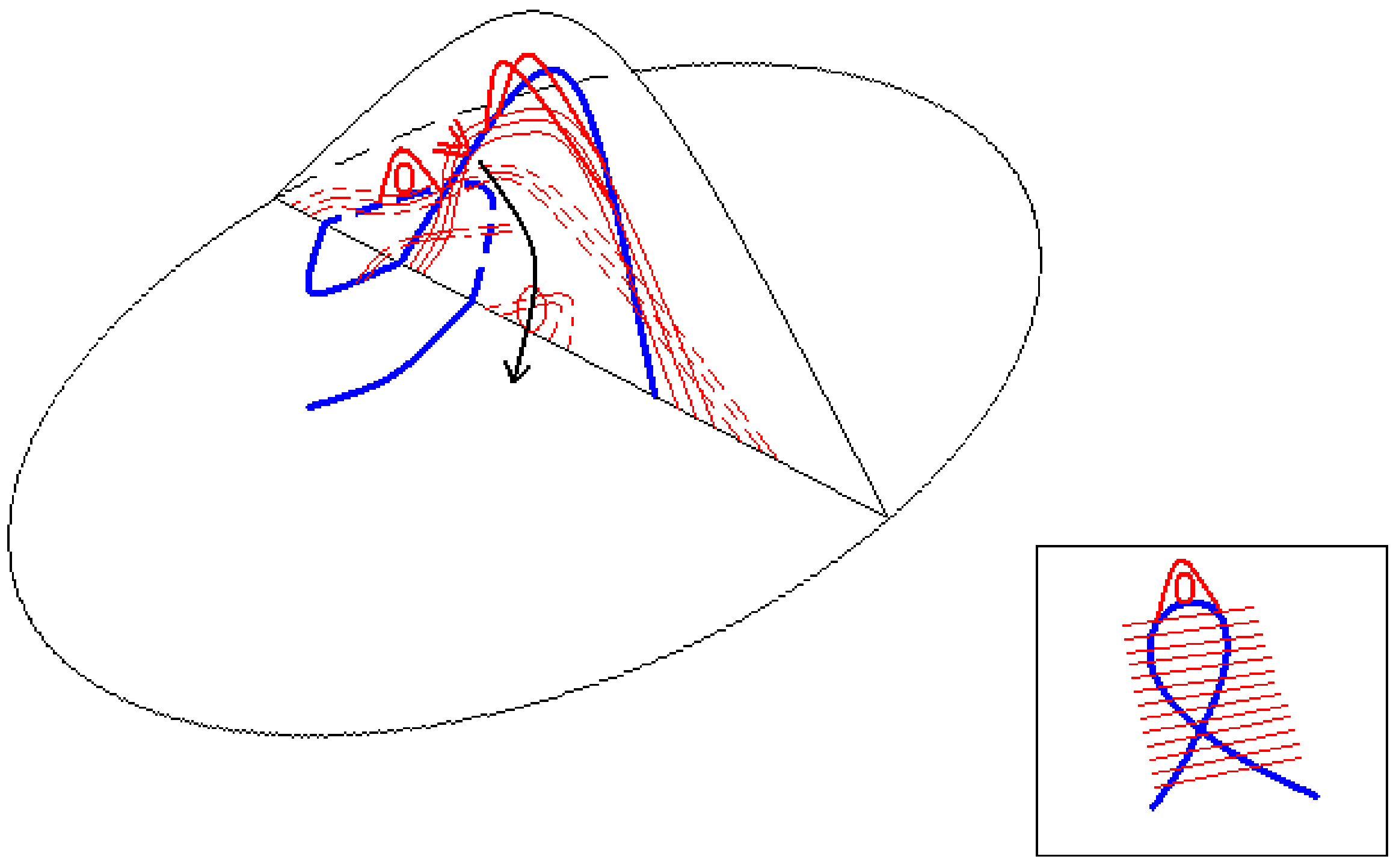

**Figure 49- the sliced leftside loop- 4**

we perform the step:

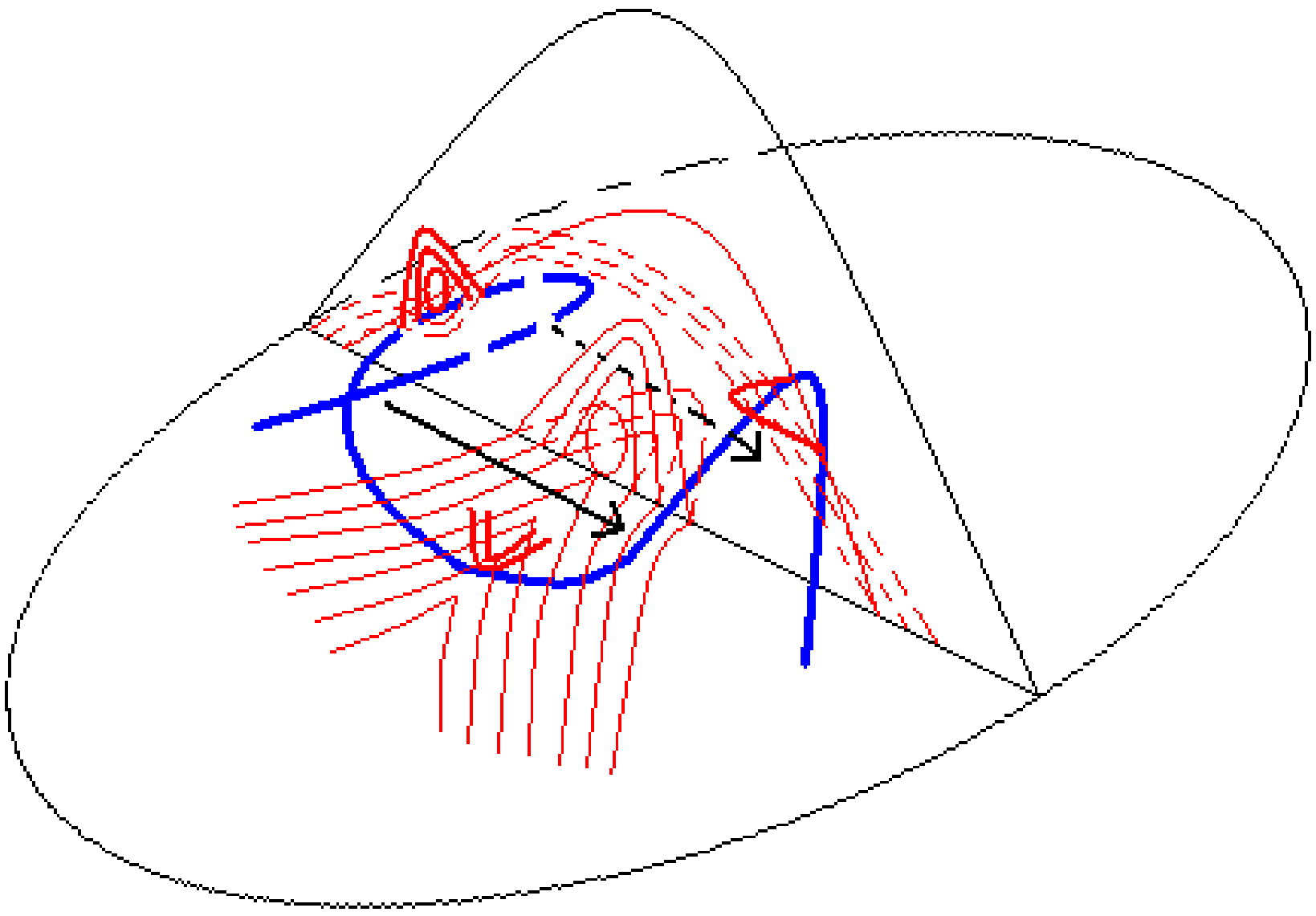

**Figure 50- the sliced leftside loop- 5**

We can perform the next step as indicated by the arrows, the description below shows the backside component:

The 2 figures on top show, that the step can be performed:
The attaching curve has a saddlepoint and this fits together with the slices of the backside component. We can slide the saddlepoint during the step to the turning point.
If there is no saddlepoint on the attaching curve, the 2 figures below show a contradiction:
Each slice of the component near the turning point of the attaching curve intersects this curve in 2 points, therefore the line segments at these points extend the slice on the attached 2-cell and these would join together at the turning point. However that is not permited by the Quinn list of local transitions:

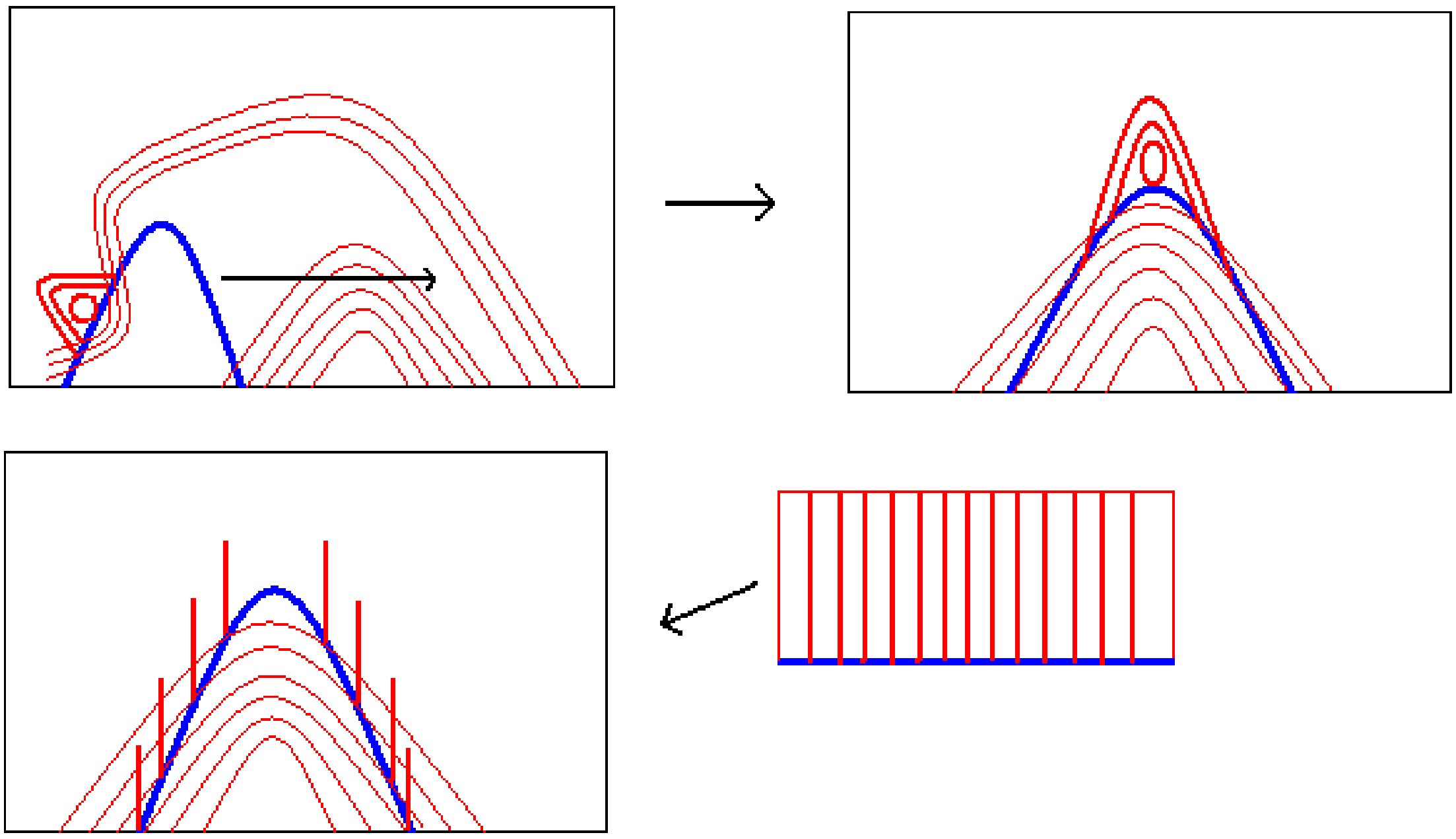

**Figure 51- the sliced leftside loop- 6**

We perform the step:

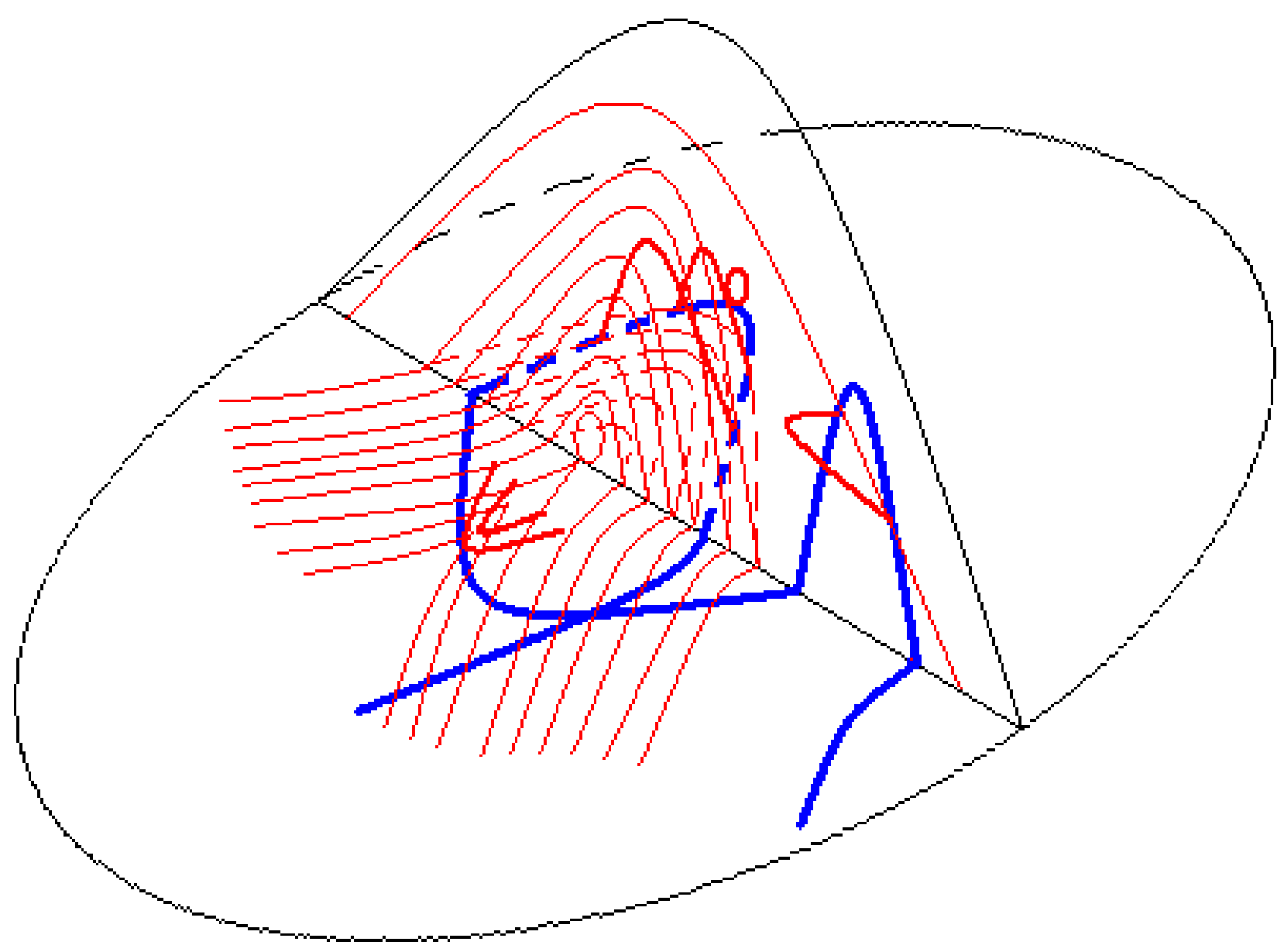

**Figure 52- the sliced leftside loop- 7**

We push the curve of the backside component across the vertex into the frontside component. We spread the slices and shrink the loop to avoid new saddlepoints and get:

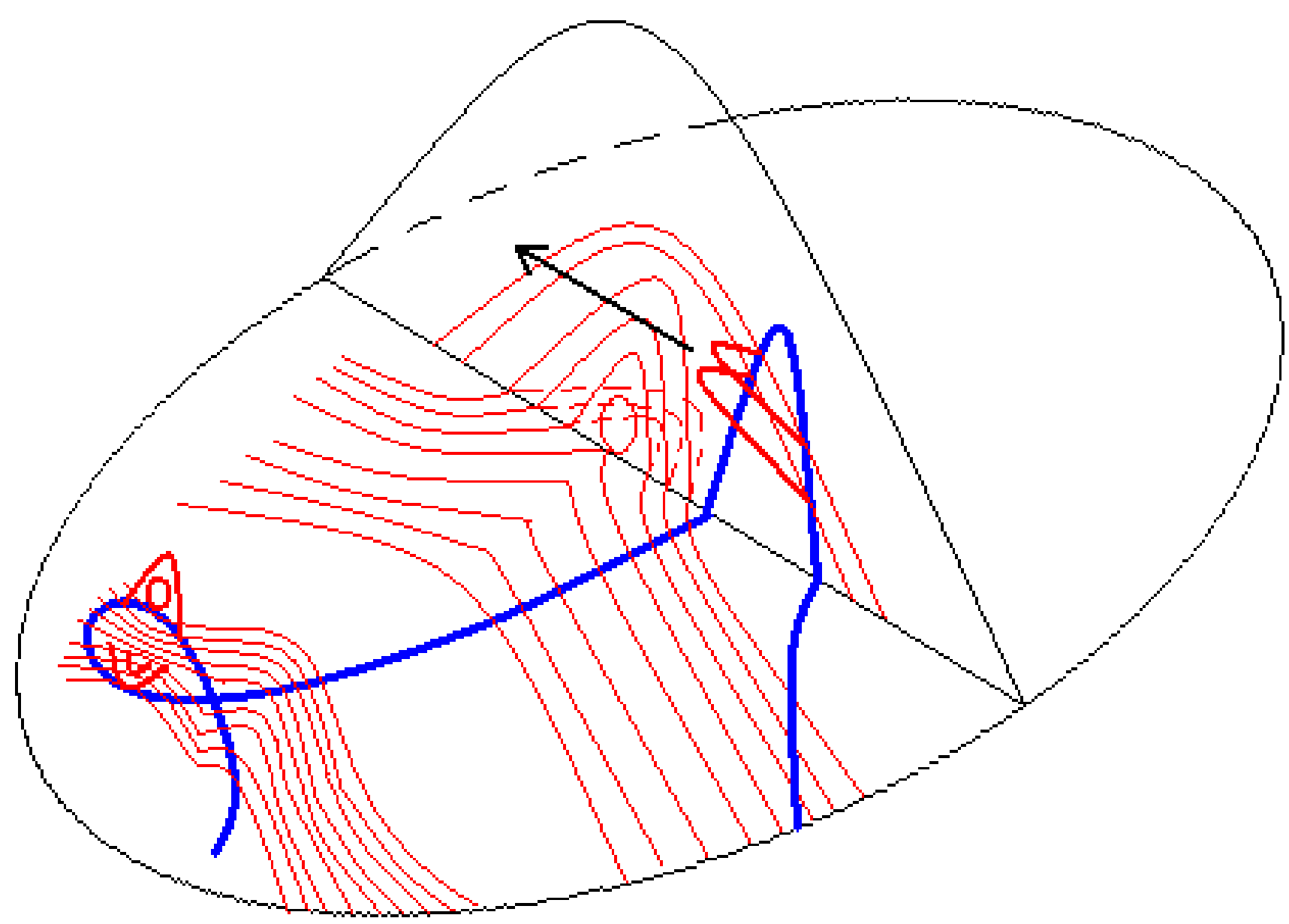

**Figure 53- the sliced leftside loop- 8**

The next step is indicated by the arrow and uses the same argument as before: We have to change the sequence of slices locally at the vertex. After performing the step, we change back to the original slices:

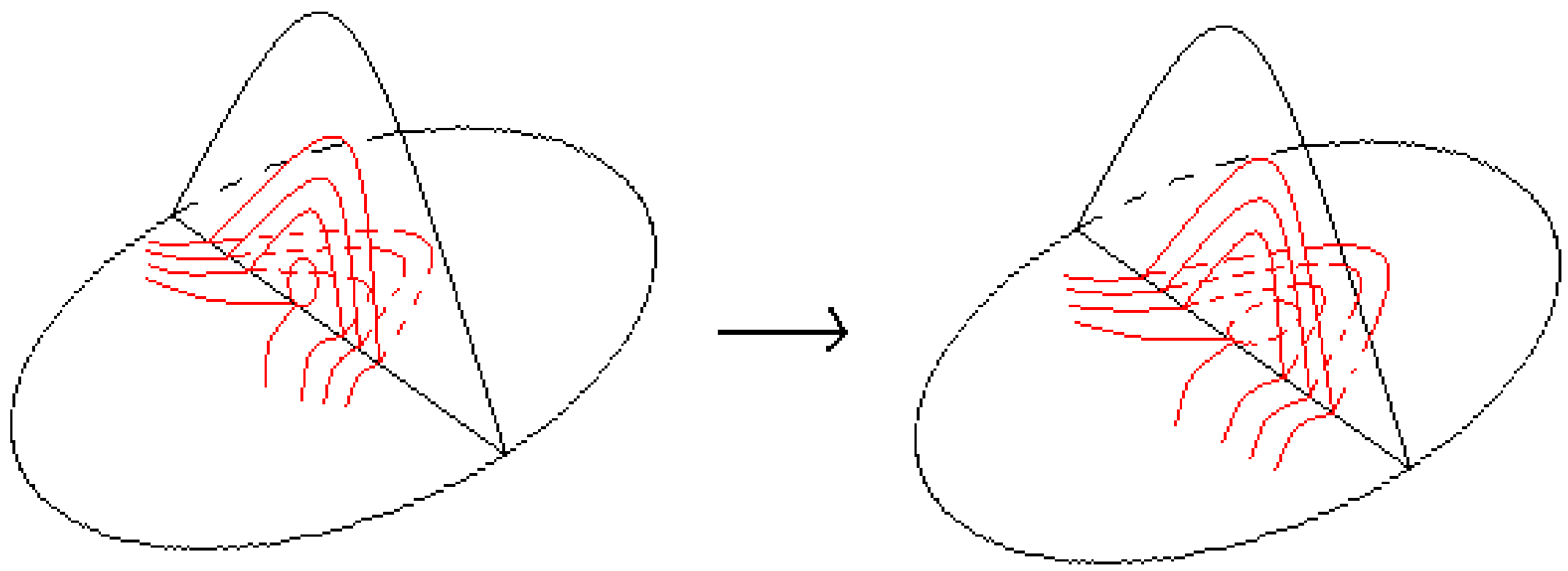

**Figure 54 - the sliced leftside loop- local change at vertex**

Result of the final sequence of slices for a leftside loop is:

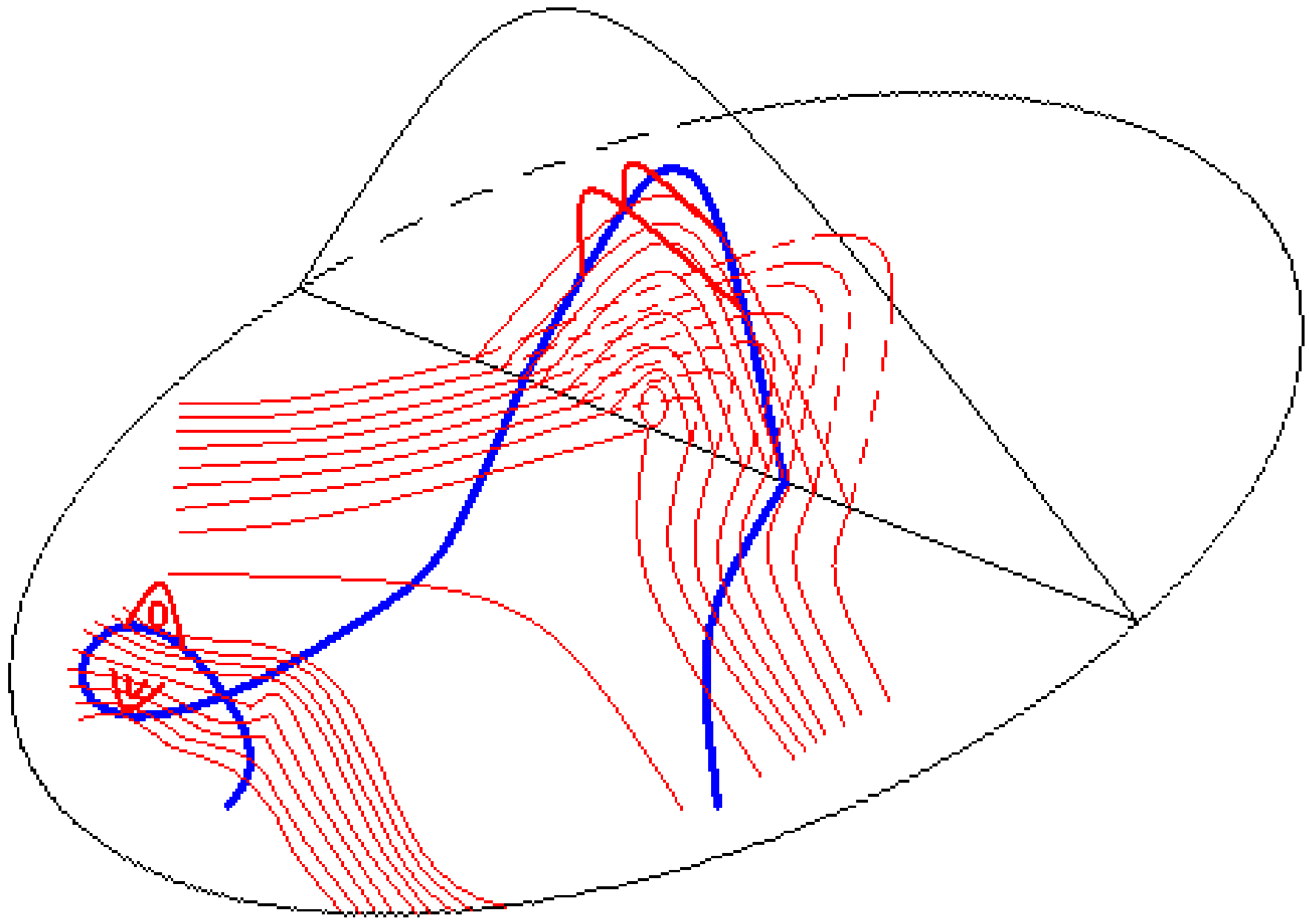

**Figure 55- the sliced leftside loop- 9 (end)**

Sometimes the attaching curve can have selfintersection without a twist. In the local vertex model the top component is crossed transversally by the bottom component. In that case we present the sequence of slices than as a “2 times flip”:

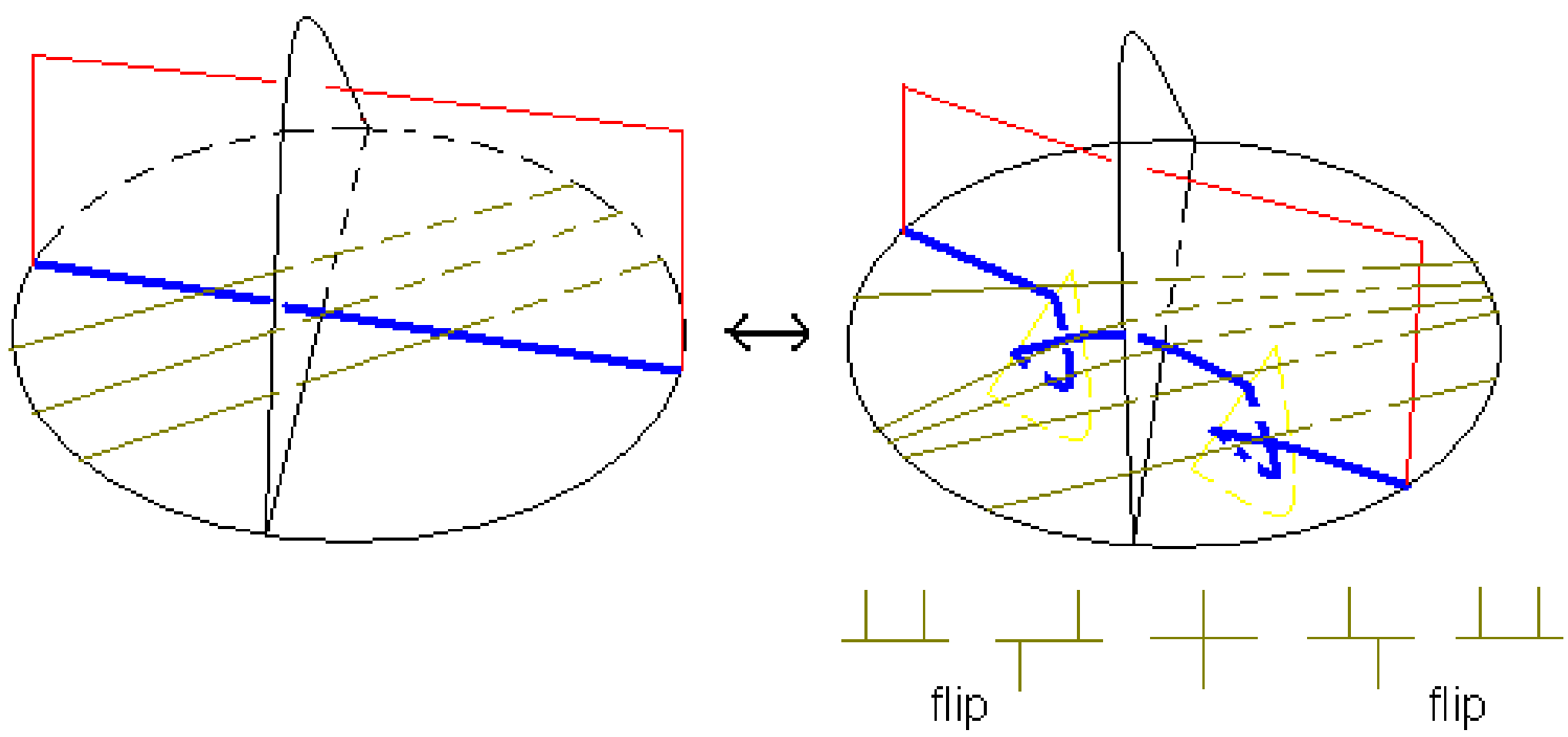

**Figure 56- the sliced leftside loop- selfintersection of attaching curve without twist present as 2 times flip**

This finishes our considerations about the construction of a loop and shows, that it is necessary to modify the slices depending on the upcoming step.

# 5 Subsolutions

## 5.1 *Cancellation of a loop pair- a special case*

First we show the process of cancelling a pair of leftside loops and present for each step the sequence of slices. Note, there appears a restriction to the variations an attaching curve can have; we consider not more than a pair of loops !!!
To loosen the restriction to leftside loops, we have to turn rightside loops into leftside loops. It is recommended to be familiar with the construction of a leftside loop and their asssociated sequence of slices.

### 5.1.1 Construct leftside and rightside loop

We start with the process of constructing a leftside respectively a rightside loop, where we omit the sequence of slices:

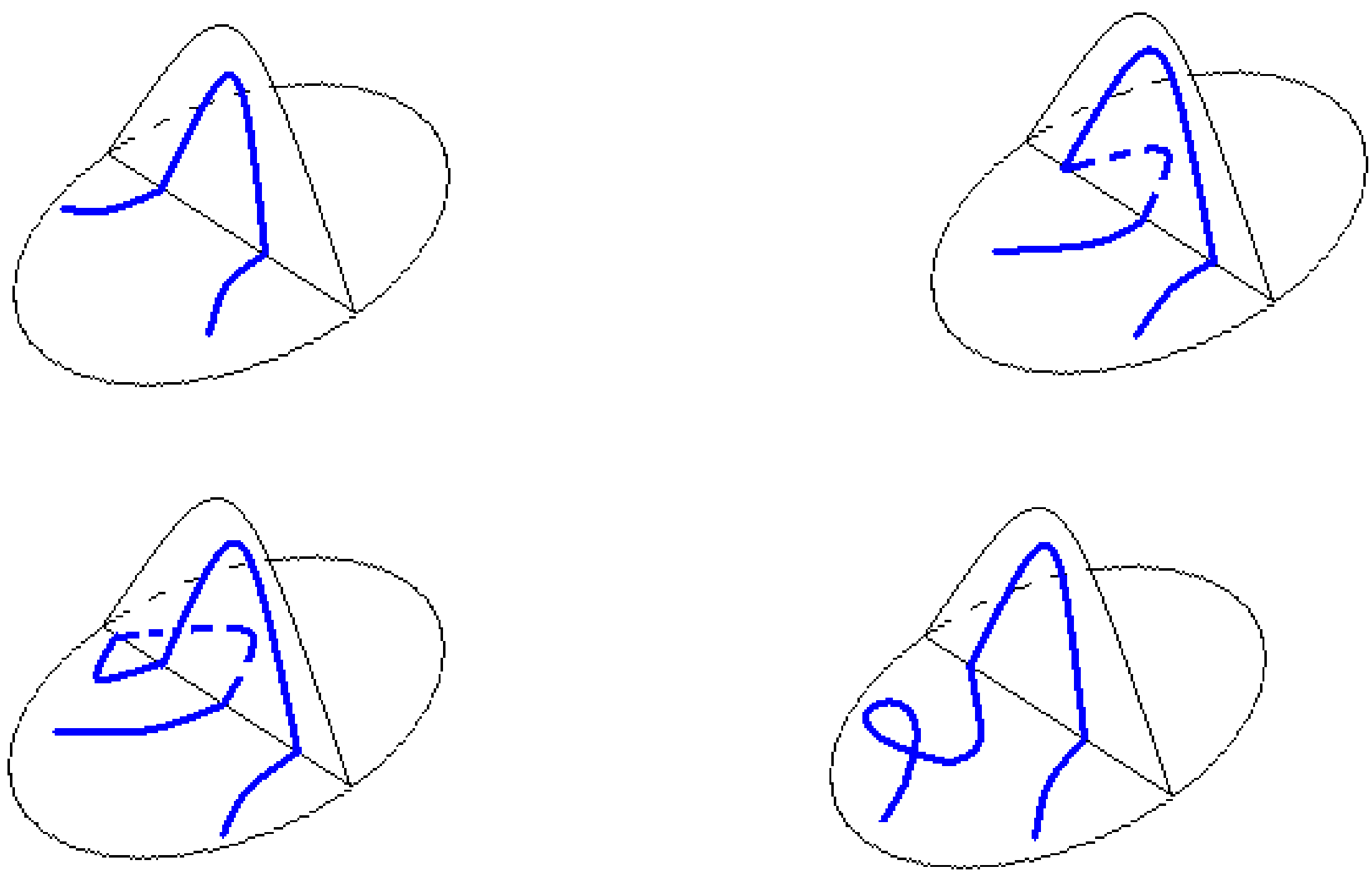

**Figure 57- Construct leftside and rightside loop- leftside loop**

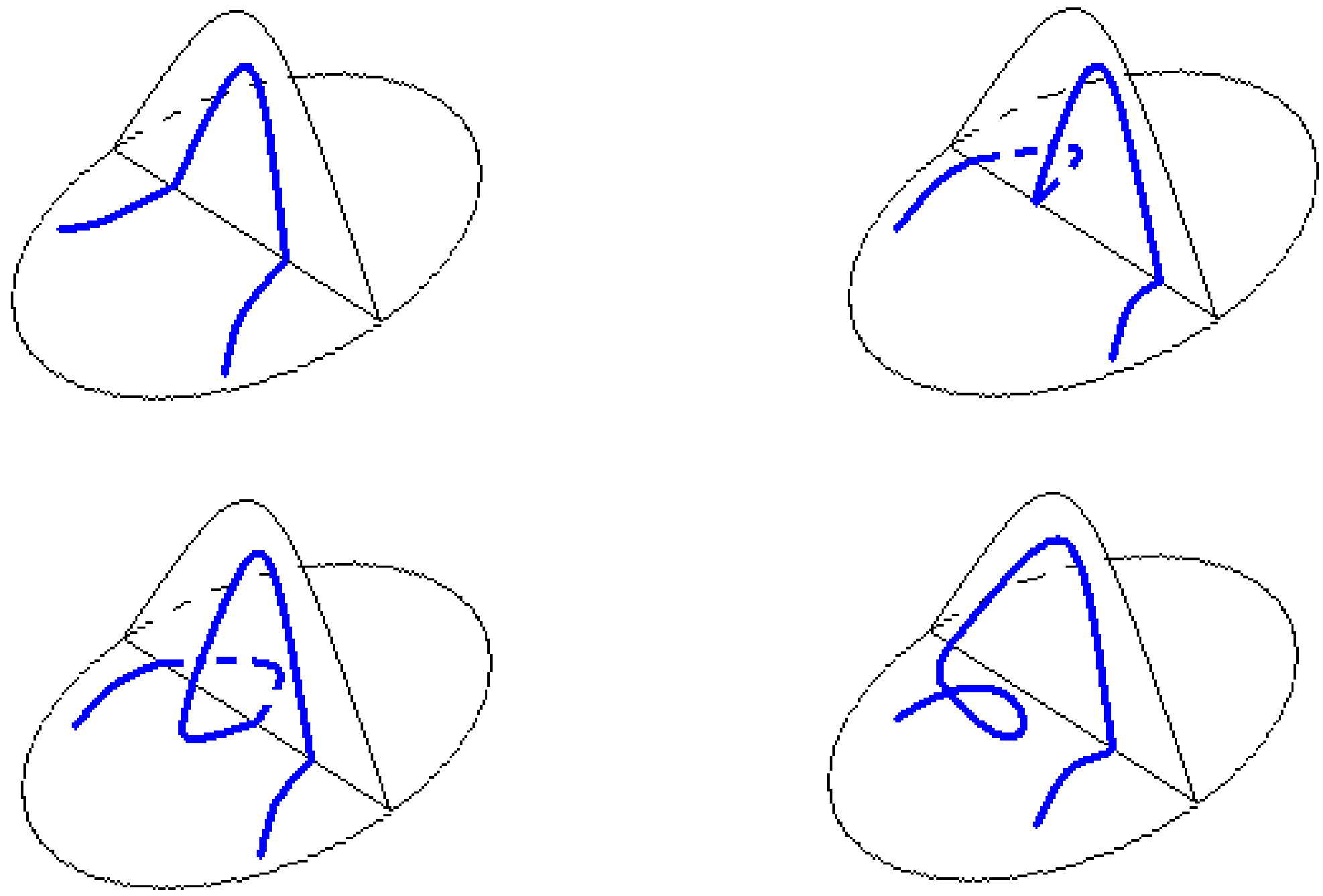

**Figure 58- Constuct leftside and rightside loop- rightside loop**

### 5.1.2 Cancel a pair of leftside loops

From now on we concentrate on the cancellation of a pair of leftside loops. We follow the idea of **[Mat]** and stack up the loops. We show, how to do that for 2 loops:

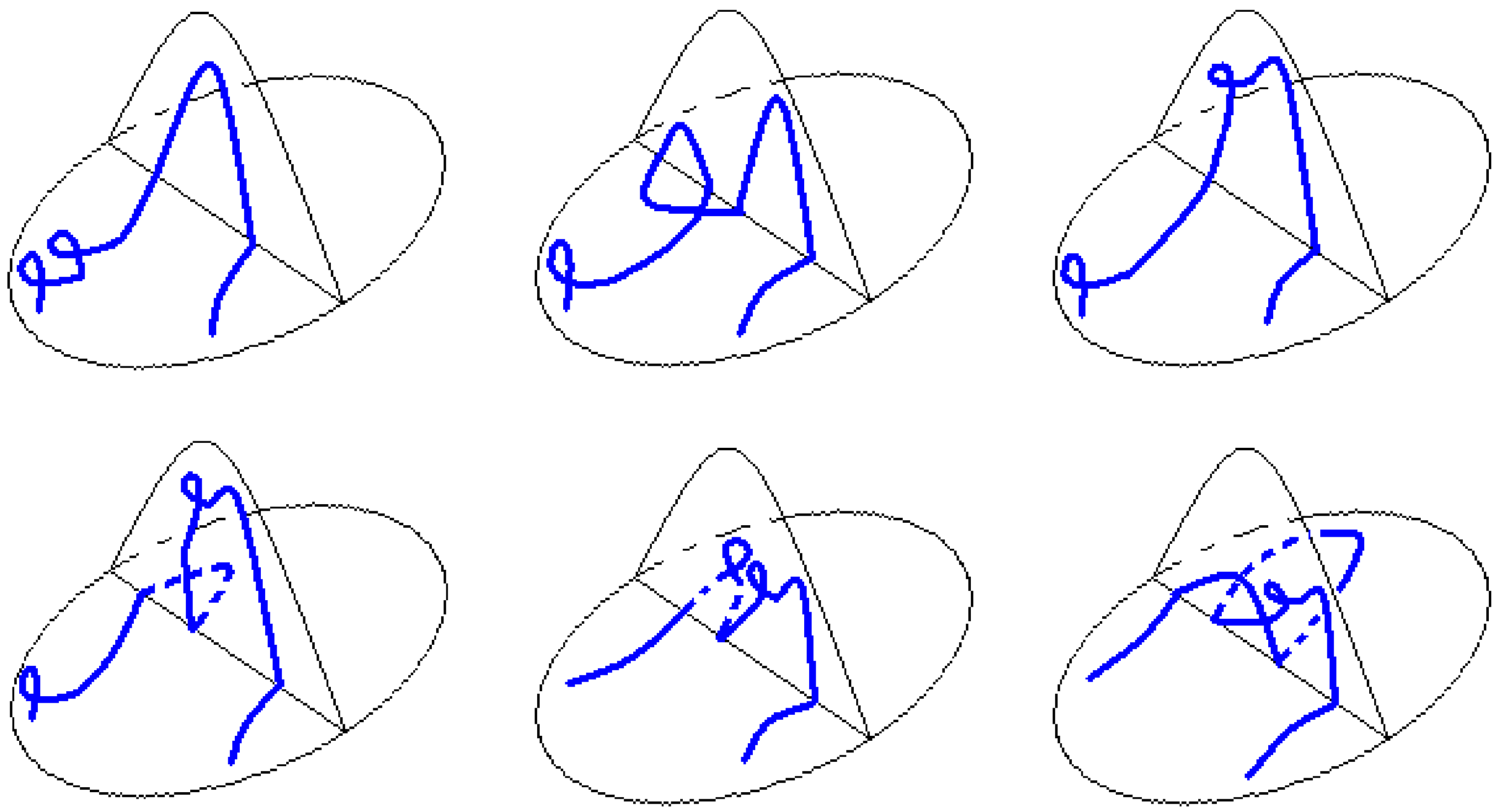

**Figure 59- Cancel a pair of leftside loop- preview part 1**

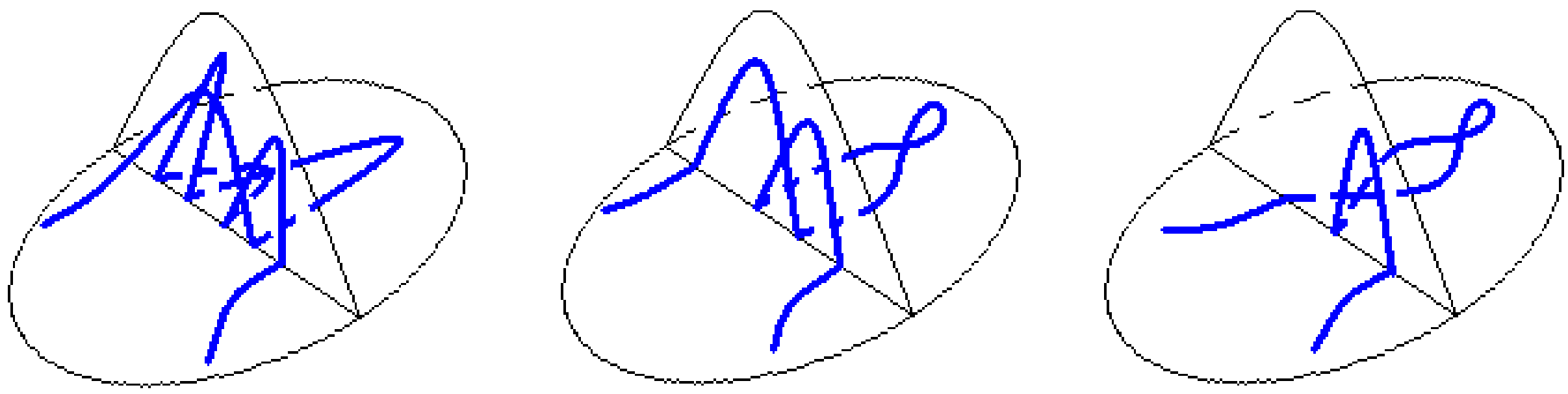

**Figure 60- Cancel a pair of leftside loops- preview part 2**

Note that perhaps we are forced to modify some steps, but we will always orientate on the given sequence. However the question is, whether there exists a sequence of slices to terminate the process ?
We will only draw for each step at most the essential parts of the sequence of slices to point out that these are the most important ones for performing the next step.
The figures contain detailed information on the sequence of slices. We do not draw arrows to avoid covering these informations.
The start figure describes the blue attaching curve on a vertex model, where the sequence of slices (red thin lines) illustrates a saddlepoint (for example the entry of a 2-cell) at the vertex. The sequence of slices of the attached 2-component is indicated by thick red lines:
For each loop we have 2 saddlepoints (not as a cancellation pair) and at turning point a saddlepoint (the exit of the attached 2-cell). Note, that the sequence of slices on the top component was constructed, such that it intersects the attaching curve in a pair of points, and one intersection at turning point. This justifies the choice of slices of the attached 2-component near that point.

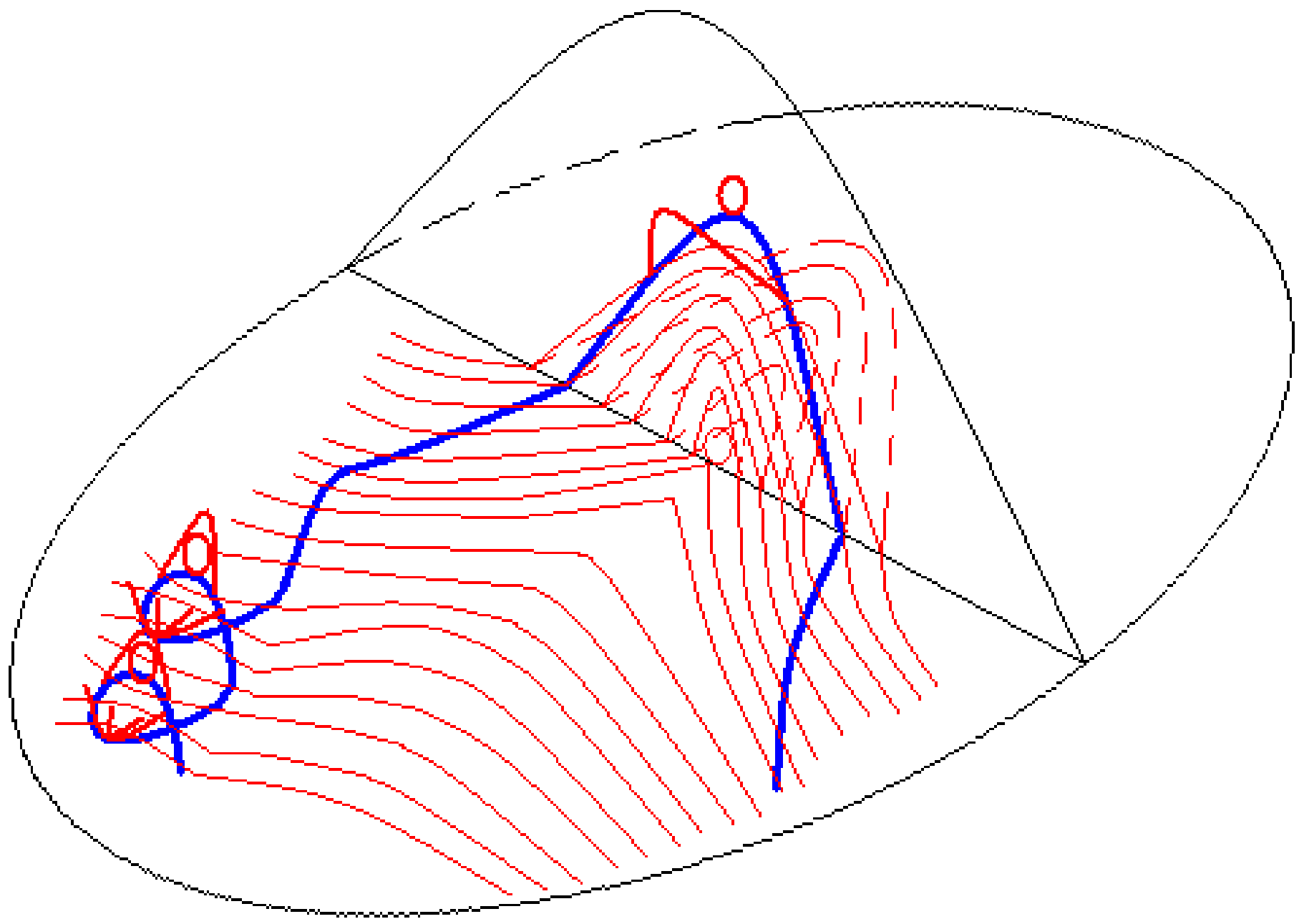

**Figure 61- Cancel a pair of leftside loops- 1 (start)**

We change the local sequence of slices at the vertex as described in chapter 4 and push the left arc of the attaching curve across the vertex in the top component and change the local sequence back:

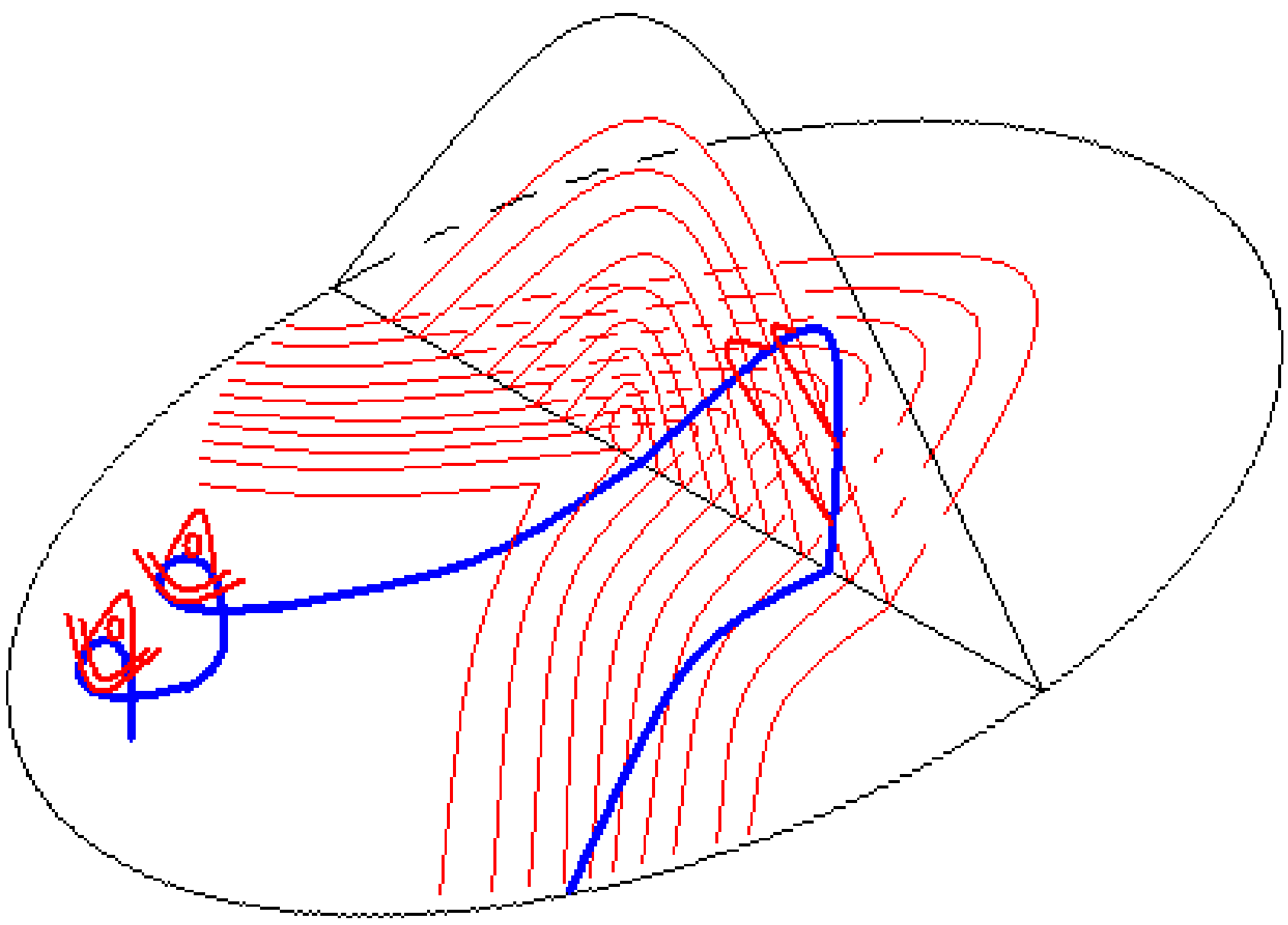

**Figure 62- Cancel a pair of leftside loops- 2**

We push one arc of the loop of the fronside component across the vertex into the top component. Note, that since the "rhythm" of the slices at the attached 2-component near the turning point agree with the "rhythm" of the slices near the center,  this step can be performed without problems:

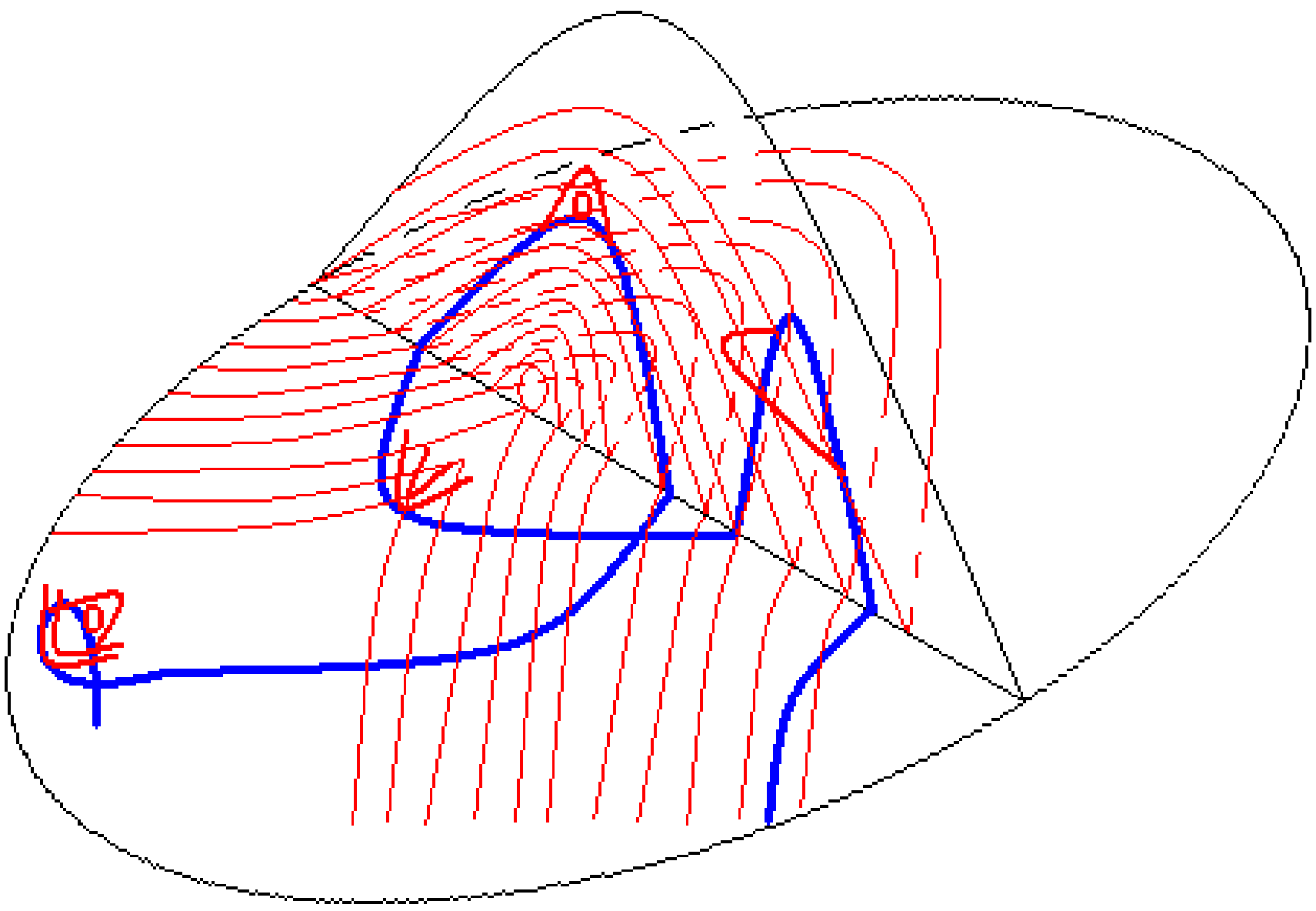

**Figure 63- Cancel a pair of leftside loops- 3**

We repeat this step with the remaining arc of the first loop in the frontside component:

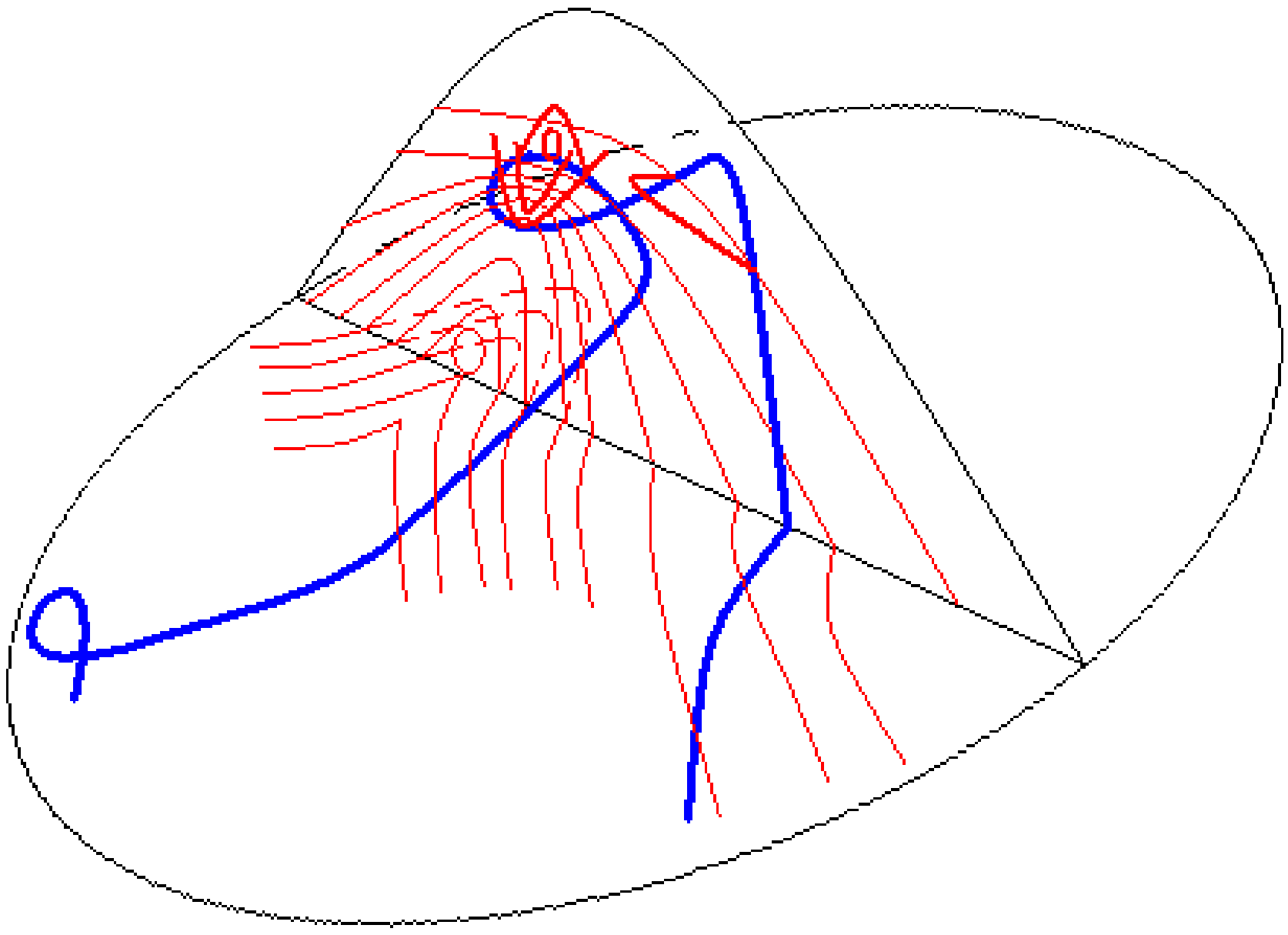

**Figure 64- Cancel a pair of leftside loops- 4**

We perform a $T_3$ move near the vertex and by observing the slices we see it results to a “good $T_3$ turn” in the backside component:

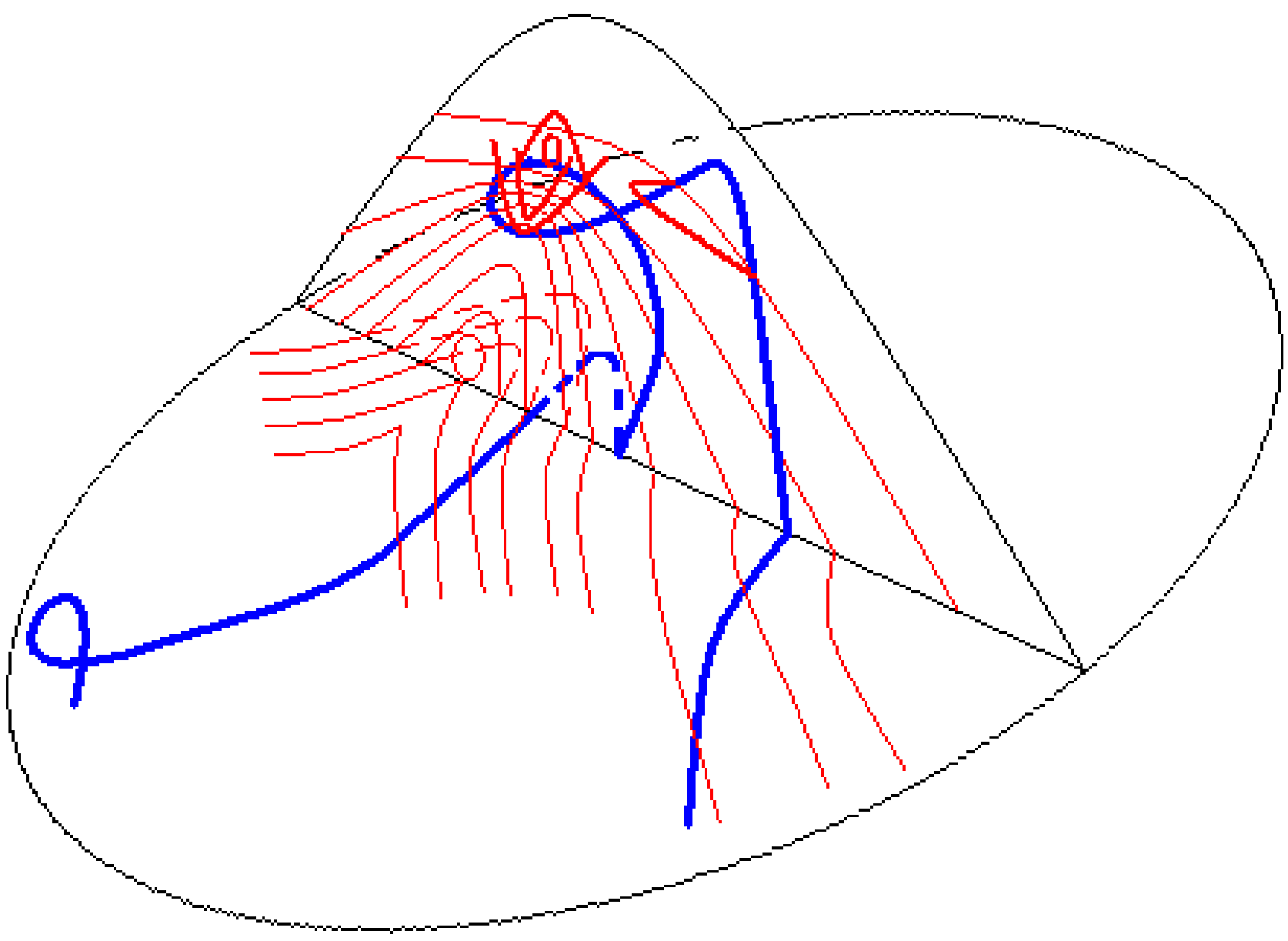

**Figure 65- Cancel a pair of leftside loops- 5**

We would like to have the second loop in the backside component:

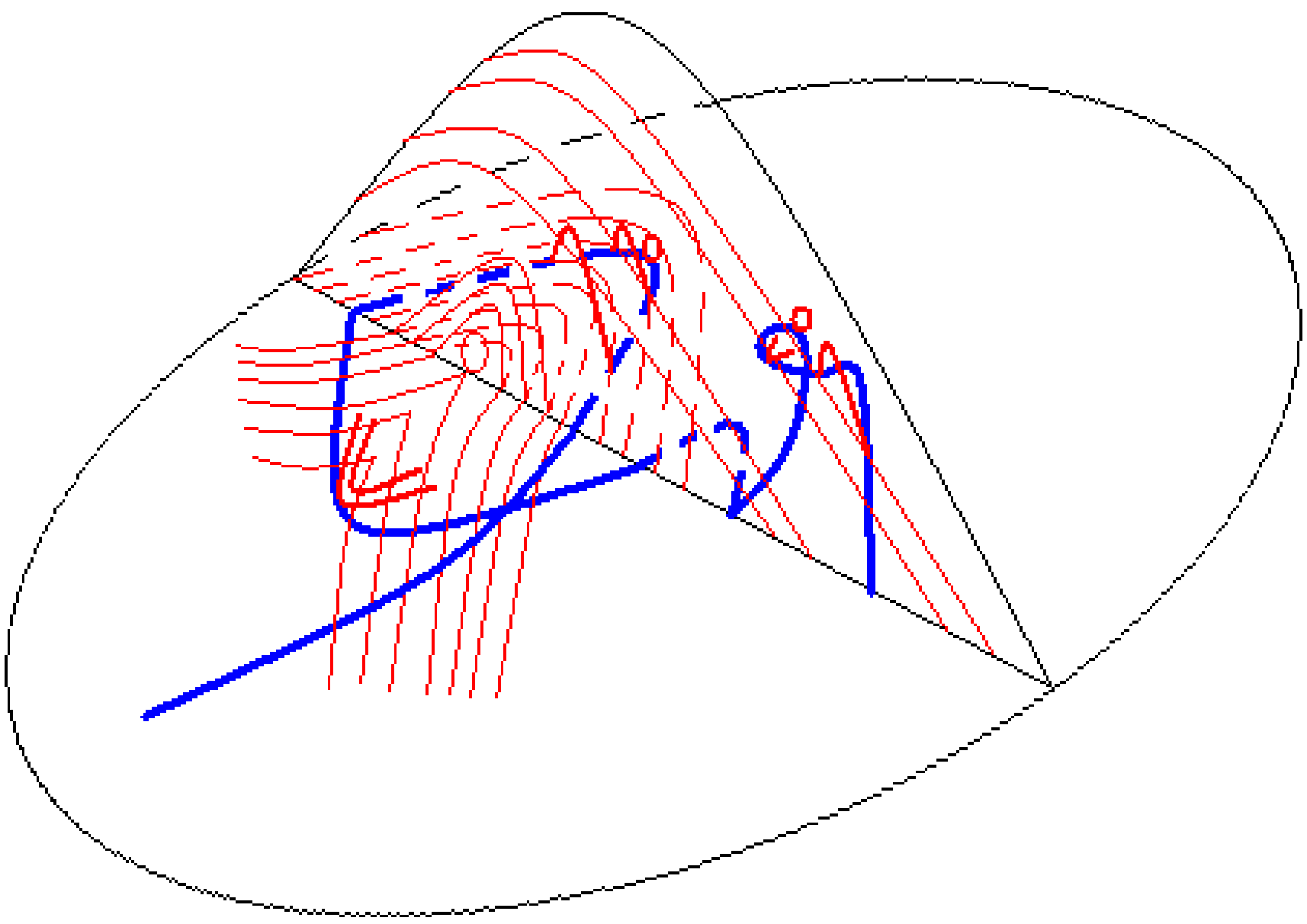

**Figure 66- Cancel a pair of leftside loops- 6**

We push the second loop completly into the backside component.
The second figure illustrates the sequence of slices on the attaced 2-cell:

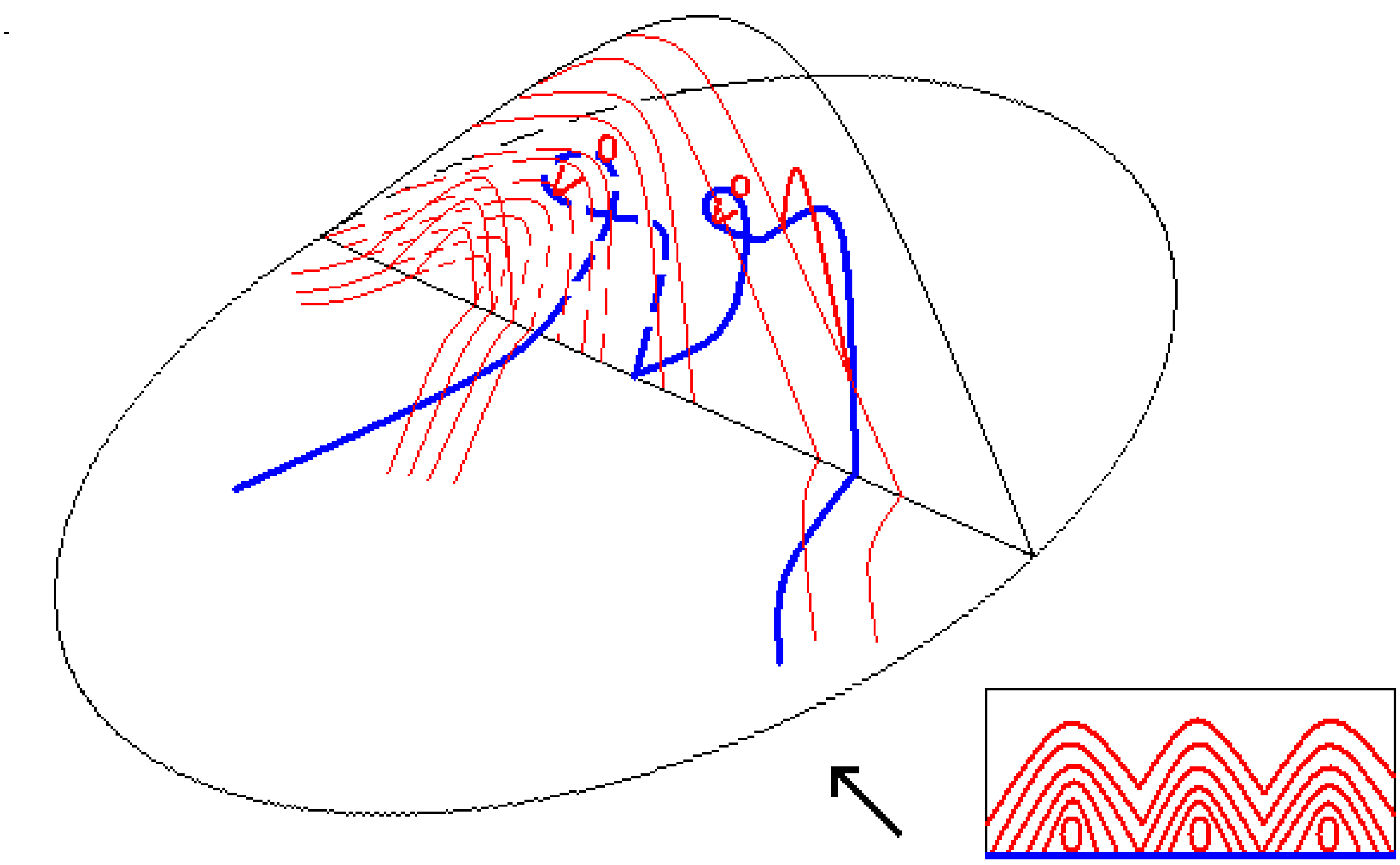

**Figure 67- Cancel a pair of leftside loops- 7**

### 5.1.2.1 Stack up the loops

The next step is a composition of Matveev moves. In the first step we push an arc of the second loop of the backside component along its crossing arc until the arc is completely on the top component. The thick black line marks the boundary that belongs to the bottom component in the corresponding vertex model:

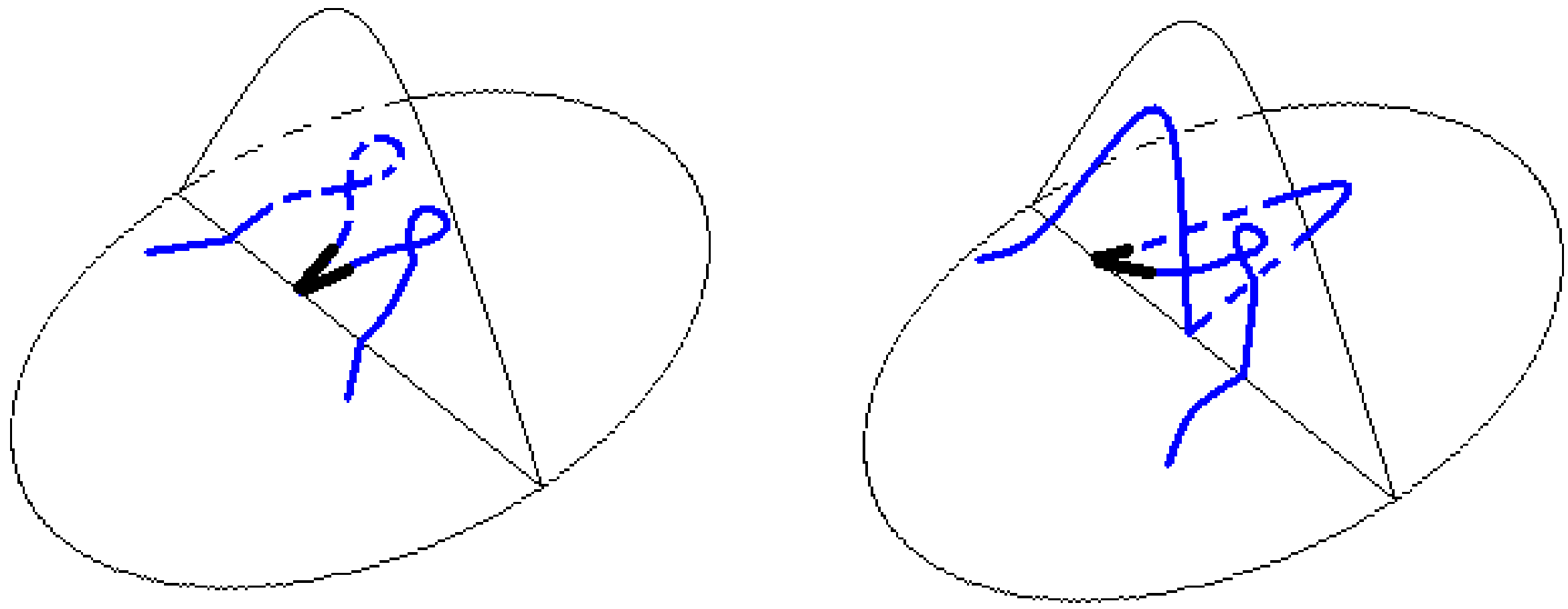

**Figure 68- Stack up loops- Preview- resolve loop in backside component**

It follows a decomposition of this move in Matveev moves, for the relevant parts in the local vertex model. In the following figures the backside component of the former figure is the frontside component:

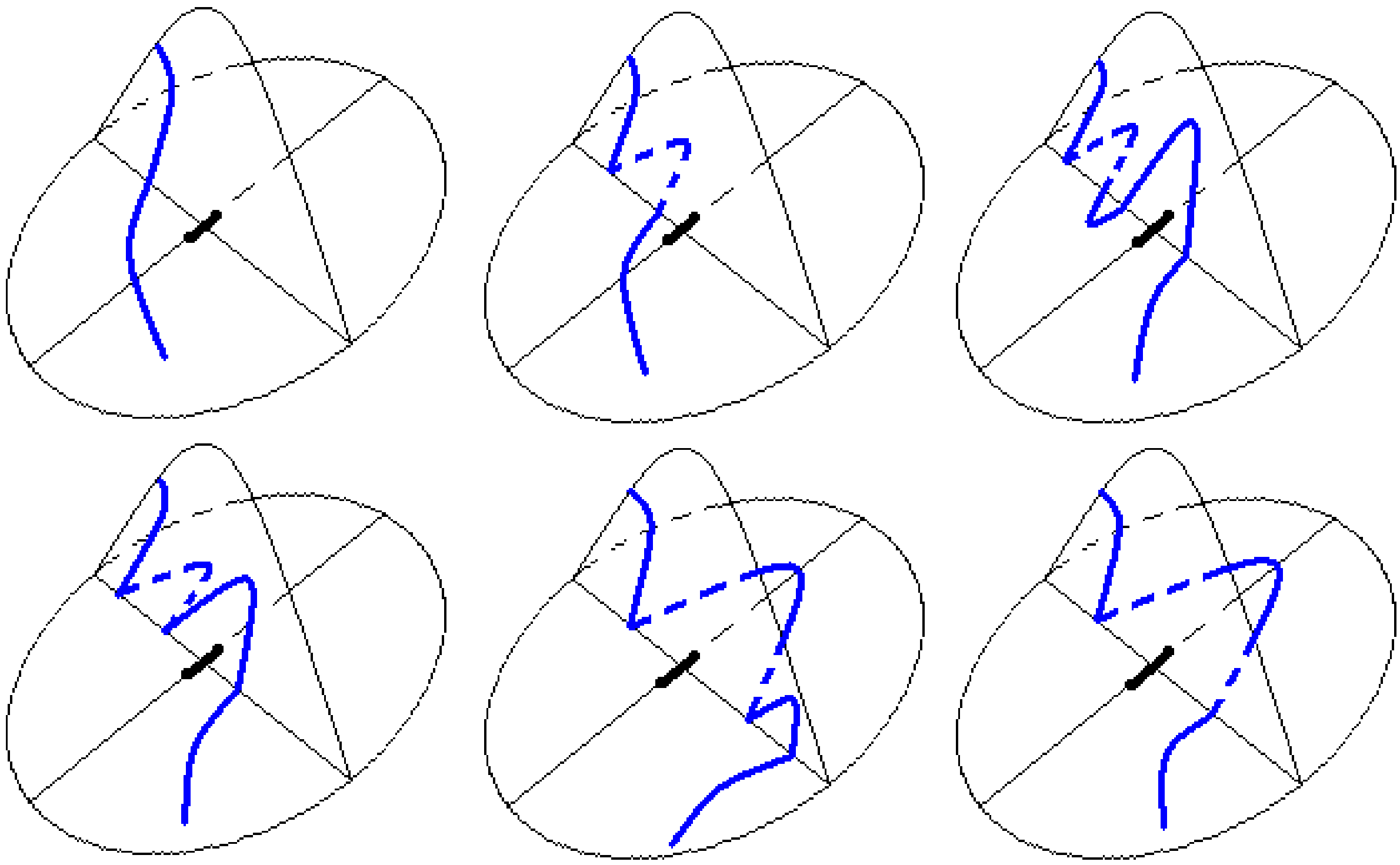

**Figure 69- Stack up loops- resolve loop in backside component- local sequence**

We embed the local sequence into a global sequence:

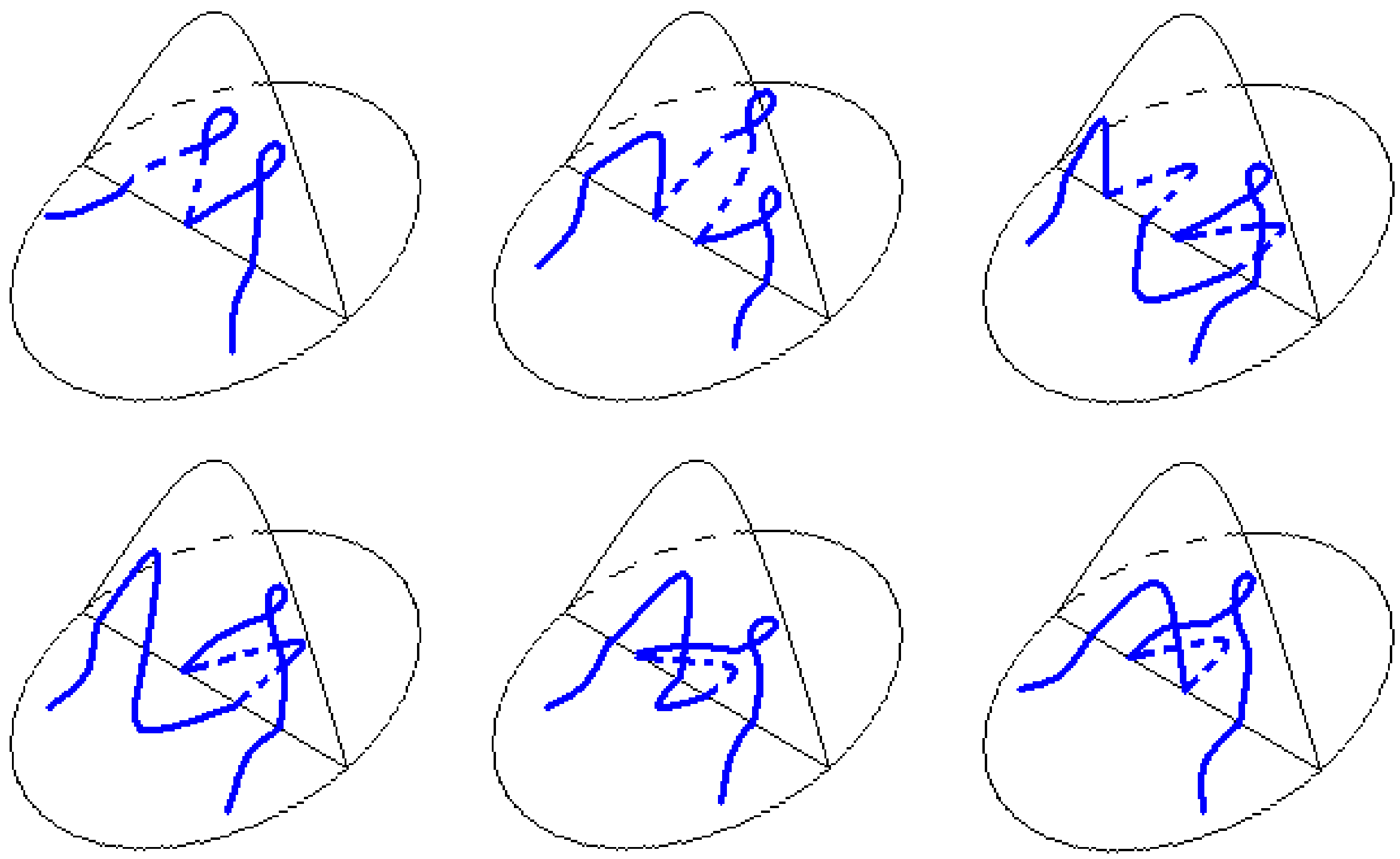

**Figure 70- Stack up the loops- resolve loop in backside component- global sequence**

We list the sequence of slices for each move. Our start figure is:

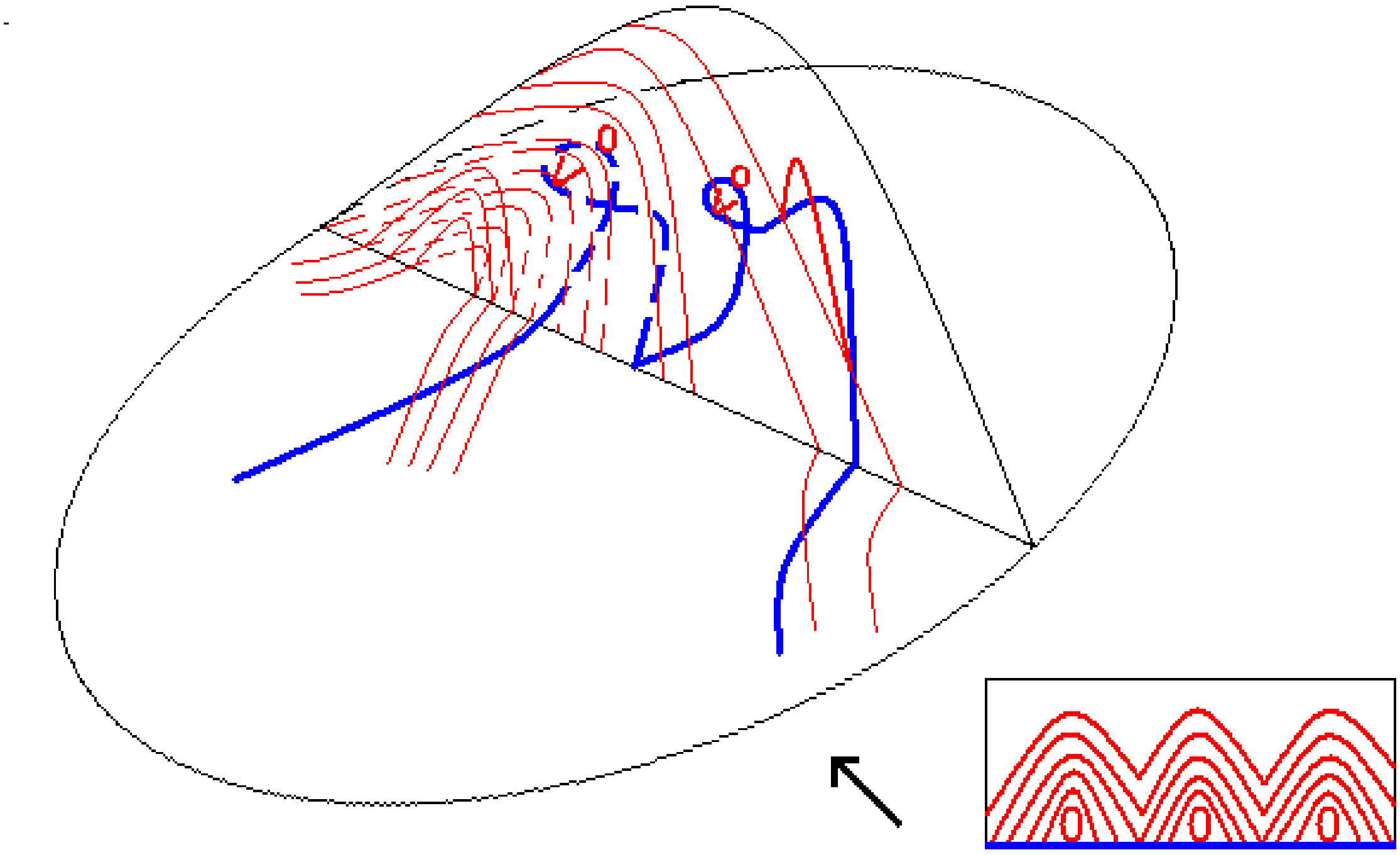

**Figure 71- Stack up loops- resolve loop in backside component- 1 (start)**

We perform the$T_3$ move and see it is a “good” one:

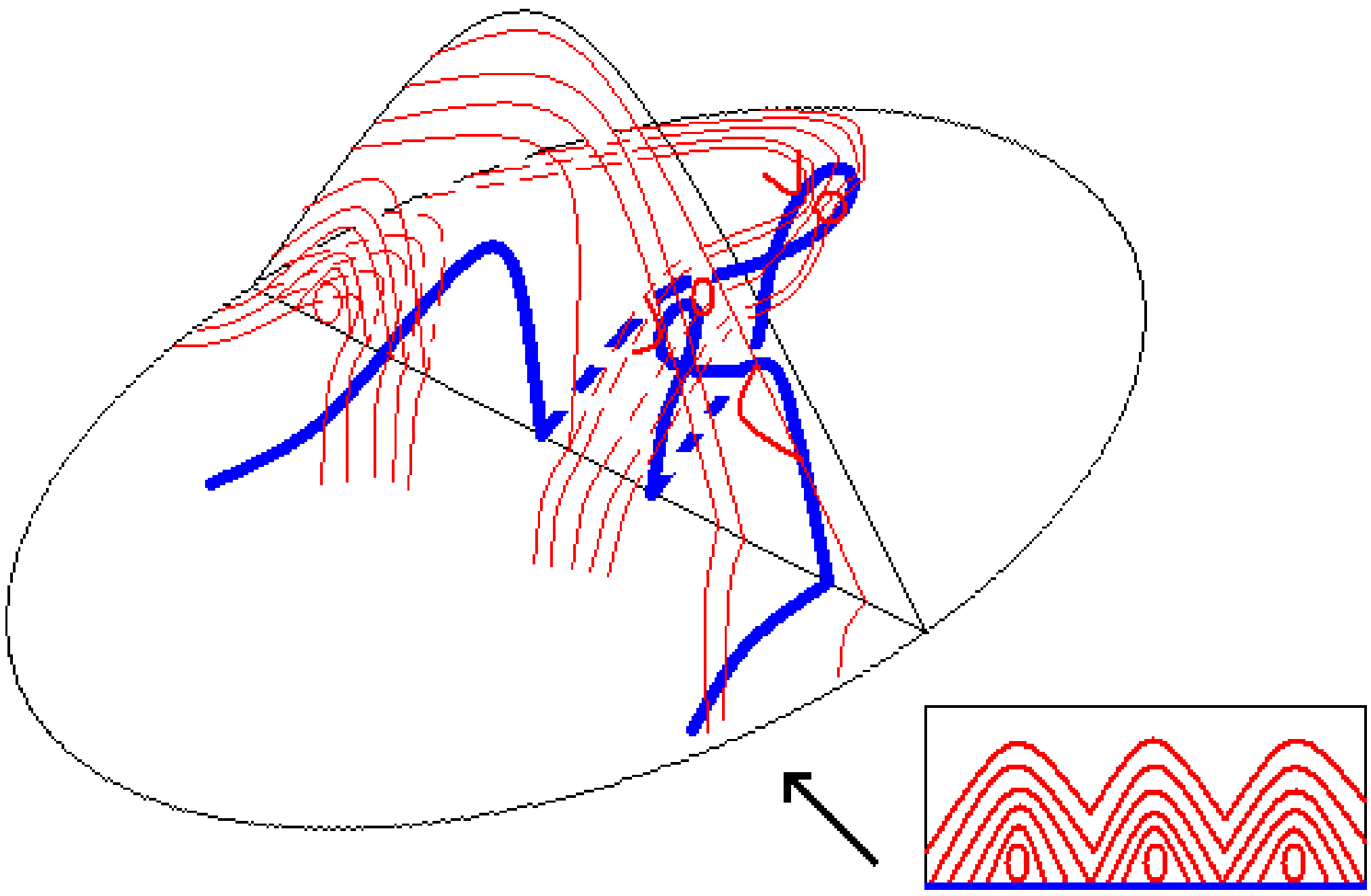

**Figure 72- Stack up the loops- resolve loop in backside component- 2**

We perform the $T_2$ move on the loop in the backside component and thereby we resolve the loop. We remark, that there arises two $T_3$ turns in the backside component:
The first one (counting from left to right) is “good” and the other one is sliced with two saddlepoints (which comes from the loop) and therefore also “good” :

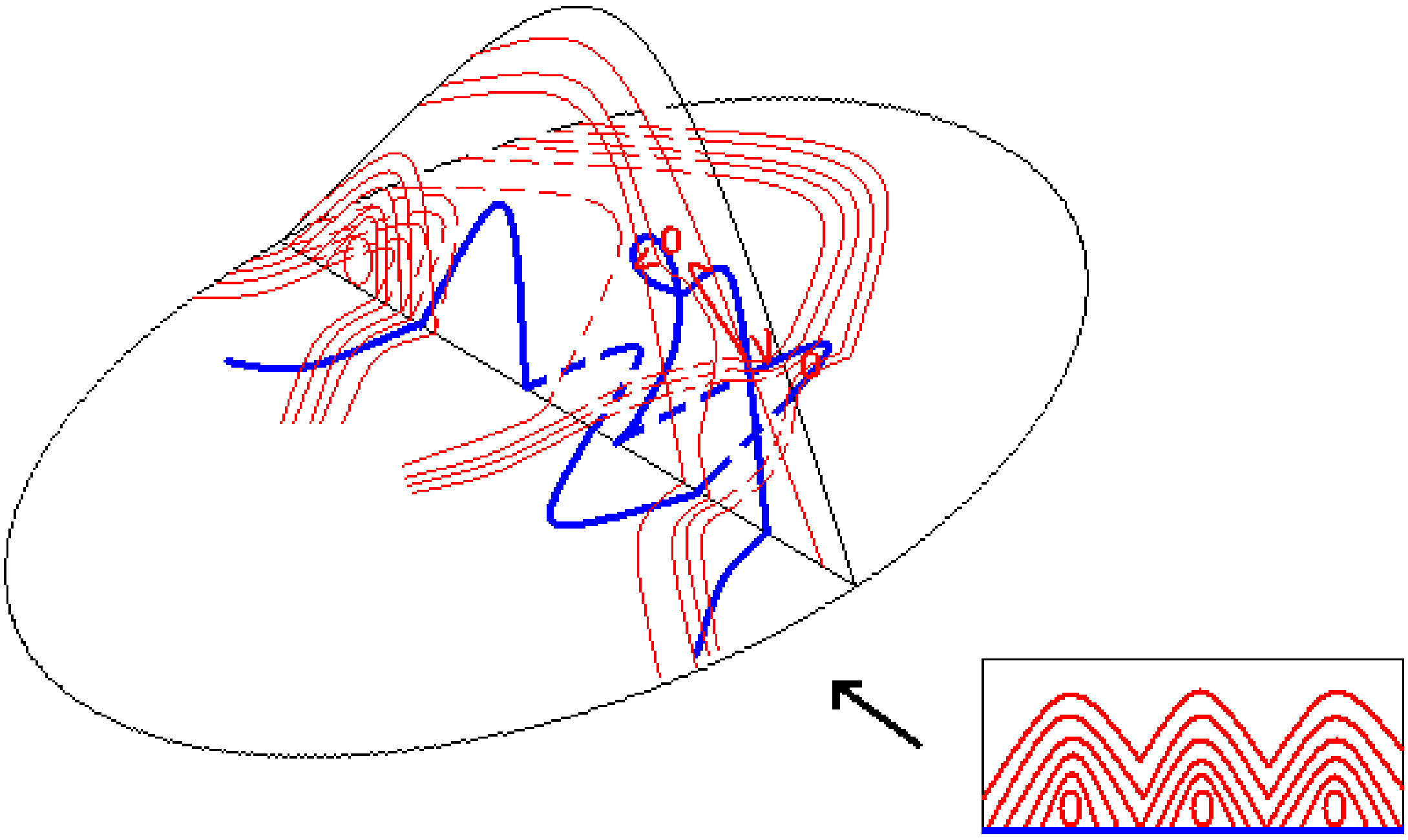

**Figure 73- Stack up loops- resolve loop in backside component- 3**

We cancel the first one of that $T_3$ turn by its inverse. Note that the new $T_3$ turn in the frontside component is also a “good” one:

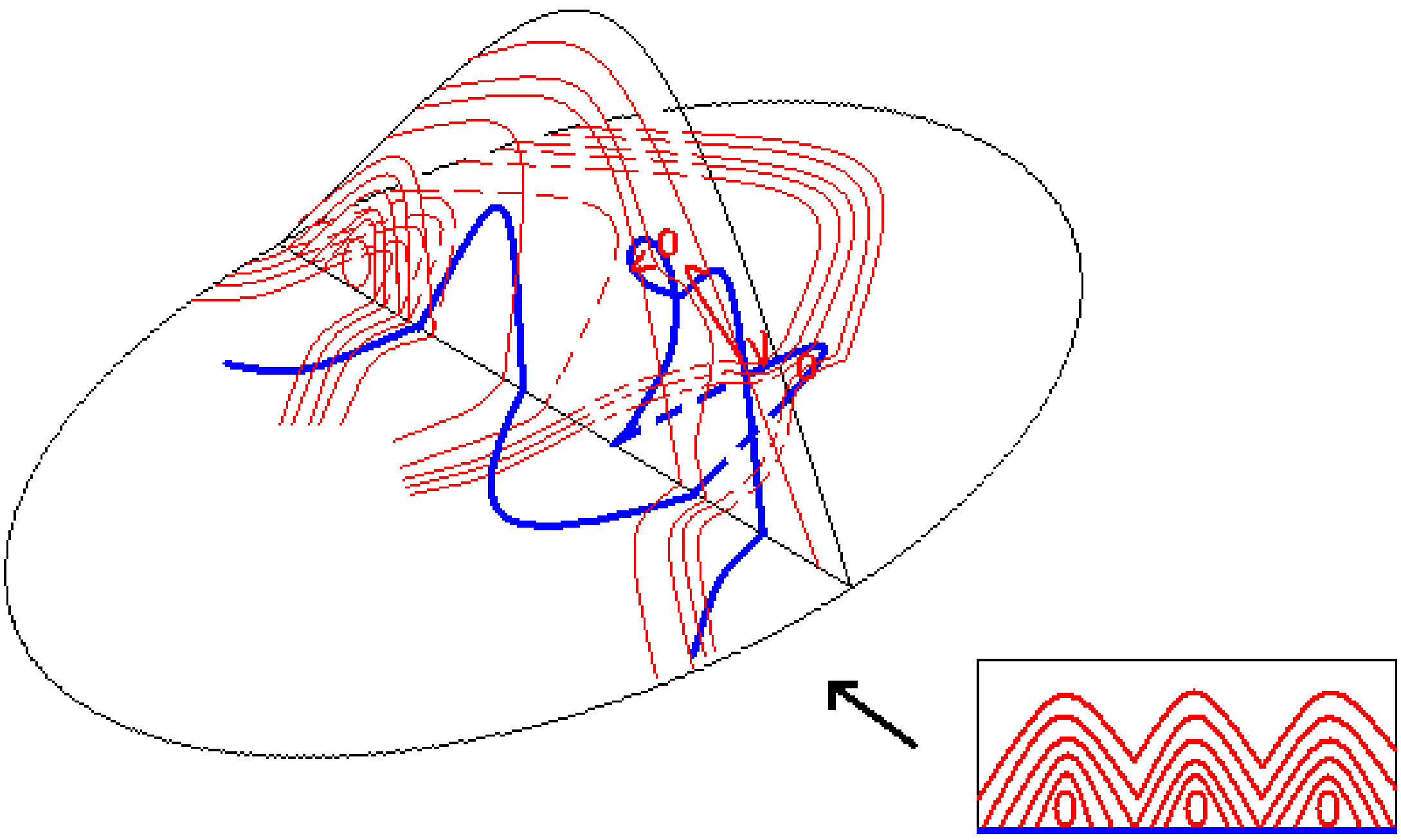

**Figure 74- Stack up loops- resolve loop in backide component- 4**

We perform T* to the $T_3$ turn in the frontside component and come to:

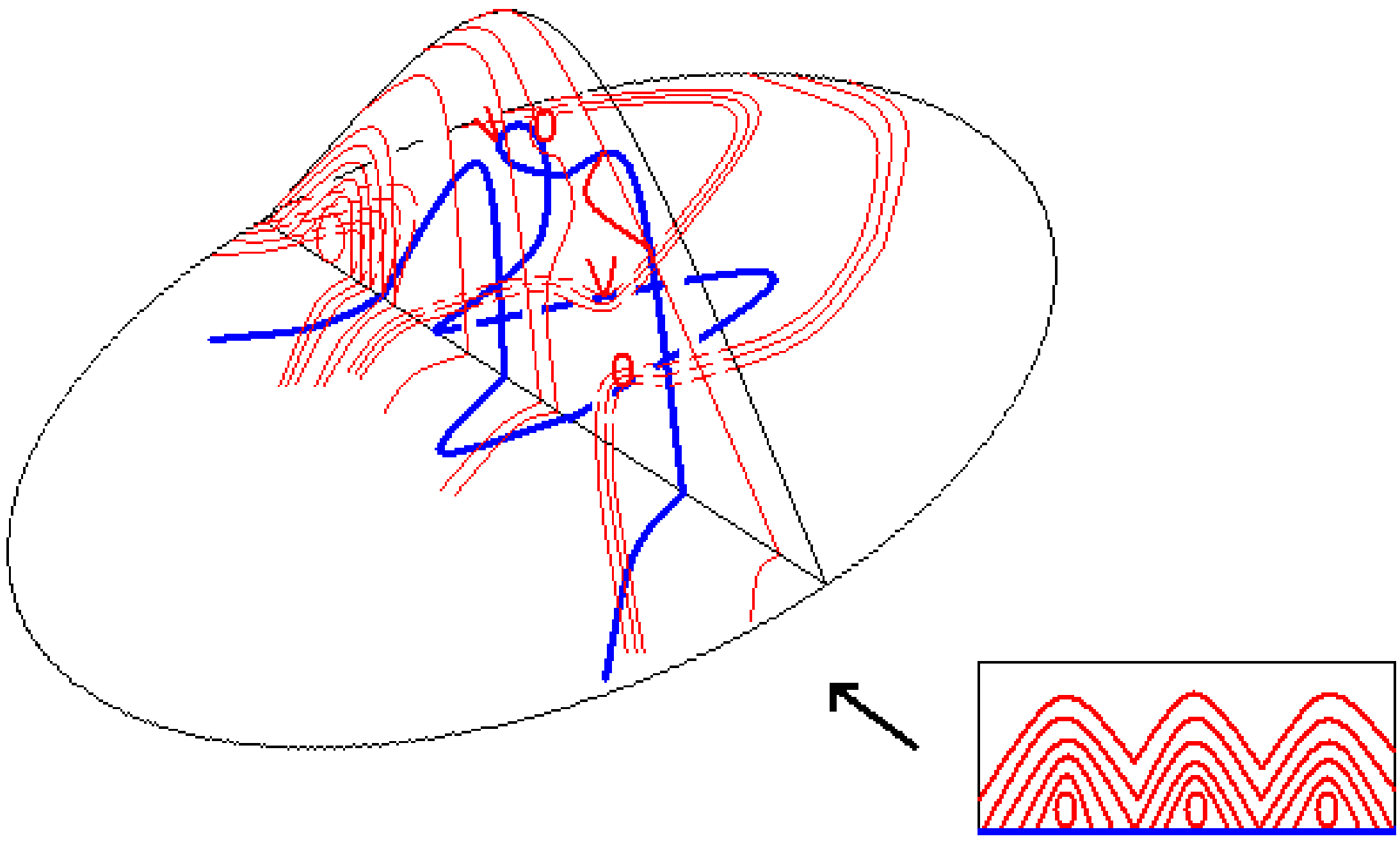

**Figure 75- Stack up the loops- resolve loop in backside componet- 5**

We annihilate the $T_3$ turn on the frontside component and reach the end of that sequence:

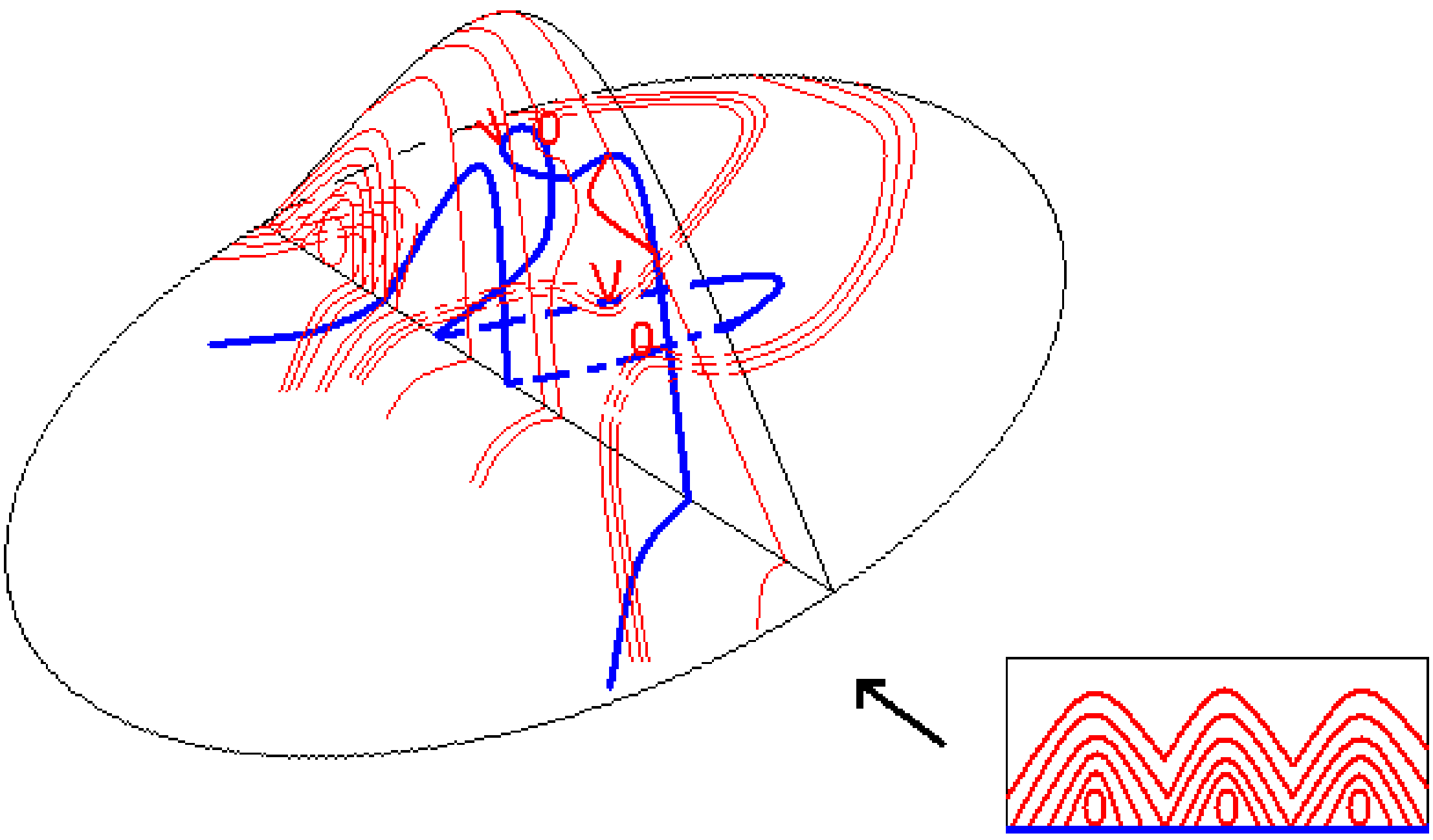

**Figure 76- Stack up the loops- resolve loop in backside component- 6 (end)**

The next step is similar to the former one. We resolve the loop in the top component, through pushing the arc along its crossing arc into the backside component:

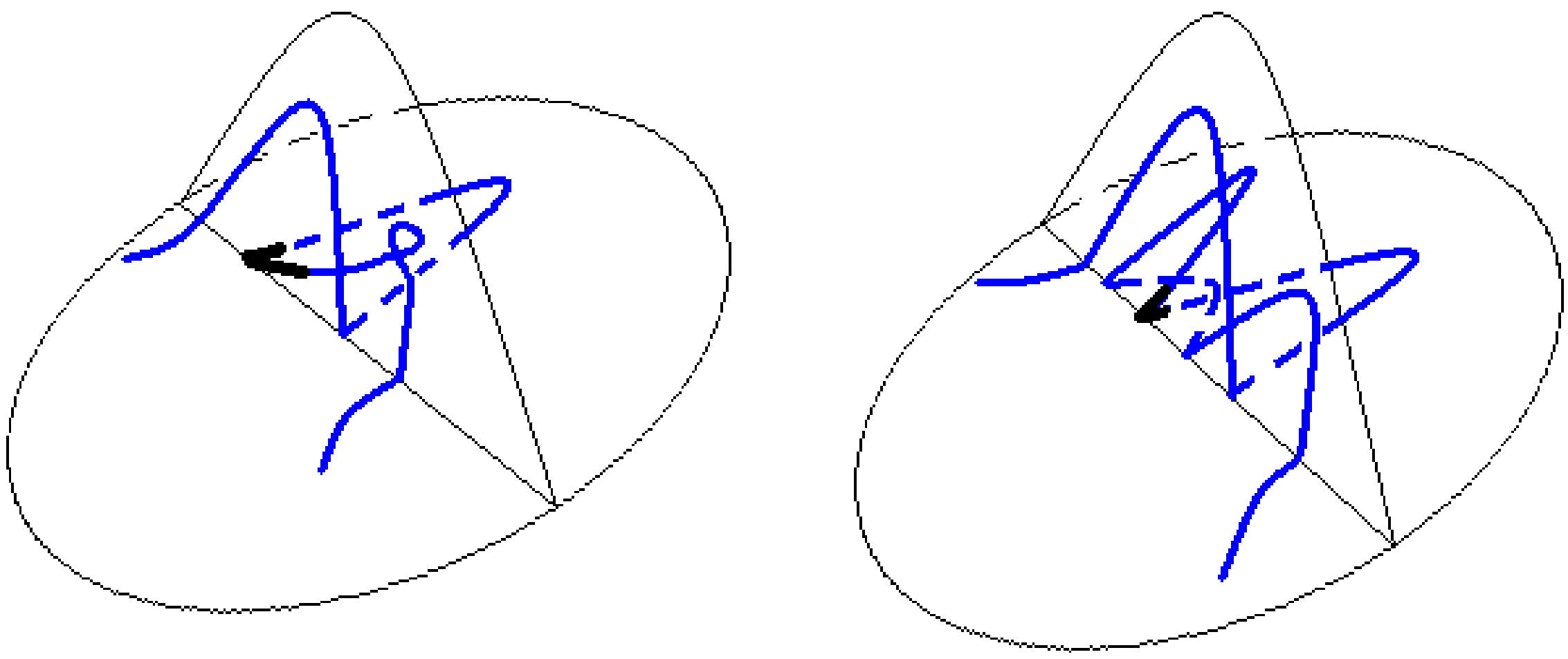

**Figure 77- Stack up the loops- preview- resolve loop in top component**

We decompose this step in Matveev moves. Note, that the sequence for the local part is the same as in the last step. Hence we only look at the global part:

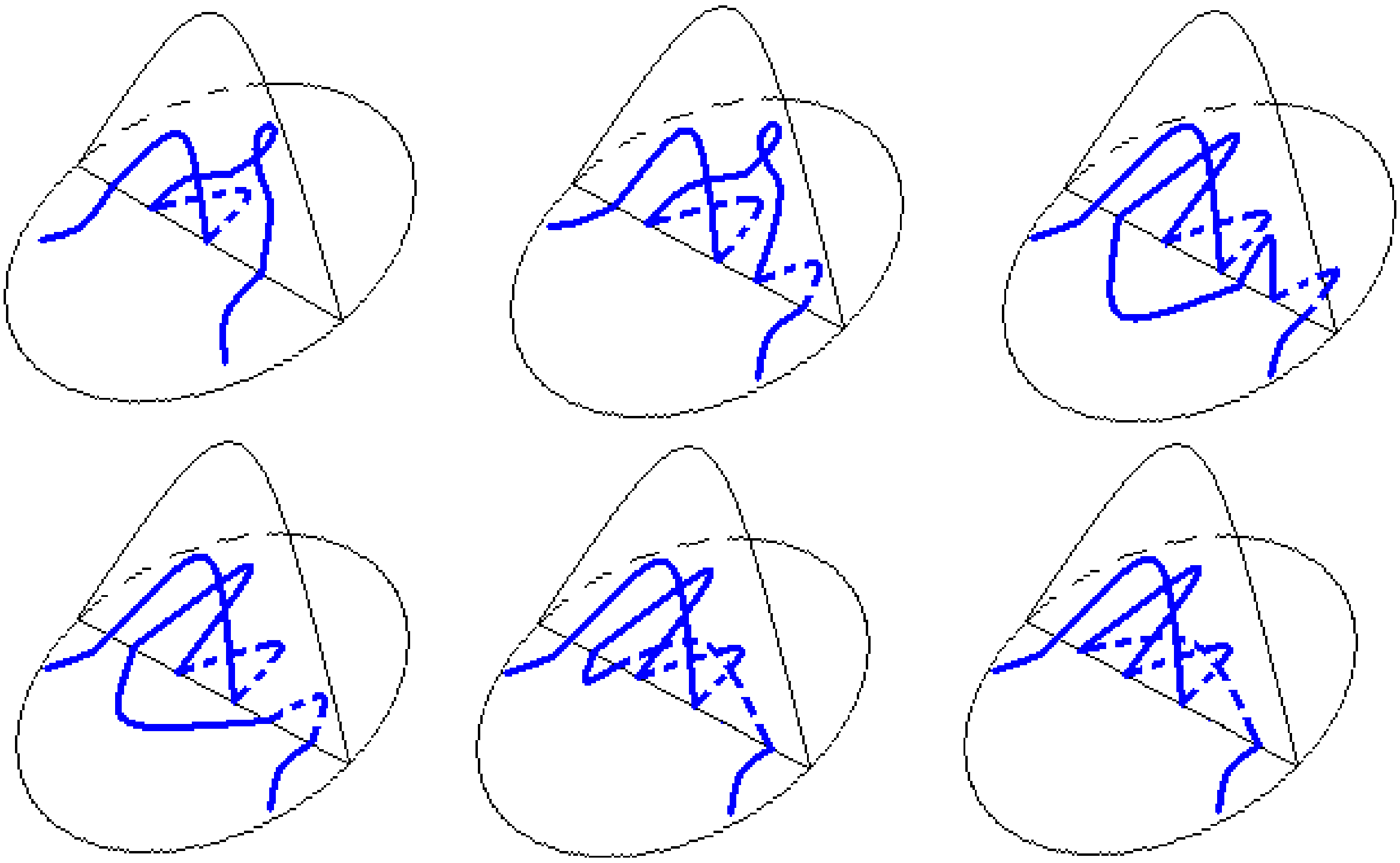

**Figure 78- Stack up the loops- resolve loop in top component- part 1**

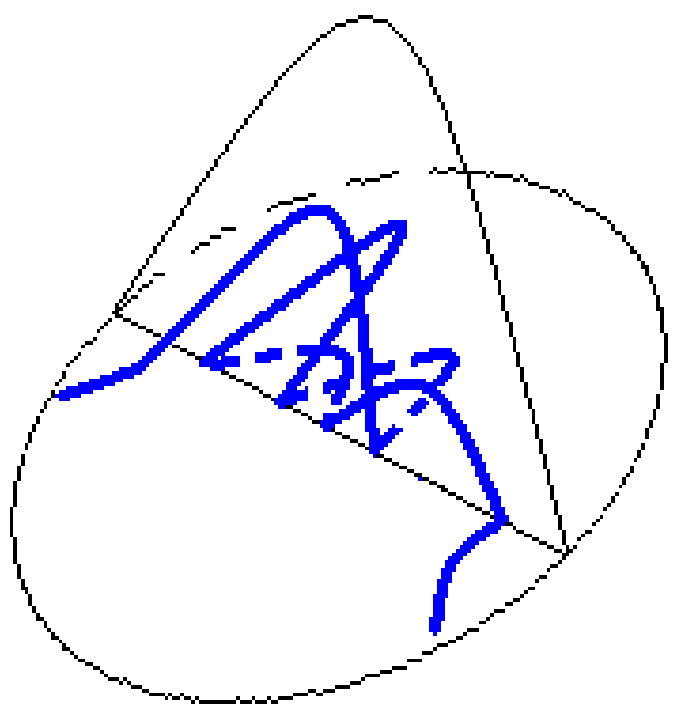

**Figure 79- Stack up the loops- resolve loop in top component- part 2**

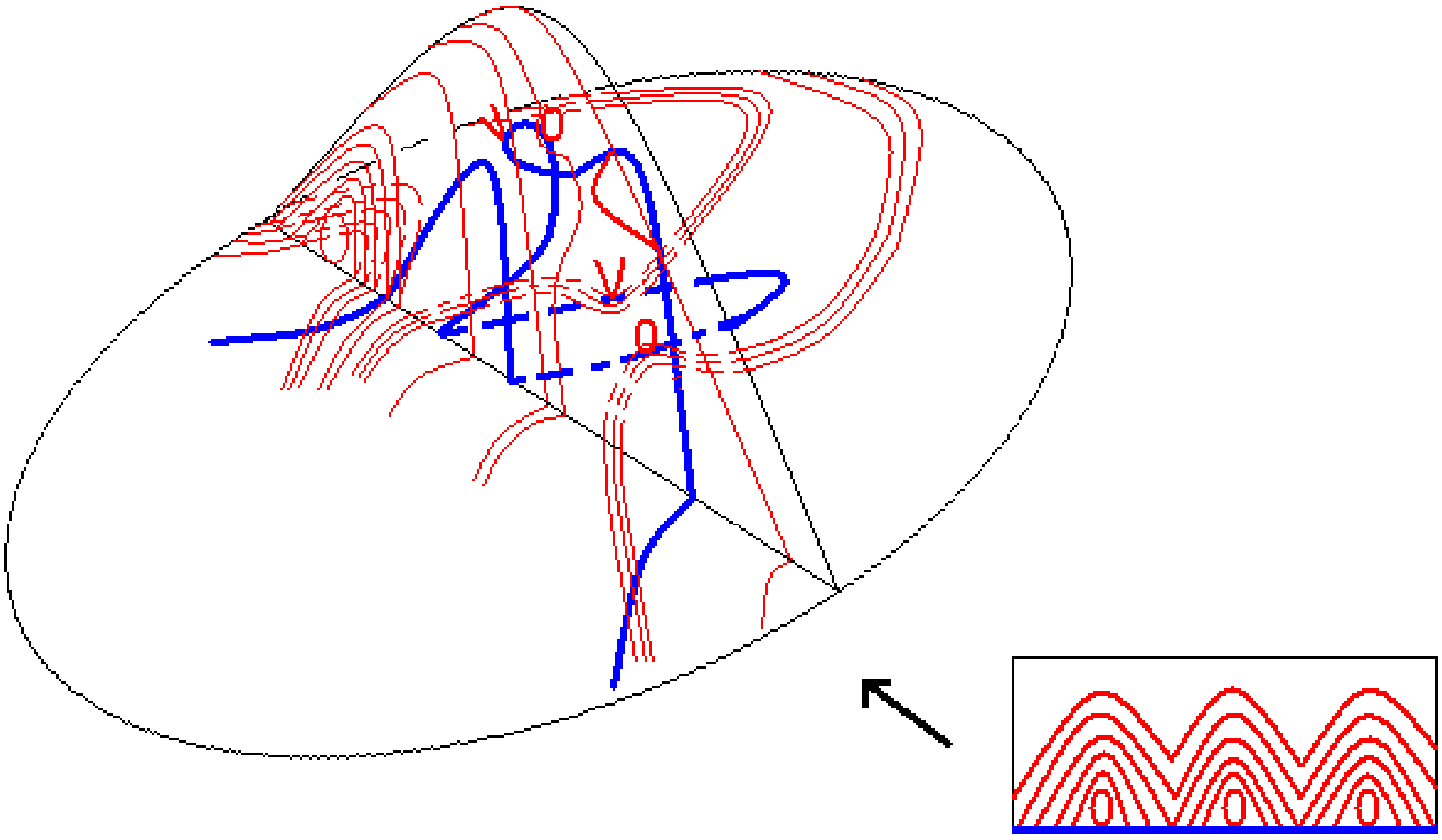

**Figure 80- Stack up the loops- resolve loop in top component- 1**

We generate a pair of saddleopoints in the frontside component on the right arc in order to prepare a $T_3$ move:

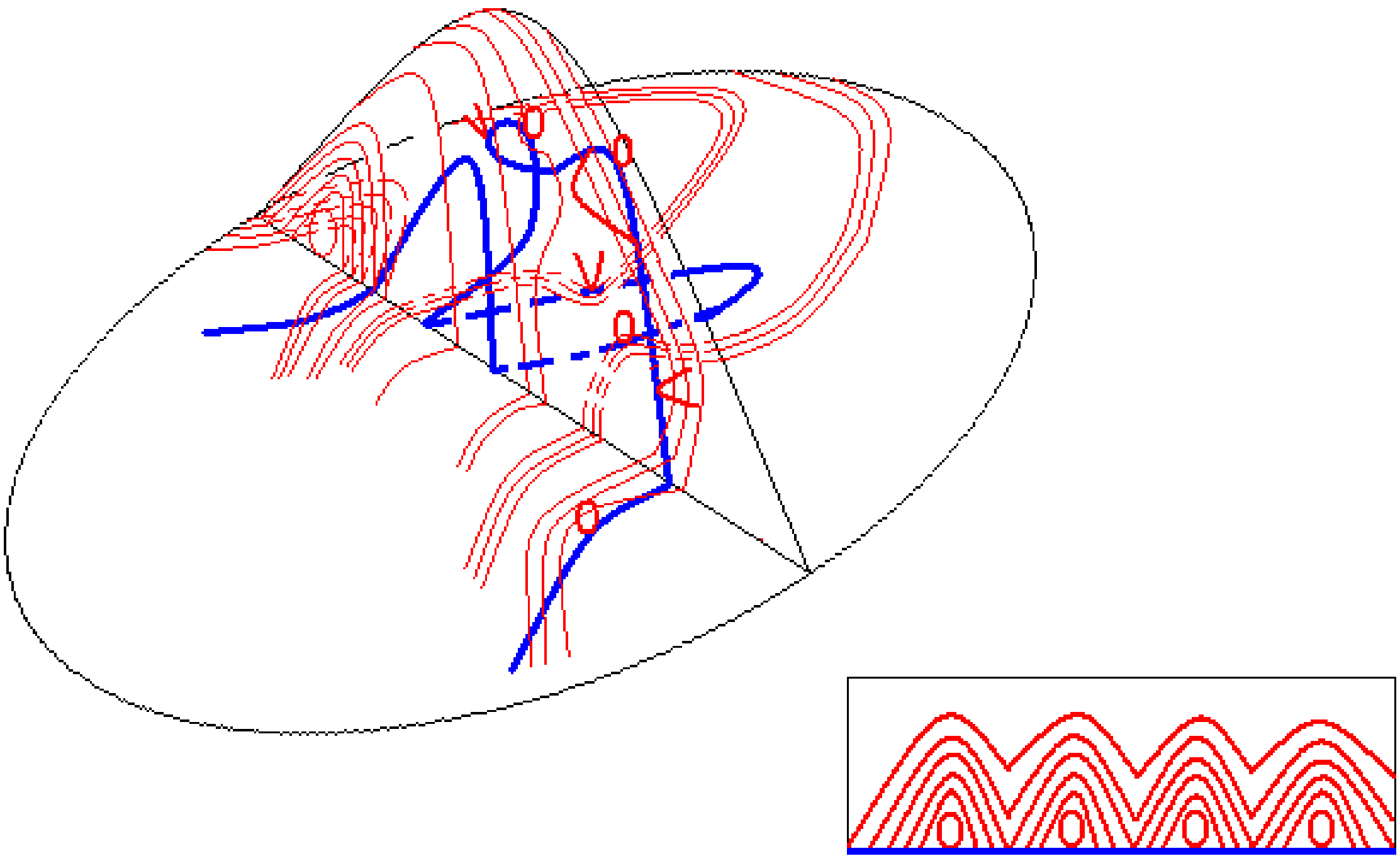

**Figure 81- Stack up the loops- resolve loop in top component- 2**

Now we can perform a “good $T_3$ turn” in the backside component for that arc:

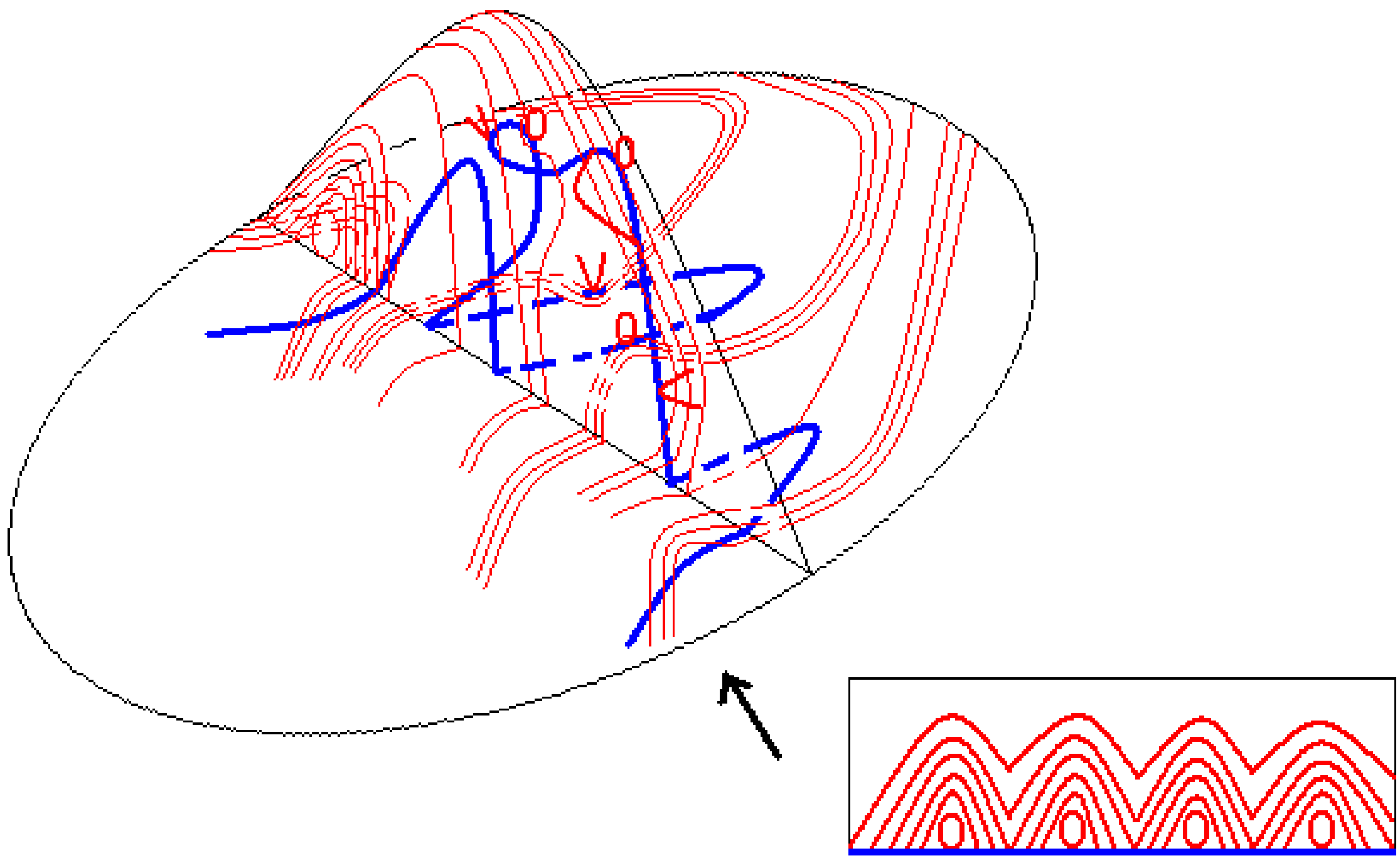

**Figure 82- Stack up the loops- resolve loop in top component- 3**

We cancel the pair of saddlepoints on the right arc in the top component and perform $T_2$ moves.
(we summarize here more than one $T_2$ move) to resolve the loop. We remark, that for the arising $T_3$ turn in top component we get from the former figure that the slices are “good”:

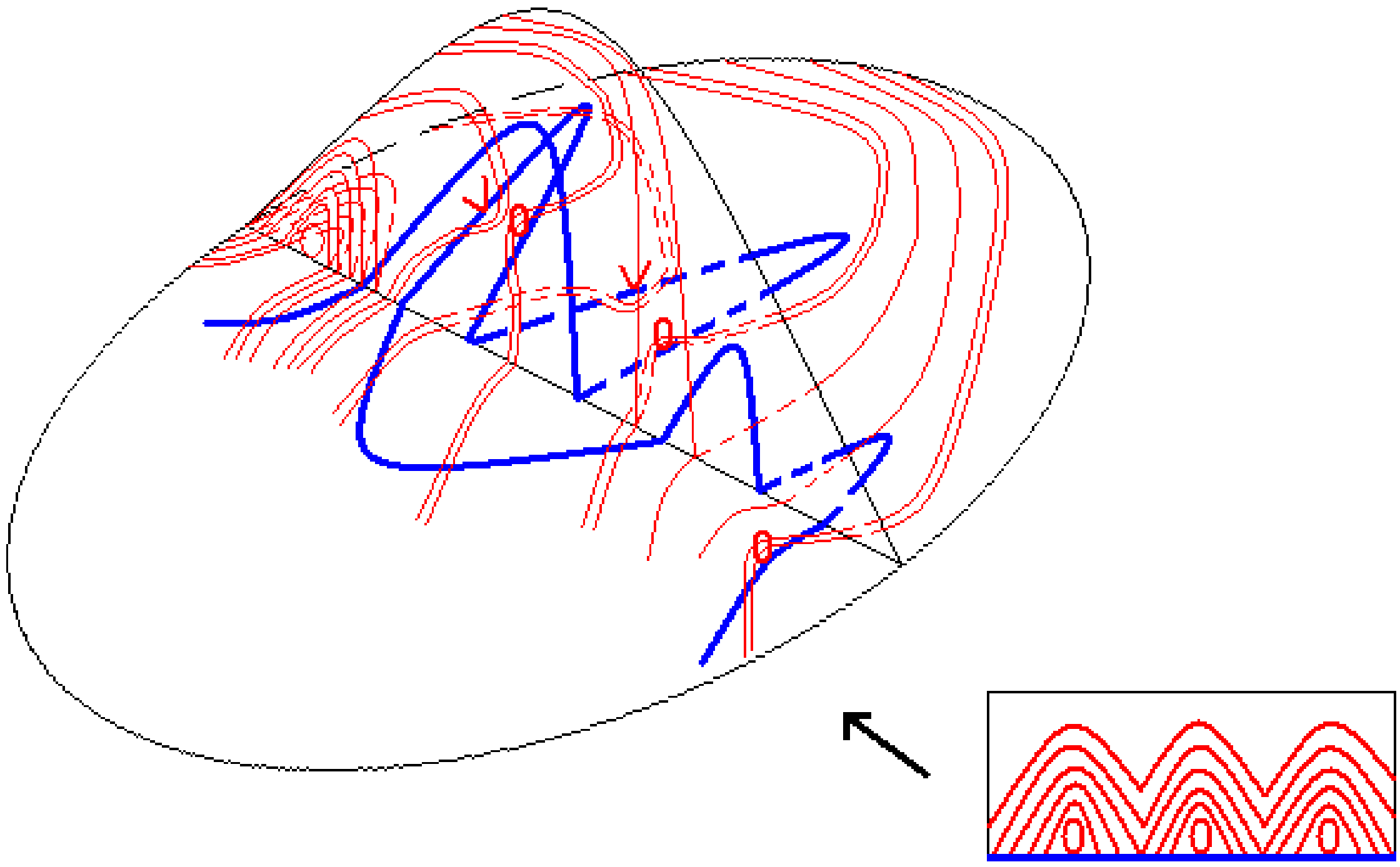

**Figure 83- Stack up the loops- resolve loop in top component- 4**

We annihilate the right $T_3$ turn on the top component by $T_3^{-1}$ and remark, that the $T_3$ turn arising during the annihilation in the frontside component by this step is a “good” one, which comes from the choice of slices in the step before:

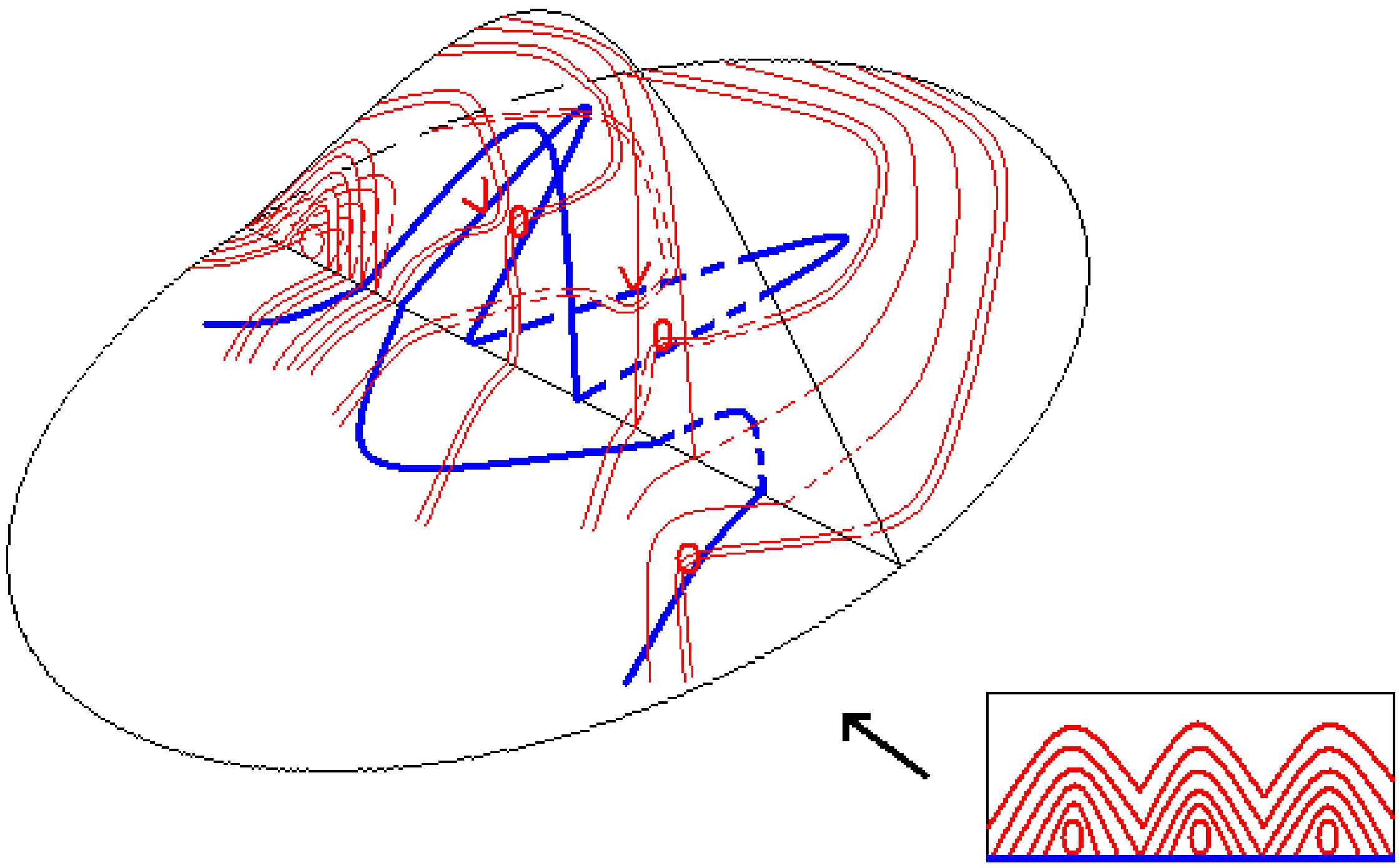

**Figure 84- Stack up the loops- resolve the top component- 5**

We perform the move T* 2 times on the $T_3$ turn in the frontside component and the “good $T_3$ turn” stays good:

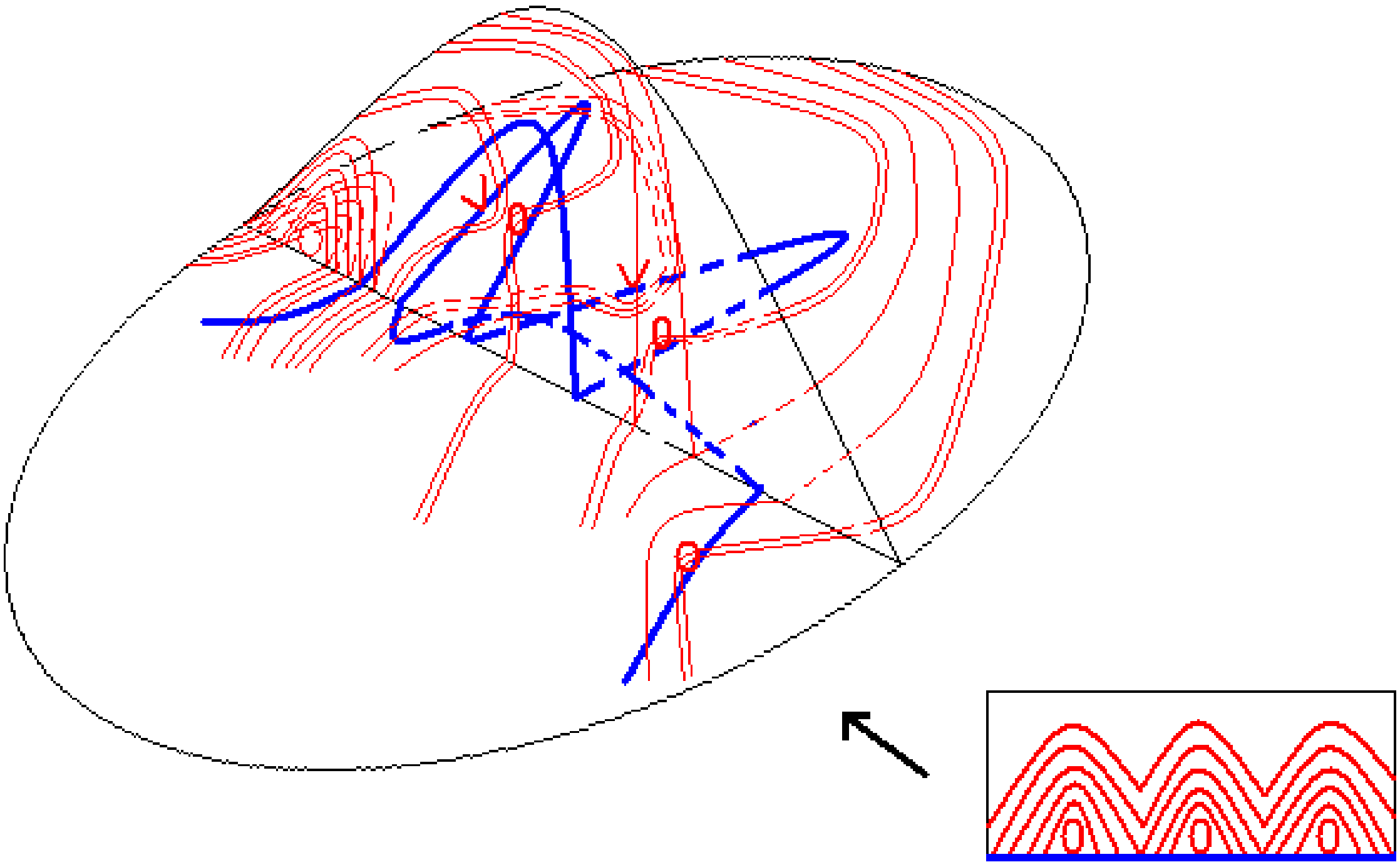

**Figure 85- Stack up the loops- resolve the top component- 6**

We annihilate that $T_3$ turn in the frontside component by $T_3^{-1}$ and remark, that the new $T_3$ turn in the backside component is “good” again:

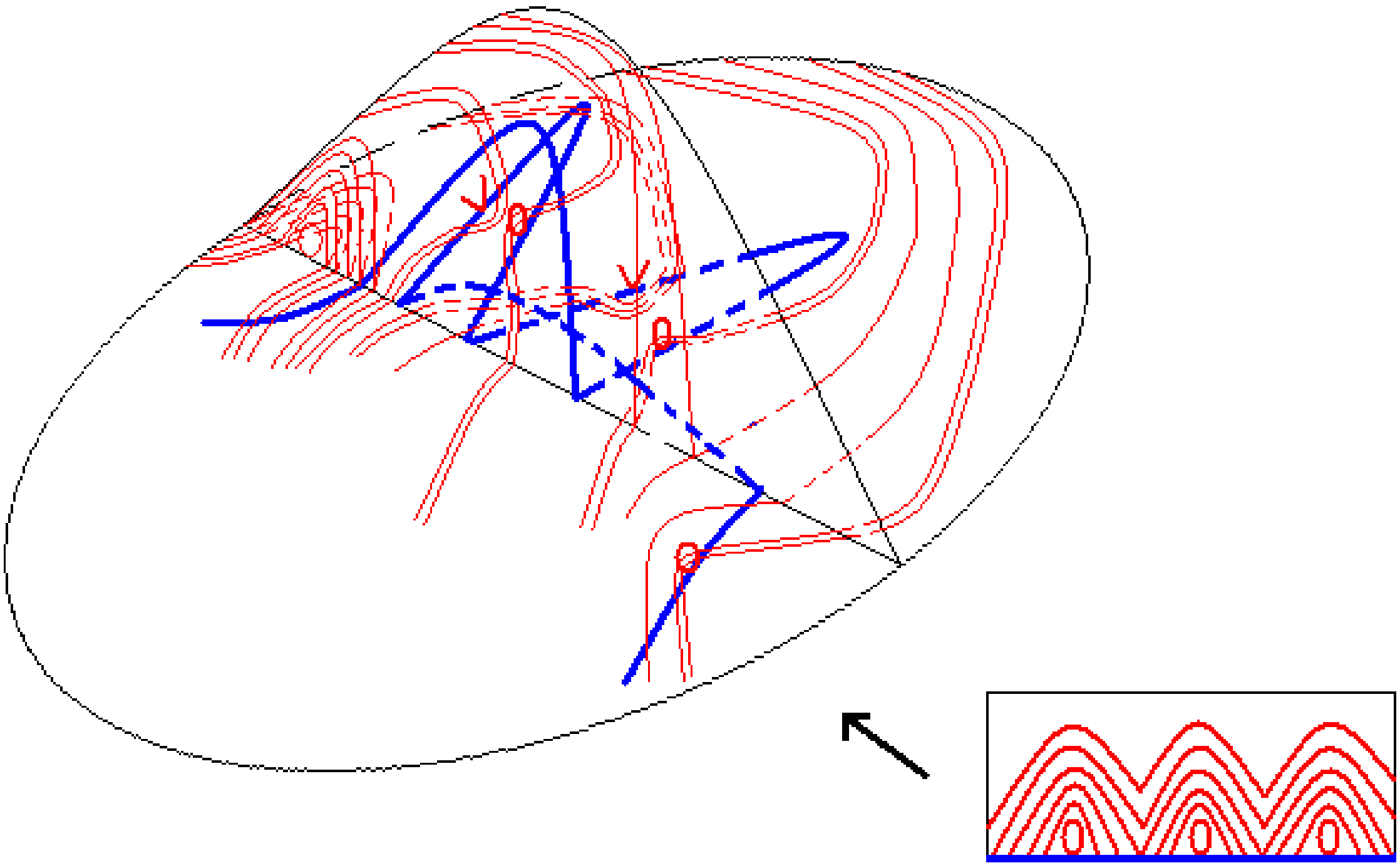

**Figure 86- Stack up the loops- resolve loop in top component- 7**

We take the arc of that $T_3$ turn in the backside component and perform a modified $T_2$ move (the modifications will be discussed in chapter 6.2). Again the new $T_3$ turn on the top component defines a “good “one, the slices for the connected arc in the backside component are unchanged:

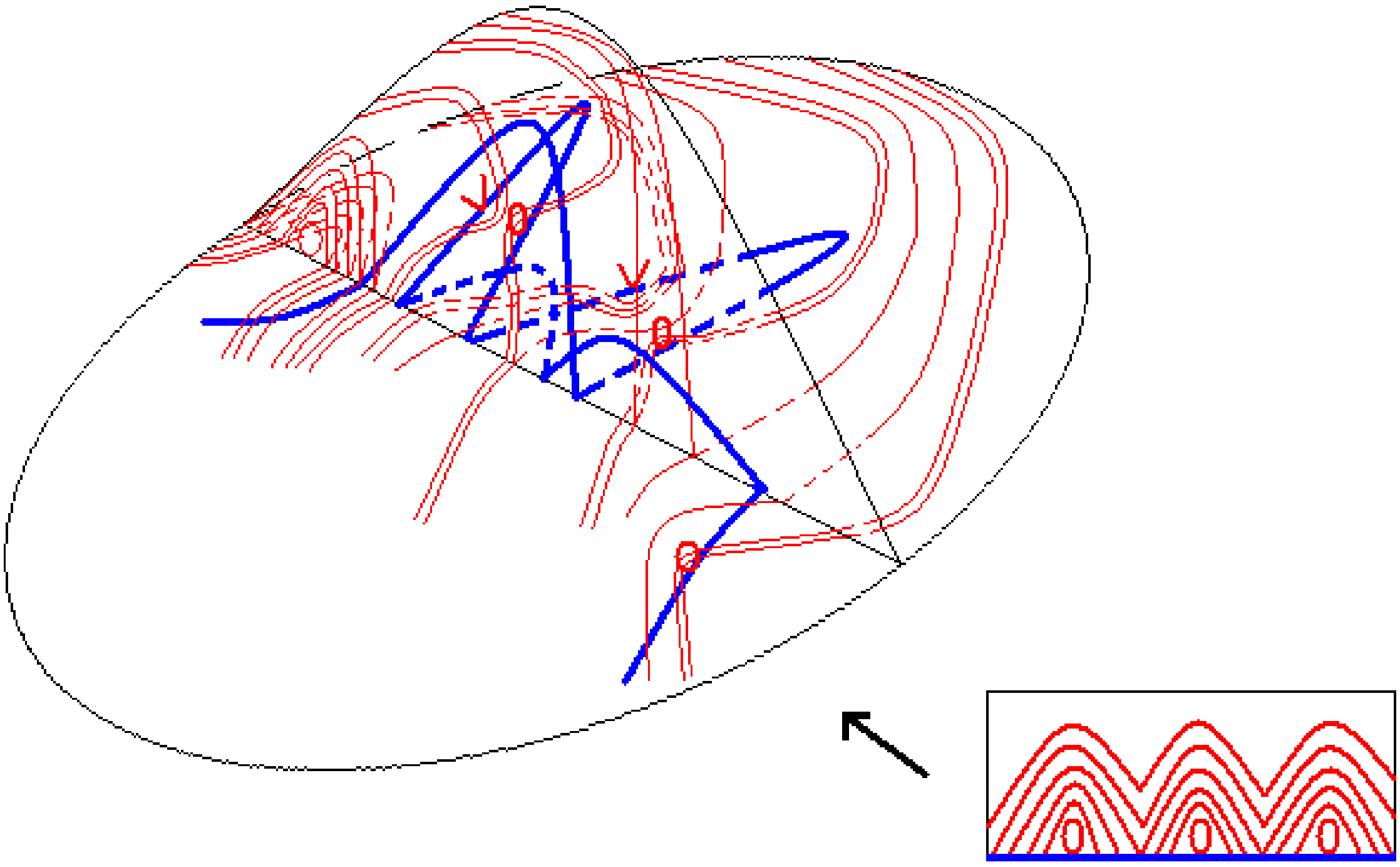

**Figure 87- Stack up the loops- resolve loop on top component- 8**

To finish that sequence, we generate a pair of saddlepoint on the right arc, one on the frontside and one on the top component

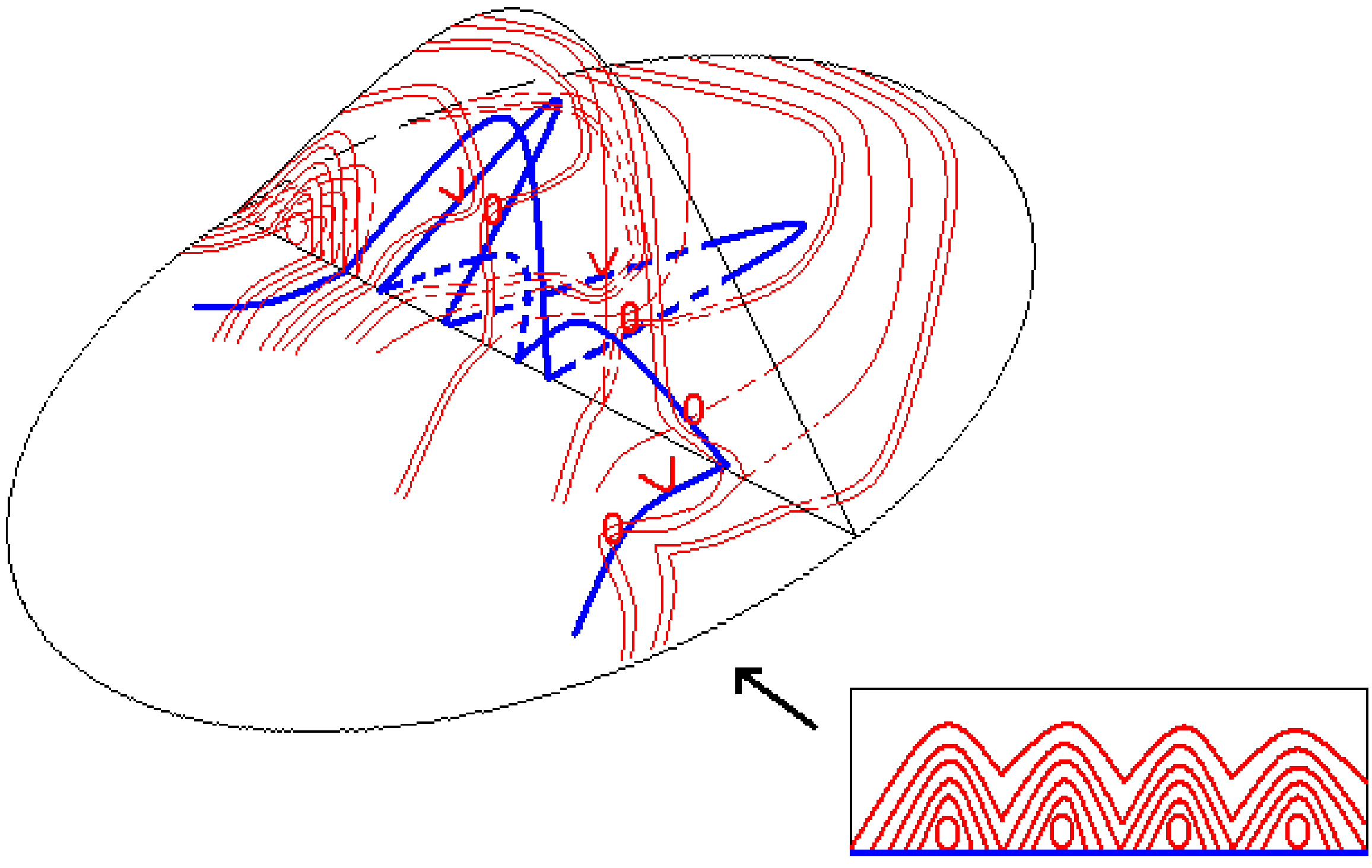

**Figure 88- Stack up the loops- resolve loop in top component- 9**

and cancel the pair of saddlepoints on the right arc in the frontside component:

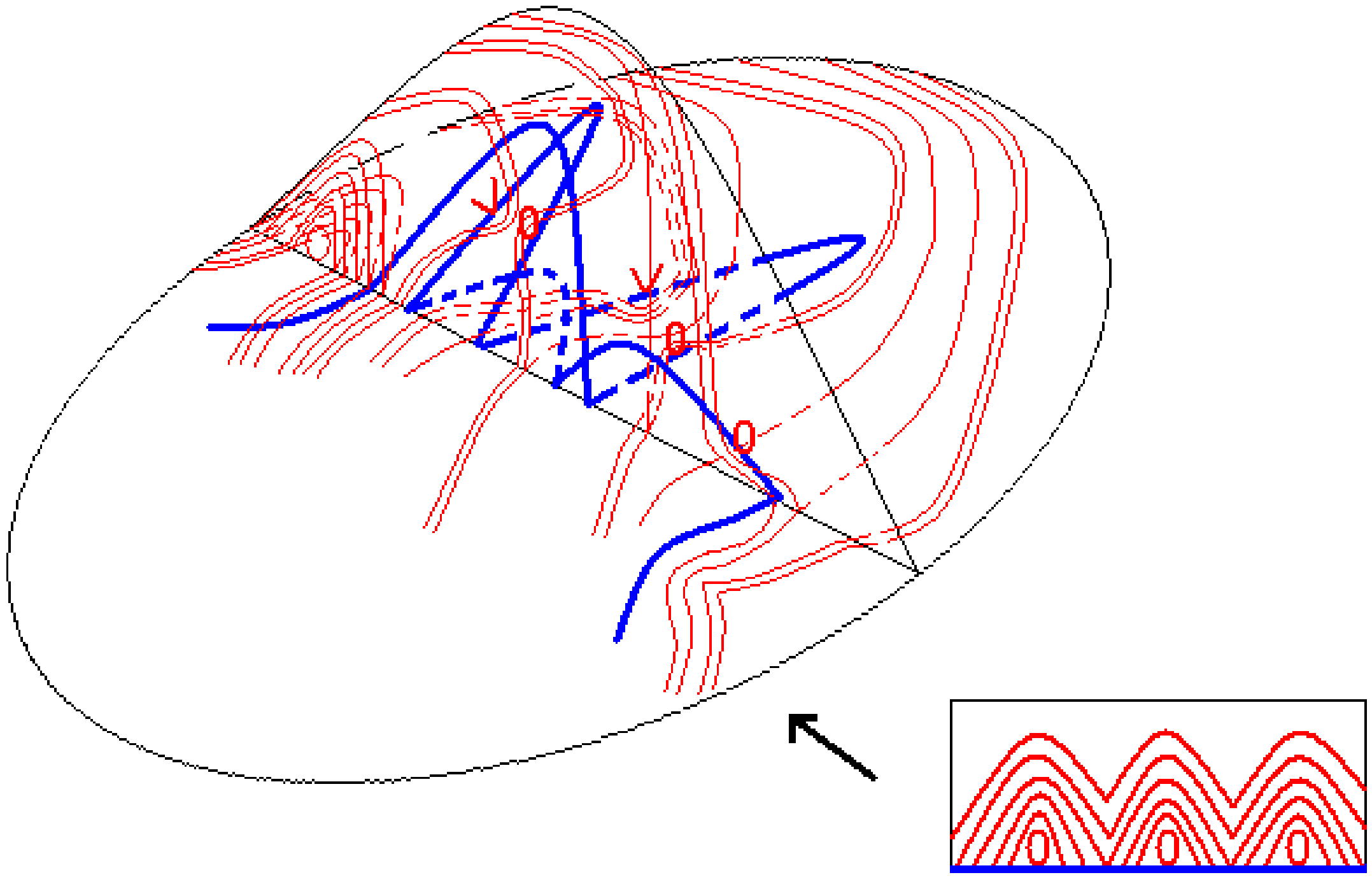

**Figure 89- Stack up the loops- resolve loop in top component- 10 (end)**

By that we have annihilate our introduced pair of saddlepoints at the beginning.

### 5.1.2.2 Resolve the stacked up loops

We give a preview to the next steps. At first push the left big $T_3$ turn in the top component from left to right as in the first three figures. In the third figure left is a loop, composed from parts in top and in the backside component. Move the part in the top component into the backside component, to have the loop entirely in that component, as in the fourth figure. Annihilate the left $T_3$ turn in the top component:

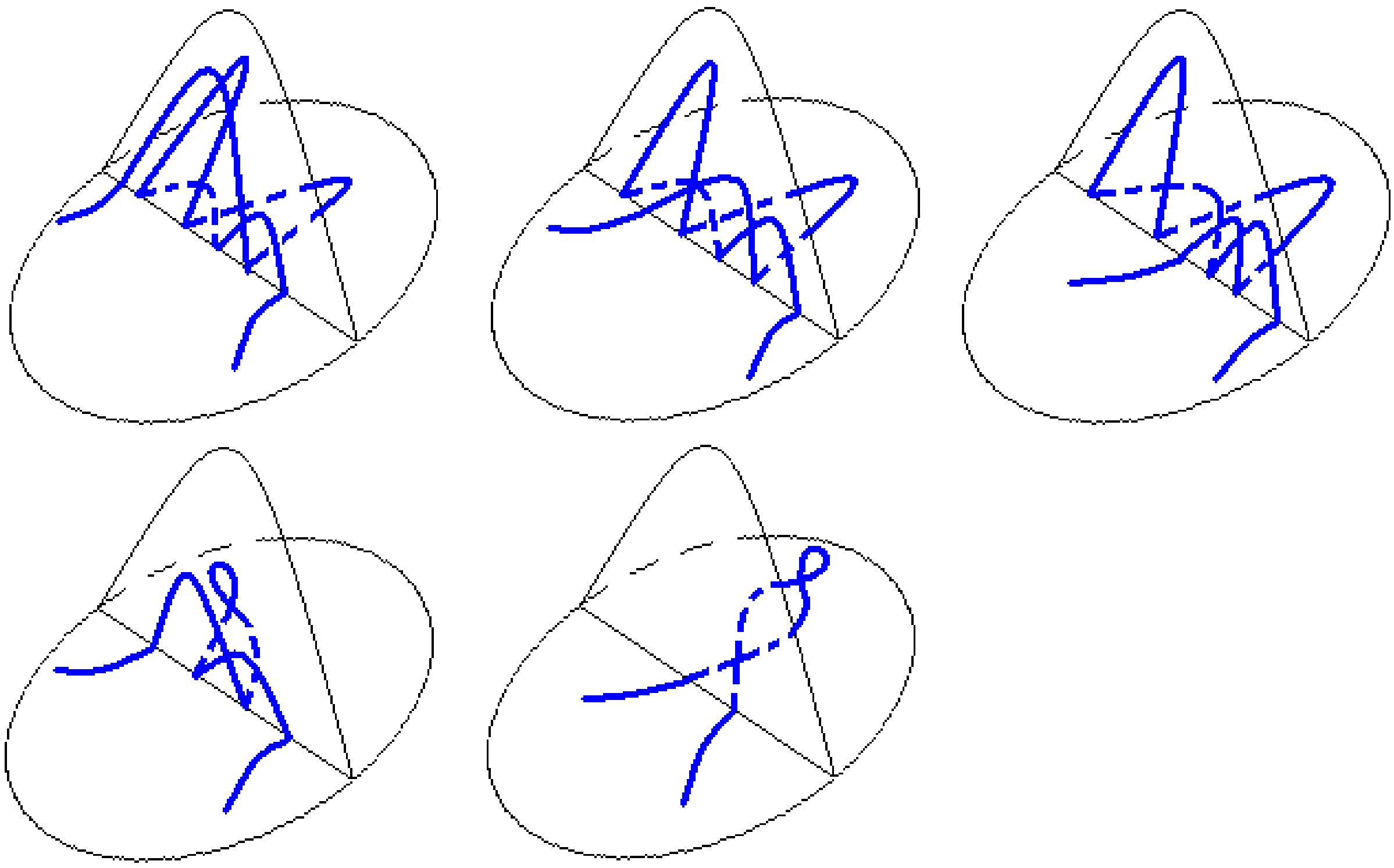

**Figure 90- Resolve the stacked up loops- preview**

We decompose the move of the left big $T_3$ turn in the top component:

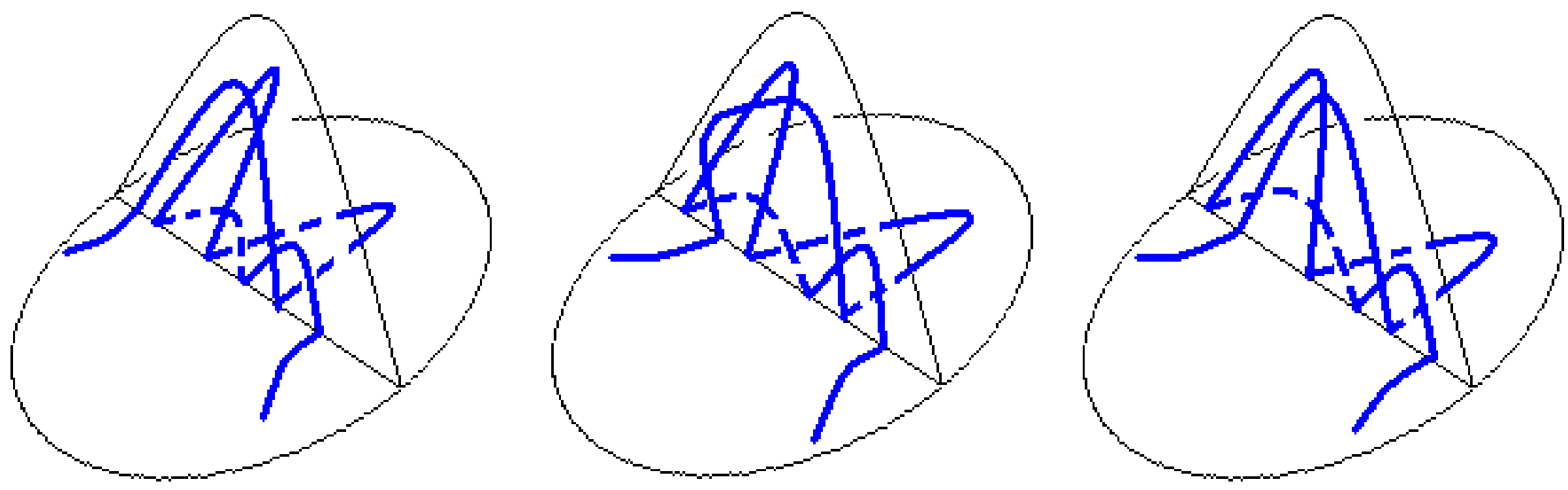

**Figure 91- Resolve the stacked up loops- preview- move of big left $T_3$ turn**

The idea to arrange the sequence of slices to perform these steps is to introduce an extremum on the arc of the big $T_3$ turn to be able to pass the given extremum of the second arc in the top component:

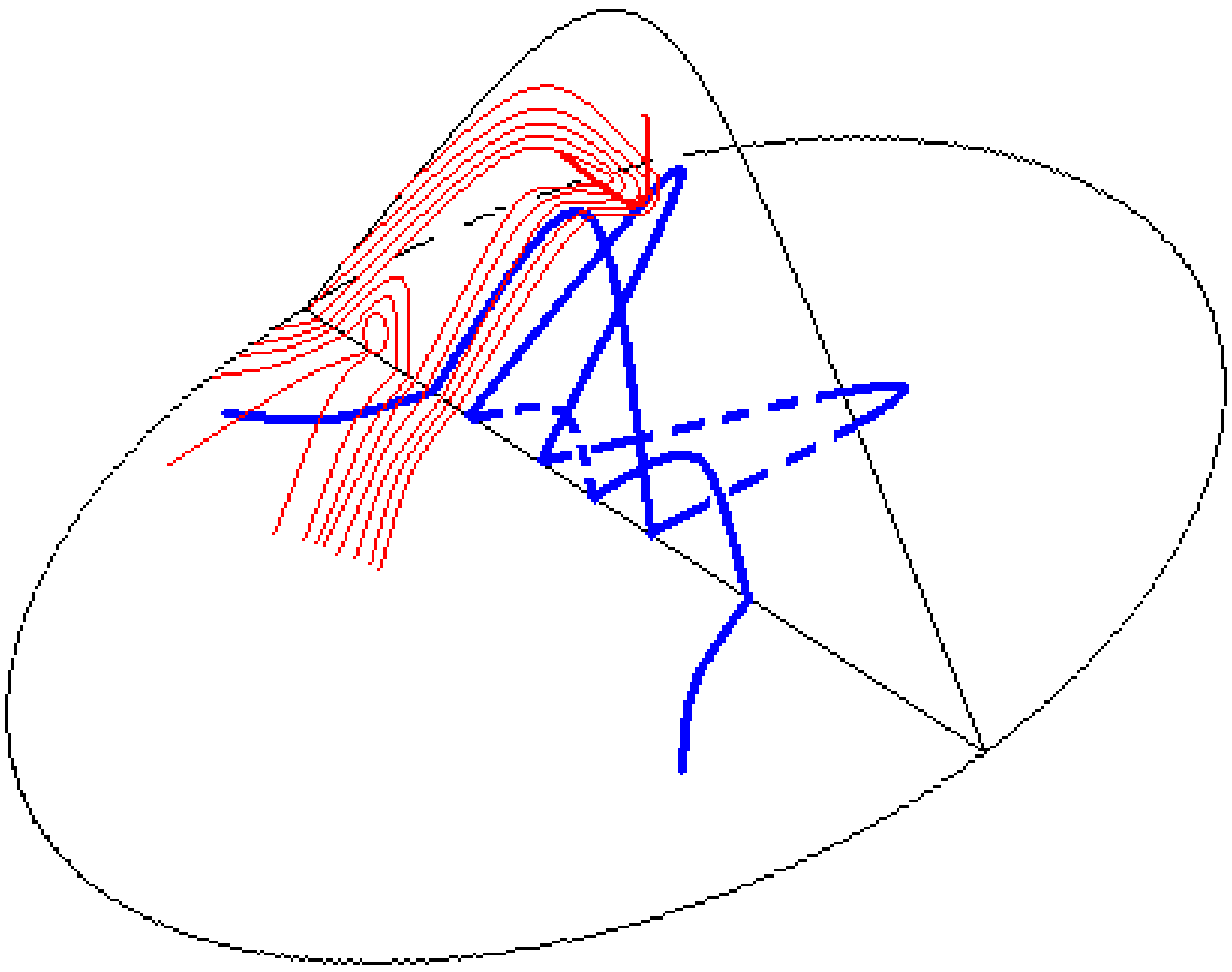

**Figure 92- Resolve the stacked up loops- move of big left $T_3$ turn- 1**

Clearly we have to generate a pair of saddlepoints on the big $T_3$ turn line:

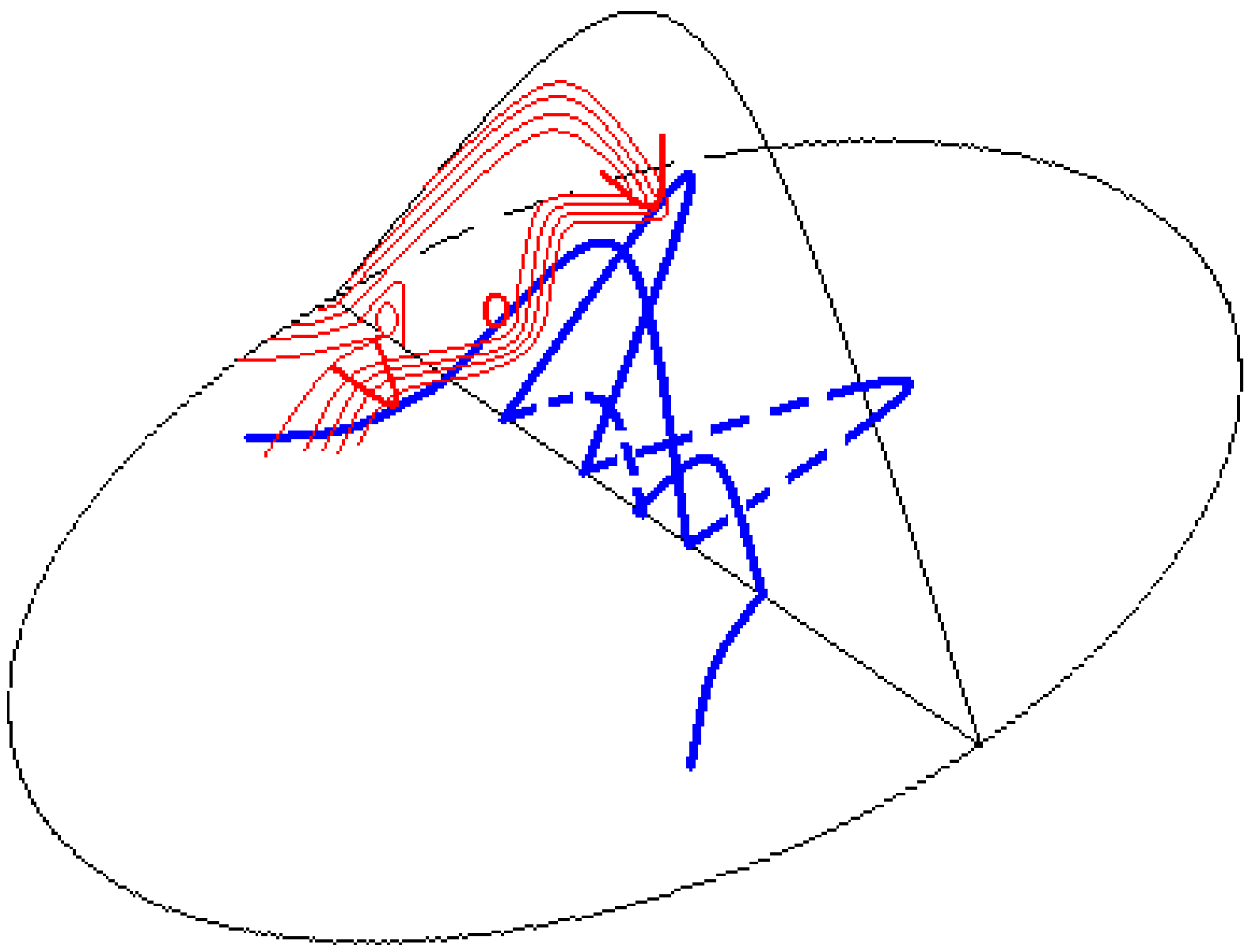

**Figure 93- Resolve of stacked up loops- move of big left $T_3$ turn- 2**

Now we can perform the move:

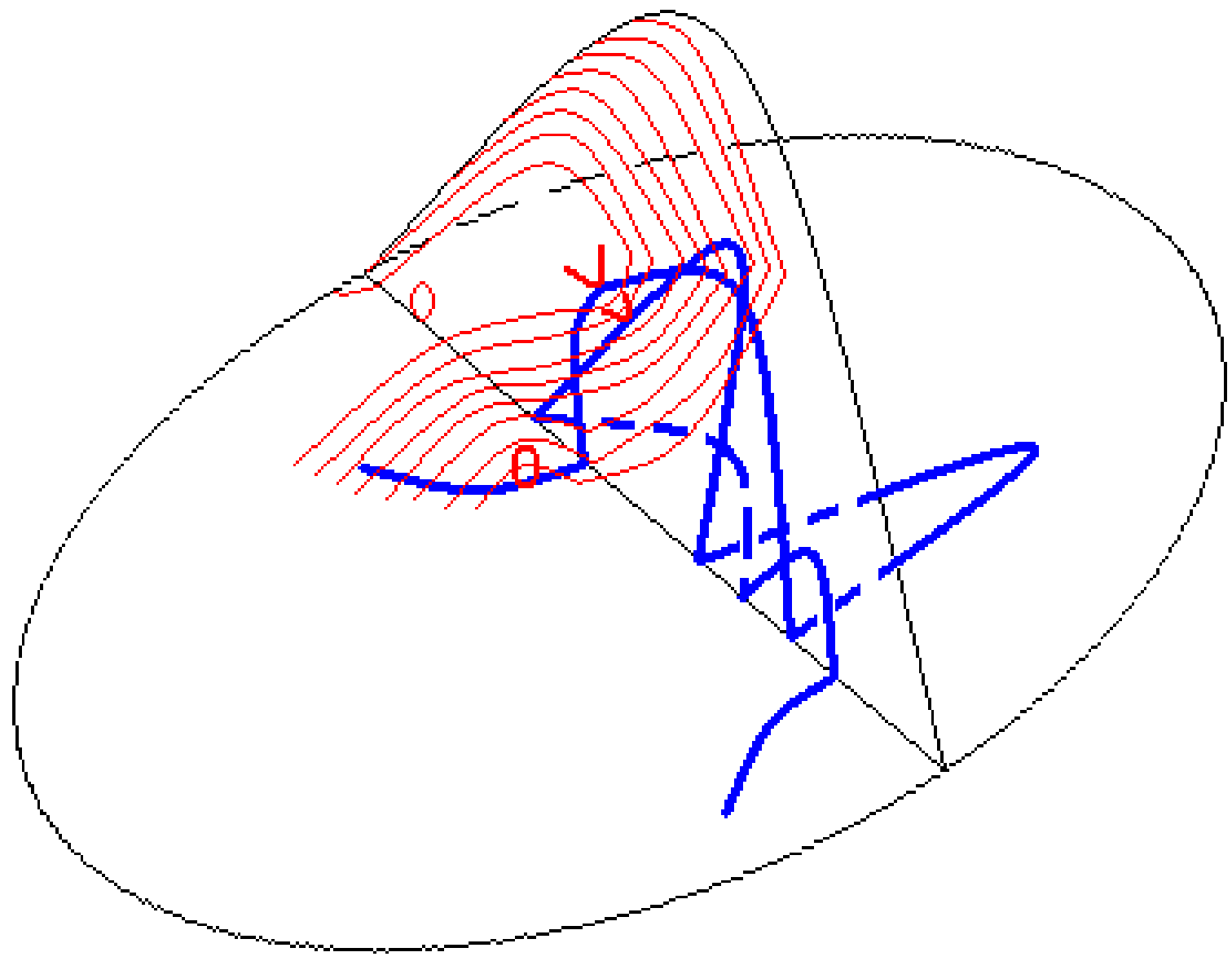

**Figure 94- Resolve the stacked up loops- move of big left $T_3$ turn- 3**

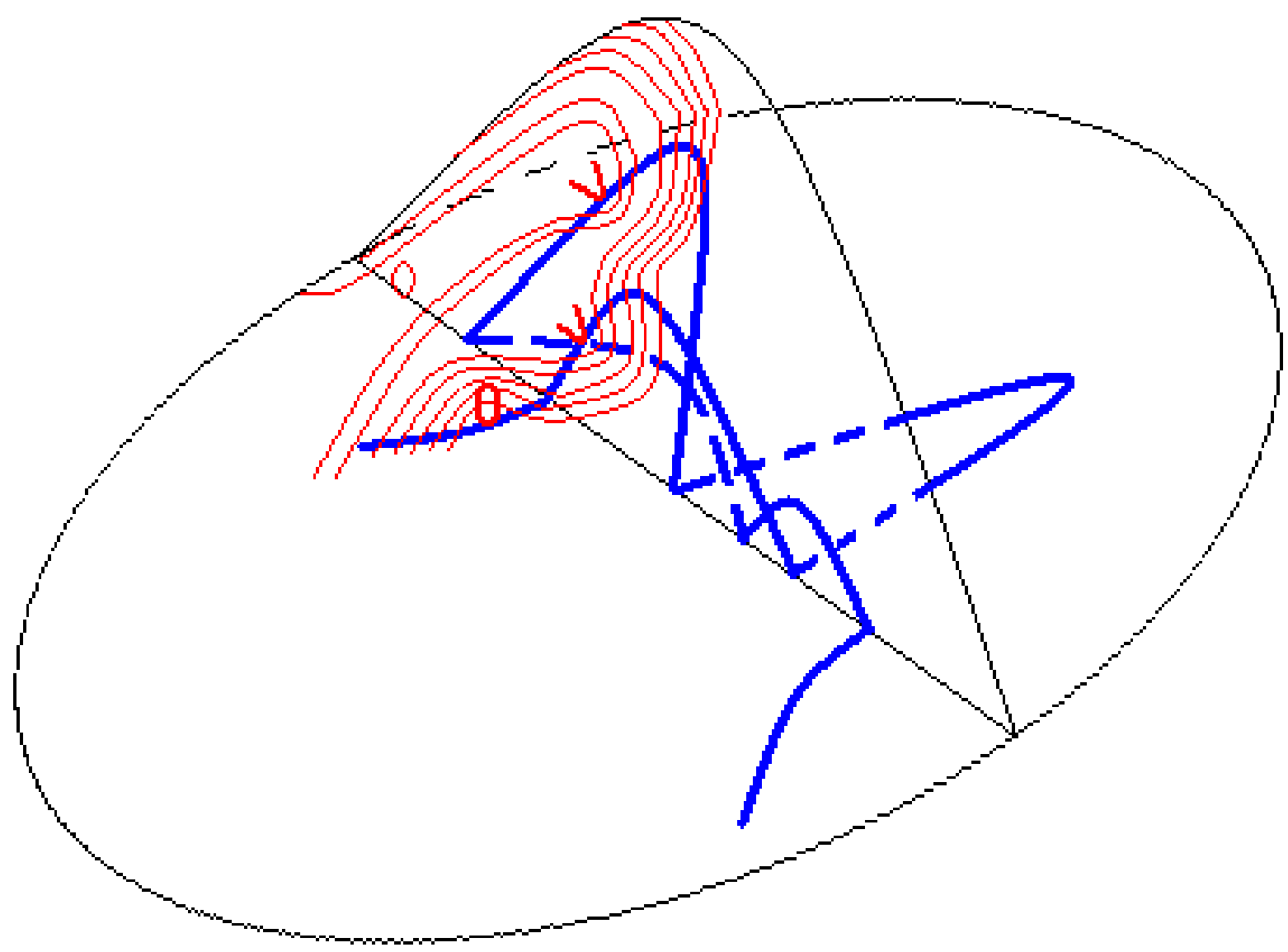

**Figure 95- Resolve the stacked up loops- move of big left $T_3$ turn- 4**

It remains to cancel the introduced pair of saddlepoints:

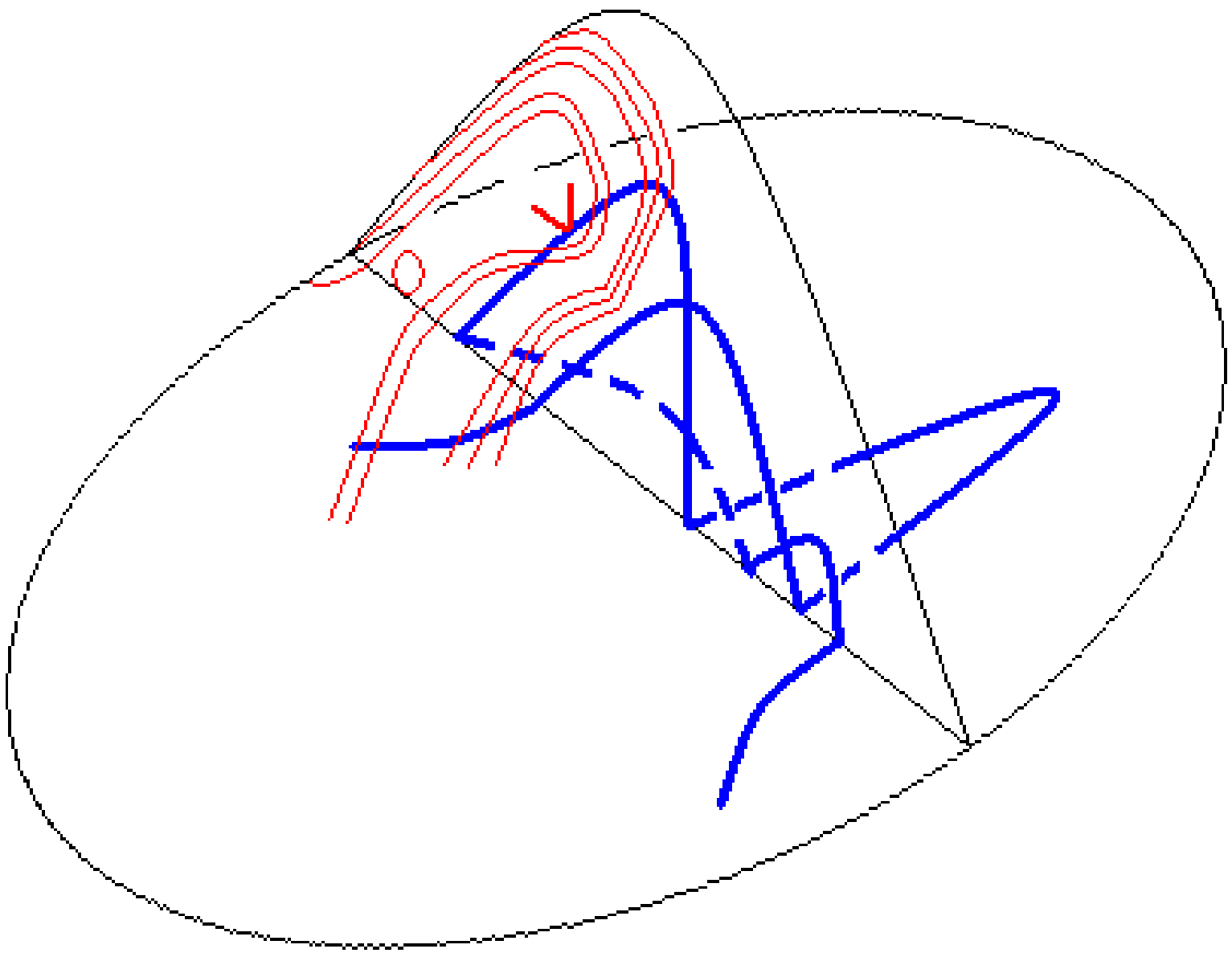

**Figure 96- Resolve of the stacked up loops- move of big left $T_3$ turn- Figure 5 (end)**

To perform the next step, repeat the construction above, but also introduce a pair of saddlepoints on the arc which will be passed and annihilate them again:

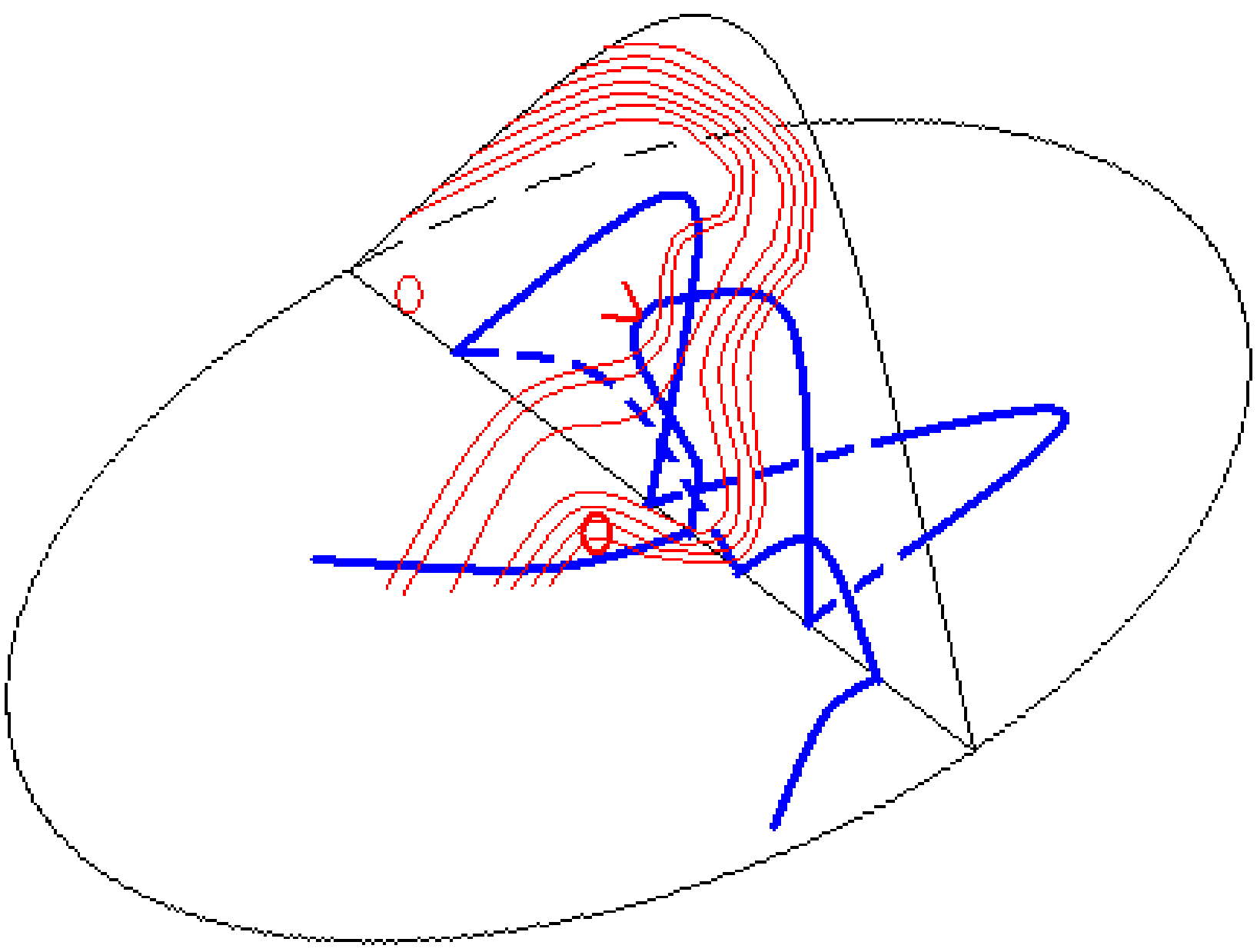

**Figure 97- Resolve the stacked up loops- continue move of big left $T_3$ turn- 1**

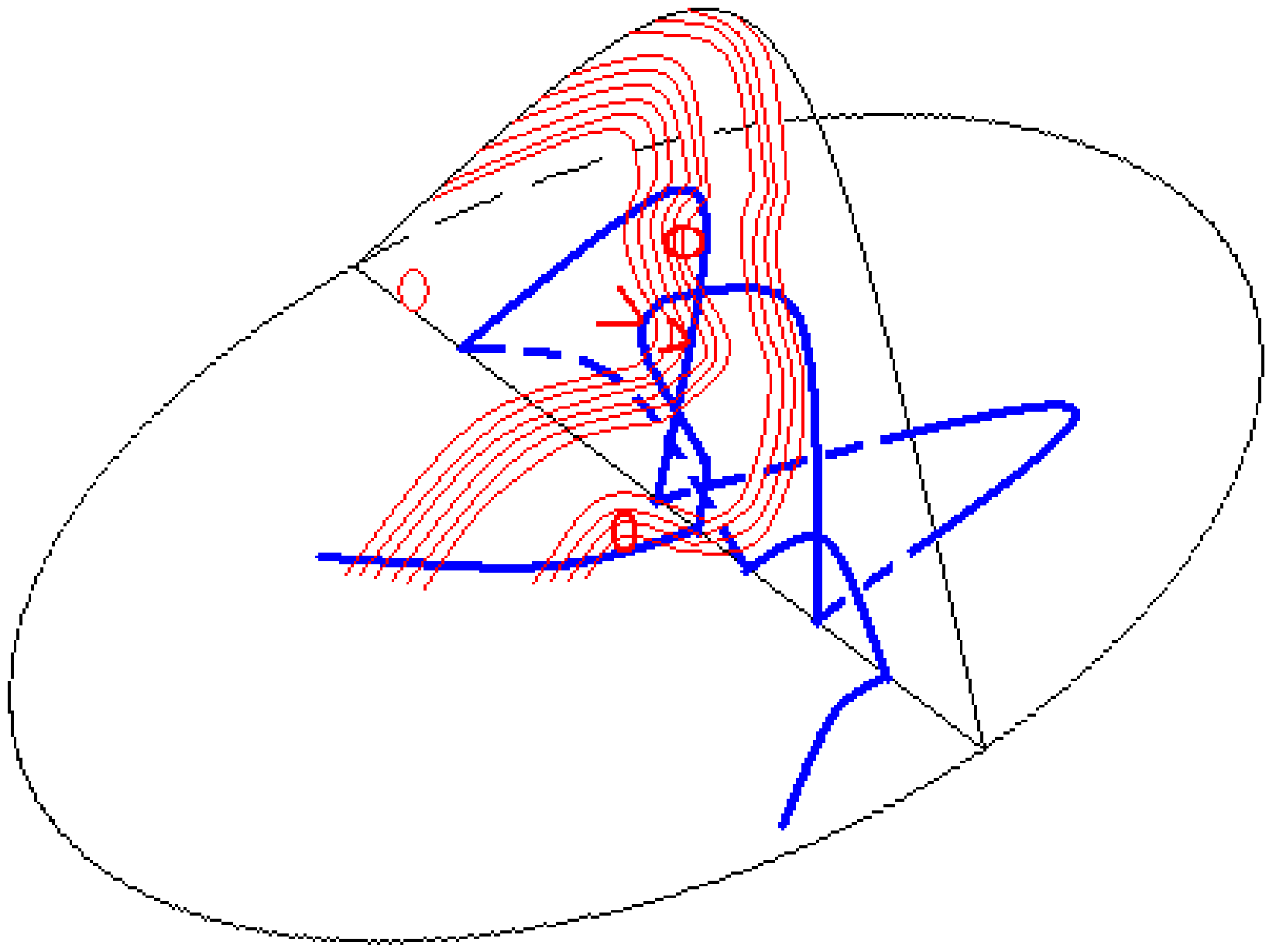

**Figure 98- Resolve the stacked up loops- continue move of big left $T_3$ turn- 2**

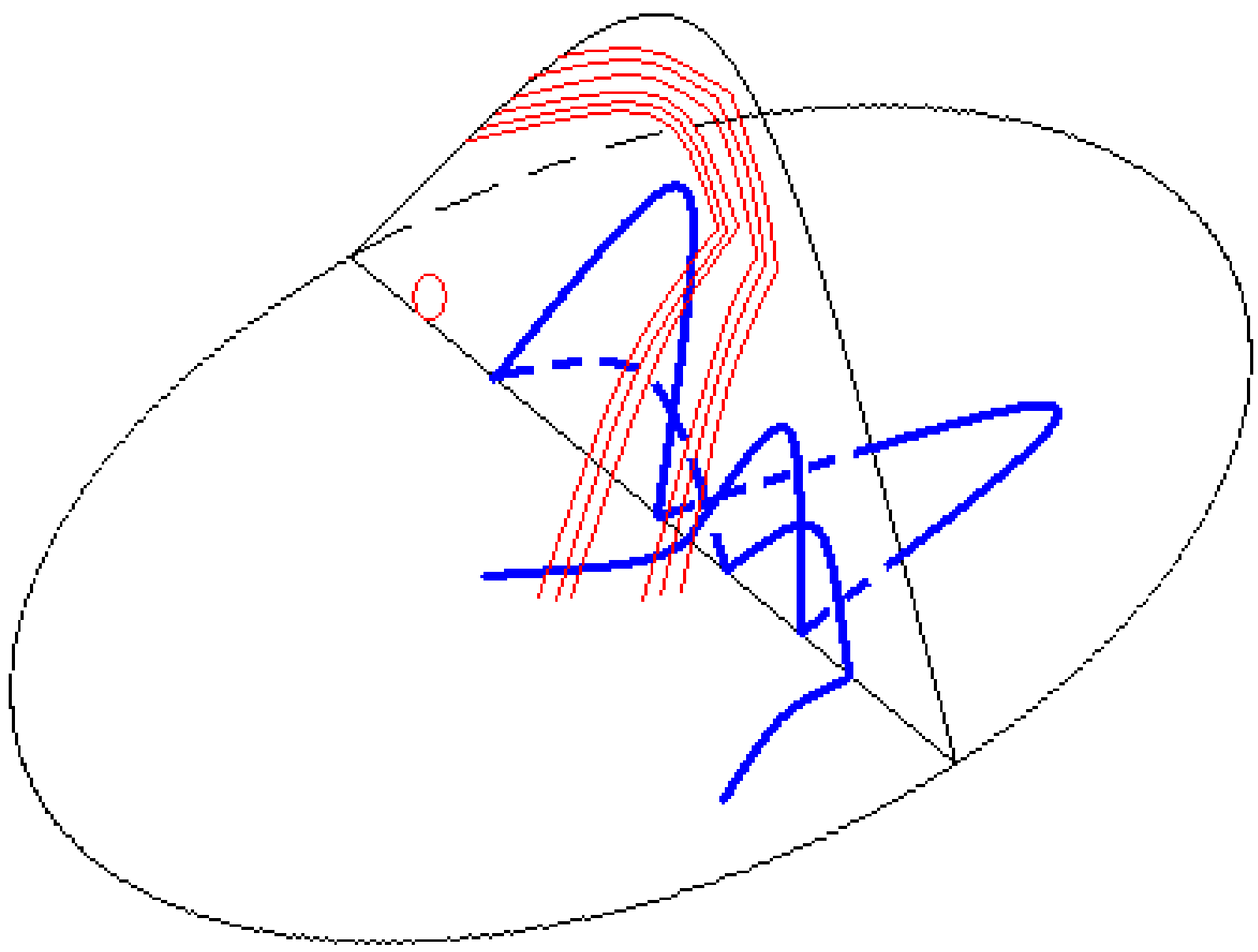

**Figure 99- Resolve the stacked up loops- continue move of big left $T_3$ turn- 3**

We give a preview of the next step but we omit further details about the sequence of slices because the construction is similar to our previous considerations . We do not perform the last transition, because we would get a new $T_3$ turn with line segment slices, so we have to be very carefully:

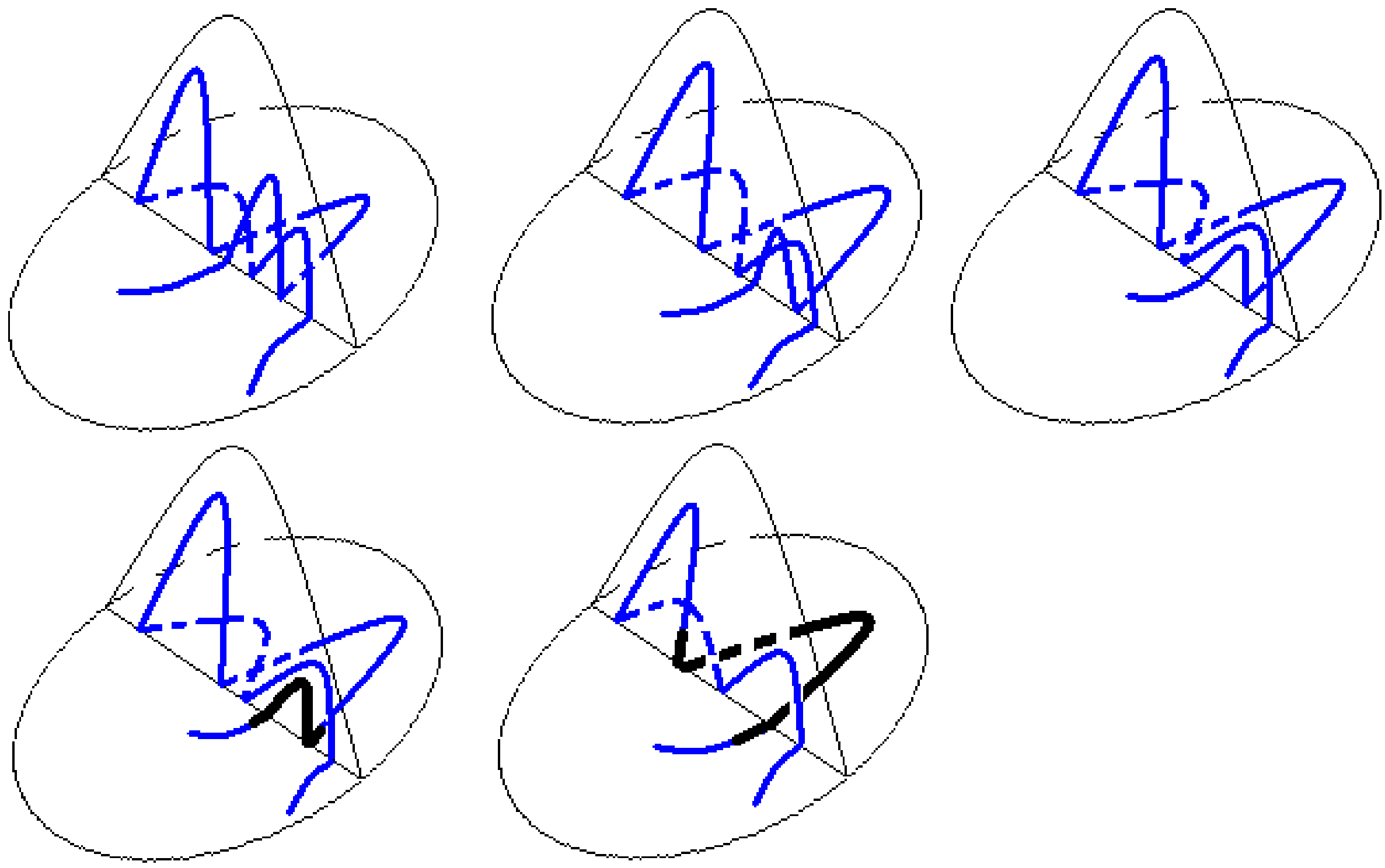

**Figure 100- Resolve the stacked up loops- preview- slide left $T_3$ turn under right $T_3$ turn**

When we perform the previous step, the result is:

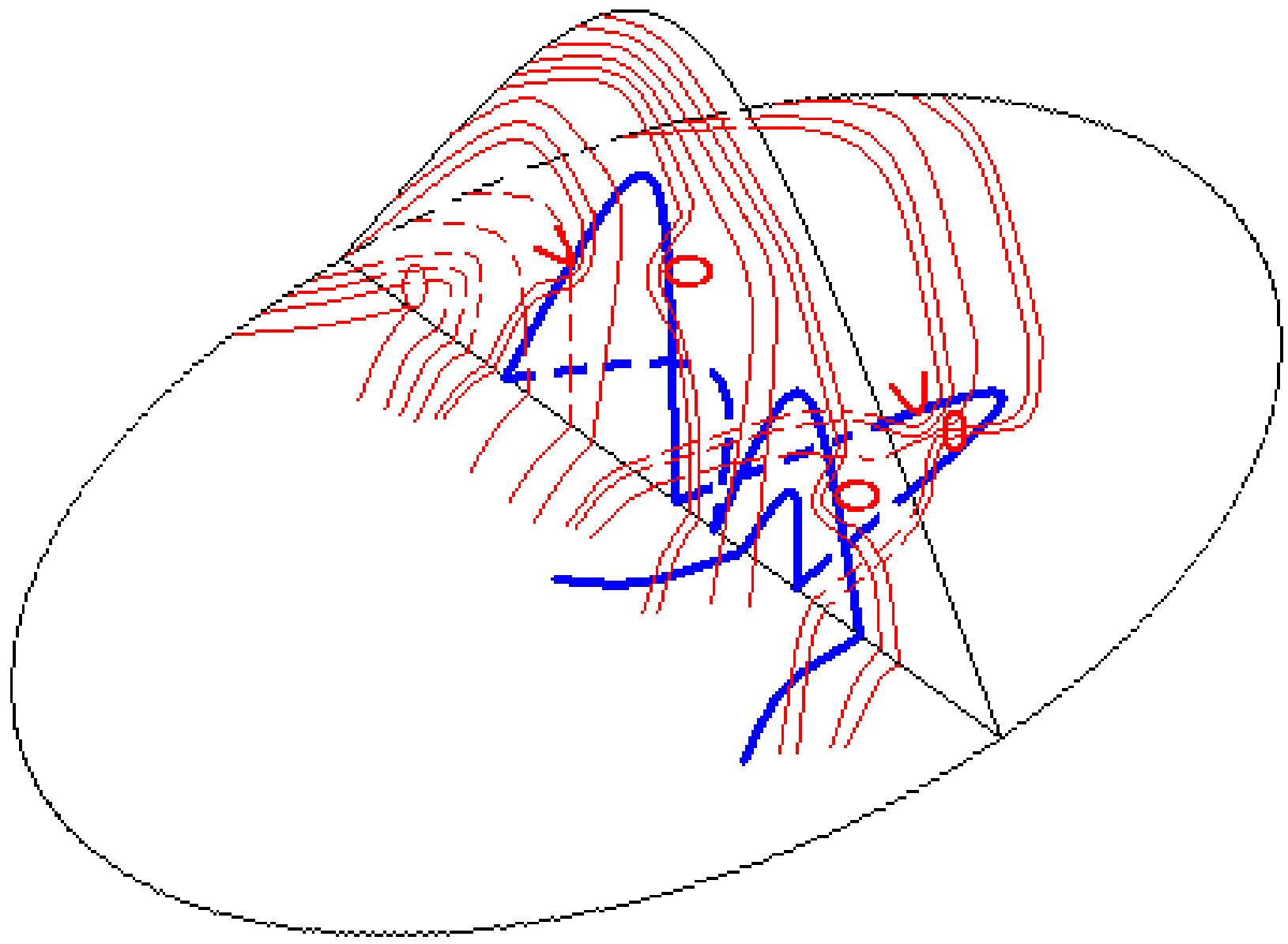

**Figure 101- Resolve the stacked up loops- slide left $T_3$ turn under right T3 turn- end**

We want to push the arc in the top component into the backside component to have the entire loop in that component. Furthermore we restrict our slices to the relevant arcs:

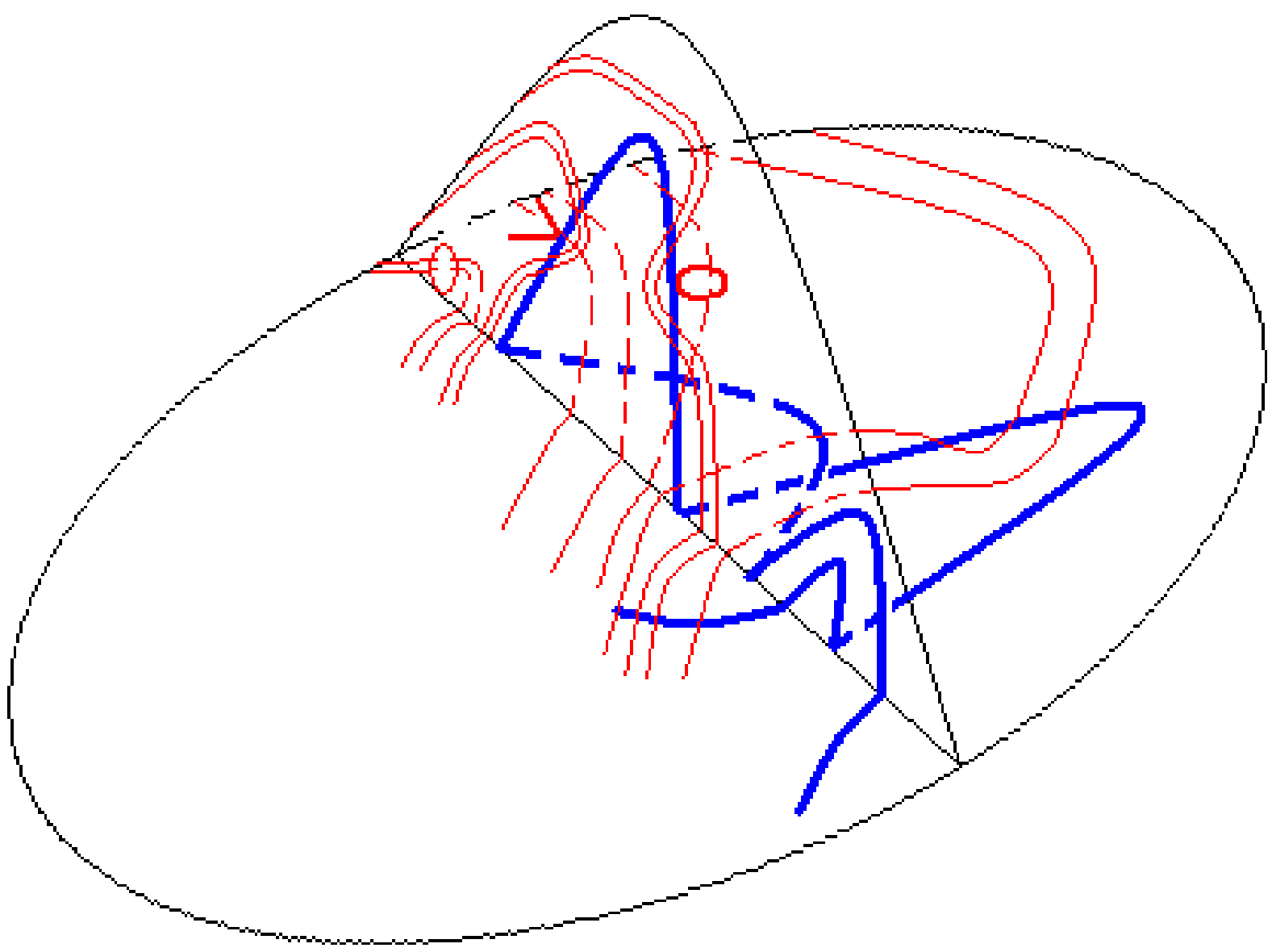

**Figure 102- Resolve the stacked up loops- push loop into backside component- 1**

The idea is to cancel the extremum on the turn in the top component, so this turn gets sliced as a line segment and then there is no obstruction to push it into the backside component. We shift the extrema into the backside component playing the game of introducing and cancelling pairs of saddlepoints. Ignore the local changes at the vertex, these belong to a more complicated idea which we do not study here:

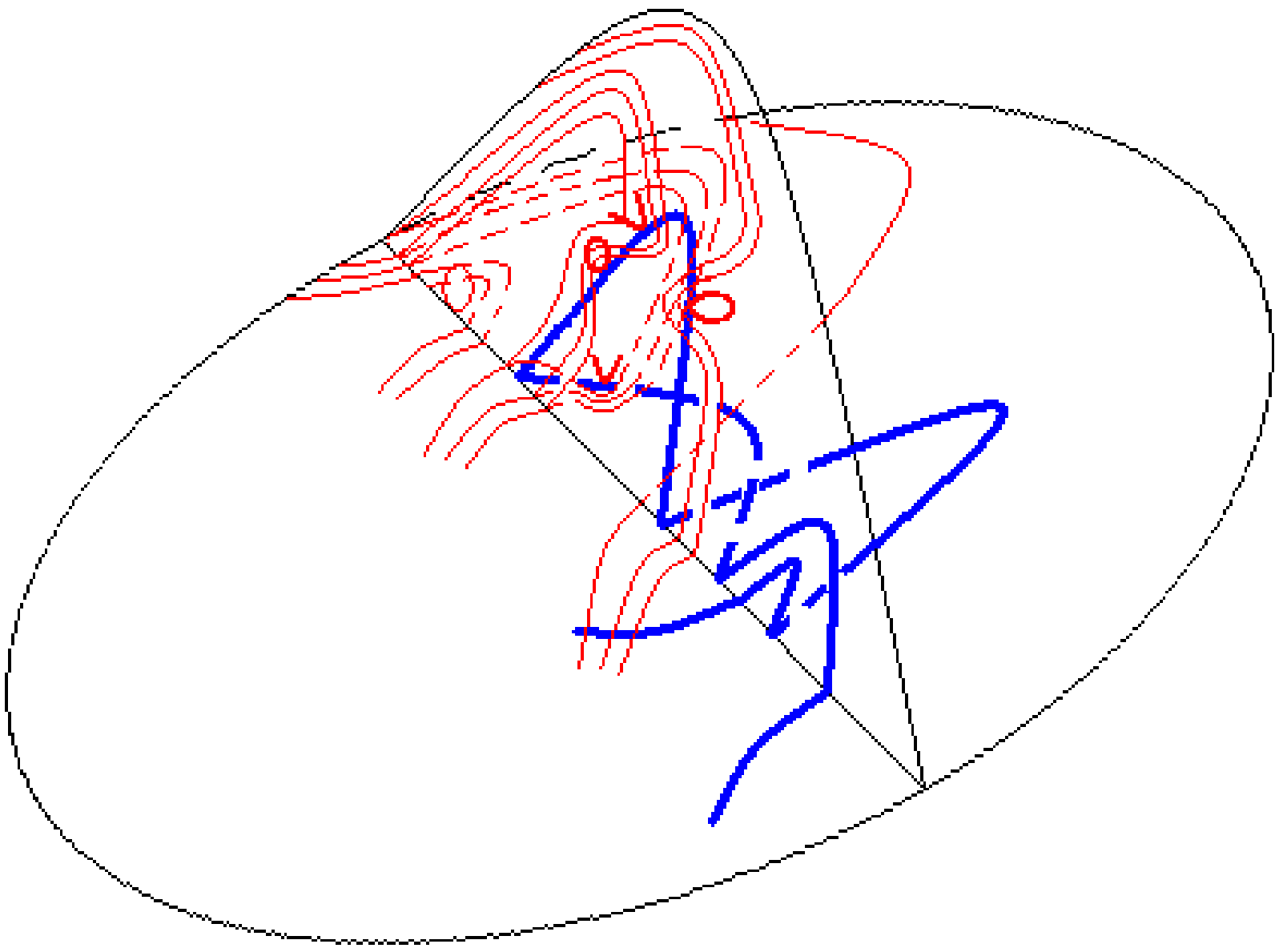

**Figure 103- Resolve the stacked up loops- push loop into backside component- 2**

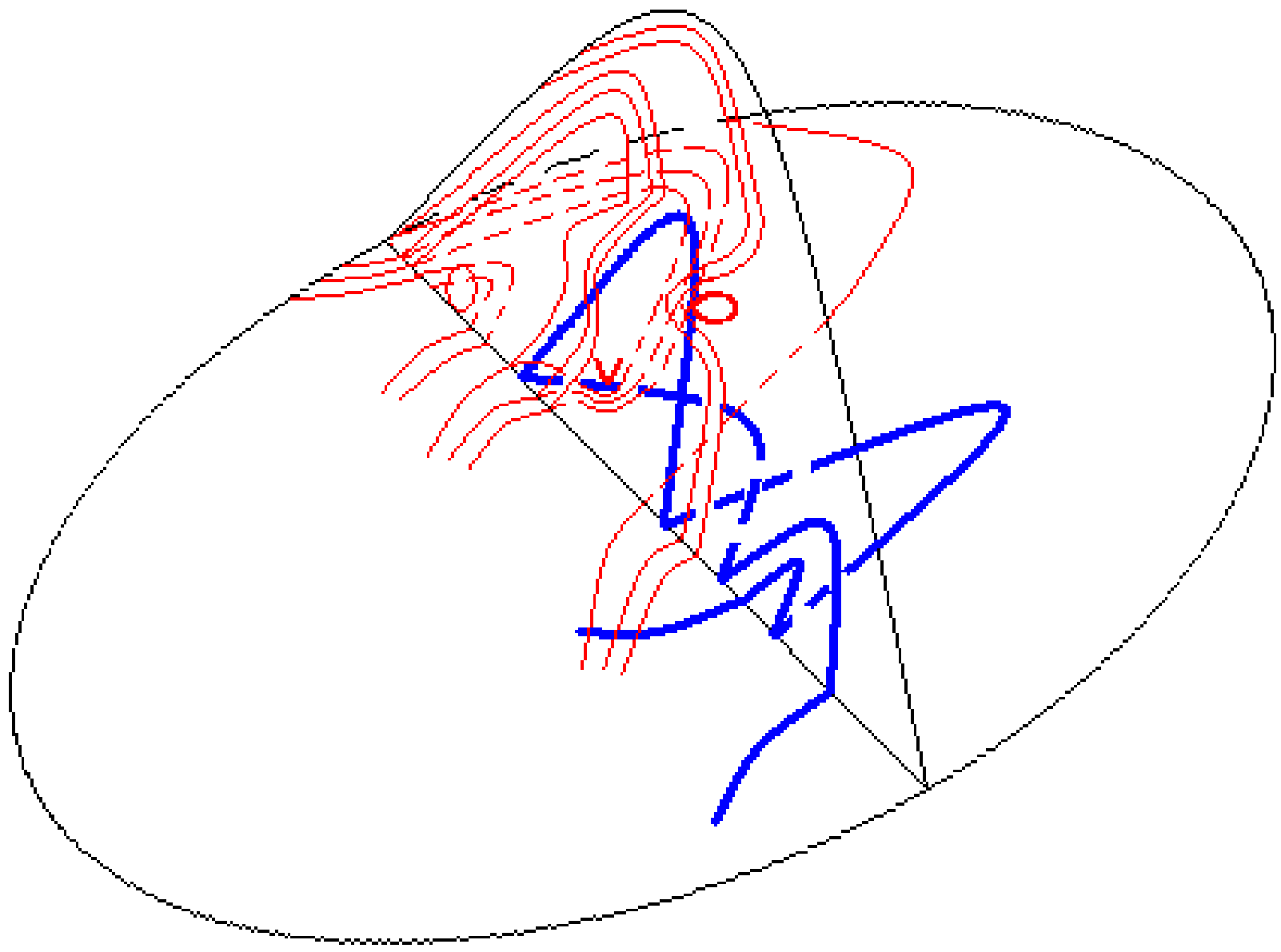

**Figure 104- Resolve the stacked up loops- push loop into backside component- 3**

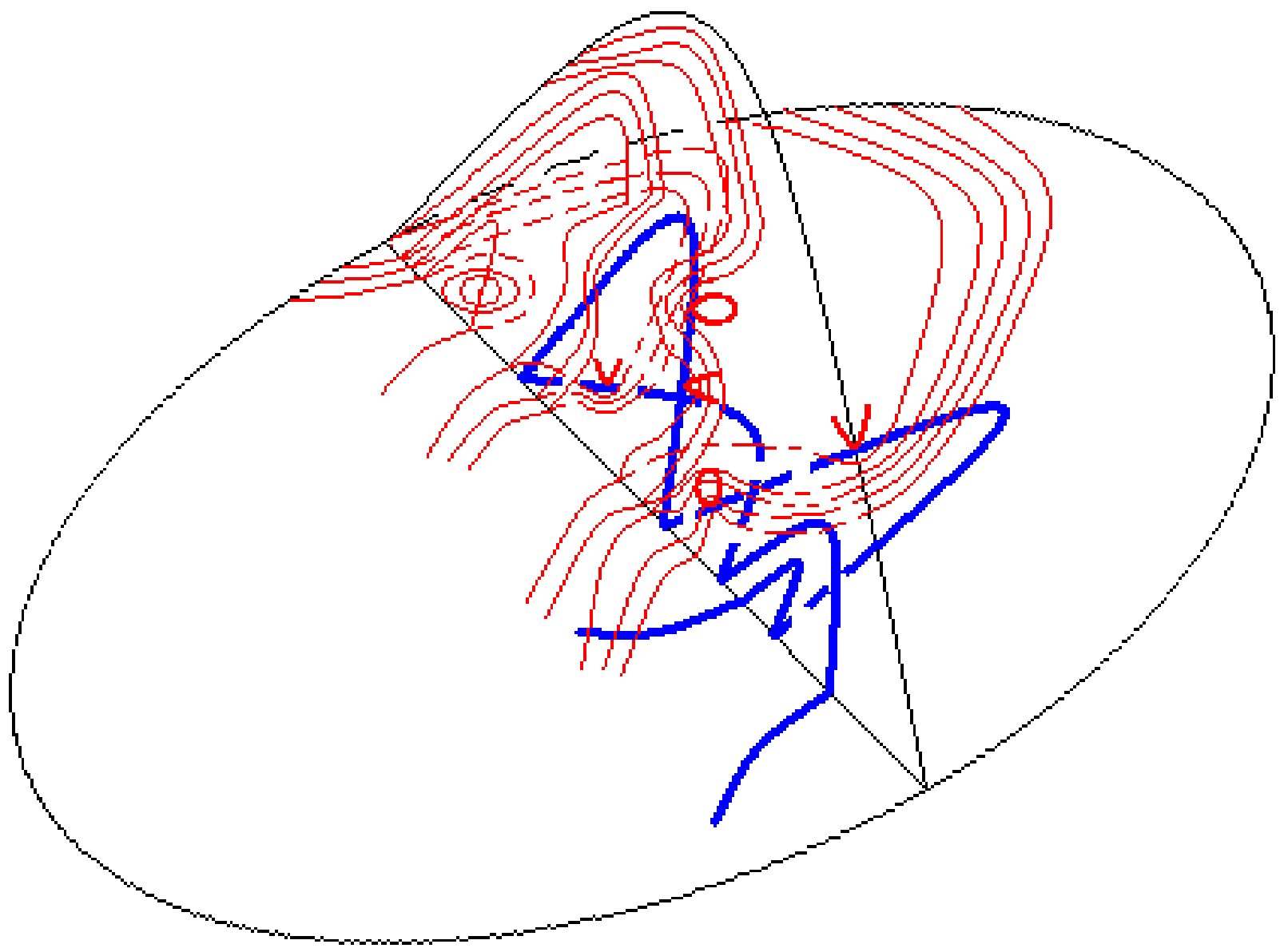

**Figure 105- Resolve the stacked up loops- push loop into bakside component- 4**

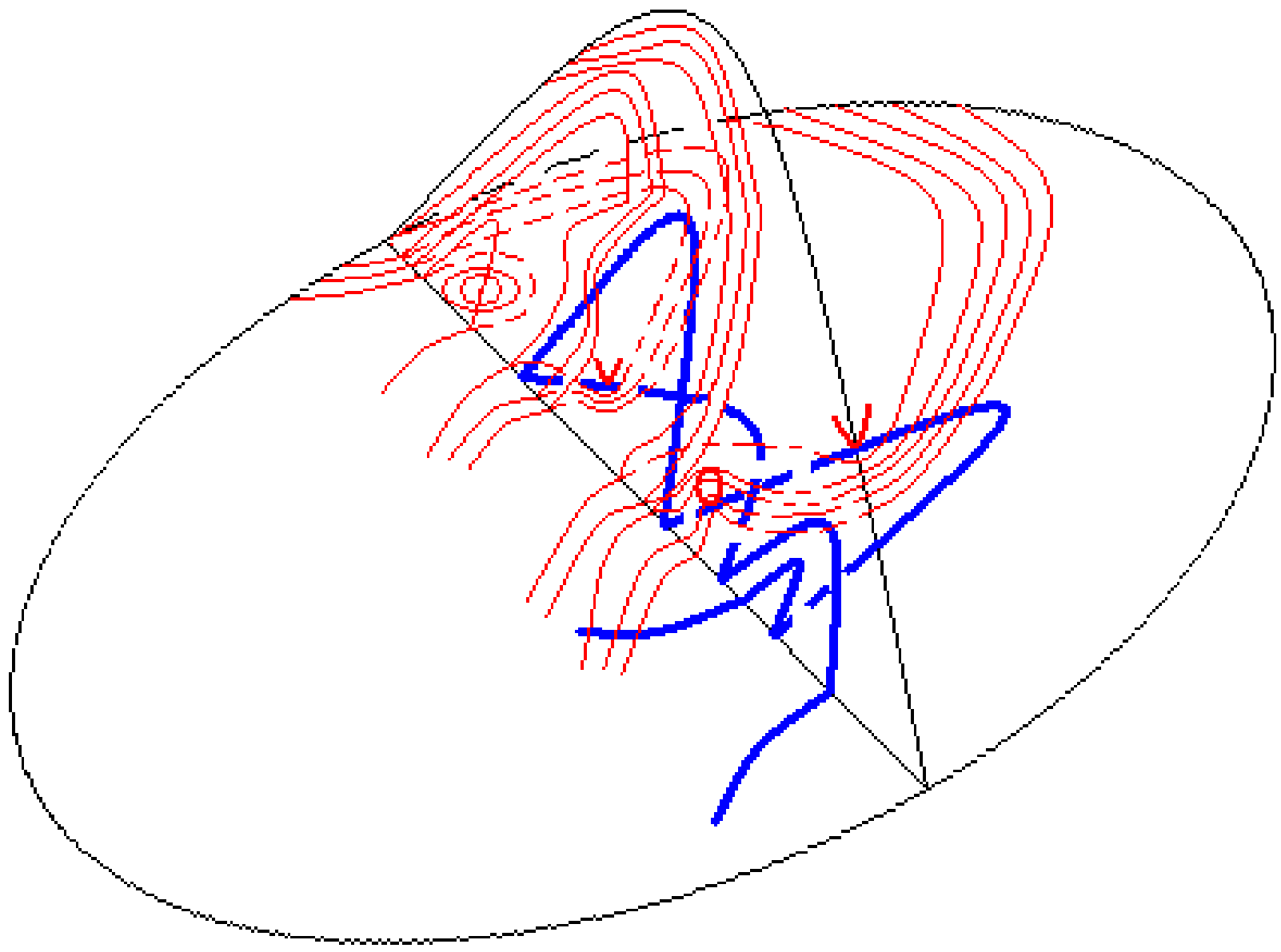

**Figure 106- Resolve the stacked up loops- push loop into backside component- 5**

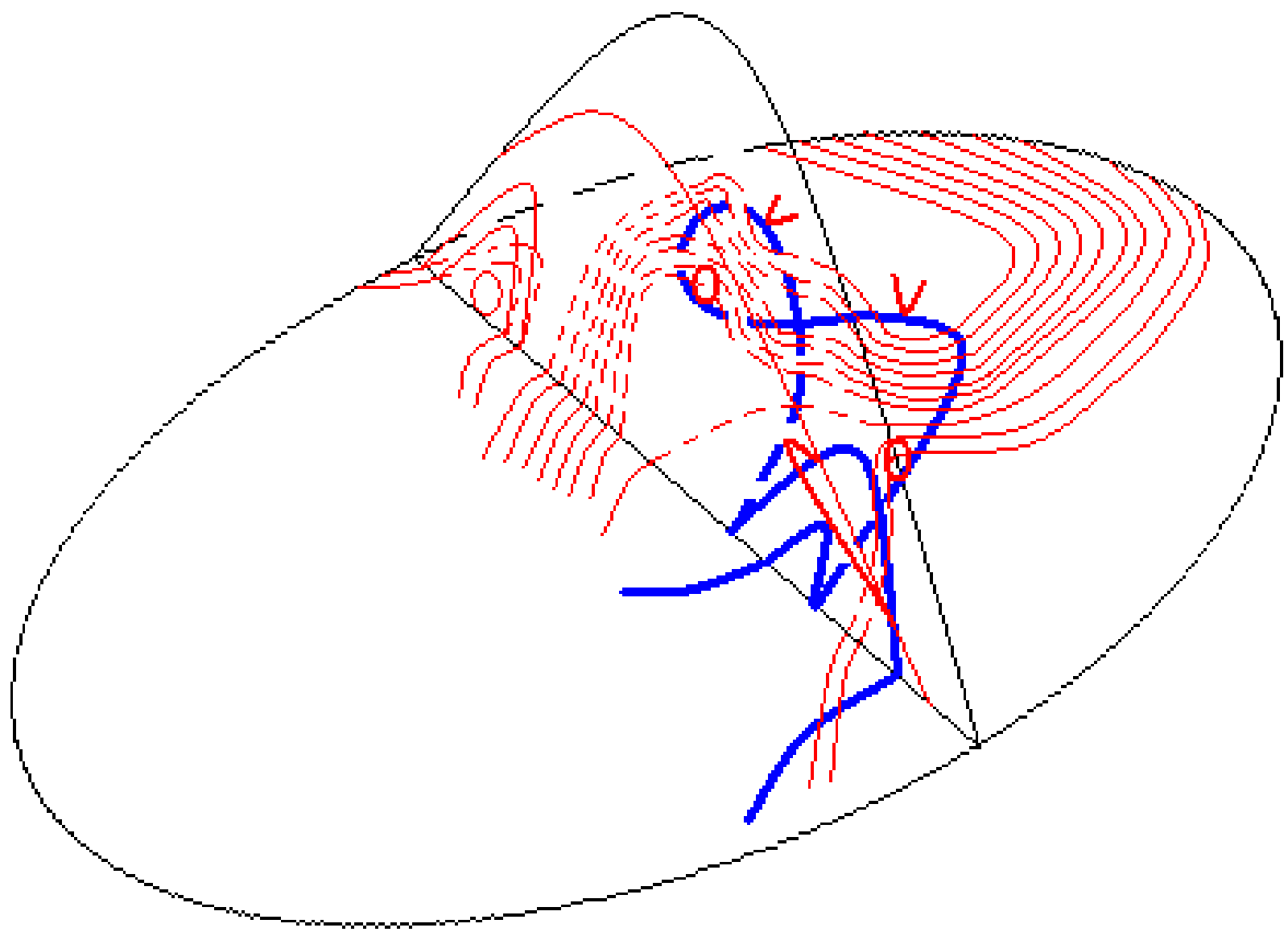

**Figure 107- Resolve the stacked up loops- push loop into backside component- 6 (end)**

We annihilate the pair of saddlepoints as shown in the next local picture

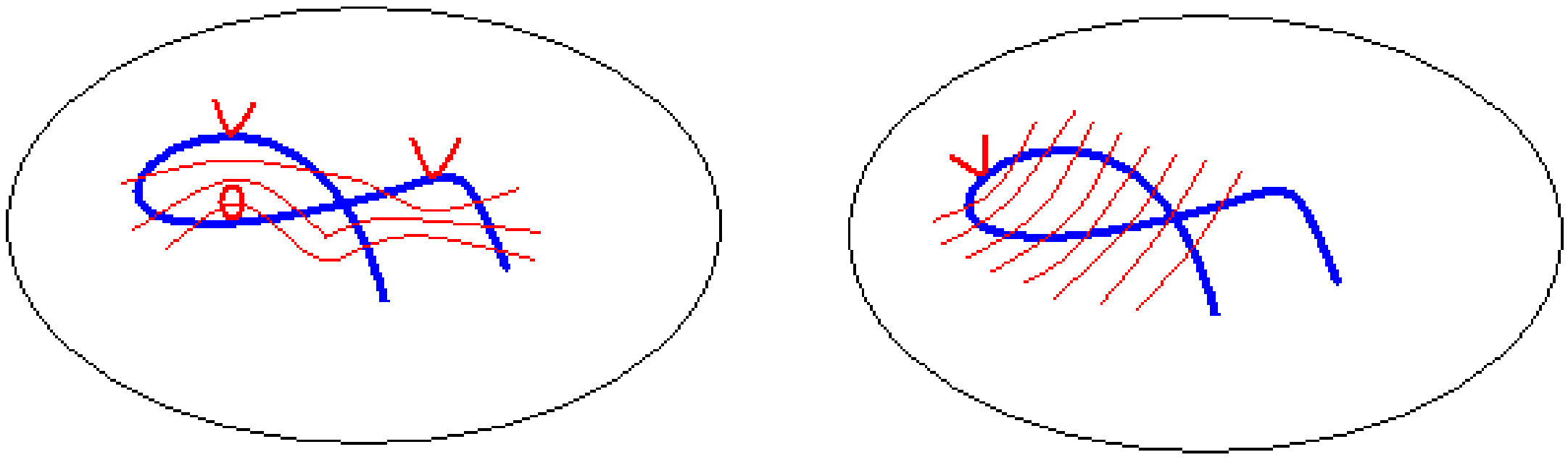

**Figure 108- Resolve the stacked up loops- annihilate a pair of saddlepoints- local**

and get:

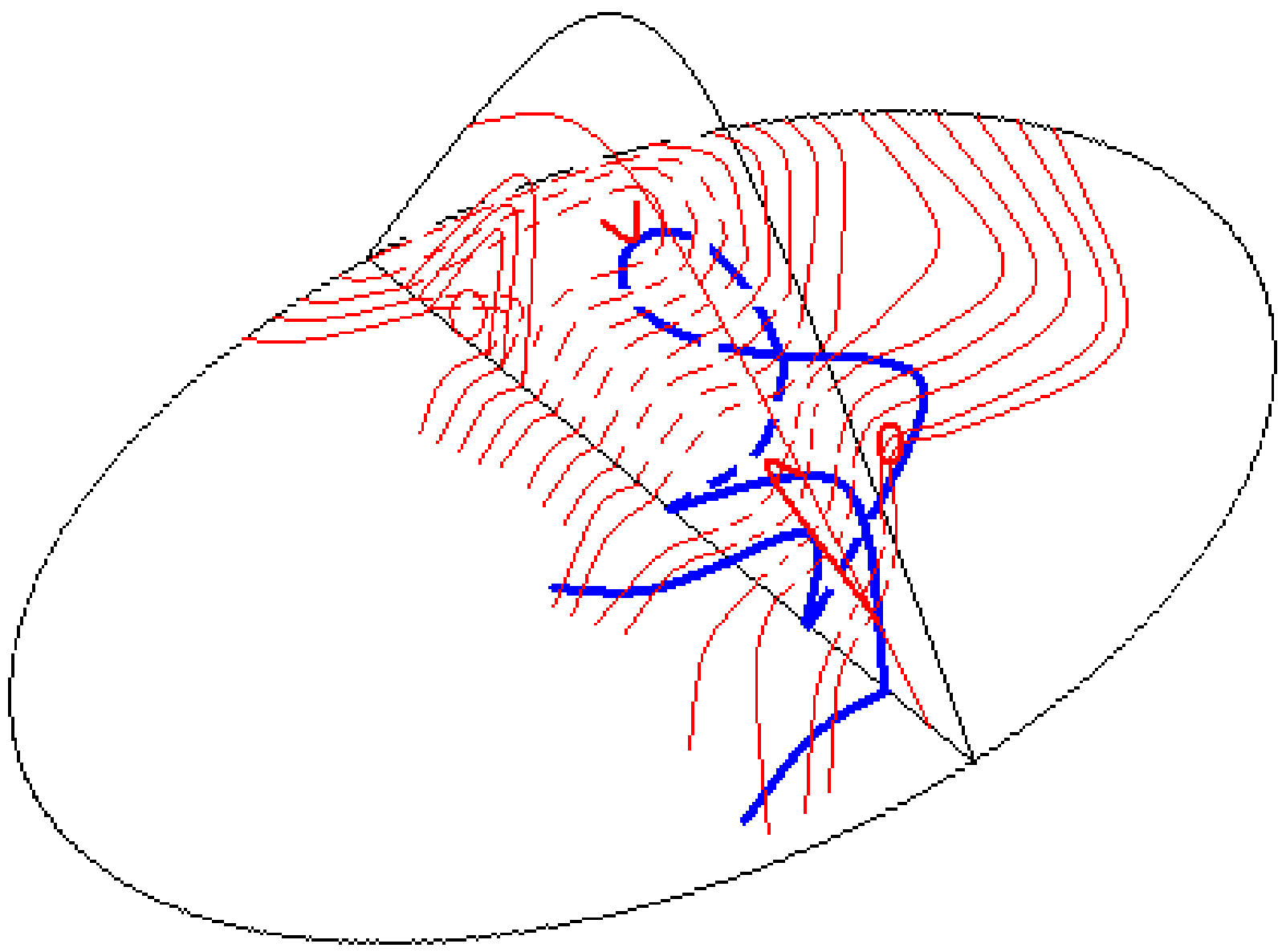

**Figure 109- Resolve the stacked up loop- annihilate a pair of saddlepoints- global**

Since this $T_3$ turn in the top component represents a ”good” one, we cancel it by the move $T_3^{-1}$:

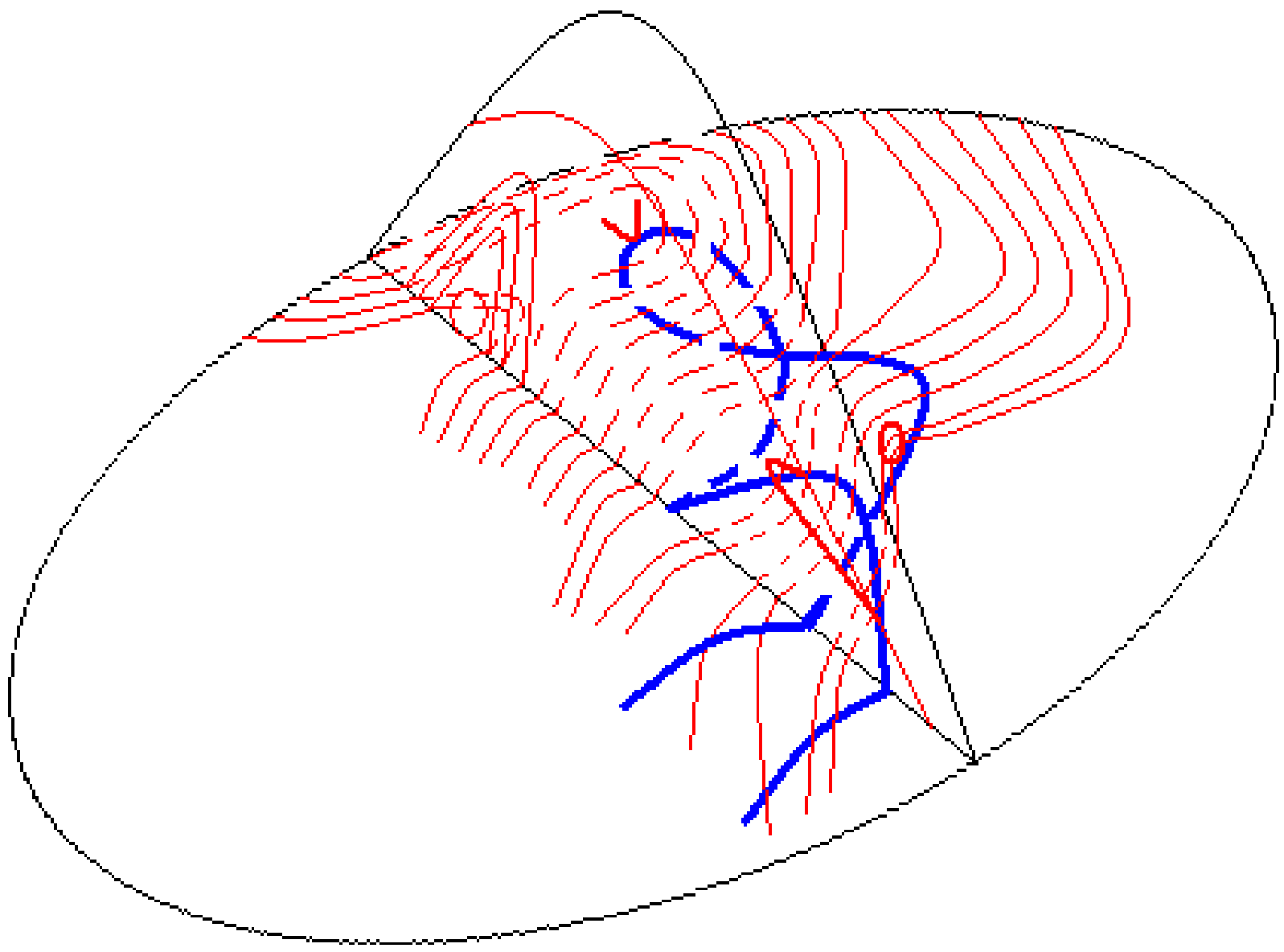

**Figure 110- Resolve the stacked up loops- annihilate $T_3$ turn**

Now we introduce a new pair of saddlepoints to be able to resolve the loop:

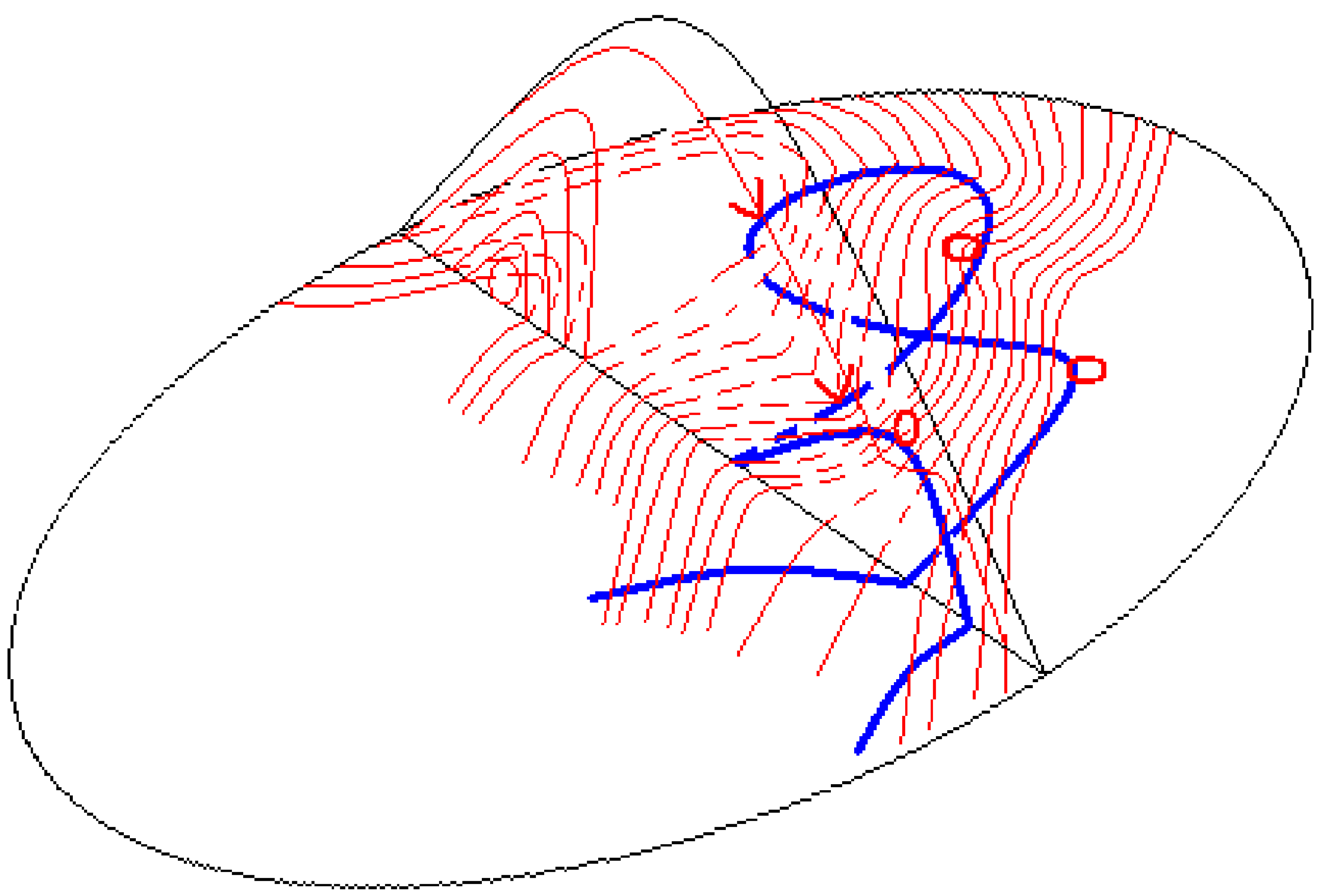

**Figure 111- Resolve the stacked up loops- prepare crossing to resolve the loop**

We consider the arc in the backside component, which connects the $T_3$ turn with the crossing point of the loop. We can pass this using T*  and thereby get a new $T_3$ turn in the backside component:

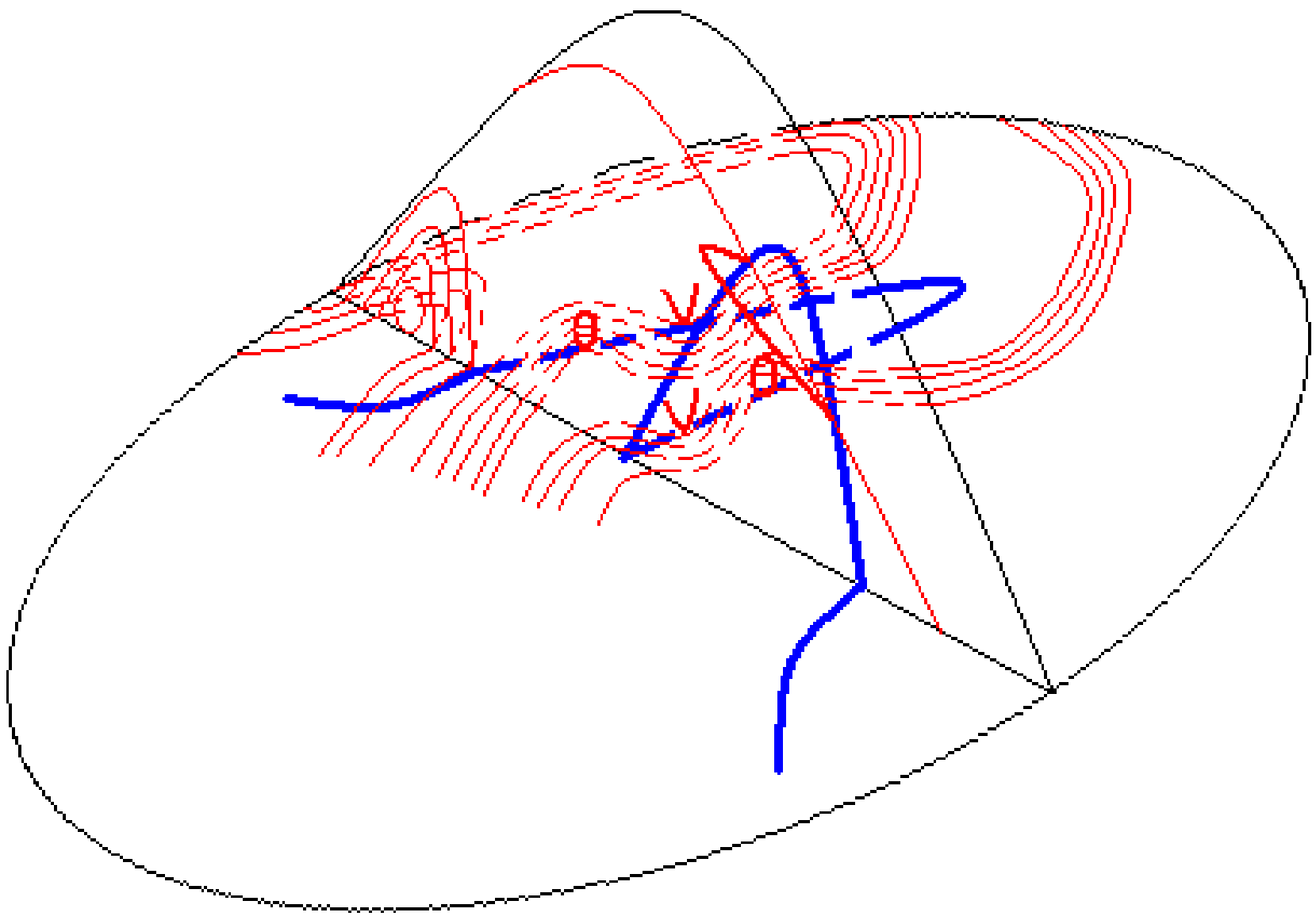

**Figure 112- Resolve the stacked up loops- new $T_3$ turn**

We annihilate the pair of saddlepoints and consider the slices on the top component to be extended. This represents a “good  $T_3$ turn”:

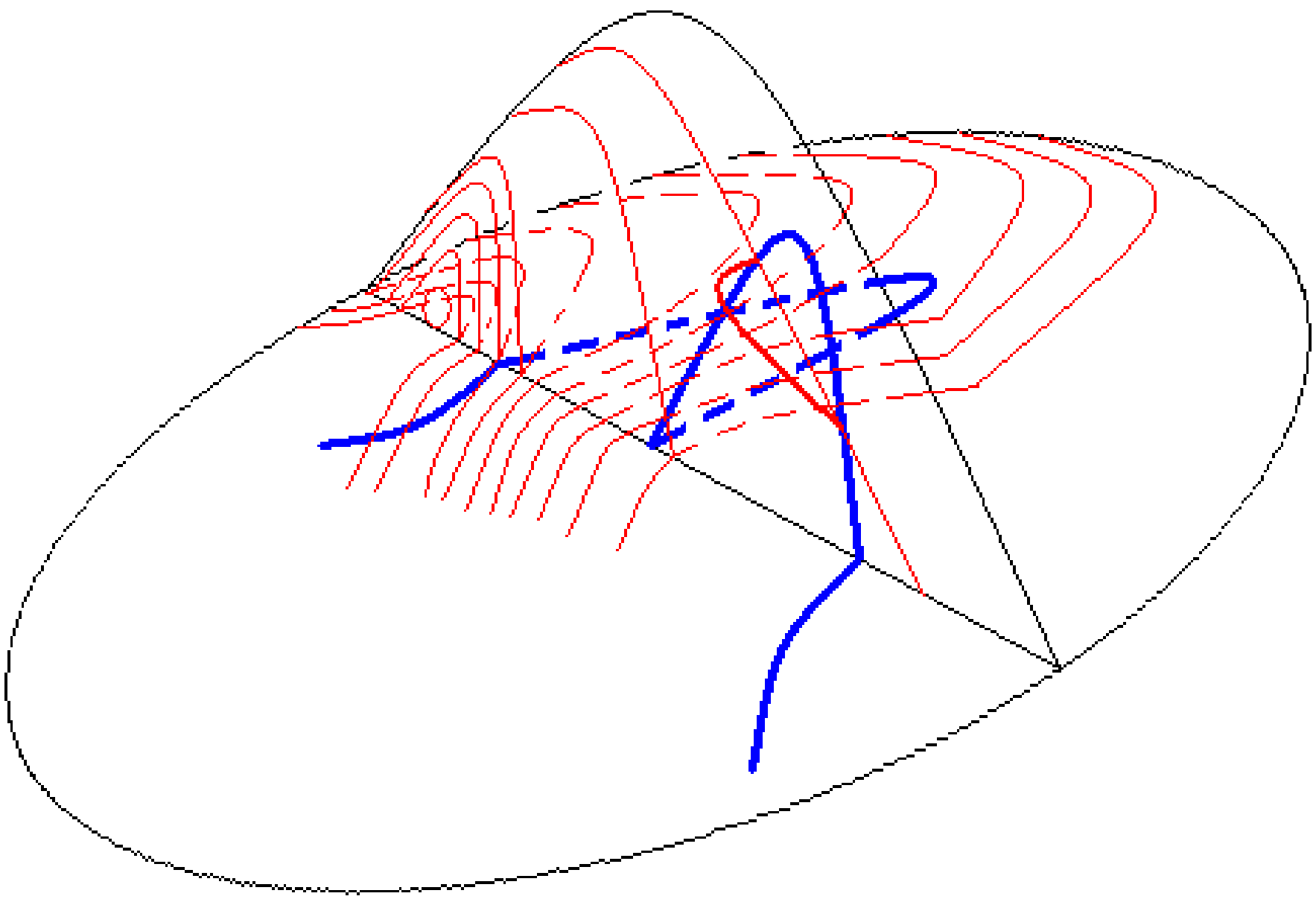

**Figure 113- Resolve the stacked up loops - good $T_3$ turn**

We cancel the $T_3$ turn in the backside component by the move $T_3^{-1}$. We can perform this move since the $T_3$ turn and the top component have the same line segment slices:

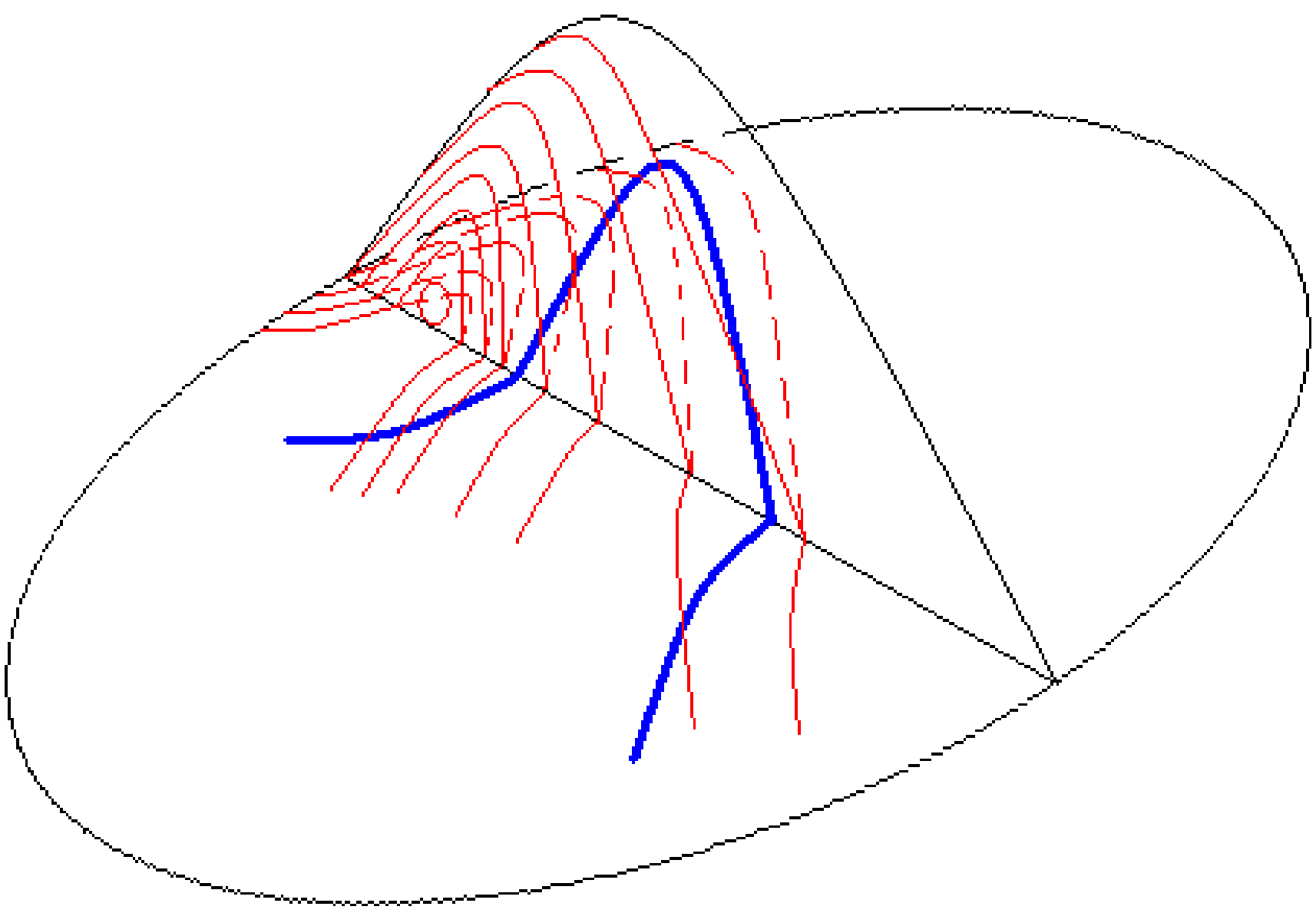

**Figure 114- Resolve the stacke up loops- startfigure without loops**

It remains to push the left arc across the vertex, using a local change of the slices. After rechanging the slices at the vertex, we get our start figure back. Hence the process of cancelling a pair of leftside loops was successful.

### 5.1.3 Transfer rightside loop to leftside loop

Of course it may be possible to repeat the whole process for a pair of rightside loops or a mixed pair of rightside and leftside loops, but we think that it is better to show the transfer of a rightside to a leftside loop:

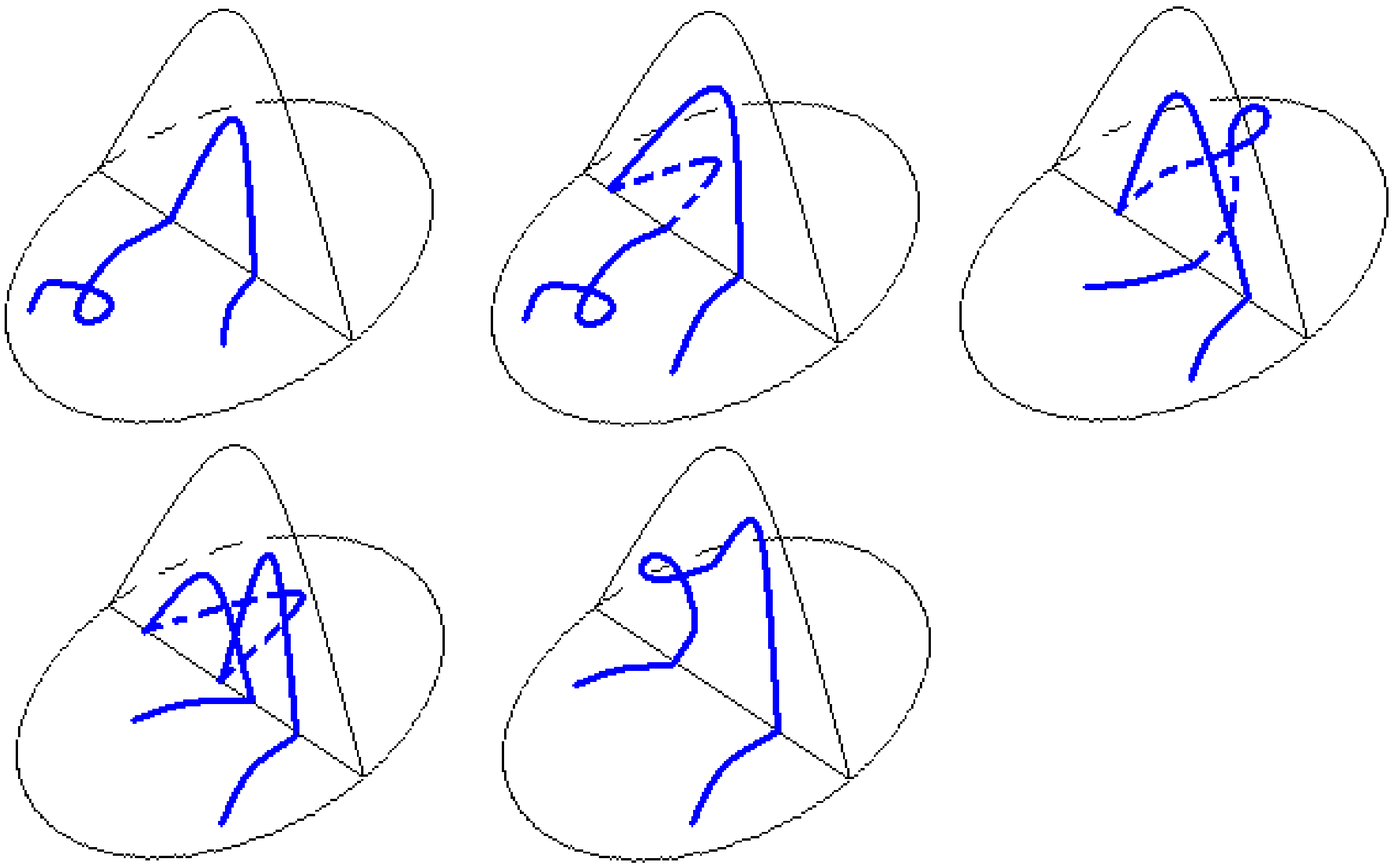

**Figure 115- Transfer rightside loop to leftside loop- preview**

The sequence of slices for the start figure looks like:

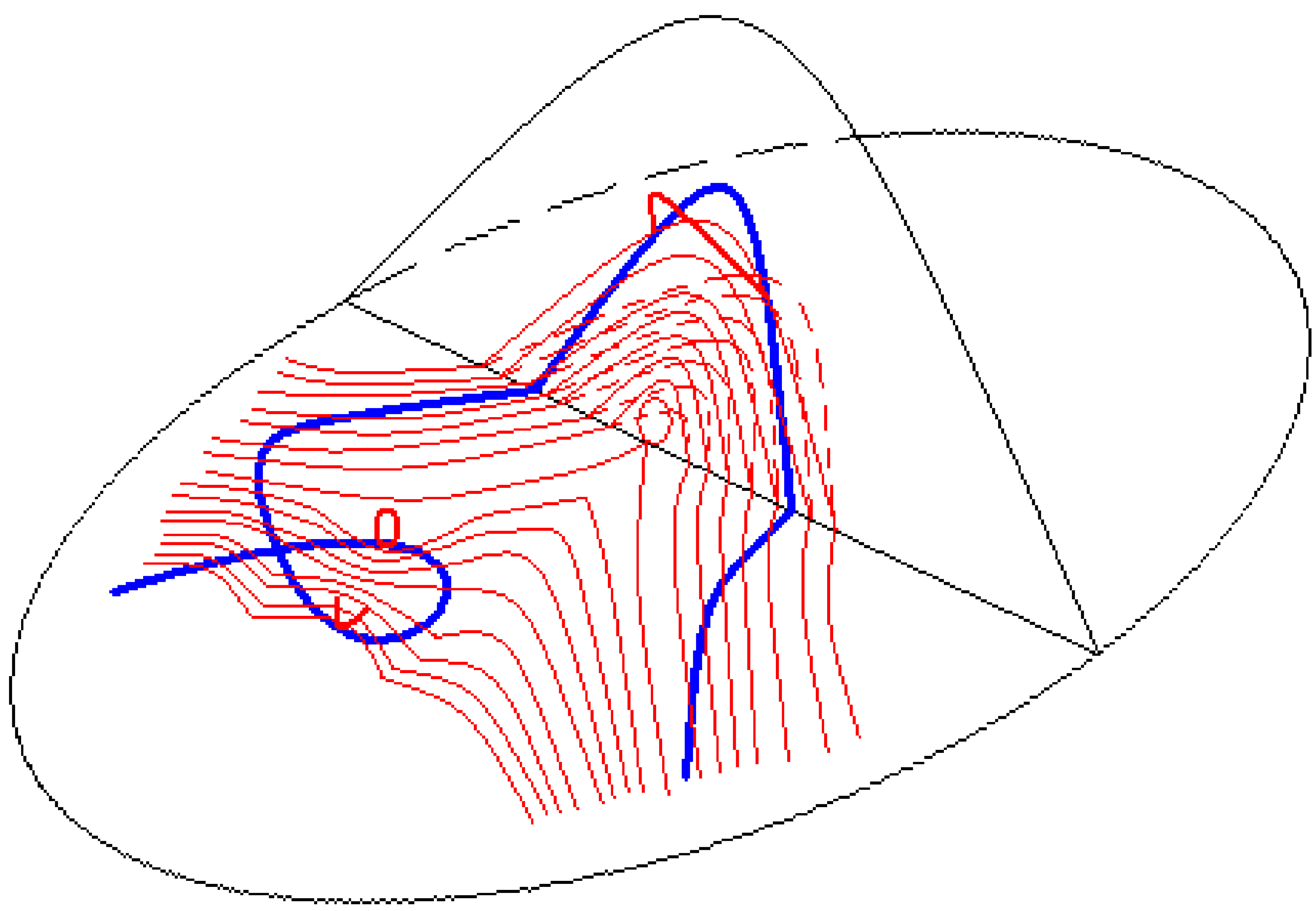

**Figure 116- Transfer rightside loop to leftside loop- 1**

We build a $T_3$ move with turn in the backside component. Note that the $T_3$ turn that arises in this component with line segment slices is a “good” one, the other $T_3$ turn in the top component has the slices of a saddlepoint, so we can not run into trouble:

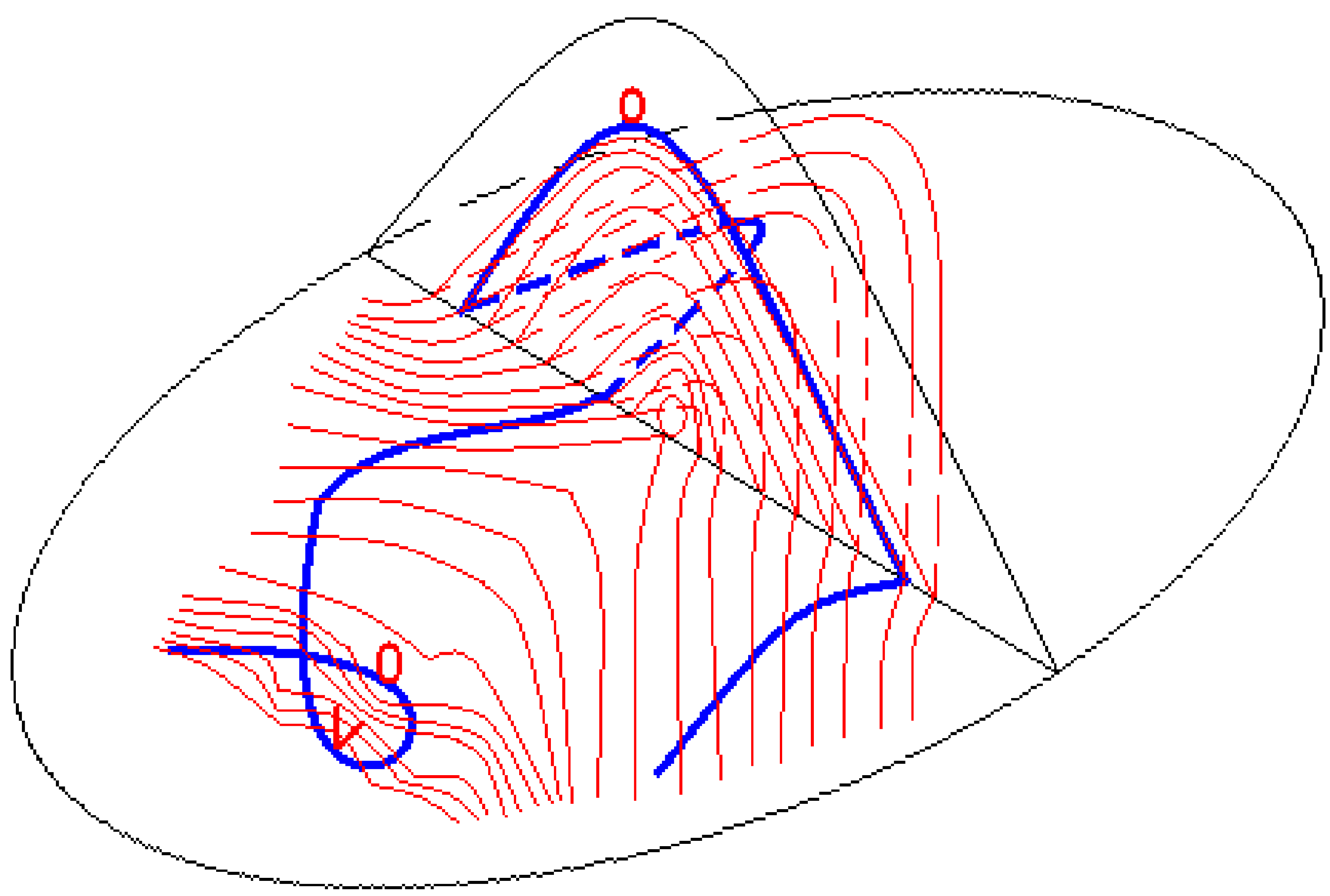

**Figure 117- Transfer rightside loop to leftsid loop- 2**

We push the loop across the vertex into the backside component:

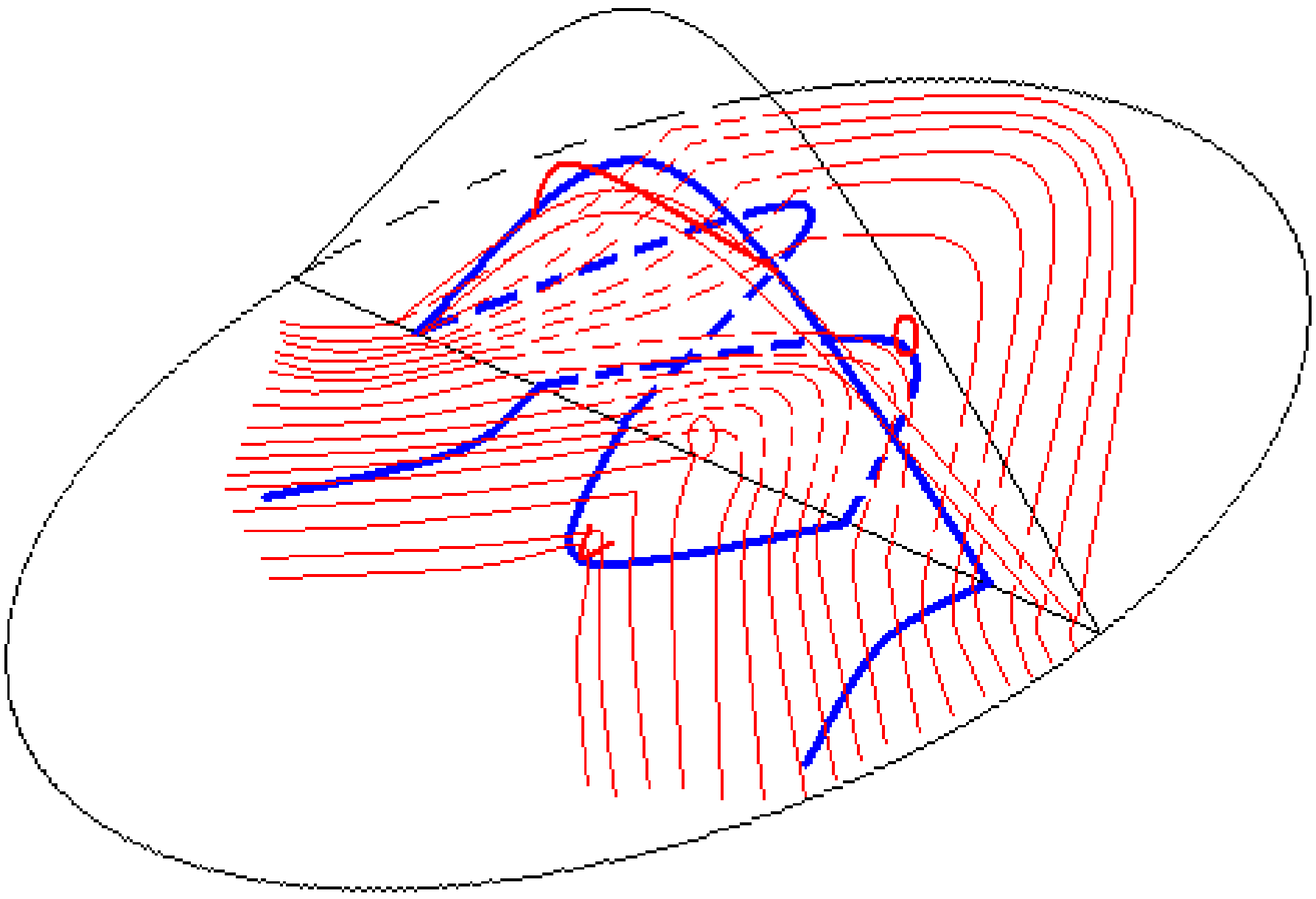

**Figure 118- Transfer rightside loop to leftside loop- 3**

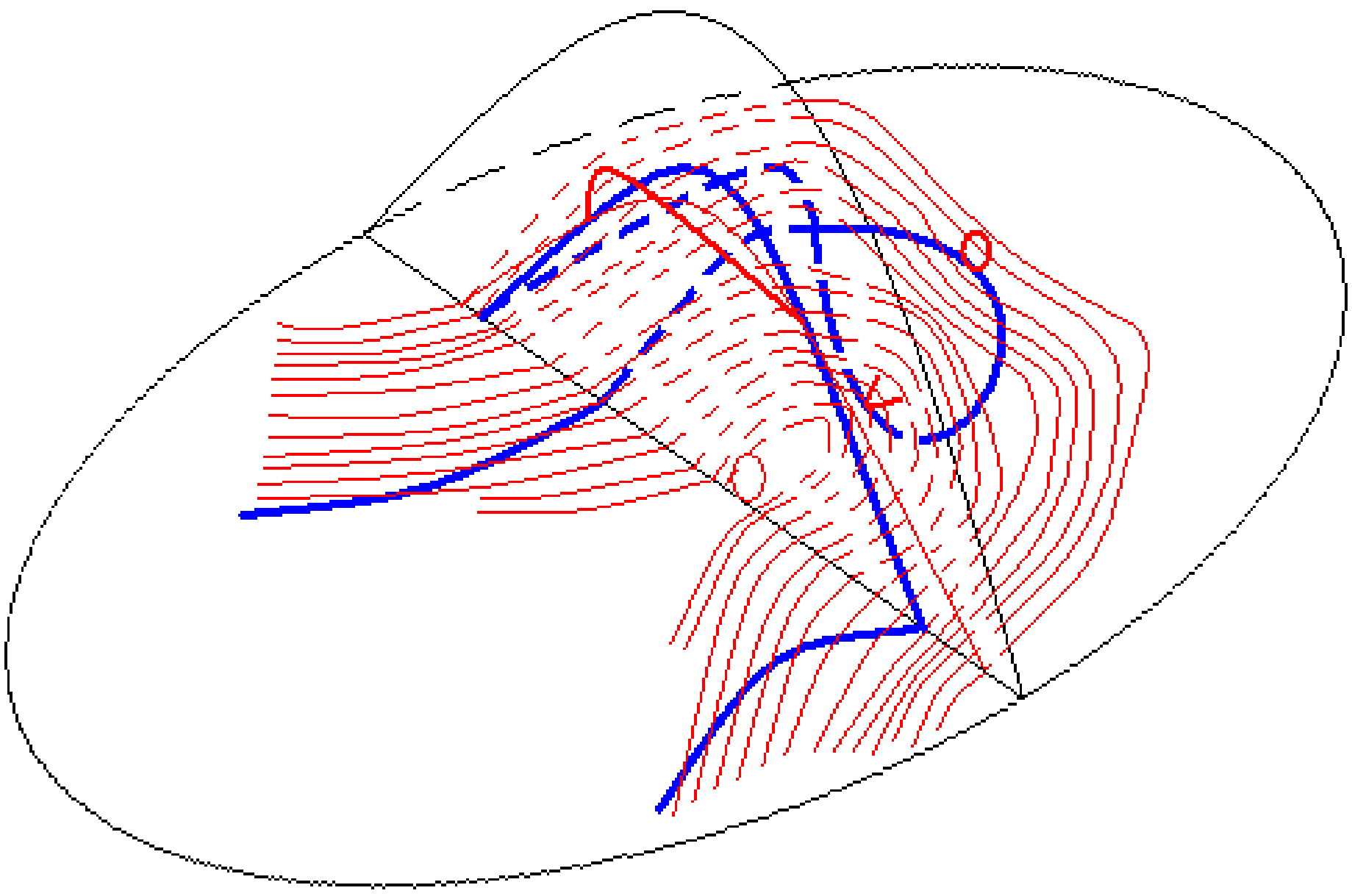

**Figure 119- Transfer rightside loop to leftside loop- 4**

The next step is a composition of moves, which we will describe in detail. The thick little line marks the vertex of the local model:

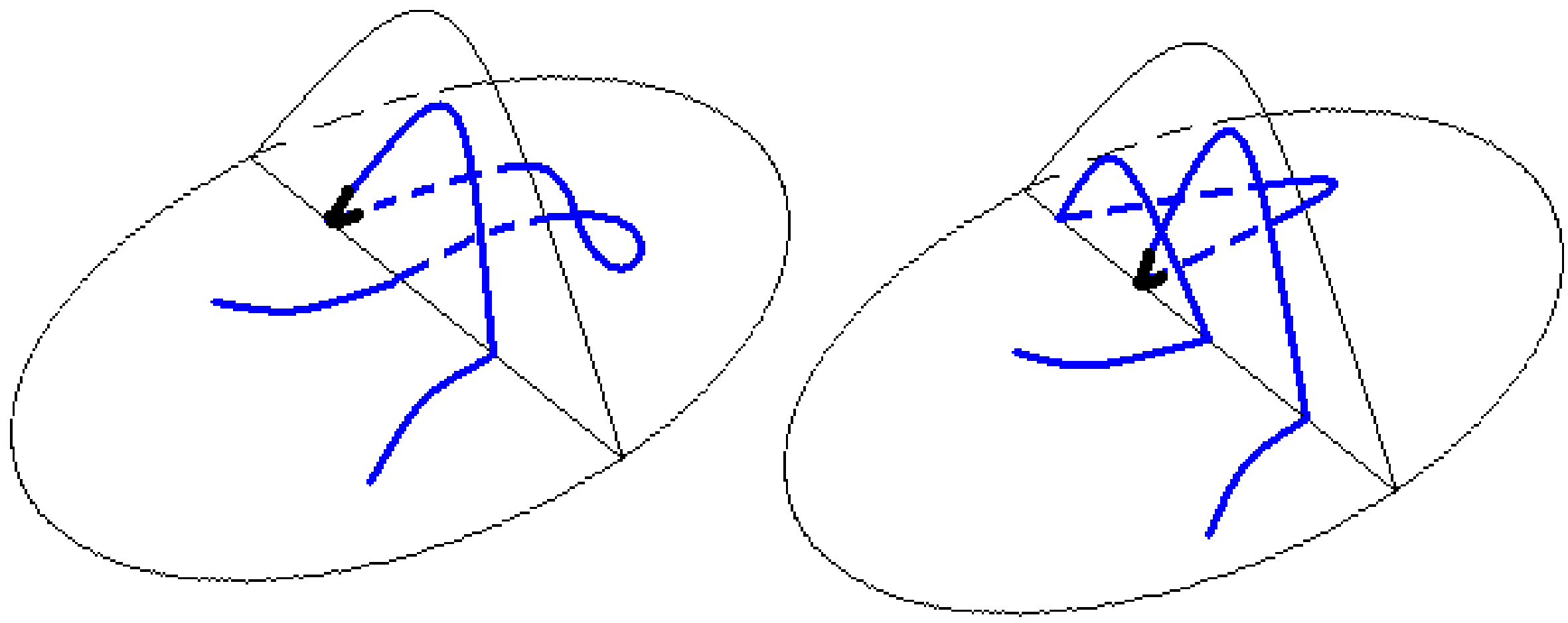

**Figure 120- Transfer rightside loop to leftside loop- preview - resolve loop in backside component**

The moves for the local model are (the backside component of the former figure corresponds to the frontside component):

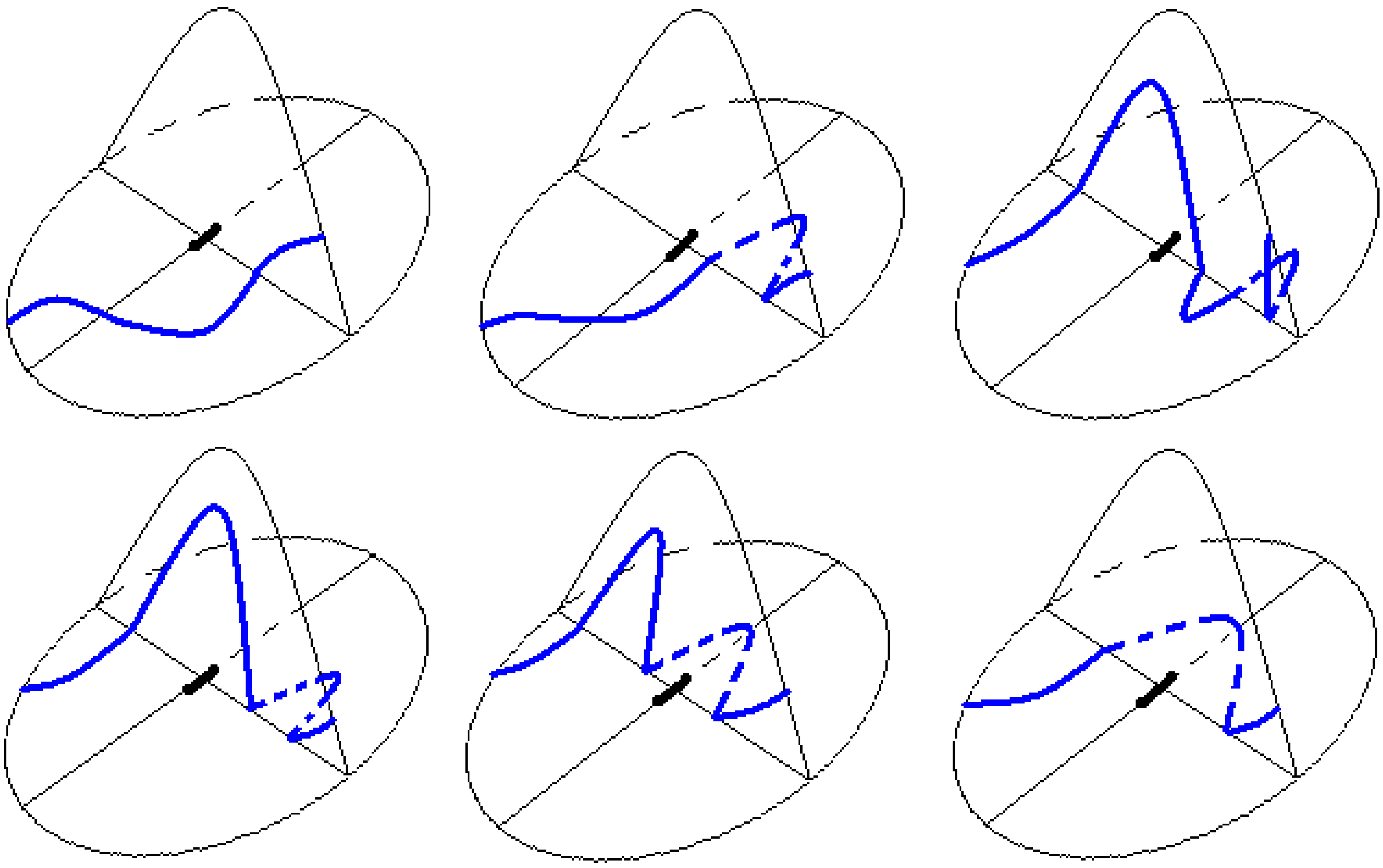

**Figure 121- Transfer rightside loop to leftside loop- resolve loop in backside component- local**

We embed this local sequence into the start figure of the current step and get:

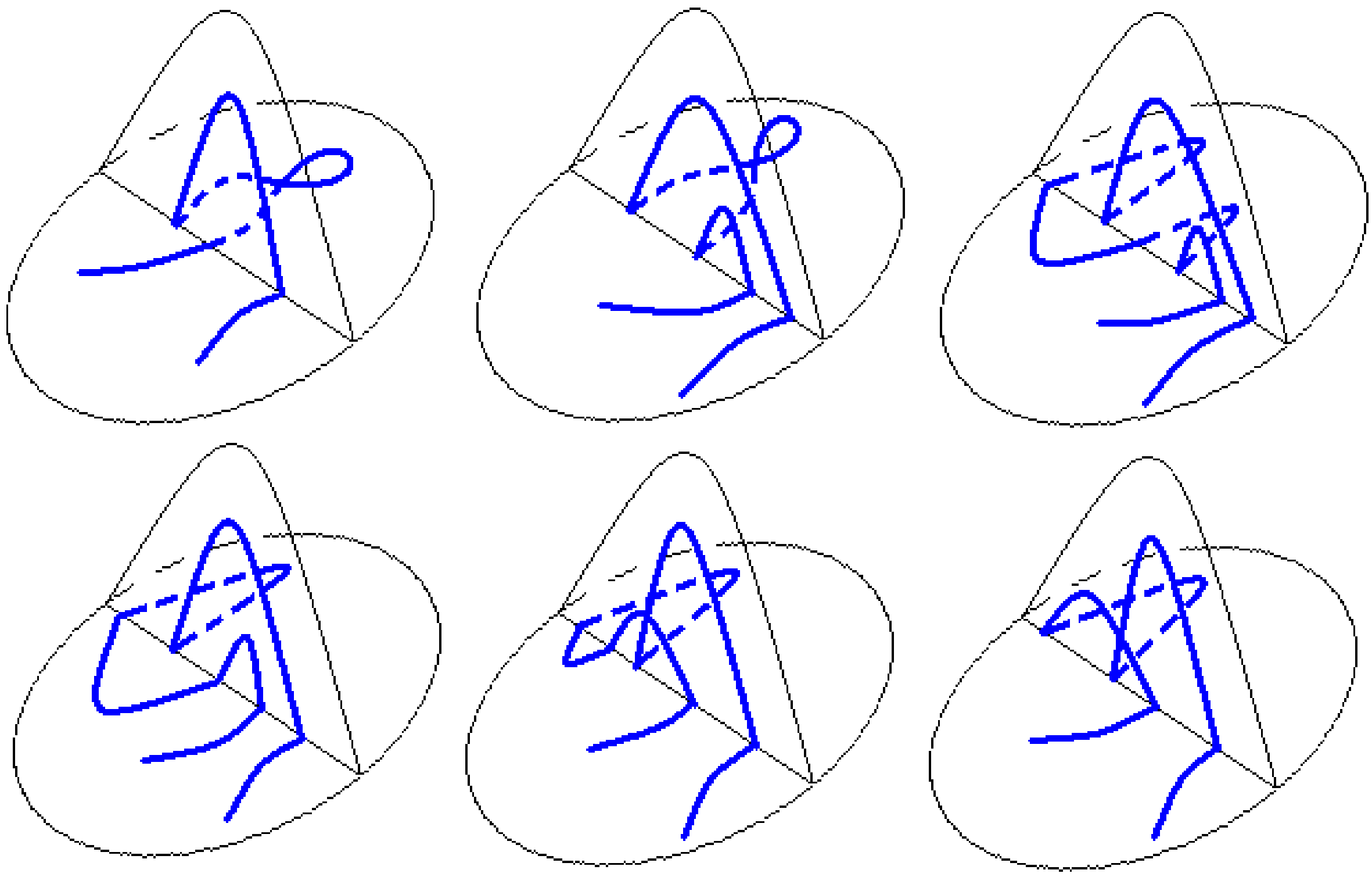

**Figure 122- Transfer rightside loop to leftside loop- resolve loop in backside component- global**

We have to ensure, that there exists a sequence of slices for each move of the step. We perform a $T_3$ move in the start figure of the sequence and the $T_3$ turn in the top component is a “good” one:

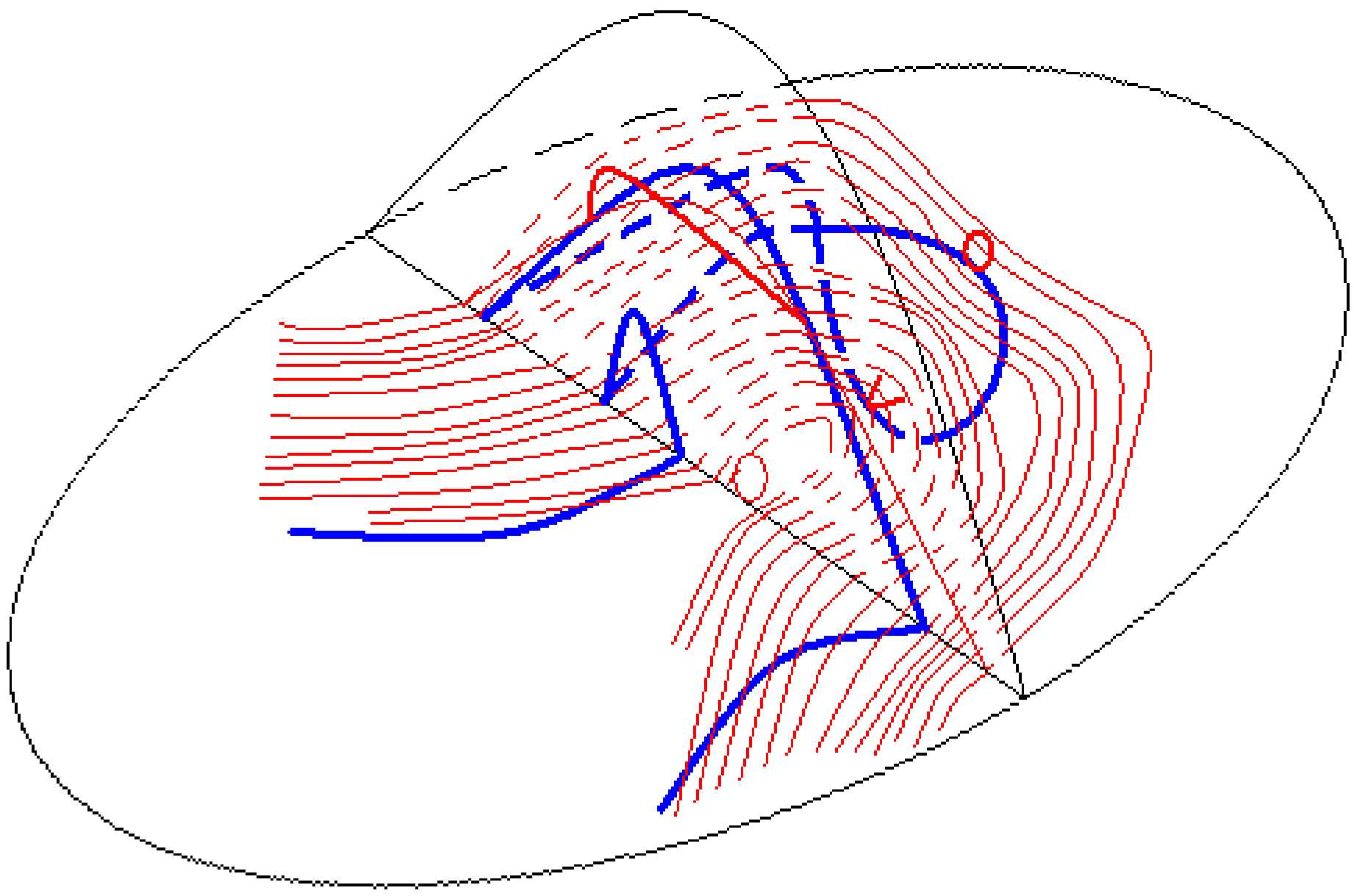

**Figure 123- Transfer rightside loop to leftside loop- resolve loop in backside- 1**

We take a part of the crossing arc in the backside component of the loop which is connected to the new $T_3$ turn. By executing the $T_2$ move on this part, we get two new $T_3$ turns. One is on the left side with a $T_3$ turn in the backside component and inherits saddlepoints from the previous slices. The other one also has its $T_3$ turn in the backside component and inherits the slice such that it results in a ”good” one:

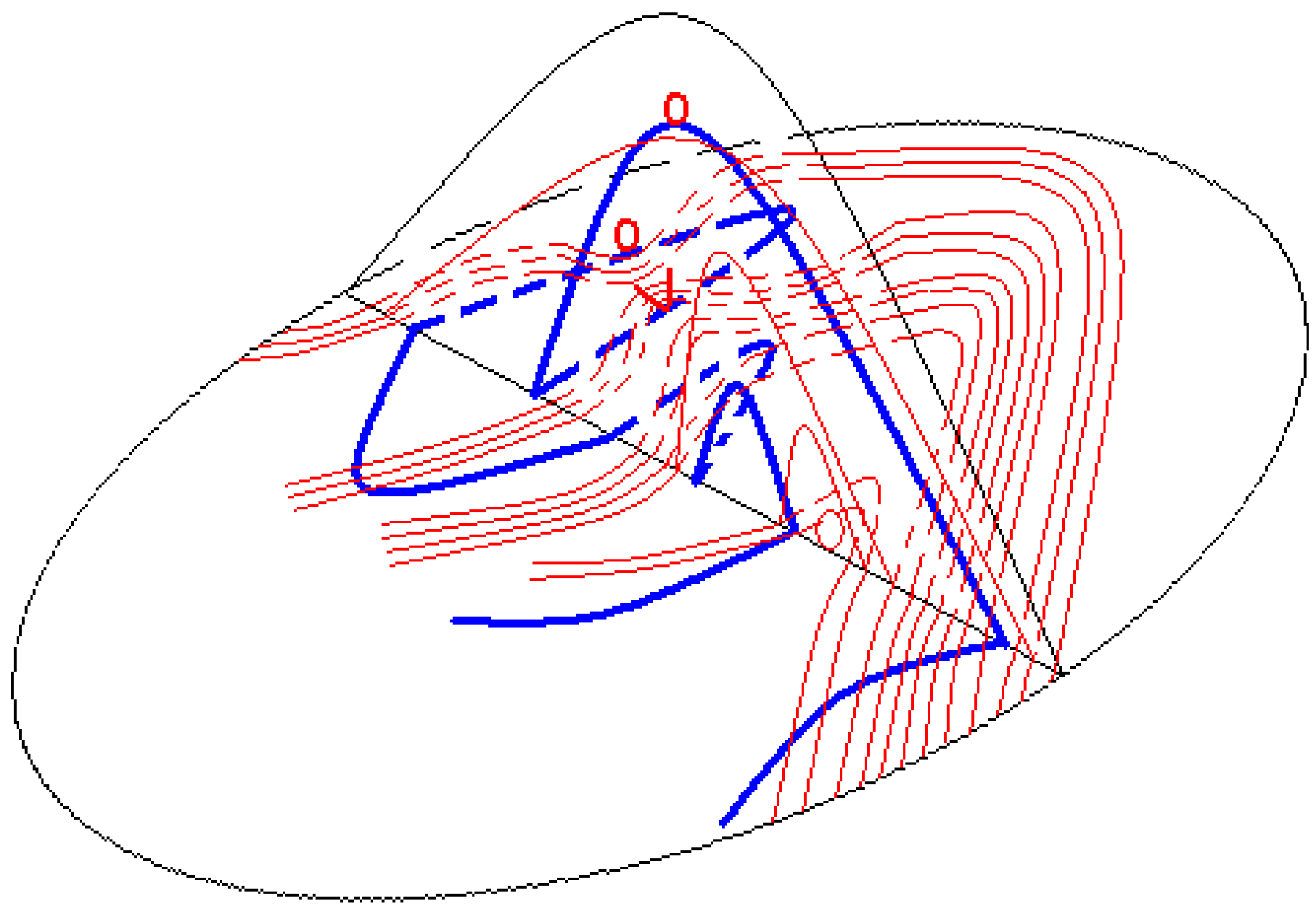

**Figure 124- Transfer rightside loop to leftside loop- resolve loop in backside- 2**

We annihilate this $T_3$ turn by $T_3^{-1}$ and produce a new “good” $T_3$ turn in the frontside component, what we can see by observing the slices of the connected arc in the top component (by continuing the slices):

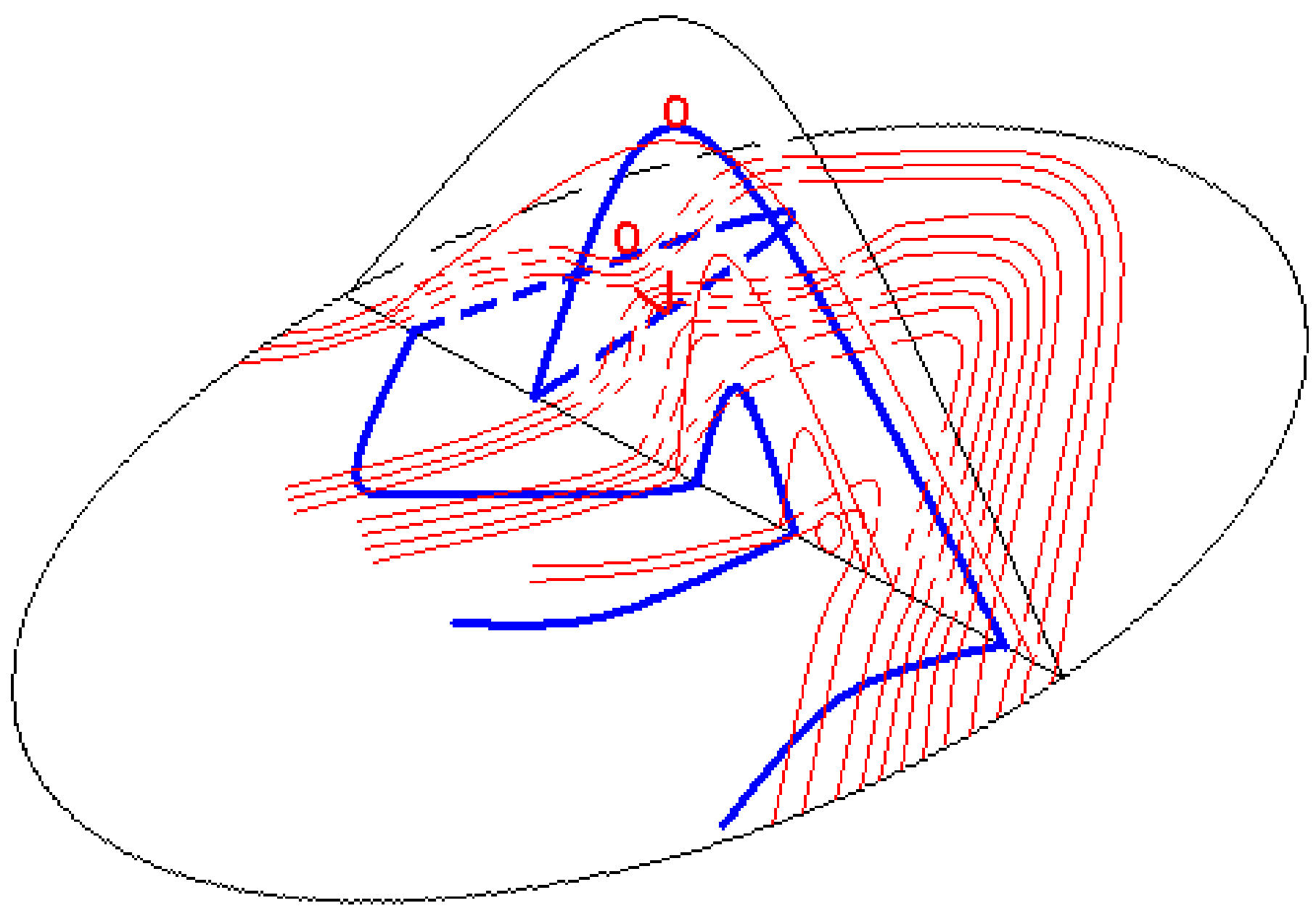

**Figure 125- Transfer rightside loop to leftside loop- resolve loop in backside- 3**

We perform T* on the $T_3$ turn in the frontside component and transfer the slices:

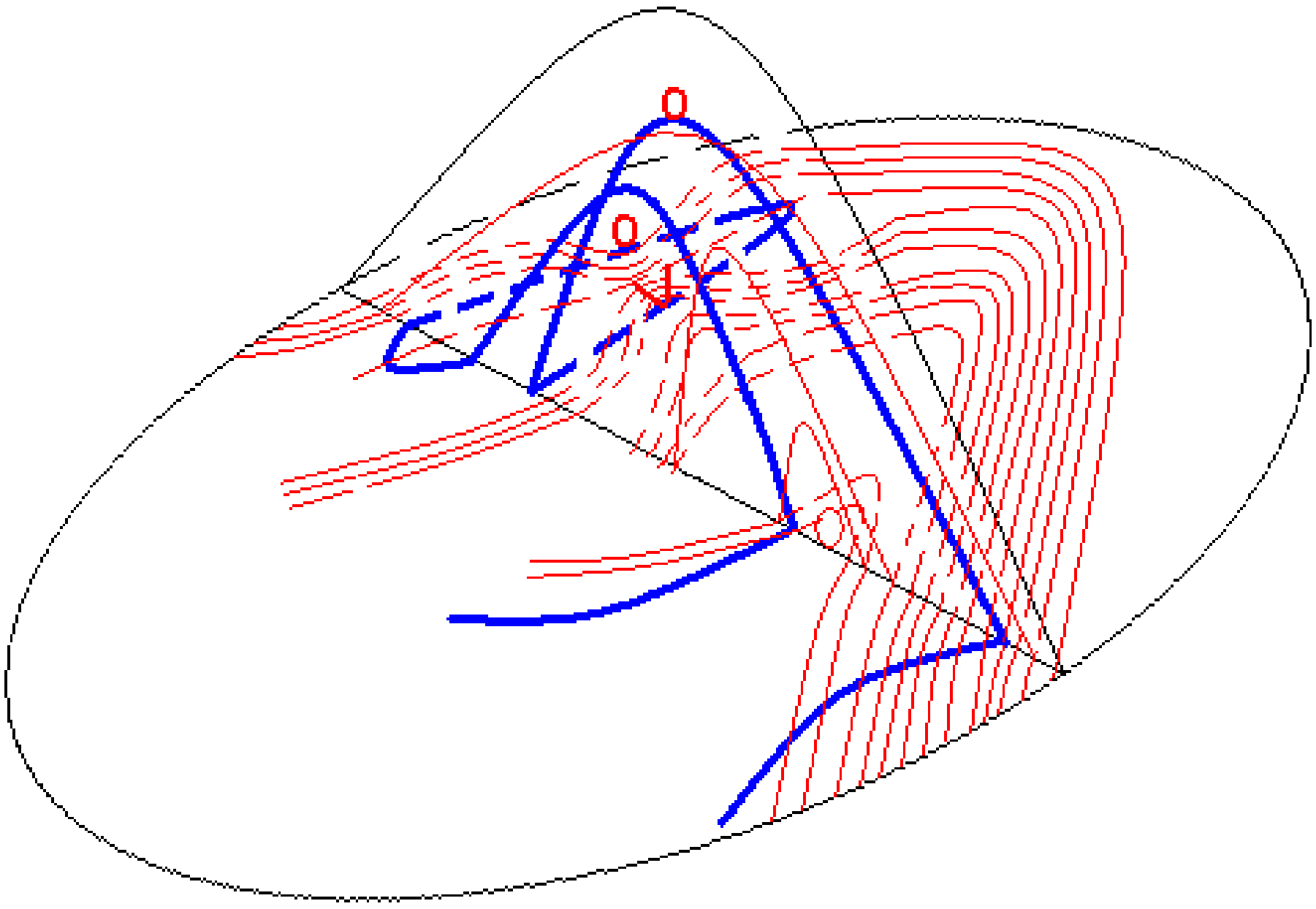

**Figure 126- Transfer rightside loop to leftside loop- resolve loop in backside- 4**

We annihilate the $T_3$ turn in the frontside component by $T_3^{-1}$ and get:

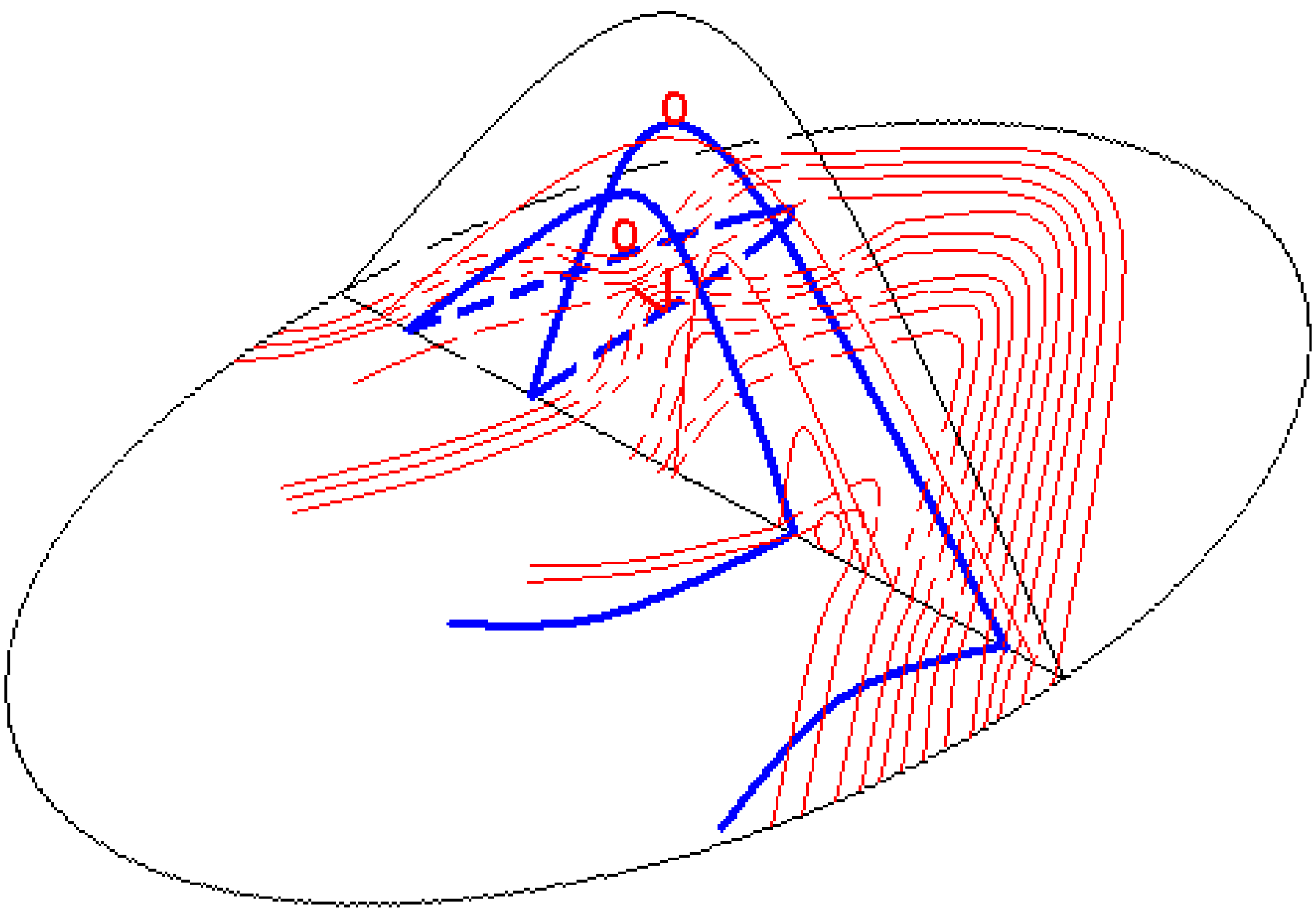

**Figure 127- Transfer rightside loop to leftside loop- resolve loop in backside- 5 (end)**

It remains to push the turn in the backside component into the top component. Similar to the case for the leftside loop we proceed by generating and annihilating pairs of saddlepoints. The idea behind this is to shift the saddlepoints into the top component to have only line segment slices for the turn in the backside component:

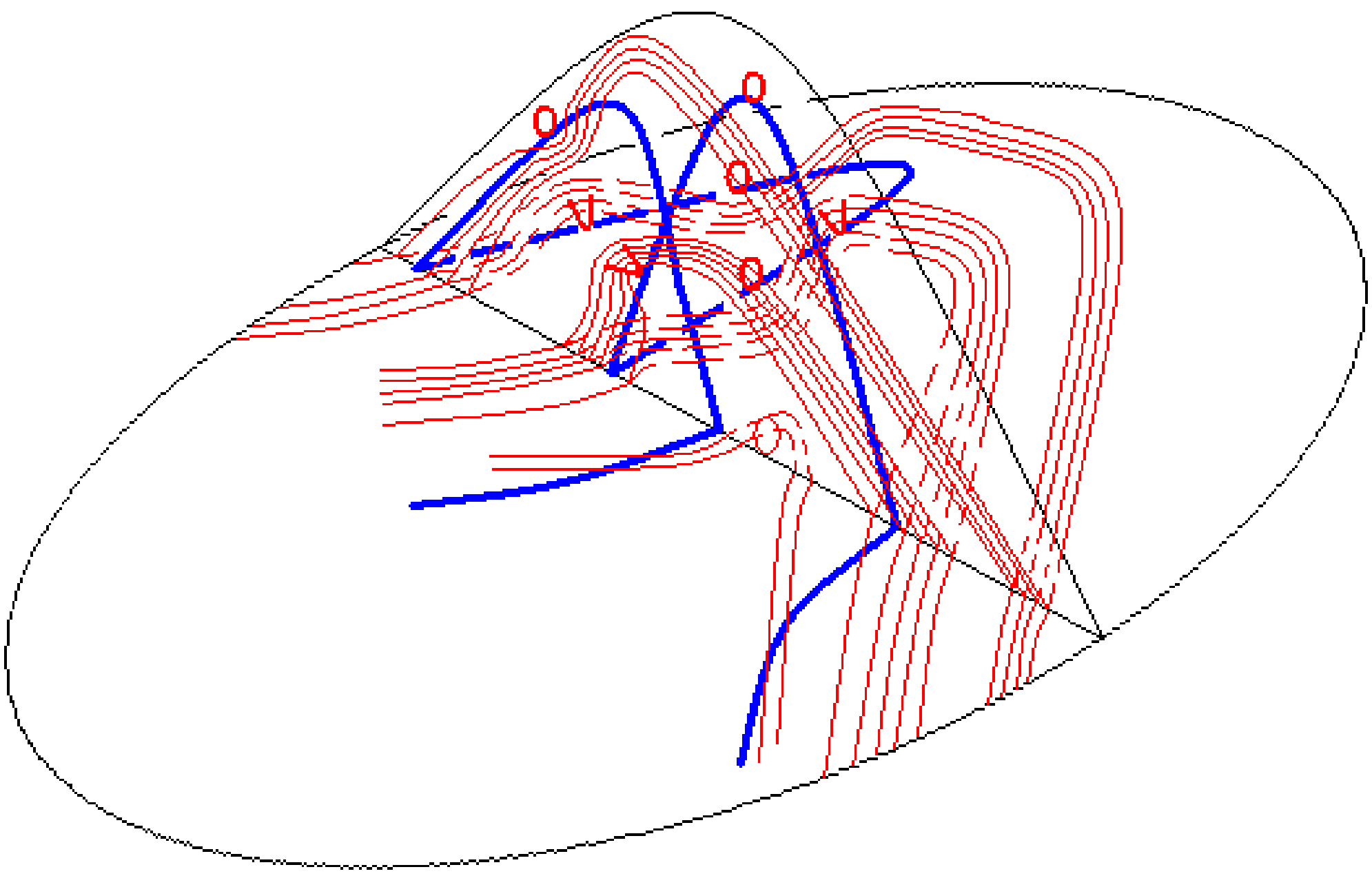

**Figure 128- Transfer rightside loop to leftside loop- push loop in top component- 1**

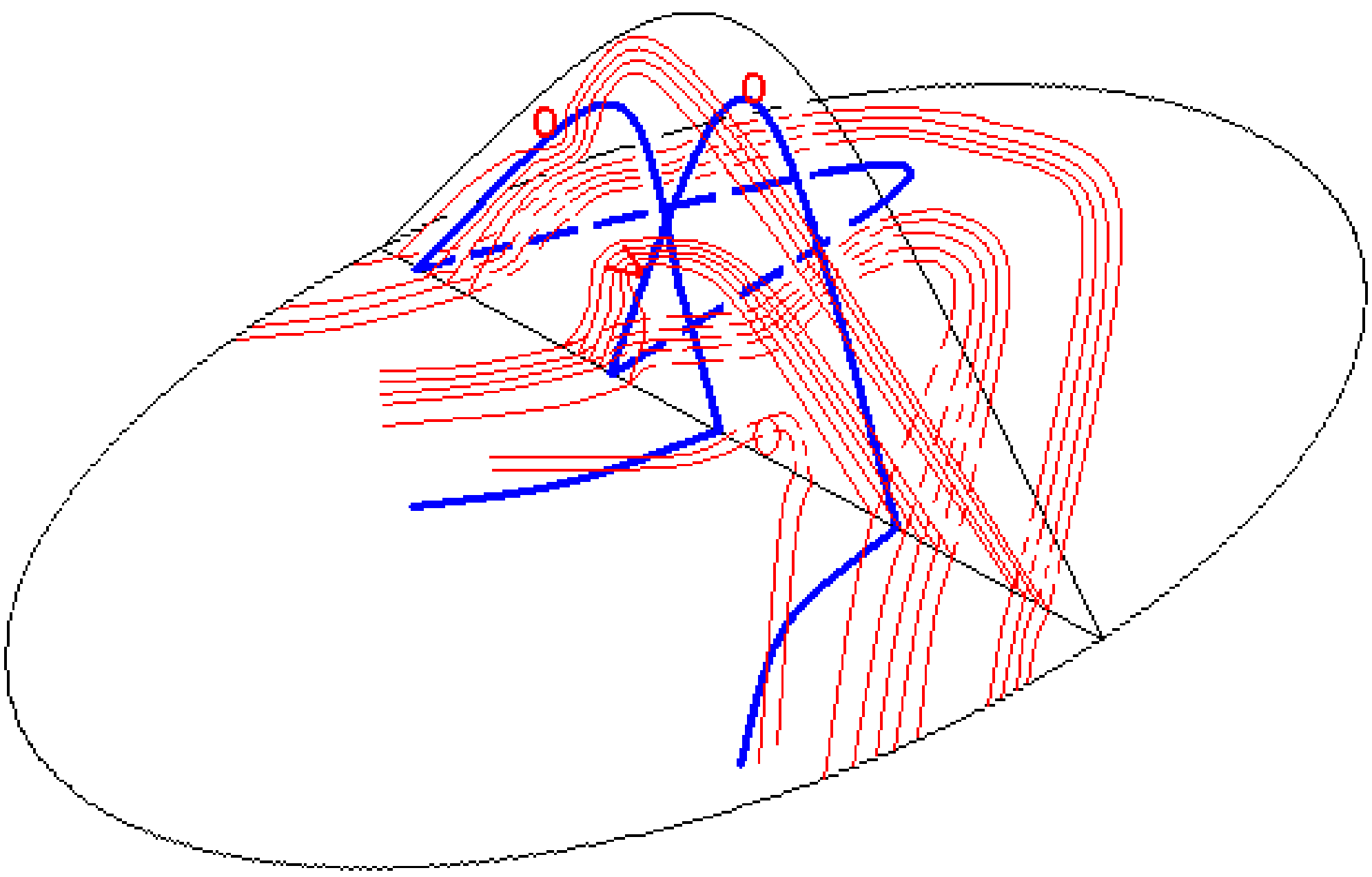

**Figure 129- Transfer rightside loop to leftside loop- push loop in top component- 2**

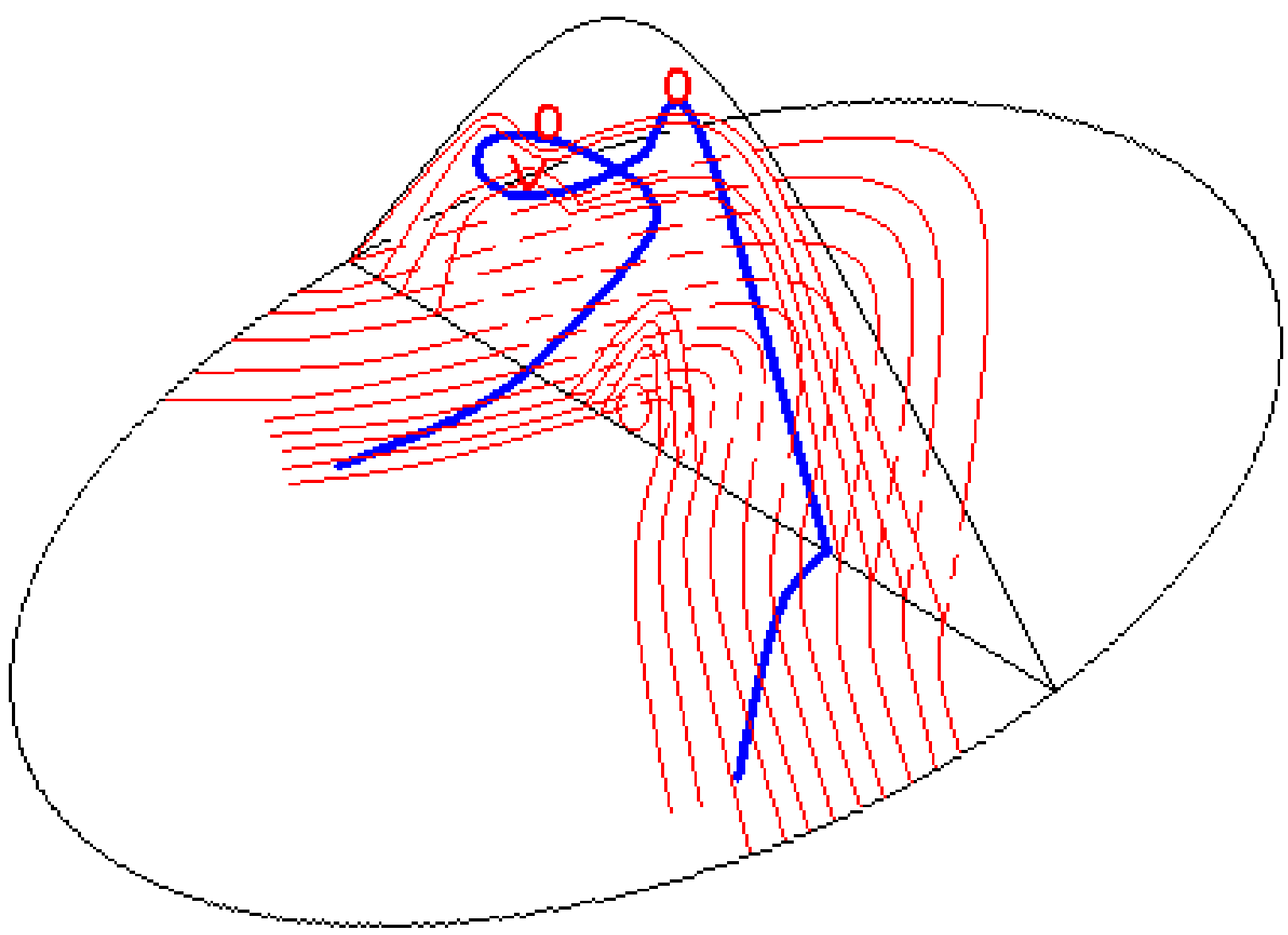

**Figure 130- Transfer rightside loop to leftside loop- push loop in top component- 3 (end)**

This picture shows the standard sequence of slices for the leftside loop and our work is done.

## *5.2 Cancellation of a loop pair- general case*

In the previous chapter we have performed the cancellation of a loop pair in the case, where the vertex model has the sequence of slices for a saddlepoint. We wish to loosen that restriction now and generalize the special case.

### 5.2.1 Pullback subsolution from special case

Start with a pair of leftside loops embedded in the Quinn model, such that the vertex has line segment slices. We change the position of the loop pair a little bit:

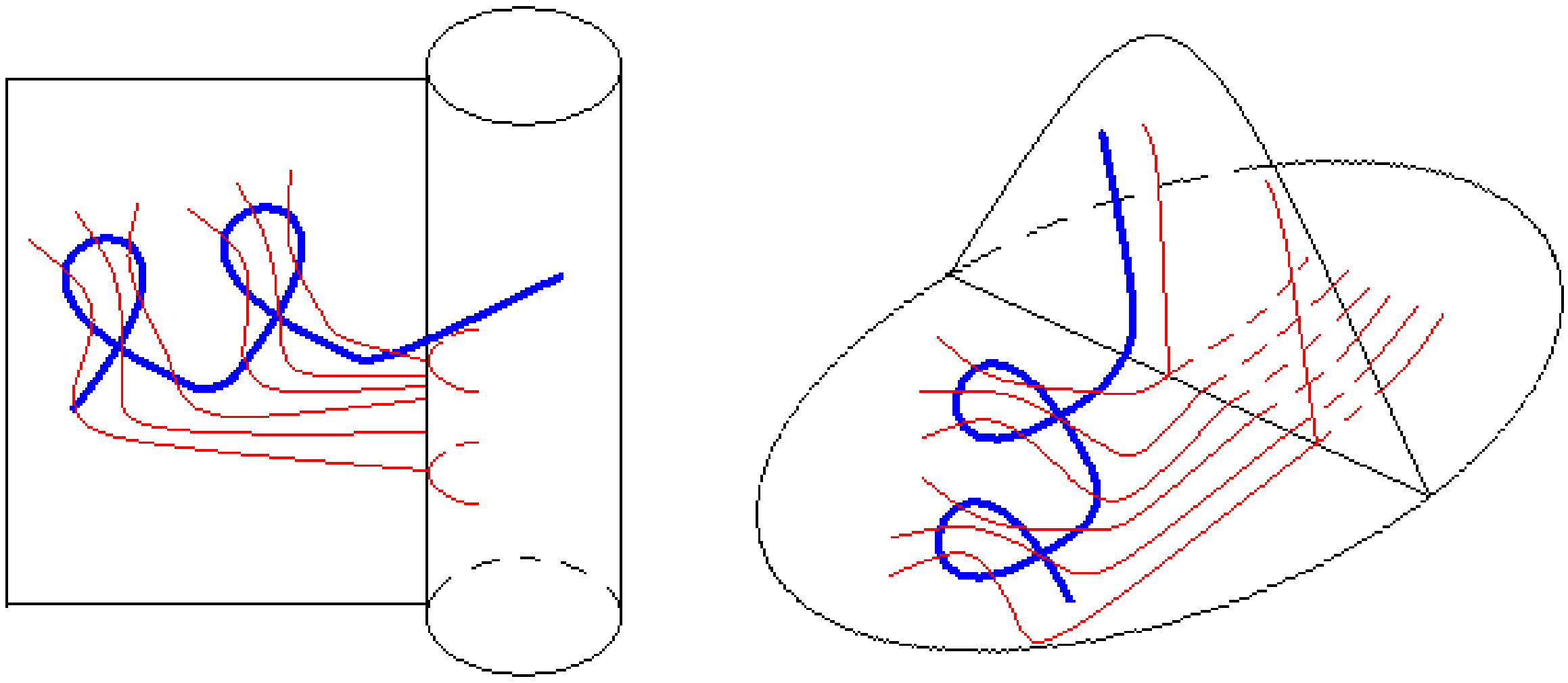

**Figure 131- Pullback subsolution from special case- a pair of leftside loop**

We want to replace the straight attaching curve by a wave curve and locally get the start figure of the special case, where we can apply the cancellation of the loop pair:

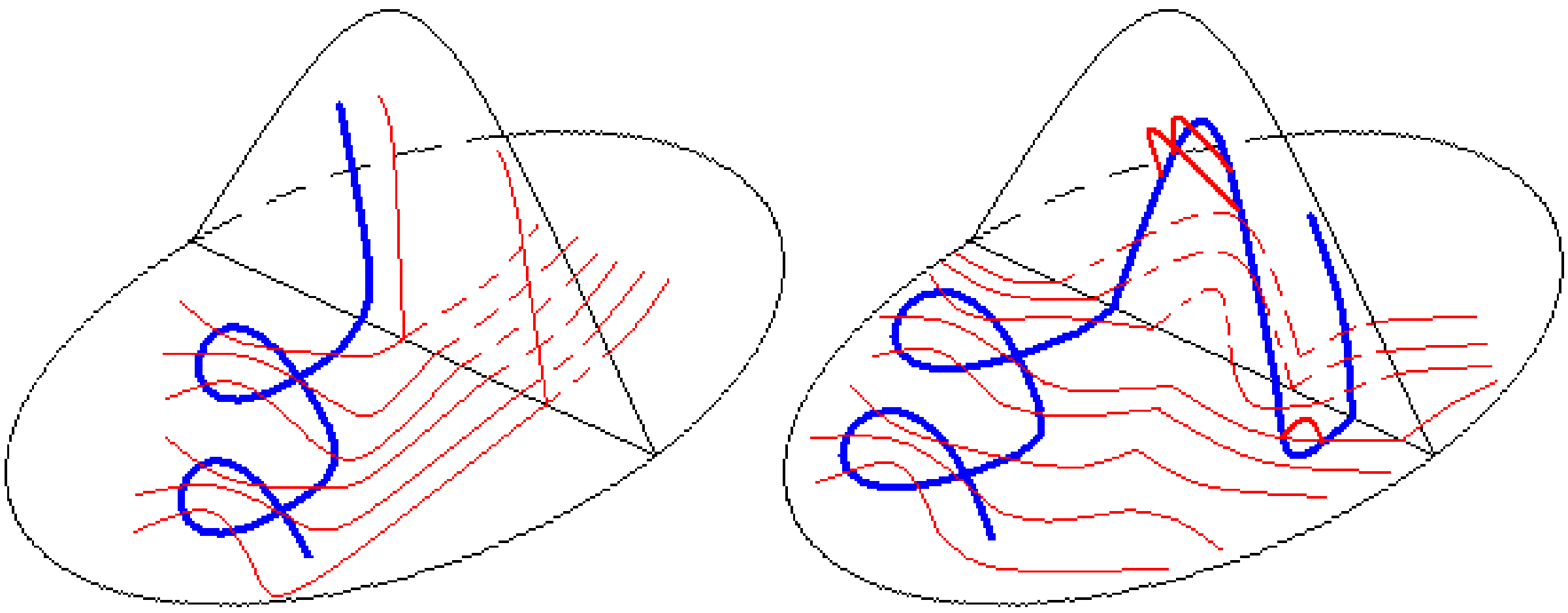

**Figure 132- Pullback subsolution from special case- preview- transfer to use starfigure from subsolution of special case**

We consider the local case:

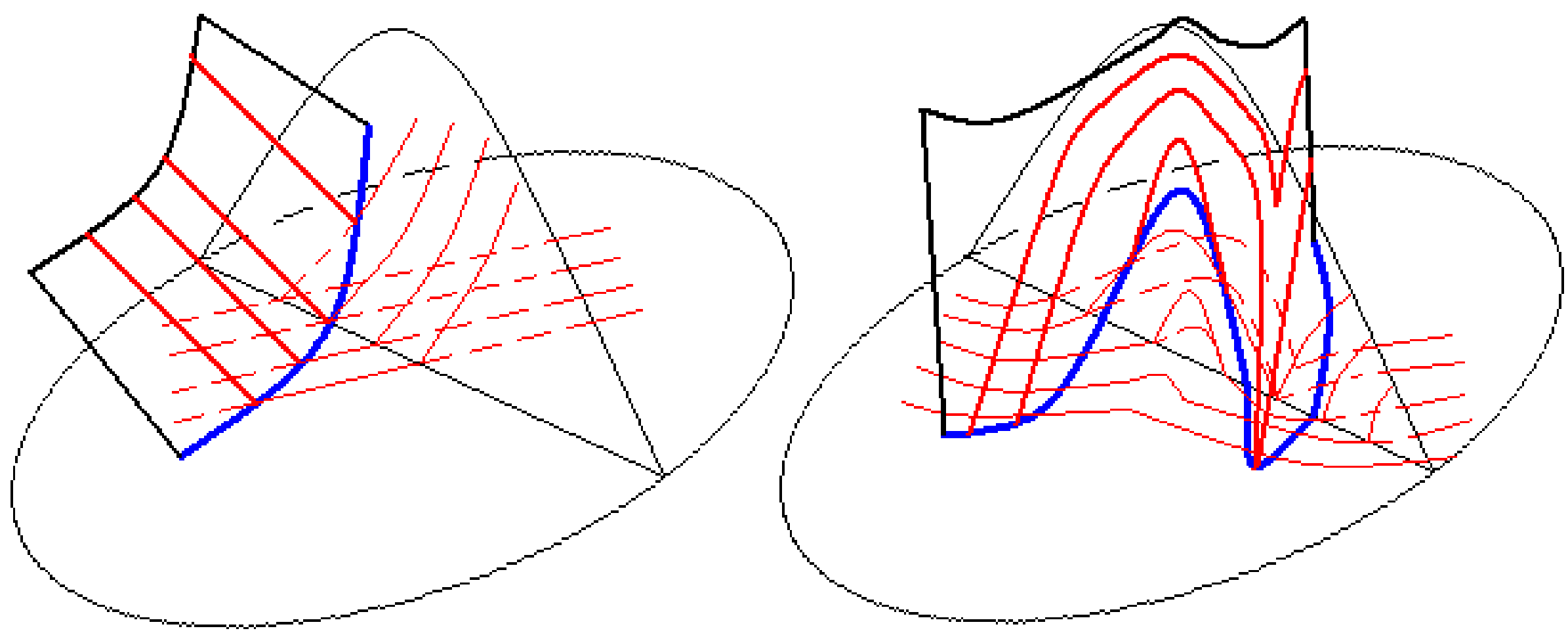

**Figure 133- Pullback subsolution from special case- preview - reduce to local case**

We arrange it as follows: First we introduce a pair of saddlepoints on the attached component:

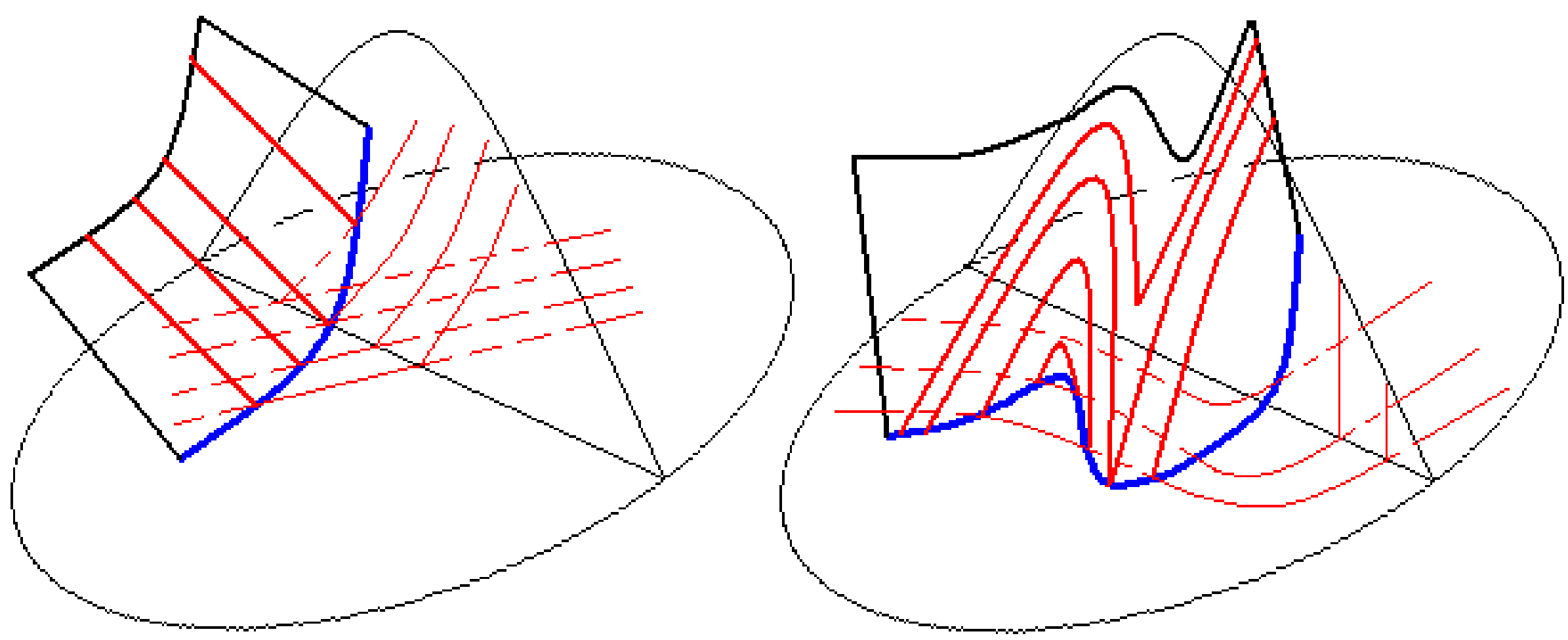

**Figure 134- Pullback subsolution from special case- reduce tansfer to local case- 1**

The next step is to change from line segment slices to slices for a pair of saddlepoints at the vertex:

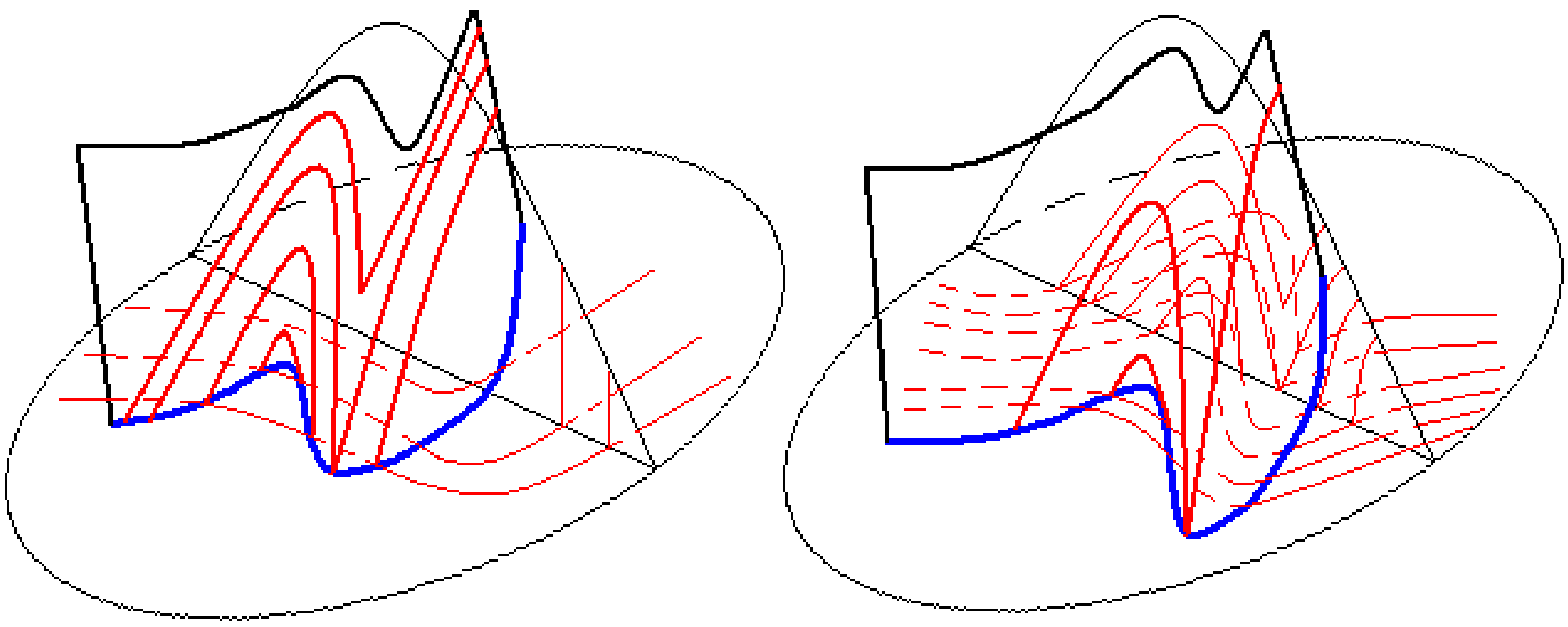

**Figure 135- Pullback subsolution from special case- reduce tansfer to local case- 2**

It remains to push a part of the curve from the frontside component across the vertex into the top component and get our desirable result. From that point on we can perform the cancellation of a loop pair and after that rechange the sequence of slices to line segment slices.

When there is a pair of rightside loops, we transfer it to a pair of leftside loops, using the previous construction and then apply the cancellation for the special case and rechange the slices:

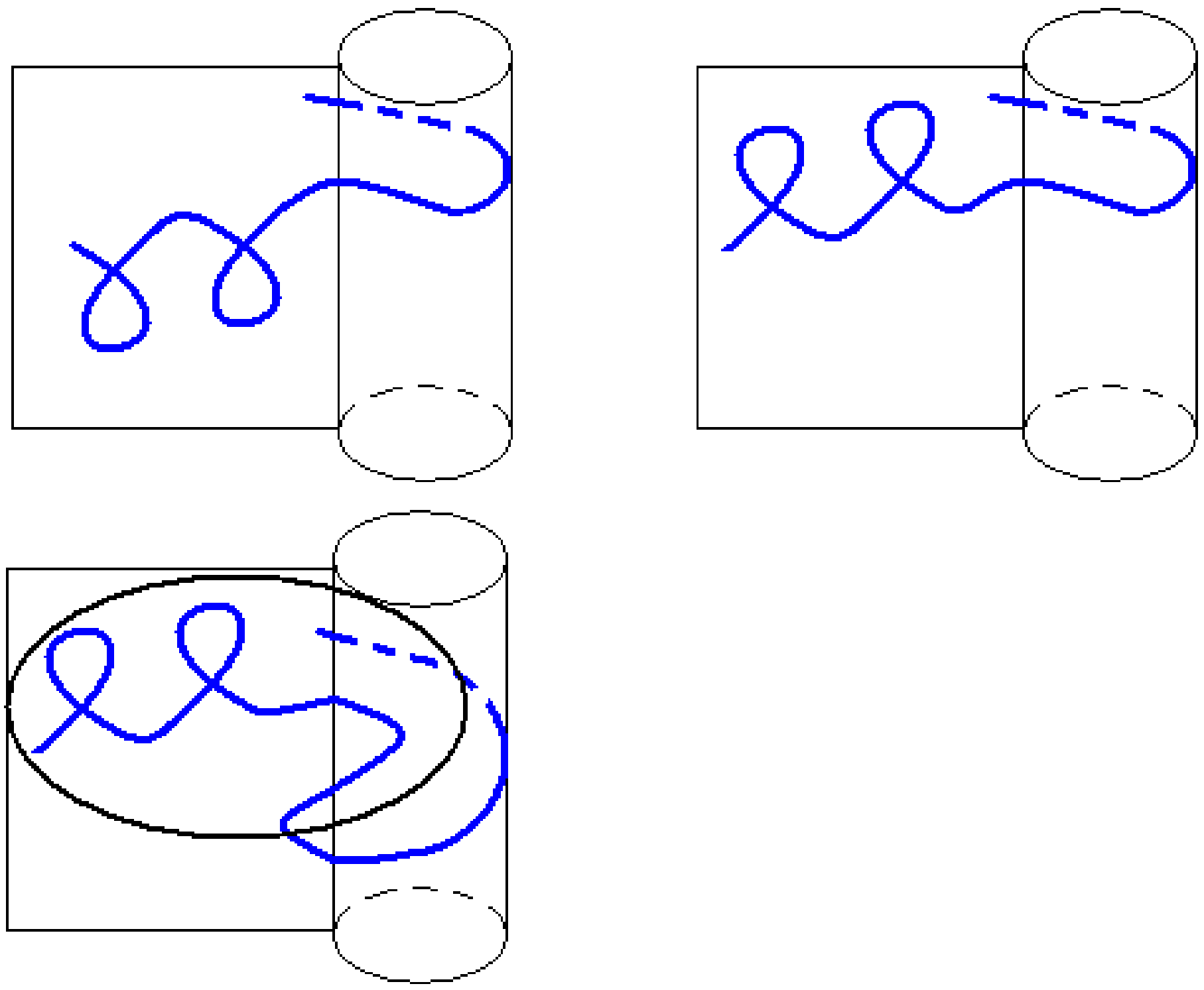

**Figure 136- Pullback subsolution from special case- a pair of rightside loops**

### 5.2.2 Shift a loop

We finish this section with the shift of a loop. Assume there are two separated loops at the attaching curve. We can shift them next to each other and then annihilate them:

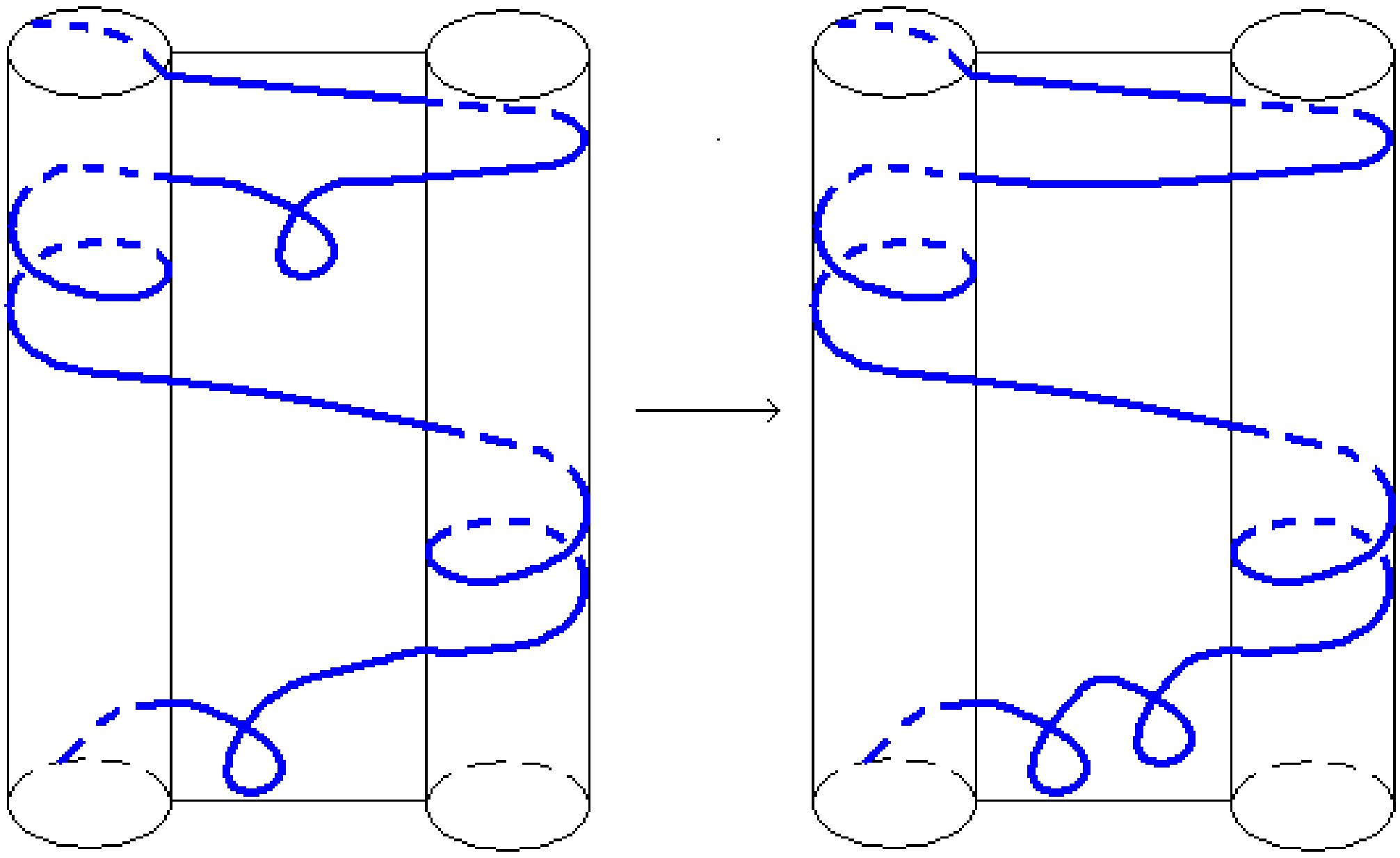

**Figure 137- Shift a loop- shift two separated loops together**

There are two cases to consider:

- passing at the generator cylinder
- crossing of the attaching curve

but we can summarize them to a single local case:
We have to change the line segment slices to slices of an introduced pair of saddlepoints and then we push the loop across the vertex:

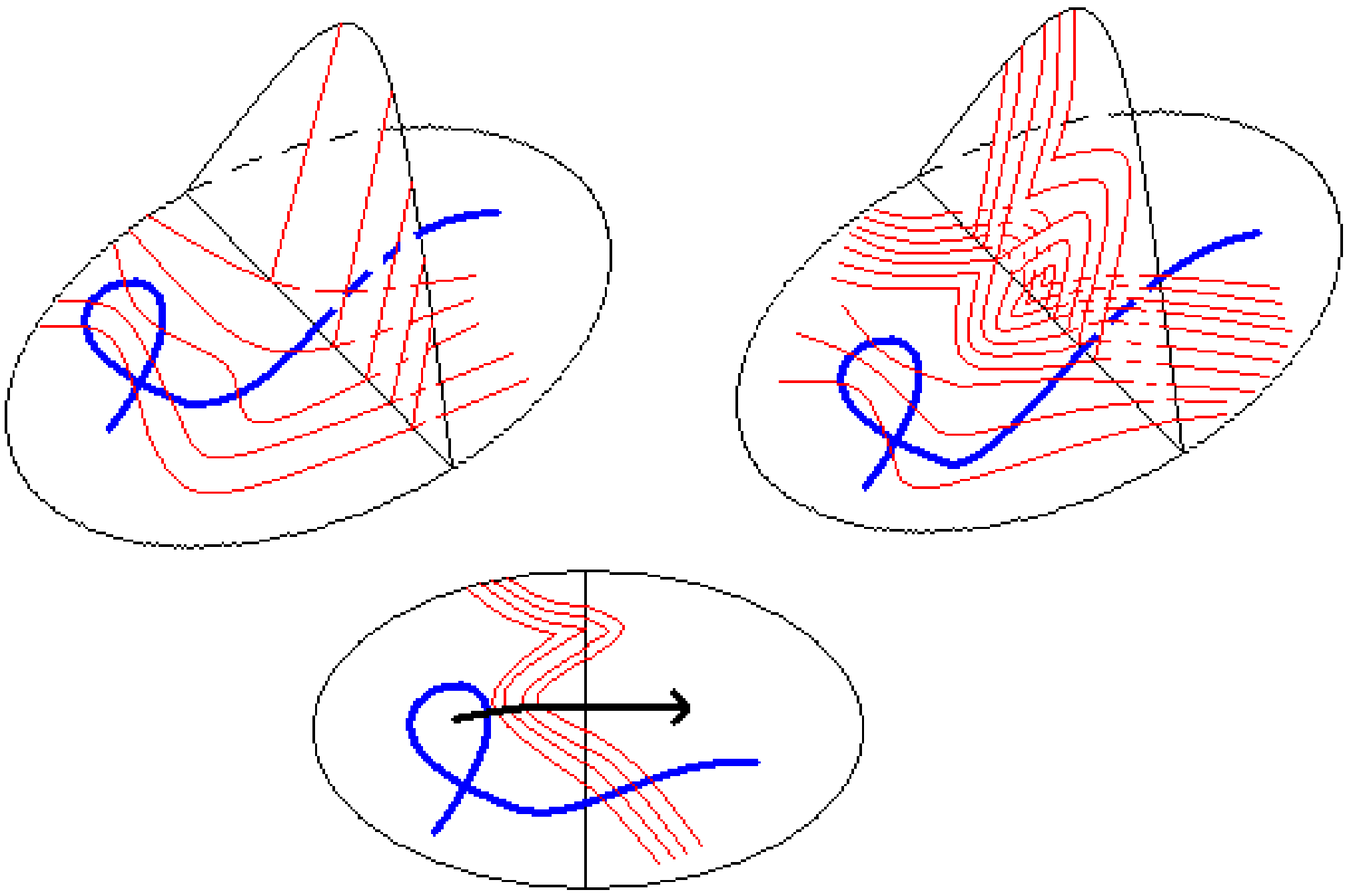

**Figure 138- Shift a loop- reduce to local case**

# 6 Different ways of performing twists

## *6.1 Twists on the generator cylinder*

In chapter 3 we dealt with the appearance of twists in two ways. Now we construct the twists at the generator cylinder and look for the composition of these. Note, that the turn around the generator cylinder will always be fixed so we omit this part by looking for the sequence of slices. We do not describe each step, but only the “way”.

a) The first case result in a leftside loop:

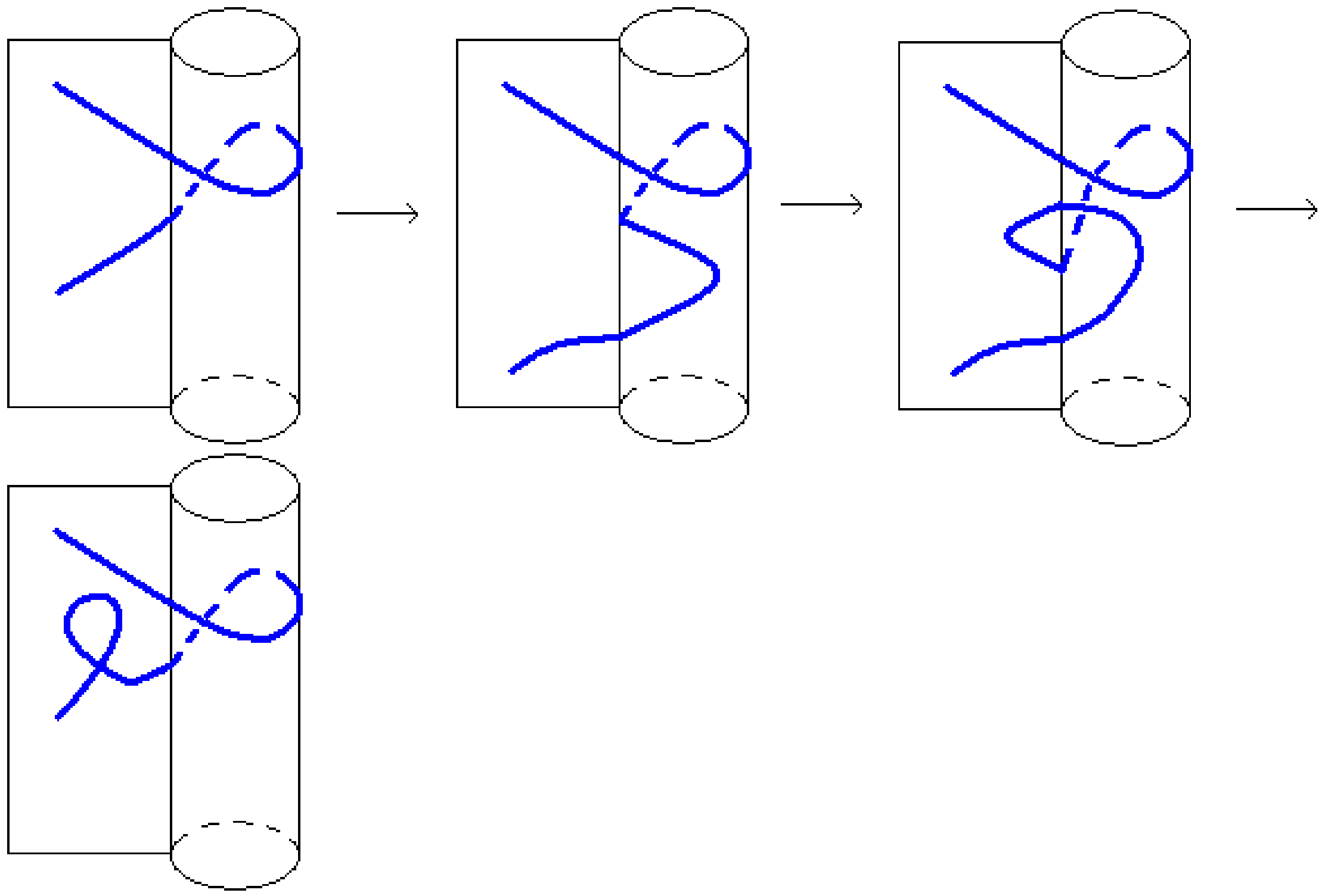

**Figure 139- Twists on the generator cylinder- Case a)**

We start with the first twist. To perform the second one, we have to introduce a pair of saddlepoints to avoid a “bad $T_3$ turn” in the next step:

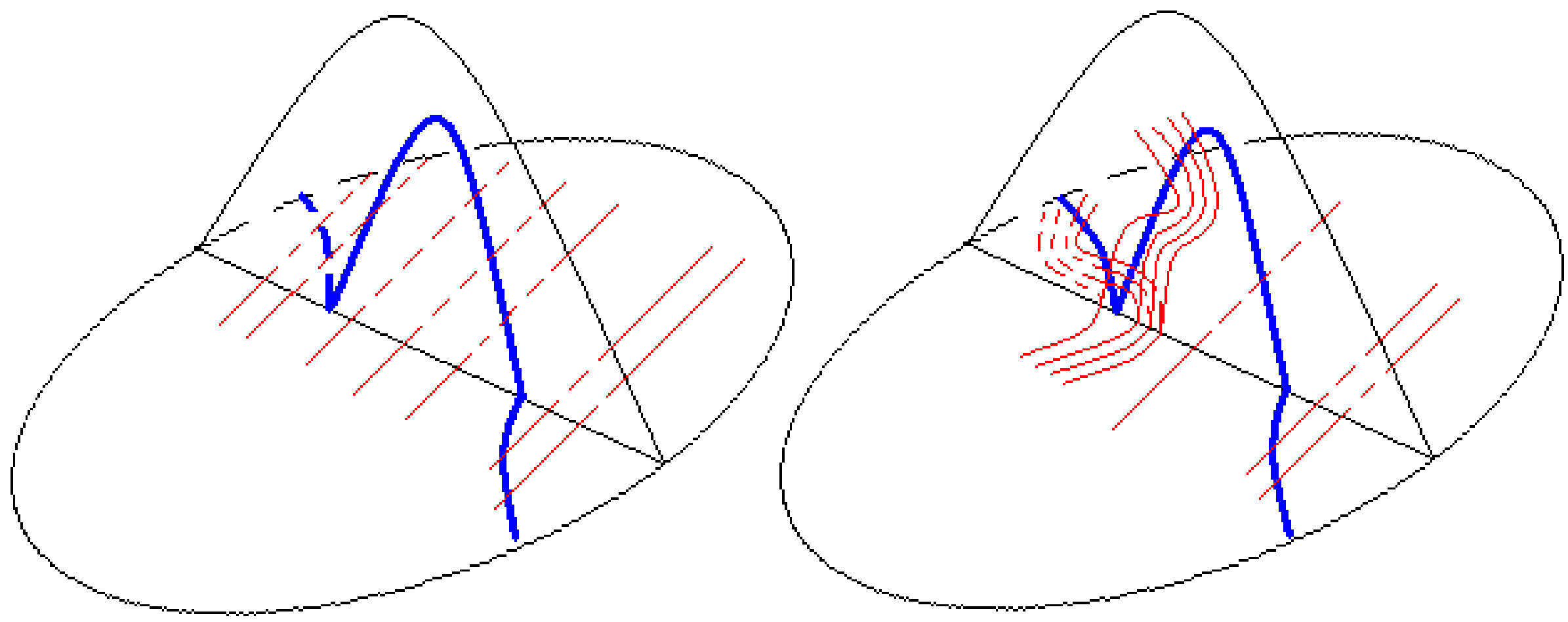

**Figure 140- Twists on the generator cylinder- Case a) local- 1**

We perform the twist and in the next steps we shift the saddlepoints to the frontside component to get a loop there:

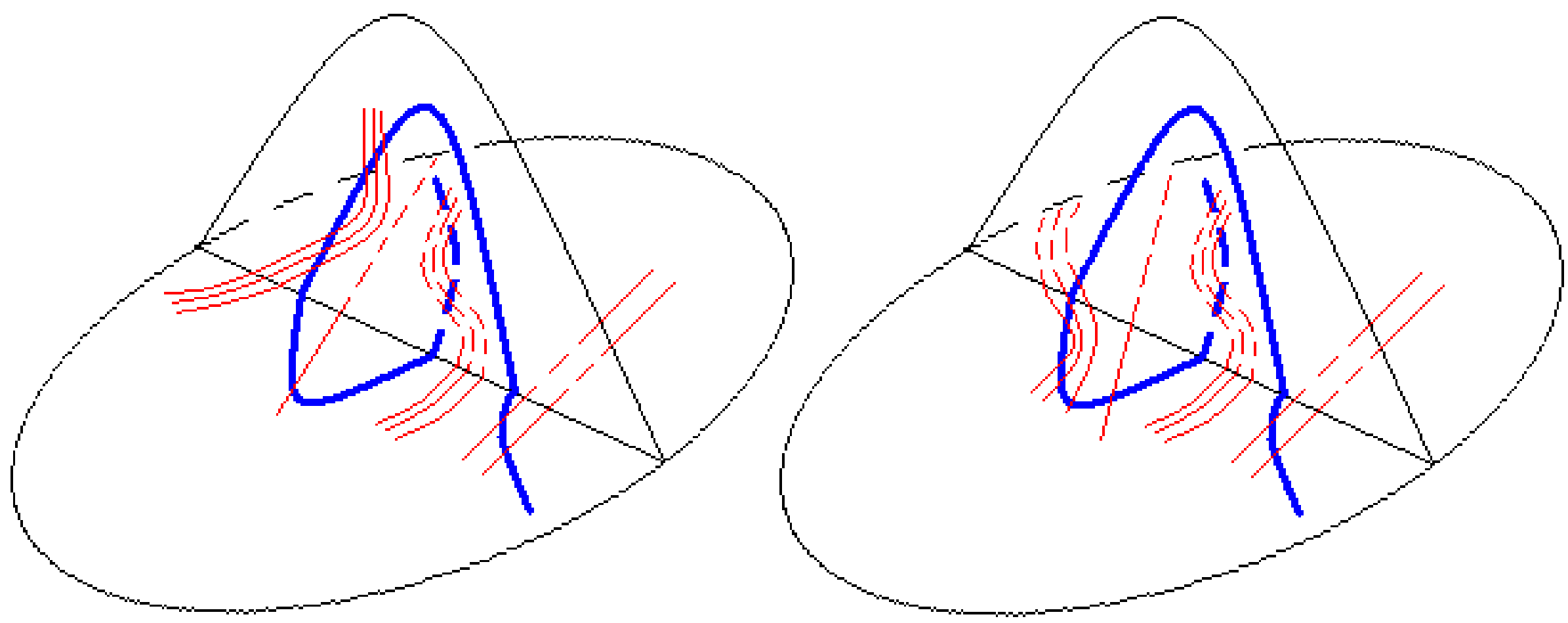

**Figure 141- Twists on the generator cylinder- Case a) local- 2**

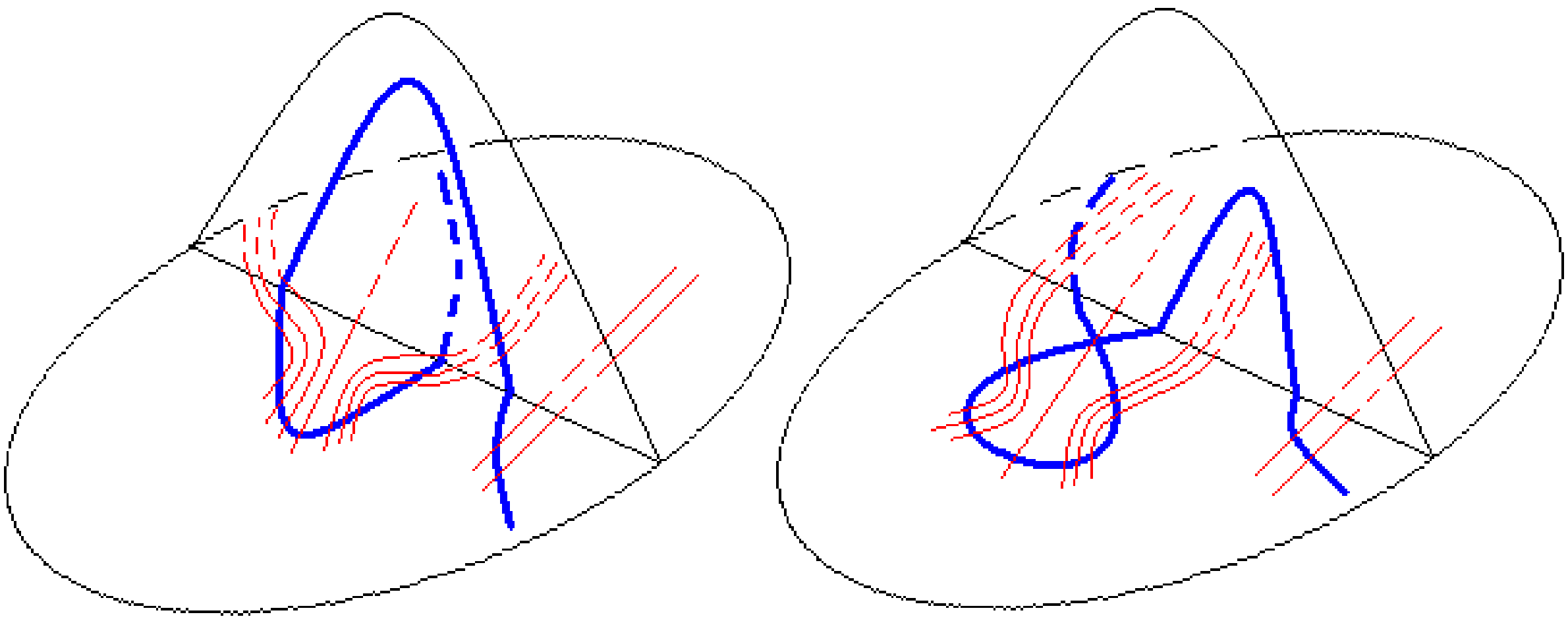

**Figure 142- Twists on the generator cylinder- Case a) local- 3**

b) The next sequence describes a rightside loop:

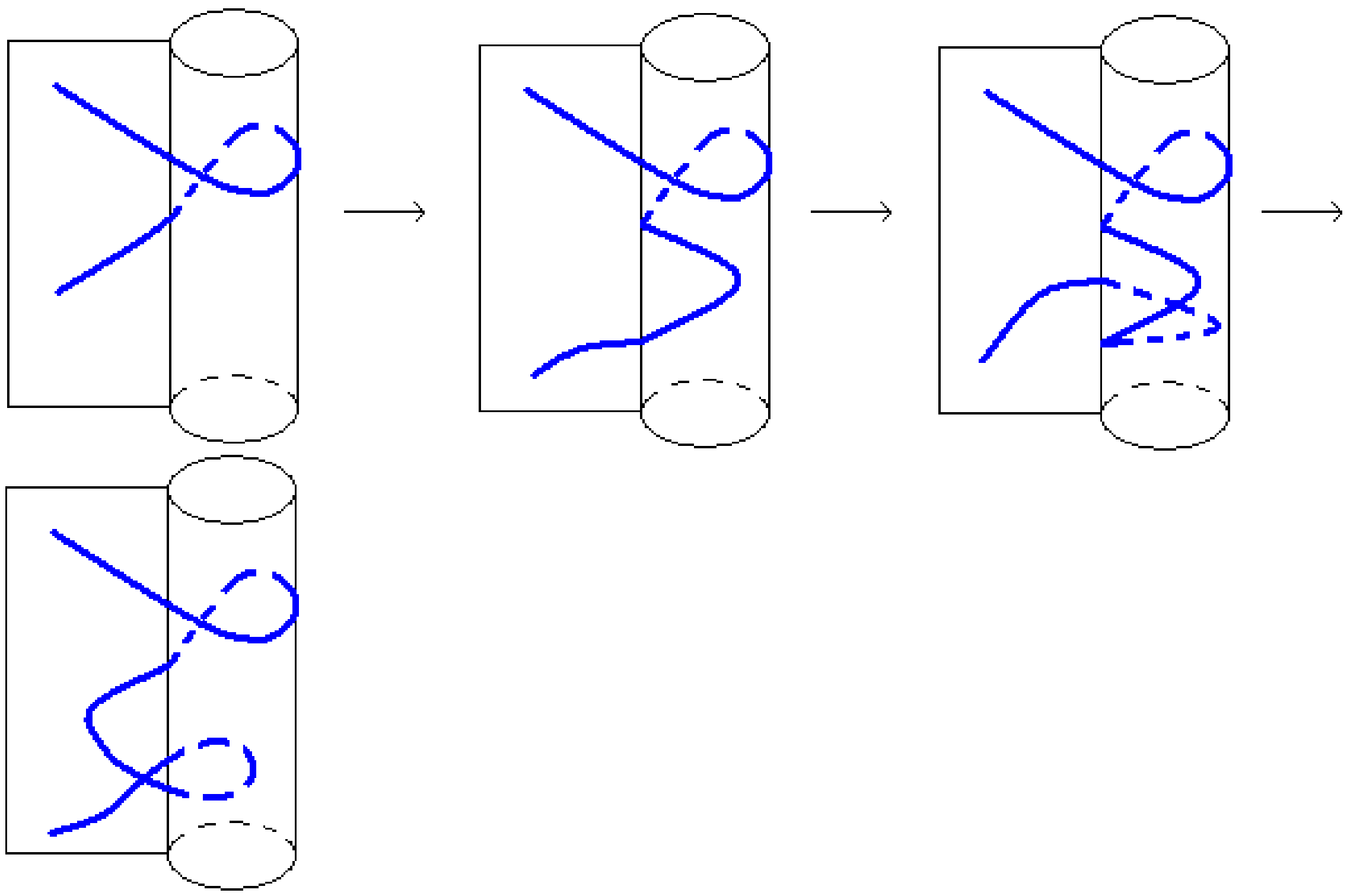

**Figure 143- Twists on the generator cylinder- Case b)**

We introduce a pair of saddlepoints to perform the next twist as a"good" one:

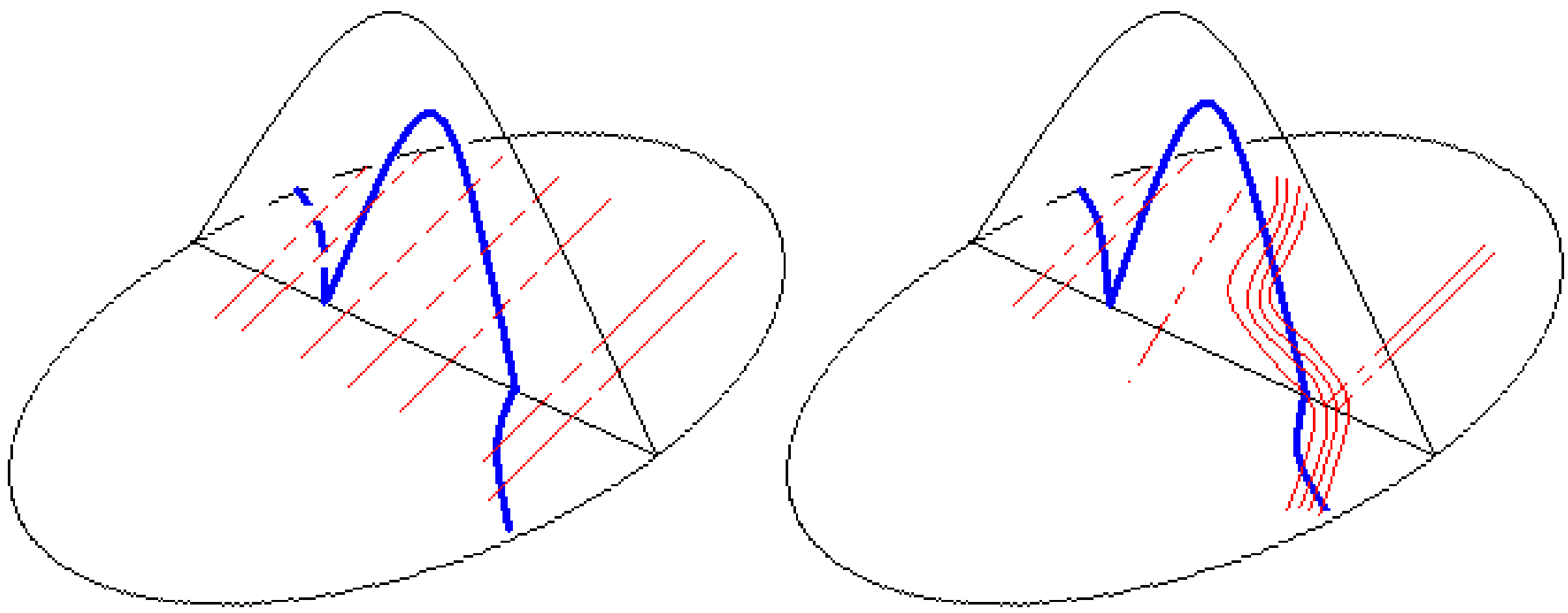

**Figure 144 Twists on the generator cylinder- Case b) local- 1**

We perform the twists and in the next steps we shift the saddlepoint into the backside component and construct the loop there. Then we shift the saddlepoints to the frontside component, so we can push the line segment sliced-turn of the backside component into the frontside component, as shown in the next figures:

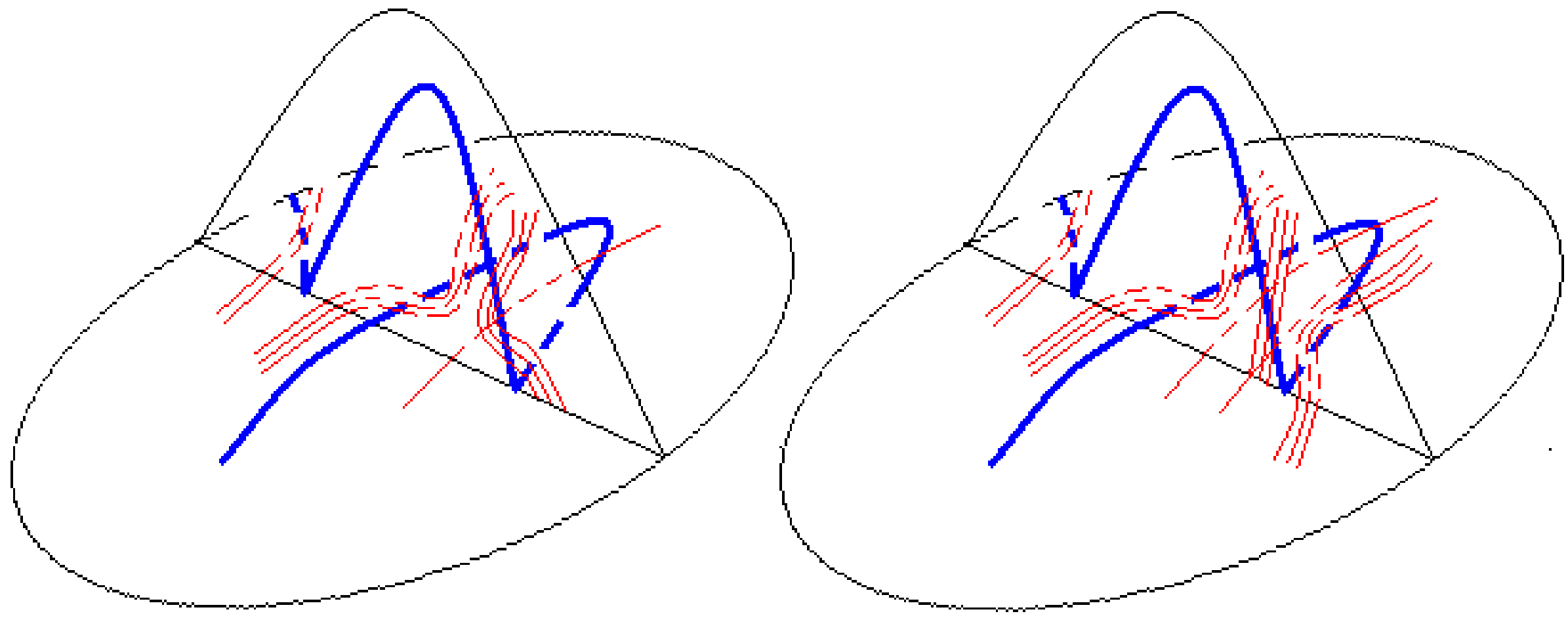

**Figure 145- Twists on the generator cylinder- Case b) local- 2**

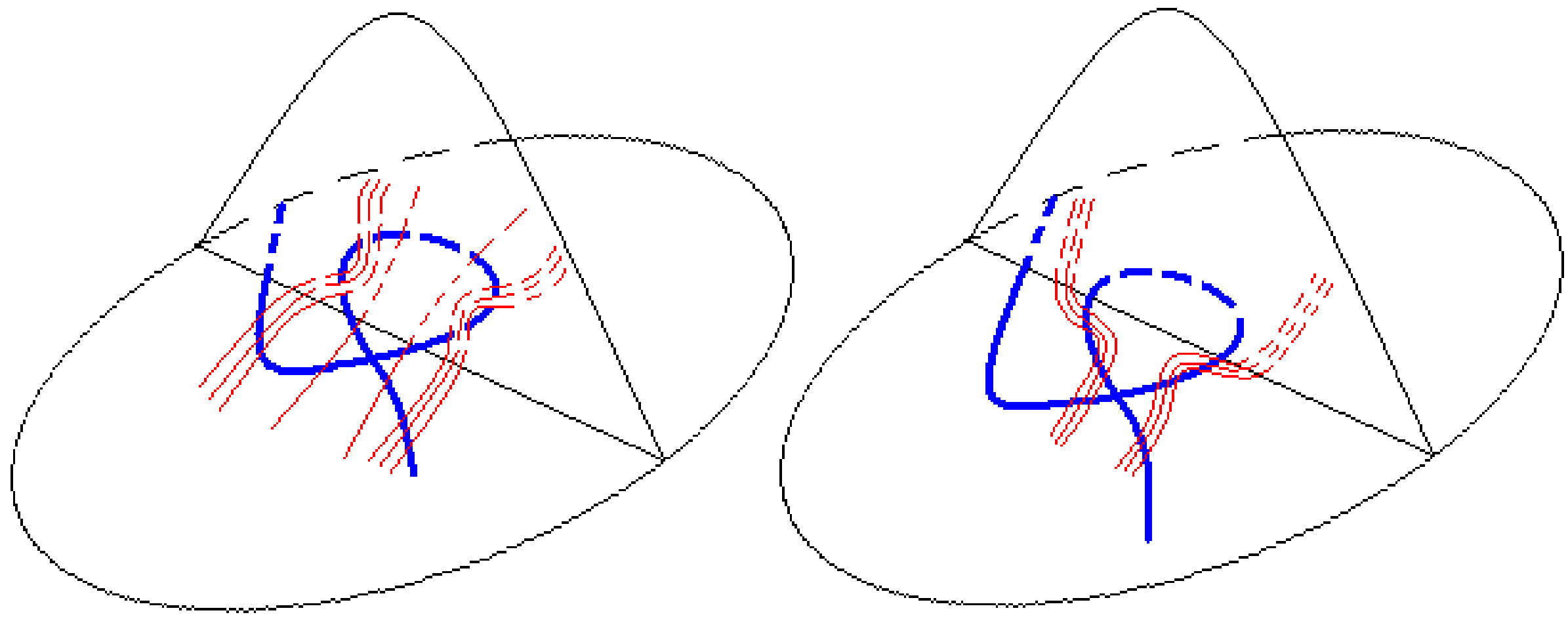

**Figure 146- Twists on the generator cylinder- Case b) local- 3**

c) The next sequence leads to a leftside loop:

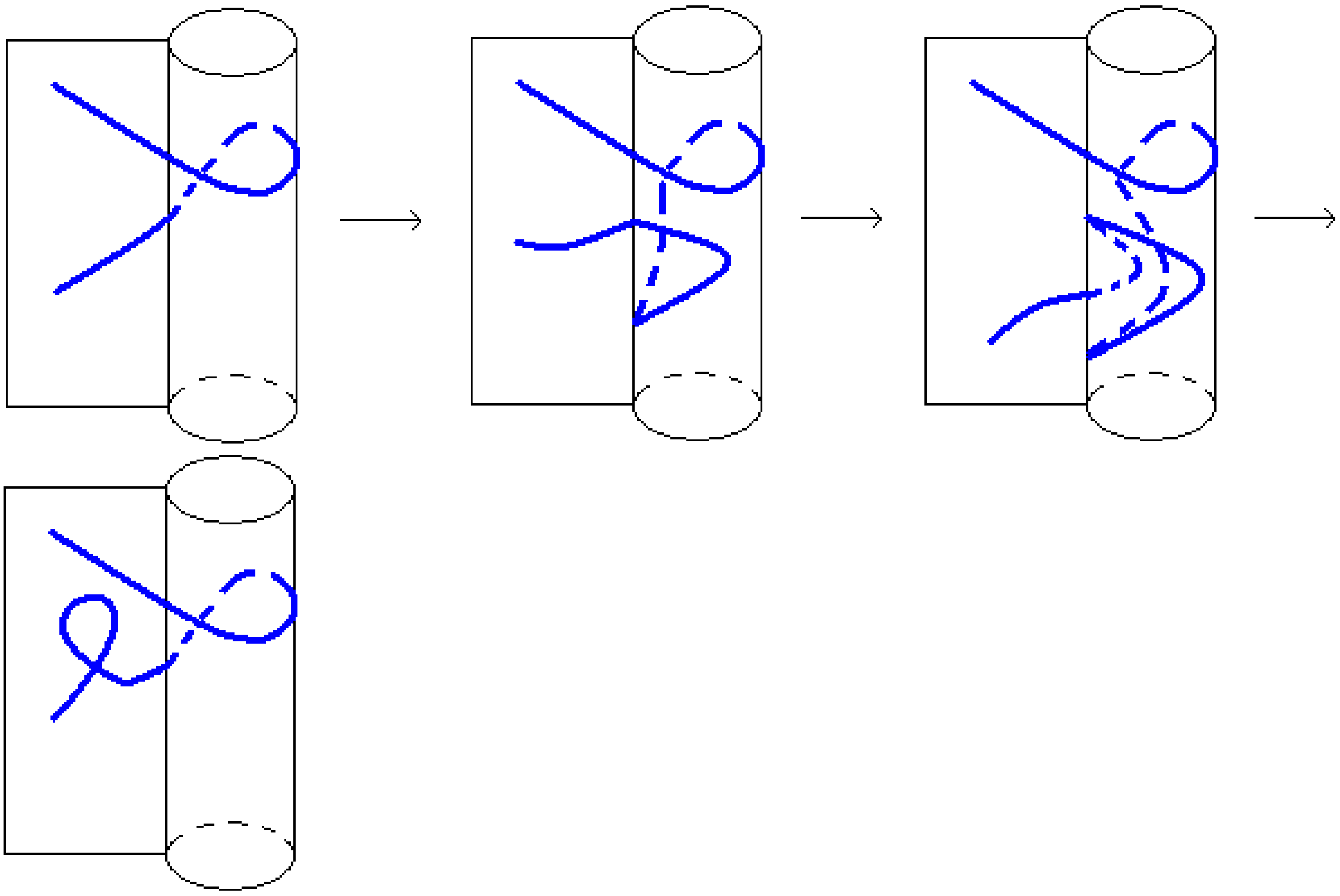

**Figure 147- Twists on the generator cylinder- Case c)**

To perform the first twist, we have to change the situation into slices with a pair of saddlepoints, otherwise it would be a “bad $T_3$ turn”. We shift one saddlepoint from the frontside component into the top component:

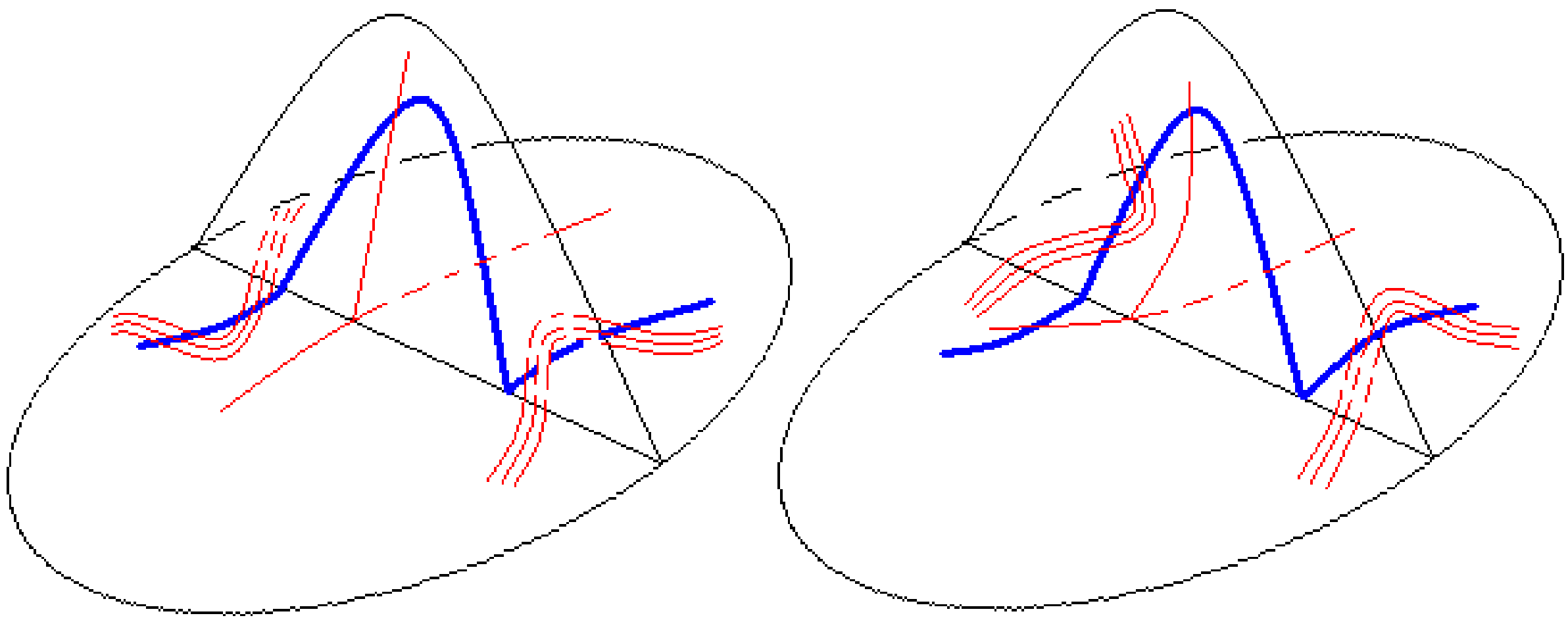

**Figure 148- Twists on the generator cylinder- Case c) local- 1**

Now we can perform the second twist as a “good $T_3$ turn”. We shift the saddlepoint from the top component into the backside component:

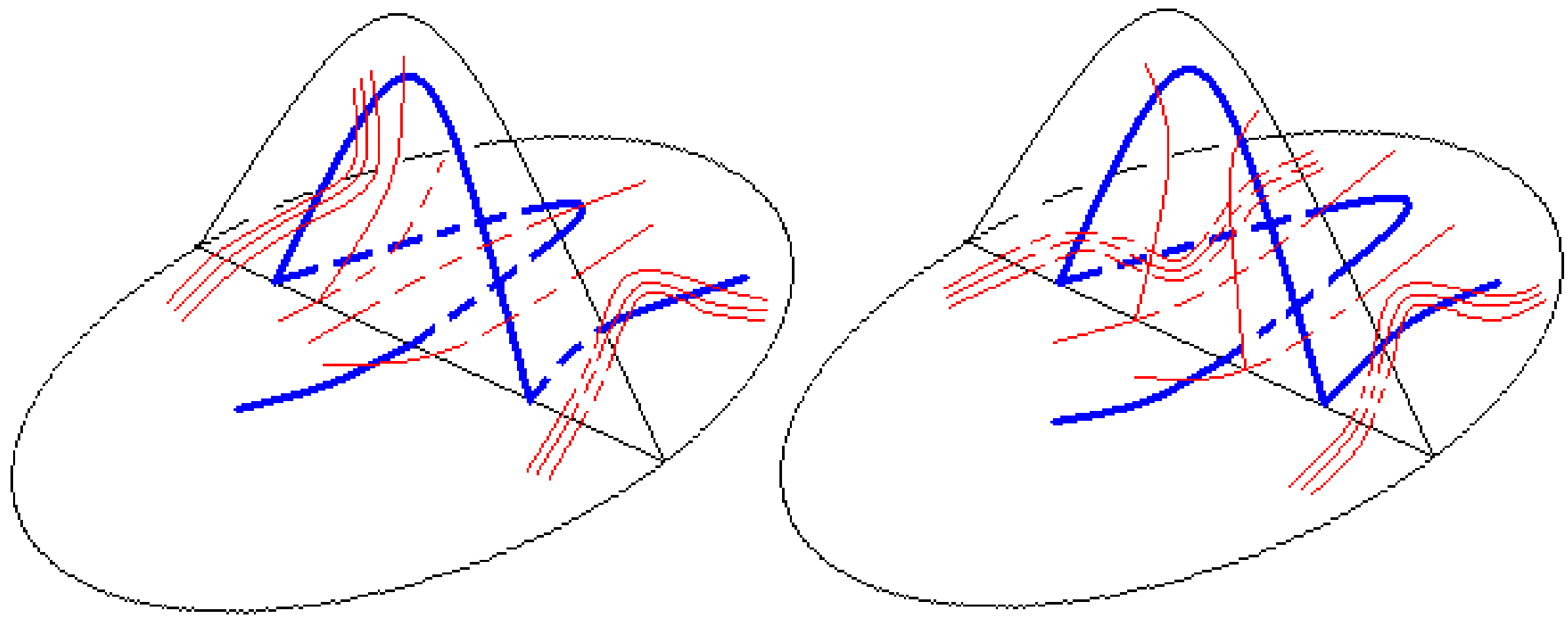

**Figure 149- Twists on the generator cylinder - Case c) local- 2**

Hence we can use the turn with line segment slices in the top component to construct the loop in the frontside component. Then we push the saddlepoint into the frontside component:

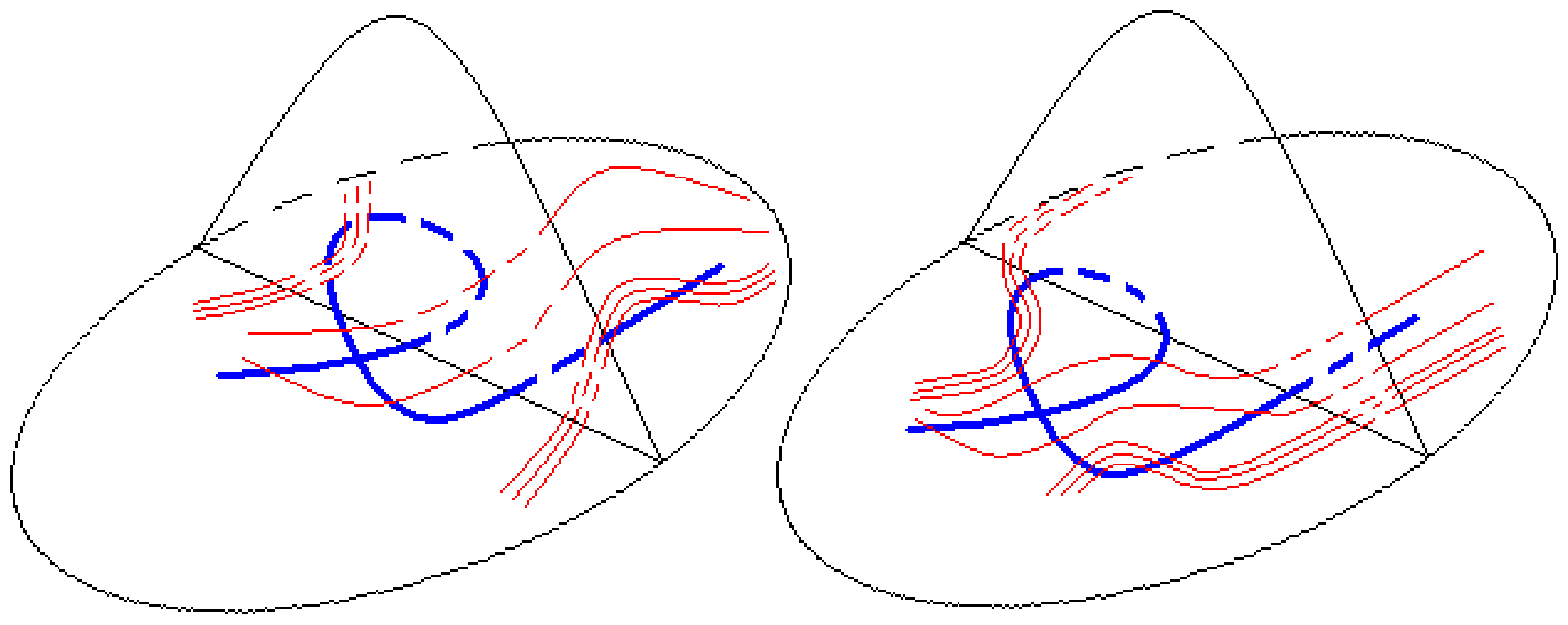

**Figure 150- Twists on the generator cylinder- Case c) local- 3**

We change the slices to a standard loop slices by introducing a pair of saddlepoints to each saddlepoint and get the result shown in the right picture:

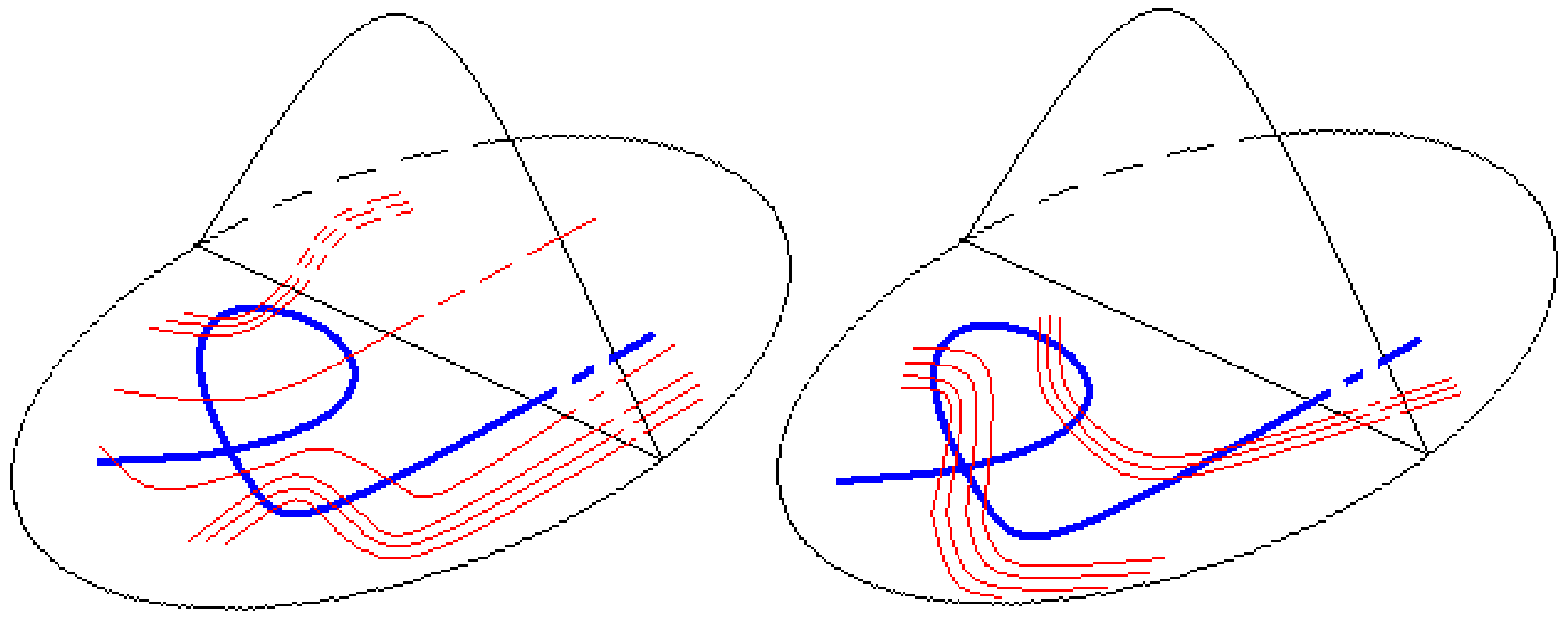

**Figure 151- Twists on the generator cylinder- Case c) local- 4**

d) The last case, which results in a loop, is the construction of a rightside loop:

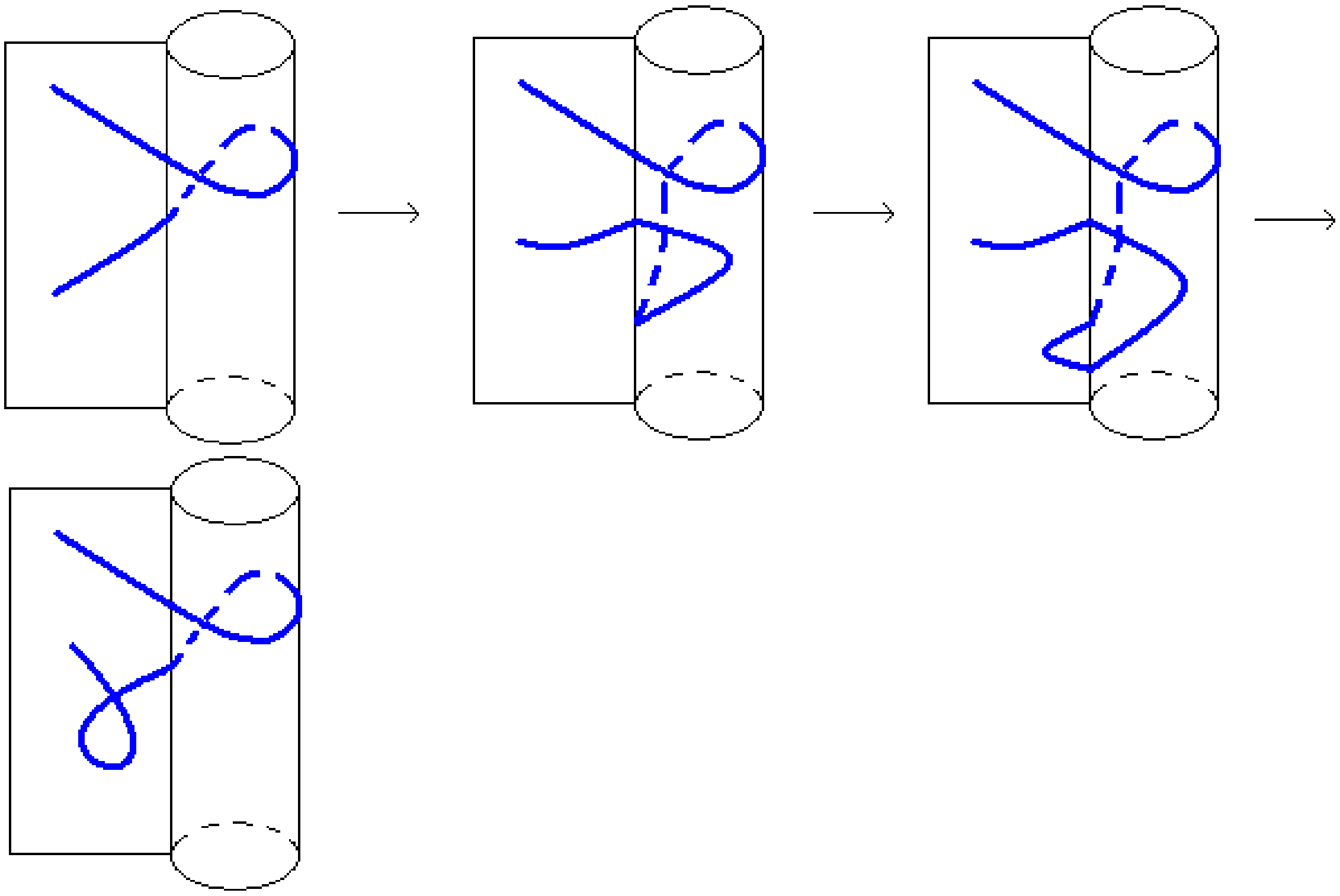

**Figure 152- Twists on the generator cylinder- Case d)**

For the sequence of slices we modify the steps in the case c) and get:

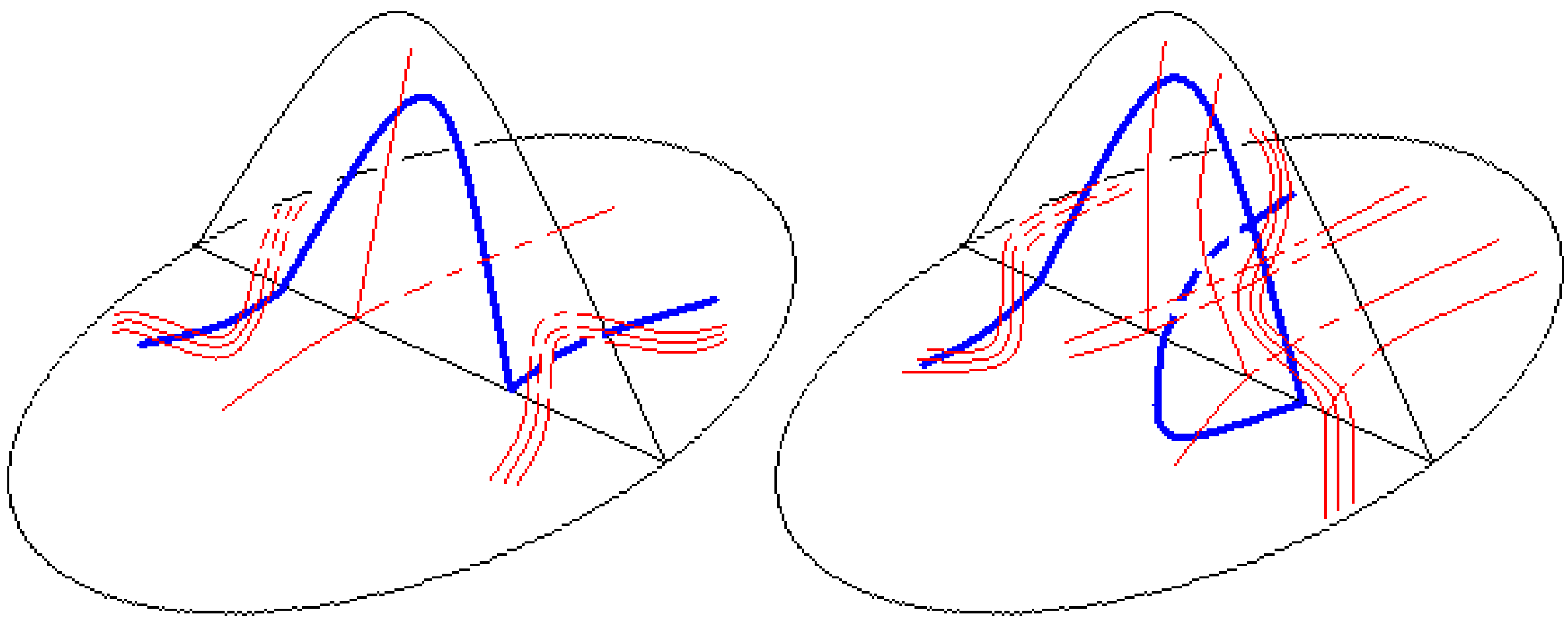

**Figure 153- Twists on the generator cylinder- Case d) local- 1**

We compare the right picture with the left picture of:

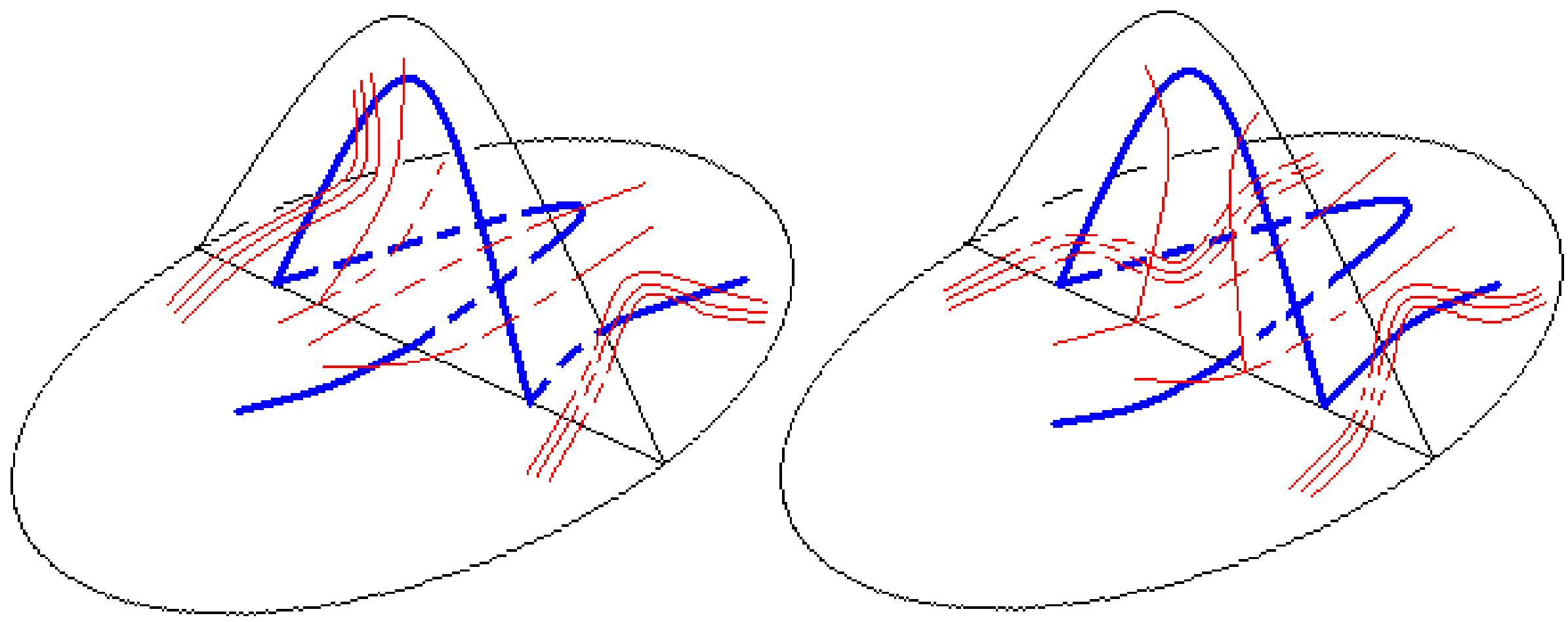

**Figure 154- Twists on the generator cylinder - Case d) local- 2**

The step from the right picture in the case c) to the left picture in the case d) is a rotation of the local model which changes the role of frontside and backside component. It remains to transfer the steps in c) by that rotation to get steps for d)

In all the other cases the twists annihilate to a line:

We only desribe one case in detail, where the slices have to change. The others are similar or easier.

e) twists annihilate to a line:

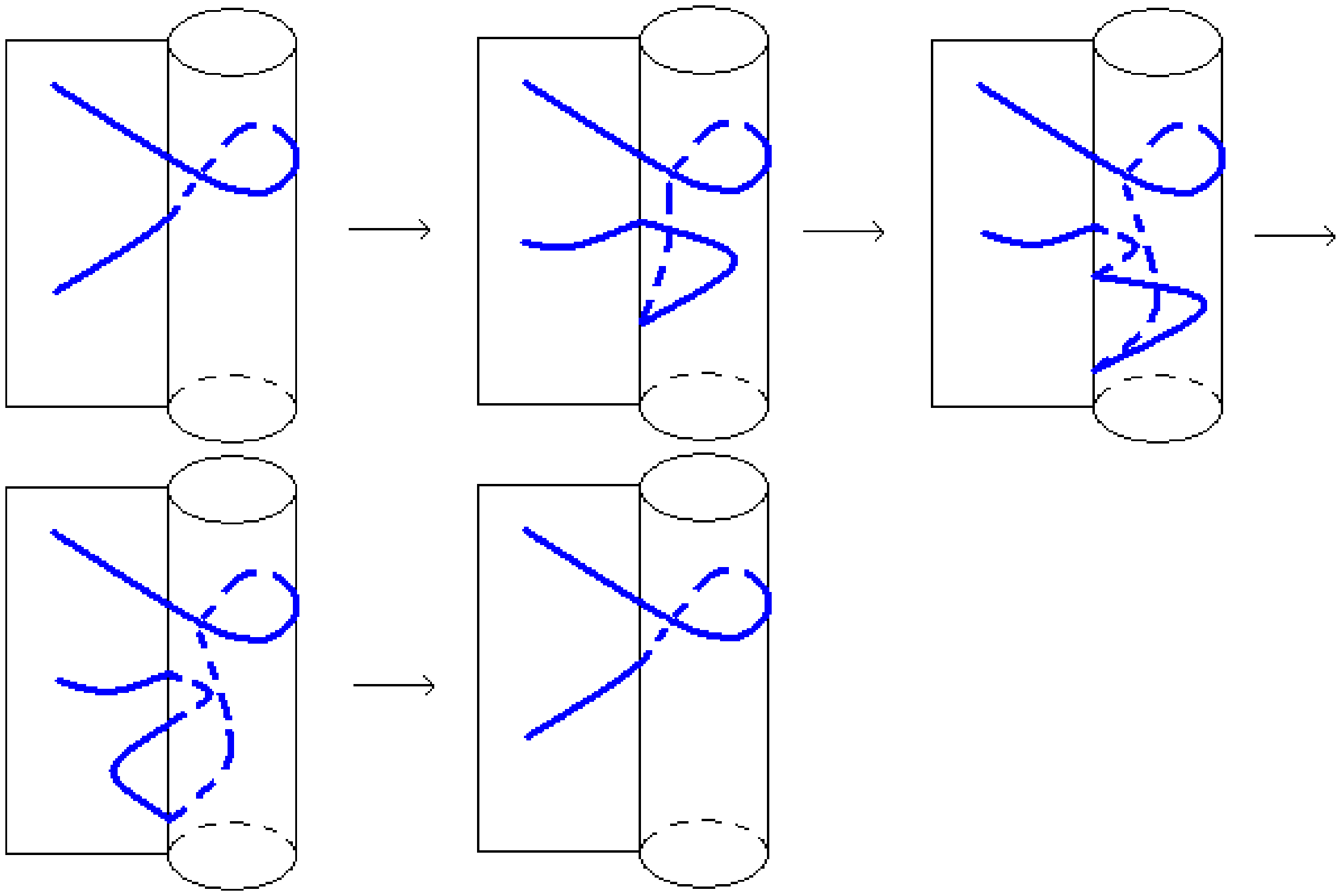

**Figure 155- Twists on the generator cylinder- Case e)**

We have to introduce a pair of saddlepoints to be able to perform the first twist. Then we perform the second one:

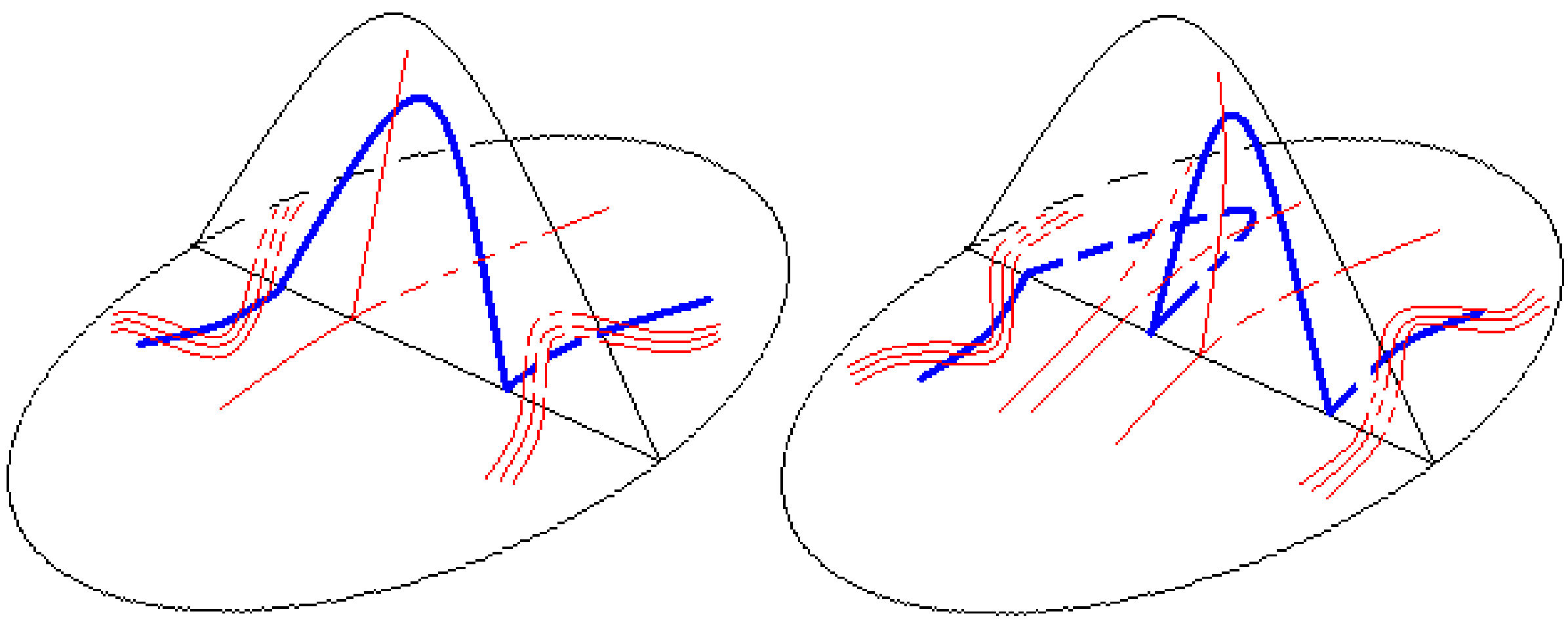

**Figure 156- Twists on the generator cylinder - Case e) local- 1**

We push the turn in the top component into the frontside component. The result is a wave:

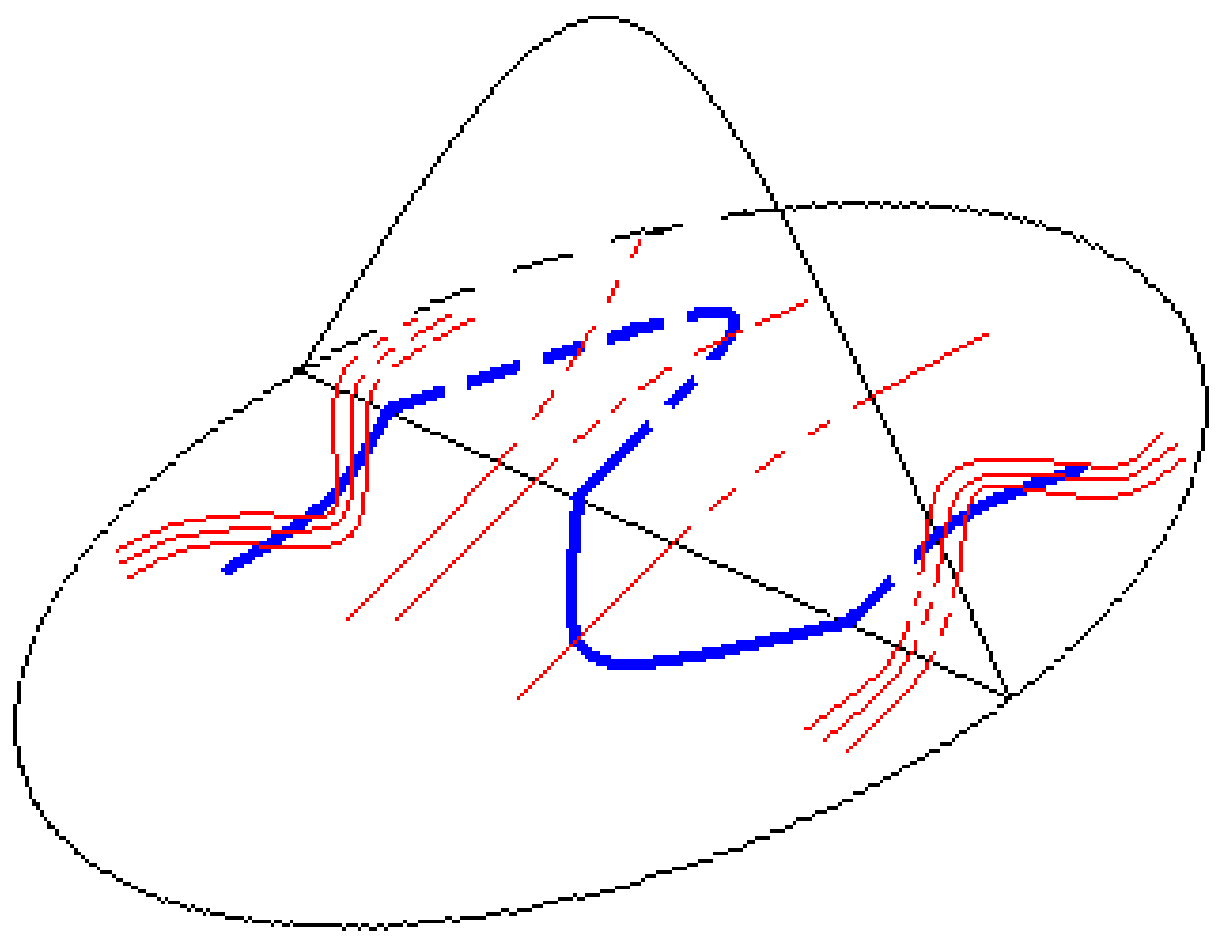

**Figure 157- Twists on the generator cylinder- Case e) local- 2**

We deform the wave to a line and annihilate the pair of saddlepoints.

The other cases for a line are:

f) twists annihilate to a line:

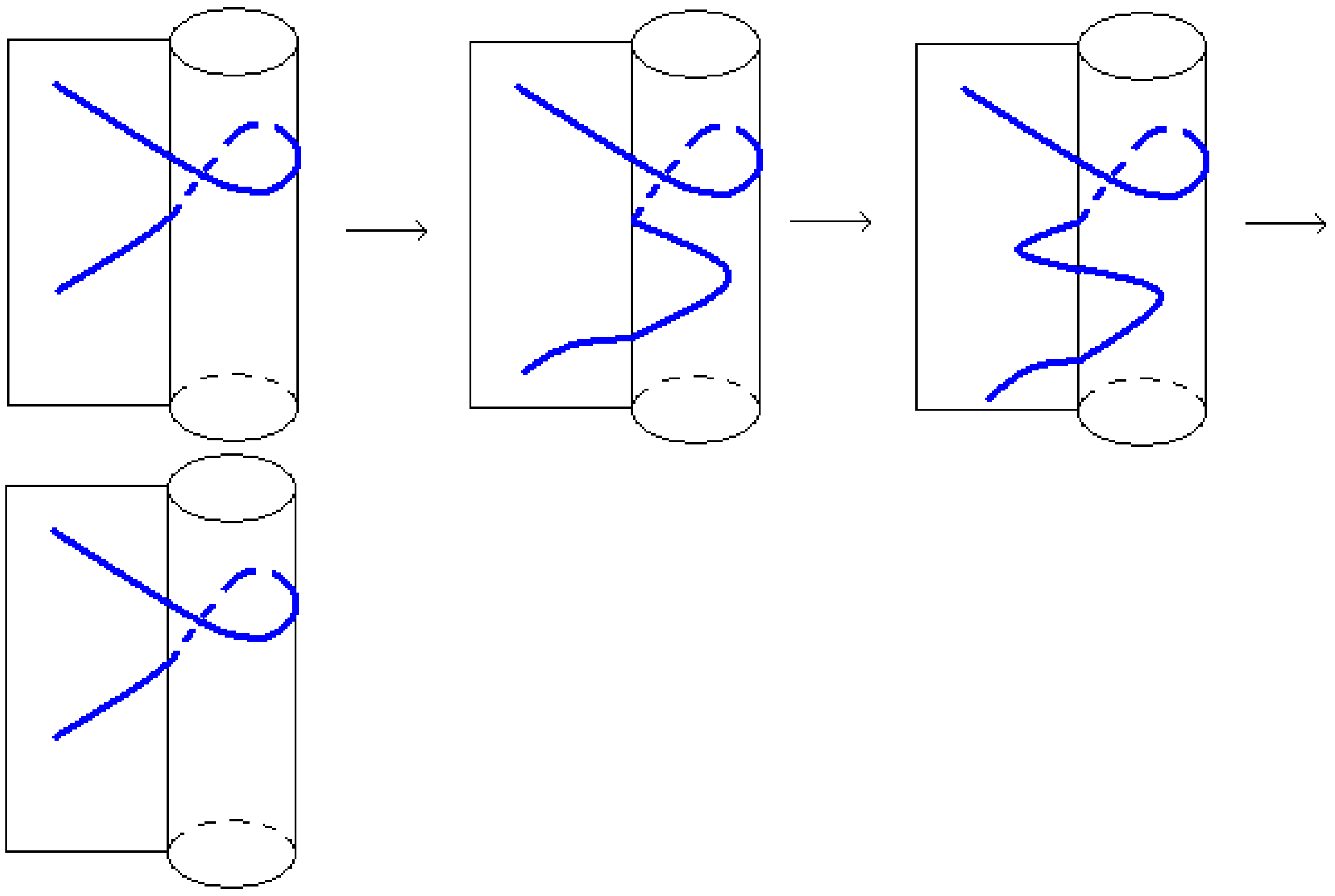

**Figure 158- Twists on the generator cylinder- Case f)**

g) twists annihilate to a line:

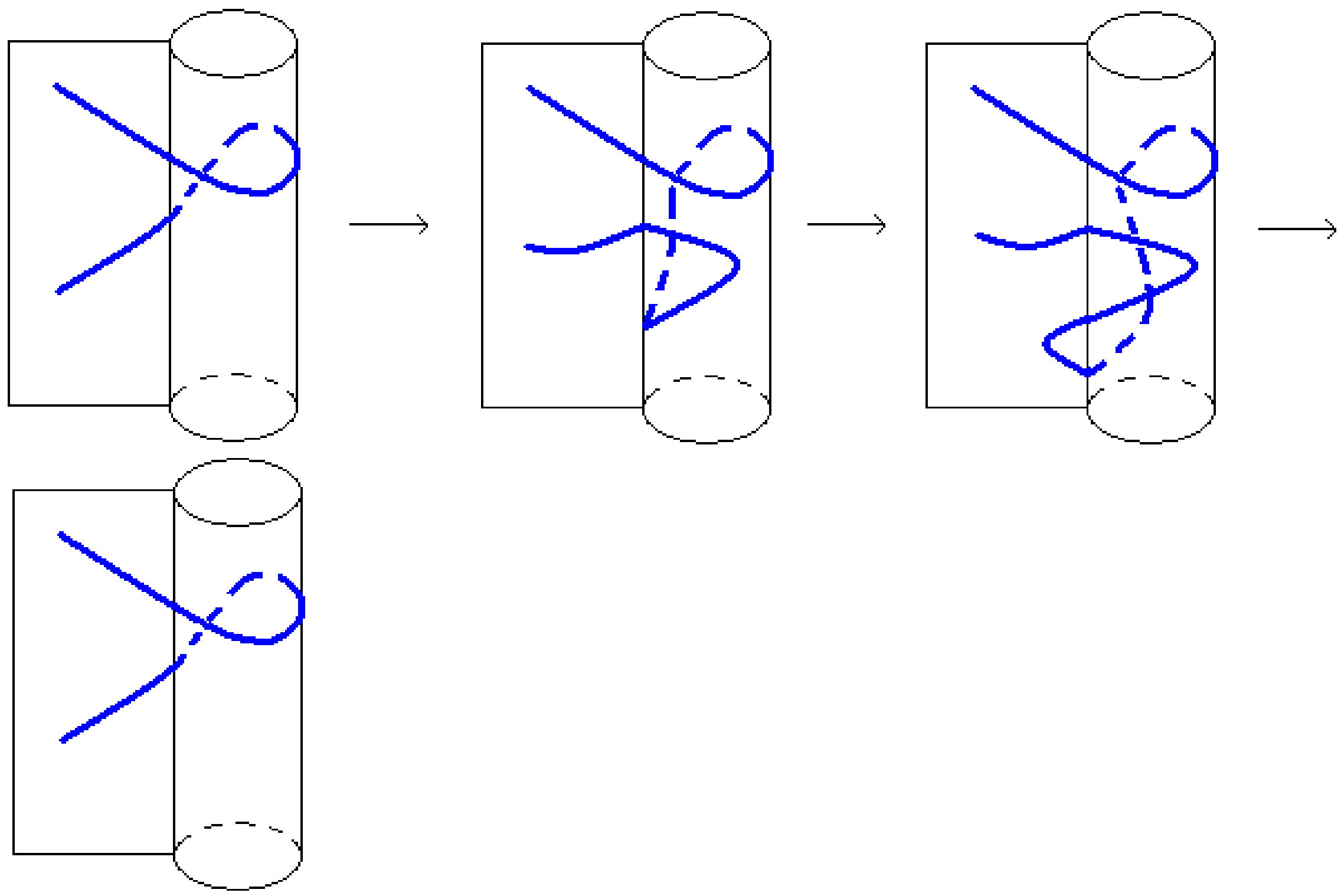

**Figure 159- Twists on the generator cylinder- Case g)**

### 6.2 *Realize twists at the 2-cell*

We consider the case, where two relations are multiplicated. This can be realized by sliding  a little arc of a relation on the other. In chapter 3 we have motivated the appearance of twists. In the previous chapter we have performed the twists at the generator cylinder:

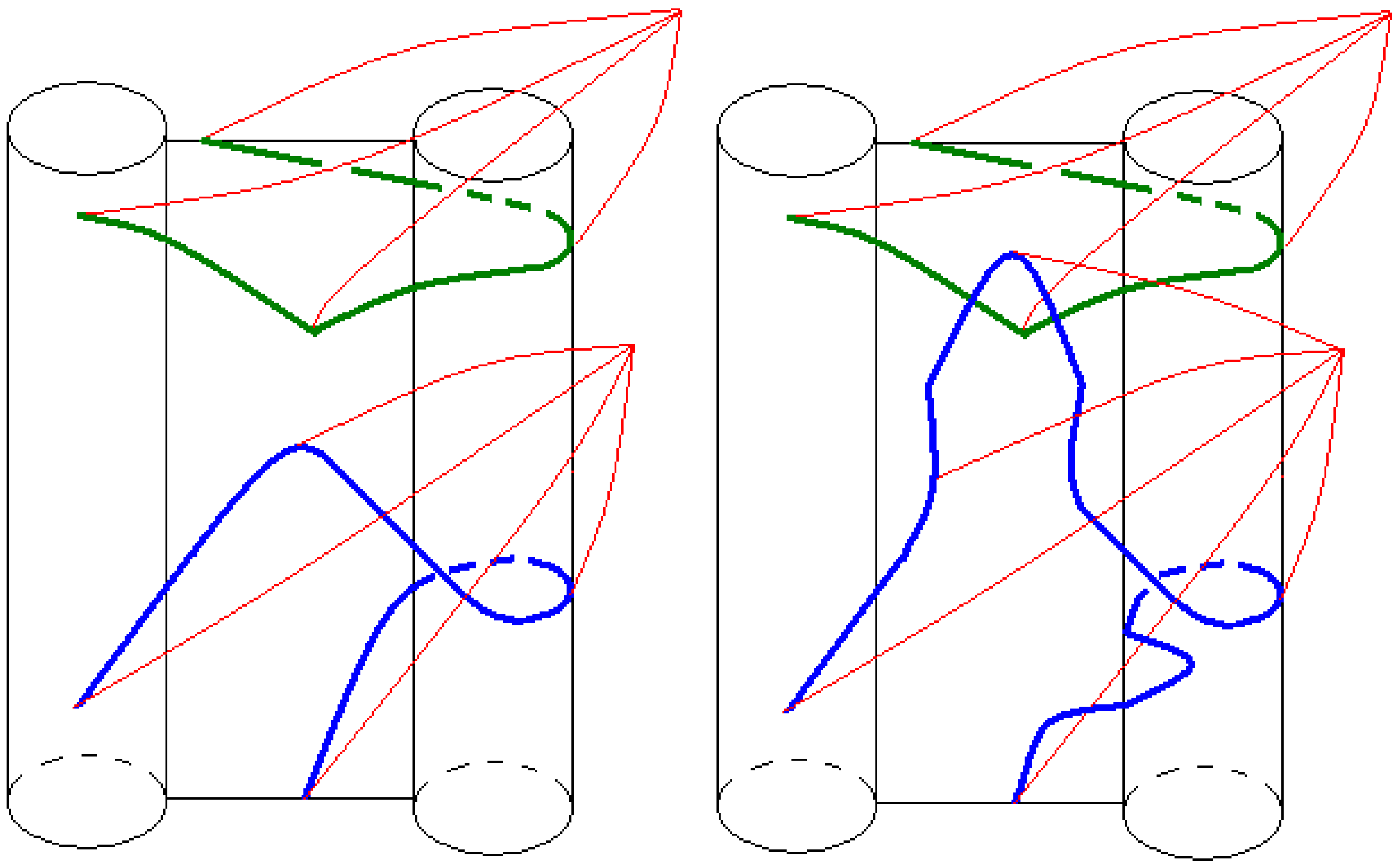

**Figure 160 - 5.2 Realize twists at the 2-cell – twist on generator cylinder**

In the next chapter on the Q-transformations, we will apply our results; but here we work out an alternative, where the twists are realized at the 2-cell. It has the disadvantage, that for each Matveev move in a Q-transformation we have to drag along the twists. Clearly that makes life difficult, nevertheless we want to study it.

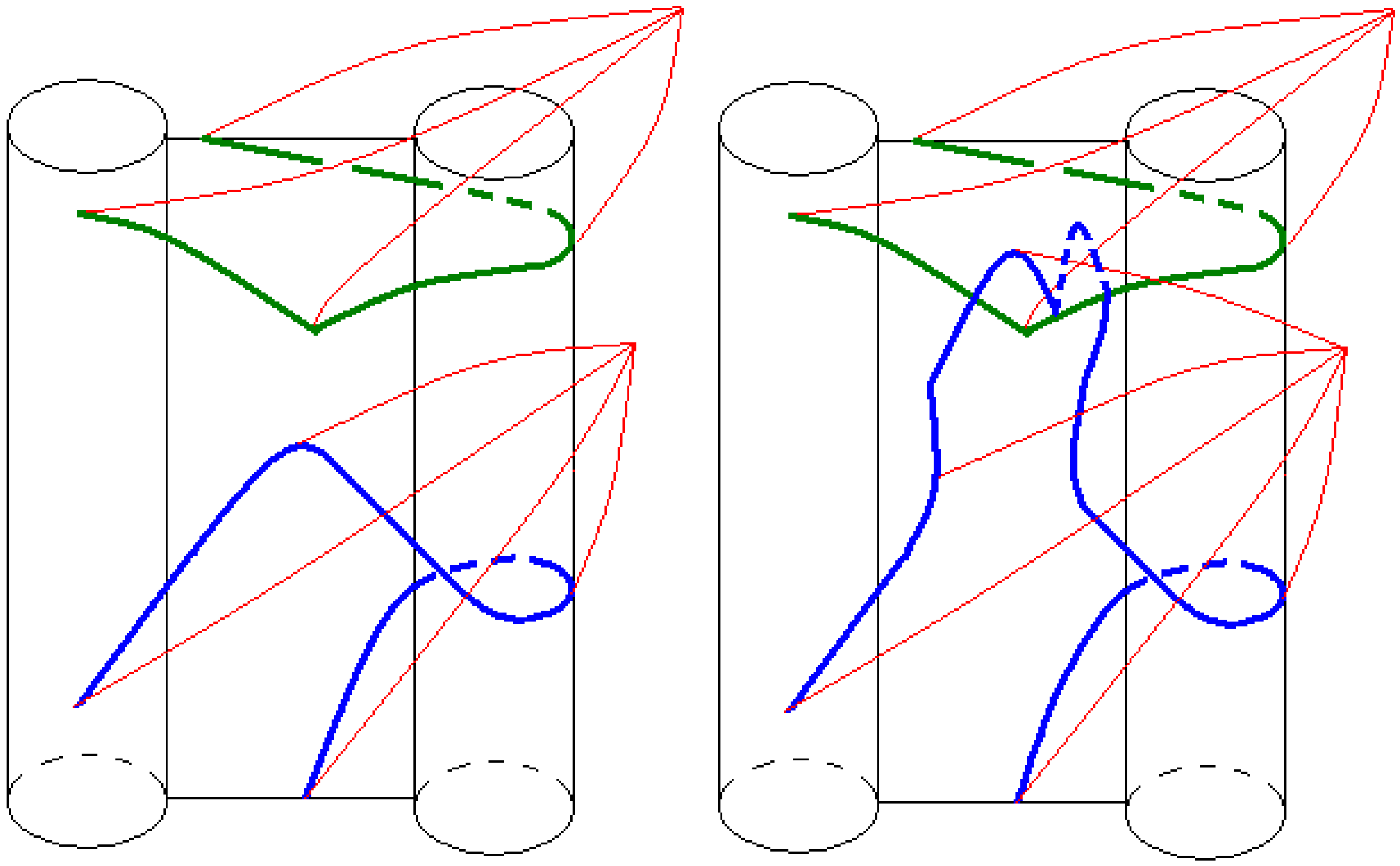

**Figure 161- Realize twists at the 2-cell- twist at 2-cell**

Some comments to the picture:
The green thick line represents the attaching curve of the relation, which is fixed under the multiplication, the thick blue line indicates the attaching curve of the other relation, which slided on the 2-cell that belongs to the green curve. Note, that the broken blue line of the twist is in the rectangle between the generator cylinders. We will develop the sequence of slices and support these by different colours for the attached 2-cells.

In the sequence of slices before we realize the twist, the situation looks like in the left picture, but to perform the twist, we have to change it by introducing a pair of saddlepoints as shown in the right picture:

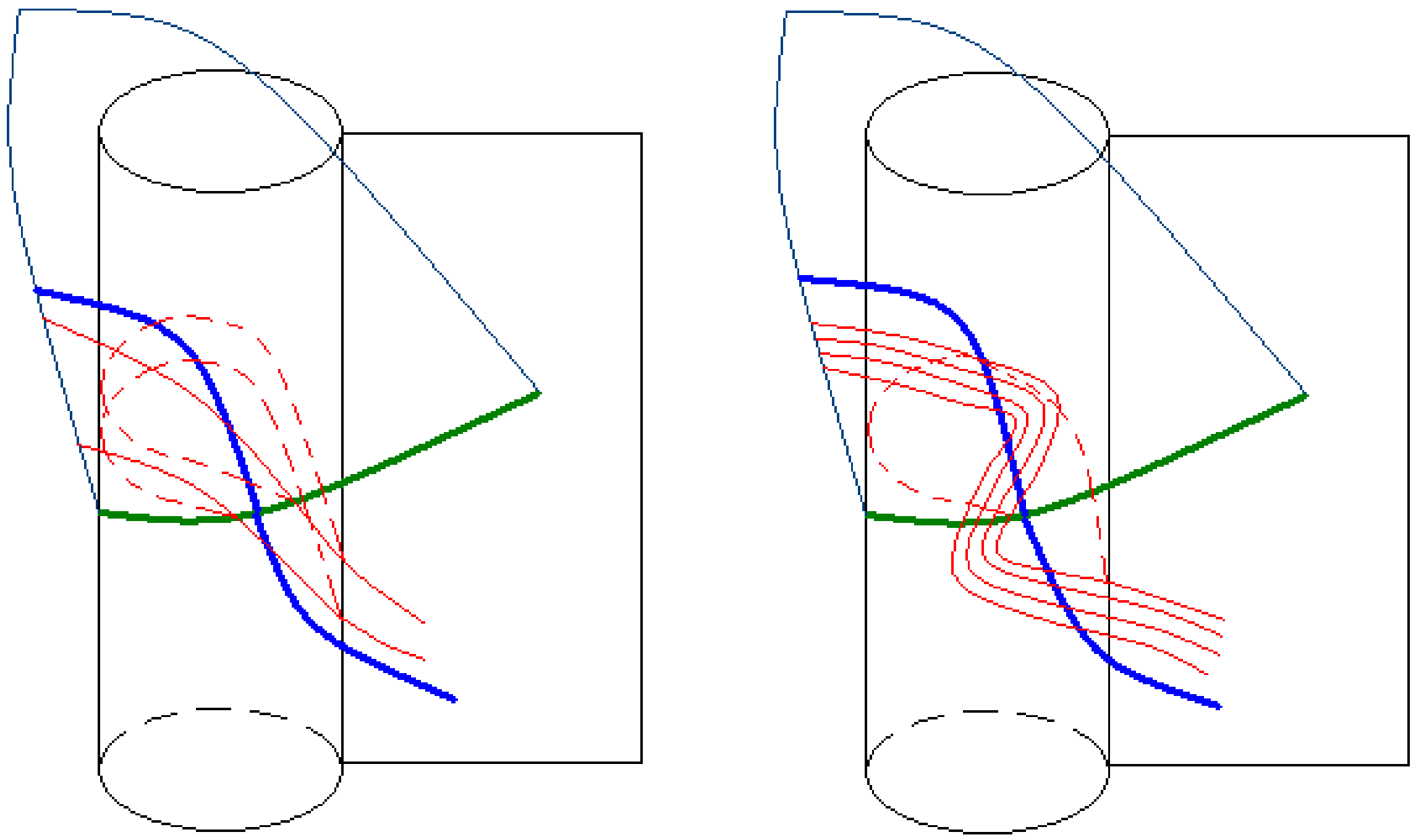

**Figure 162- Realize twists at the 2-cell– prepare twist**

1) Pass from the generator cylinder to the rectangle

We give a preview for that sliding:

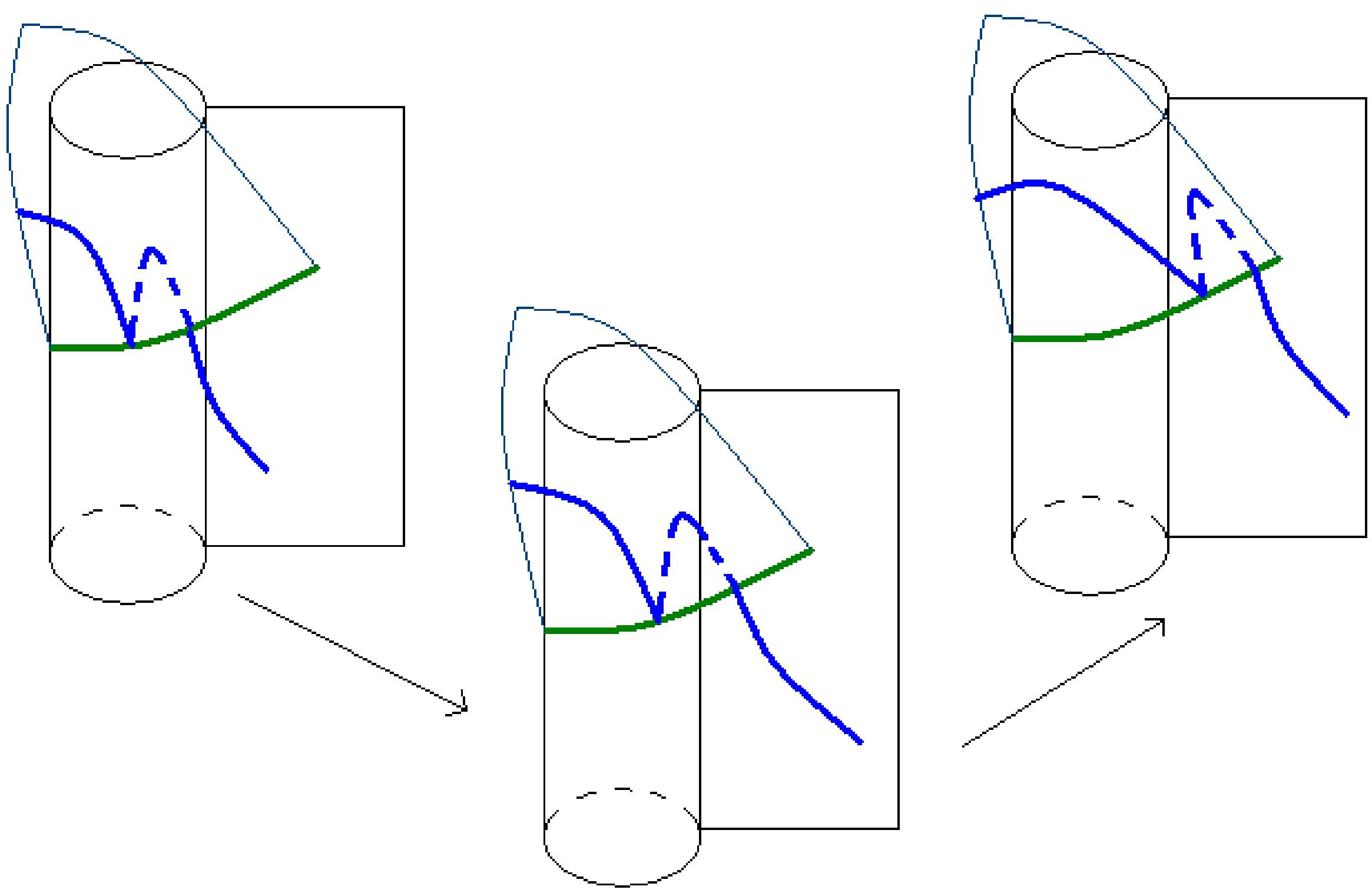

**Figure 163- Realize twists at the 2-cell- preview- slide twist from generator cylinder to rectangle**

We draw the slices of the twisted curve in the Quinn model and assign it to the local model, the colours support the correspondence:

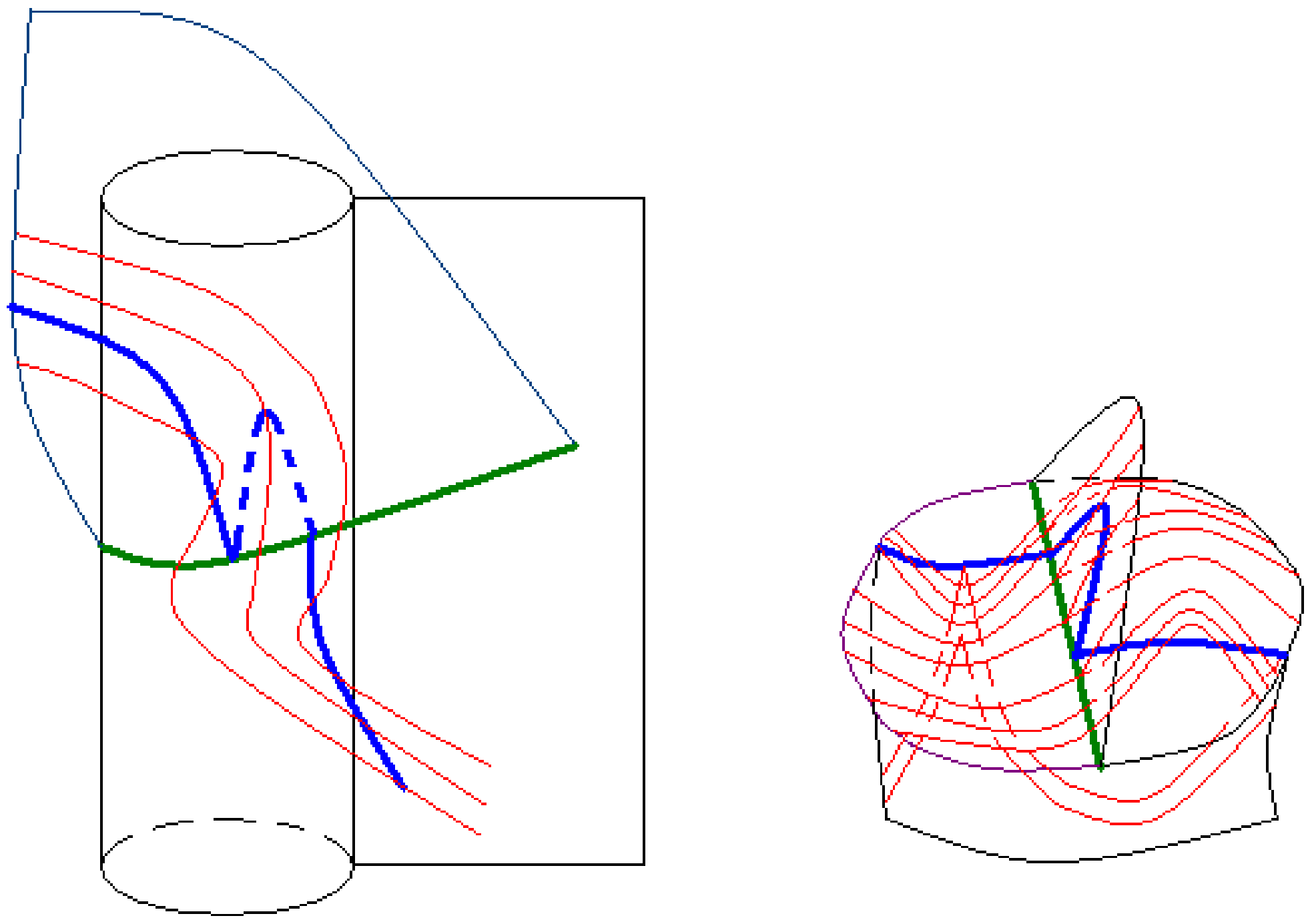

**Figure 164- Realize twists at the 2-cell- $T_3$ turn on generator cylinder**

The first step in the preview leads to the local move:

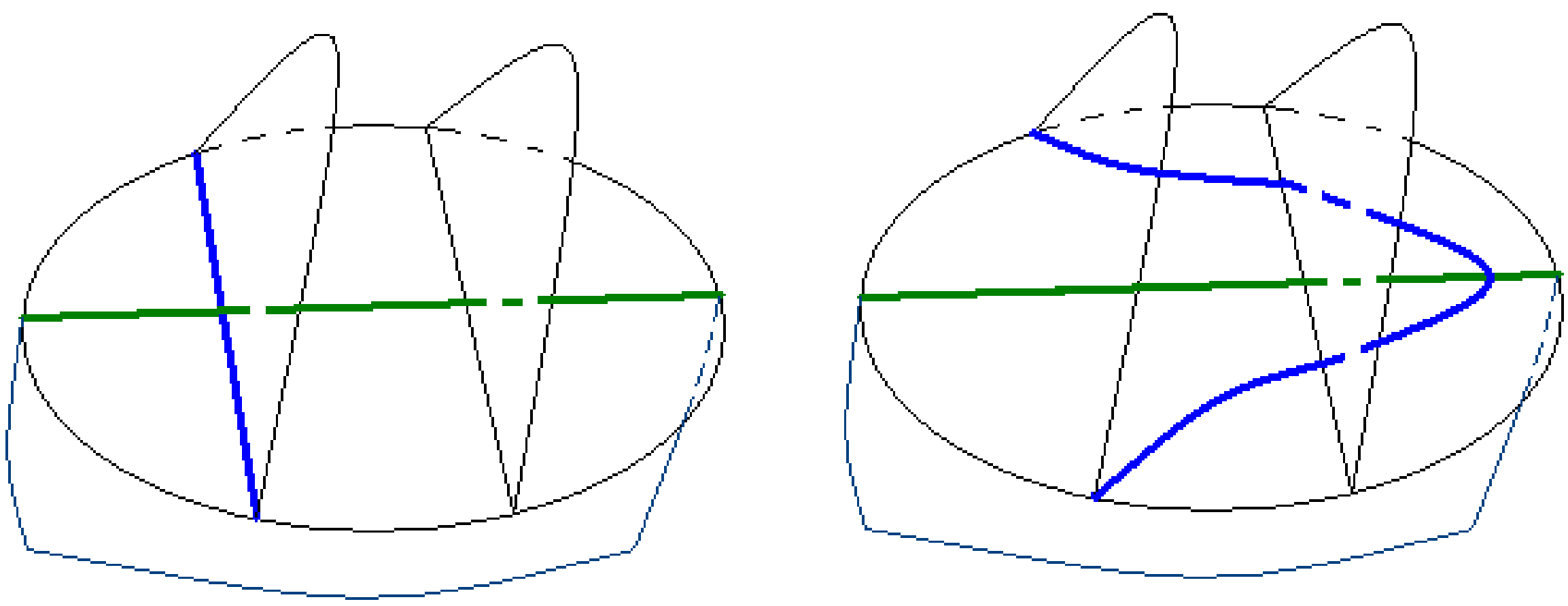

**Figure 165- Realize twists at the 2-cell- first step in slide- modified $T_2$ move**

This move appears in chapter 5, where it is called a “modified $T_2$ move”. We decompose it into a sequence of Matveev moves:

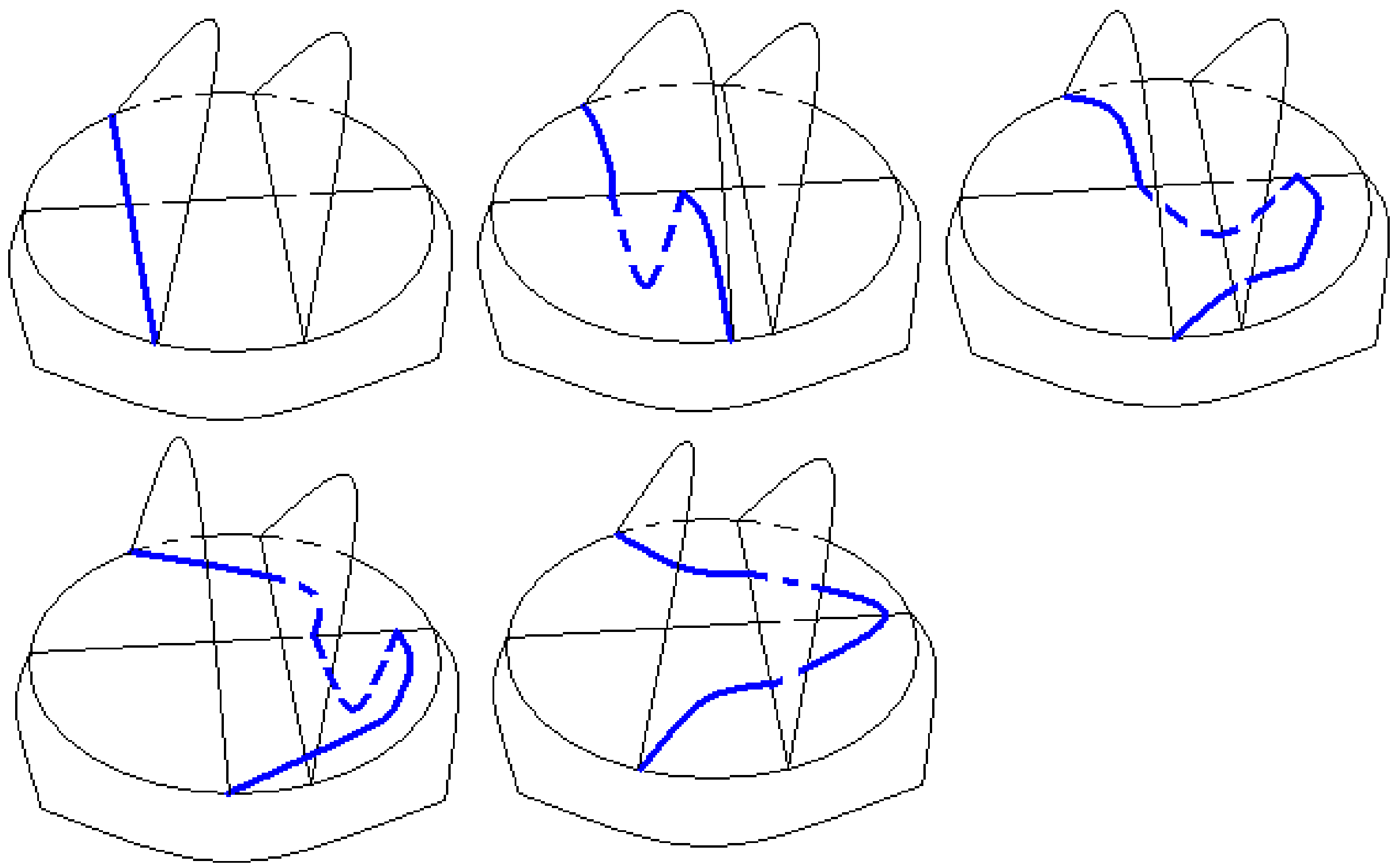

**Figure 166- Realize twists at the 2-cell- preview- modified $T_2$ move as composition of Matveev moves**

To develop the sequence of slices for each step, it is sufficient only to concentrate on the main slices to see what’s going on, otherwise the figures would lose their essential message. We locally change the slices to prepare “good $T_3$ turns” and deform it so slightly, that we do not change the “rhythm” of the sequence of slices.

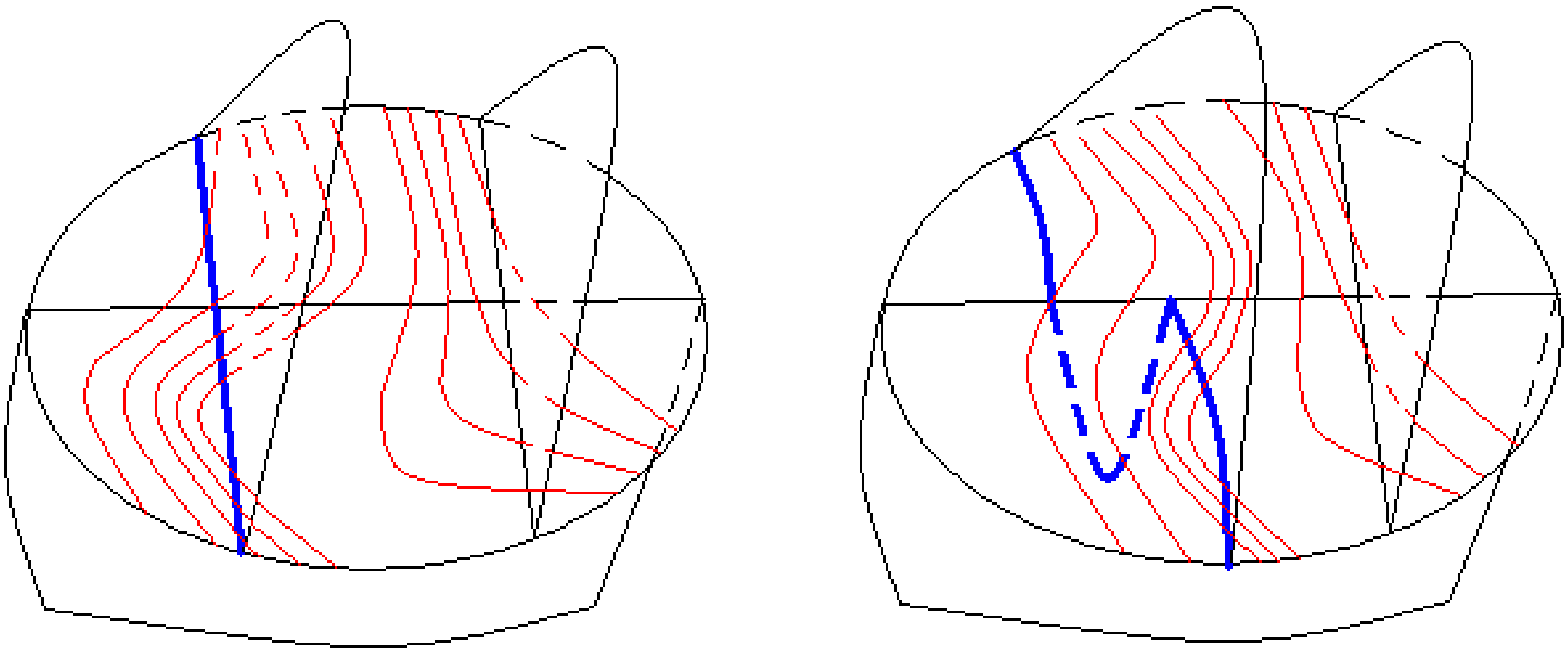

**Figure 167- Realize twists at the 2-cell- modified $T_2$ move reduced- 1**

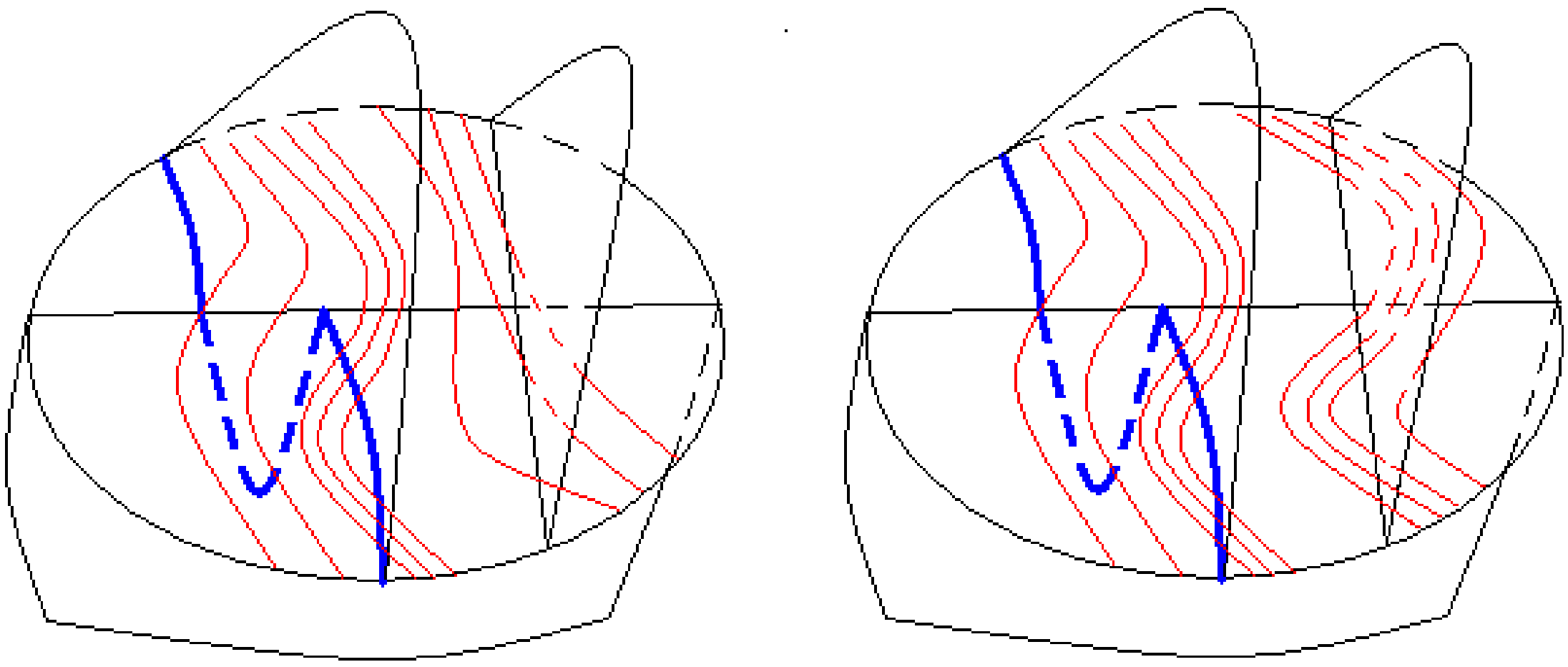

**Figure 168- Realize twists at the 2-cell- modified $T_2$ move reduced- 2**

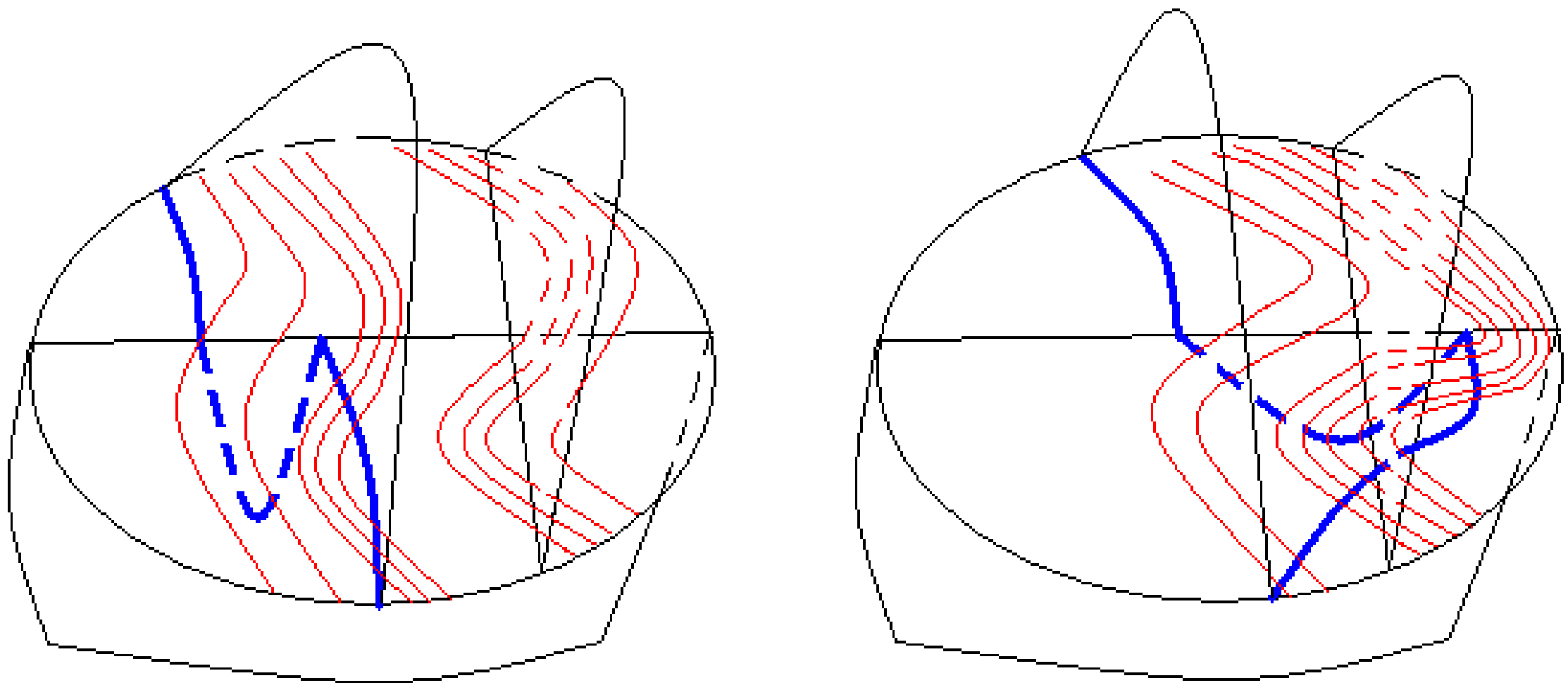

**Figure 169- Realize twists at the 2-cell- modified $T_2$ move reduced- 3**

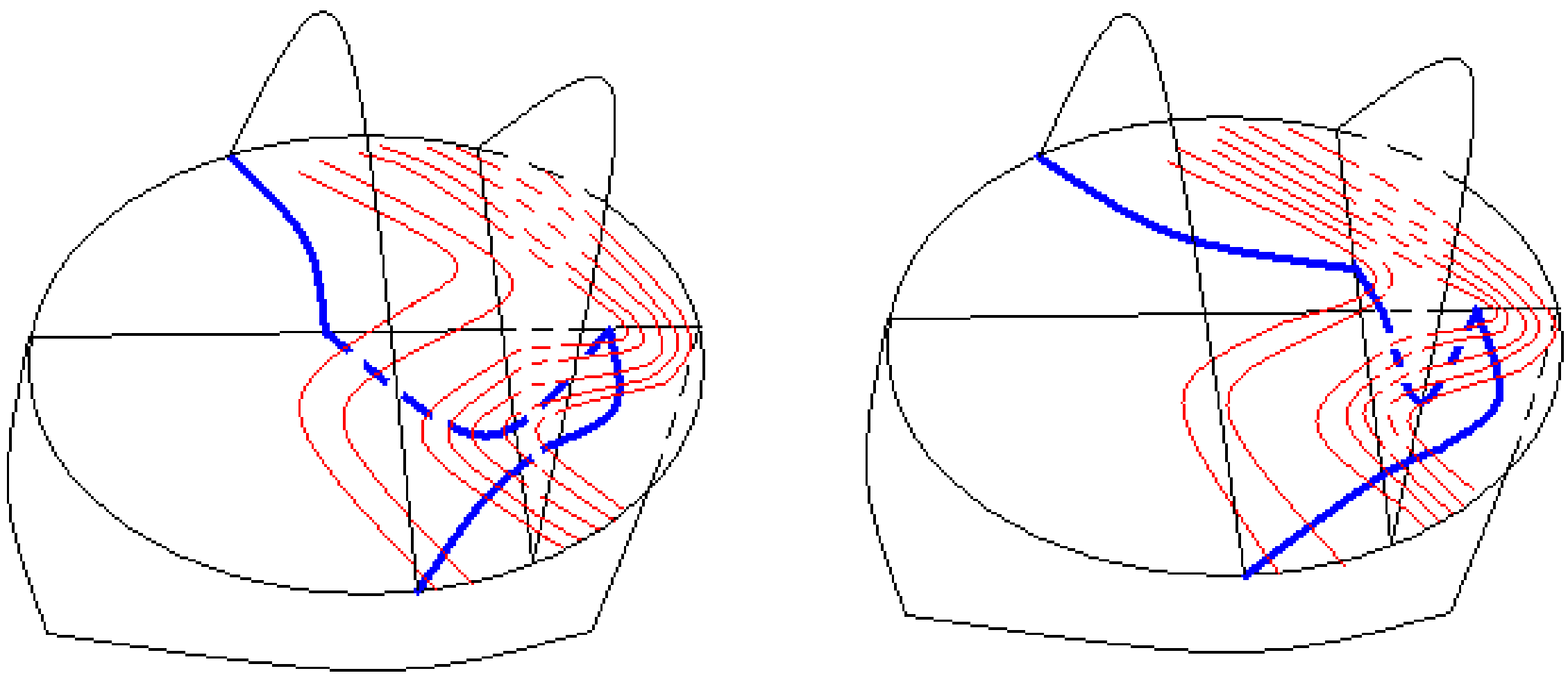

**Figure 170- Realize twists at the 2-cell- modified $T_2$ move reduced- 4**

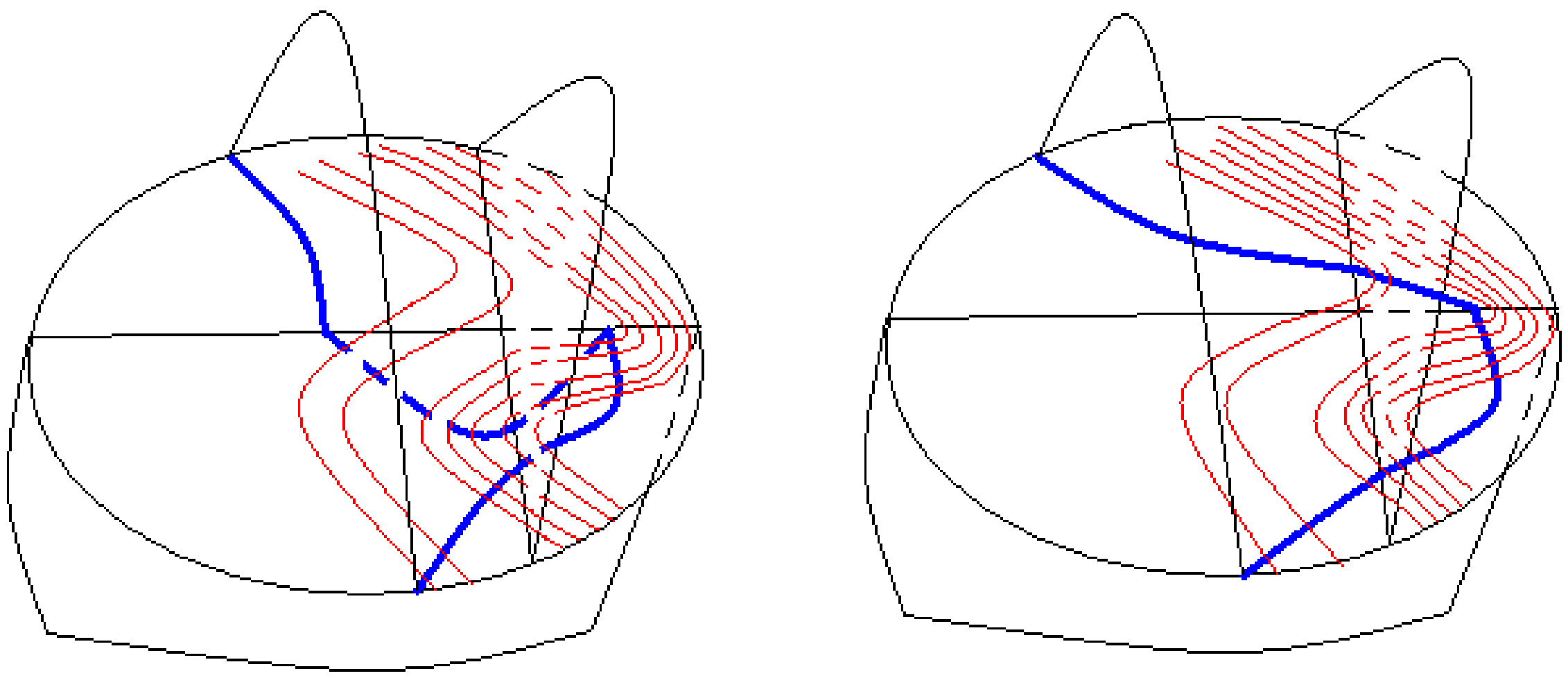

**Figure 171- Realize twists at the 2-cell- modified $T_2$ move reduced- 5 (end)**

Now we have finished the move and translate it back to the Quinn model:

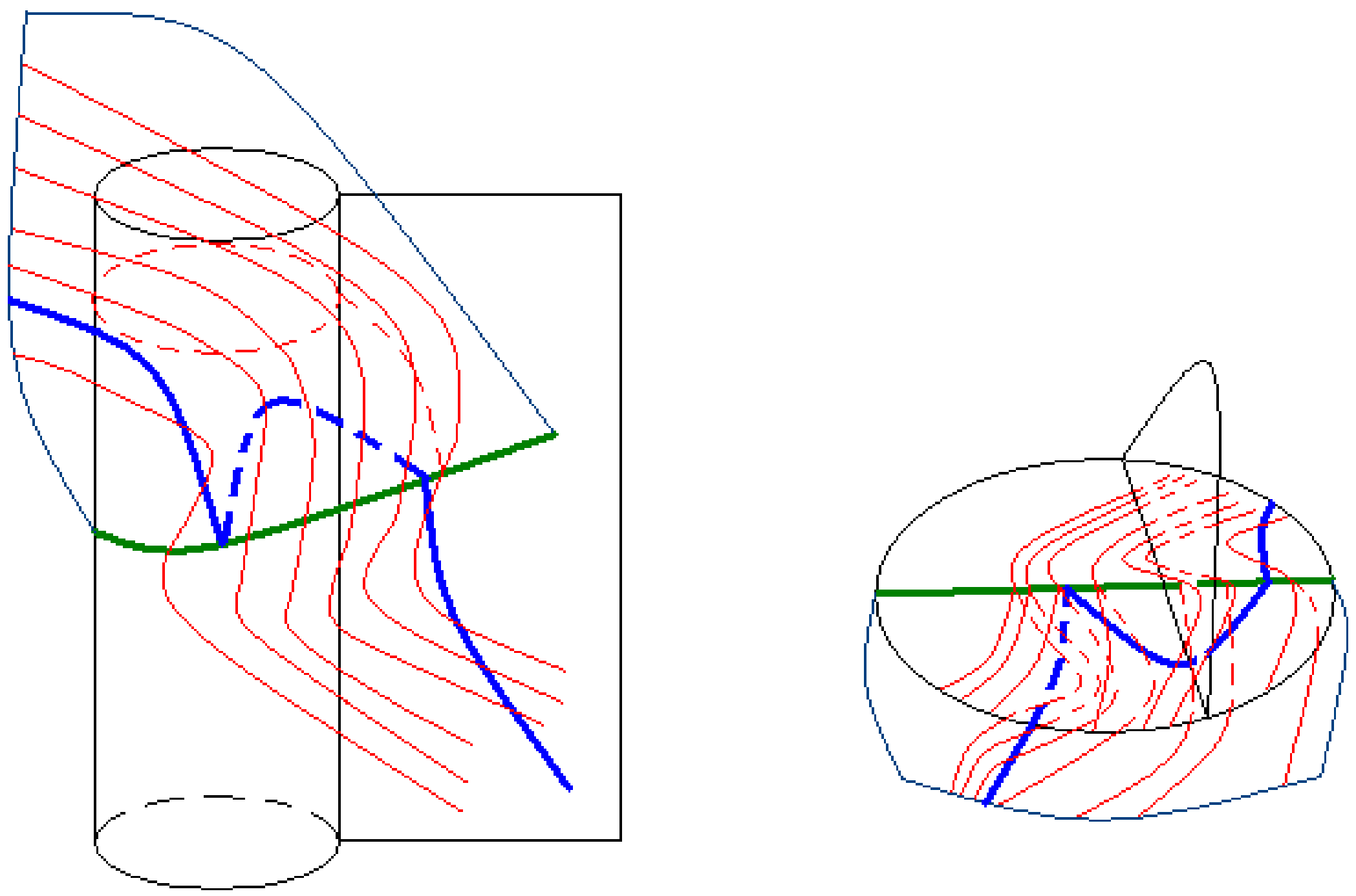

**Figure 172- Realize twists at the 2-cell- $T_3$ turn on (generator cylinder and rectangle)**

We perform the second (and last) step in our preview. Note that we have to draw more slices on the bottom component to clearify that we have “good $T_3$ turns”:

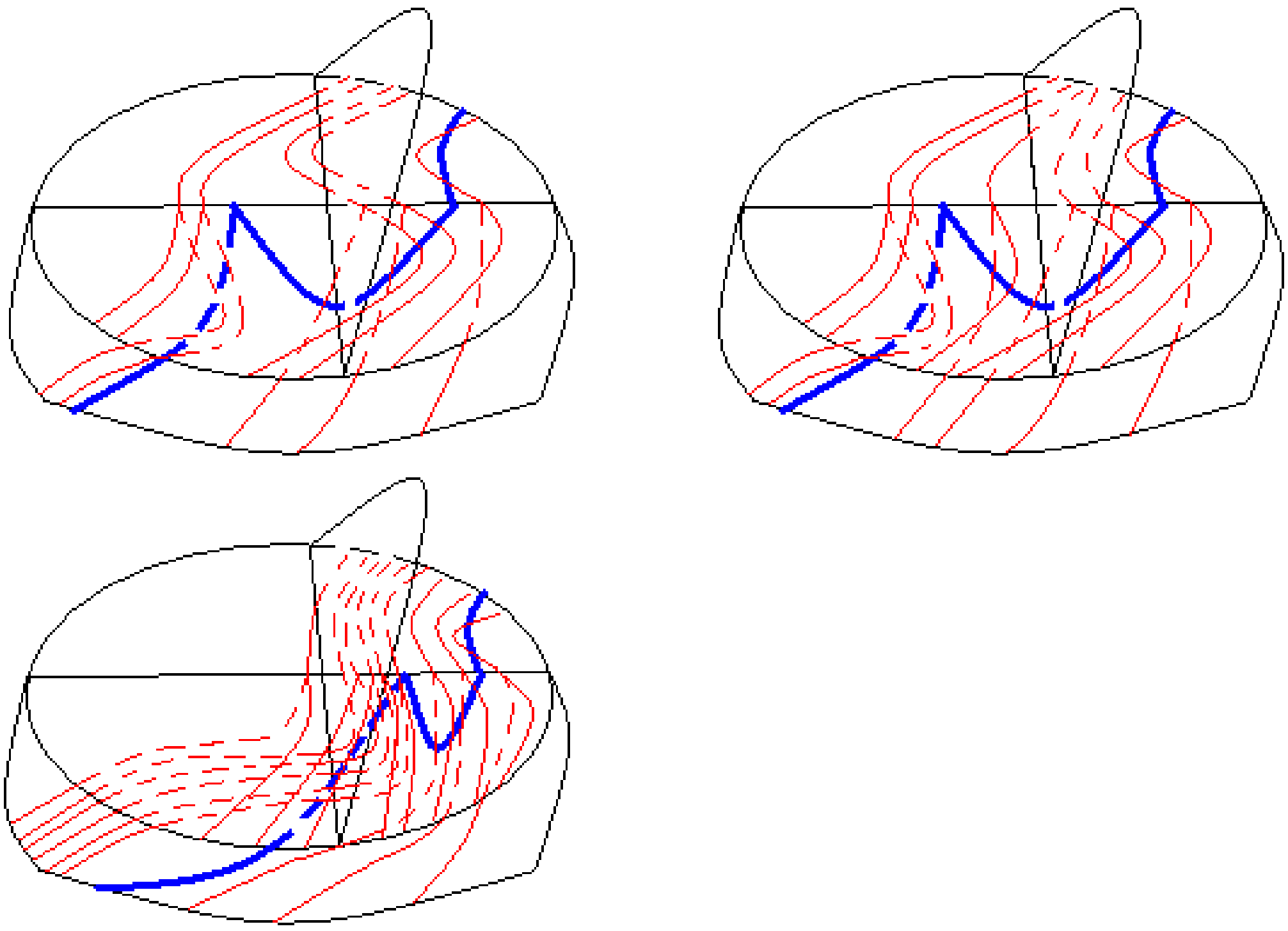

**Figure 173- Realize twists at the 2-cell- slide to rectangle**

We transfer the result to the Quinn model:

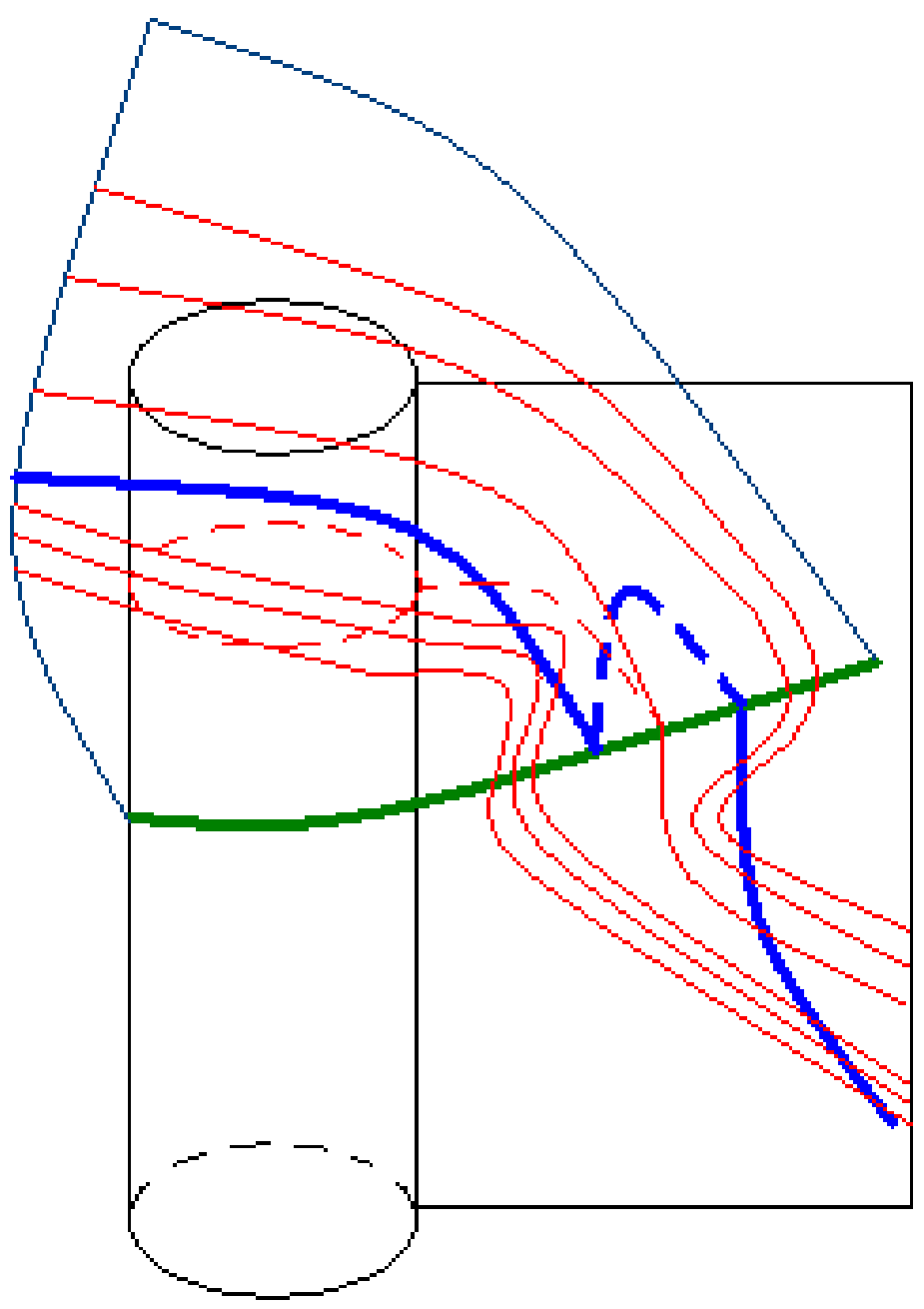

**Figure 174- Realize twists at the 2-cell- $T_3$ turn on rectangle**

We see, that we have not made essential changes to the sequence of slices.

We only list the other cases for sliding a twist, but compared to the previous one there is no difference in the moves and the slices:

2) Pass from rectangle to generator cylinder:

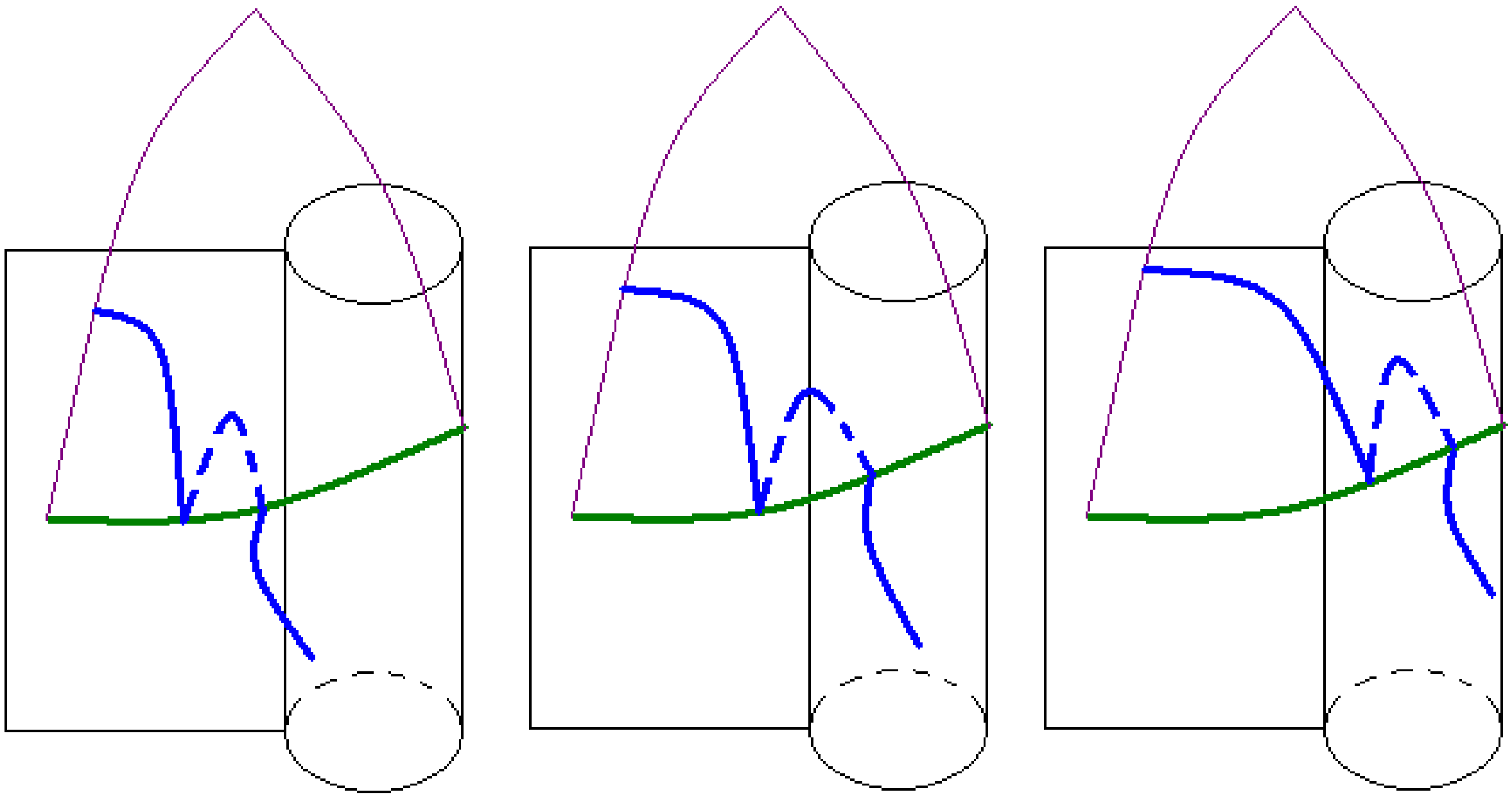

**Figure 175- Realize twists at the 2-cell- slide twist from rectangle to generator cylinder**

3) Pass a selfintersection of the attaching curve:

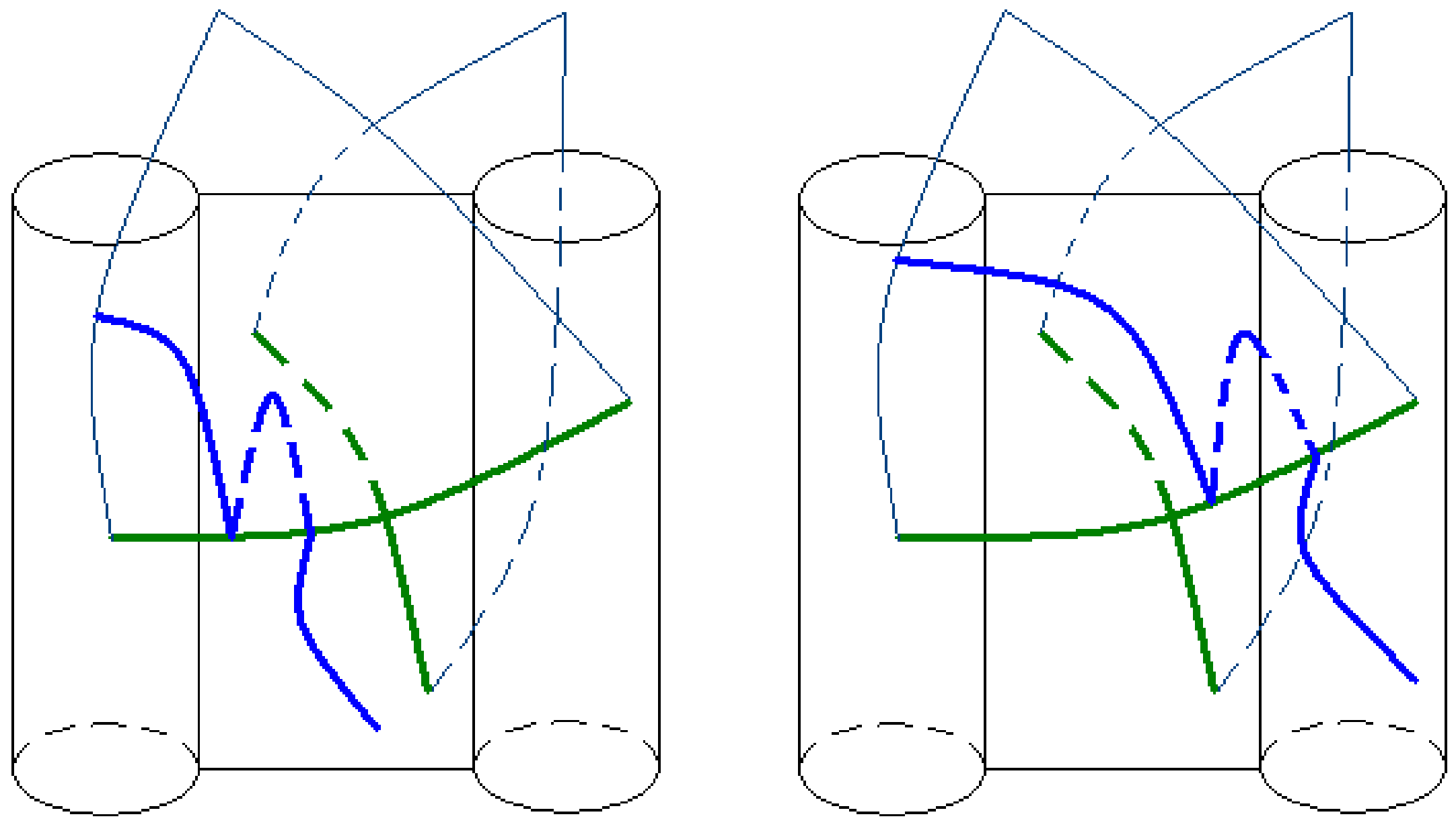

**Figure 176- Realize twists at the 2-cell- twist pass a selfintersection**

# 7 Q-transformations and 2-deformations in the Quinn model

We consider the Q-transformations in the Quinn model and translate them into a homotopy of attaching maps and then decompose it into a sequence of Matveev moves. The idea to get the Matveev moves is, to perform single steps of the homotopy by crossing vertices of that model with arcs of the attaching curve. We have to analyse all these crossings.
We develop the sequence of slices and compare the starting sequence with the resulting one. We use the results of the former chapters as follows:
We realize twists, compositions of twists or loops at the generator cylinder far away from all Matveev moves and simplify or annihilate them when possible.
Therefore we can exclude them in our considerations about Q-transformations !!!
We list the Q-transformation (see chapter 3), restricted to a 2-complex with two generators and relations:
Let P = ⟨a,b / R,S ⟩ be a presentation of the 2-complex, represent in the Quinn model
Then there are three types of transformations for the relations:

1) multiplication:
   R → SR
   S → S

2) conjugation:
   R→ $wRw^{-1}$ w is a word in the generators a, b

3) inverse:
   R → $R^{-1}$

The other transformation about generalized prolongation belongs to the 2-deformation which is easier to handle.

## 7.1 The multiplication

The multiplication corresponds to the slide of the 2-cell R on the 2-cell S, so the attaching curve of R changes to that of SR

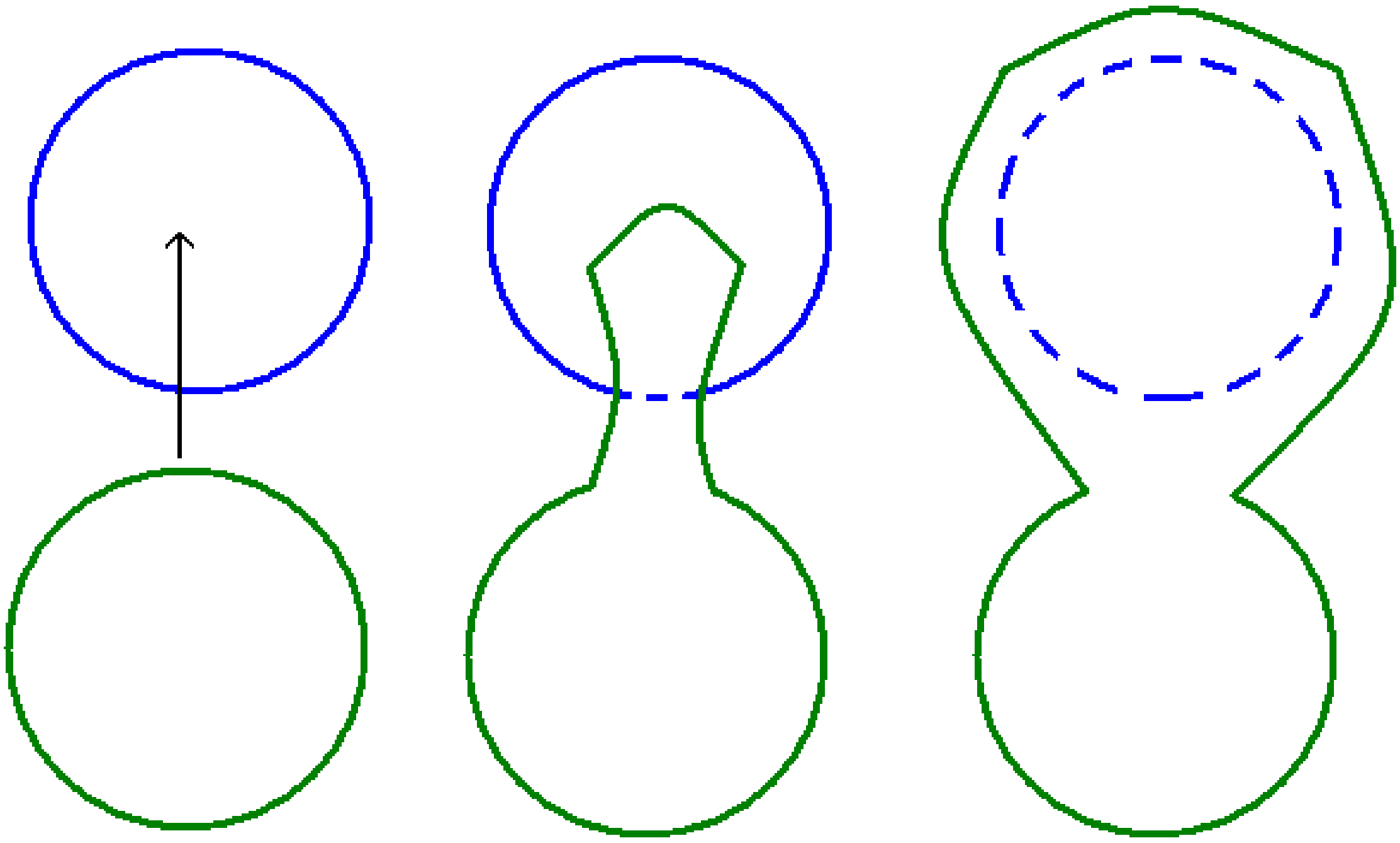

**Figure 177- The multiplication- slide of the attaching curve**

We see, that it is sufficient to slide a small arc of the attaching curve of R onto S. We transfer the three phases into the Quinn model:

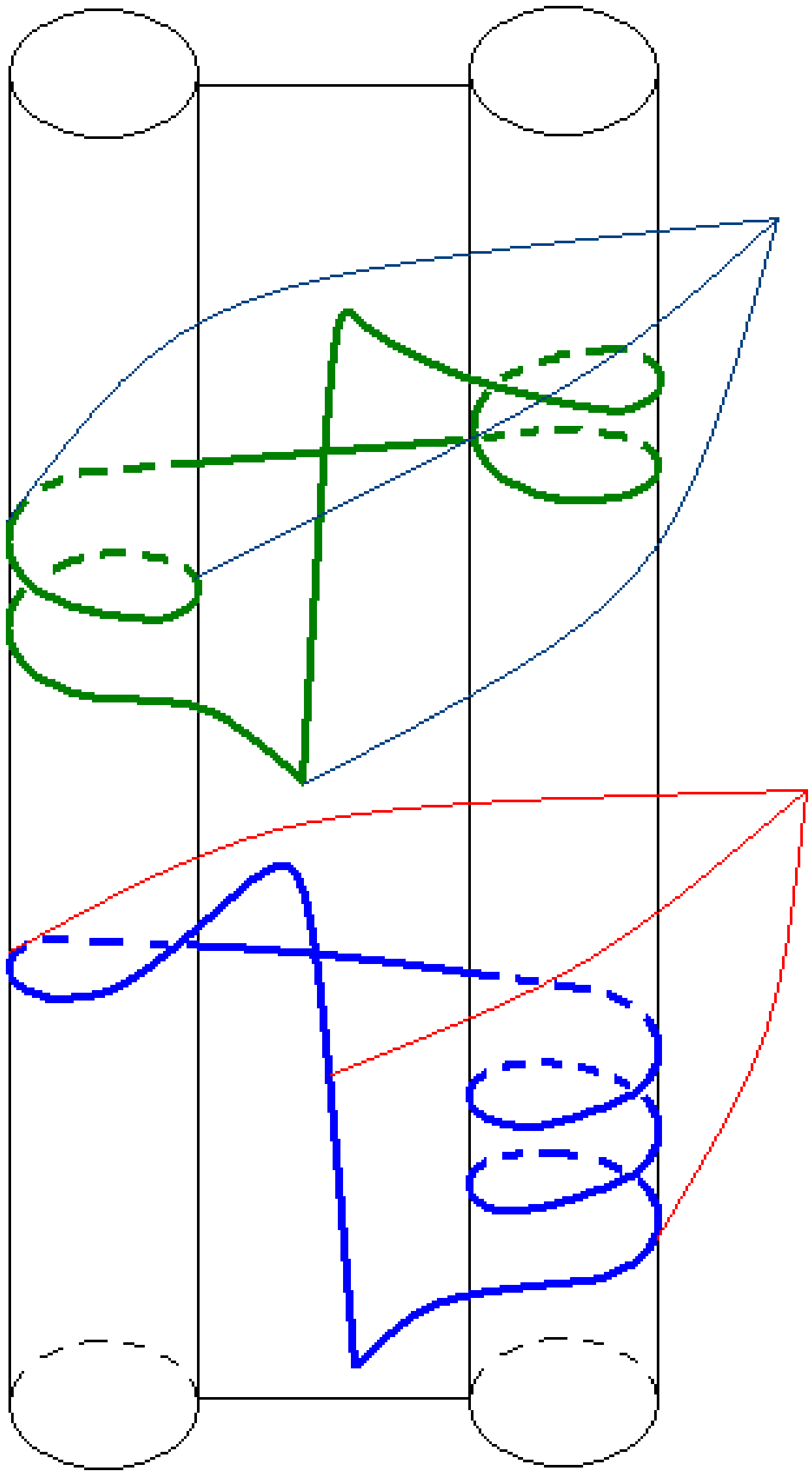

**Figure 178- The multiplication- disjoint attached 2-cells in Quinn model**

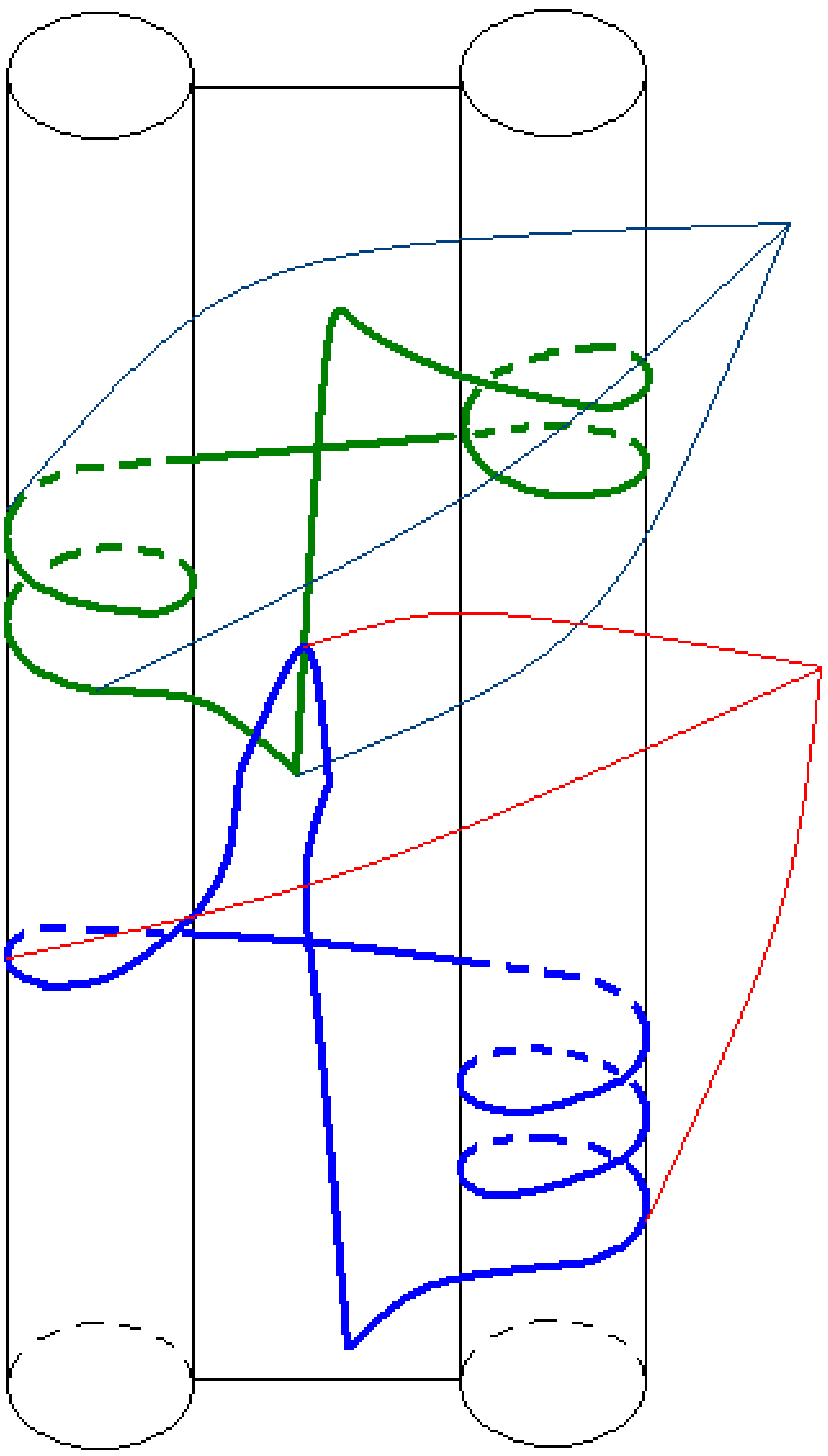

**Figure 179- The multiplication- entry of the slide in Quinn model**

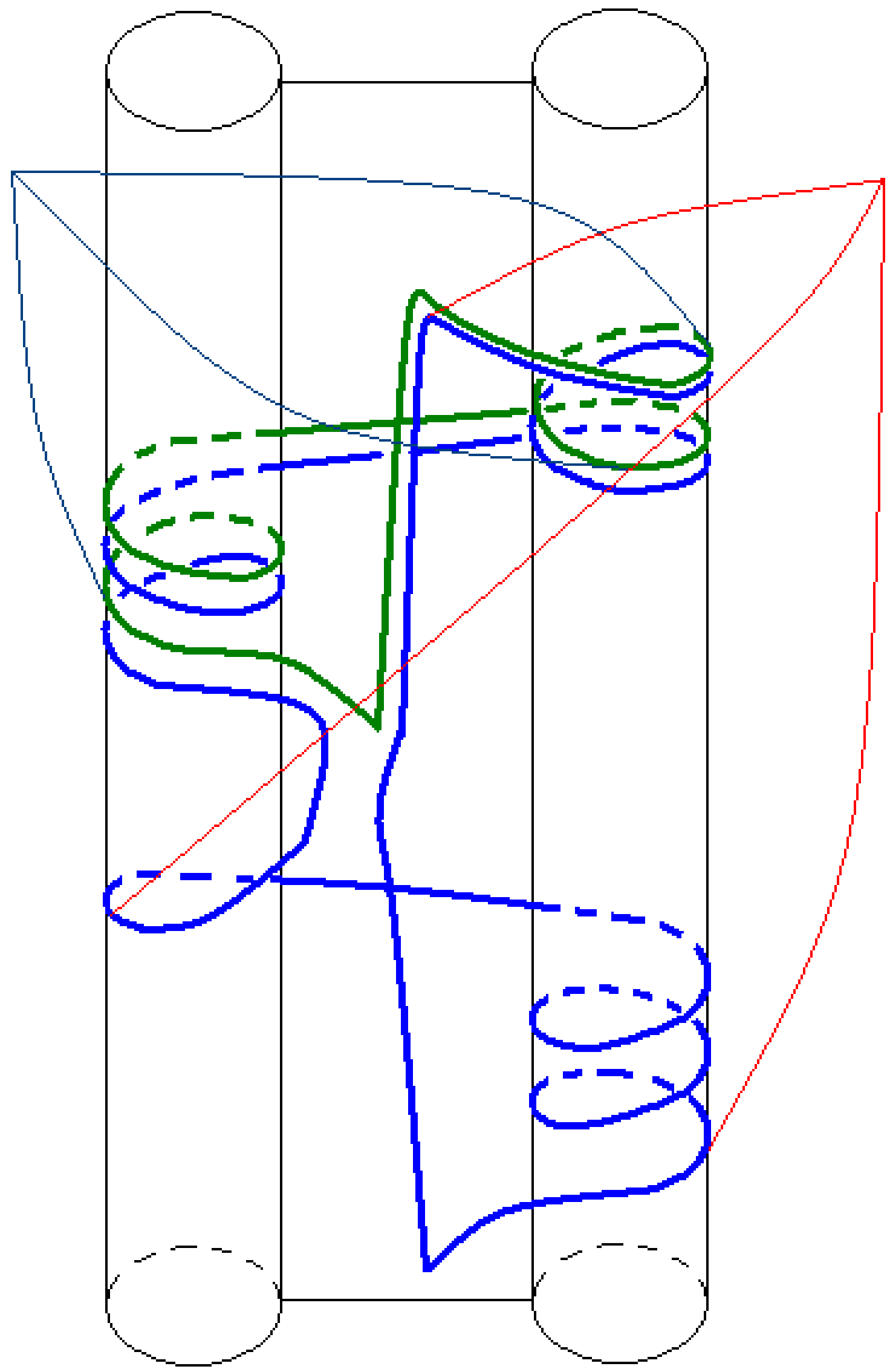

**Figure 180- The multiplication- almost parallel attaching curves of both 2-cells in Quinn model**

Note that now the slided blue arc is almost parallel to the green attaching curve. We have to drop down this part to the unchanged rest of the blue curve. Then and not before that, the blue attaching curve corresponds in the Quinn model to the multiplication of two relations.

We explain the process by the given example, where all crossing situations appear. Clearly that leads to a decomposition in Matveev moves. We drop down the blue line and our strategy will be to keep the level order which means:
If h is a height function, x drops to x', y drops to y'
$h(x) < h(y) \Rightarrow h(x') < h(y')$
When there is no crossing with the green curve left, we can drop the blue line down along the generator cylinders and rectangle. Hence the process is finished. We omit the 2-cells in our figures. The start figure is:

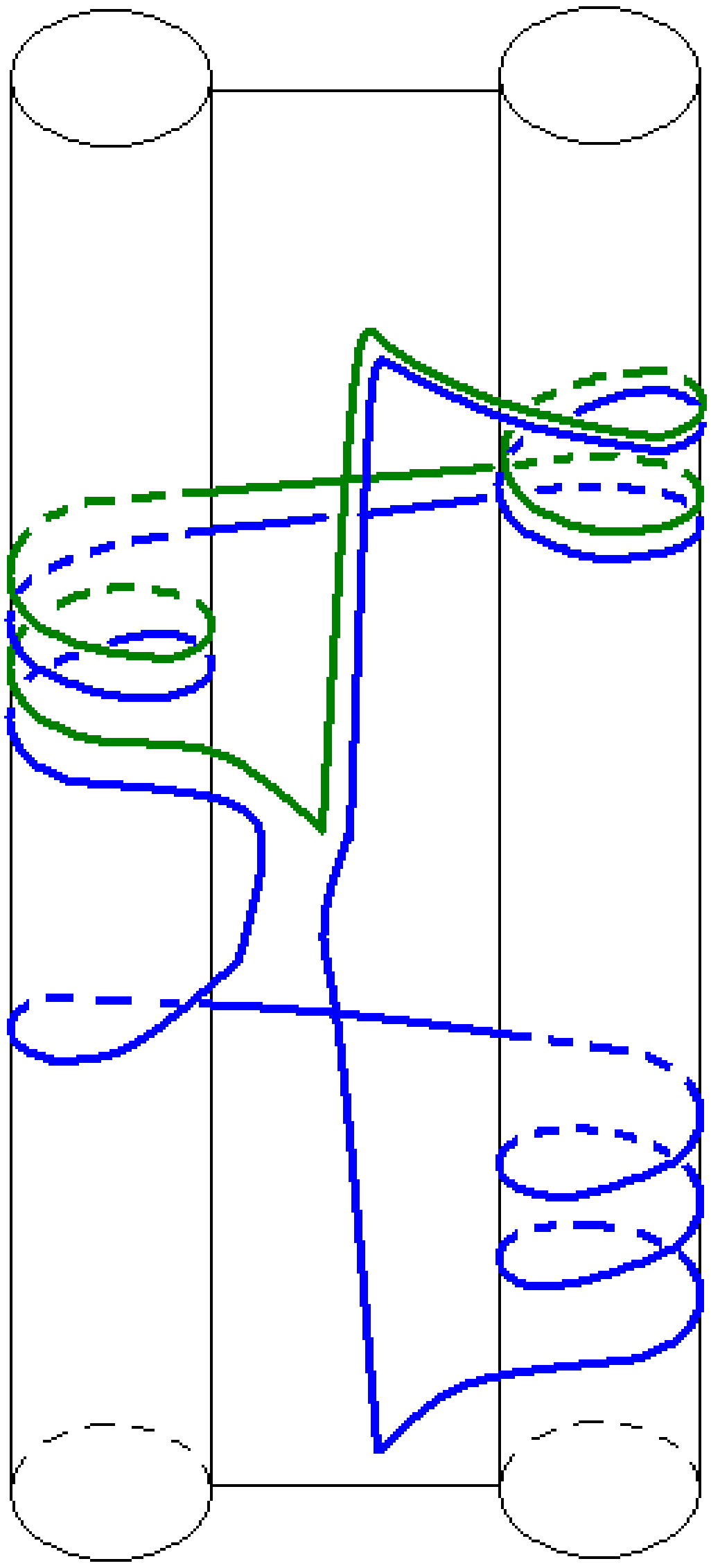

**Figure 181- The multiplication- drop down slided attaching curve- 1**

We start to drop down on the left generator cylinder. We cross the green line twice and prepare the move $T_2^{-1}$ by using the move T*:

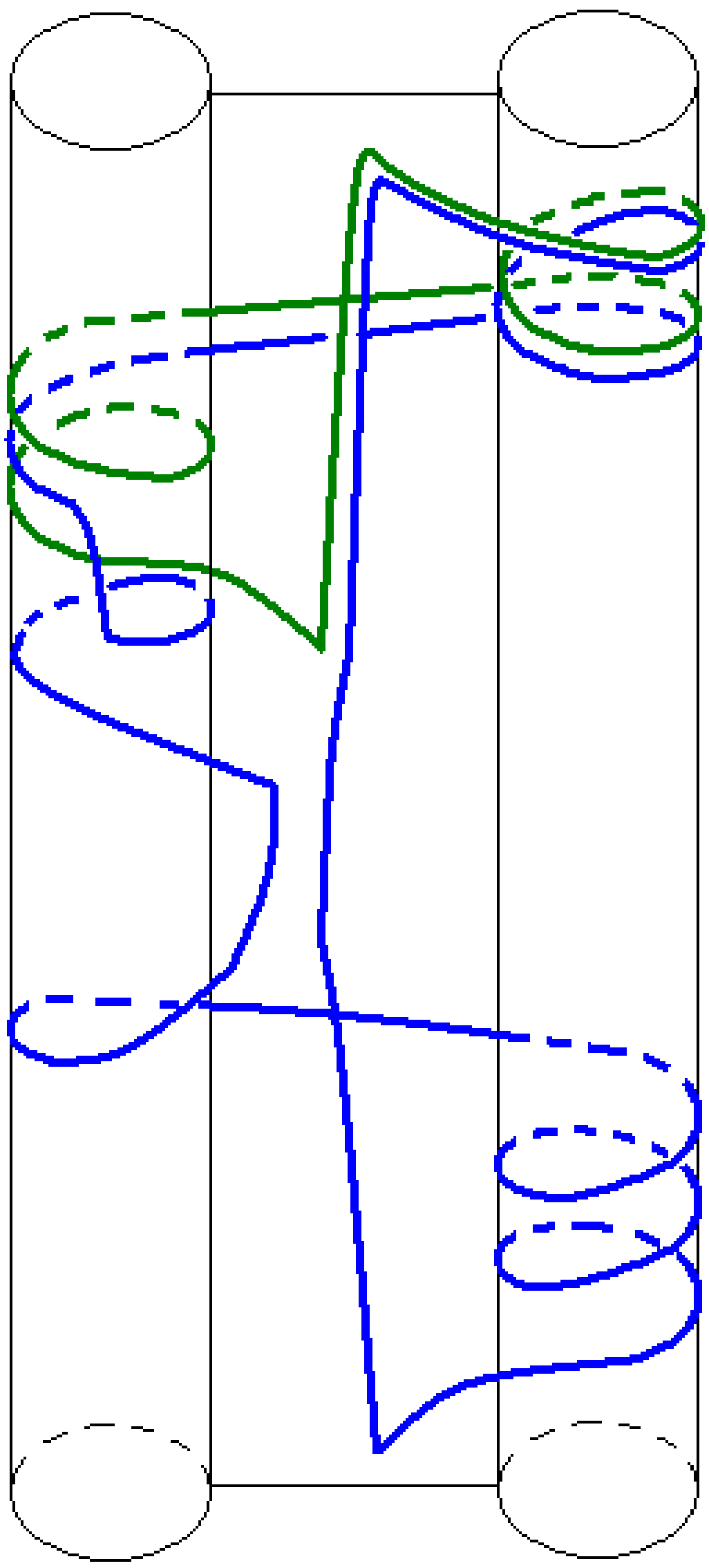

**Figure 182- The multiplication- drop down slided attaching curve- 2**

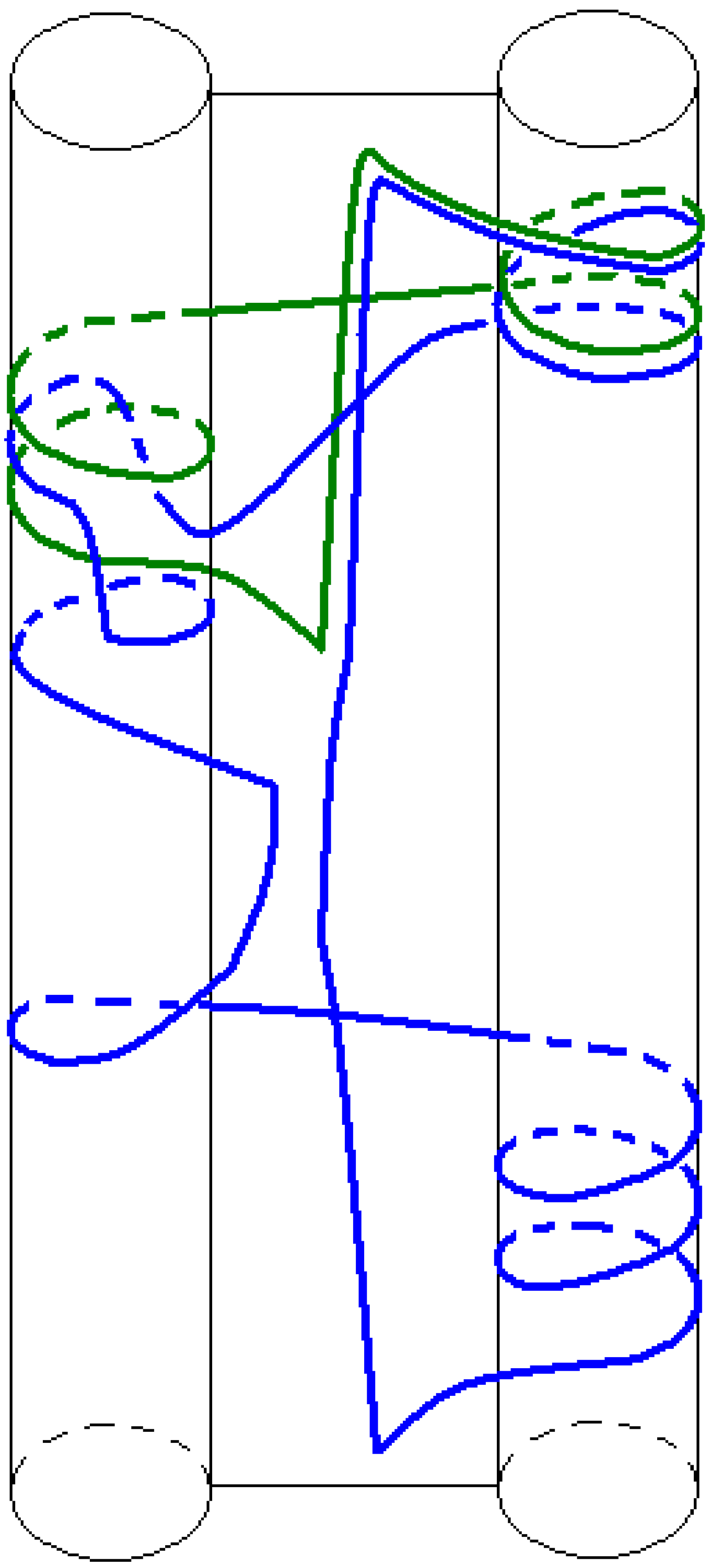

**Figure 183- The multiplication- drop down slided attaching curve- 3**

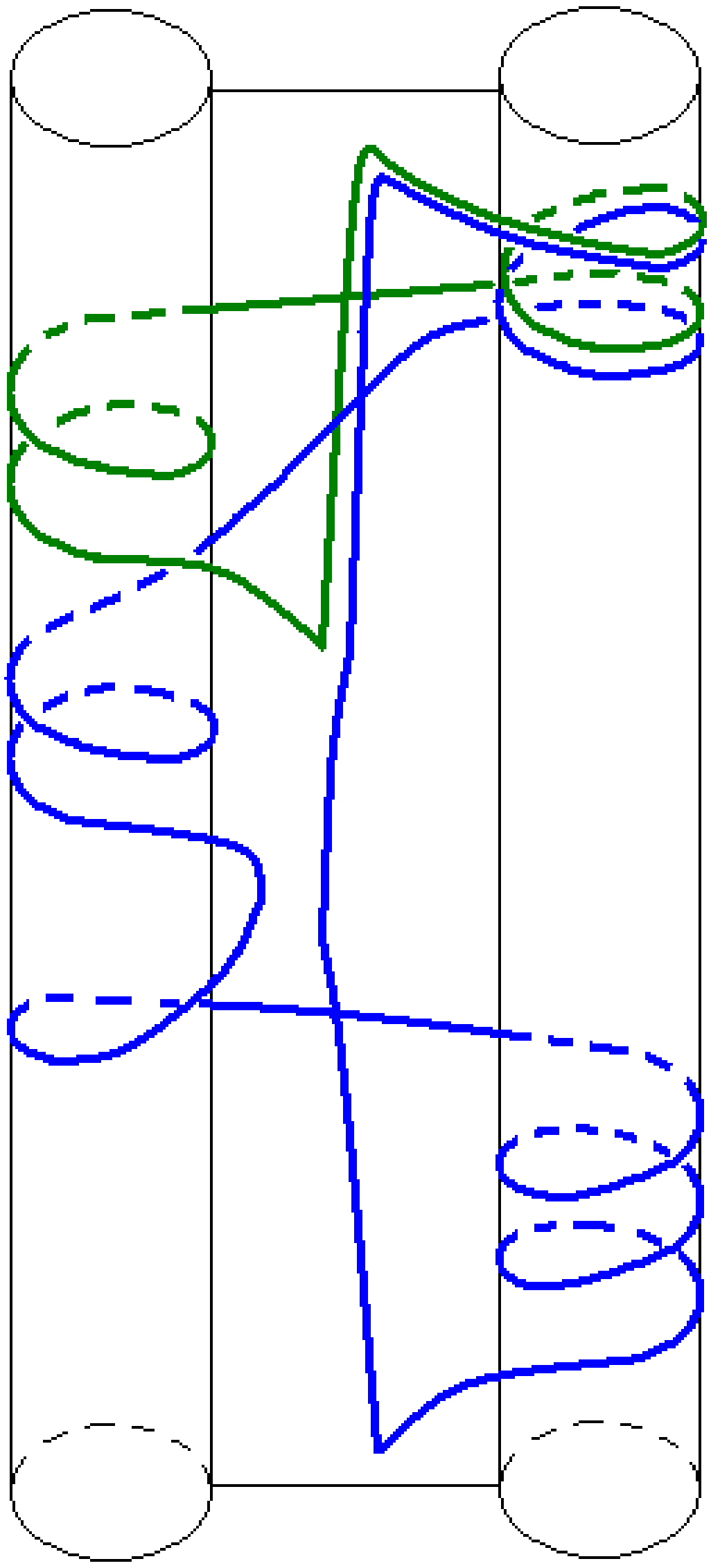

**Figure 184- The multiplication- drop down slided attaching curve- 4**

Now we can perform the move and hence we get the blue curve free on the left generator cylinder:

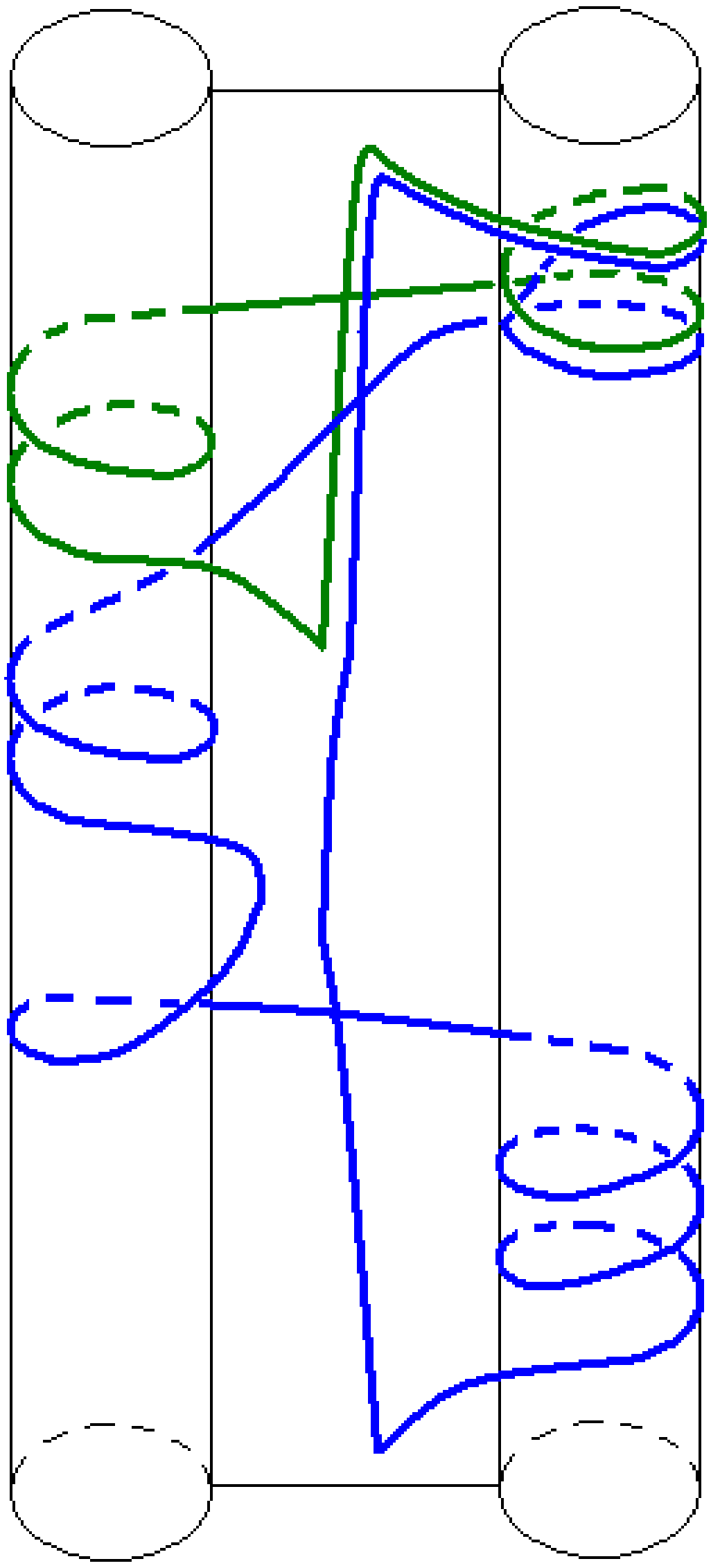

**Figure 185- The multiplication- drop down slided attaching curve- 5**

We repeat this process on the right generator cylinder:

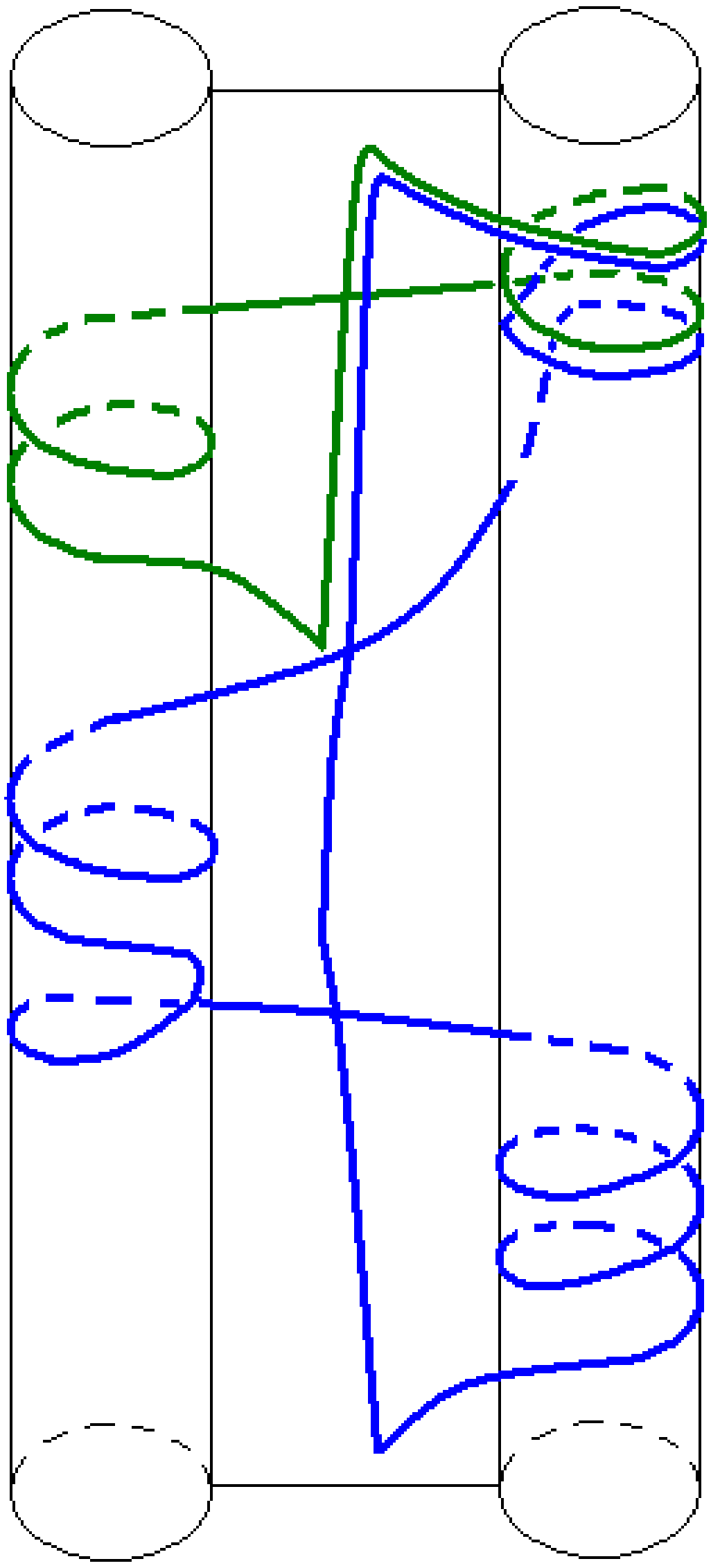

**Figure 186- The multiplication- drop down slided attaching curve- 6**

(The next  picture describes the same situation as the previous one, but it is easier to see the next step)

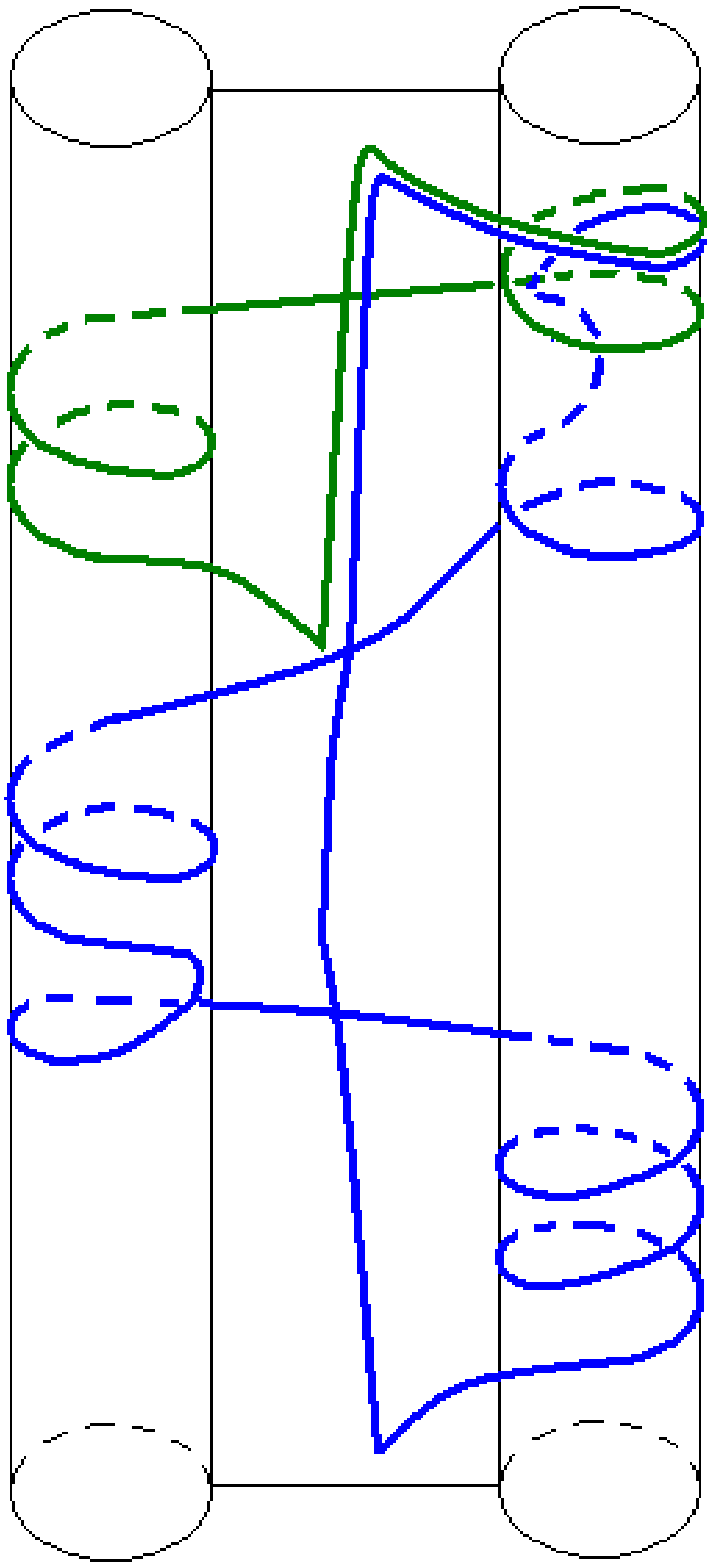

**Figure 187- The multiplication- drop down slided attaching curve- 7**

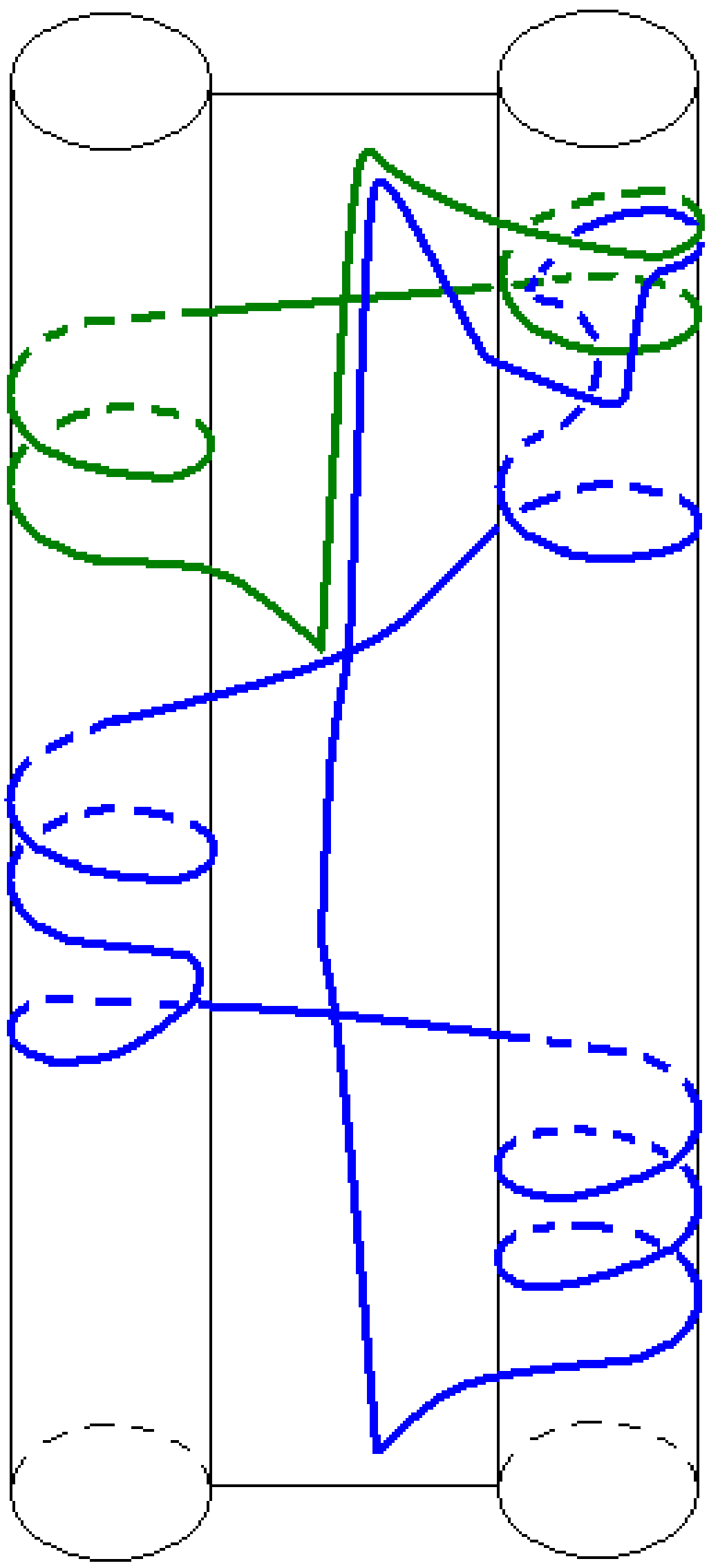

**Figure 188- The multiplication- drop down slided attaching curve- 8**

Now we can perform $T_2^{-1}$, so that the blue curve has no intersections with the green curve on the right generator cylinder (we say the blue curve is free on the right generator cylinder):

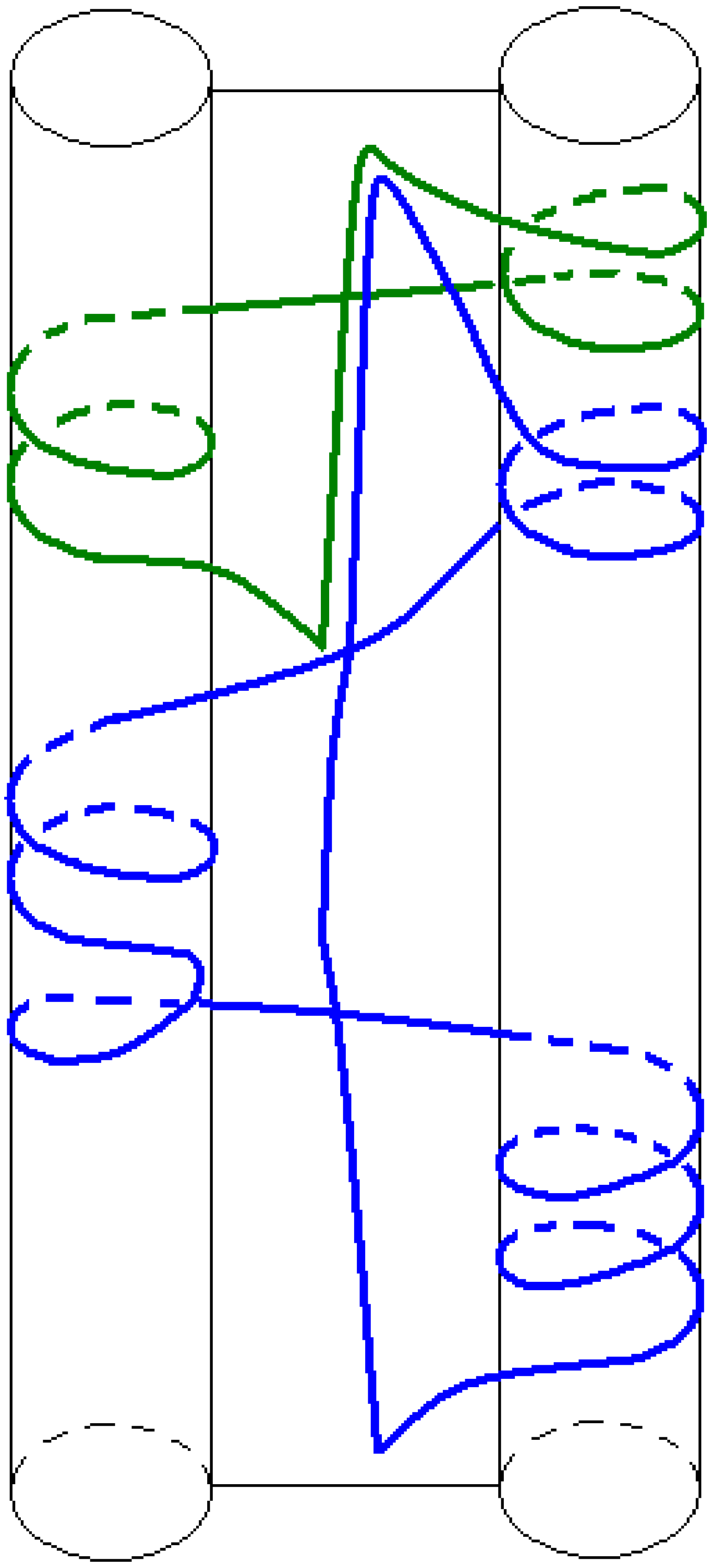

**Figure 189- The multiplication- drop down slided attaching curve- 9**

It remains to make the blue curve free on the rectangle:

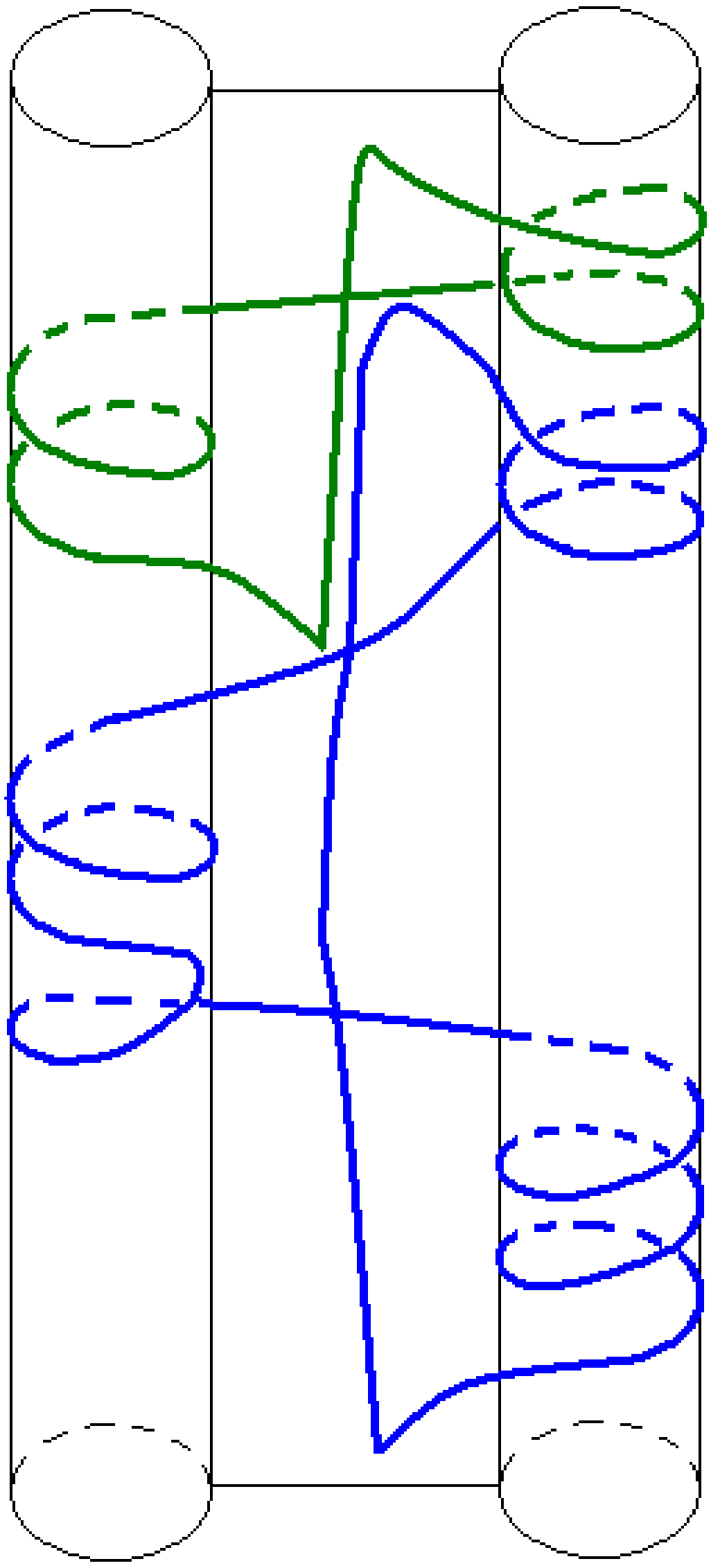

**Figure 190- The multiplication- drop down slided attaching curve- 10 (end)**

Now we can drop down the line to the rest without any crossing with the green curve.

We will consider the Matveev moves and the sequence of slices in more detail, so we have to study all possible crossings of the blue arc along the green curve. Note that in general we add flanges to get the local transitions from the Quinn list, but these are only drawn by the sequence of slices.

The entry of the slide:
We pass the entry of the 2-cell R (the green curve) with the blue arc and the slices are shown too.

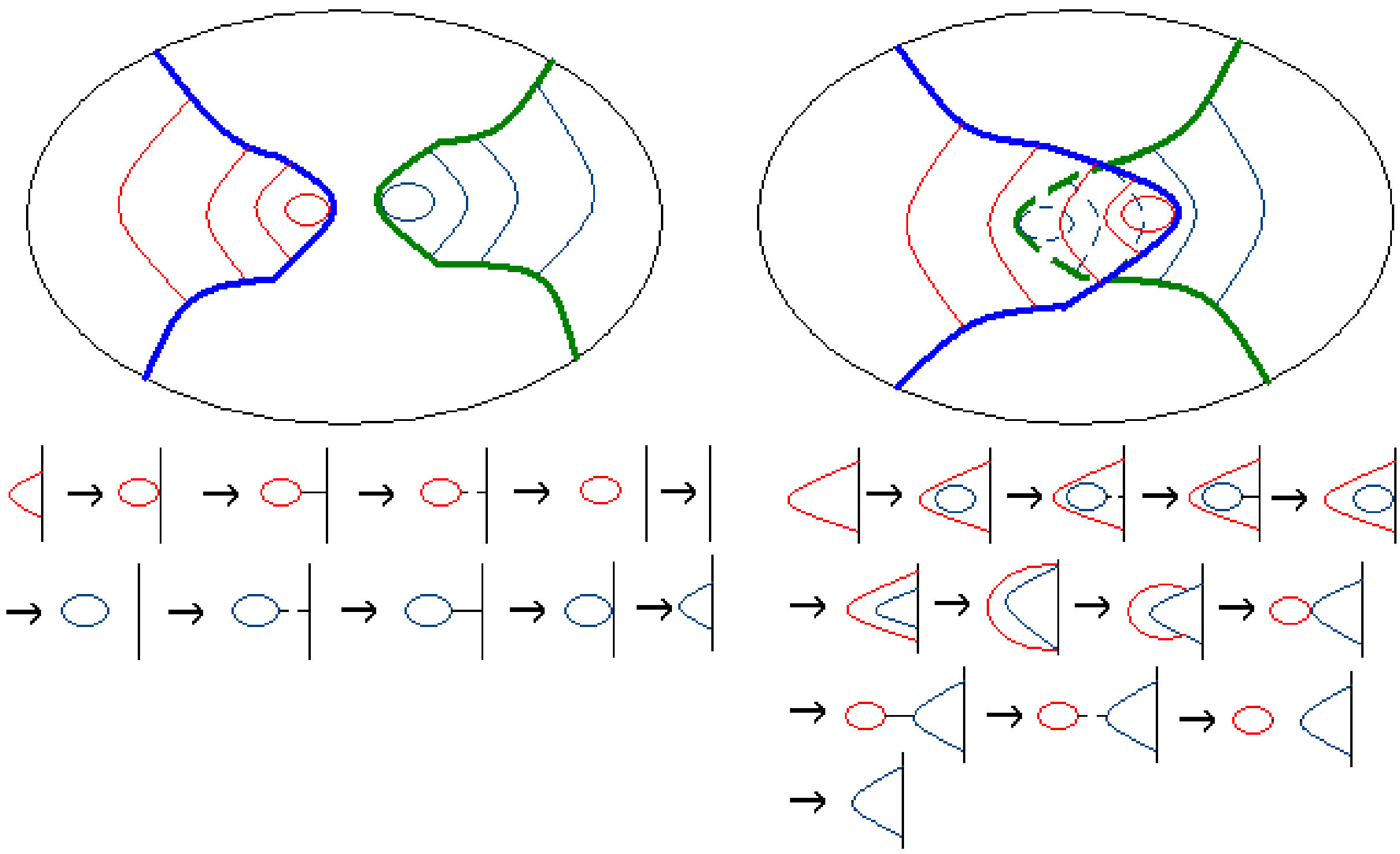

**Figure 191- The multiplication- entry of  the slide- sequence of slices**

Then we study all the crossings of the blue slided arc along  the green curve:

a) crossing from rectangle to generator cylinder:

This corresponds to the Matveev move T*. Note that some parts of the blue arc are on the 2-cell R and on the rectangle:

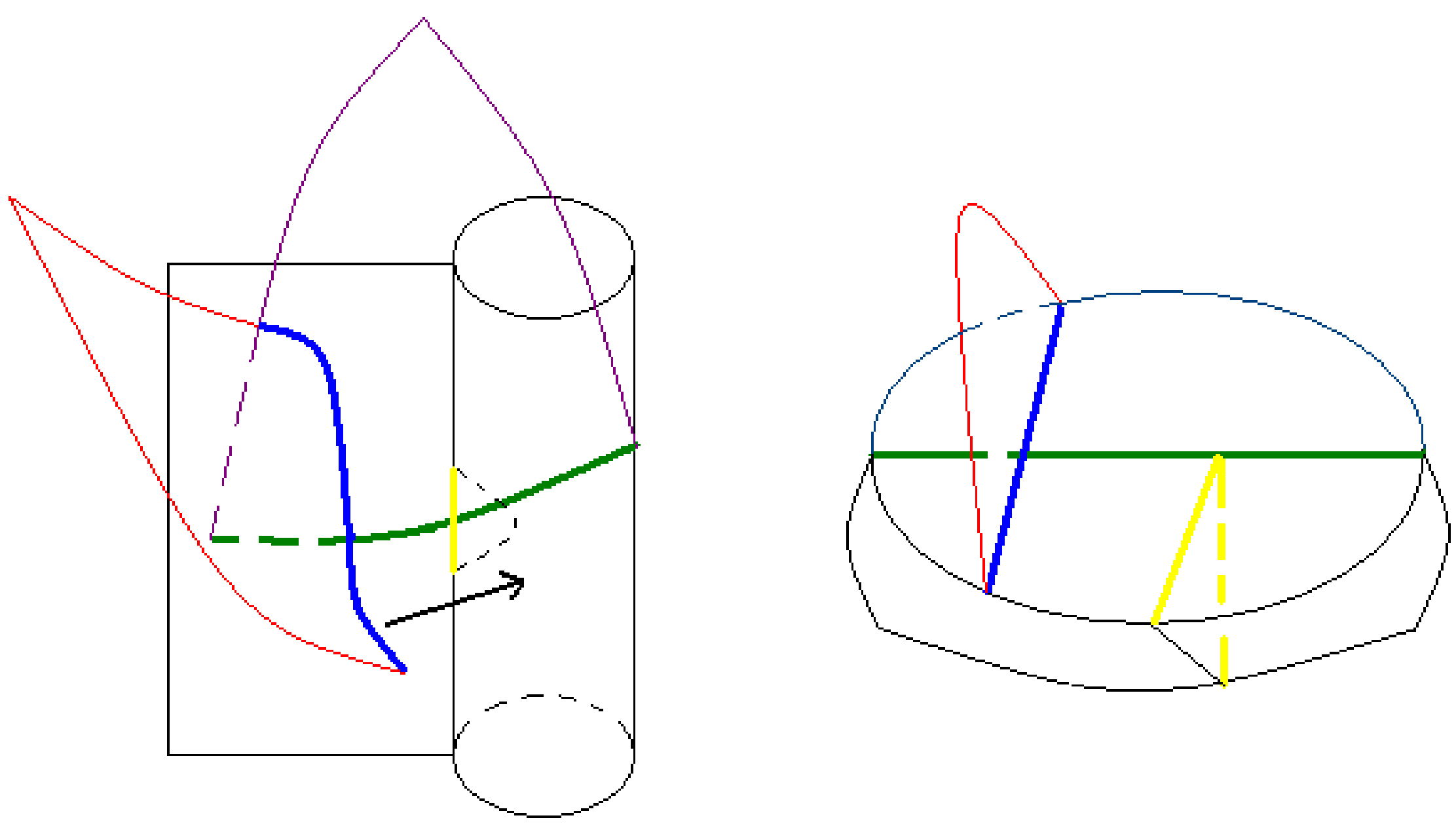

**Figure 192- The multiplication- a) crossing from rectangle to generator cylinder**

b) crossing from generator cylinder to rectangle:
This corresponds to the Matveev move T*. Note that some parts of the blue arc are on the 2-cell R and on the generator cylinder:

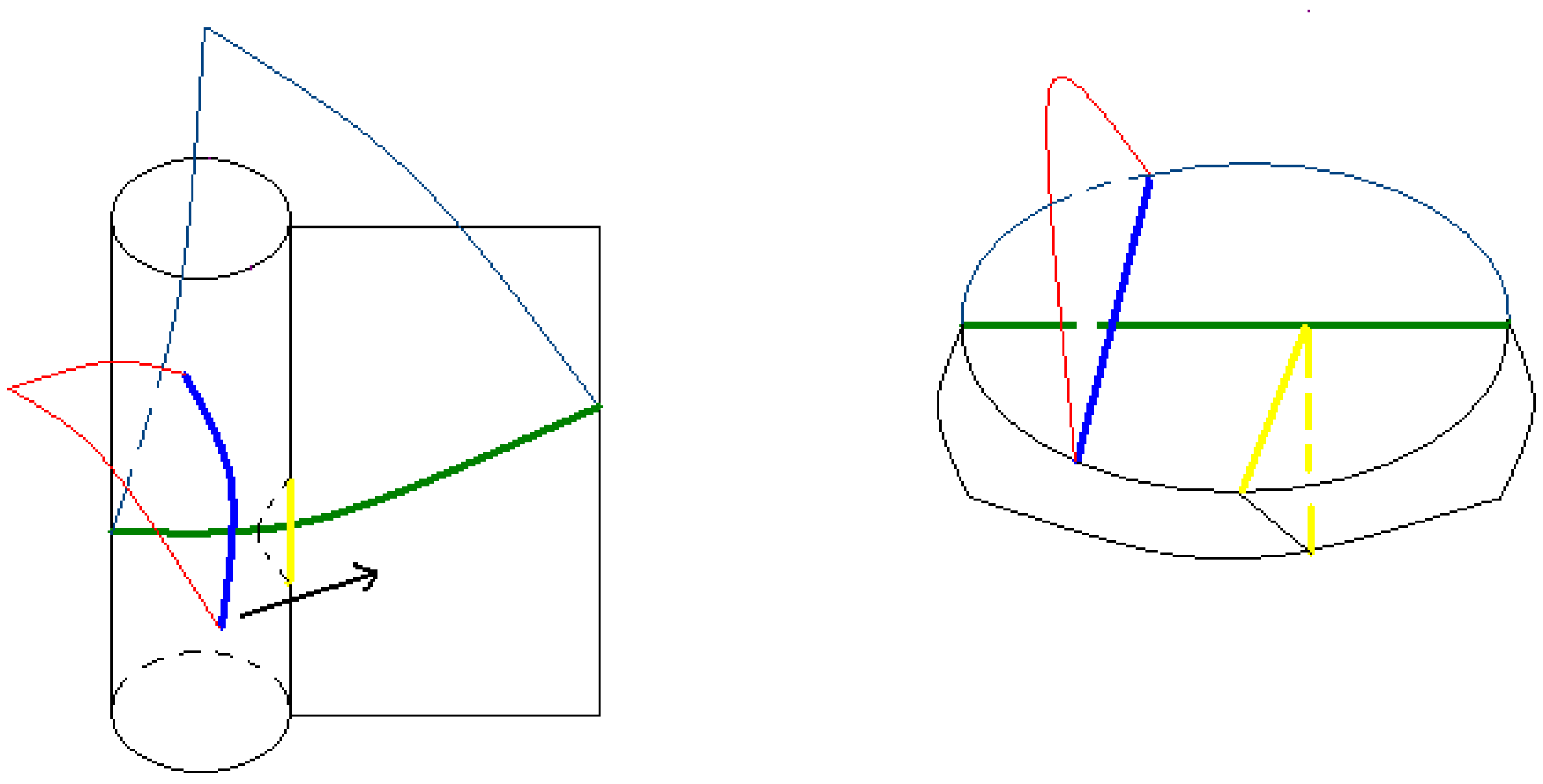

**Figure 193- The multiplication- b) crossing from generator cylinder to rectangle**

c) crossing a selfintersection of the green curve:
This corresponds to the Matveev move T*. Note that some parts of the blue arc are on the 2-cell R and on the rectangle:

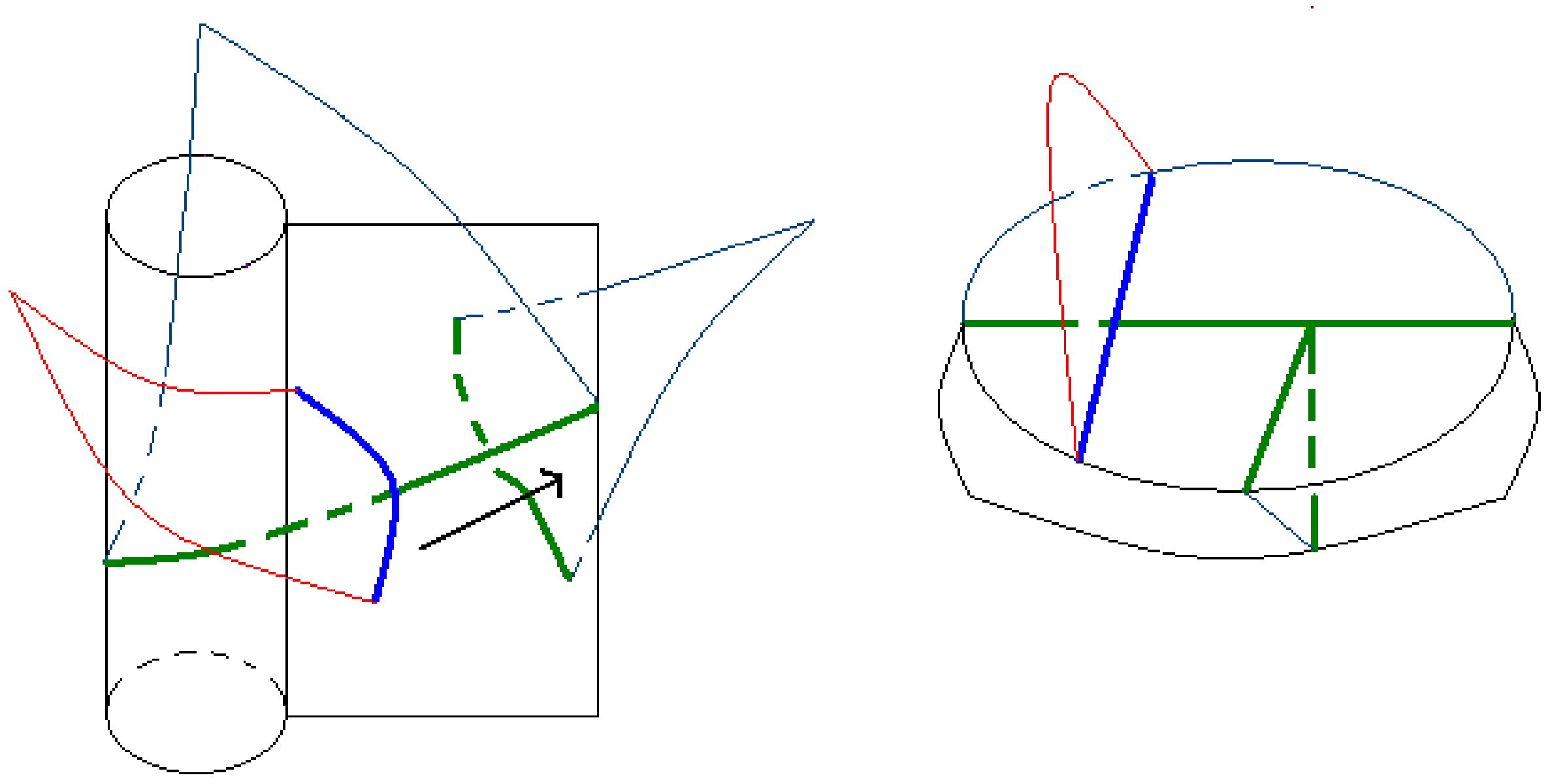

**Figure 194- The multiplication- c) crossing a selfintersection- version 1**

d) crossing a selfintersection of the green curve:
This corresponds to the Matveev move T*. Note that some parts of the blue arc are on the 2-cell R (drawn in backside) and on the rectangle:

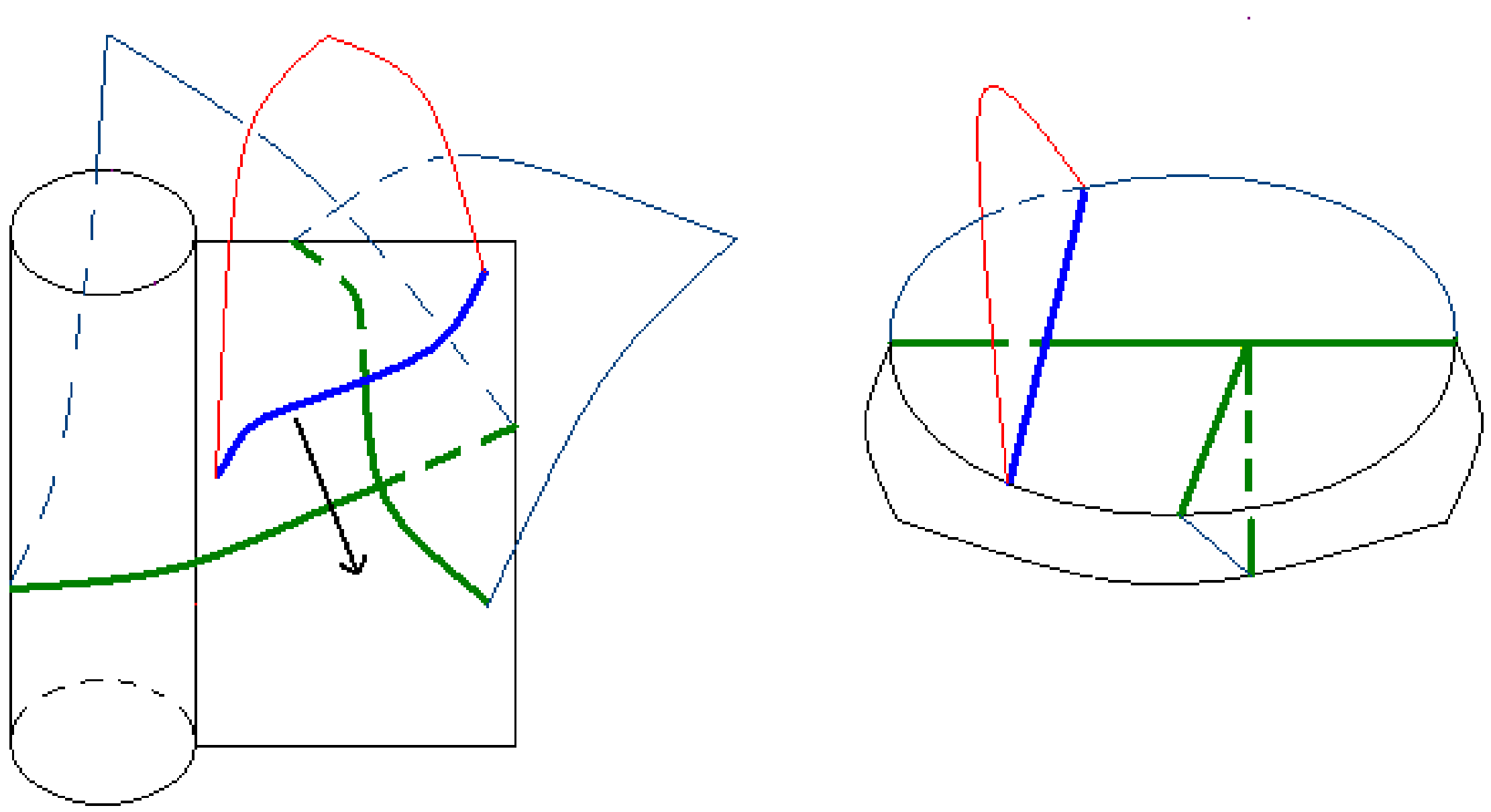

**Figure 195- The multiplication- d) crossing a selfintersection- 2**

The exit from the 2-cell R (we use the figure for entry of the slide)
Move the blue curve away from the 2-cell R and the slices are shown too.

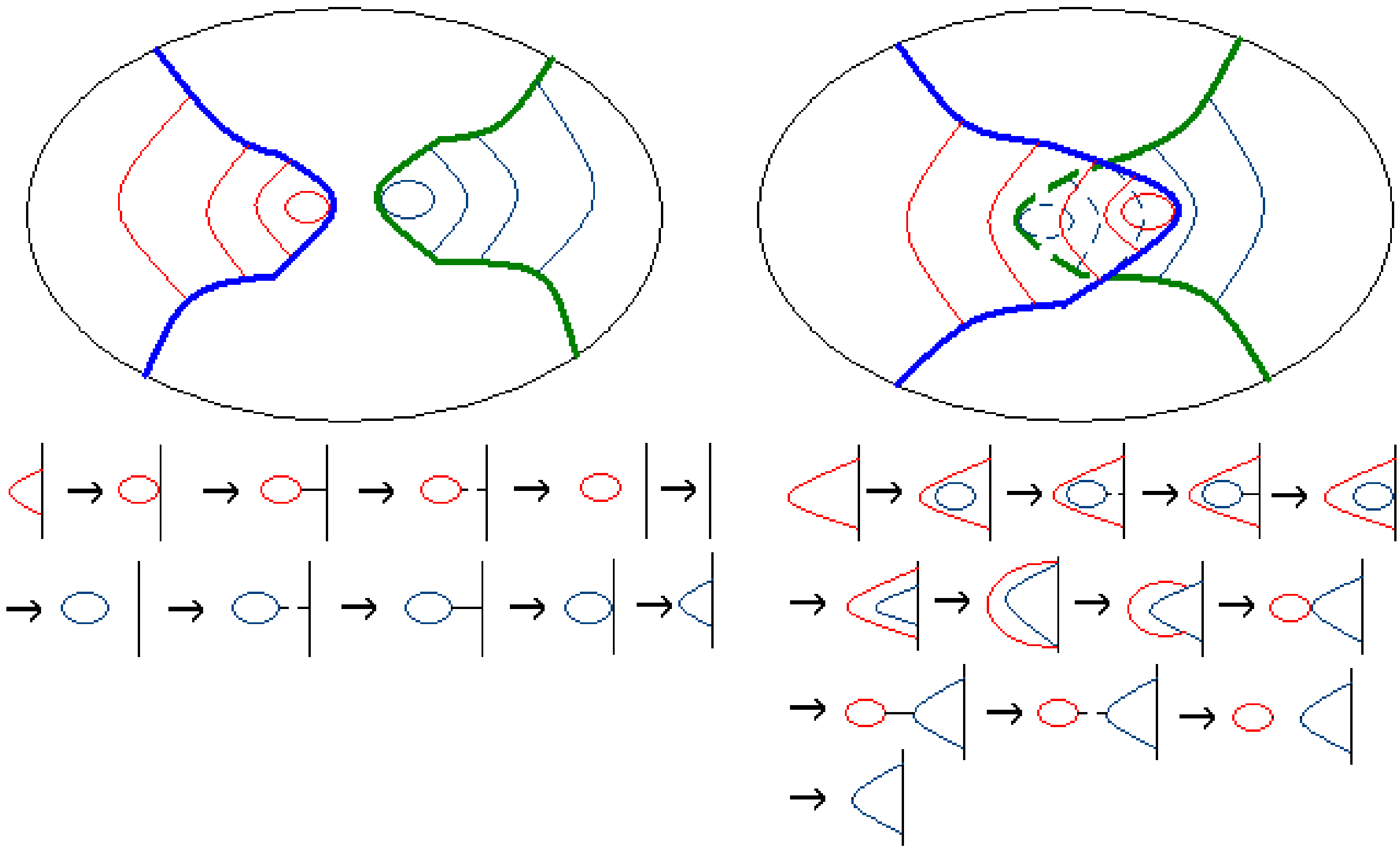

**Figure 196- The multiplication- exit from the 2-cell- sequence of slices**

We have developed the cases a) - d) for the sequence of slices for T*. Since the slices are line segment slices, we list the two cases:

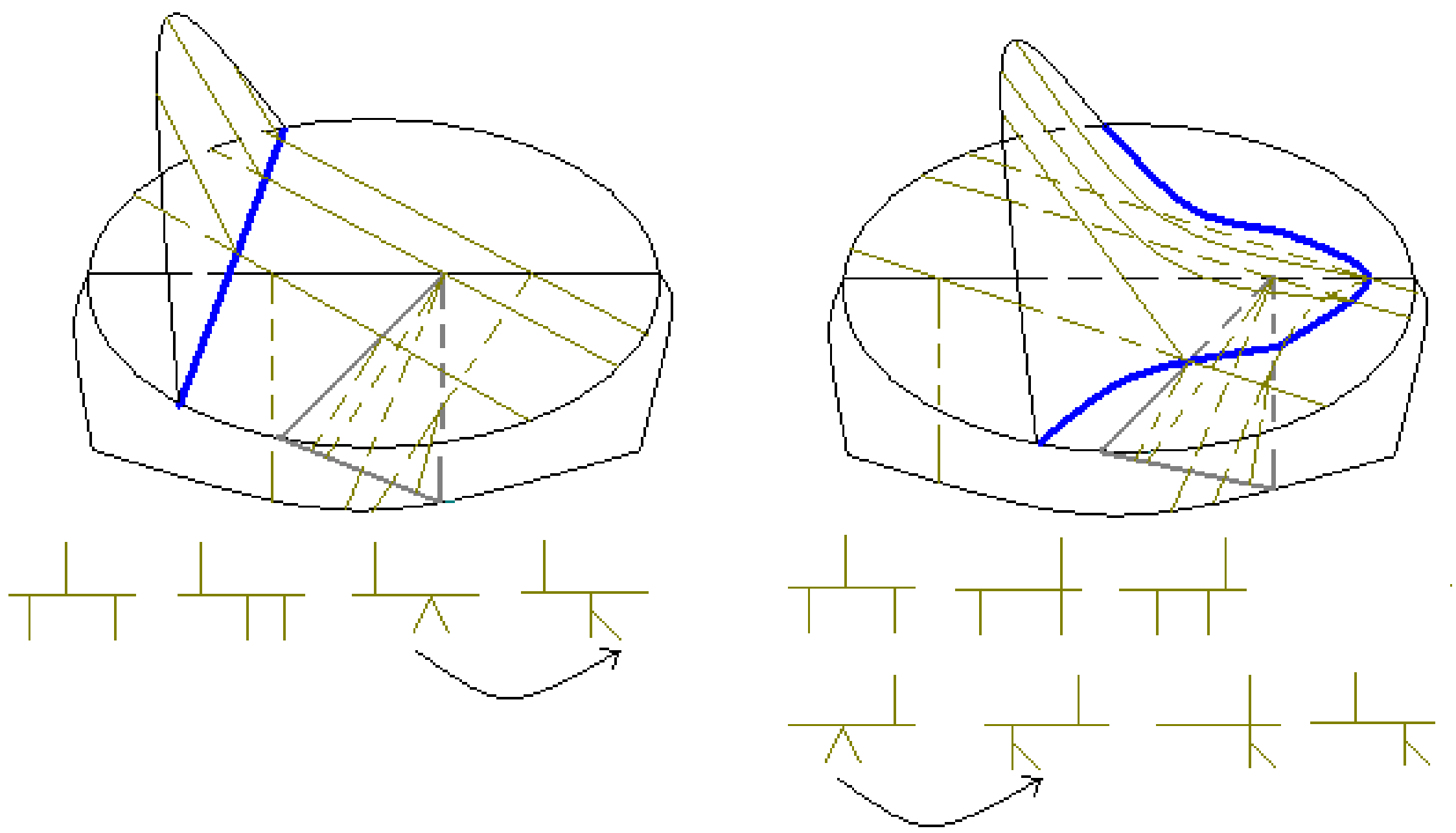

**Figure 197- The multiplication- T* move- sequence of slices- version 1**

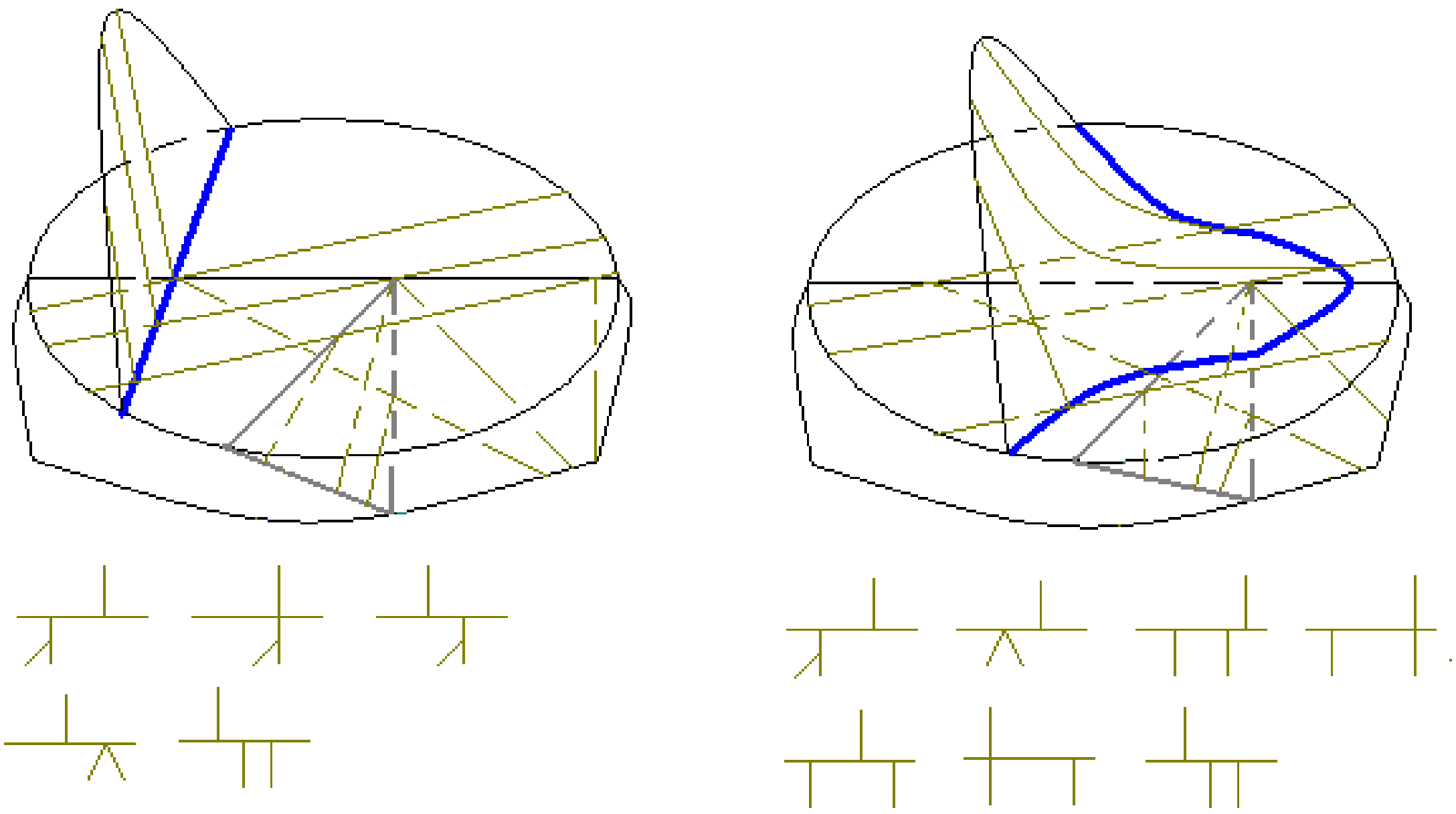

**Figure 198- The multiplication- T* move- sequence of slices- version 2**

From the process of dropping down the blue curve we will pick up three standard cases:

(1) leads to the Matveev move T*, so the sequence of slices is known:

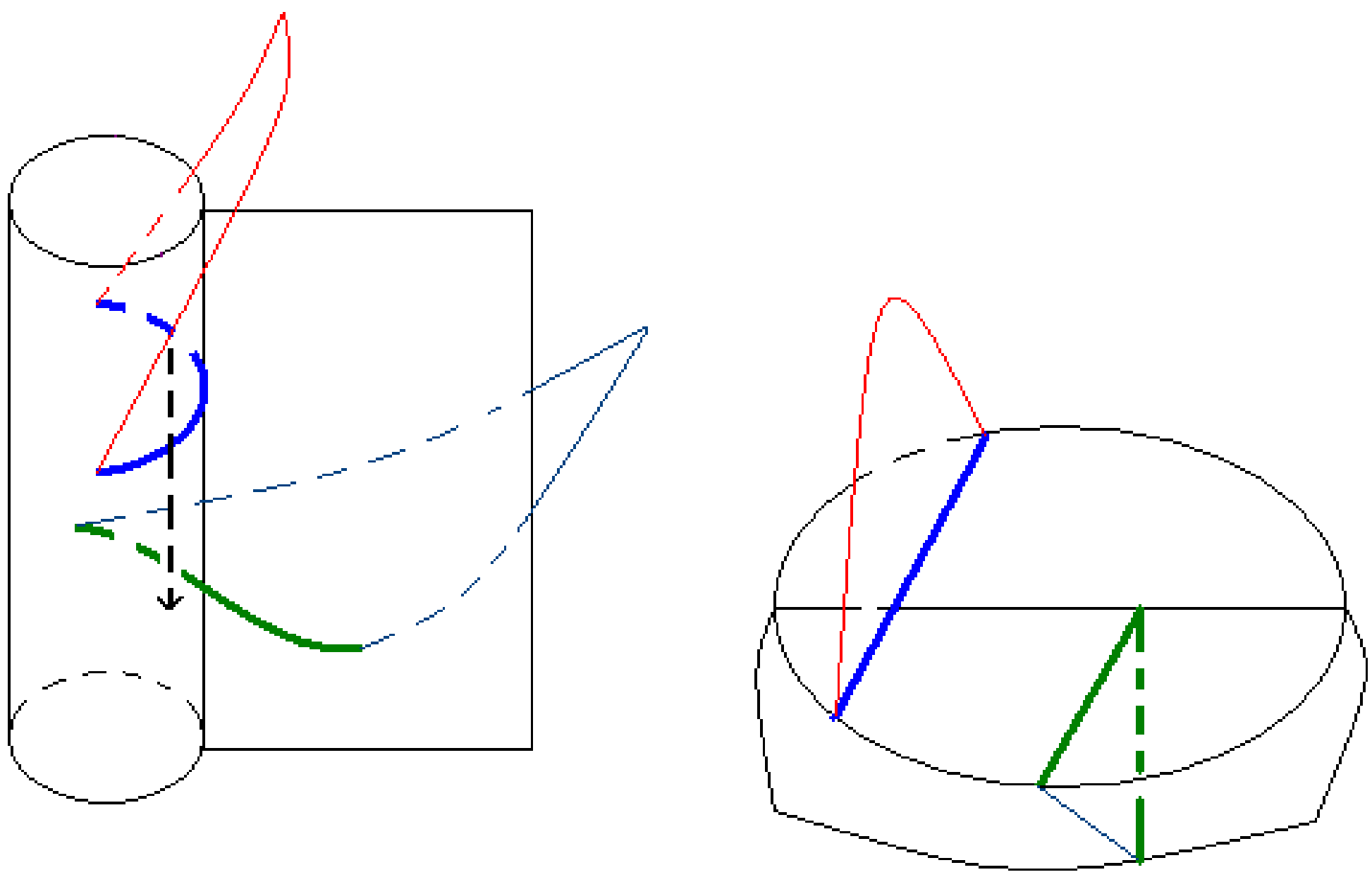

**Figure 199- The multiplication- drop down slided attaching curve- T* move**

(2) leads to the Matveev move $T_2^{-1}$:

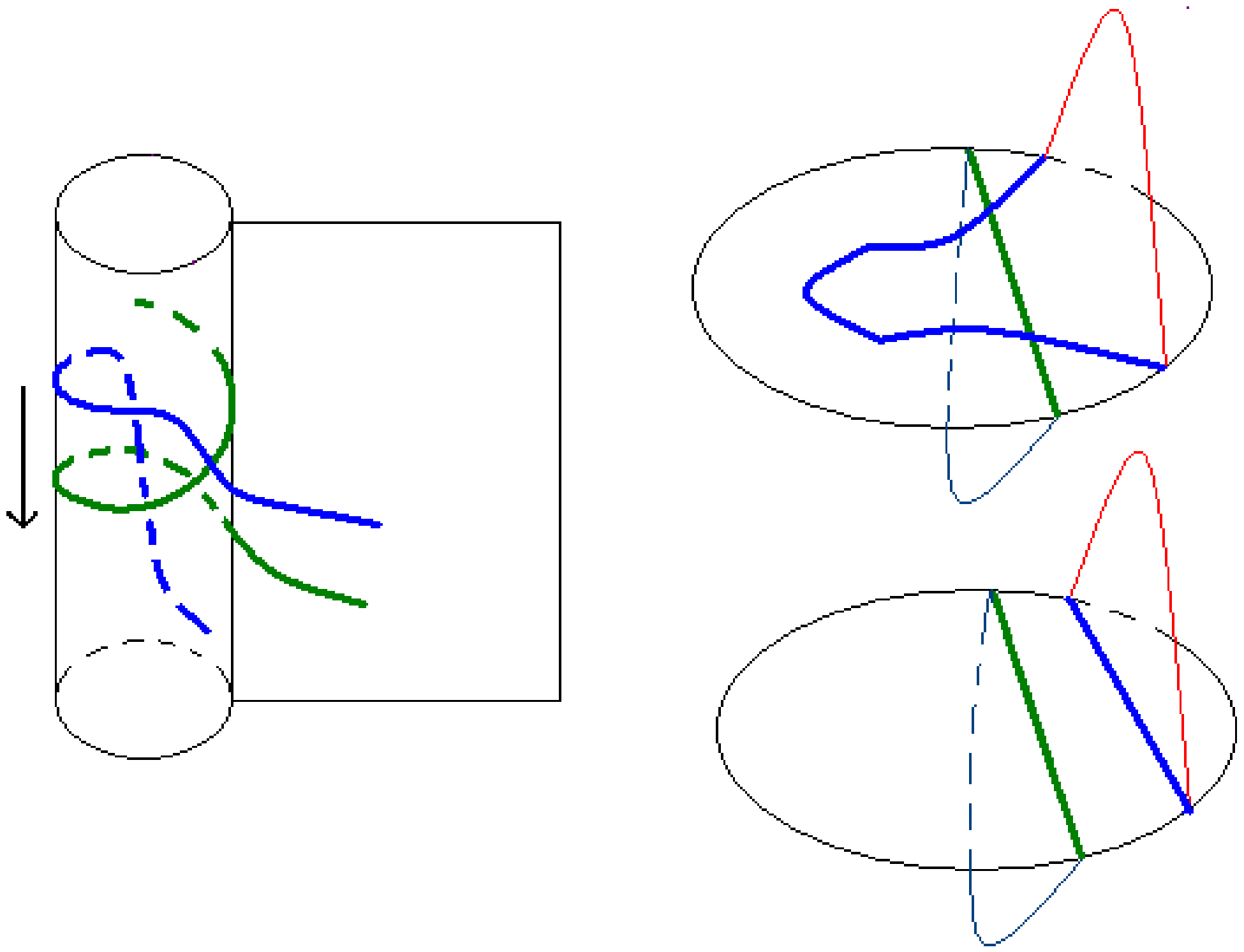

**Figure 200- The multiplication- drop down slided attaching curve- $T_2^{-1}$ move**

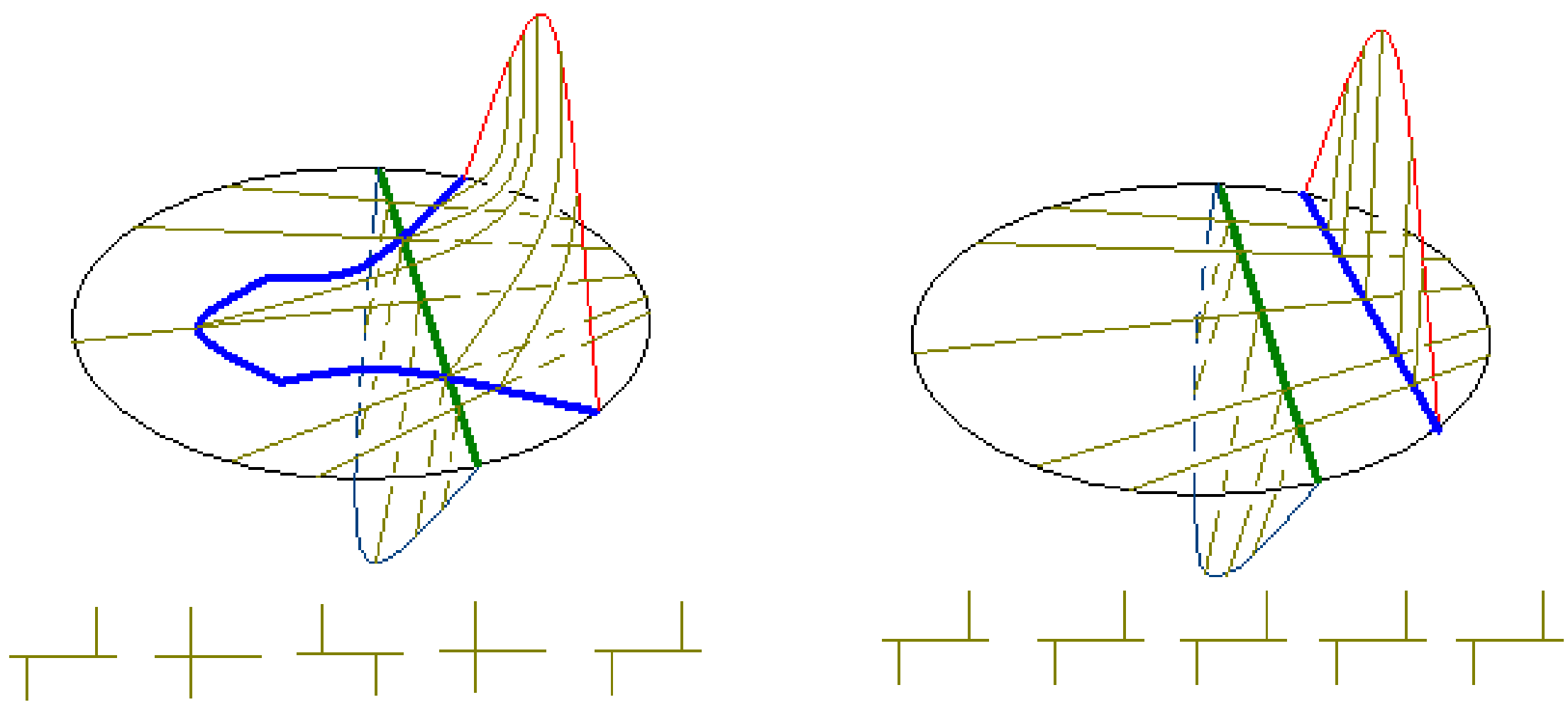

**Figure 201- The multiplication- drop down slided attaching curve- $T_2^{-1}$ move- sequence of slices**

(3) is is a composition of Matveev moves, we called it in chapter 5.1 “modified $T_2$ move”:

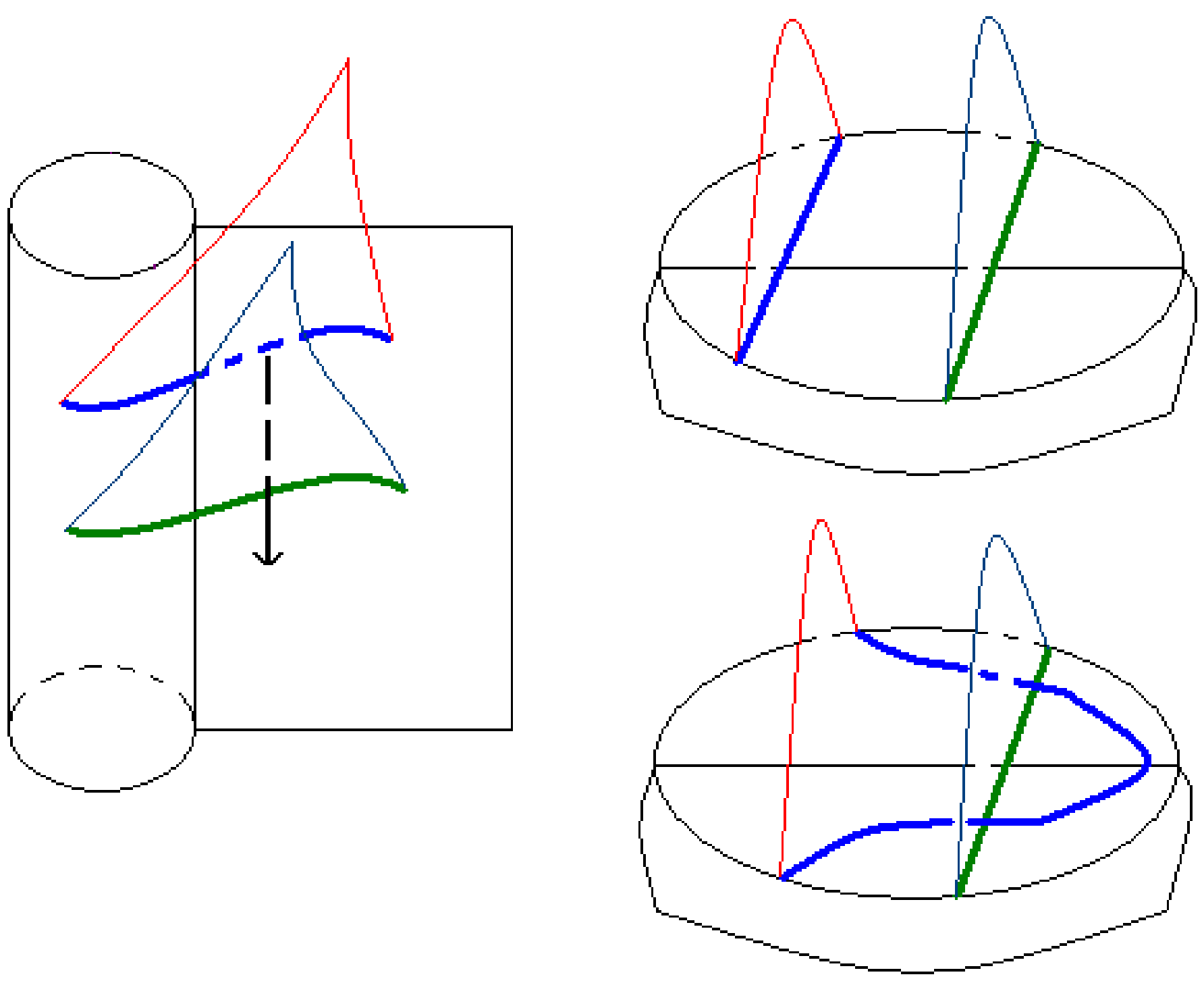

**Figure 202- The multiplication- drop down slided attaching curve- modified $T_2$ move**

We decompose it into its single Matveev moves:

i. $T_3$ move

Remark: The $T_3$ turn was chosen to be a “good” one:

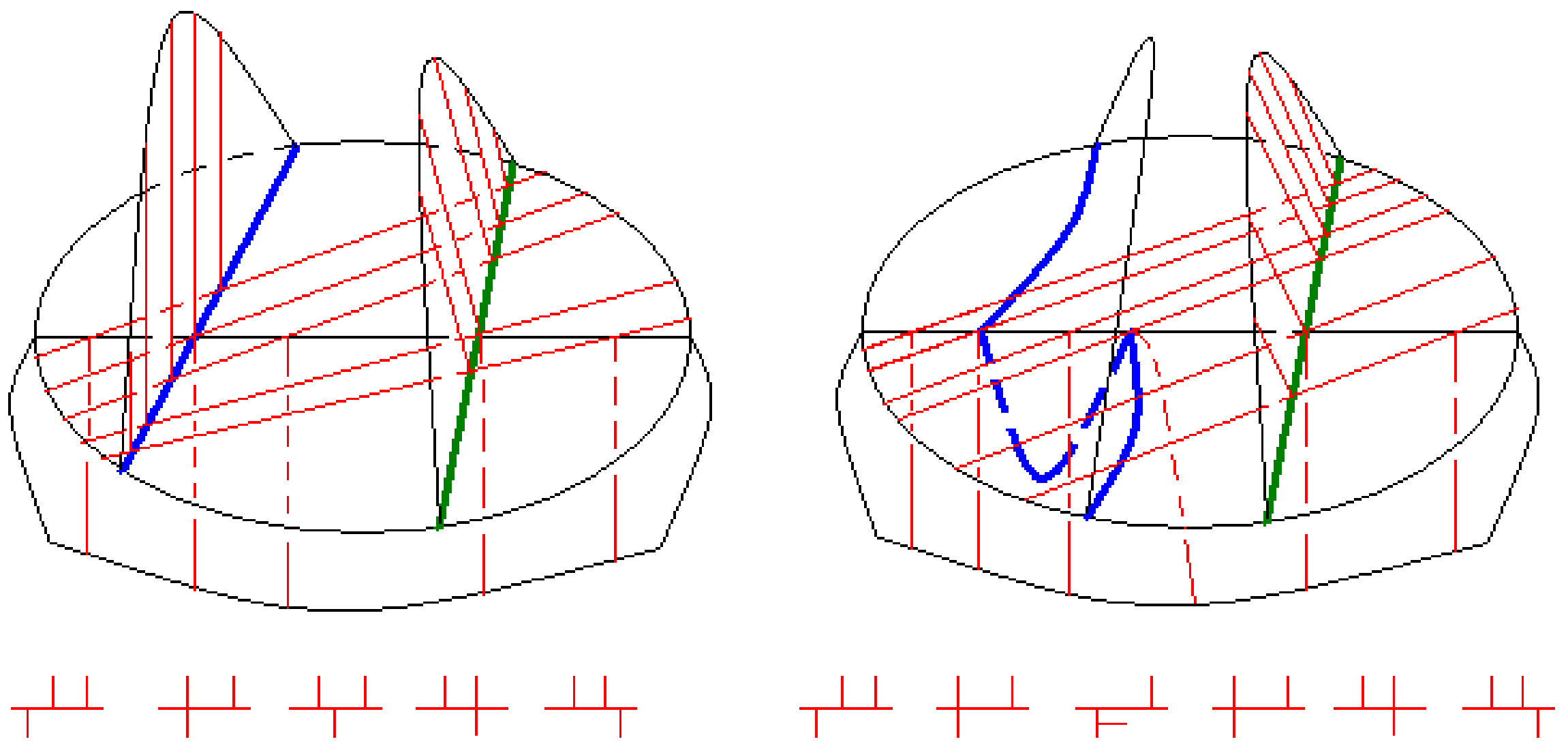

**Figure 203- The multiplication- modified $T_2$ move- sequence of slices- 1**

ii. T* move:

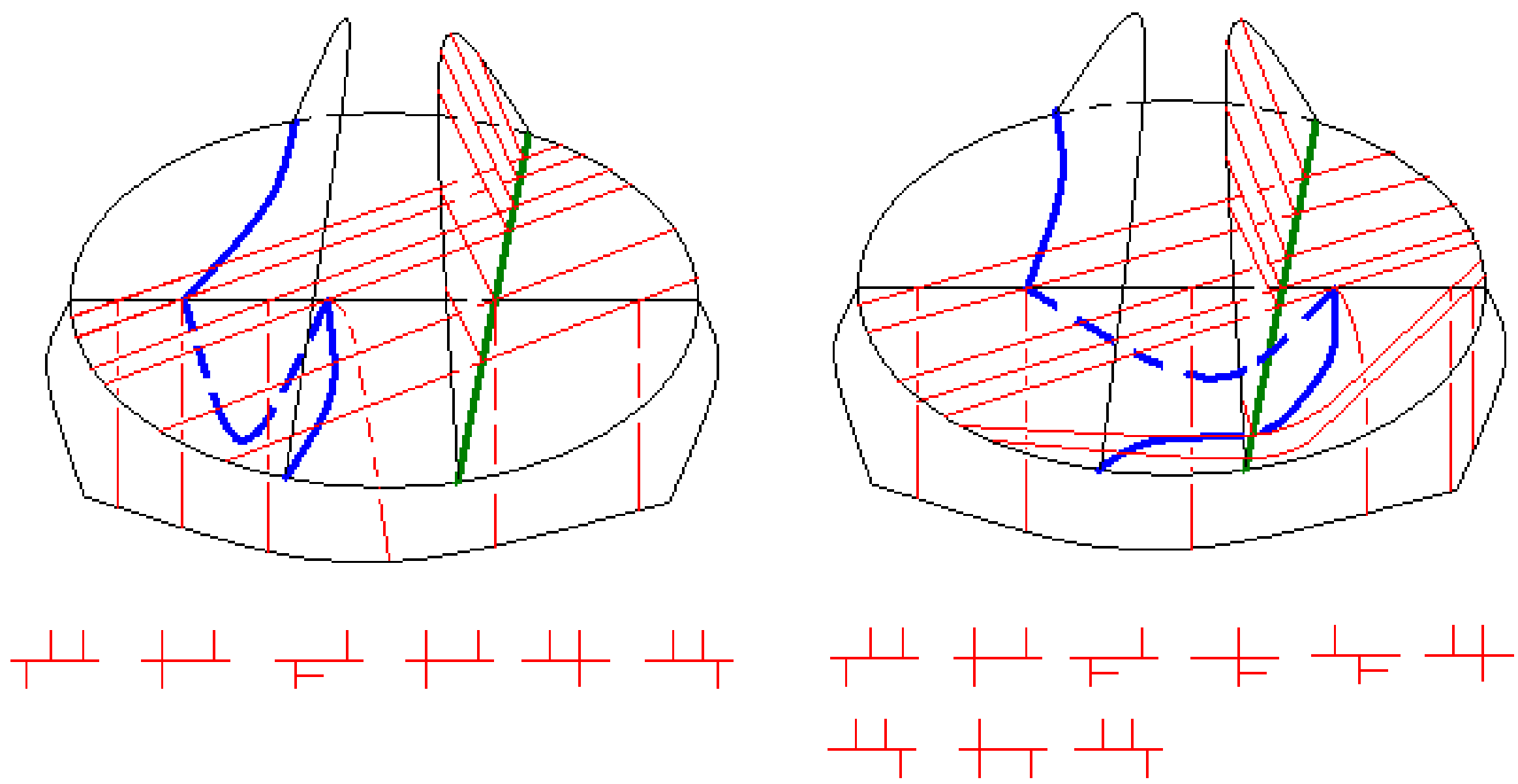

**Figure 204- The multiplication- modified $T_2$ move- sequence of slices- 2**

iii. T* move:

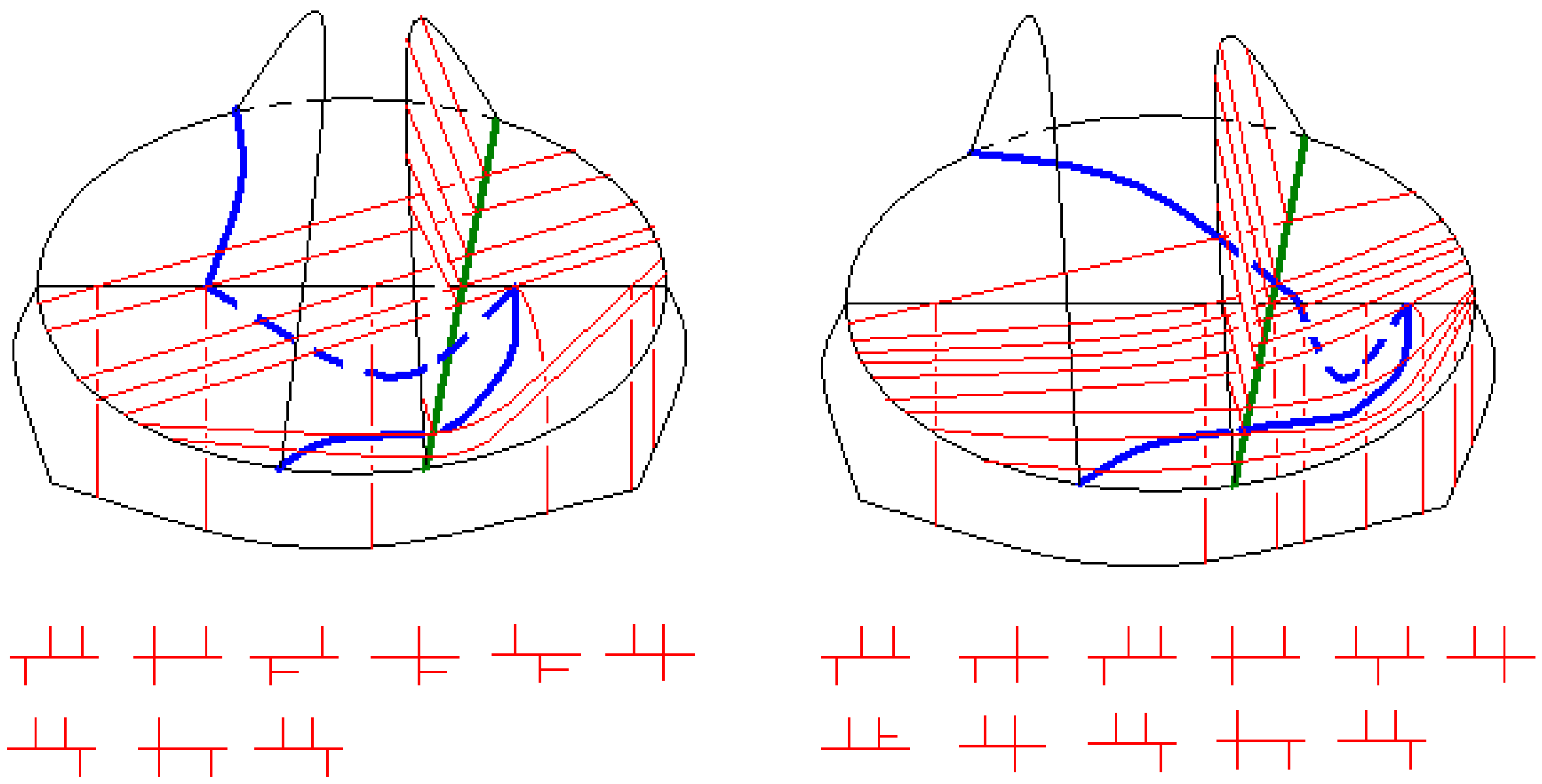

**Figure 205- The multiplication- modified $T_2$ move- sequence of slices- 3**

iv. $T_3^{-1}$ move:

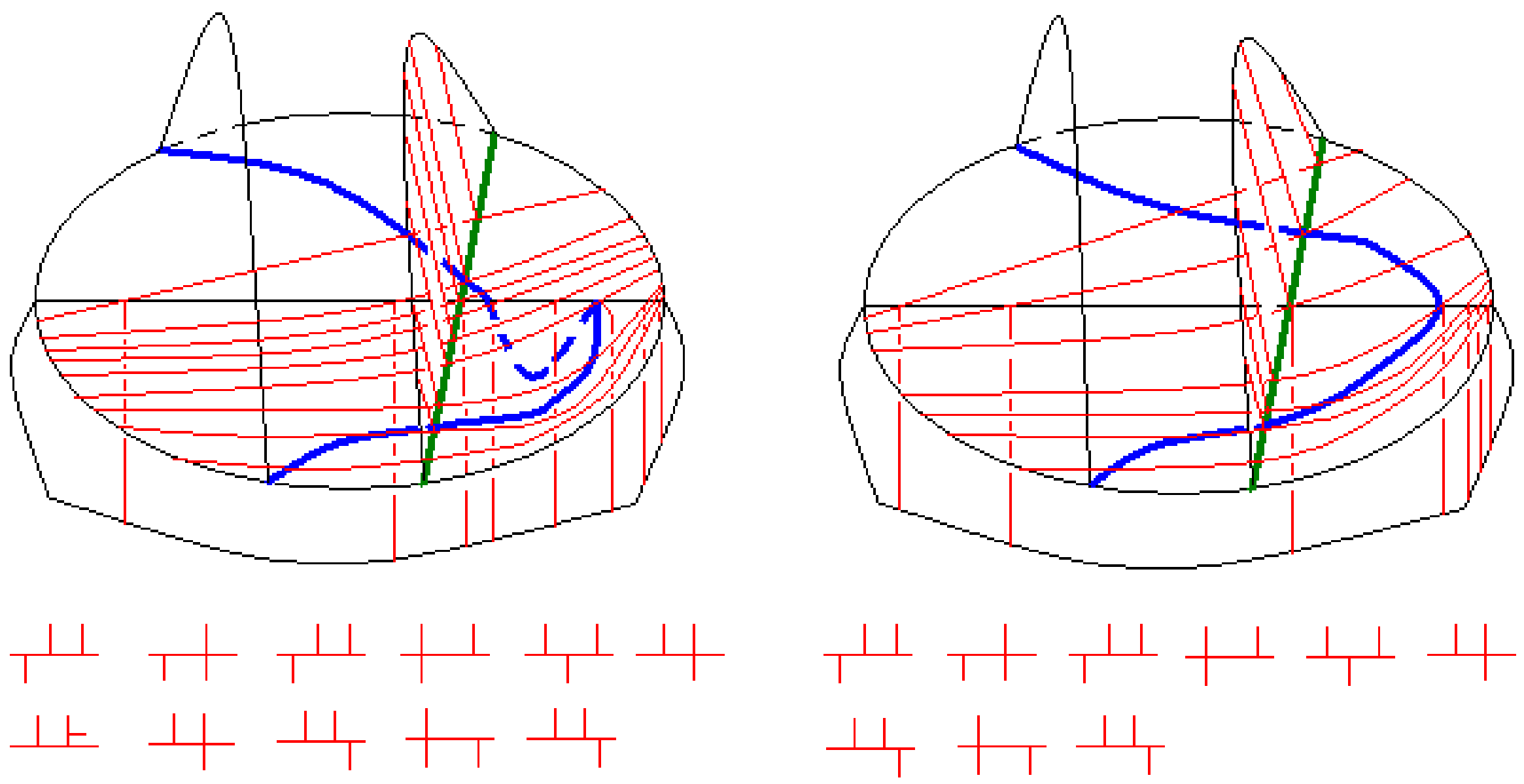

**Figure 206- The multiplication- modified $T_2$ move- sequence of slices- 4**

We do not want to forget the drop down on the rectangle, where we have to change the slices to push the extremum through the S-cell and afterwards rechange the slices:

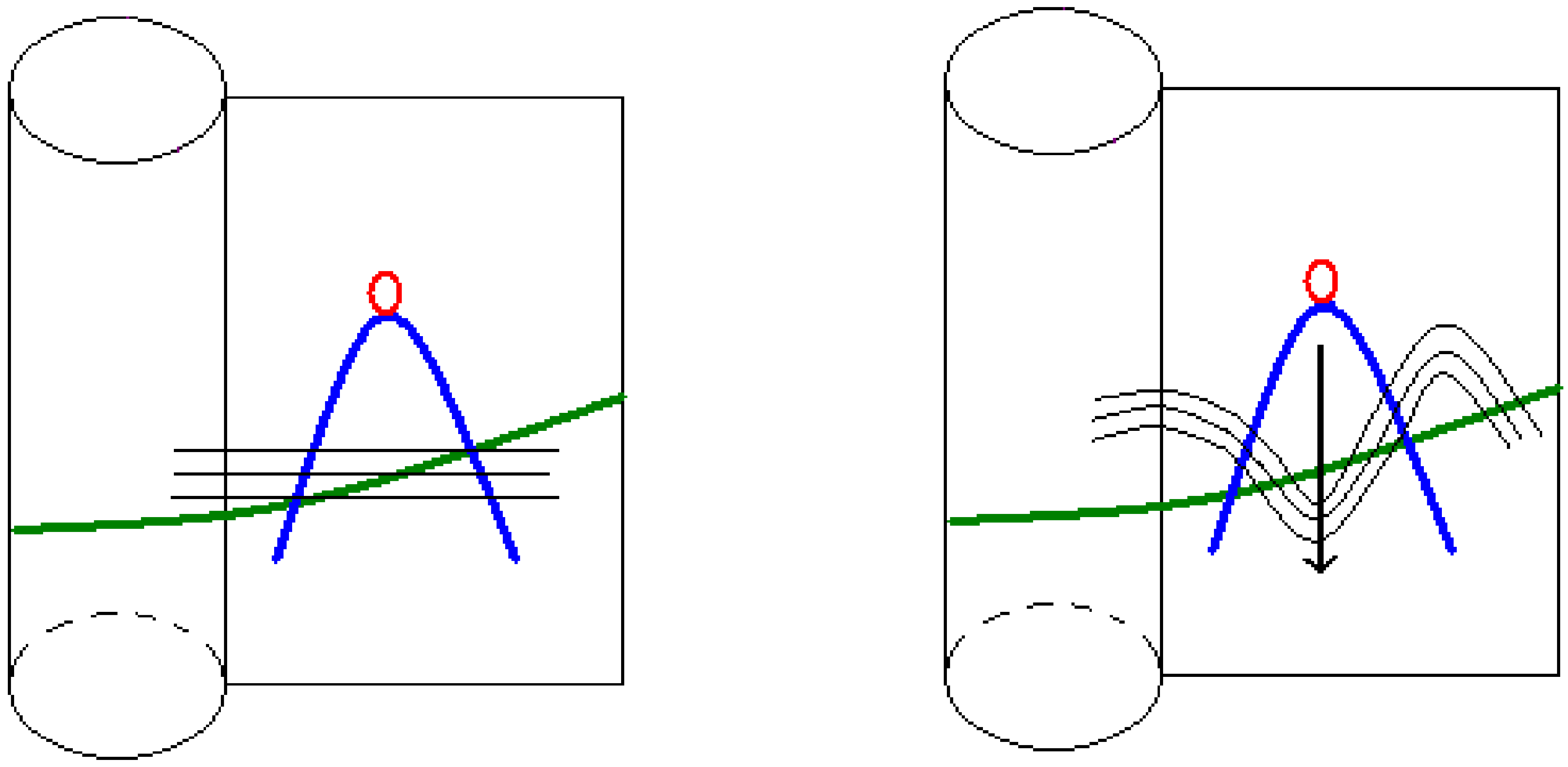

**Figure 207- The multiplication- drop down attaching curve on rectangle**

## *7.2 The conjugation*

We start with a relation R in the Quinn model:

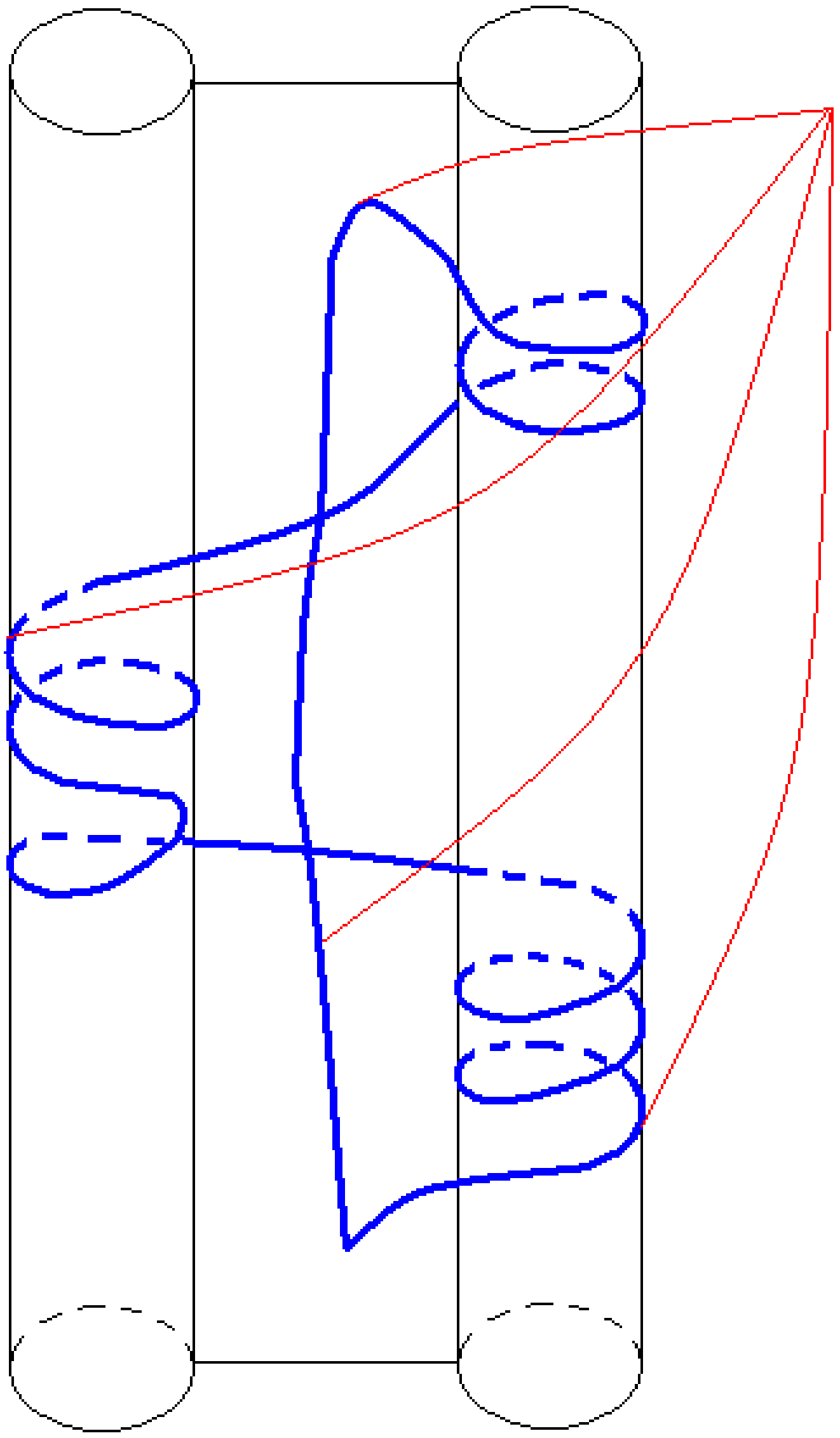

**Figure 208- The conjugation- a relation in Quinn model**

We conjugate the relation R with a single generator w and regard the attaching curve to the relation $wRw^{-1}$. We come back to the attaching curve of R, if we push the connecting arc (from level zero to level 1) around the generator cylinder. Doing this process many crossings with the R-curve occur and these result in Matveev moves.

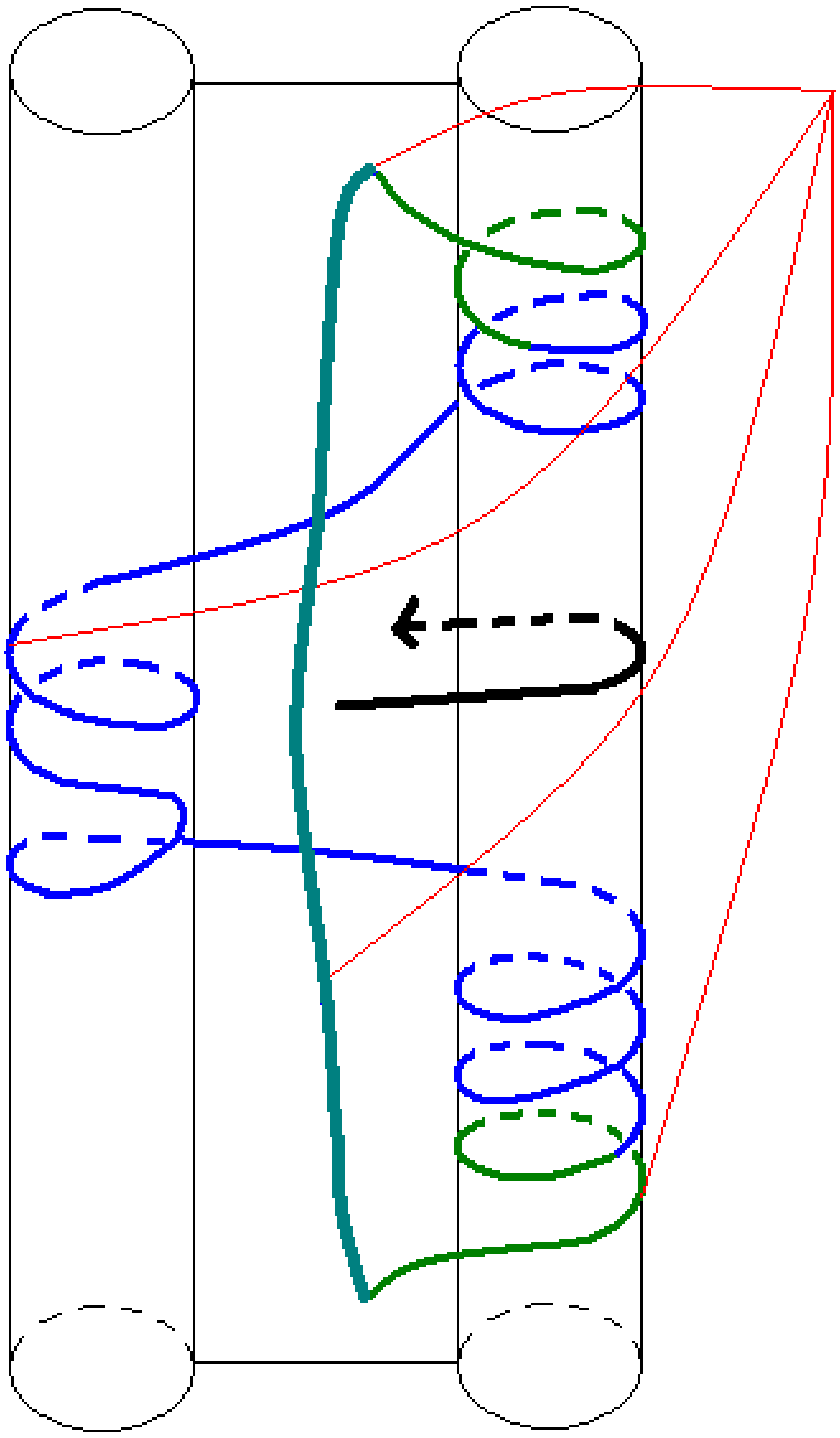

**Figure 209- The conjugation- transfer from conjugate of a relation to the relation in Quinn model**

Note, that it is a good advice to keep the level order of the attaching curve while defining a height function on it, (compare with chapter 7.1 the action to drop down the curve), so we have to drag the extrema with the connecting arc around the generator cylinder, too.
It is easy to shift extrema; it is not more than the game of introducing and cancelling pairs of saddlepoints:

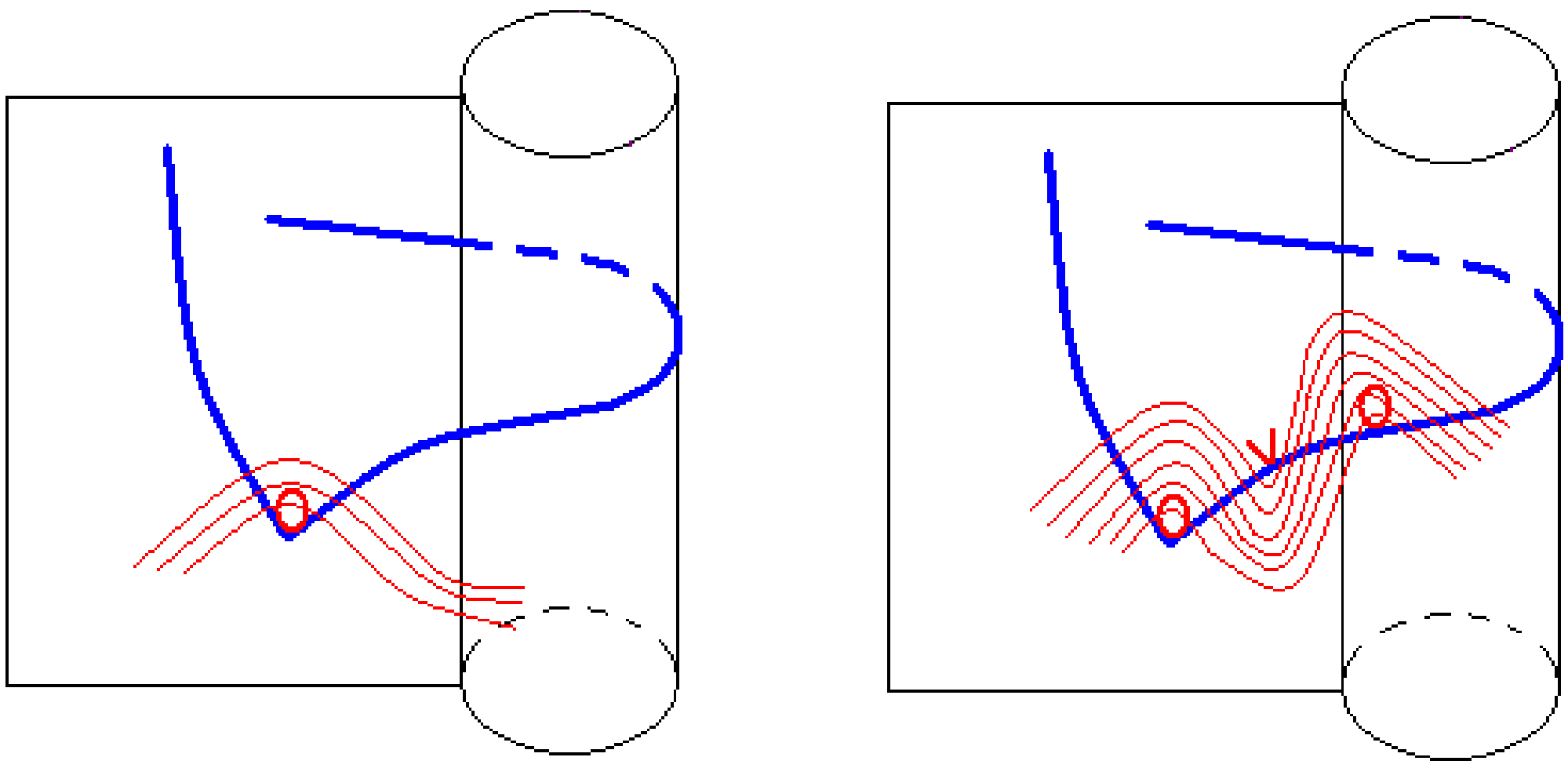

**Figure 210- The conjugation- shift a extrema in Quinn model**

We analyse the several variations how the attaching curve can turn around the generator cylinder and study the arising Matveev moves, when the connecting arc passes the vertices:

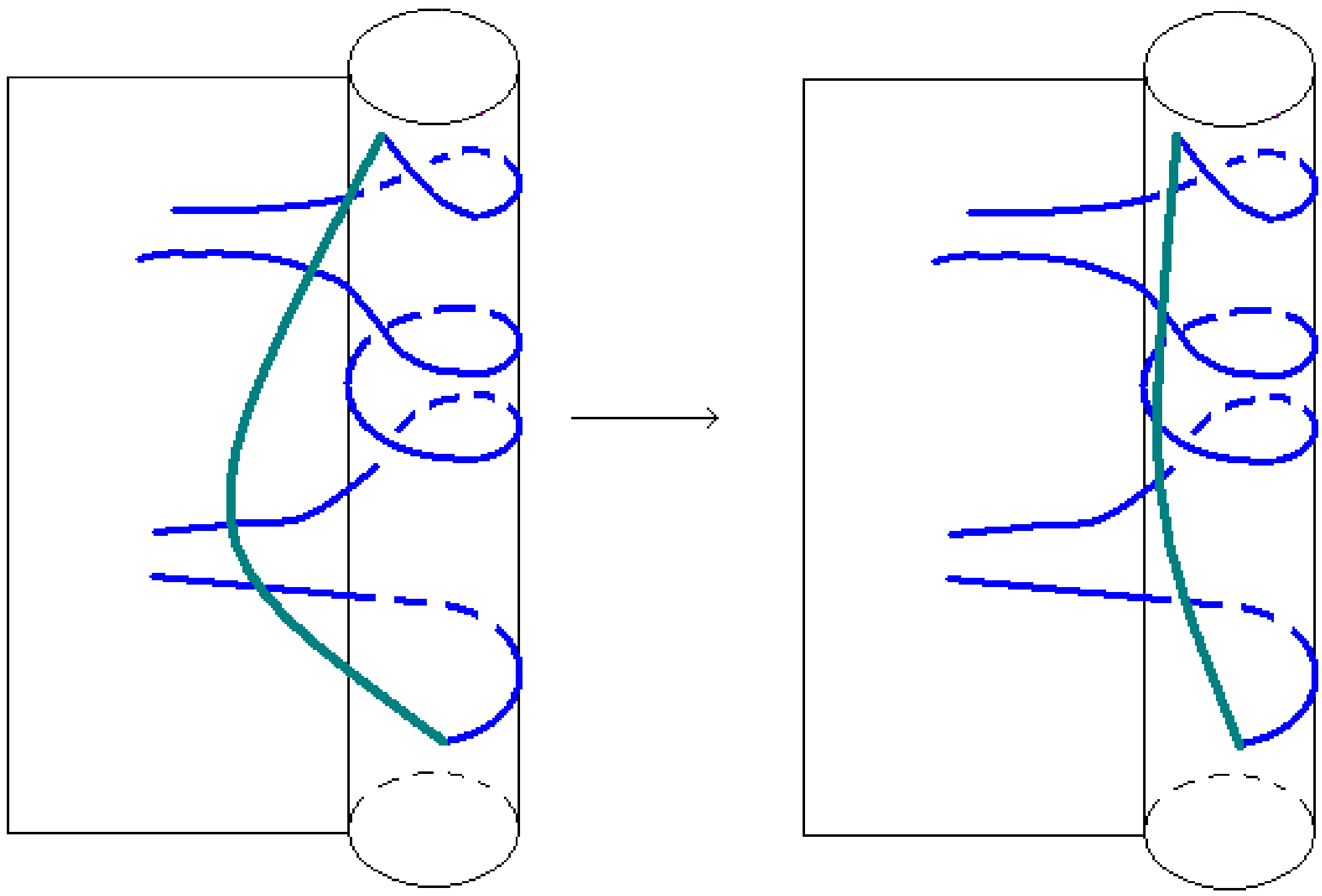

**Figure 211- The conjugation- list of crossings for the slided arc inQuinn model**

a) The curve wraps around the generator counterclockwise:

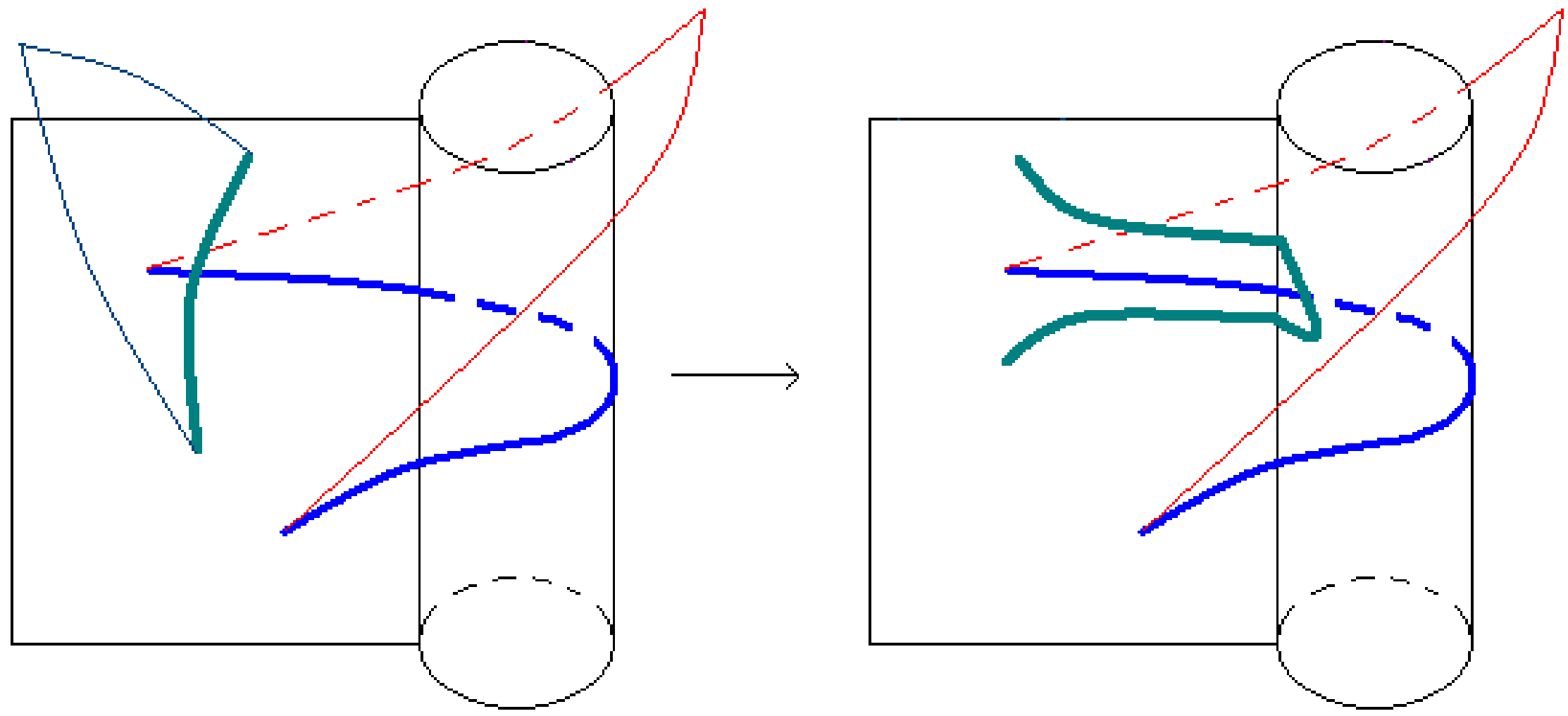

**Figure 212- The conjugation- list of crossing for slided arc- a) curve counterclockwise**

The corresponding Matveev move is $T_2$:

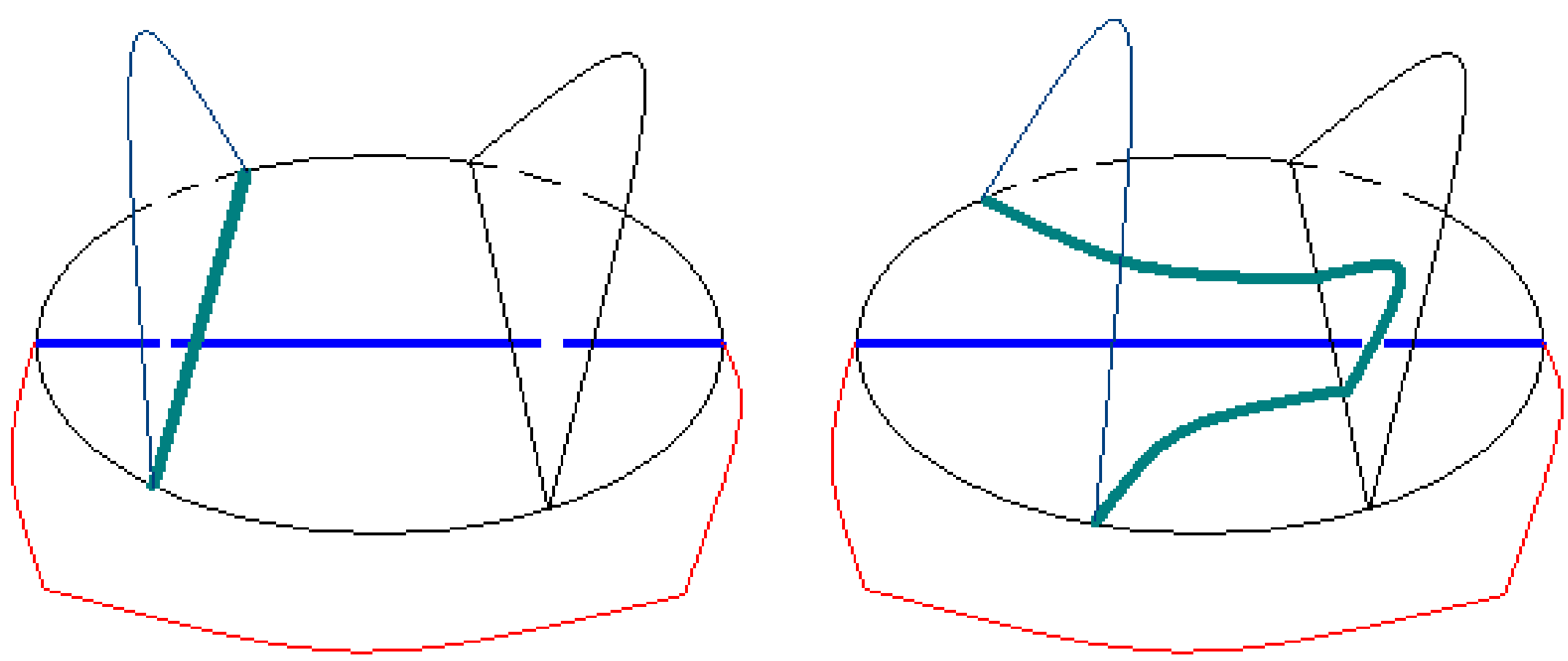

**Figure 213- The conjugation- a) curve counterclockwise- $T_2$ move**

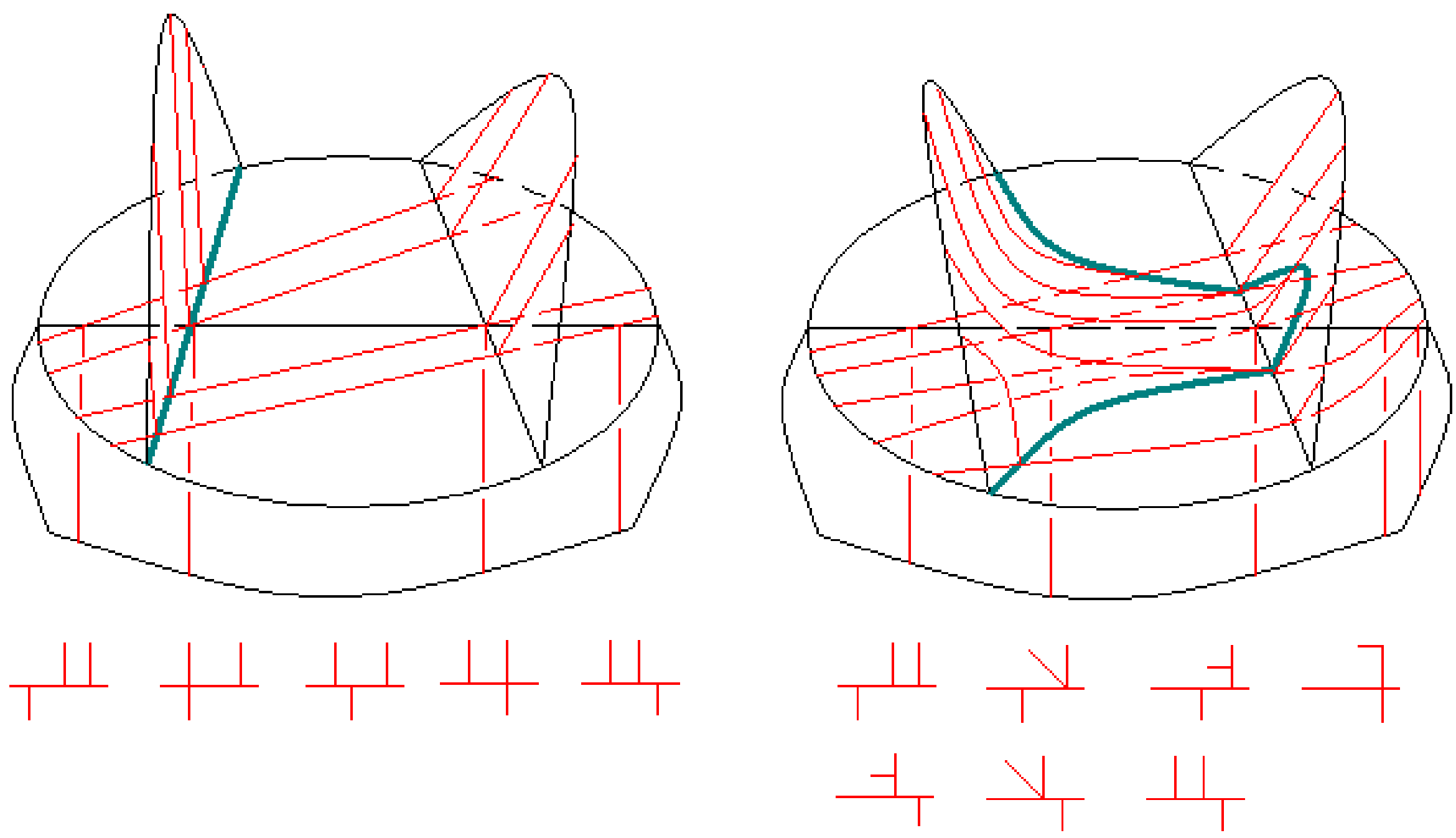

**Figure 214- The conjugation– a) curve counterclockwise- $T_2$ move– sequence of slices**

b) The curve wraps around the generator clockwise:

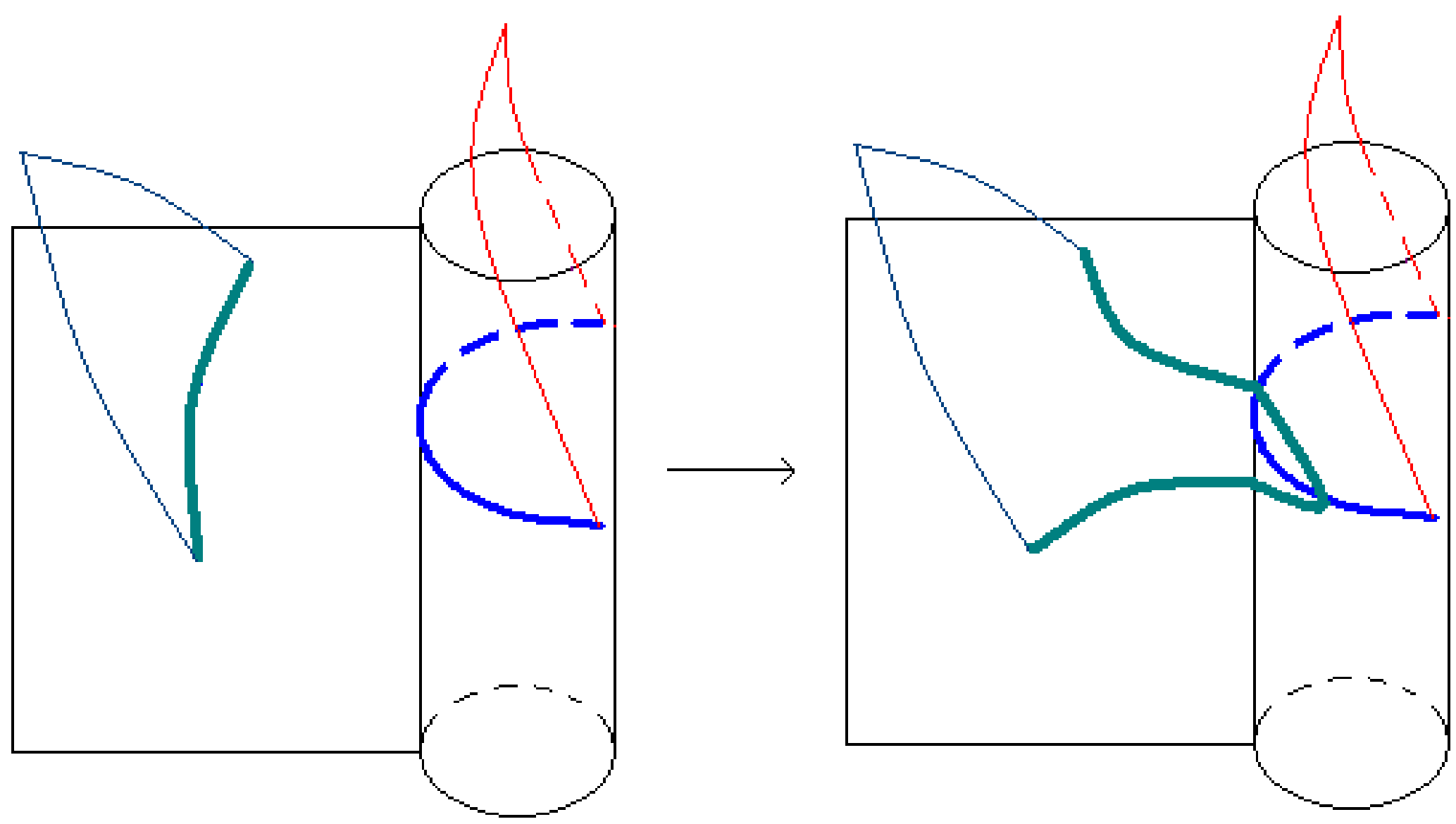

**Figure 215- The conjugation- list of crossing for slided arc- b) curve in clockwise orientation**

The move is a composition of $T_2$ and T*:

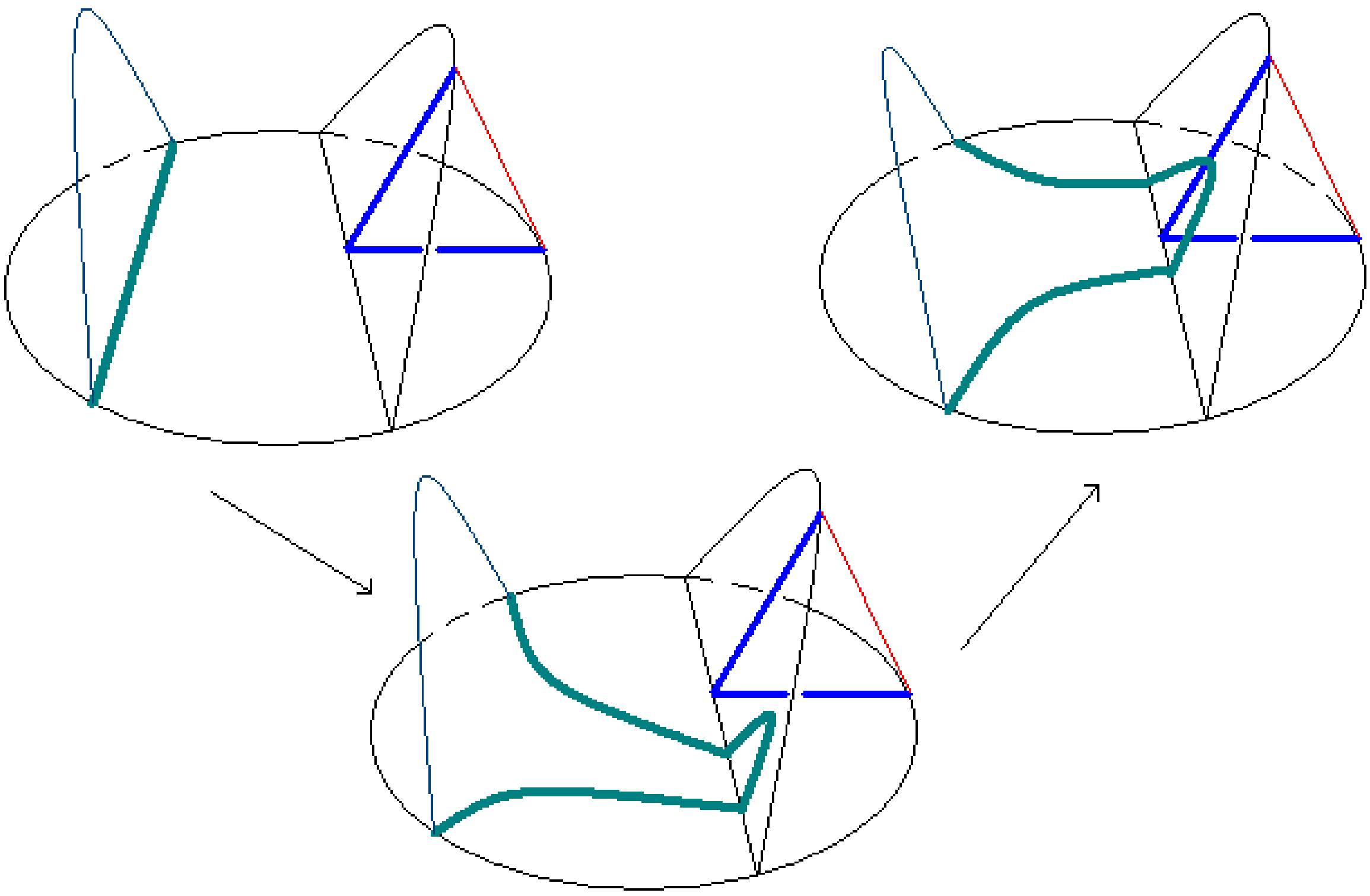

**Figure 216- The conjugation- list of crossing for slided arc- b) curve opposite to clockwise orientation- composition of $T_2$ and T***

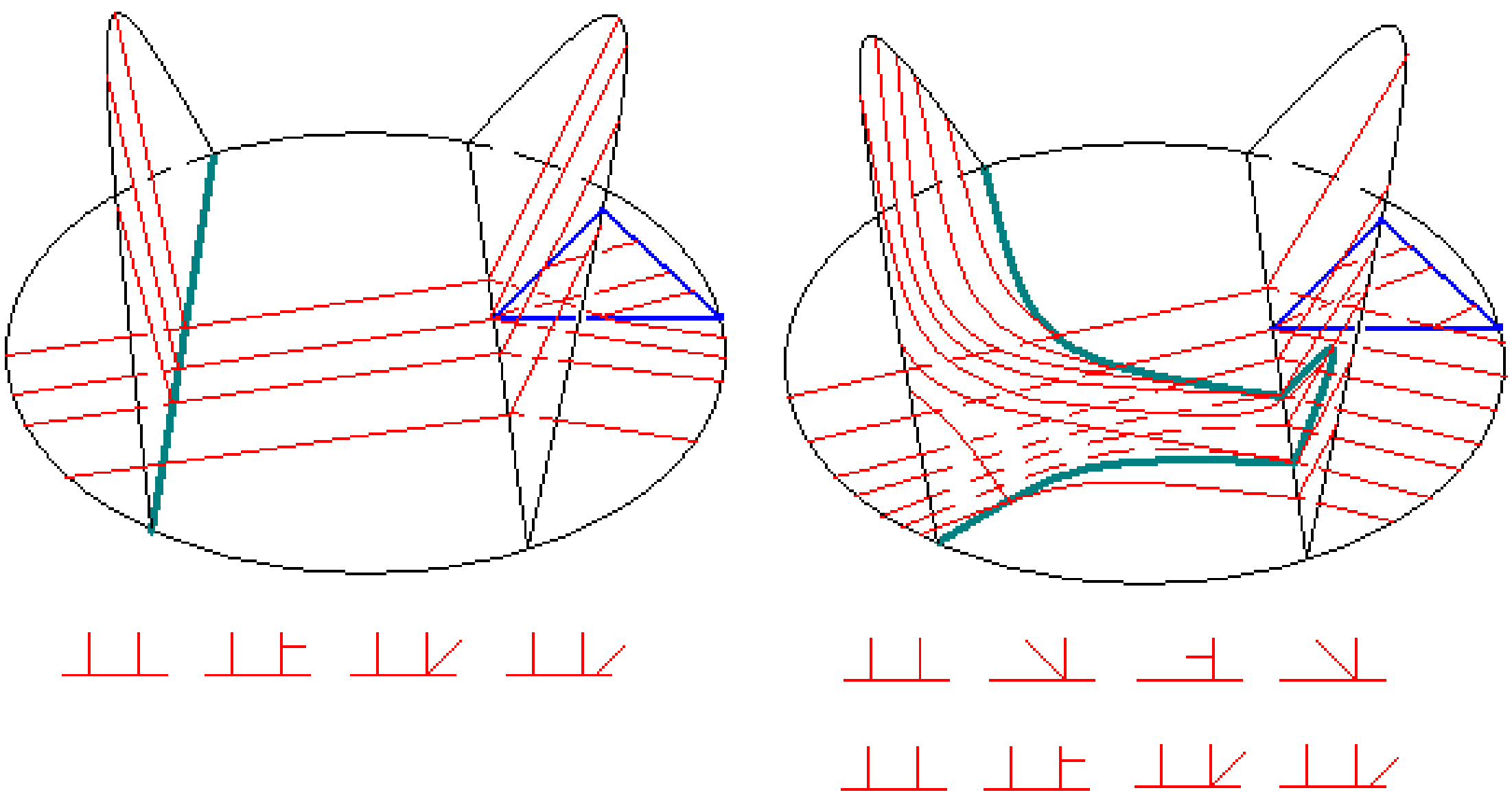

**Figure 217- The conjugation- b) curve in clockwise orientation- composition of $T_2$ and T*- sequence of slices- part 1**

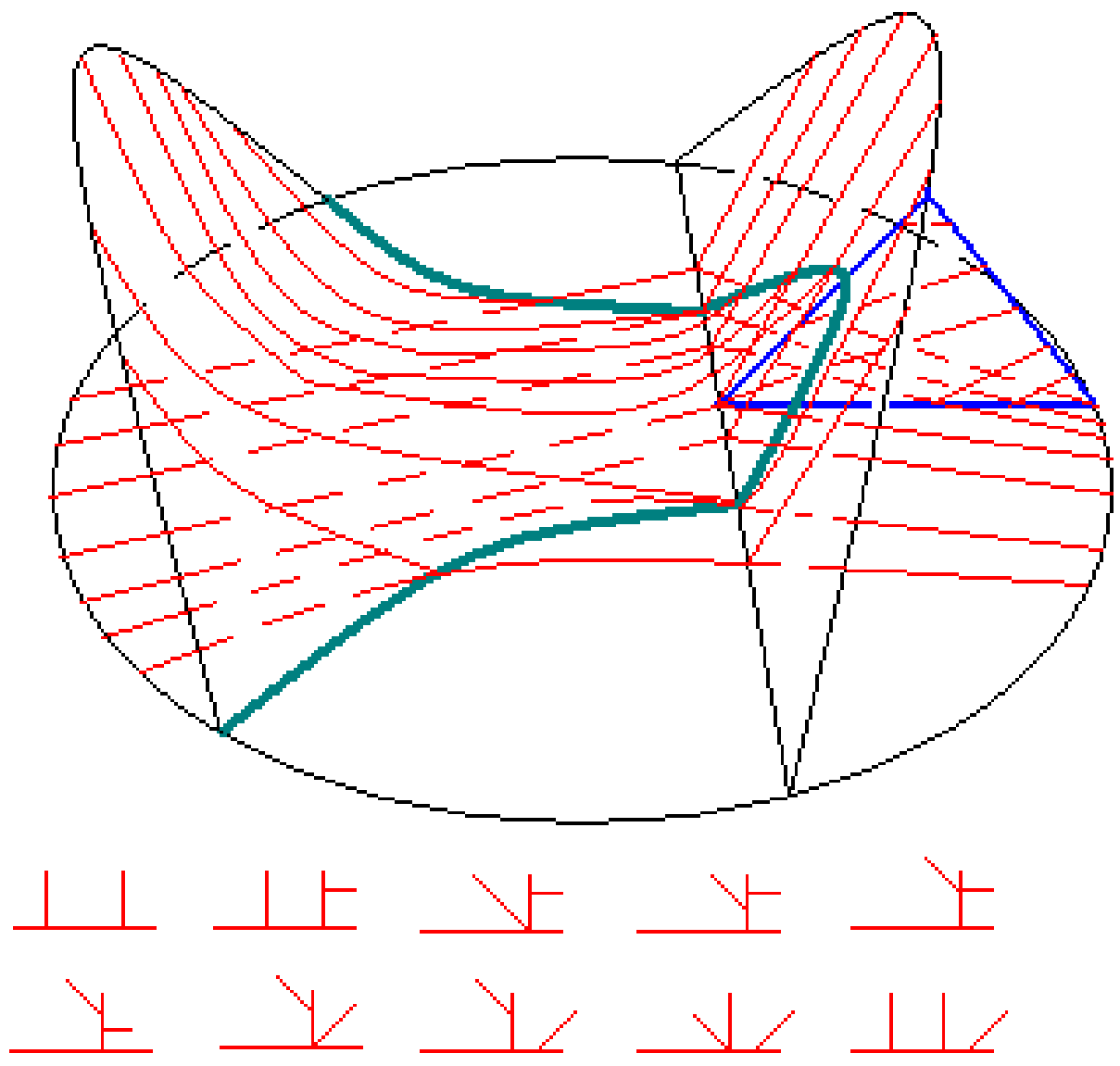

**Figure 218- The conjugation- b) curve in clockwise orientation- composition of $T_2$ and T*- sequence of slices- part 2**

c) The arc passes from the generator cylinder to the rectangle:

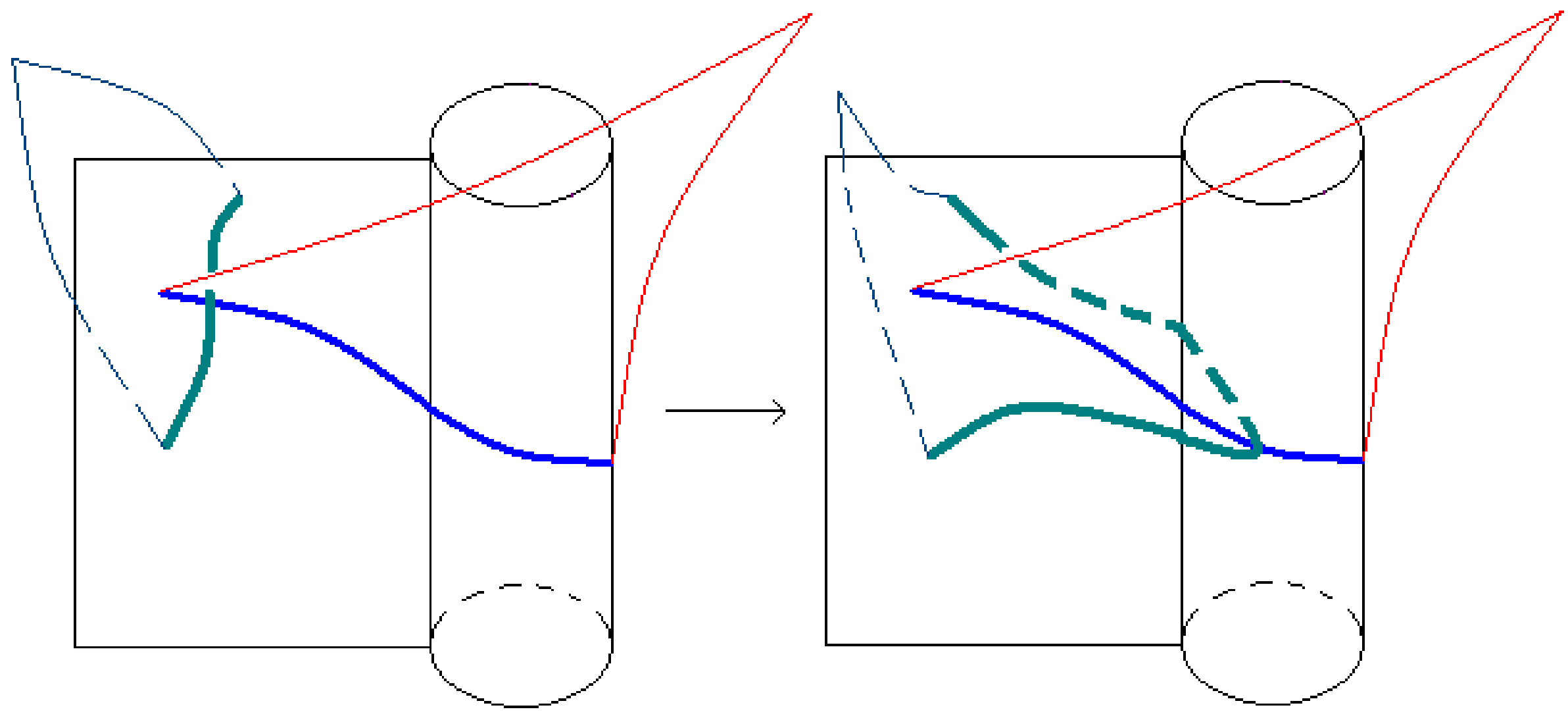

**Figure 219- The conjugation- c) arc pass from generator clinder to rectangle**

We get a modified $T_2$ move:

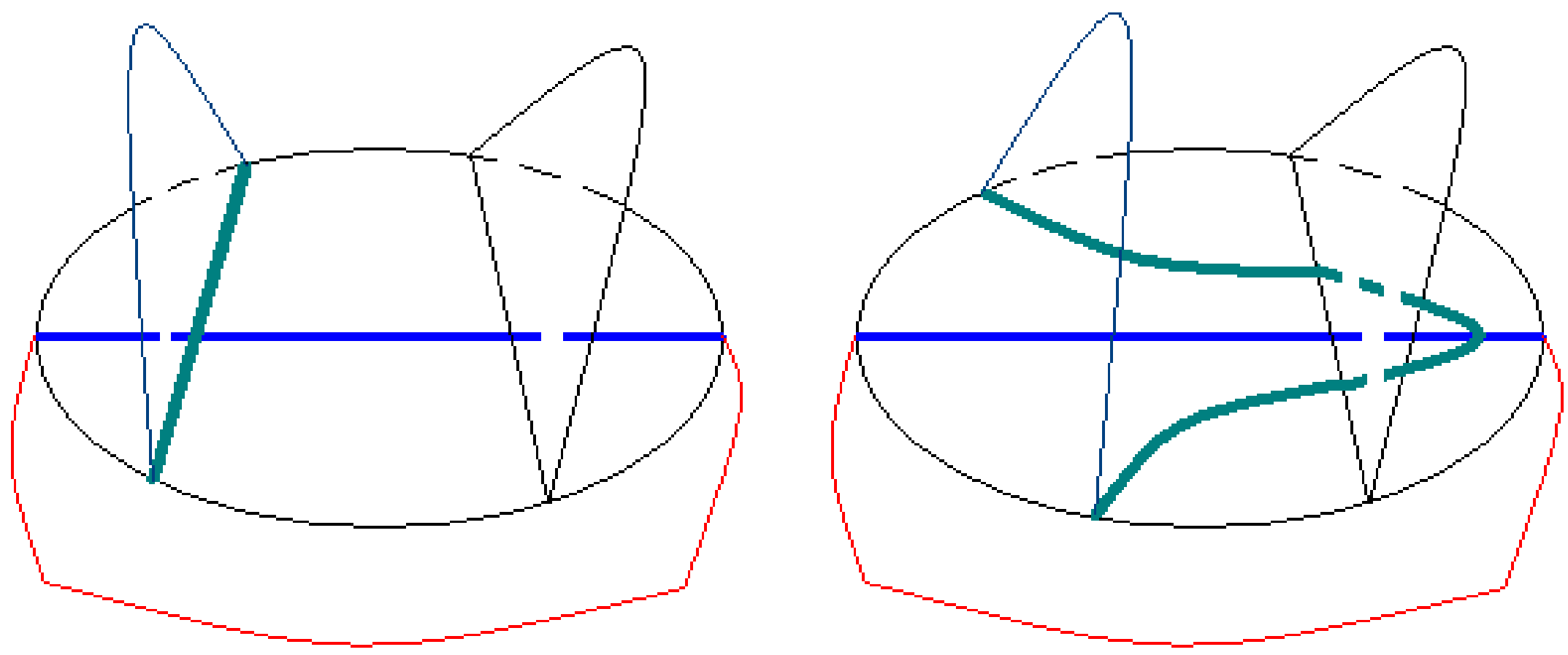

**Figure 220- The conjugation- c) arc pass from generator clinder to rectangle- modified $T_2$ move**

## *7.3 The inverse*

Assume we have a relation R in a 2-complex. We could change to the relation $R^{-1}$ by taking the same attaching map of R with reverse orientation of the 2-disk $D^2$, but this is not useful for our considerations. We read the relation R in the Quinn model from bottom to top:

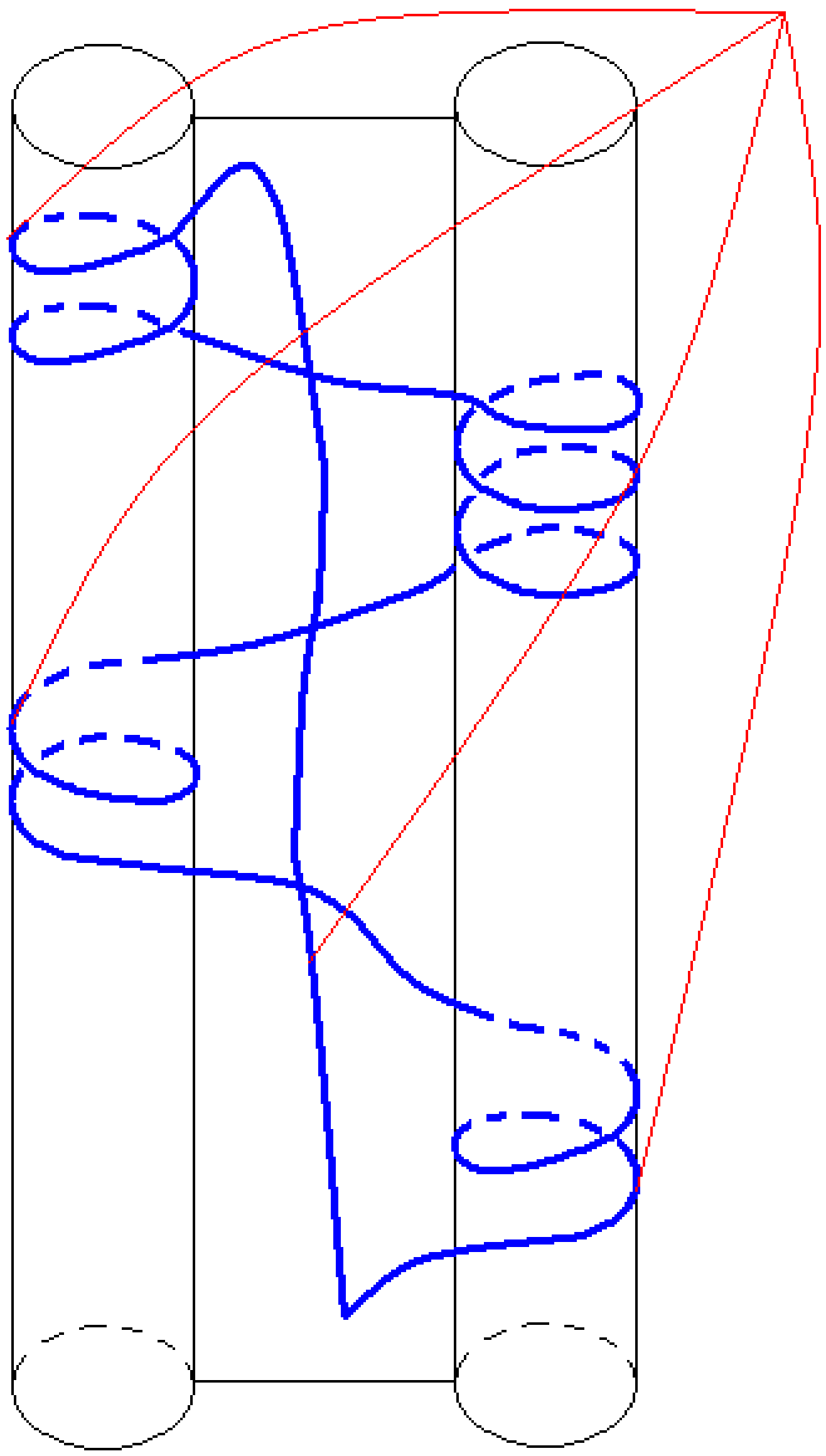

**Figure 221- The inverse- relation in Quinn model read from bottom to top**

If we read from top to bottom, we get $R^{-1}$:
Note, that for each generator in R we would read the inverse of the generator, hence we construct our homotopy by dropping down the relation R from top to bottom according to the appearance of the generators in the attaching map. The self crossings of the attaching map during these processes lead to Matveev moves again. The next sequence of  pictures indicates the homotopy for the relation given above:

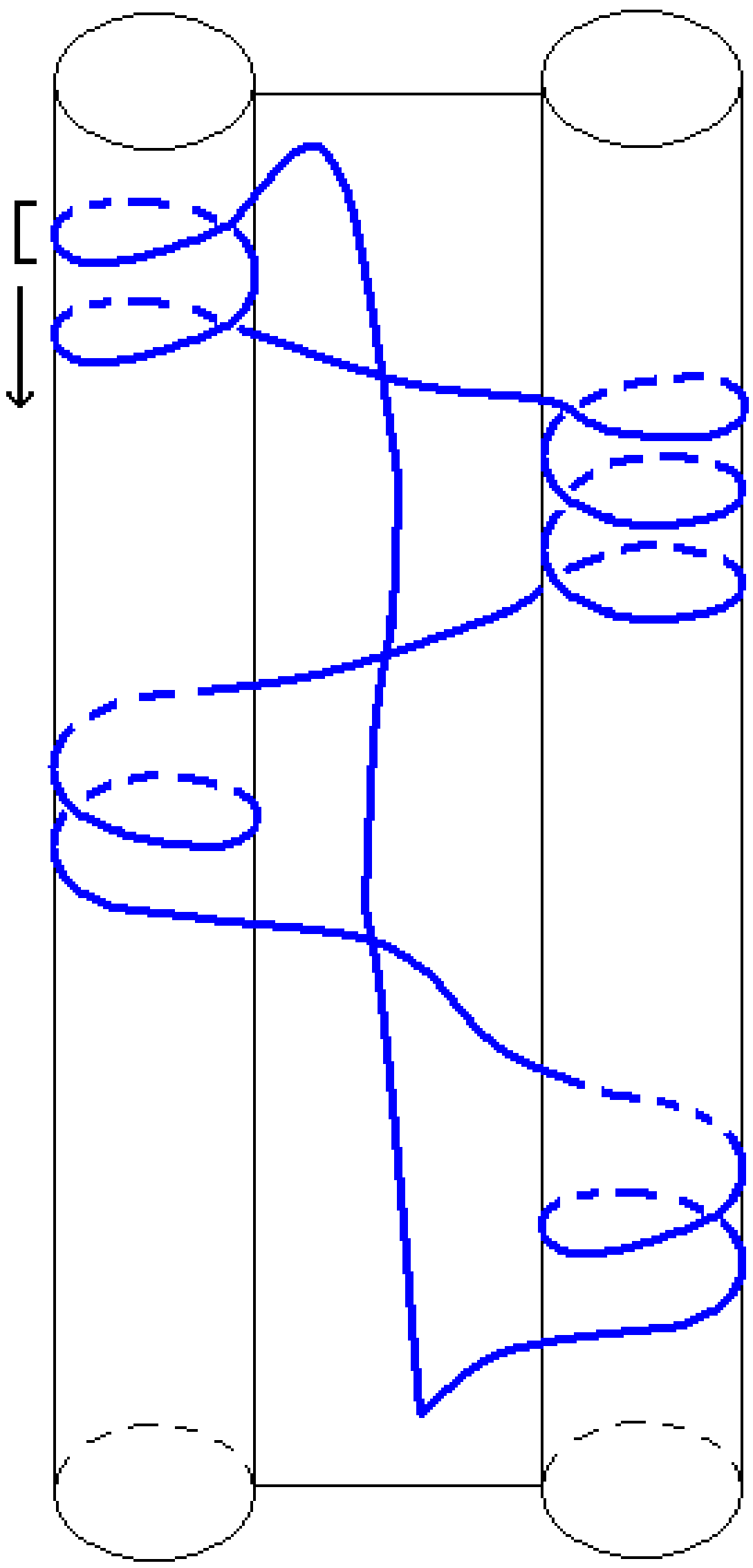

**Figure 222- The inverse- transfer relation to the inverse relation- 1**

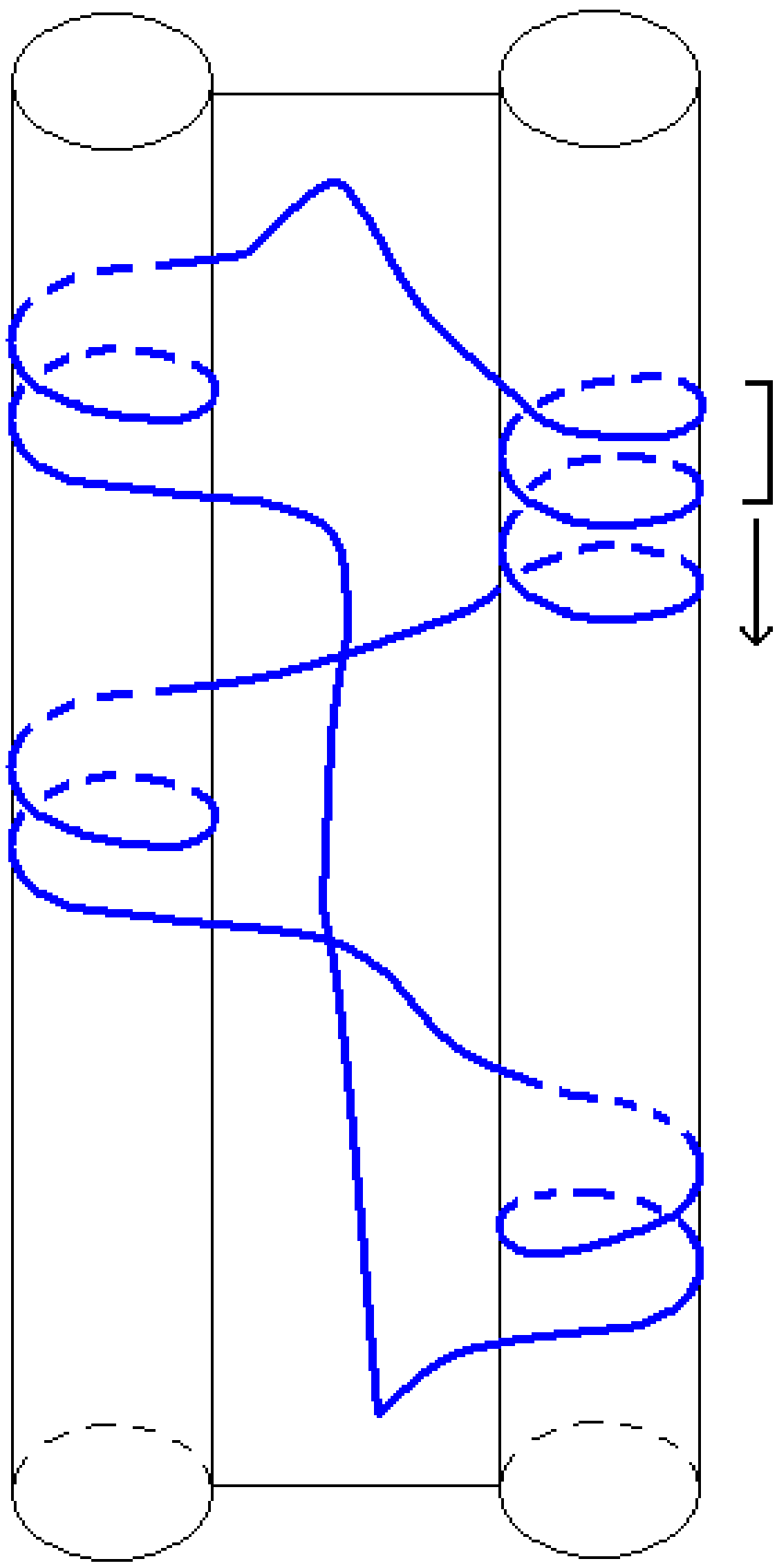

**Figure 223- The inverse- transfer relation to the inverse relation- 2**

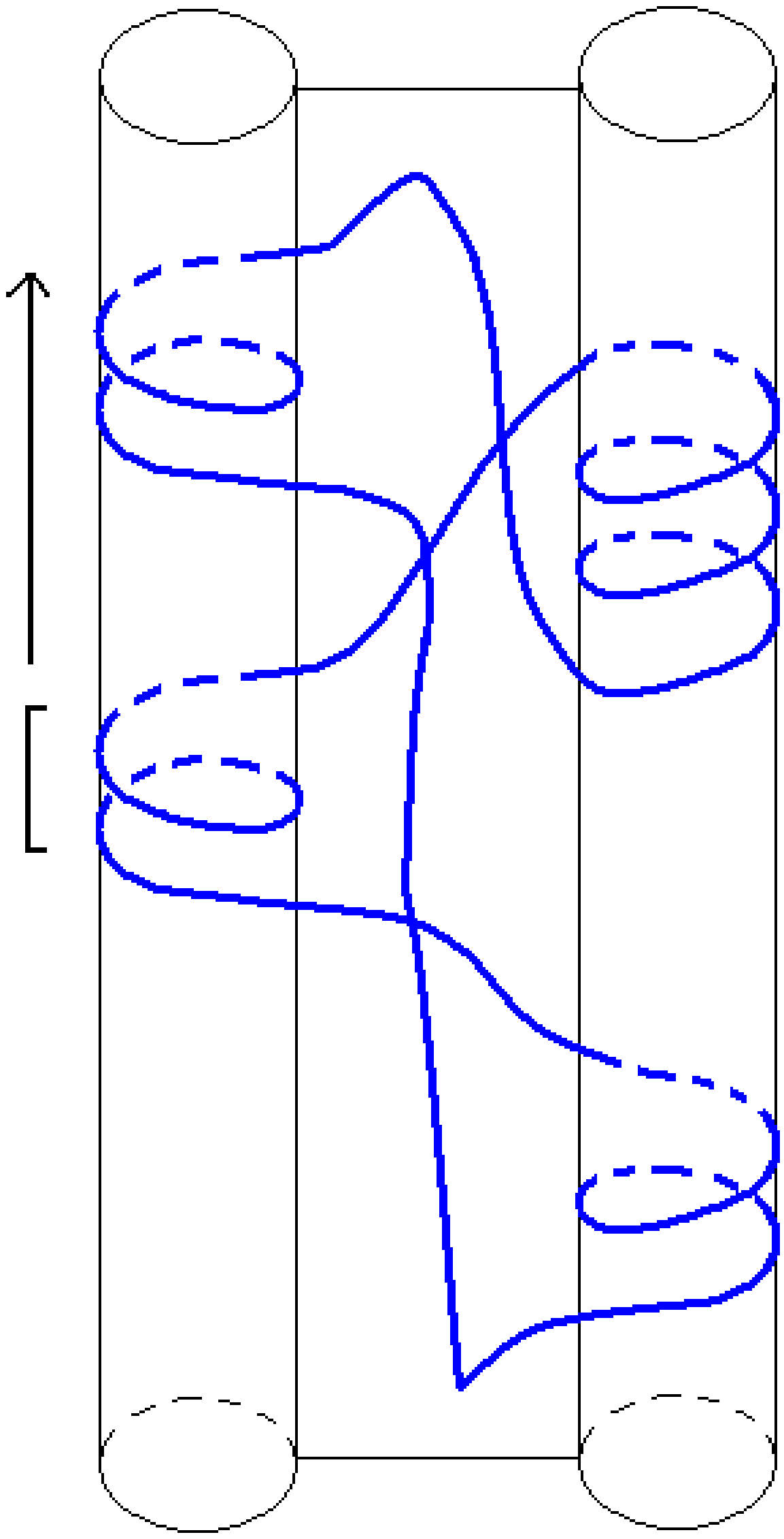

**Figure 224- The inverse- transfer relation to the inverse relation- 3**

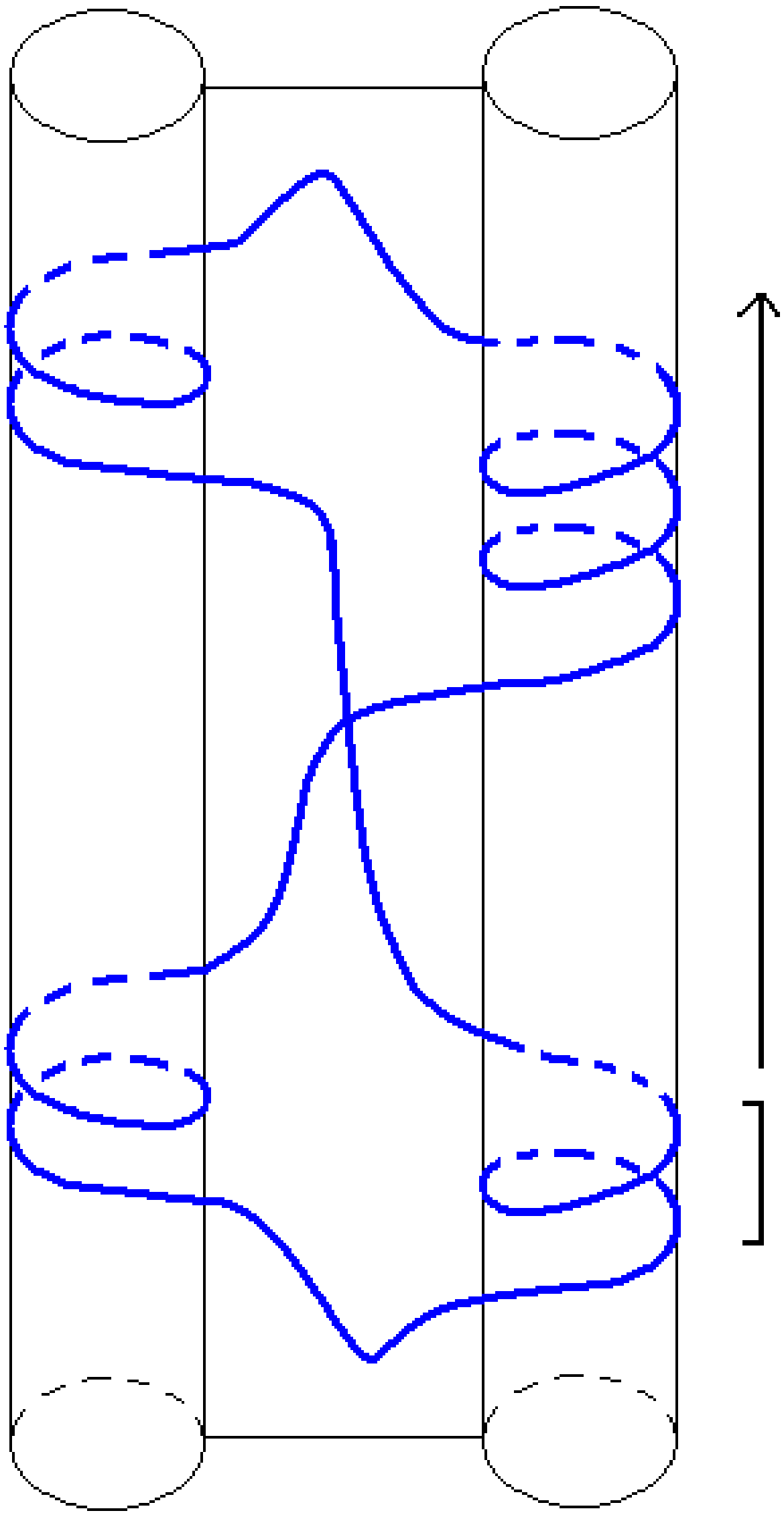

**Figure 225- The inverse- transfer relation to the inverse relation- 4**

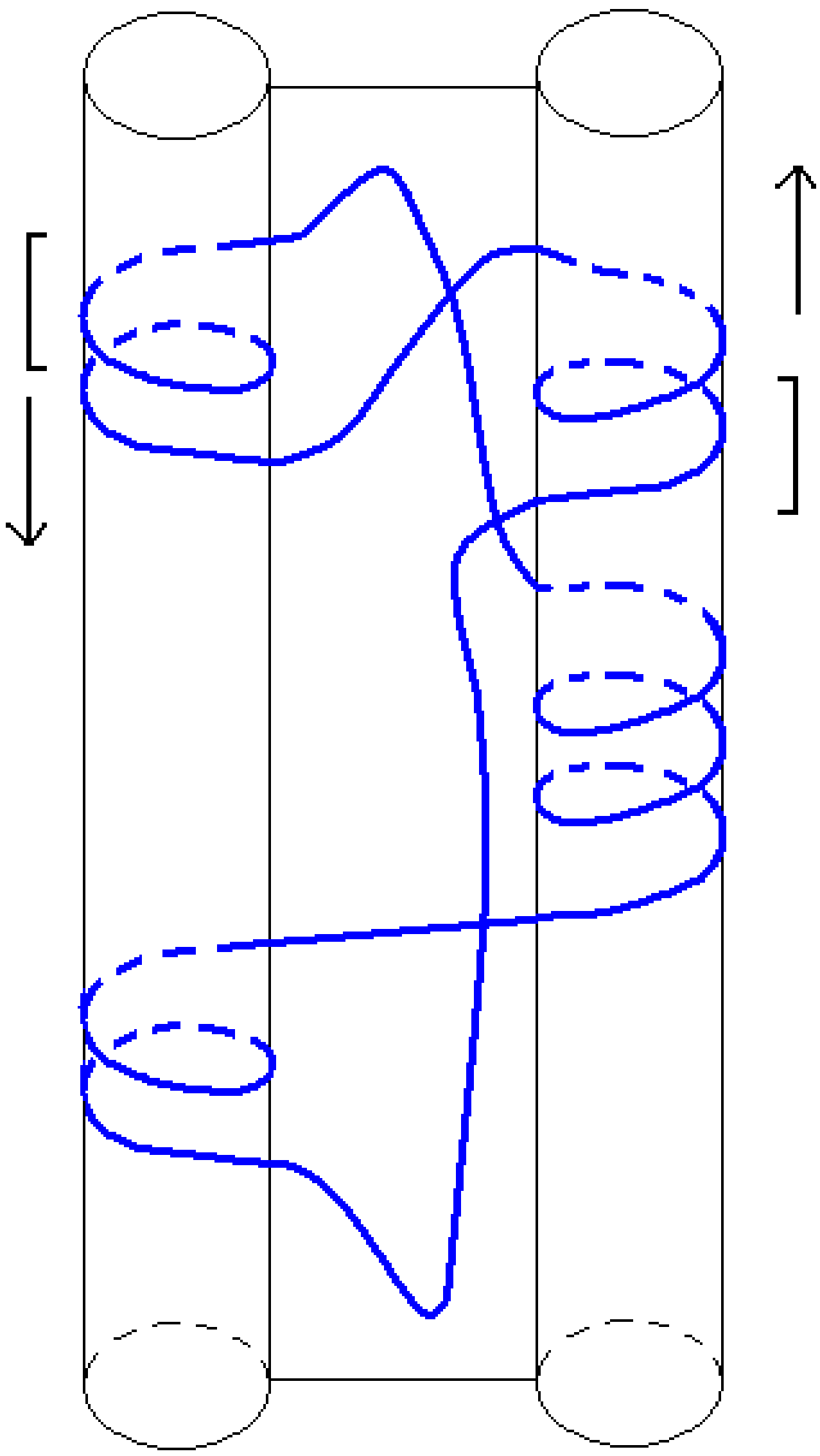

**Figure 226- The inverse- transfer relation to the inverse relation- 5**

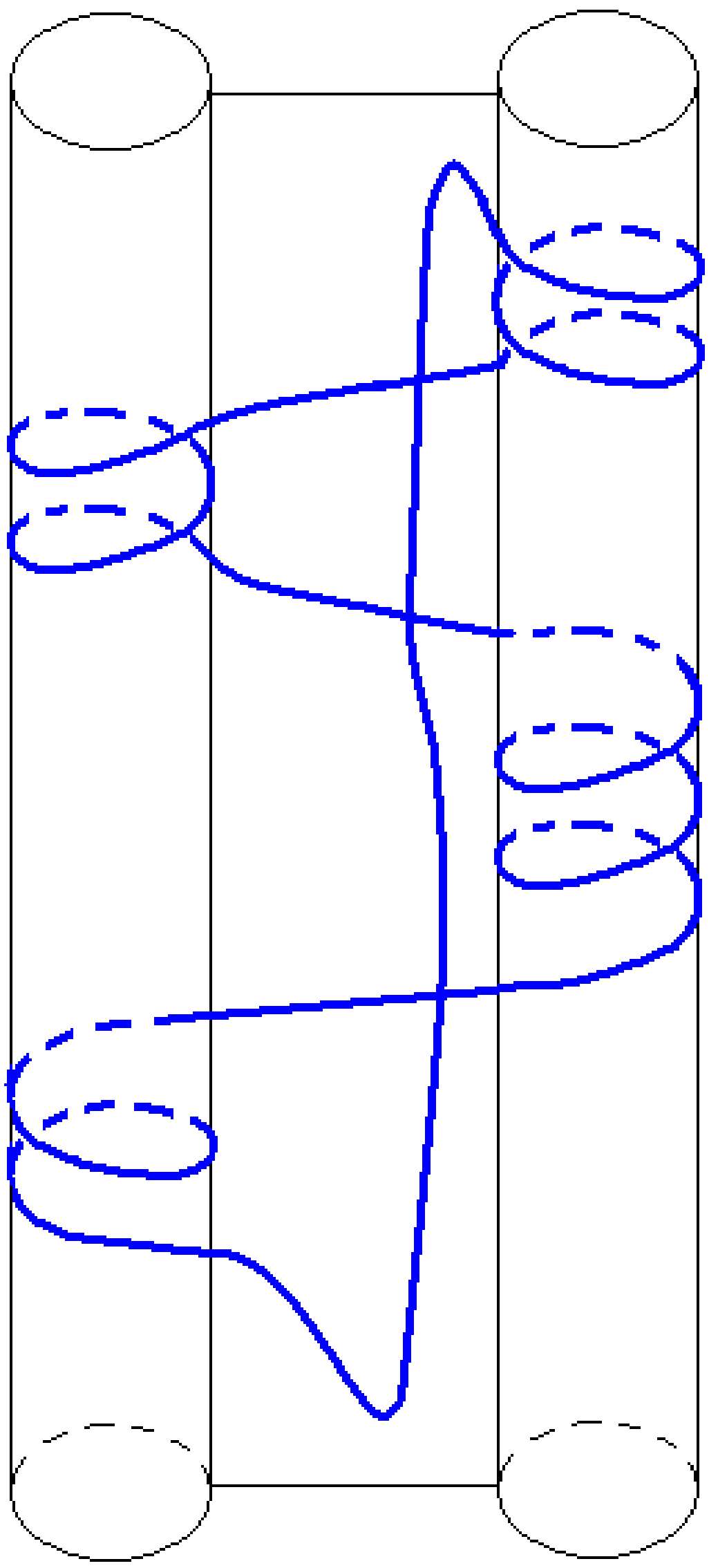

**Figure 227- The inverse- transfer relation to the inverse relation- 6**

We pick some steps and illustrate, how to compose these as a sequence of Matveev moves.
First we present the basic example. The whole process is a multiple application of it:

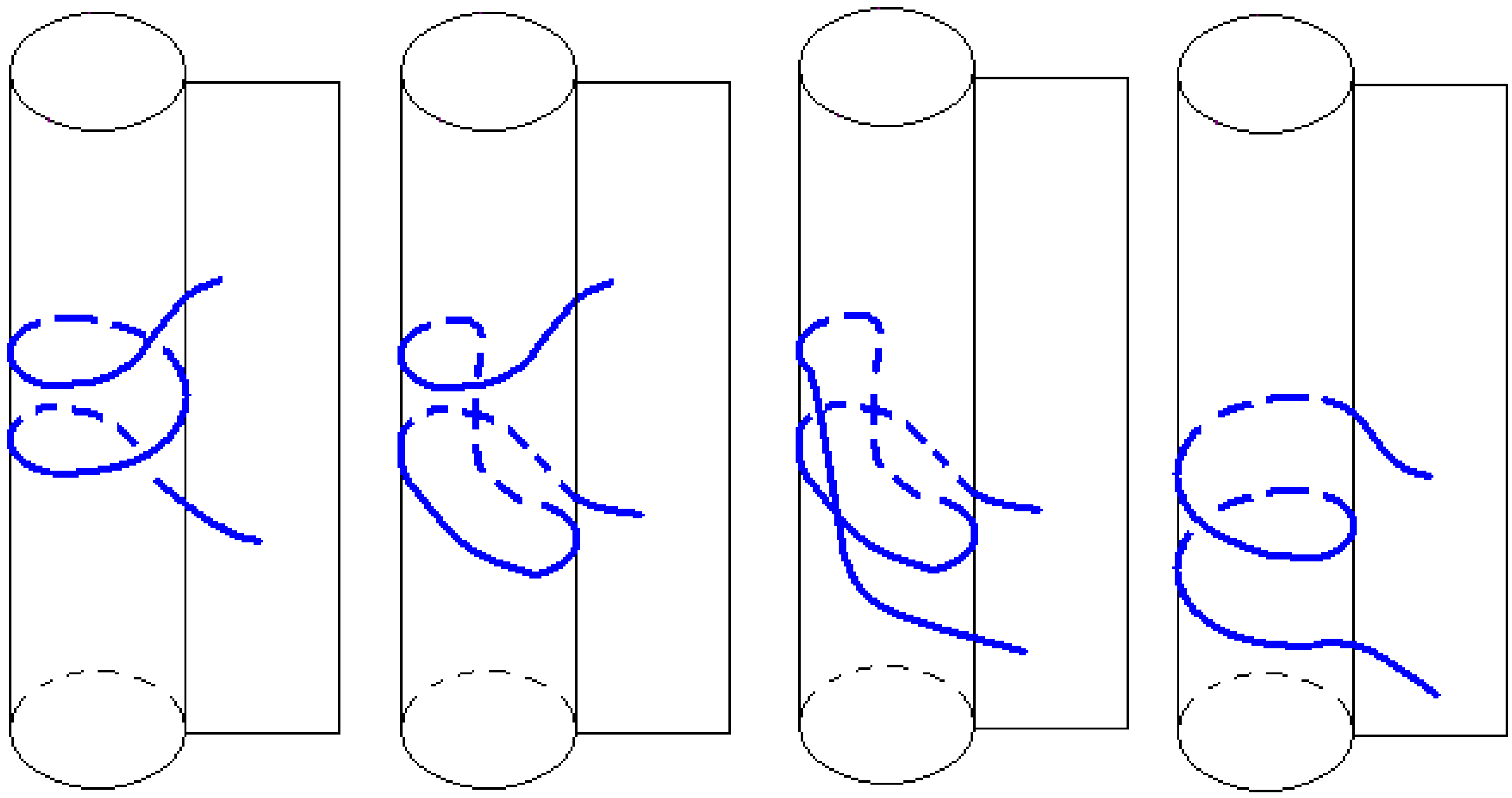

**Figure 228- The inverse- the basic example**

It is sufficient to study this example in detail, but before that, we consider a more complicated example:
We want to exchange the blue curve with the green curve.

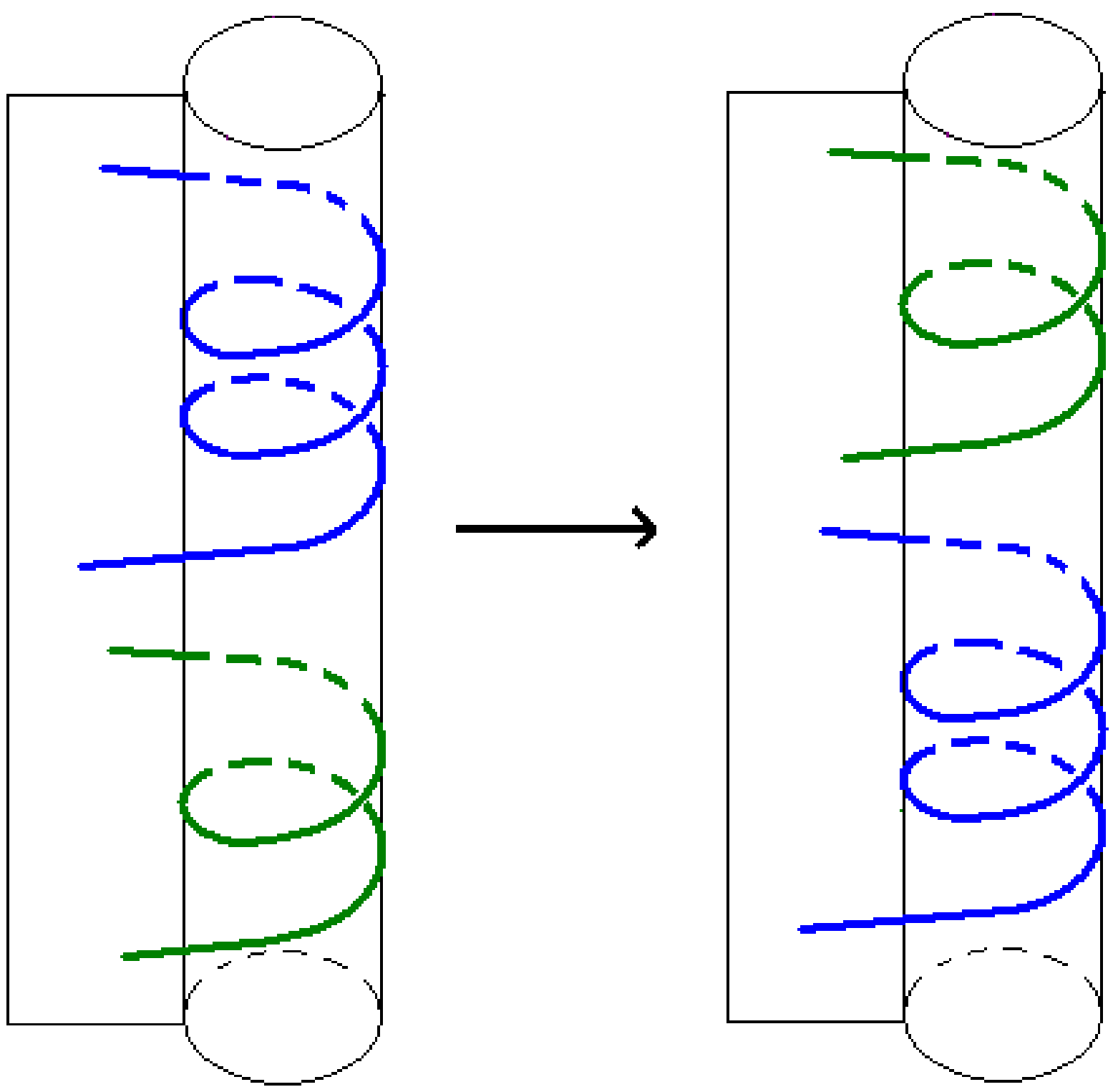

**Figure 229- The inverse- preview- a more complicated example**

One way to do that is shown in the following sequence:

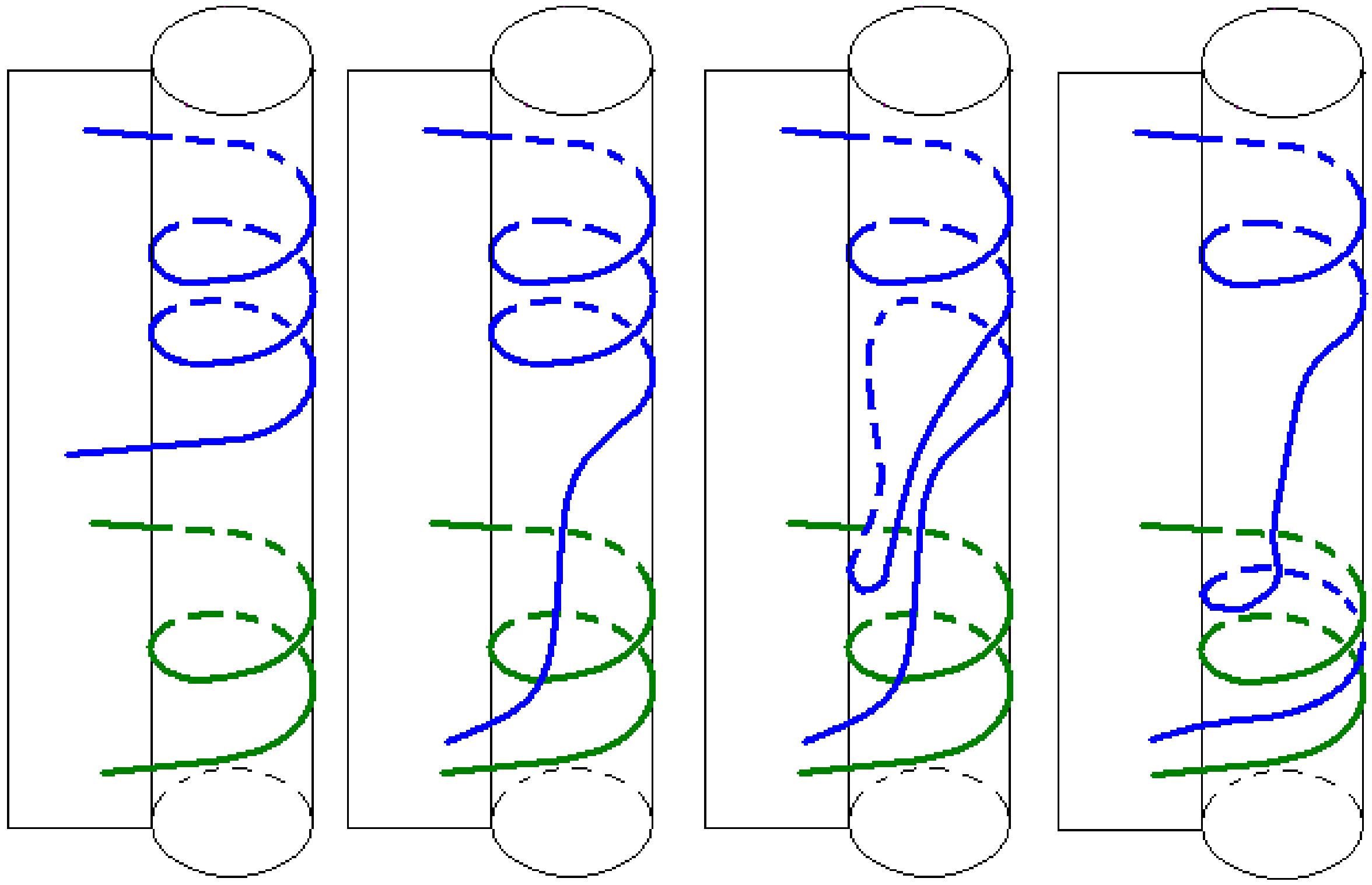

**Figure 230- The inverse- a more complicated example- part 1**

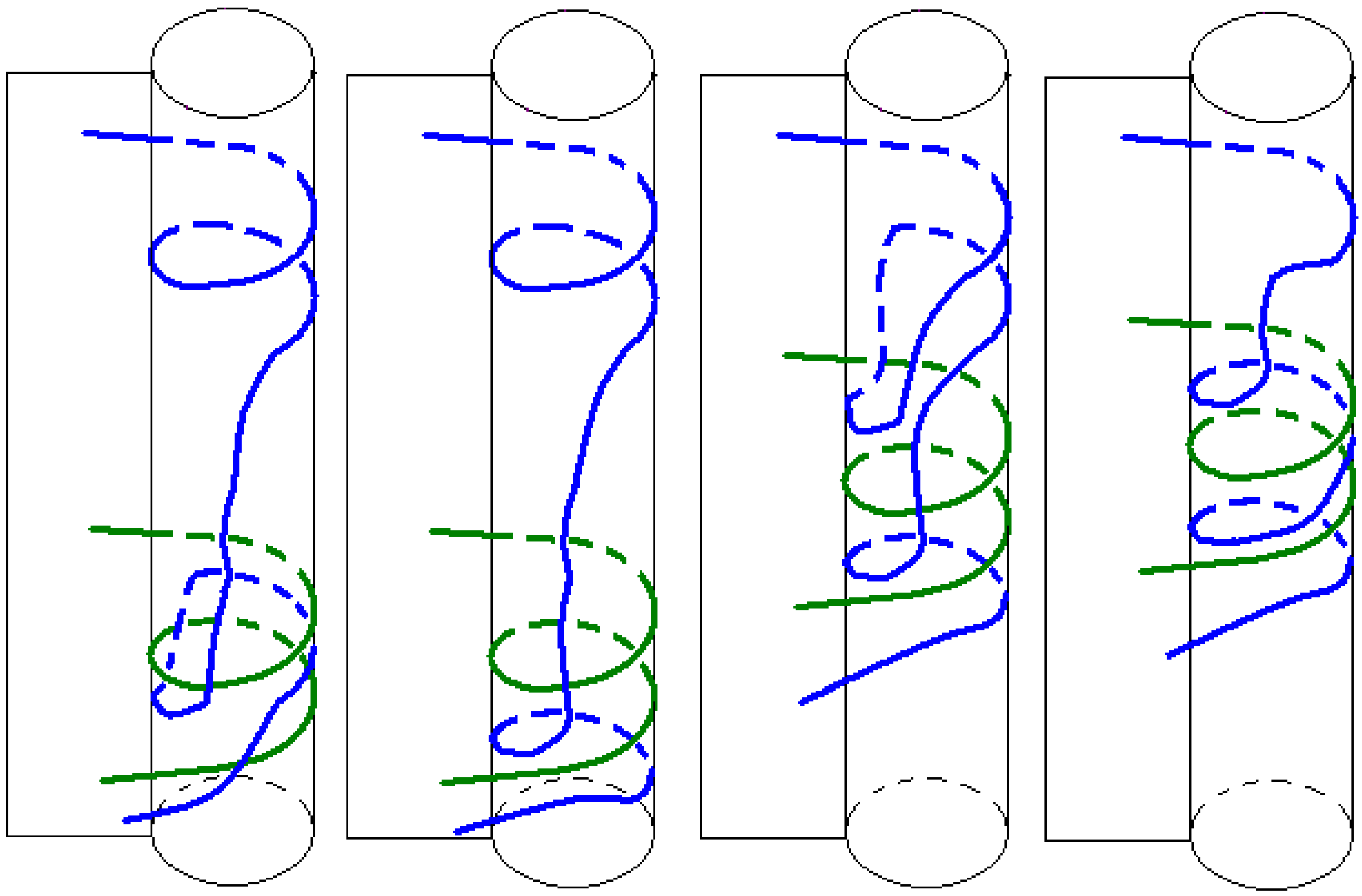

**Figure 231- The inverse- a more complicated example- part 2**

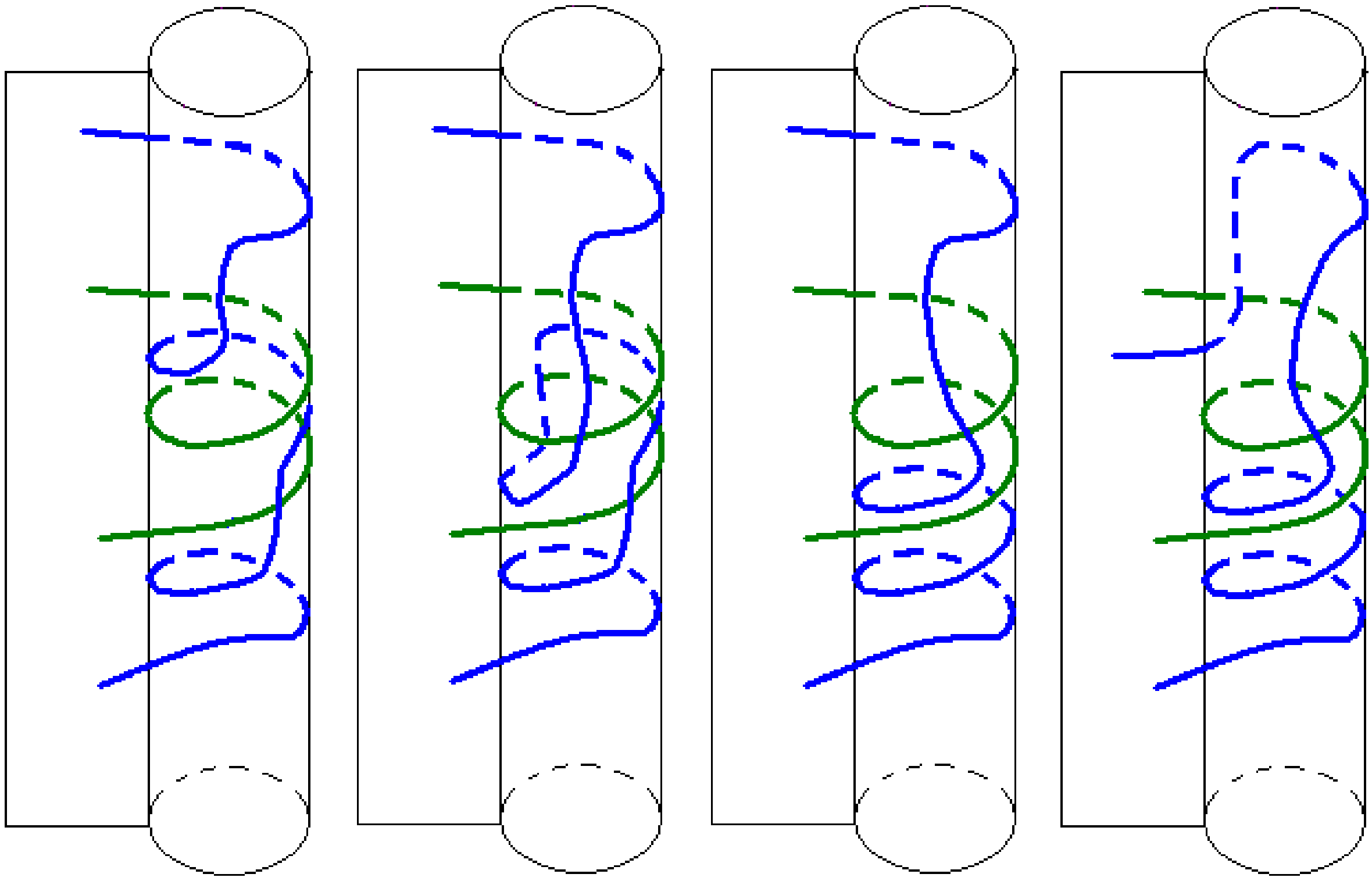

Figure 232- The inverse- a more complicated example- part 3

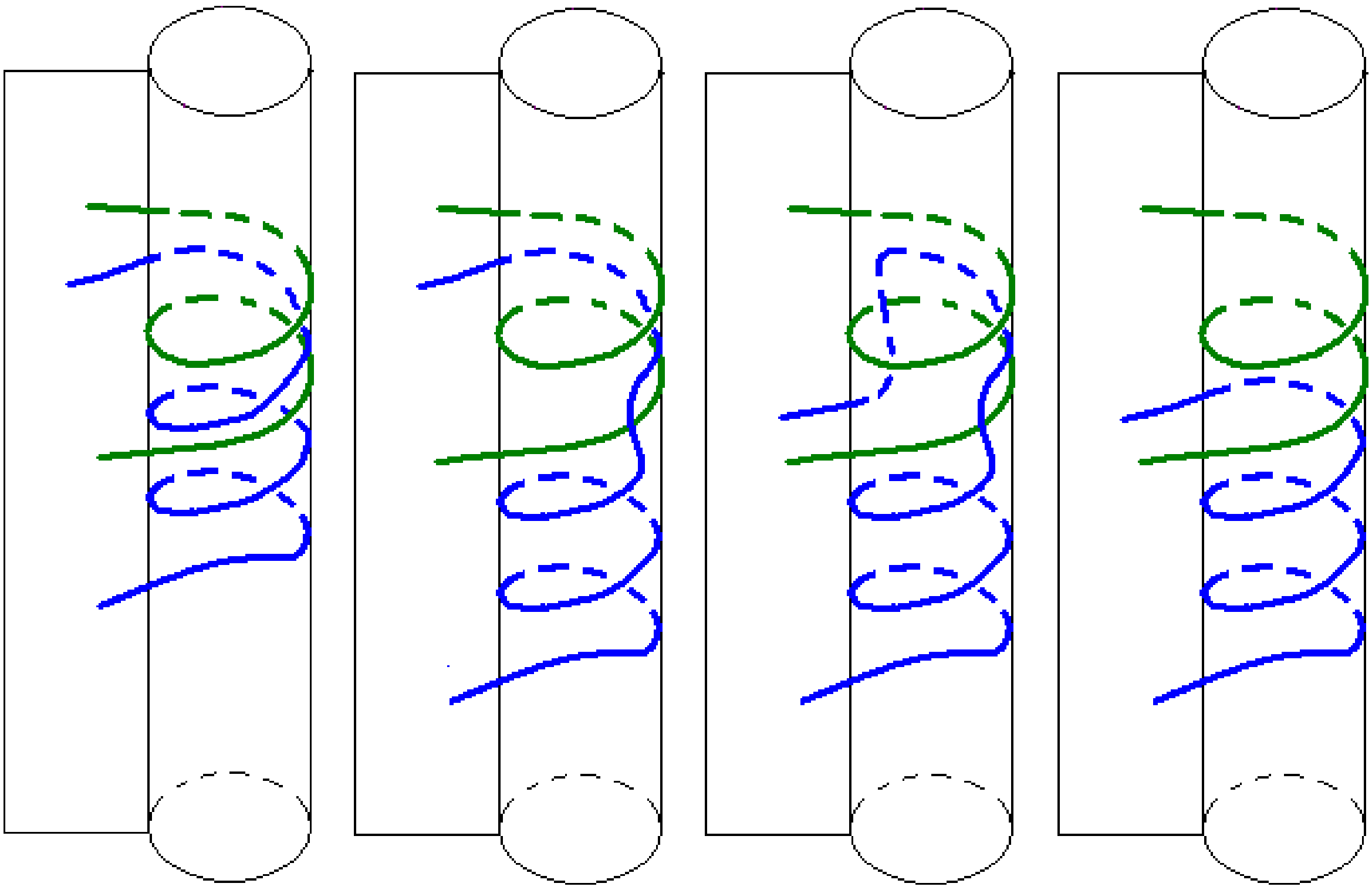

Figure 233- The inverse- a more complicated example- part 4 (end)

We will concentrate on the basic example. First it should be pointed out, that local extremas arise:

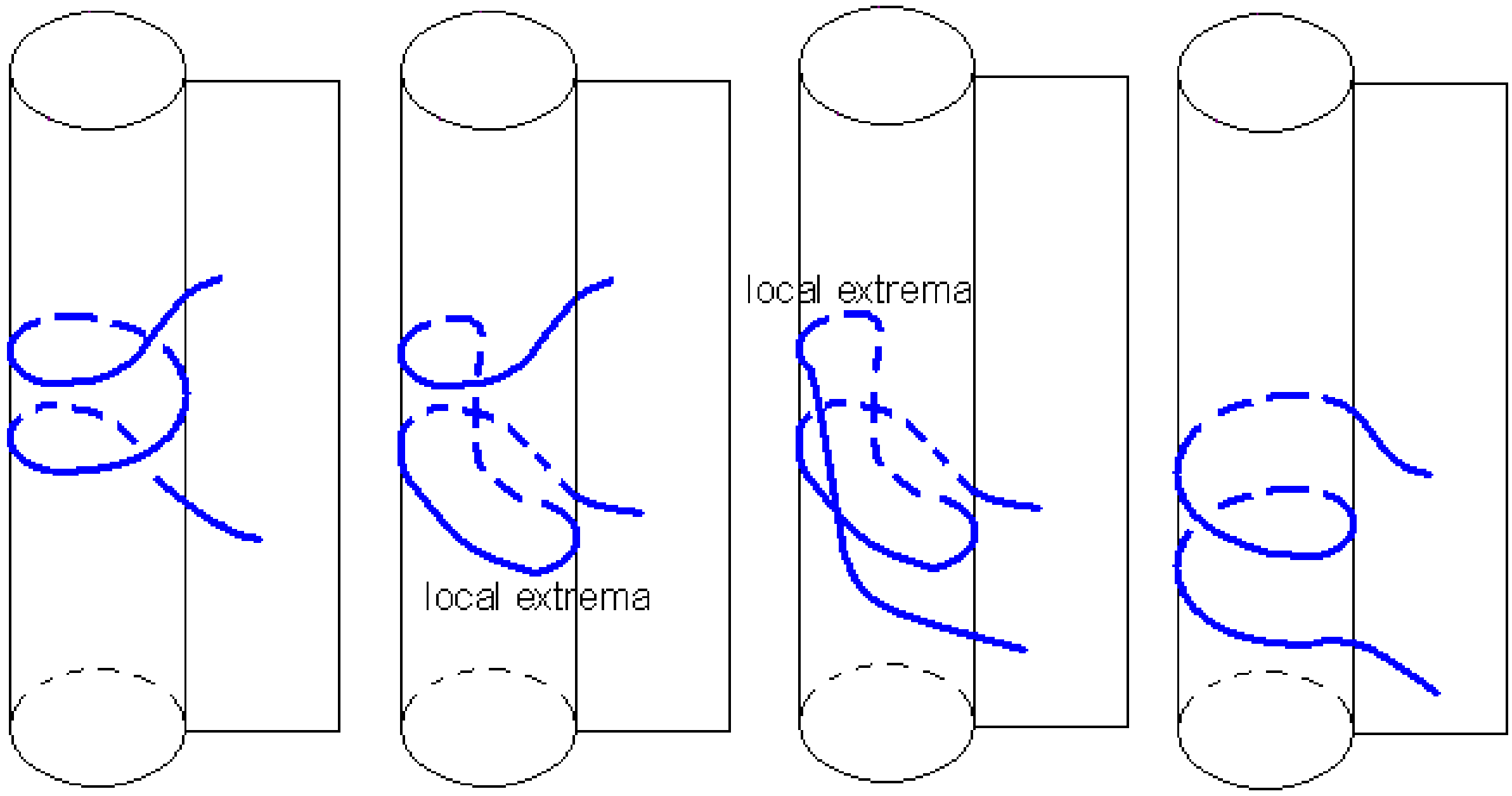

**Figure 234- The inverse- basic example with viewpoint on the arising local extrema**

We divide the process into two parts, the “top” part which describes the changes at the top loop and similar the “bottom” part. Note that we have to prepare the sequence of slices for the appearence of the extrema. It is better to do this separately and then compose the pieces together.

i. bottom part:

To prepare the slices for extrema, we introduce a pair of saddlepoints:

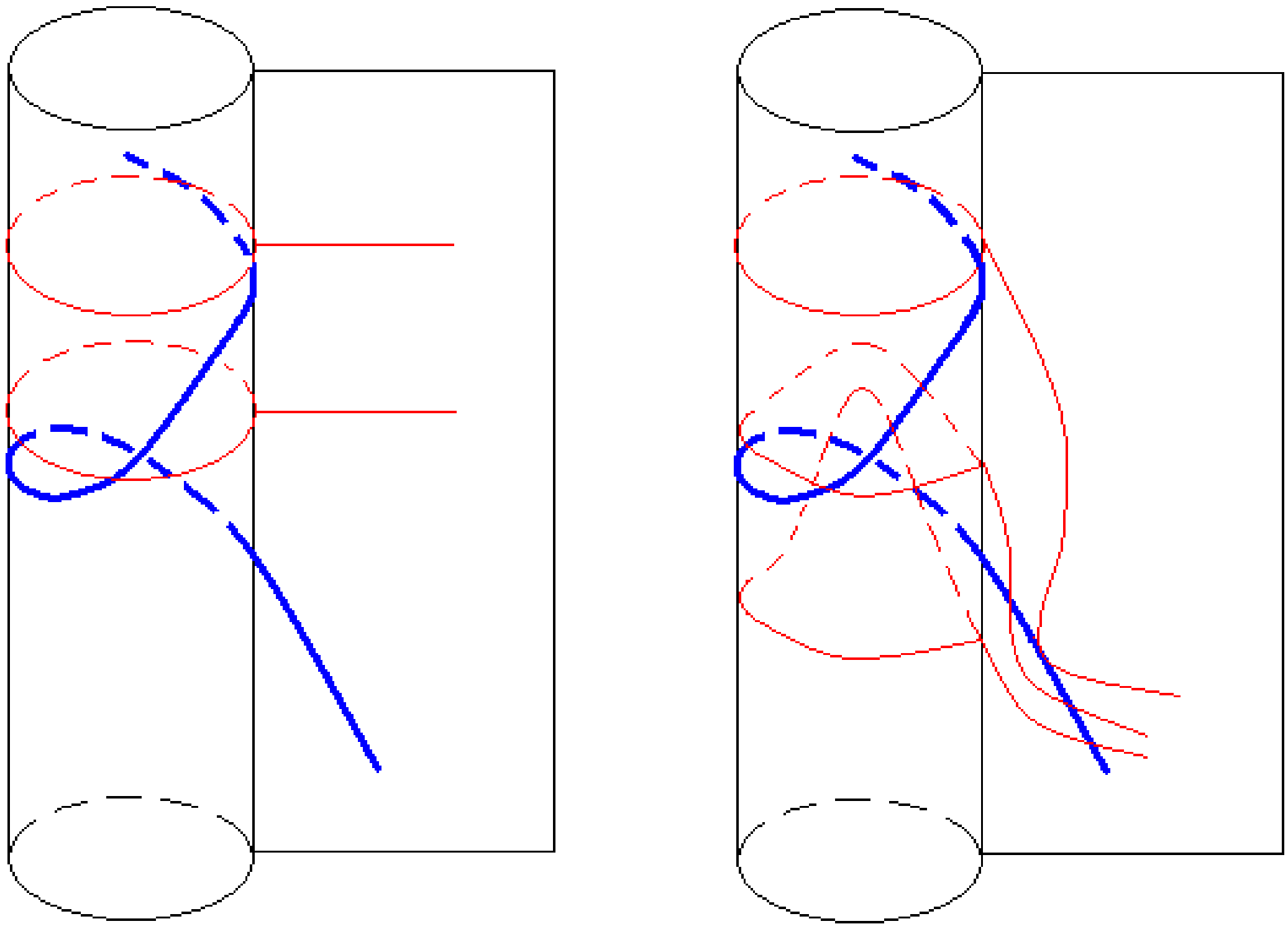

**Figure 235- The inverse- basic example- bottom part- prepare extrema**

Then we drop down the arc and get the extrema:

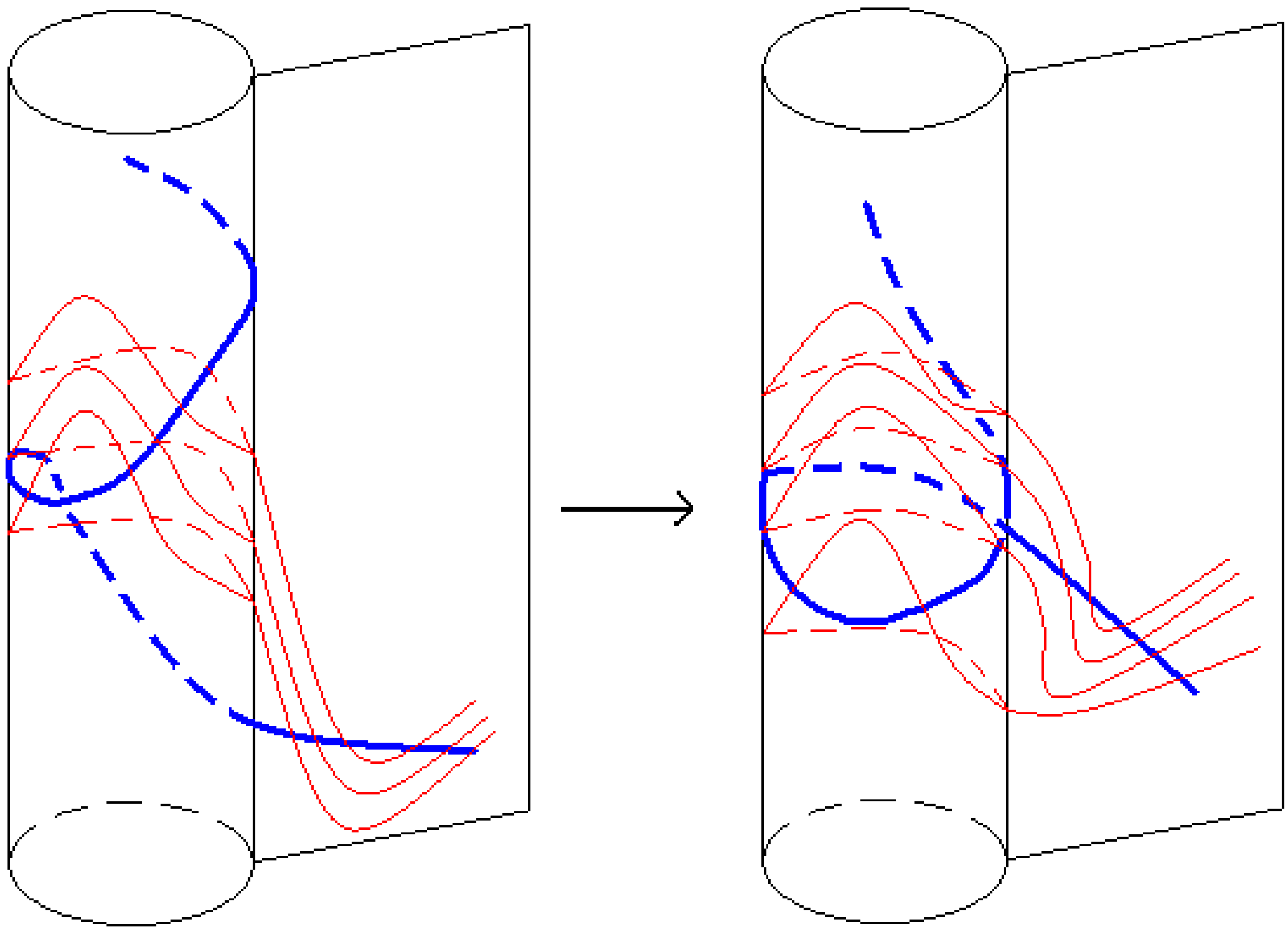

**Figure 236- The inverse- basic example- bottom part- perform step 1**

ii. top part:

Introduce a pair of saddlepoints to prepare the extrema:

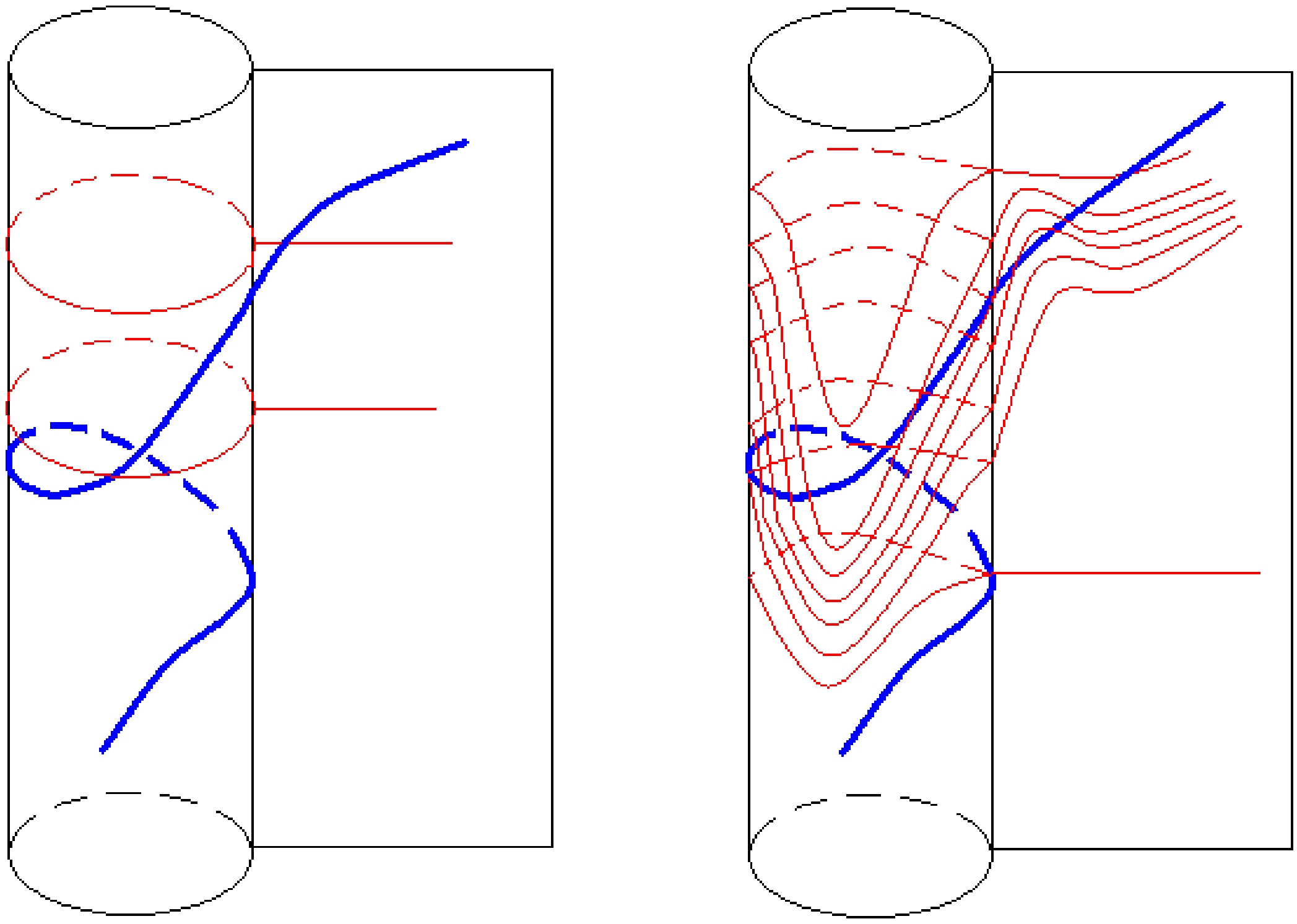

**Figure 237- The inverse- basic example- top part- prepare extrema**

Drop down the arc:

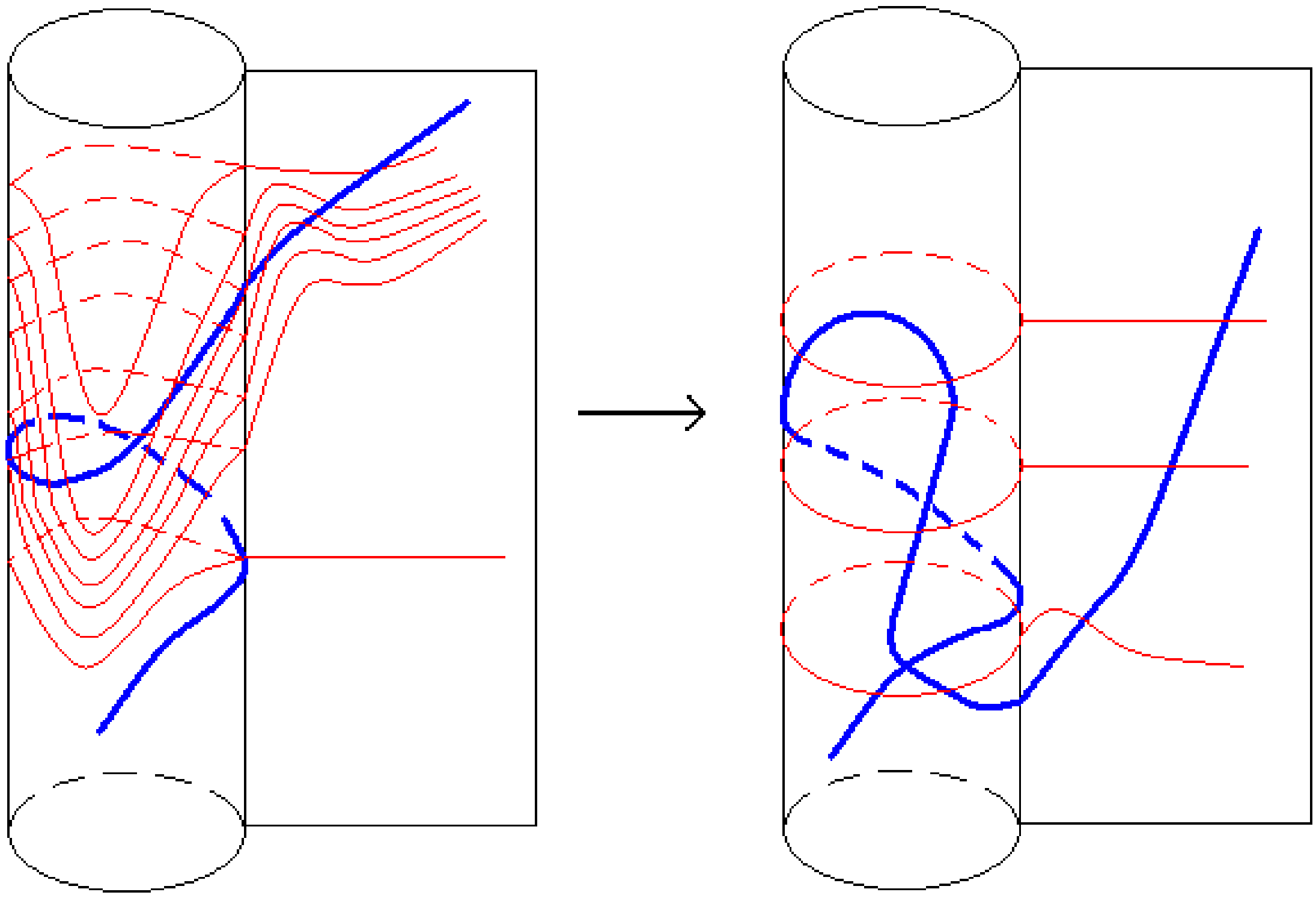

**Figure 238- The inverse- basic example- top part- perform step 2**

We compose the resulting pieces to get a complete picture:

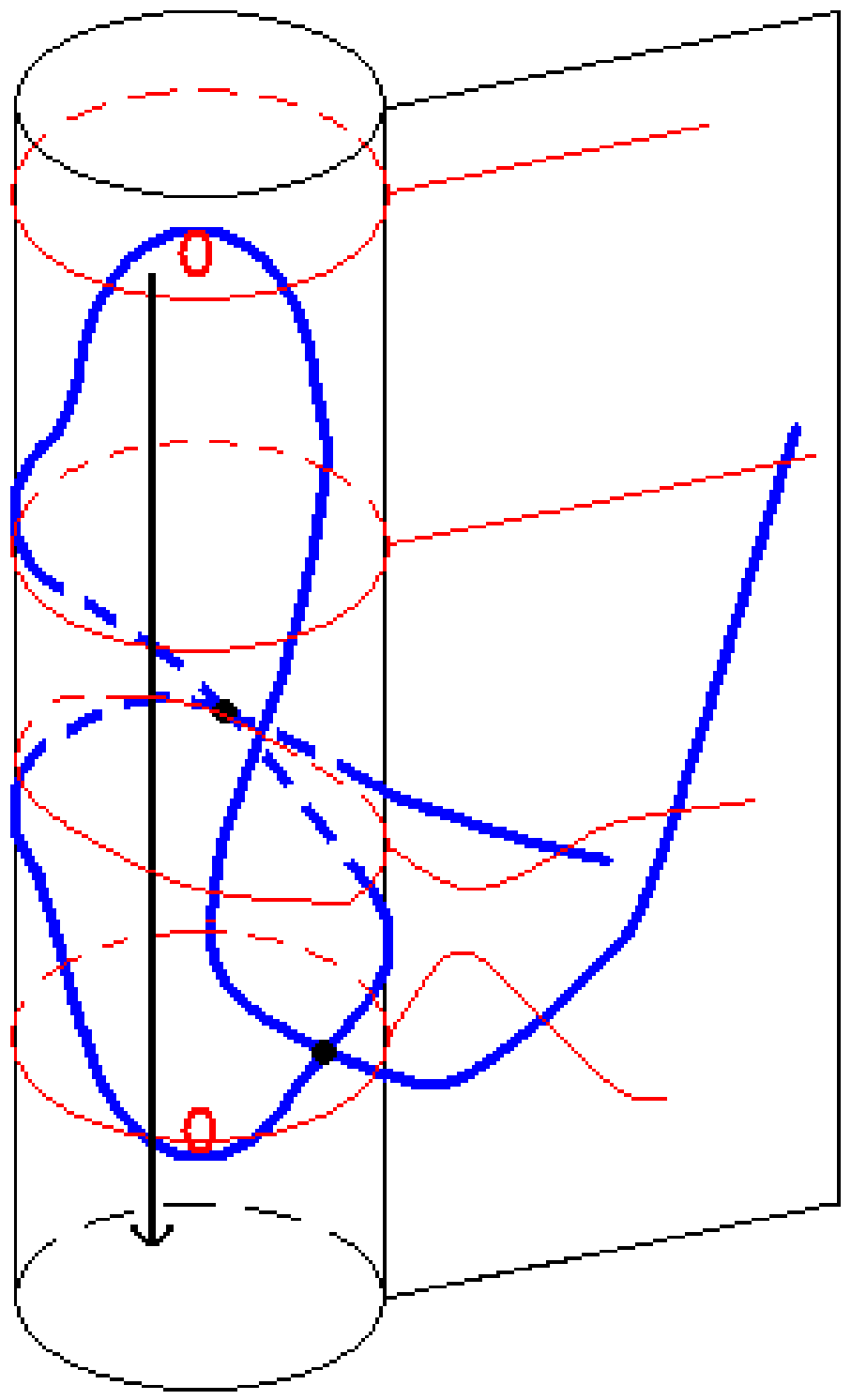

**Figure 239- The inverse- basic example- summmarize top and bottom part**

By performing the move as indicated by the arrow we get the right picture:

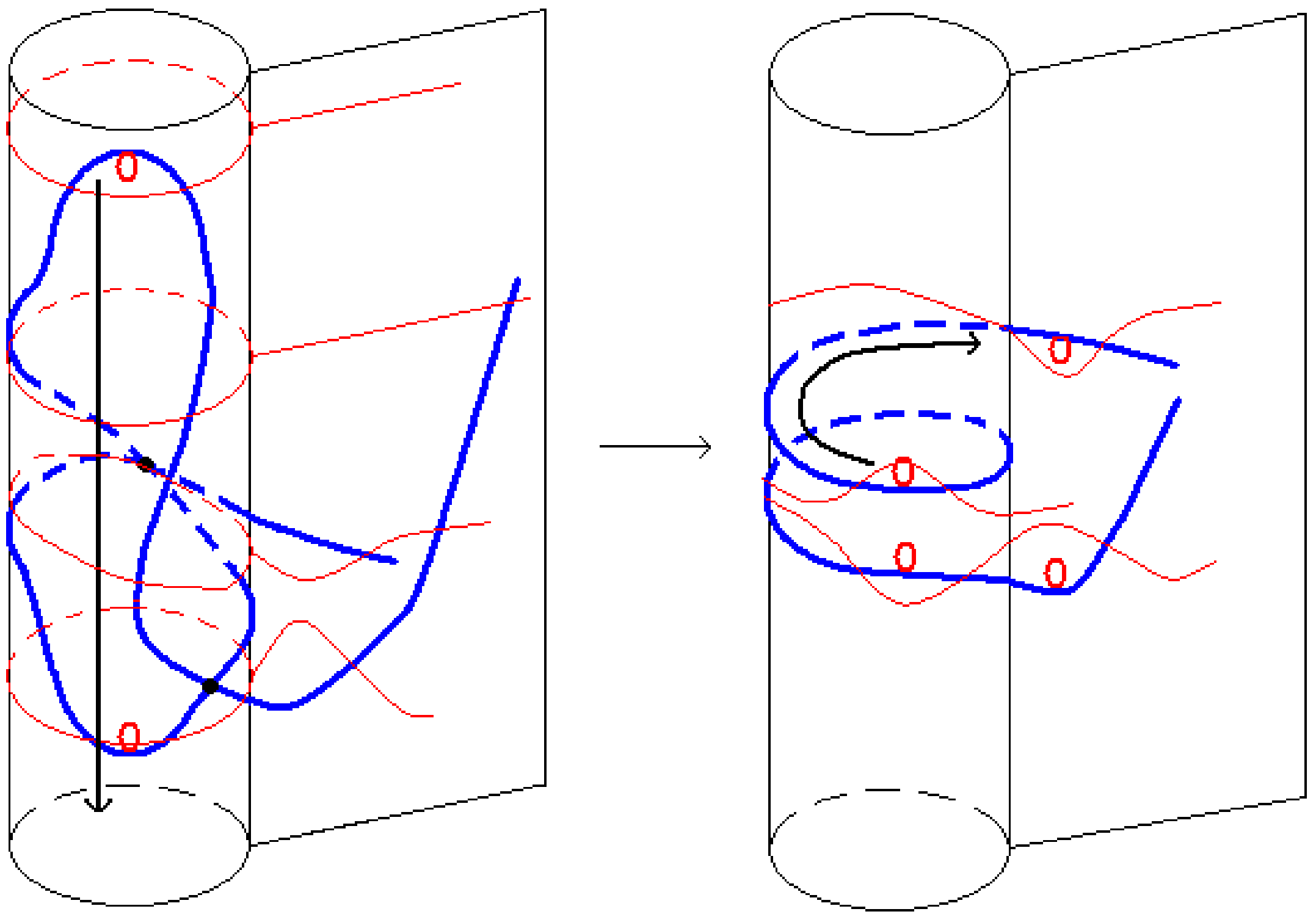

**Figure 240- The inverse- basic example- perform the move on the whole figure**

We draw the extrema and indicate the slices. Furthermore it shows, that the two pairs of saddlepoints can be annihilated, so no extrema appears after doing that step. Note, that the blue curve is strictly increasing on the generator cylinder.

It remains to look for the Matveev moves and the slices, which we work out separately for the top and bottom part, and then the resulting sum.

i. bottom part (local):

The Matveev move is T*:

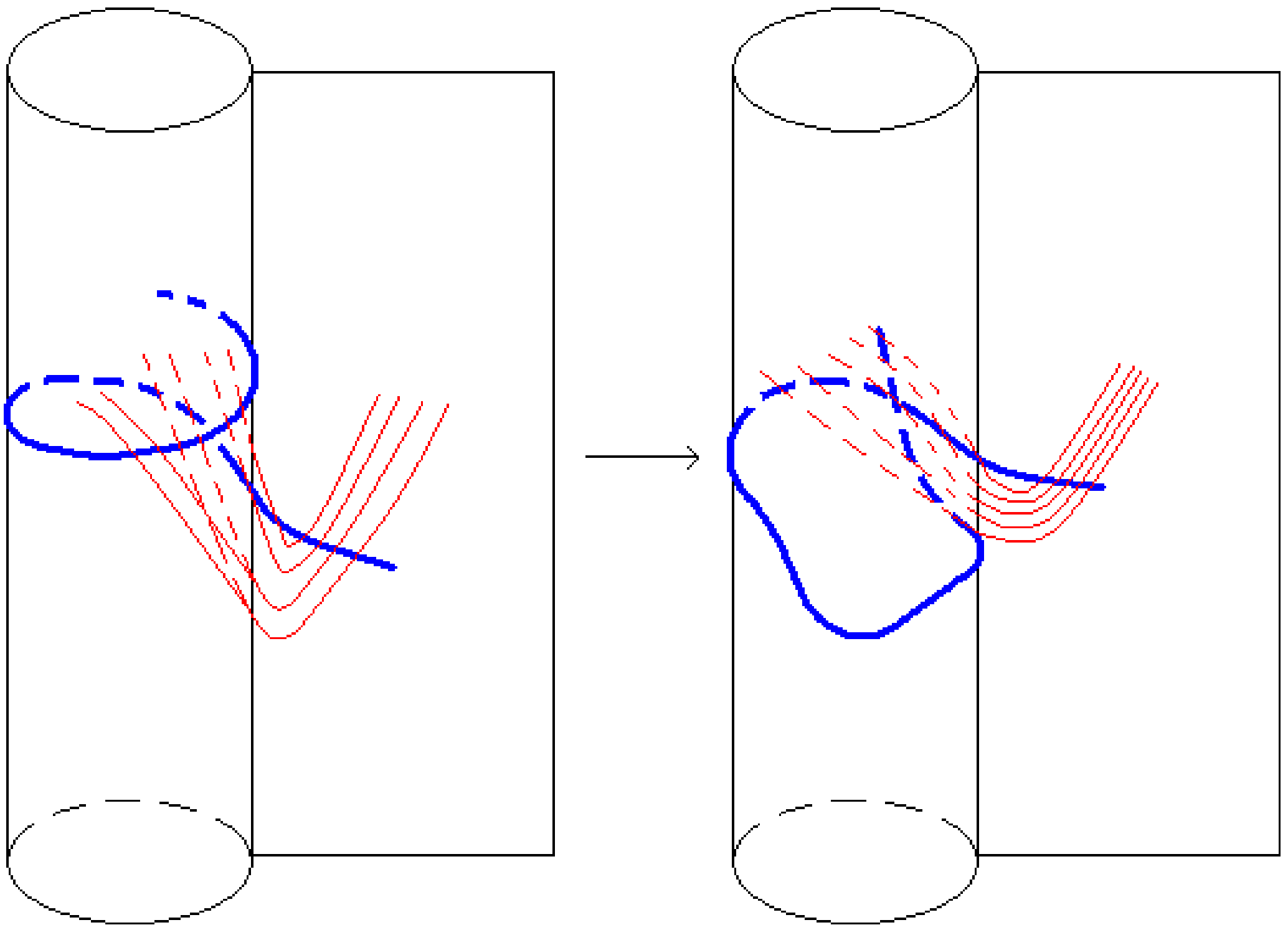

**Figure 241- The inverse- basic example- bottom part local- the move T* in Quinn model**

The sequence of slices:

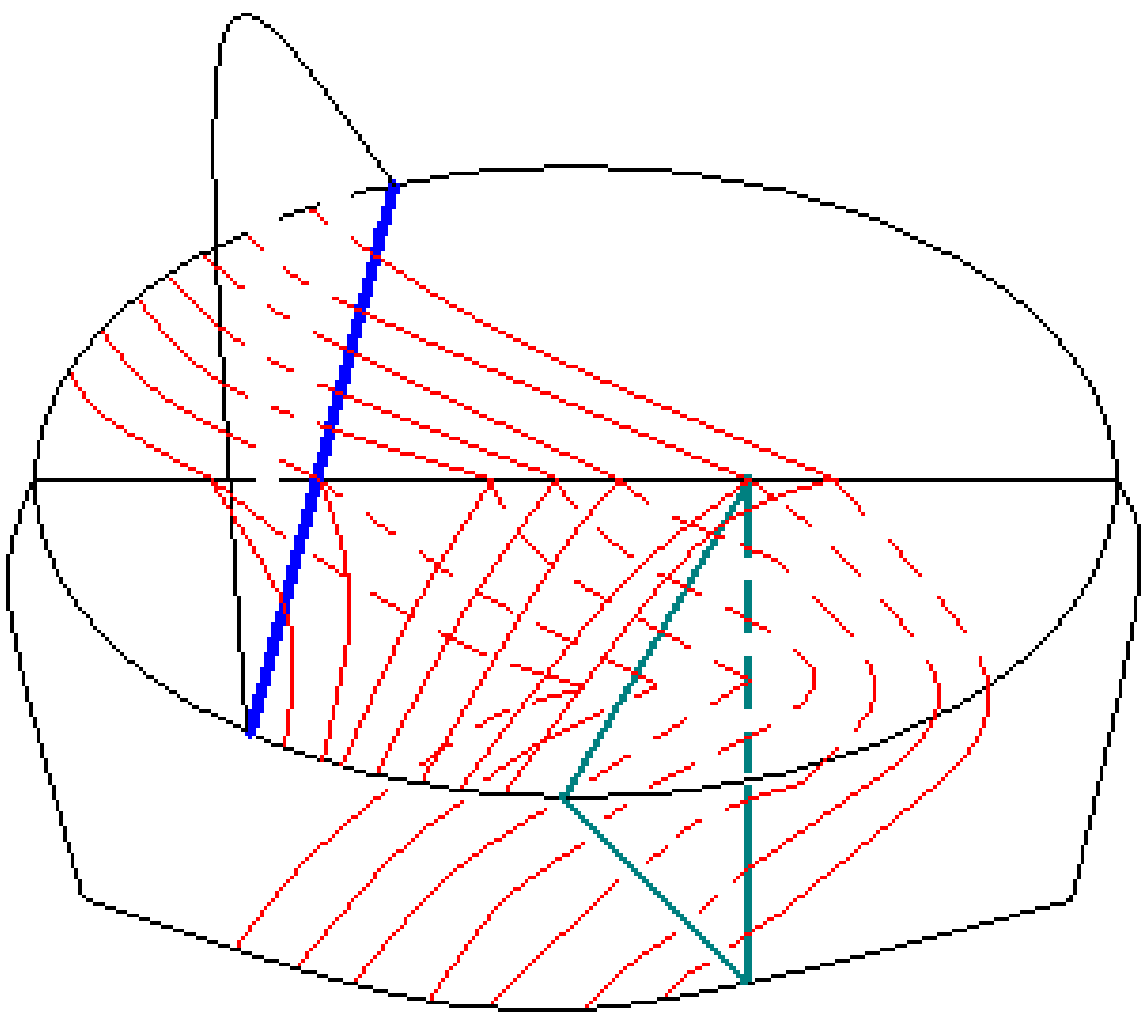

**Figure 242- The inverse- basic example- bottom part local- the move T* in local model**

ii. top part (local):

The Matveev move is T*:

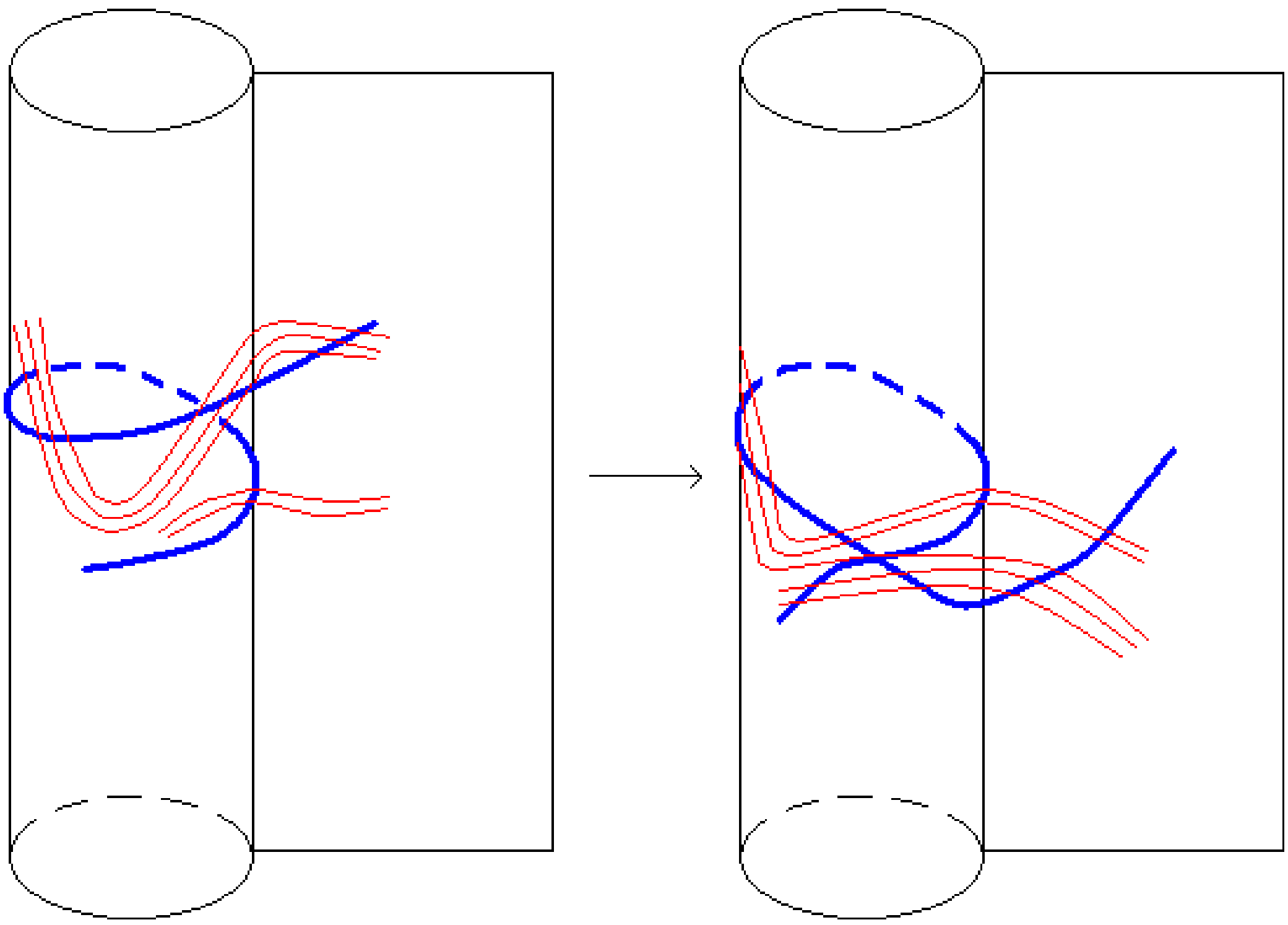

**Figure 243- The inverse- basic example- top part local- the move T* in Quinn model**

The sequence of slices:

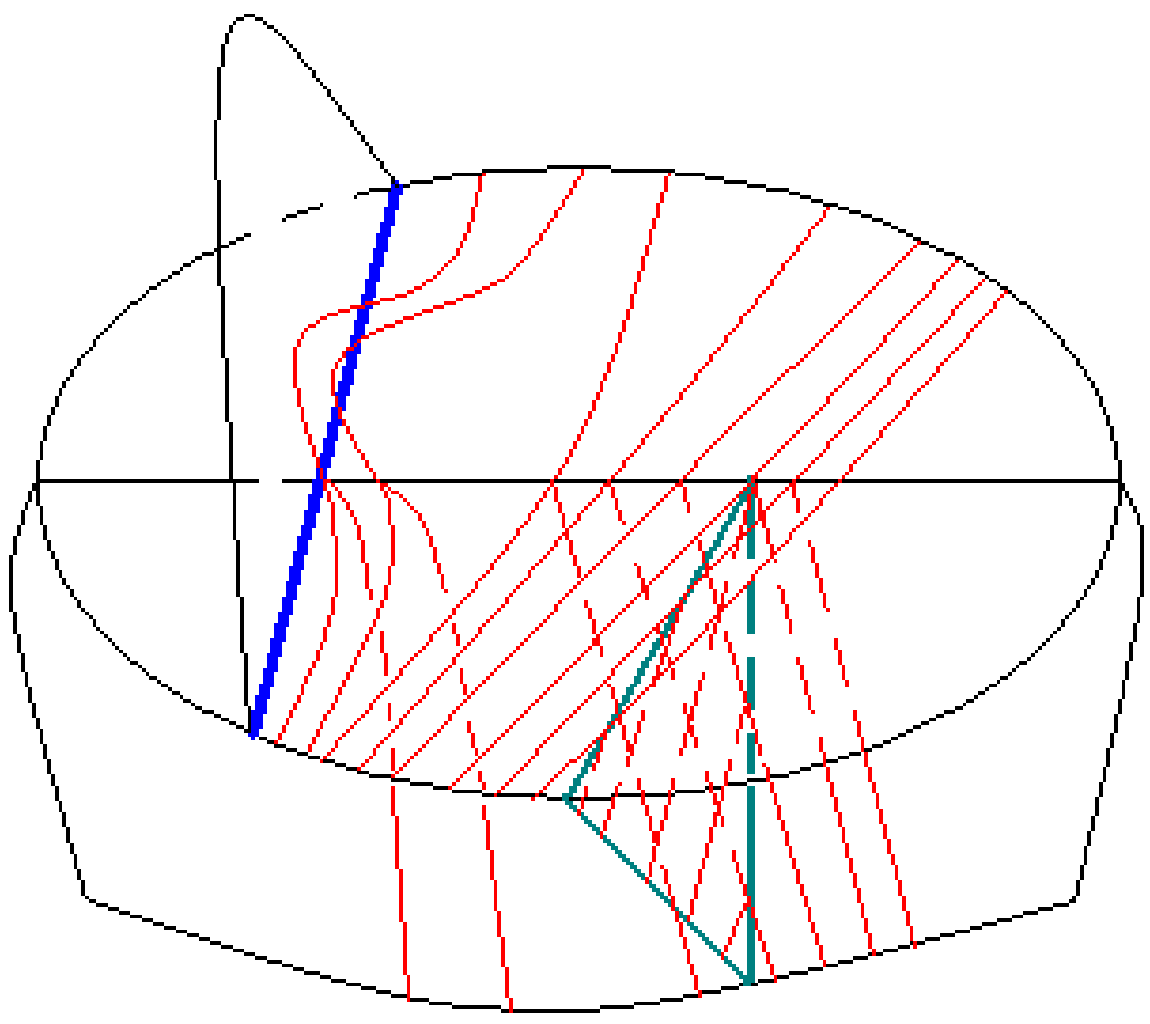

**Figure 244- The inverse- basic example- top part local- the move T* in local model**

iii. sum of both parts (local):

The Matveev move is $T_2^{-1}$:

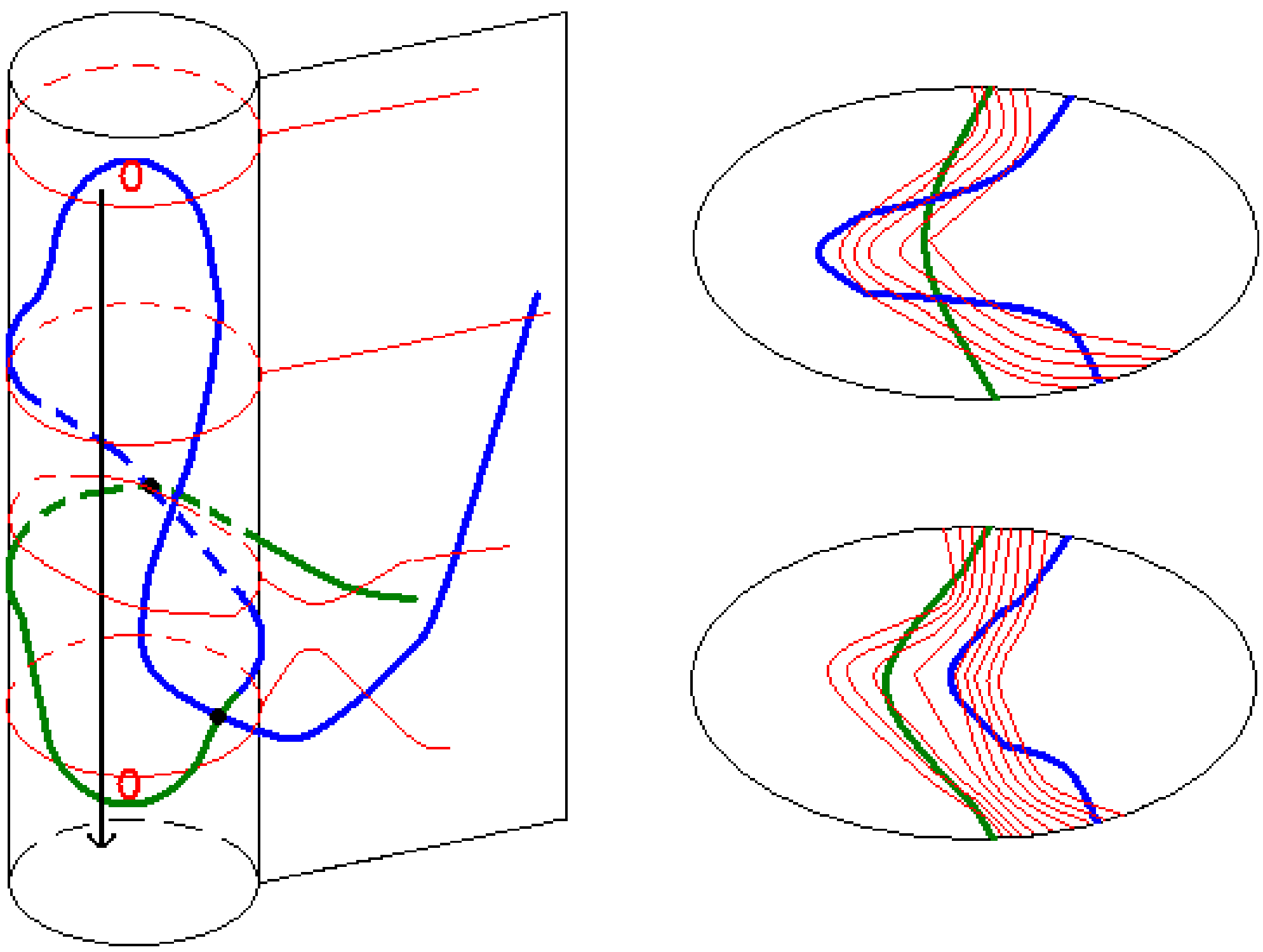

**Figure 245- The inverse- basic example- sum part local- the move $T_2^{-1}$ in Quinn and local model**

This shows, that we have not produced further topological relations during the steps in the basic example.

We remark, that P = ⟨a,b / R,S ⟩ can be transformed to Q = ⟨a,b / S,R ⟩, but we do not work that out.
We have to exchange the top relation with the bottom relation. This corresponds to dropping down the attaching curve along the generator cylinder and the rectangle. We already worked out this process in detail in the former section for the so called basic example and therefore we do not get new topological relations.

## *7.4 The 2-deformation*

In this chapter we speak about the extended prolongation. Let P = ⟨a,b / R,S ⟩ be a presentation of the 2-complex and P → Q = ⟨a,b / R,S,T = $cw^{-1}$ ⟩, where w is in F(a,b). We present Q in the Quinn model, where the hole c corresponds to the third generator cylinder:

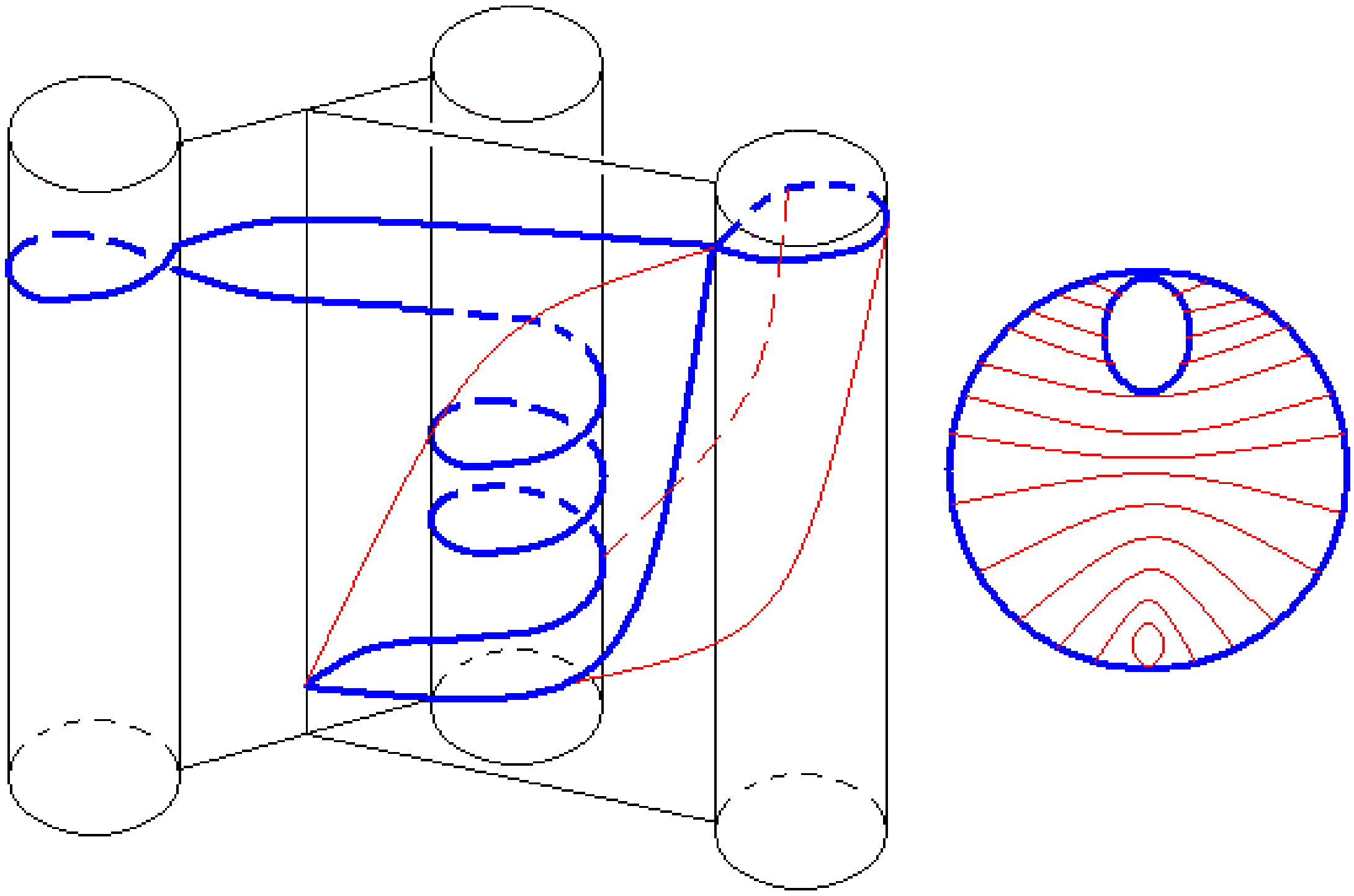

**Figure 246- The 2-deformation- extended prolongation in Quinn model**

The transformation Q → P:
First we collaps the (new) generator cylinder c, then we blow up the hole and retract c to (the boundary of T without c) = $w^{-1}$ along the slices:
There appear two types of slices:
The line segments are far from the entry of the cell T and the saddlepoint slice is near the entry. The collaps along the line segments corresponds to a relation in the Quinn list, we refer to chapter 4 . We can arrange the collaps, so that we also have line segments near the saddlepoint. By that the whole collaps is a composition of the same topological relation in the Quinn list and similar for the inverse (P → Q):

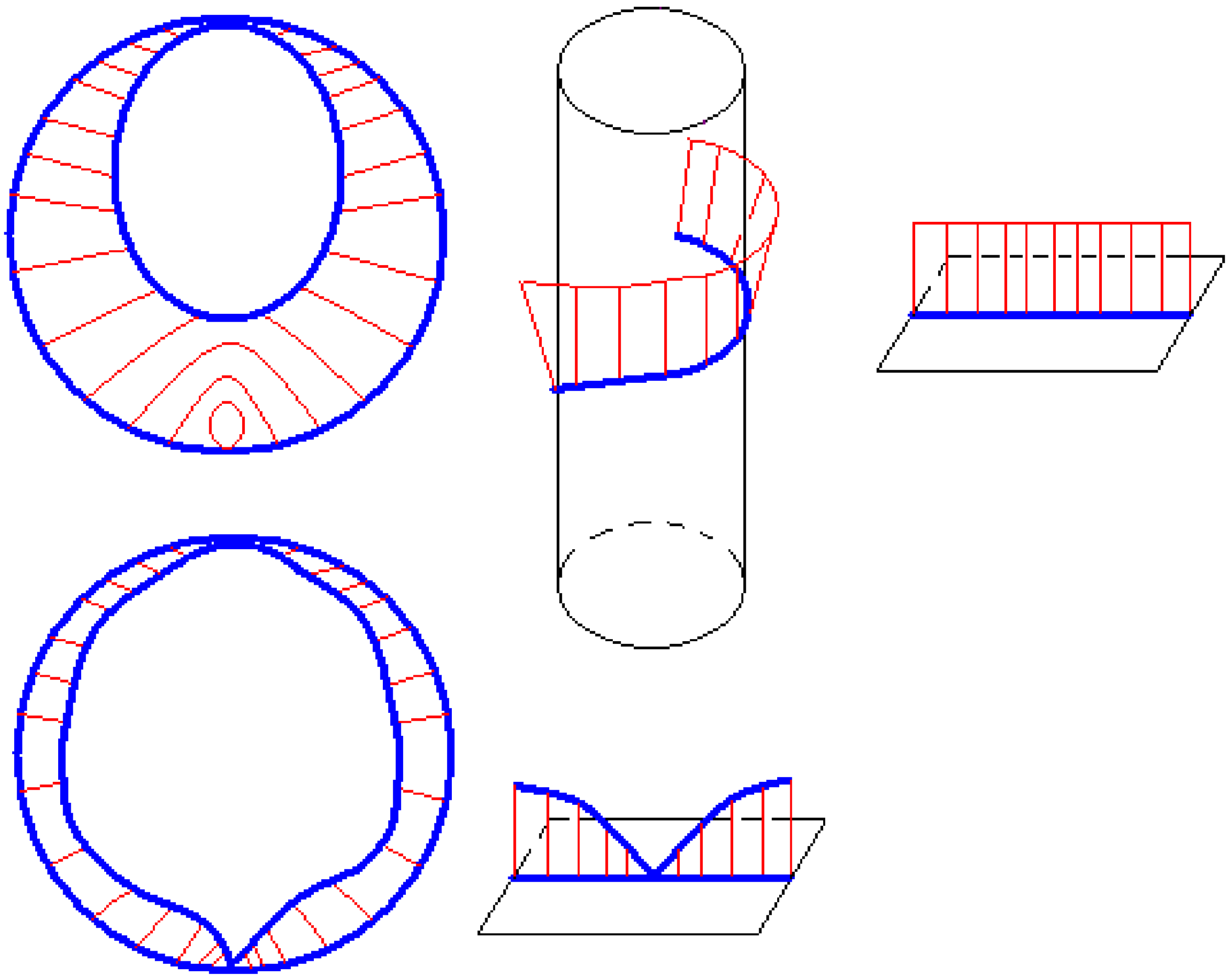

**Figure 247- The 2-deformation- collaps the perforated 2-cell**

# 8 Calculus – a TQFT example for the new sequence

## 8.1 TQFT - Compare vertex model with new sequence

In chapter 4 we deformed an attaching curve by an elementary Matveev move and applied our solution for the “good $T_3$ turn” to the “bad $T_3$ turn” case, i.e. (**before** perform the $T_3$ move) we change the sequence of slices near the vertex of the vertex model and get a new sequence.
This sequence can be extended to a topological relation. In chapter 4 we remarked, that it is still open, if this topological relation also leads to an algebraic relation in TQFT. We solve this question for a chosen tensor category:
We have to evaluate the new sequence and the sequence of passing a vertex in the context of roottrees and compare their associated homomorphism. The computations show that they are the same, hence for this example the answer is yes.
First we argue in the Quinn model, that a change to the new sequence is required:

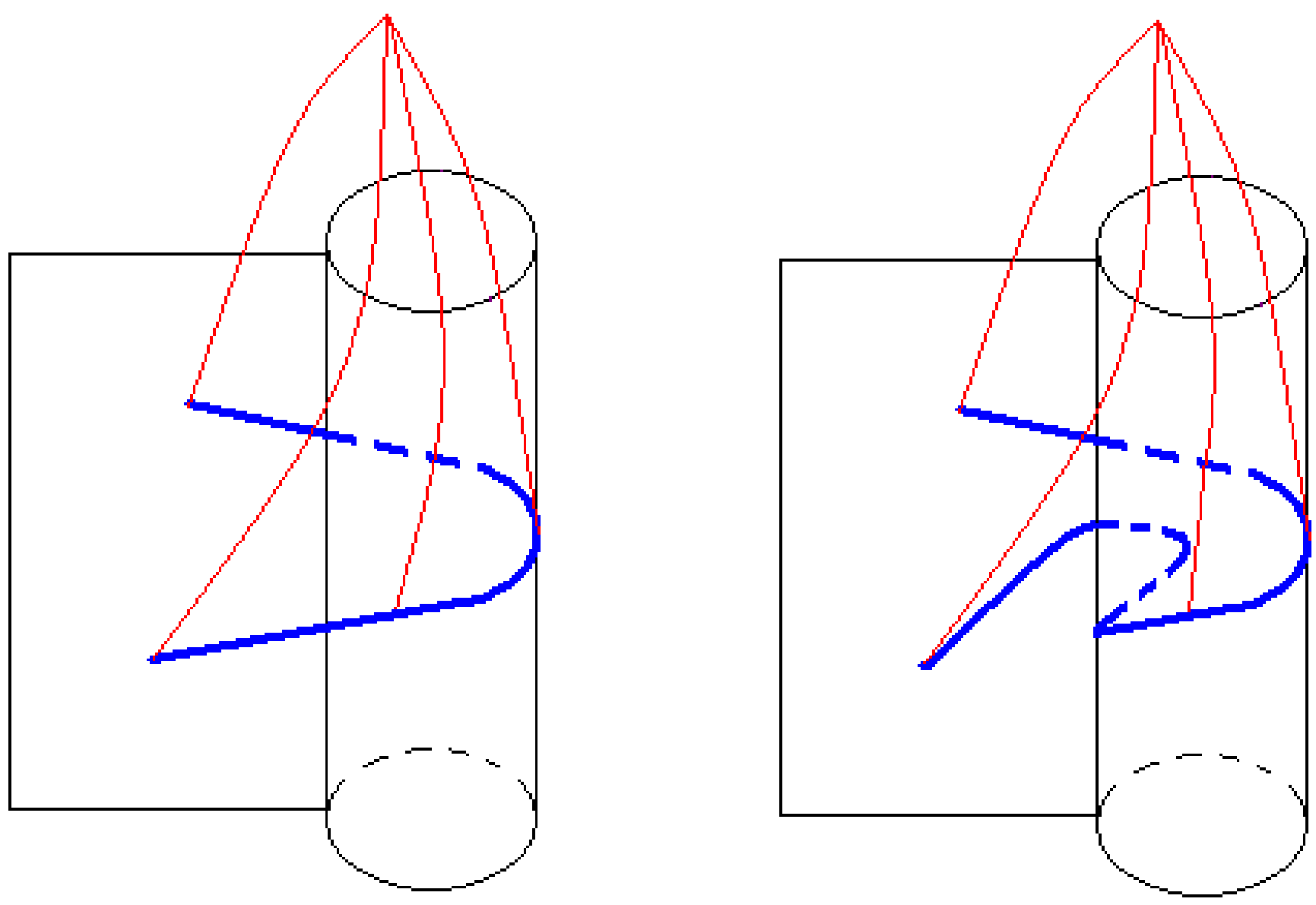

**Figure 248- TQFT- Compare vertex model with new sequence- apply $T_3$ turn**

To see that the move realizes a “ bad $T_3$ turn “, we look at the corresponding slices:

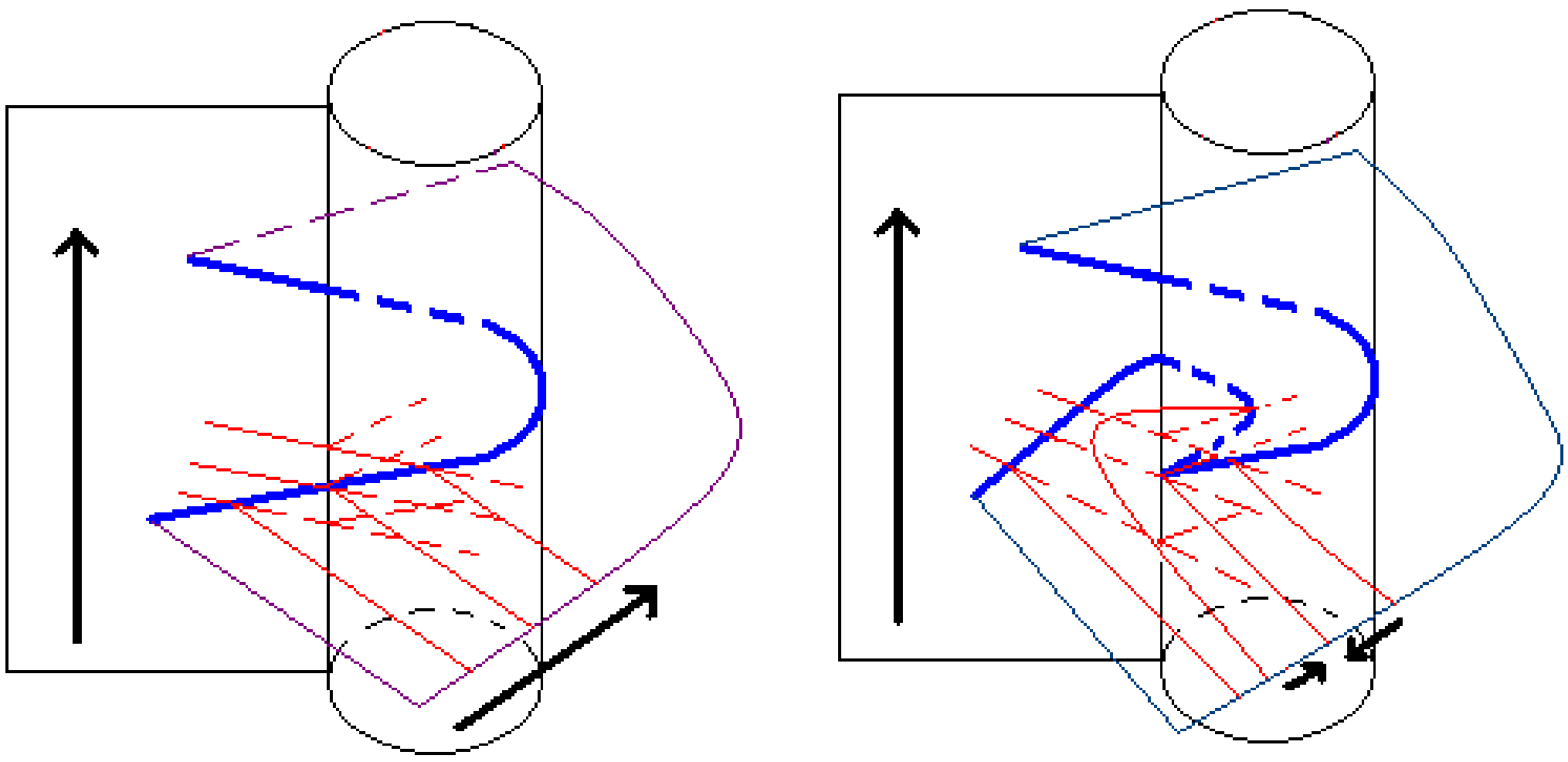

**Figure 249- TQFT- Compare vertex model with new sequence- apply $T_3$ turn- slices**

In the next figure we explain that, if we do not change the slices, we get a contradiction (the arrows in the picture indicates the height function on the bottom component):

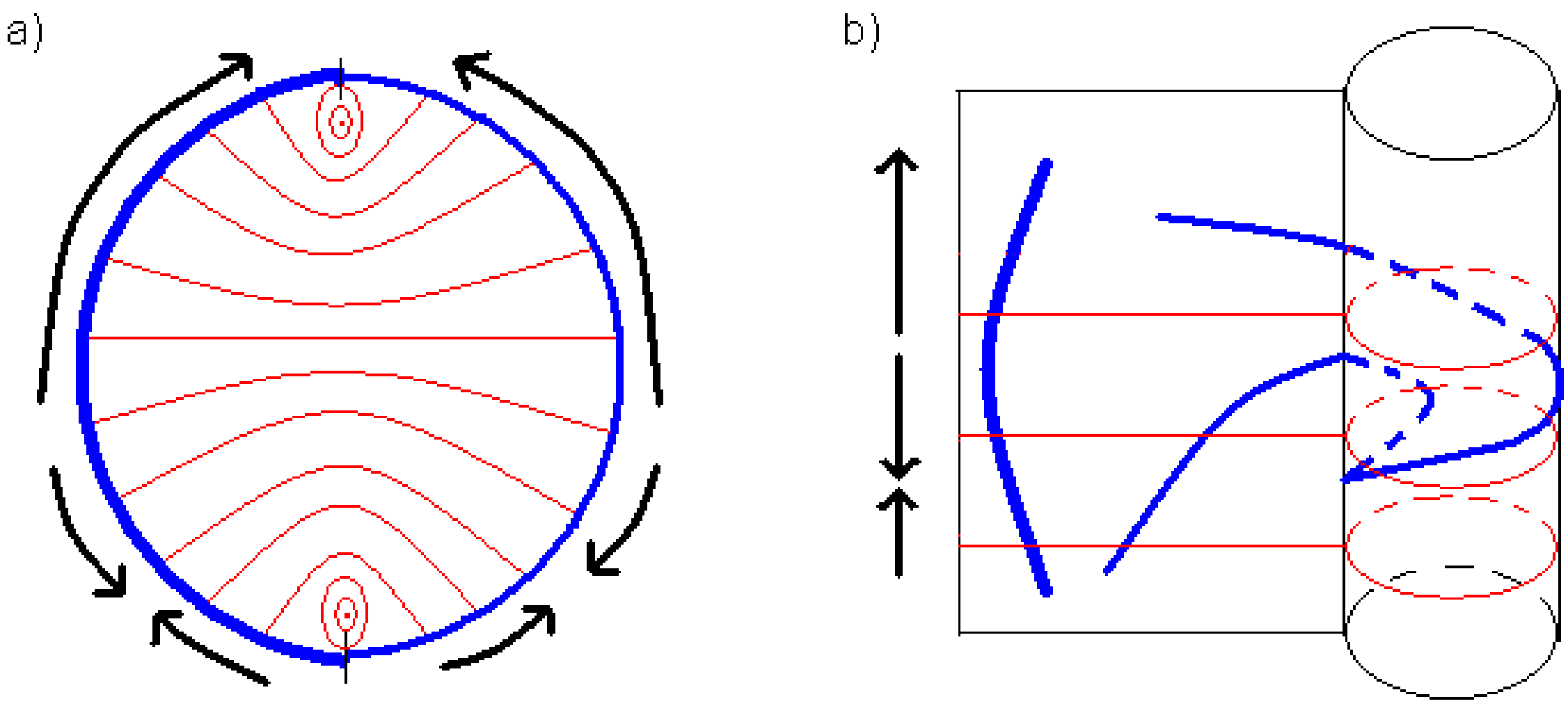

**Figure 250 - TQFT- Compare vertex model with new sequence- apply $T_3$ turn- slice do not change**

a)
the arrows indicate the corresponding height function:
Following the attaching curve from bottom to top, it increases to the first vertex, then it decreases to the second vertex and after that it increases until it reaches the top. This describes the height function for the right (thin blue) arc of the 2-cell, but the curve has to return to the bottom, as the (thick blue) line indicates. Now we see, that each slice (exept for the both circles) of the 2-cell connect both blue arcs and therefore the value of the height function of the slice in the 2-cell inherits to its boundary points, hence the value of the height function on the thick blue arc is the same.

b)
Transfer the height function to the thick blue line on the rectangle, then it intersects the slices of the retangle in points, so these also inherit the value of the height function. But this contradicts the fact, that the Quinn model requires a strictly increasing height function on the slices of generator cylinders and their connecting rectangle.

Hence a change of slices before performing the $T_3$ move is necessary !!!

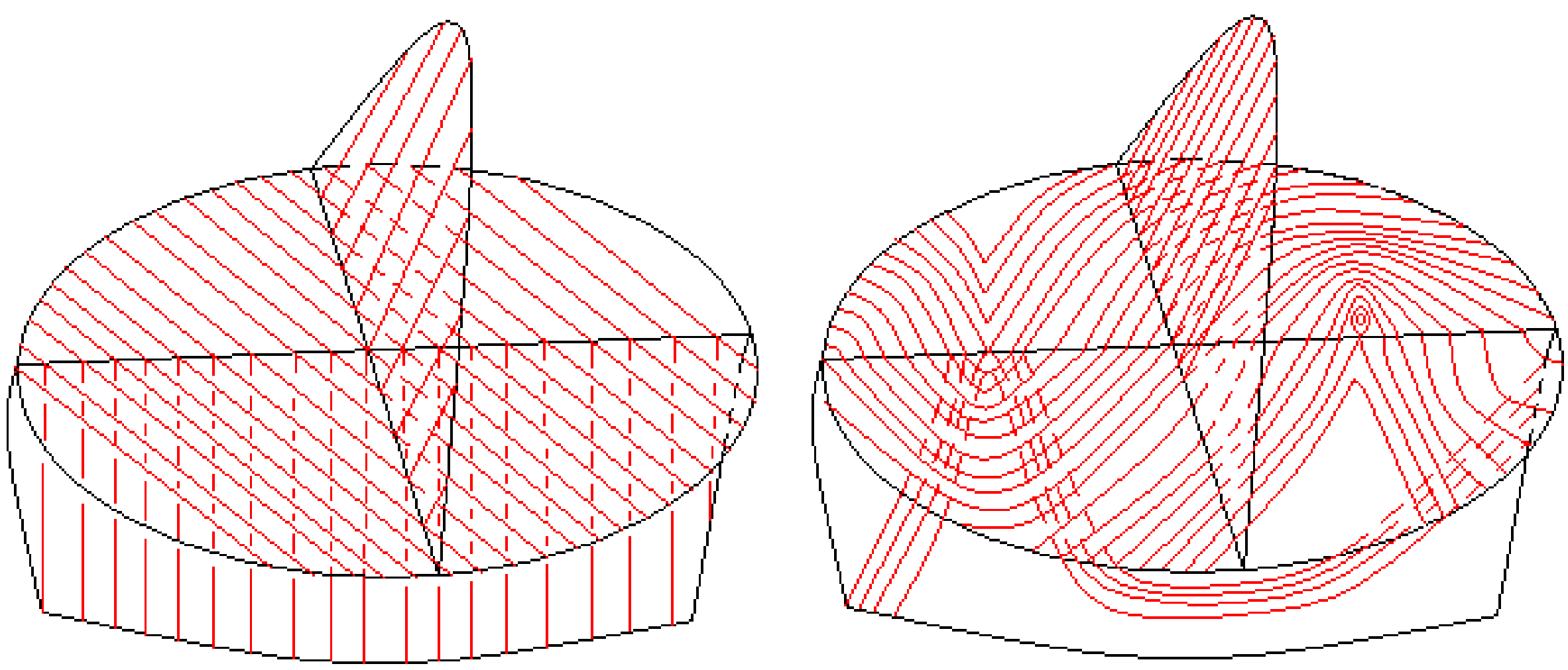

**Figure 251- TQFT- Compare vertex model with new sequence- change slice near vertex**

Transfer that local description into the Quinn model:

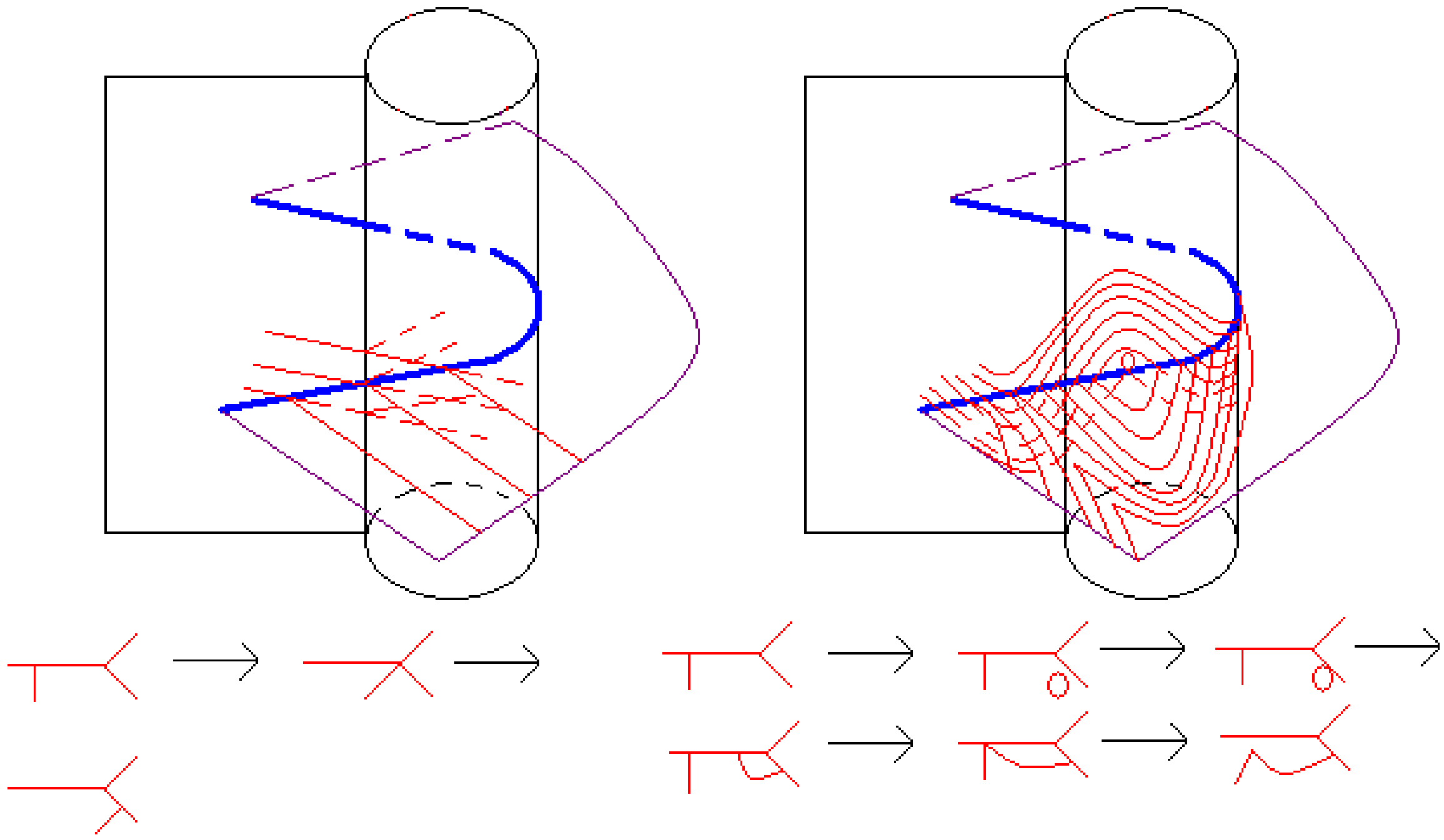

**Figure 252- TQFT- Compare vertex model with new sequence- change slices**

At the end we indicate, why that local change is sufficient for solving our problem: In a) we show the slices of the attached 2-cell, and we show in b) the slices at the rectangle and the generator cylinder. Again the arrows indicates the height function and show that it is strictly increasing on the slices of generator cylinder and their connecting rectangle:

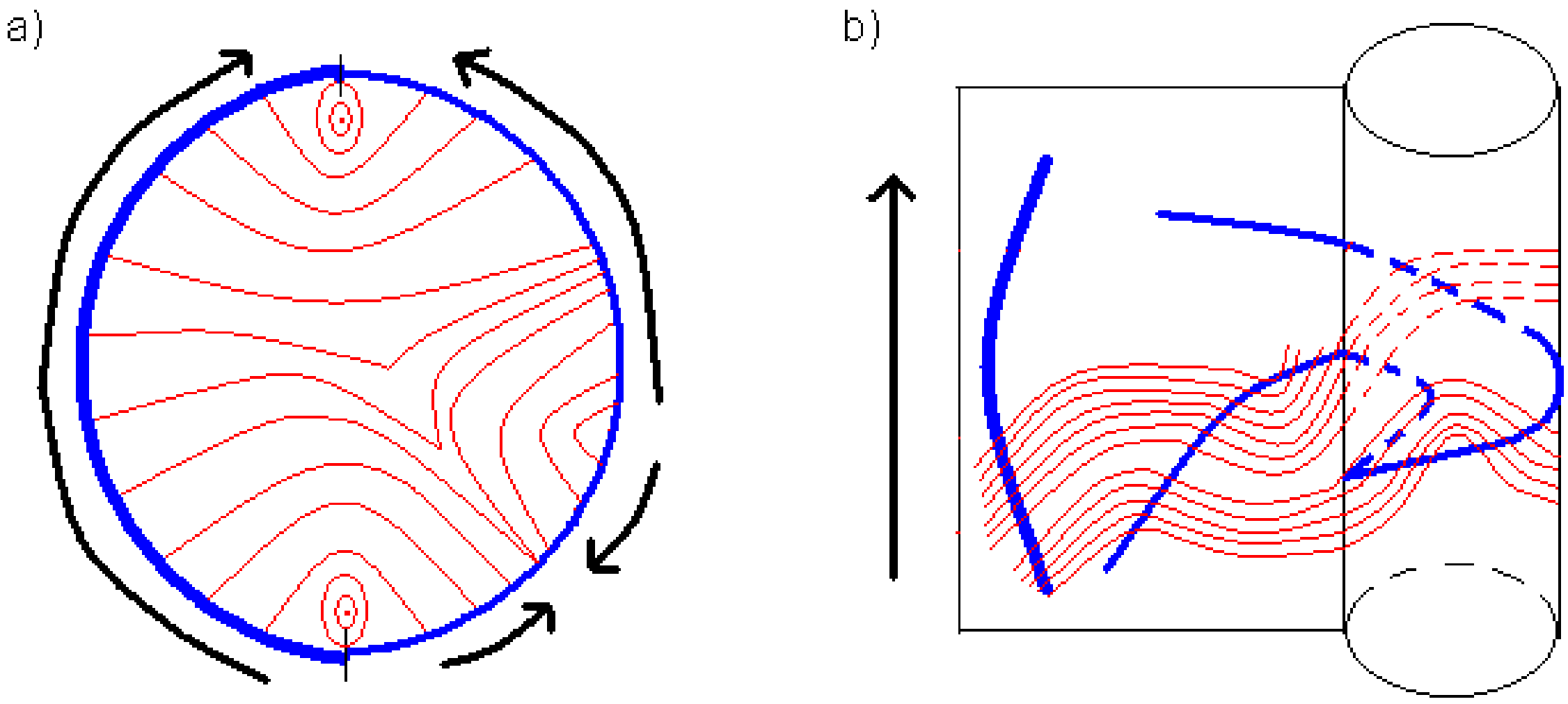

**Figure 253- TQFT- Compare vertex model with new sequence- slices when peform $T_3$**

## *8.2 Tensor category and roottrees*

We translate the sequences into the context of roottrees. Its boundary points are assigned to objects in a semisimple tensor category. We do not provide its general definition, but we present every property we use. The tensor category in our example comes from the representation of SL(2) mod $Z_5$ . It is generated by the simple objects {1,A}

$1 \otimes 1 = 1$
$1 \otimes A = A = A \otimes 1$ (1 is the unit)
$A \otimes A = 1 \oplus A$ ➔ $A = \bar{A}$ ($\bar{A}$ is the dual object)

General identities (pentagon, hexagon) in tensor categories are used to determine the matrices for associativity and commutativity, where the entries are in $Z_5$:

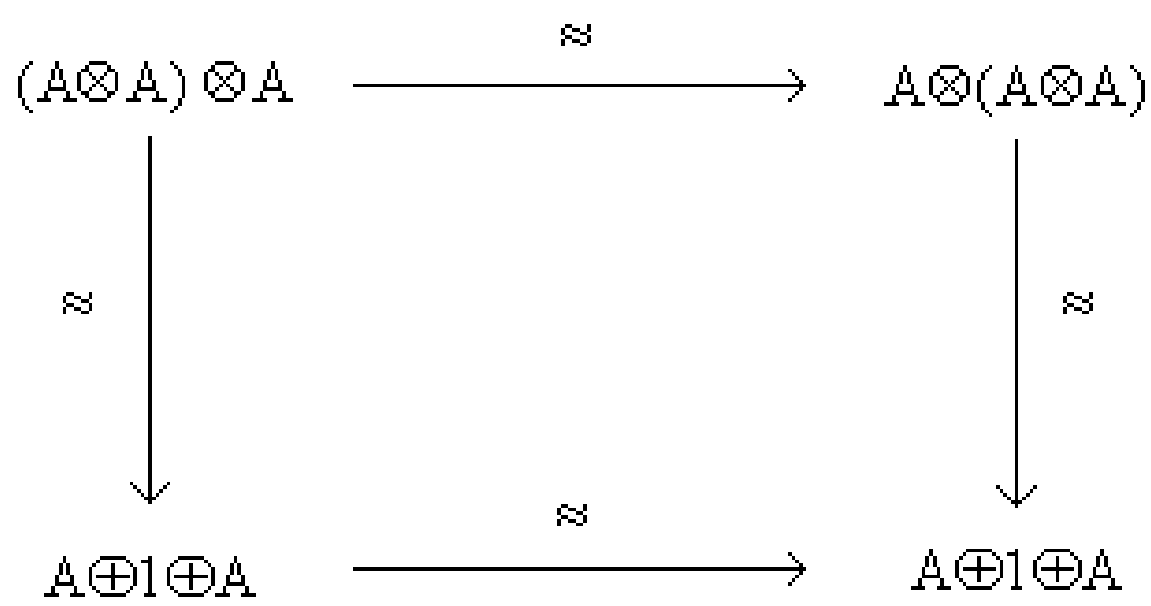

| | | |
|---|---|---|
| 2 | 0 | 4 |
| 0 | 1 | 0 |
| 3 | 0 | 3 |

and for commutativity

| | |
|---|---|
| 1 | 0 |
| 0 | 4 |

**Figure 254- Tensor category and roottrees- associativity diagram**

In general the roottree of a graph Y defines the state modul Z(Y) in TQFT. A local transition from a graph $Y_1$ to another graph $Y_2$ defines for the elementary bordism X between $Y_1$ and $Y_{2,}$ the homomorphism:

$Z_X: Z(Y_1) \rightarrow Z(Y_2)$

By composition of these we get the homomorphism assigned to the whole 2-complex, see **[Q2]** or **[Mül]** .

We will see, that the computation is deeply connected to the associativity diagramm. To make this more precise, we start with the translation of the sequence of passing a vertex:

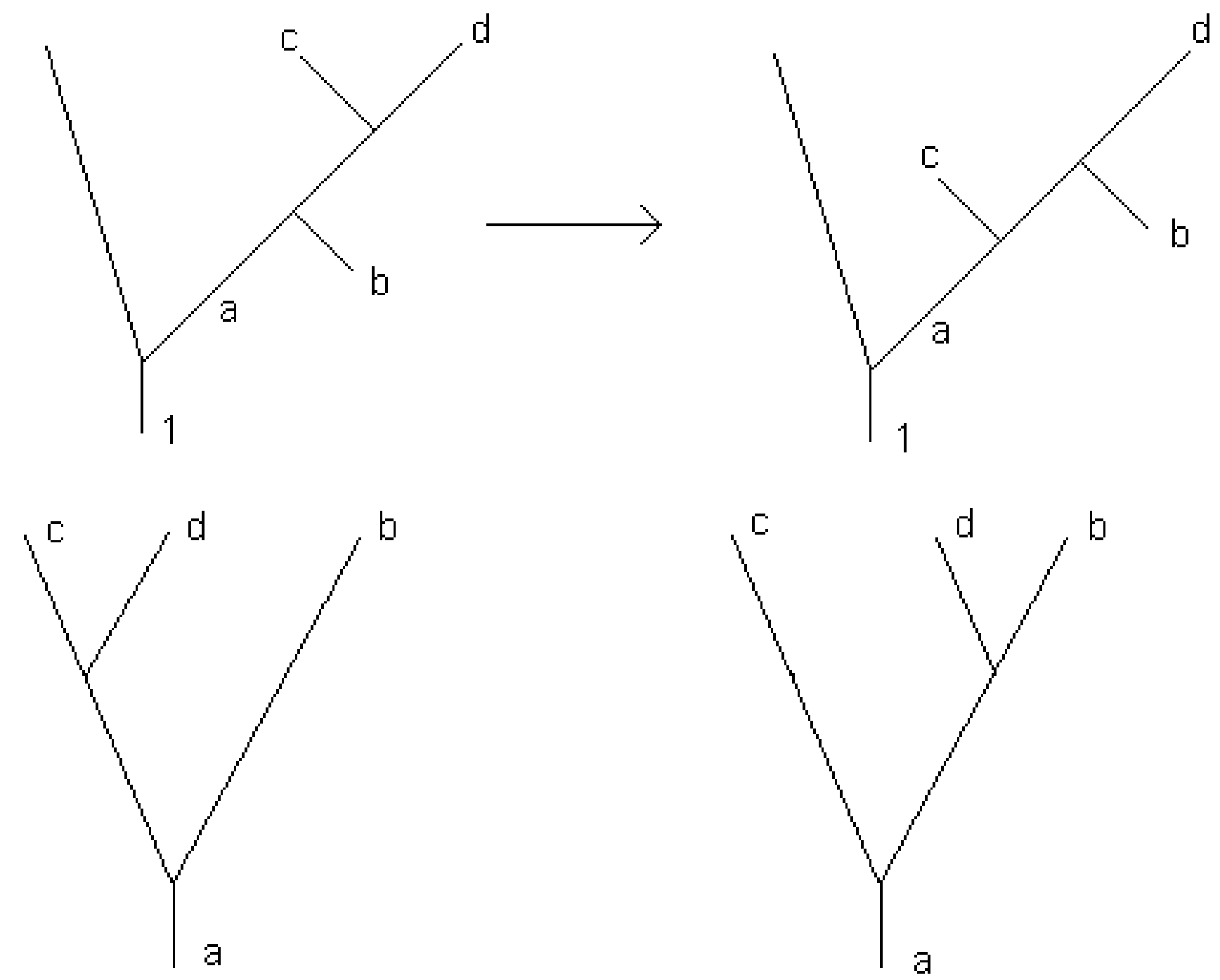

**Figure 255- Tensor category and roottrees- passing a vertex**

The first row shows the roottrees, which are labeled trees, the labels are chosen to be simple objects of the tensor category; here a, b, c, d and the root is always 1. In the first row we recognize the change of the c and b branches, which indicates the sequence of passing a vertex.
The second row shows the section consisting of the whole right branch. The geometry of the trees comes from the brackets of the tensor product respectively to associativity, for example:

Left :
a → (c ⊗ d) ⊗ b

Right:
a → c ⊗ (d ⊗ b)

## *8.3 New sequence as roottrees*

We describe the sequence of slices of the new sequence in terms of roottrees:

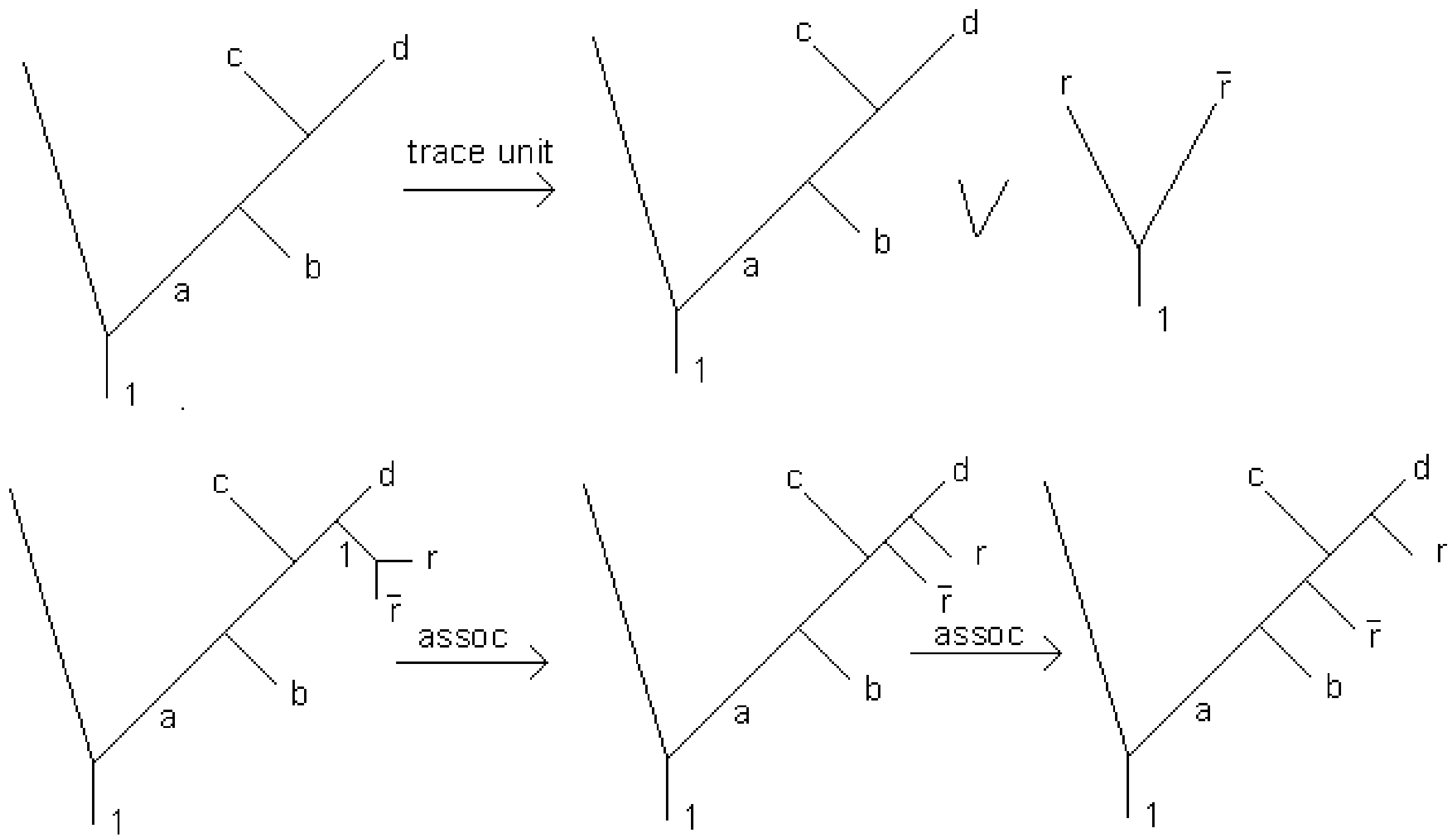

**Figure 256- new sequence as roottrees- general sequence- 1**

The trace unit describe the $S^1$. r and $\bar{r}$ stand for the algebraic identification of the boundary points. The next step moves the $\bar{r}$ branch onto the branch b by associativity. The definition of a (semisimple) tensor category requires that for each simple object r there exists a unique simple (called dual) object $\bar{r}$ with the property, that the sum decomposition of $r \otimes \bar{r}$ contains (unique) the unit 1.

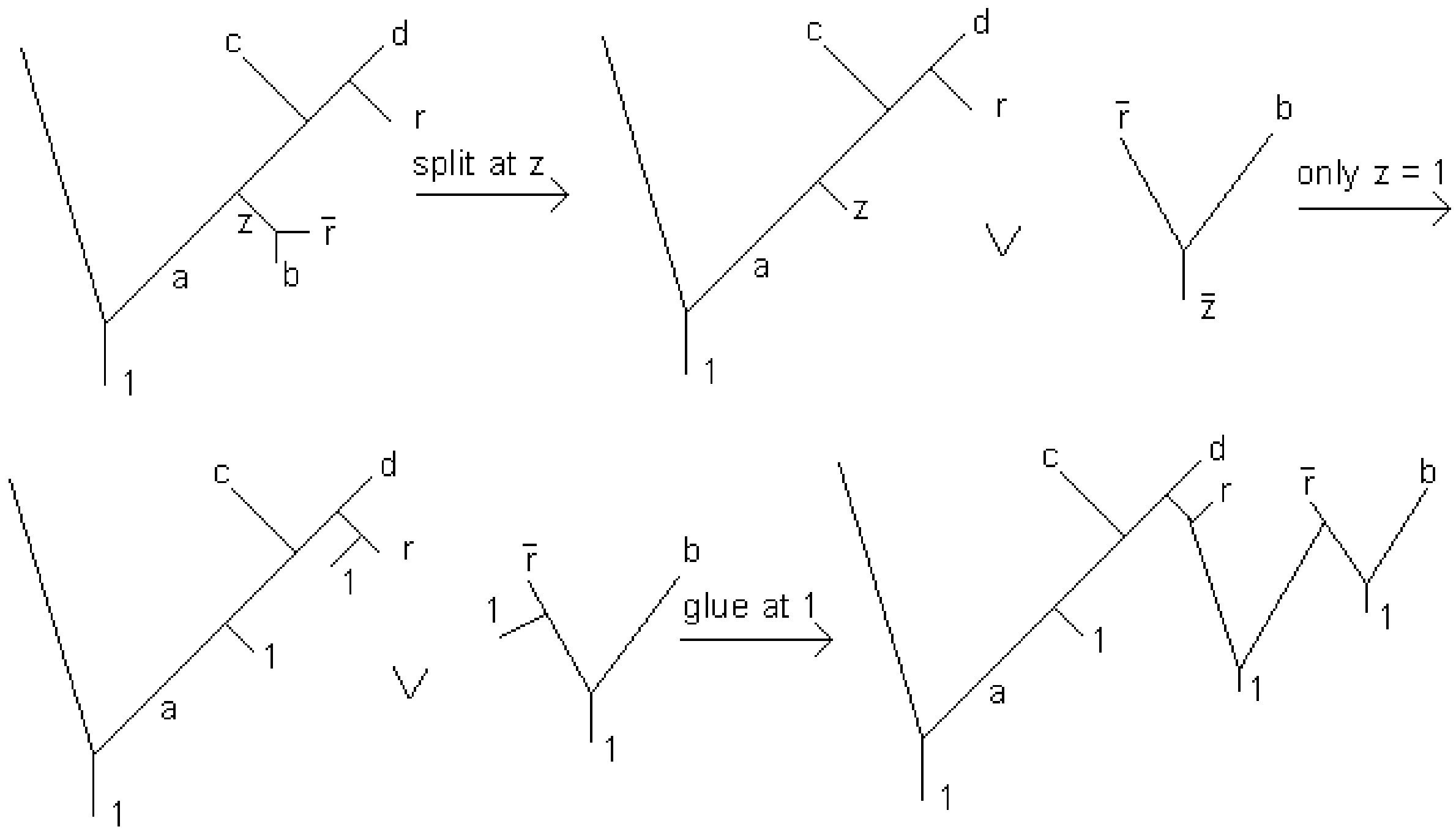

**Figure 257- new sequence as roottrees- general sequence- 2**

We remember (see chapter 4), that on the topological sequence, we split the arc arising from the $S^1$ at z and collaps the splitted branches to points (these branches correspond to additional flanges at saddlepoint) and consider the b branch as a prolongation of the $\bar{r}$ branch. But the only branch which can disappear without changing the roottree is the 1 branch (not the root) because 1 is the unit of the tensor product. To algebraically identify the both boundaries r, $\bar{r}$, we add 1 branches and connect these as a bridge for the r, $\bar{r}$ branches.
The next step simplifies by omitting the roottree some 1 branches.

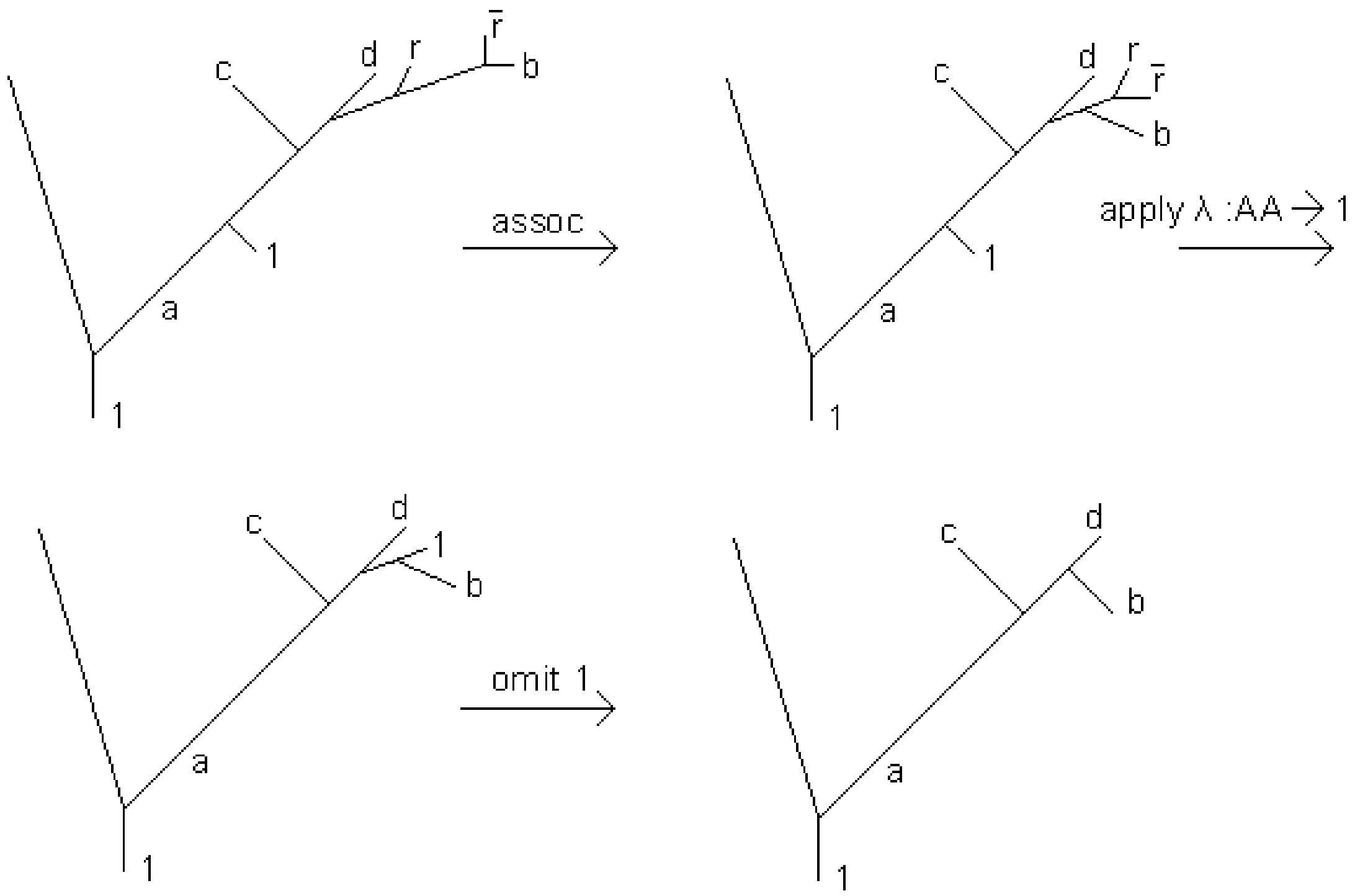

**Figure 258- new sequence as roottrees- general sequence- 3**

We use associativity to move the $\bar{r}$ branch onto the r branch. To annihilate these branches (which correponds to identify these by glueing the boundary points together), we apply a nondegenerate form $\lambda$: r ⊗ $\bar{r}$ → 1, which is an essential part in TQFT.

Furthermore we only consider semisimple tensor categories, i.e.:

a) the tensor product of 2 simple objects is a finite sum of simple objects
b) Hom (a,b) = Hom (1,1) = R•id, for a = b and zero otherwise, R = ring in a tensor product, in our example is R = $Z_5$

This leads to 2 cases. If the root to the r-$\bar{r}$ fork is assigned to 1, then the form $\lambda$: r ⊗ $\bar{r}$ → 1 is not zero, otherwise it is zero and the complete roottree disappears.

Finally omitting the arising 1 branch leads to the wanted roottree.
Now the question is, whether we get the same linear combination of the end roottrees when we apply both sequences to the start roottree ???

## *8.4 The trace unit of an ambialgebra*

To motivate, how to determine the trace unit, we first have to explain the coproduct Δ and the product m to make the rootree with 2 branches and root assigned to 1 to an ambialgebra.
An ambialgebra, where the trace unit exists is called special. The trace unit is defined as the solution c of the equation mΨΔ (c) = e.
Note, that e (unit) and c are linear combinations of roottrees with two branches, where 1 respectively A is assigned to the endpoints.

We start to explain the idea behind the trace unit. It is a fundamental concept in TQFT to associate the structure maps of the ambialgebra to elementary bordism with chosen incoming and outgoing boundaries. The coproduct Δ has two incoming (one has always the value e) and 2 outgoing boundaries, the product m has two incoming and one outgoing boundary. By exchanging the outgoing boundary of Δ by a horizontal reflection we get ΨΔ. Compose the pieces to build mΨΔ (c) = e, then we get an annulli with 2 incoming boundaries c and e and one outgoing boundary e. If we glue a disc in the hole with boundary c, the result is a pinched cylinder with incoming and outgoing boundary e:

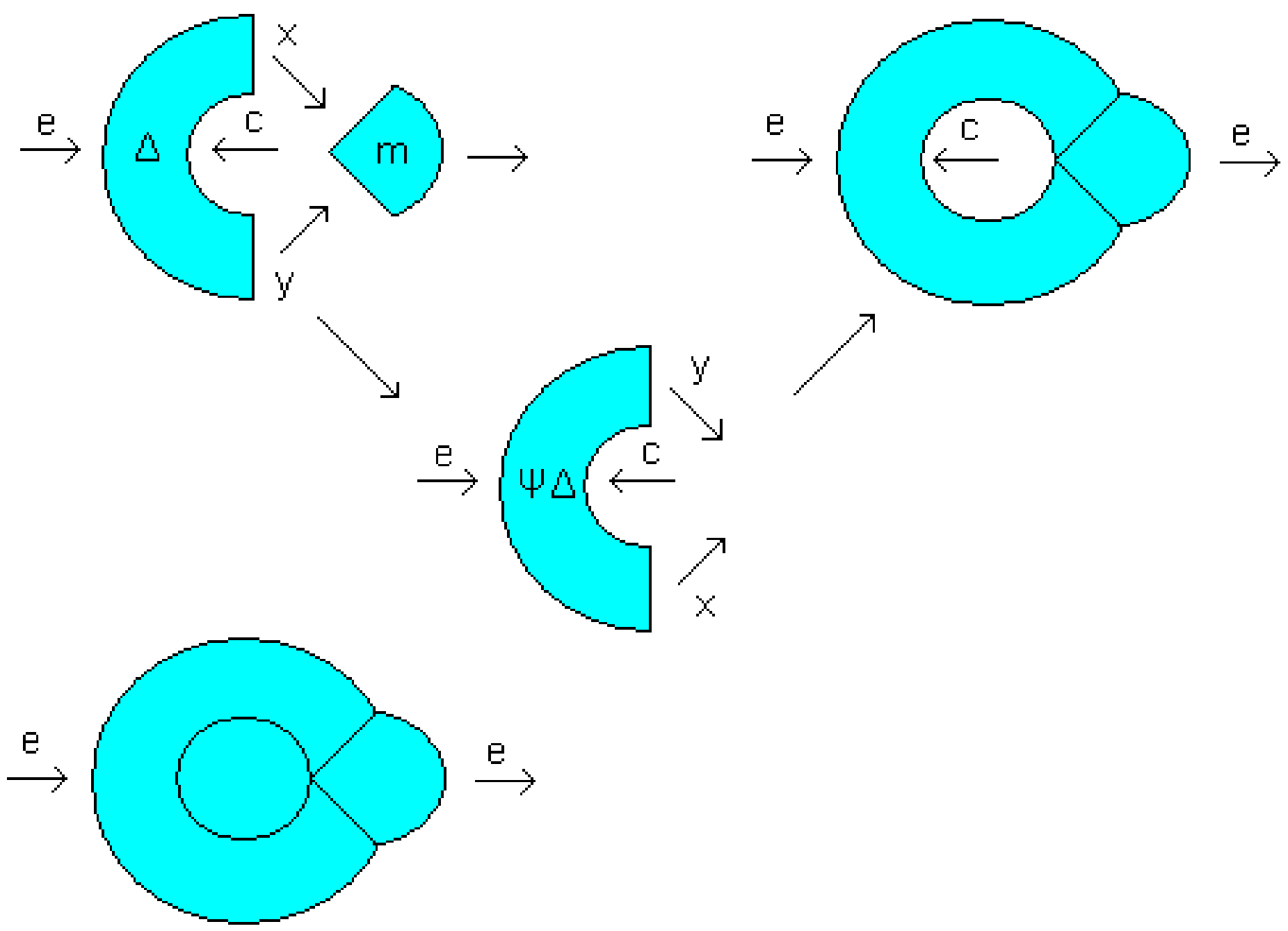

**Figure 259- trace unit– the idea of construction**

A pinched cylinder is a product bordism, so it induces the identity as homomorphism in TQFT. We conclude:
The definition of the trace unit corresponds to the existence of a disc in a TQFT.

To compute the trace unit, we formulate the structure maps of the ambialgebra in terms of rootrees. We discuss, which pair of simple objects are suitable for the endpoints of the roottrees, however the main part is to determine the coefficient for the linear combination of the roottrees.

We start with the coform Λ: 1 → A⊗A = 1⊕ A

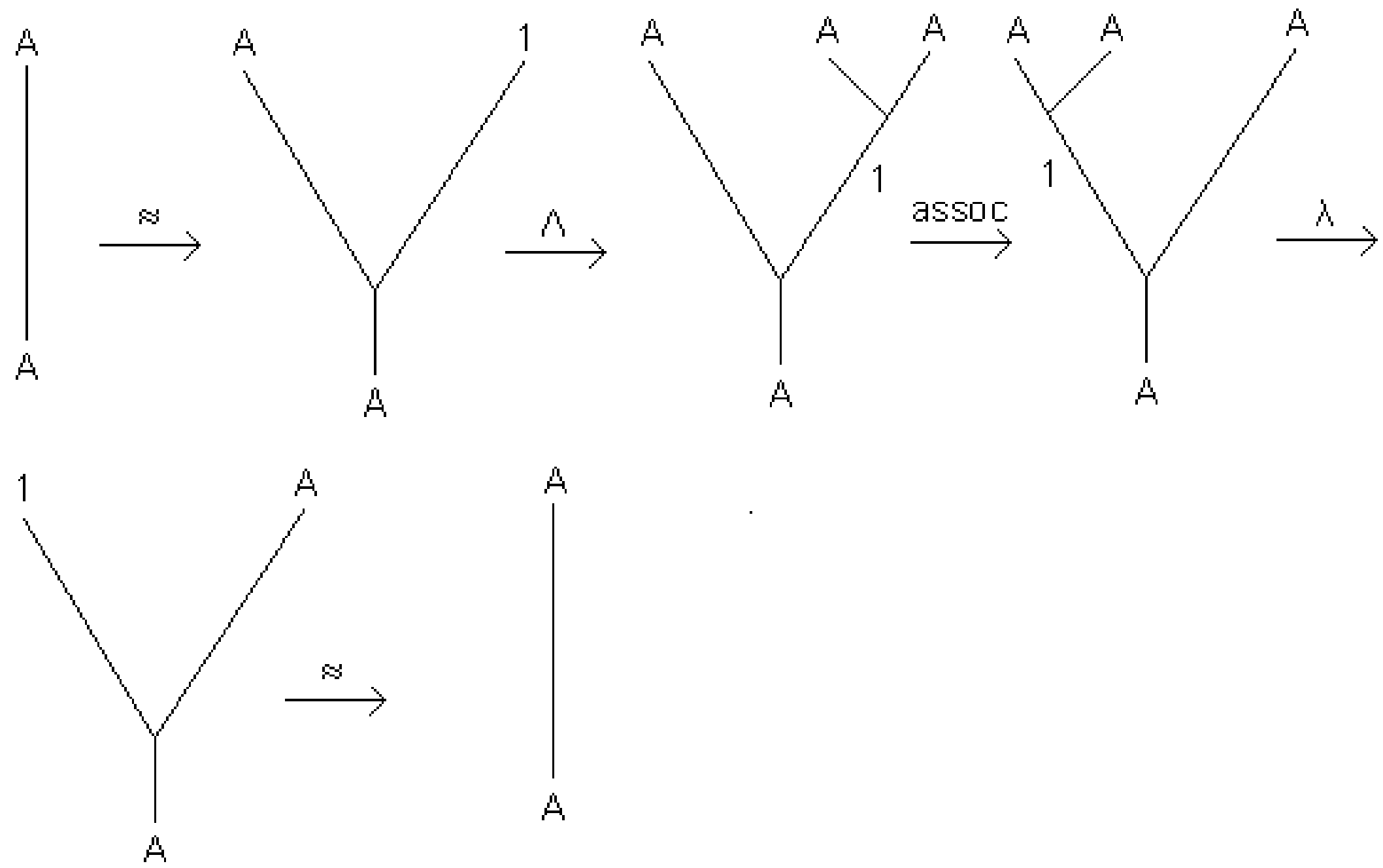

**Figure 260- trace unit– the coform**

The coform has to fulfill the requirement, that it admits a nondegenerate pairing with the form λ: A⊗A = 1⊕ A →1, hence the definition implies that the composition above has to be the identity. We remark, that after associativity we only consider the root = 1 for the left A-A fork, because this is the only case where λ is nonzero (see properties of a semisimple tensor category). Using the associativity matrix we see (more details to that type of computation in chapter 8.5) that the factor is 2, hence the factor α of the coform
Λ: 1 → A⊗A = α (1⊕ A) must be $2^{-1} = 3$ in $Z_5$ , since the composition has to be the identity.

We look at Δ:
The first observation from the separated rootrees in the last figure is, that: b = $\bar{a}$ .
The factor arises from the coform and assoc1; assoc 2 is over root 1 and therefore the identity (see associativity matrix).

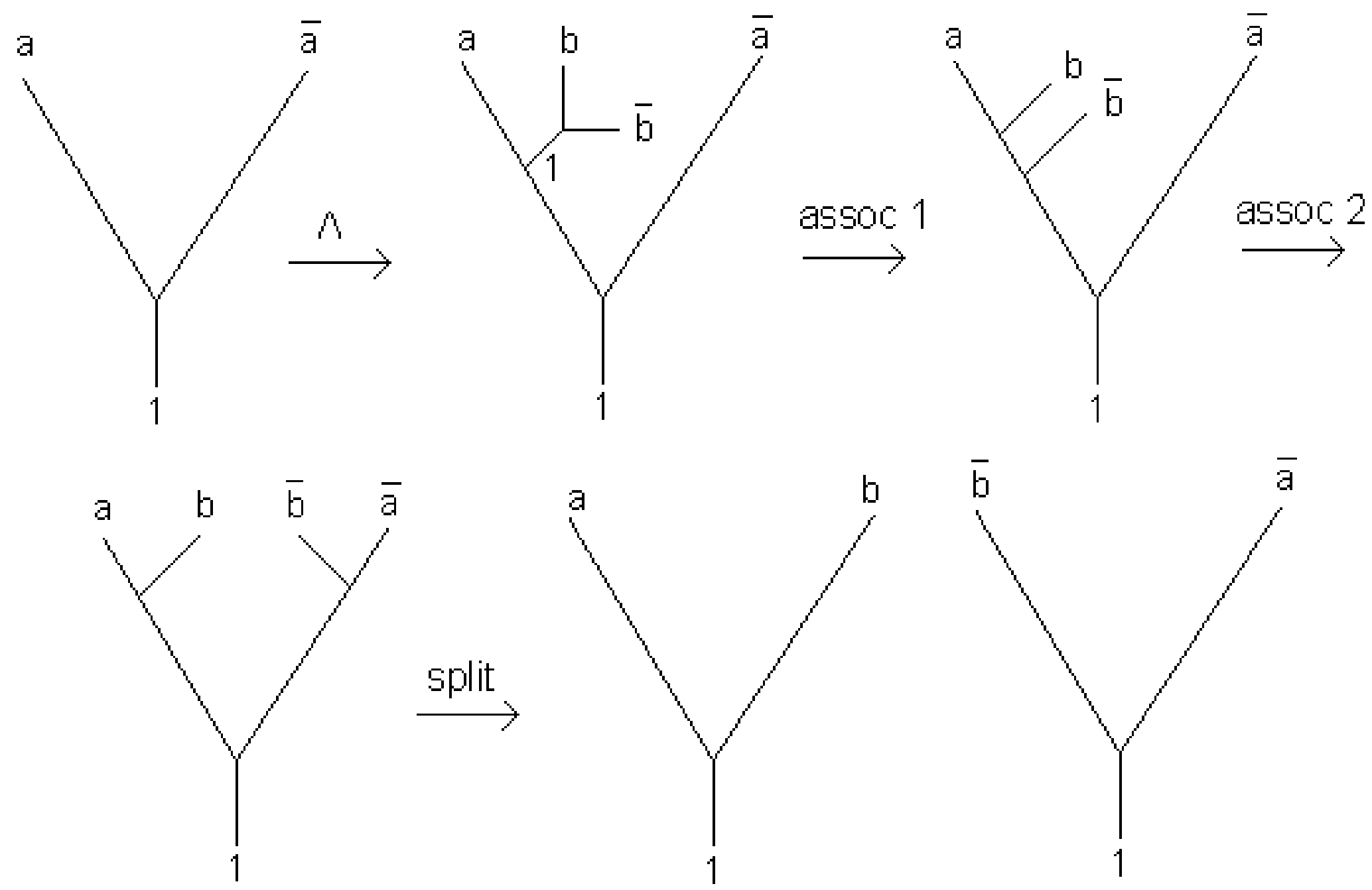

**Figure 261- trace unit– the coproduct**

From this matrix we also develop the factor for assoc 1, which is 2. We use the fact, that in the second last roottree the branches directly parting from the root must be assign to 1, other cases get lost after splitting, since we only consider roottrees with root assign to 1.

We consider the product m:
First we observe from the last roottree that: b = a
Again assoc 1 is over 1 hence the identity, assoc 2 provides the factor 2, we can repeat our argument for λ (the root for the $\bar{a}$ - b fork has to be one) as in the computation for Λ:

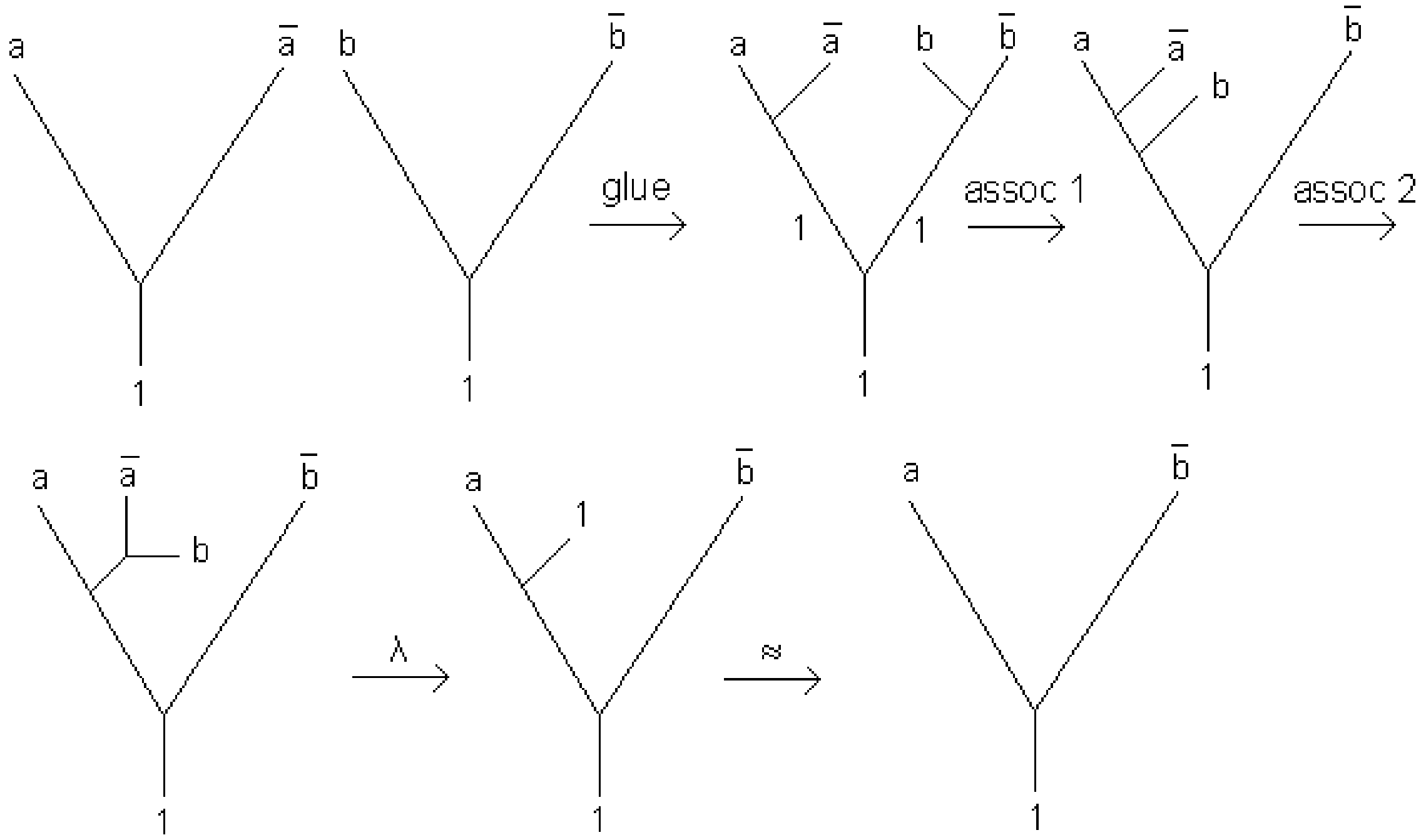

**Figure 262- trace unit– the product**

We compute the composition mΨΔ, where Ψ only changes the positions of the roottrees (It is not the Ψ from the commutativity matrix) and it does not change the factor:

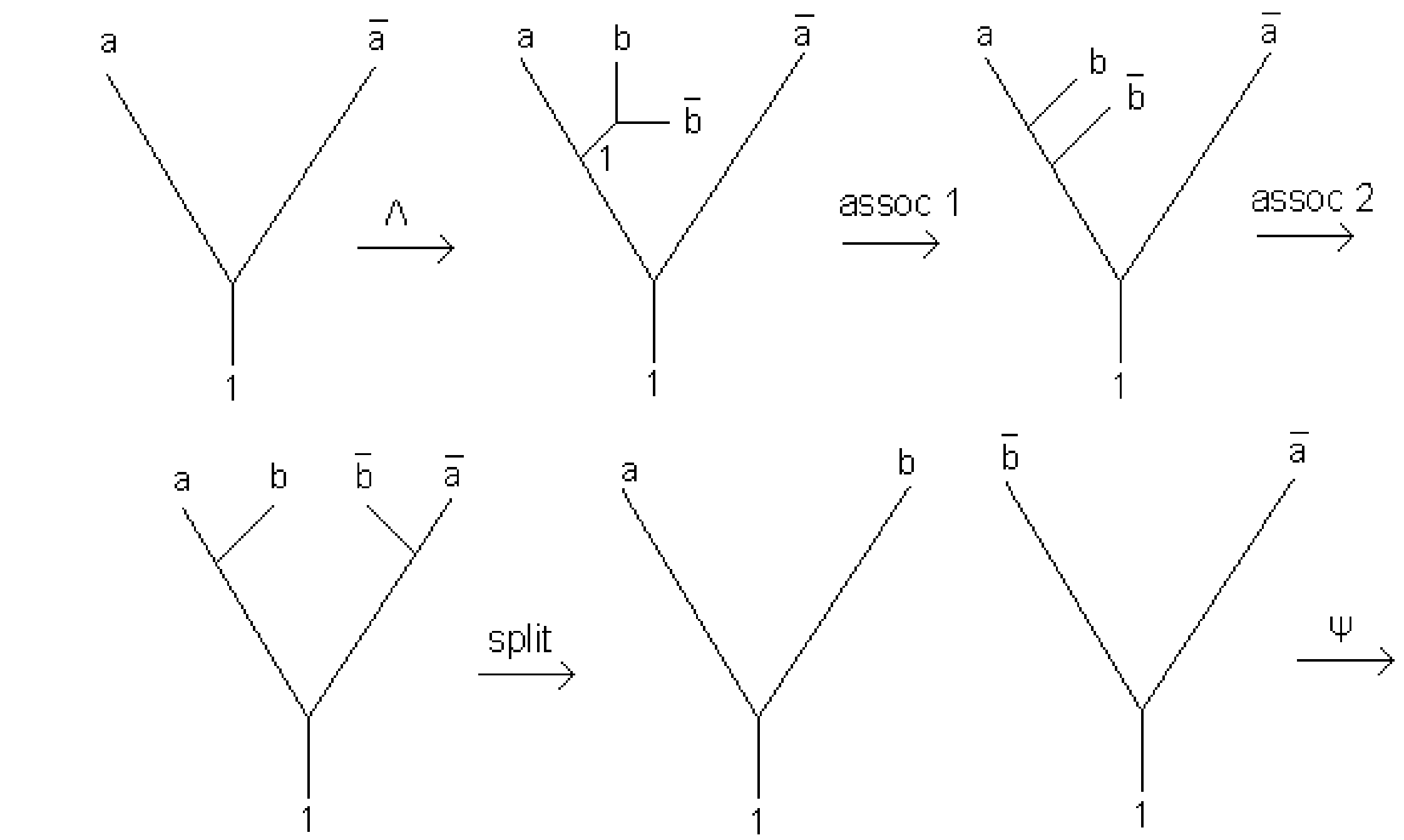

**Figure 263- trace unit– the composition 1**

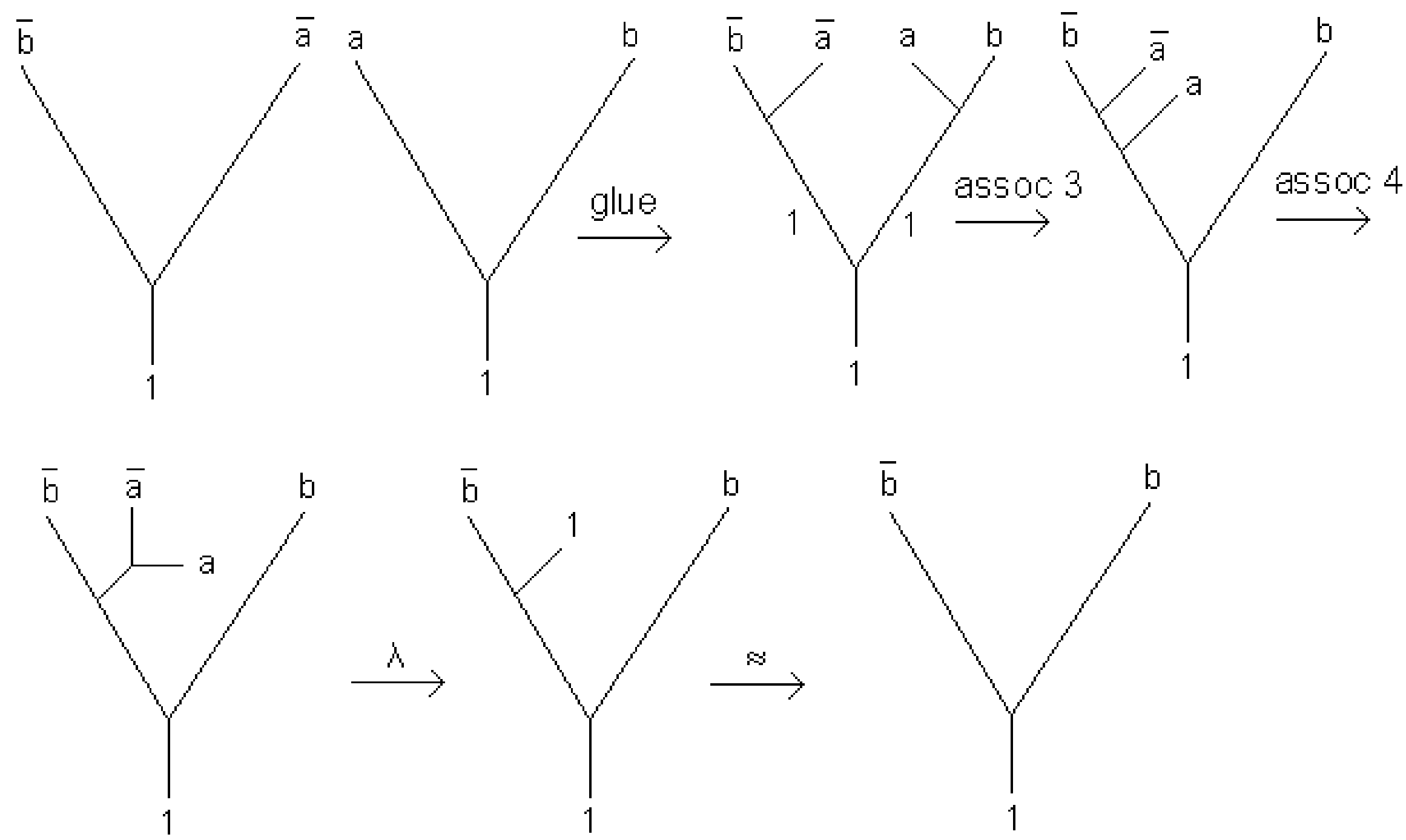

**Figure 264- trace unit– the composition 2**

We summarize the factors in the composition mod $Z_5$:
$\Delta \rightarrow 3 \cdot 2$
$m \rightarrow 2$

$m\Psi\Delta \rightarrow 3 \cdot 2 \cdot 2 = 2$

since m deliver factor 2 and if we denote the rootree with 1 objects assign to the fork as I and the other (with A objects assign to the fork) with A, the unit $e = I + e_1 A$ for the product is defined by $m (e \otimes x) = x$, where:
$x = I + x_1 A$, then:

$m (e \otimes x) = I + 2e_1 x_1 A = I + x_1 A = x$

$\Rightarrow e = I + 2^{-1} A = I + 3A$:

We solve the trace unit $c = I + c_1 A$ in $Z_5$:

$m\Psi\Delta(I + c_1 A) = I + 2c_1 A = I + 3A$

$\Rightarrow c_1 = 4$

### 8.4.1 Evaluate the trace unit in general

In general (i.e. in our case for a semisimple tensor category) we can evaluate c by using the methods introduced above. The idea is to balance the appearing factors for each simple object a as follows:

We can not assume that in general assoc over root 1 induce the identity.
Note (use $\bar{a}$ = b), that assoc 3 is inverse to assoc 2, hence their factors annihilate by composition and can be omited in our further consideration.

Nondegenerate pairing of λ and Λ provides the relation:
factor (Λ) = factor $^{-1}$ (assoc)

For the product m we have one times the associativity:
factor (m) = factor (assoc)

The coproduct gets his factor from the coform and associativity:
factor (Δ) = factor (Λ) • factor (assoc)

The unit has the factor inverse to the product m:
factor (e) = factor $^{-1}$ (assoc)

For each simple object now we can evaluate the equation mΨΔ (c) = e for the factors:

factor (assoc) • factor $^{-1}$ (assoc) • factor (assoc) • factor (c) = factor $^{-1}$ (assoc)
$\Rightarrow$ factor (c) = factor $^{-2}$ (assoc)
= factor $^{2}$ (Λ)

But if we consider the sequence of roottrees for mΨΔ (c), the factor is an element from the composition of Λ and λ , in more detail (use $\bar{a}$ = b):
$\Lambda: 1 \rightarrow \bar{a} \otimes a$

$\lambda: a \otimes \bar{a} \rightarrow 1$

Therefore the corresponding element for factor ($\Lambda$) in Hom(1,1):
$\lambda\Phi\Lambda \Rightarrow c = (\lambda\Phi\Lambda)^2$ where $\Phi: \bar{a} \otimes a \rightarrow a \otimes \bar{a}$

To clarify that the result can be different depending on the simple objects, we add an index a for each simple object a:

$c_a = (\lambda\Phi\Lambda)_a^2$

Now it only remains to construct the sum over all these objects.

In chapter 8.3 we have represented the sequence in general which we calculate now for the examples. For each example we chose different start roottrees, such for each subtree with root x and y-z fork:
$Hom(x, y \otimes z) \neq 0$

## 8.5 Computations – example 1

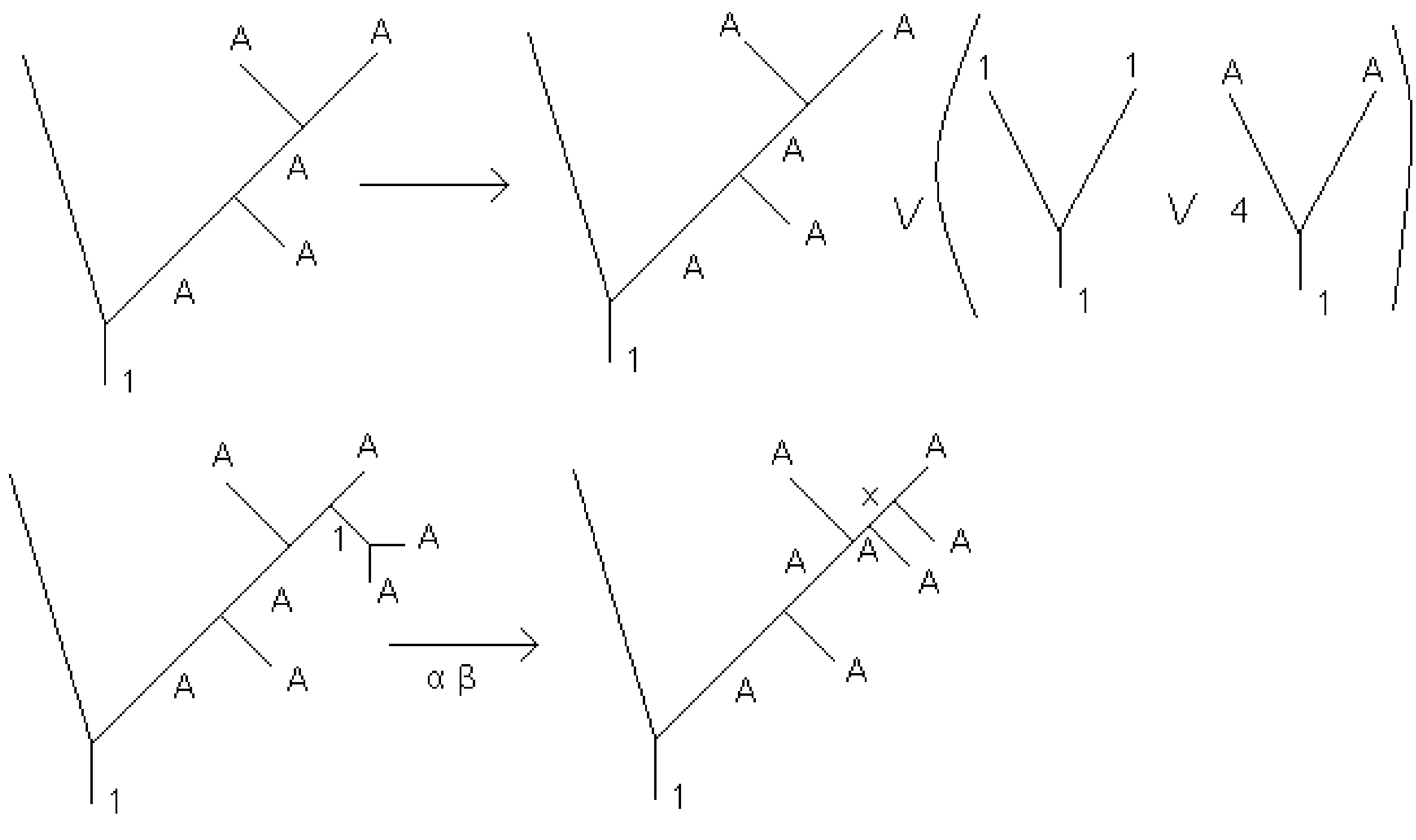

**Figure 265- Computations– example 1- trace unit**

In the first row we apply the trace unit to the start roottree, in the second row we only glue the A-A fork onto the roottree and then we apply associativity. x depends on the resulting linear combination of roottrees with coefficients α, β, that we will determine now, by considering the part, where the roottree changes:

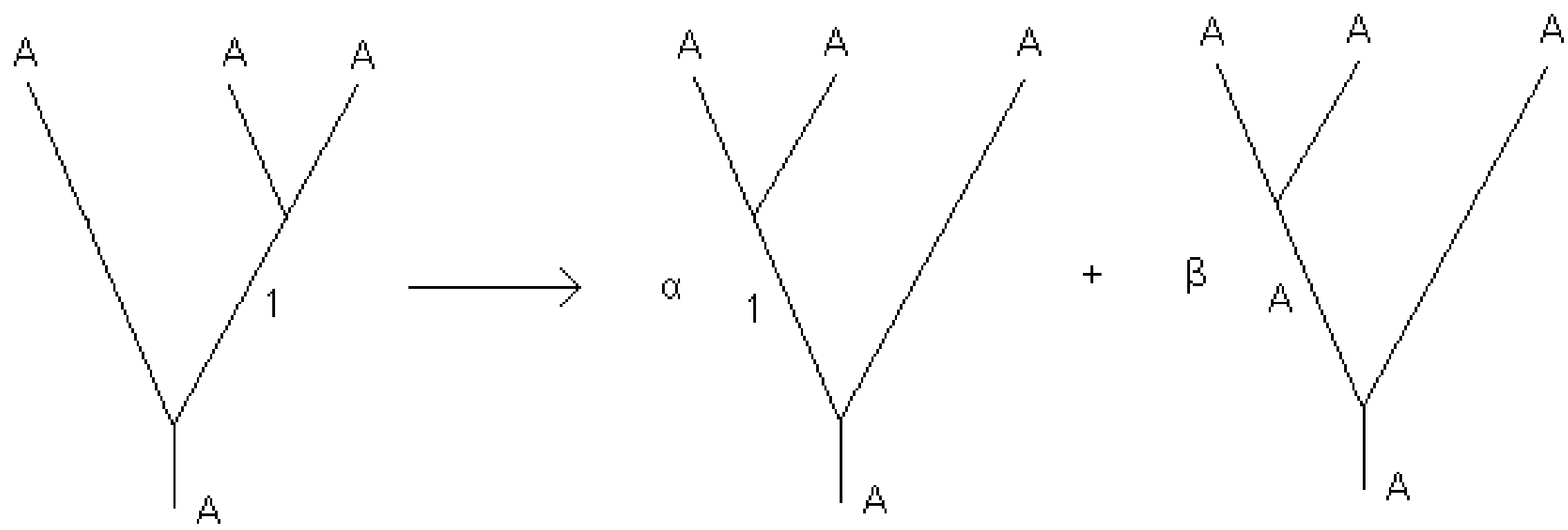

**Figure 266- Computations- example 1- determine associativity 1**

To determine the coefficients α, β, we observe the change of the 1 object in the left roottree under associativity in the right roottrees. We double underscore the relevant objects:

$A \otimes (A \otimes A) \longrightarrow (A \otimes A) \otimes A \qquad (A \otimes A) \otimes A$

$A \otimes (\underline{\underline{1}} \oplus A) \longrightarrow (\underline{\underline{1}} \oplus A) \otimes A \qquad (1 \oplus \underline{\underline{A}}) \otimes A$

$\underline{\underline{A}} \oplus 1 \oplus A \longrightarrow \underline{\underline{A}} \oplus 1 \oplus A \qquad A \oplus 1 \oplus \underline{\underline{A}}$

The associativity matrix tells us, that α = 2 and β = 3:

| | | |
|---|---|---|
| **2** | 0 | 4 |
| 0 | 1 | 0 |
| **3** | 0 | 3 |

So we get a linear combination of roottrees, where we consider each roottree separately:

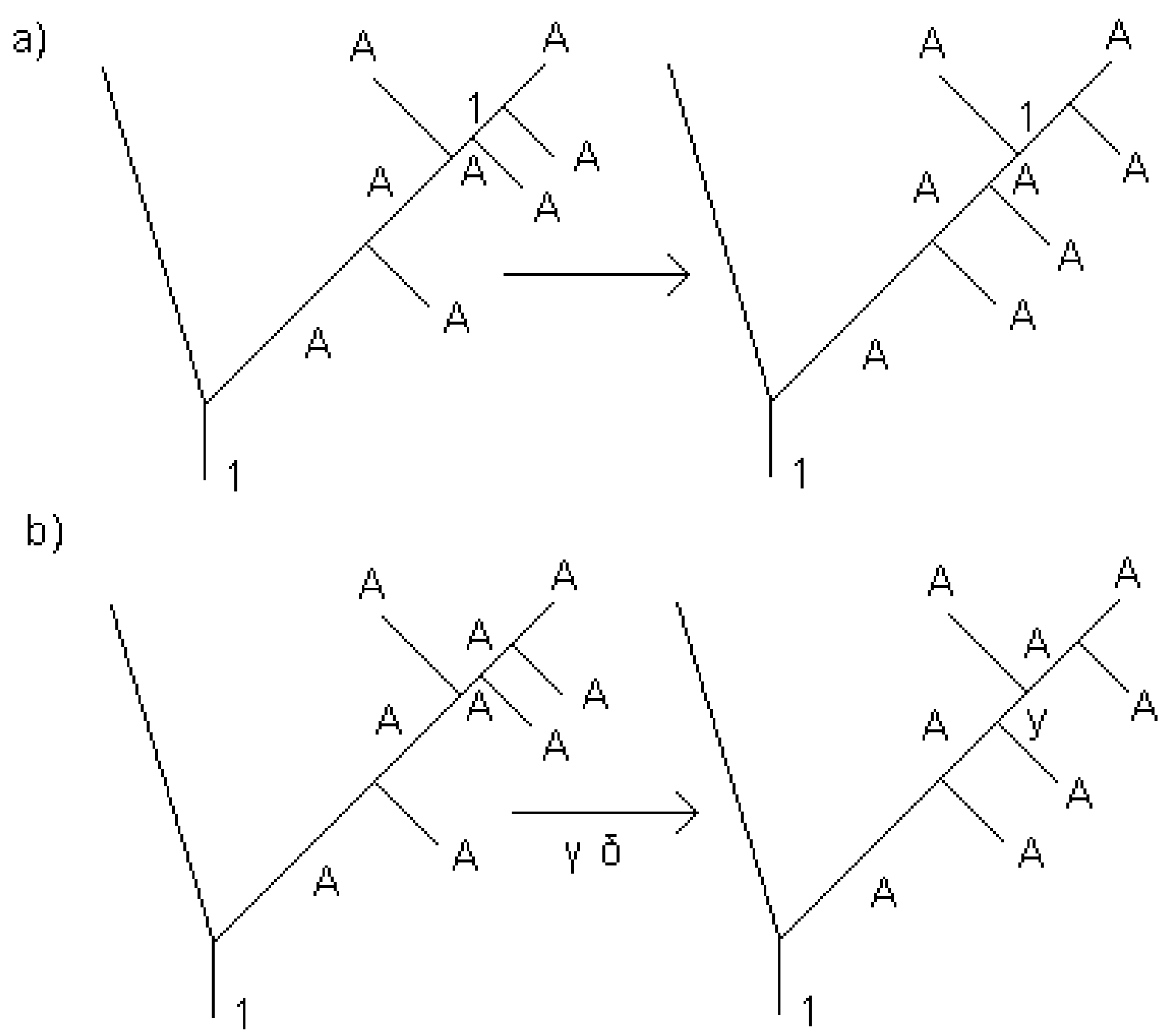

**Figure 267- Computations- example 1- determine associativity 2**

In case a) we get a unique value A  for the variable of the resulting roottree, however, in case b) we get a linear combination with variable y and coefficients γ and δ.

The change of the subtree is:

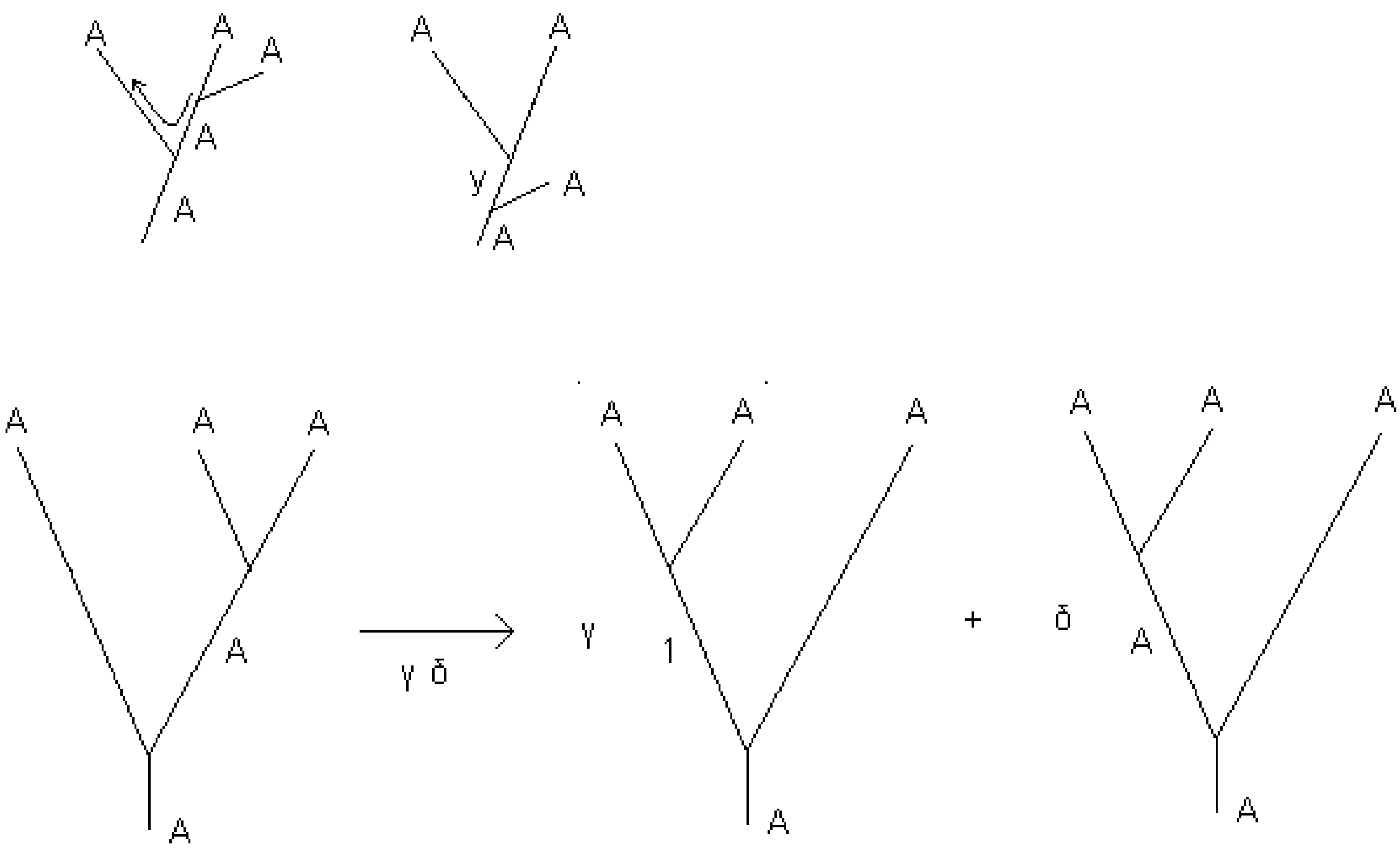

**Figure 268- Computations- example 1- determine associativity 2b**

We have to follow the A object in the left roottree under associativity:

$A \otimes (A \otimes A) \longrightarrow (A \otimes A) \otimes A \qquad (A \otimes A) \otimes A$

$A \otimes (1 \oplus \underline{\underline{A}}) \longrightarrow (\underline{\underline{1}} \oplus A) \otimes A \qquad (1 \oplus \underline{\underline{A}}) \otimes A$

$A \oplus 1 \oplus \underline{\underline{A}} \longrightarrow \underline{\underline{A}} \oplus 1 \oplus A \qquad A \oplus 1 \oplus \underline{\underline{A}}$

The associativity matrix tells us, that $\gamma = 4$ and $\delta = 3$:

| | | |
|---|---|---|
| 2 | 0 | **4** |
| 0 | 1 | 0 |
| 3 | 0 | **3** |

So after the second associativity we have 3 different roottrees, and apply the next move to prepare the splitting:

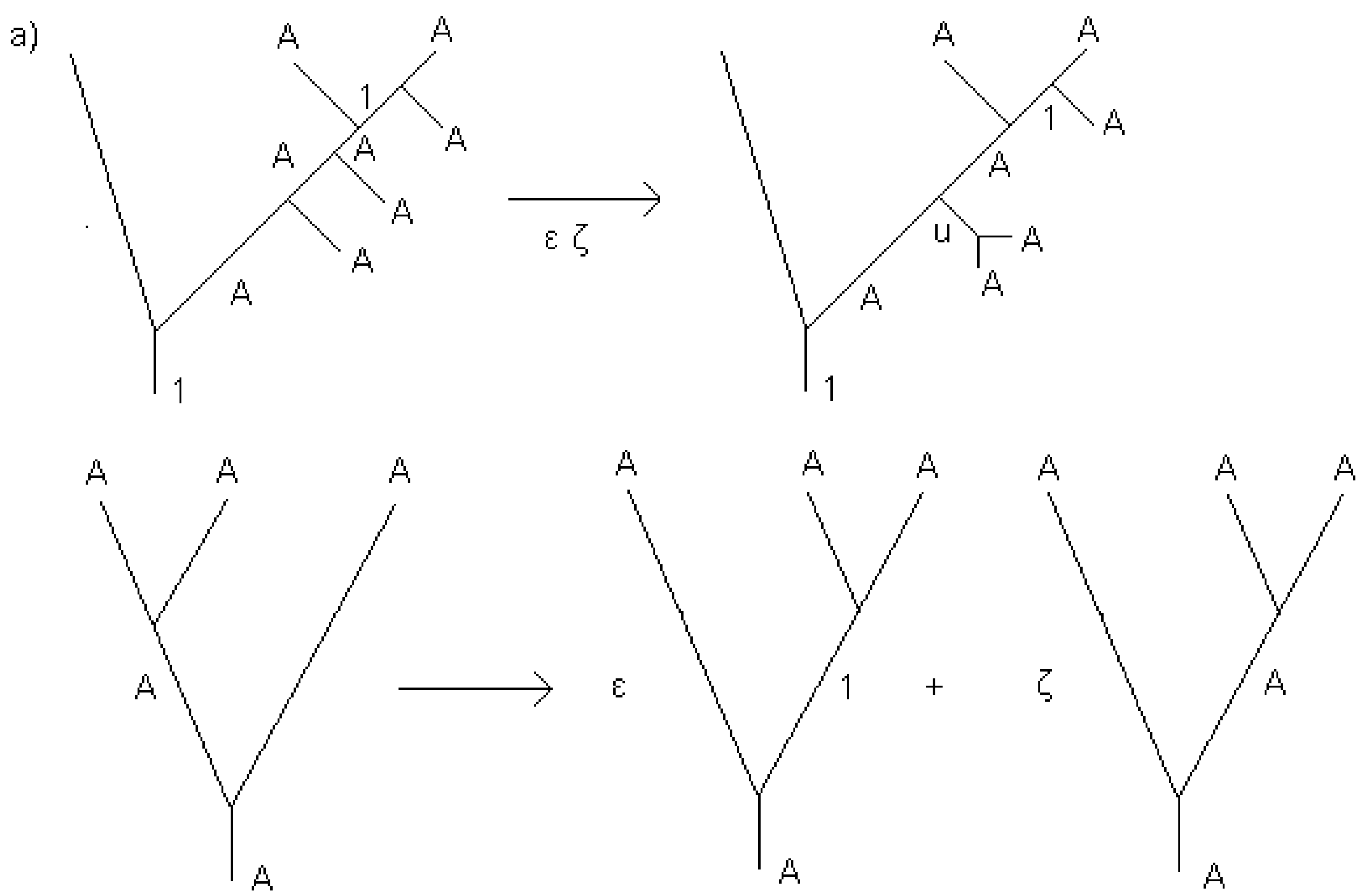

**Figure 269- Computations- example 1- determine associativity 3a**

We have to follow the A object in the left roottree under associativity. Note, that the matrix is inverse to itself:

$(A \otimes A) \otimes A \longrightarrow A \otimes (A \otimes A) \qquad A \otimes (A \otimes A)$

$(1 \oplus \underline{\underline{A}}) \otimes A \longrightarrow A \otimes (\underline{\underline{1}} \oplus A) \qquad A \otimes (\underline{\underline{1}} \oplus A)$

$A \oplus 1 \oplus \underline{\underline{A}} \longrightarrow \underline{\underline{A}} \oplus 1 \oplus A \qquad A \oplus 1 \oplus \underline{\underline{A}}$

The associativity matrix tells us, that us ε = 4 and ζ = 3.

| | | |
|---|---|---|
| 2 | 0 | **4** |
| 0 | 1 | 0 |
| 3 | 0 | **3** |

The next case b1) provides a unique simple object and the other case b2) a linear combination:

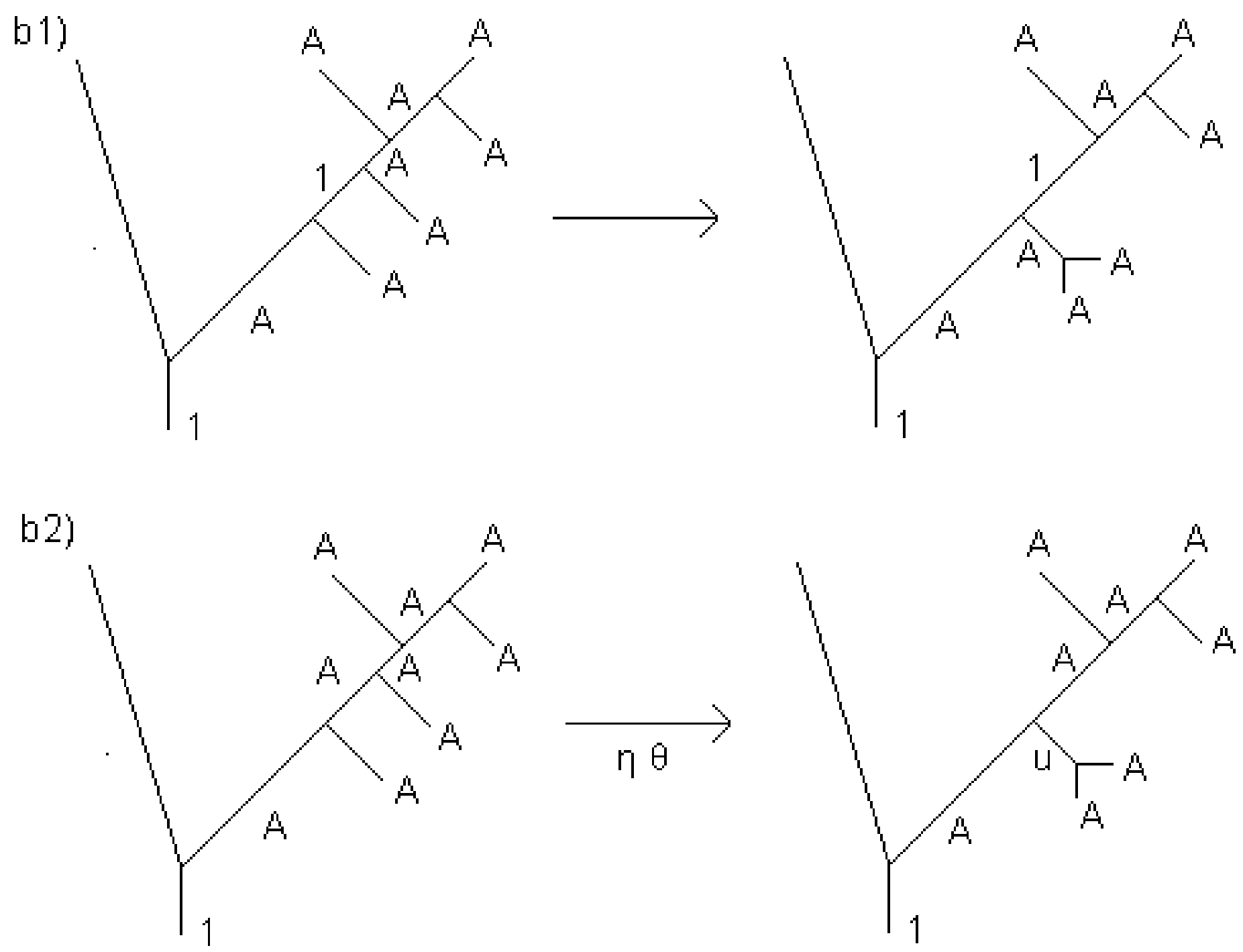

**Figure 270- Computations- example 1- determine associativity 3b**

The computation of case b2) equals case a) according to which
u = 1 → η = 4 and u = A → θ = 3

We collect the cases:

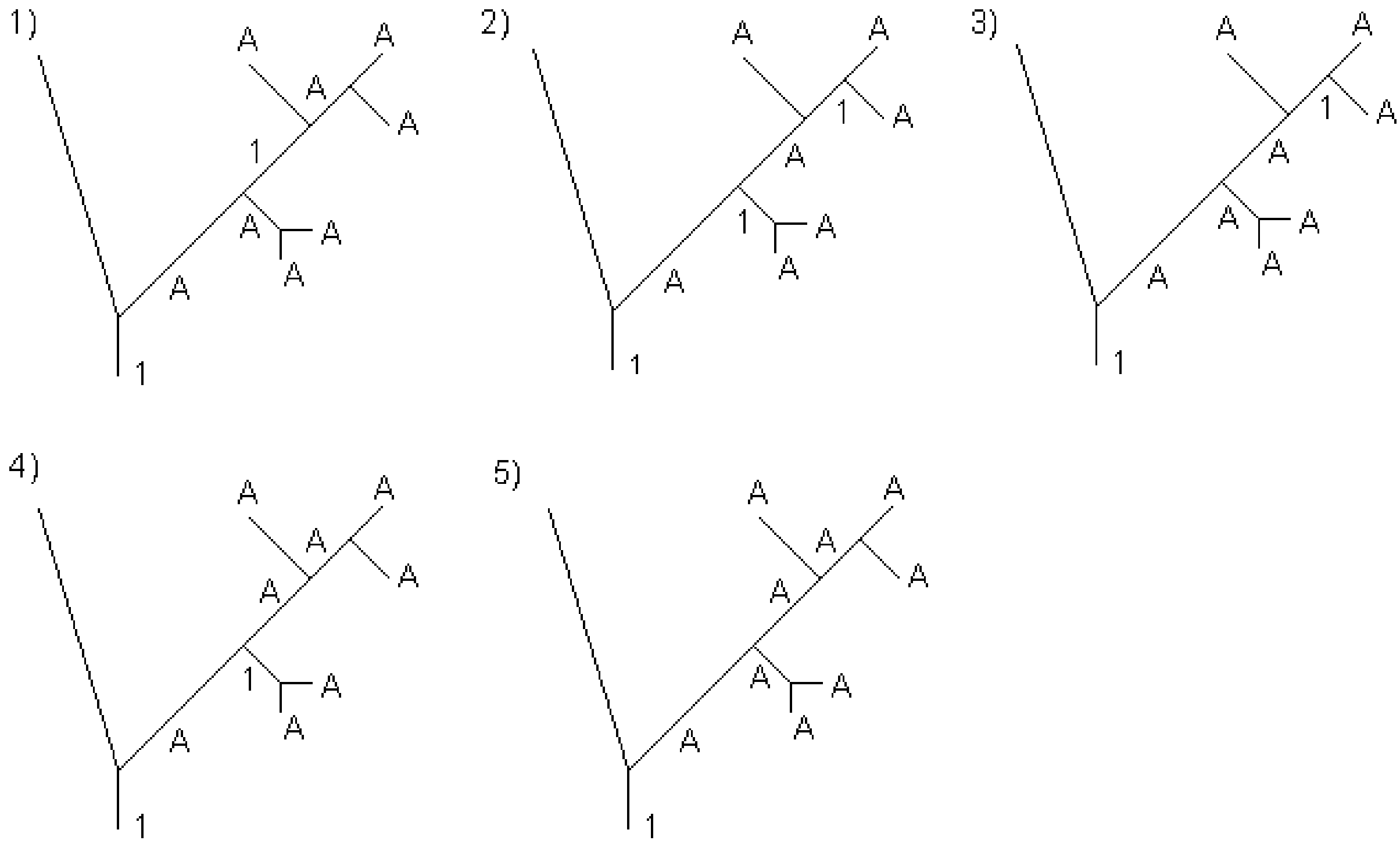

**Figure 271- Computations- example 1- collect cases before split**

The general introduction in chapter 8.3 provides that only splitting at 1 makes sense, so only 2) and 4) remain:

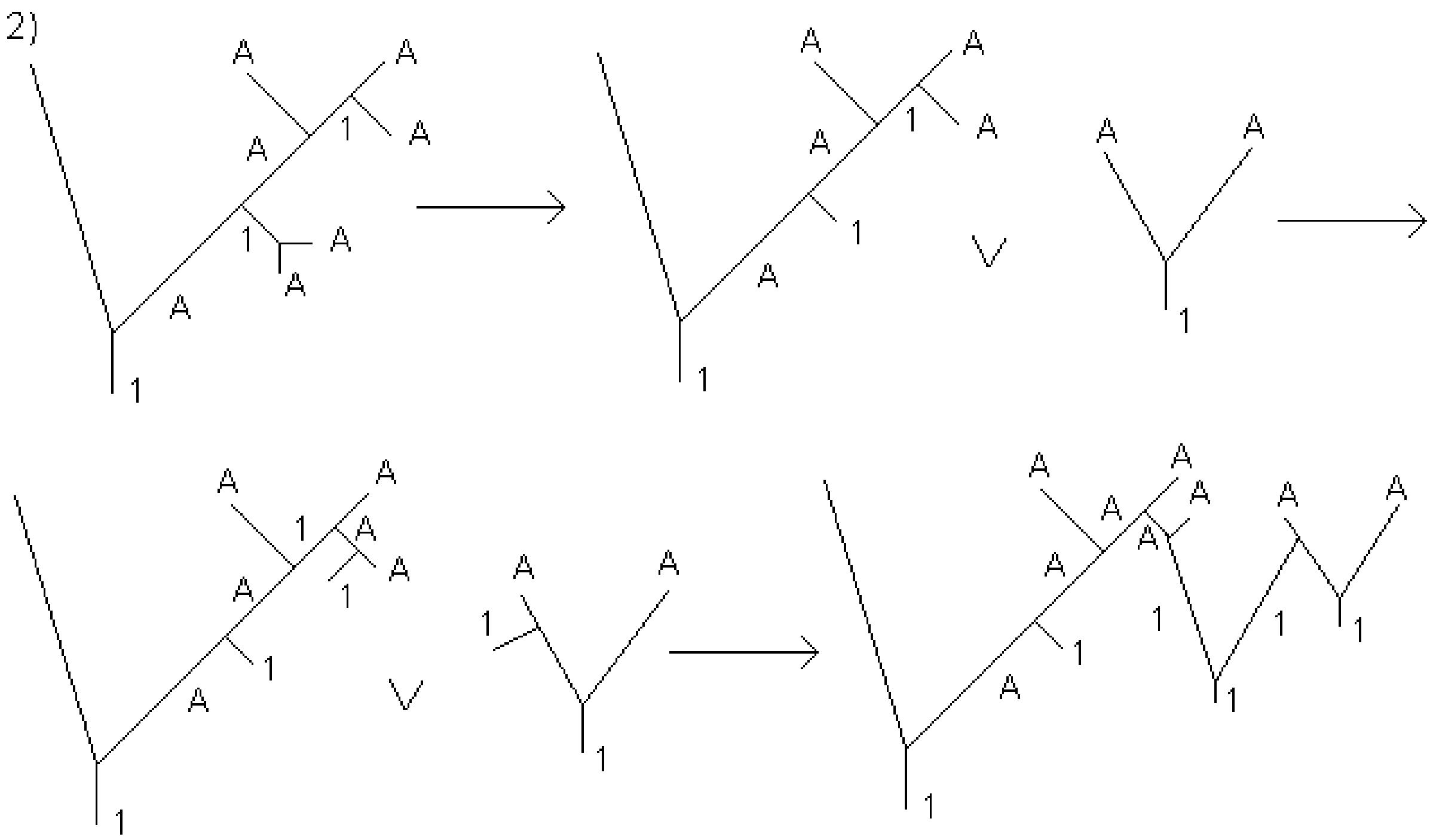

**Figure 272- Computations- example 1- case 2- glue splitted roottrees**

We split the roottree and glue both parts together at the boundary points corresponding to the trace unit, and do the same for case 4):

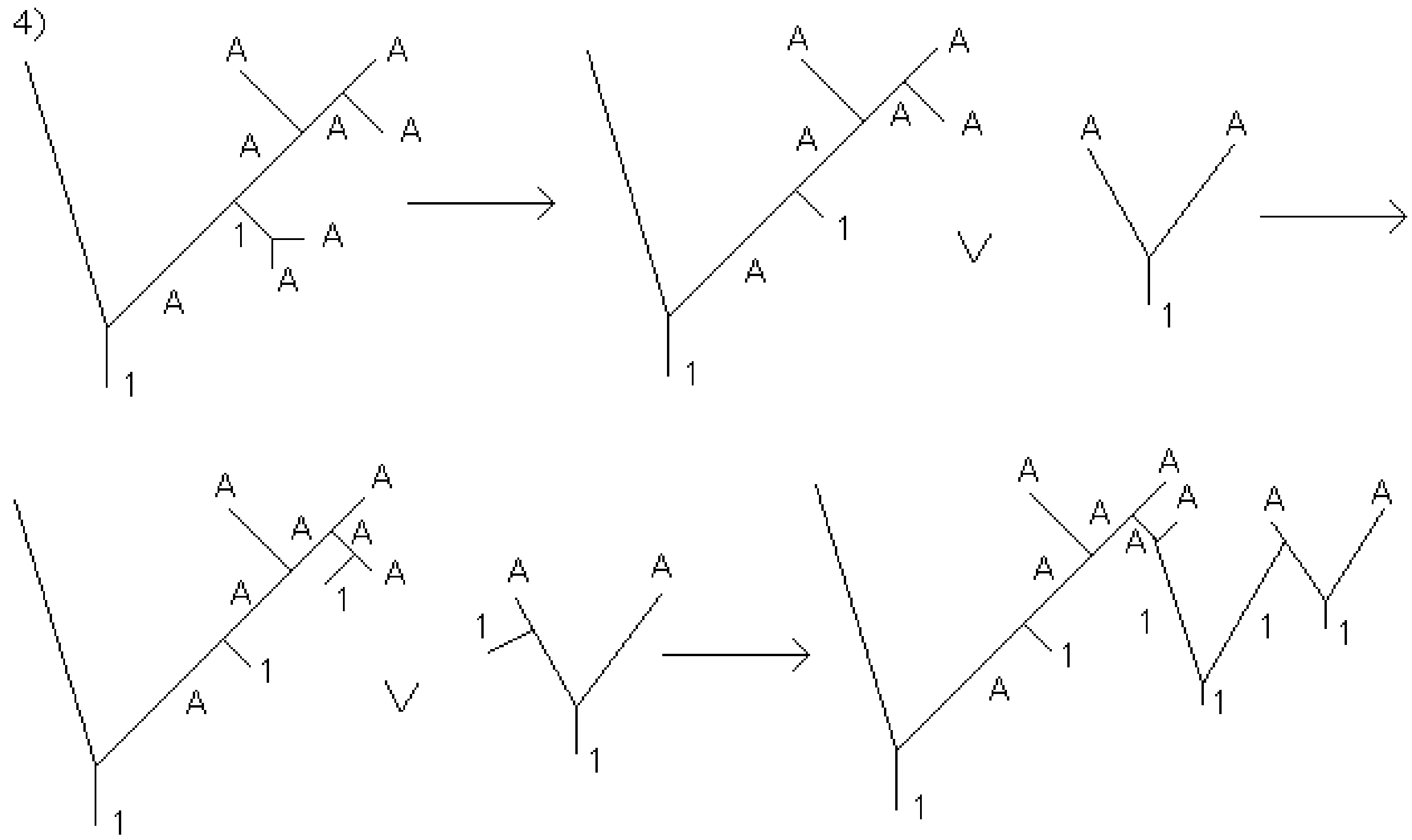

**Figure 273- Computations- example 1- case 4- glue splitted roottrees**

We finish case 2):
First we move one A branch onto another by associativity with variable object v and then apply the form $\lambda$ : v → 1.
This is nonzero, if v = 1 and zero if v = A, since the objects are in a semisimple tensor category:

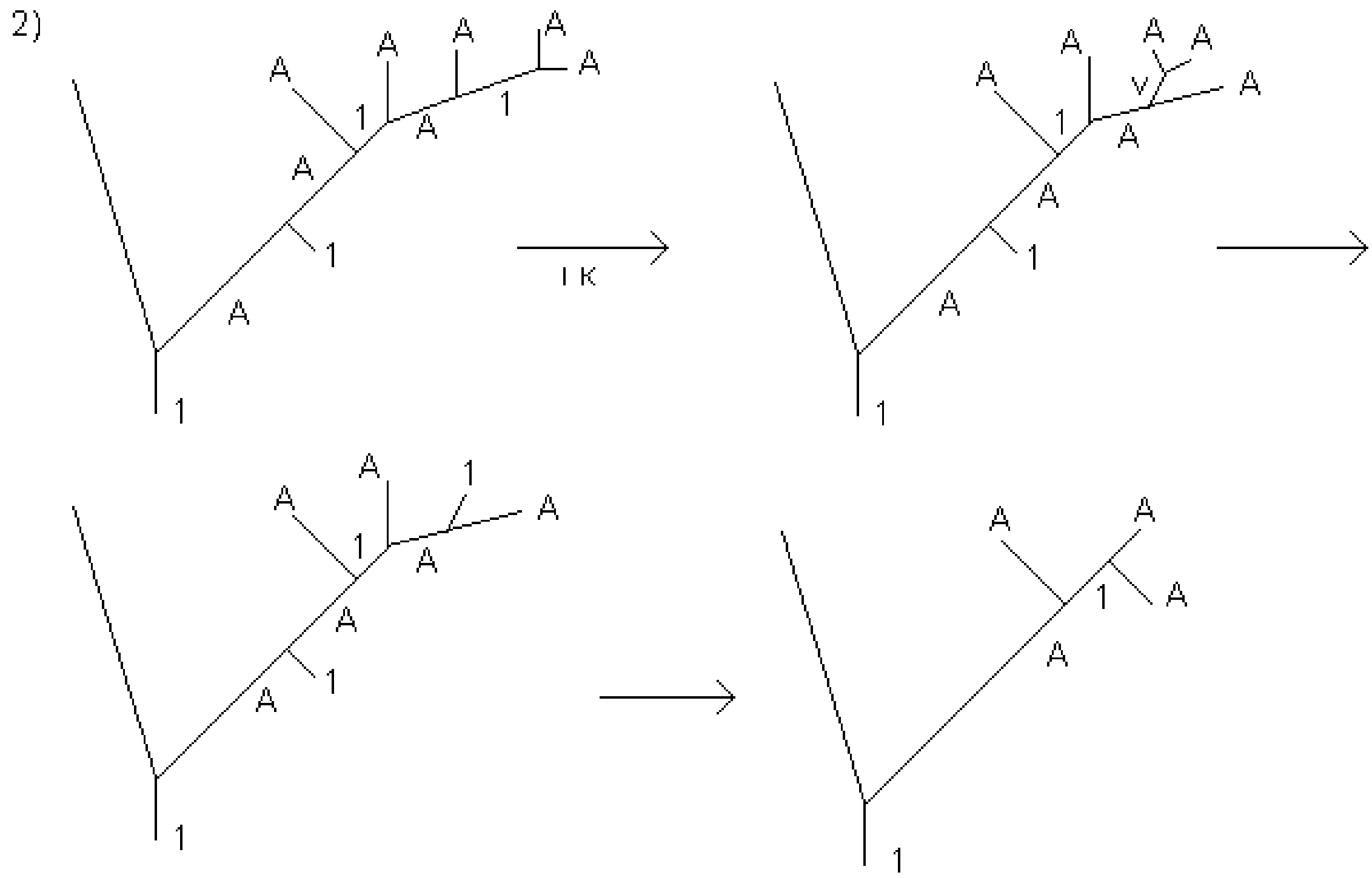

**Figure 274- Computations- example 1- case 2- finish**

We determine the coefficients ι and κ by repeating our computation, supported by:

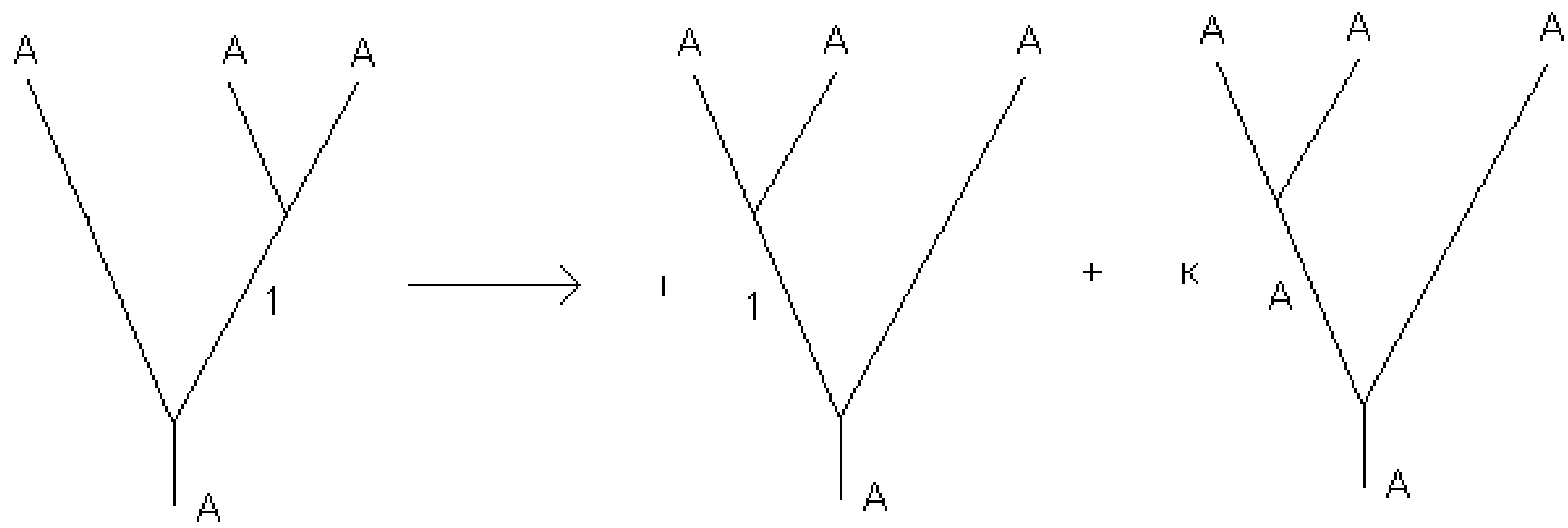

**Figure 275- Computations- example 1- case 2- associativity before apply the form**

By former computations we get ı = 2 and κ = 3. Note that only ı = 2 is relevant.

Similar we finish case 4:

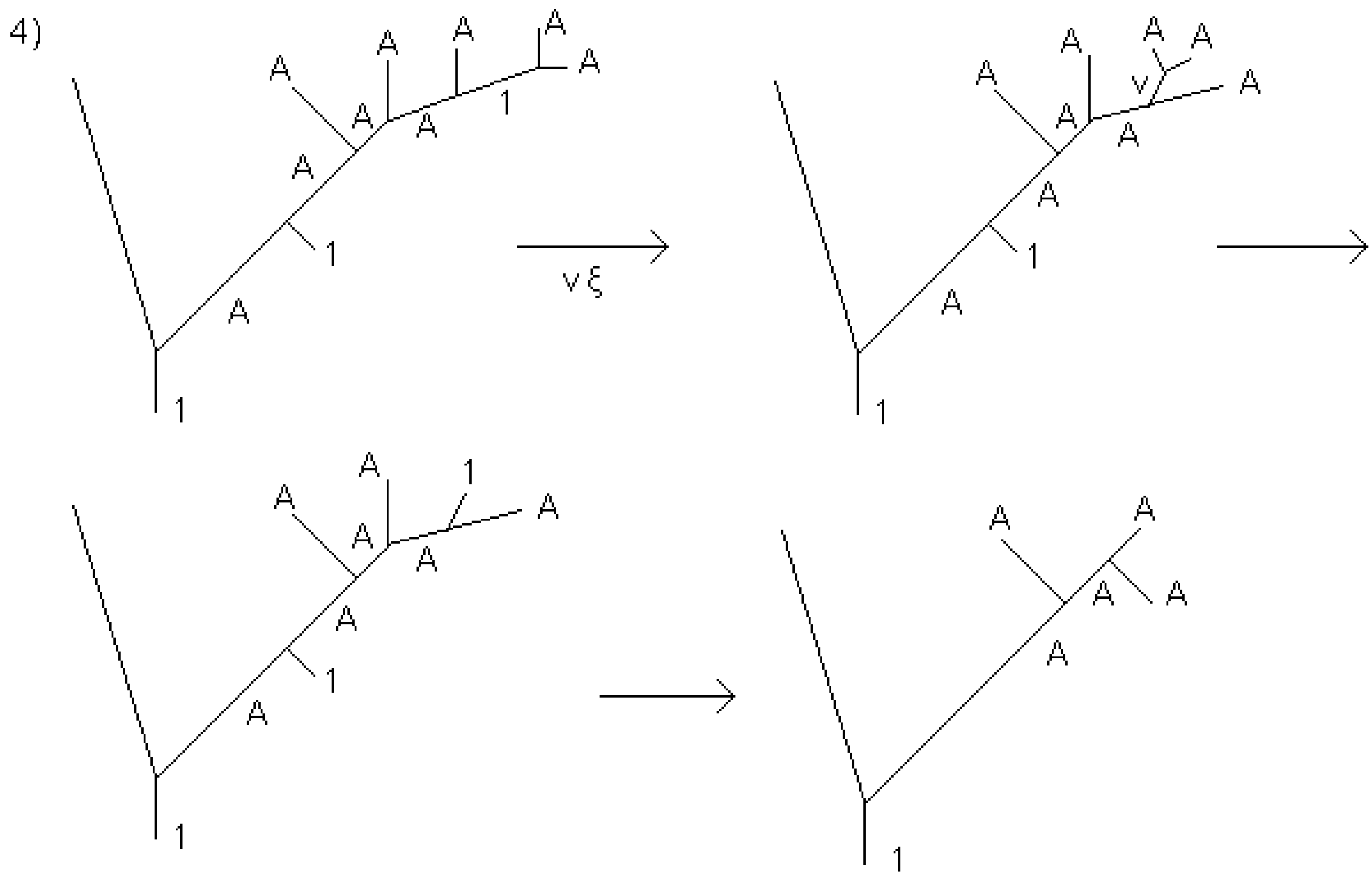

**Figure 276- Computations- example 1- case 4- finish**

As in case 2) v = 2 and ξ = 3.

We computed all details, but it remains to determine the coefficients of the end roottrees. A graphical overview of the relevant steps is useful:

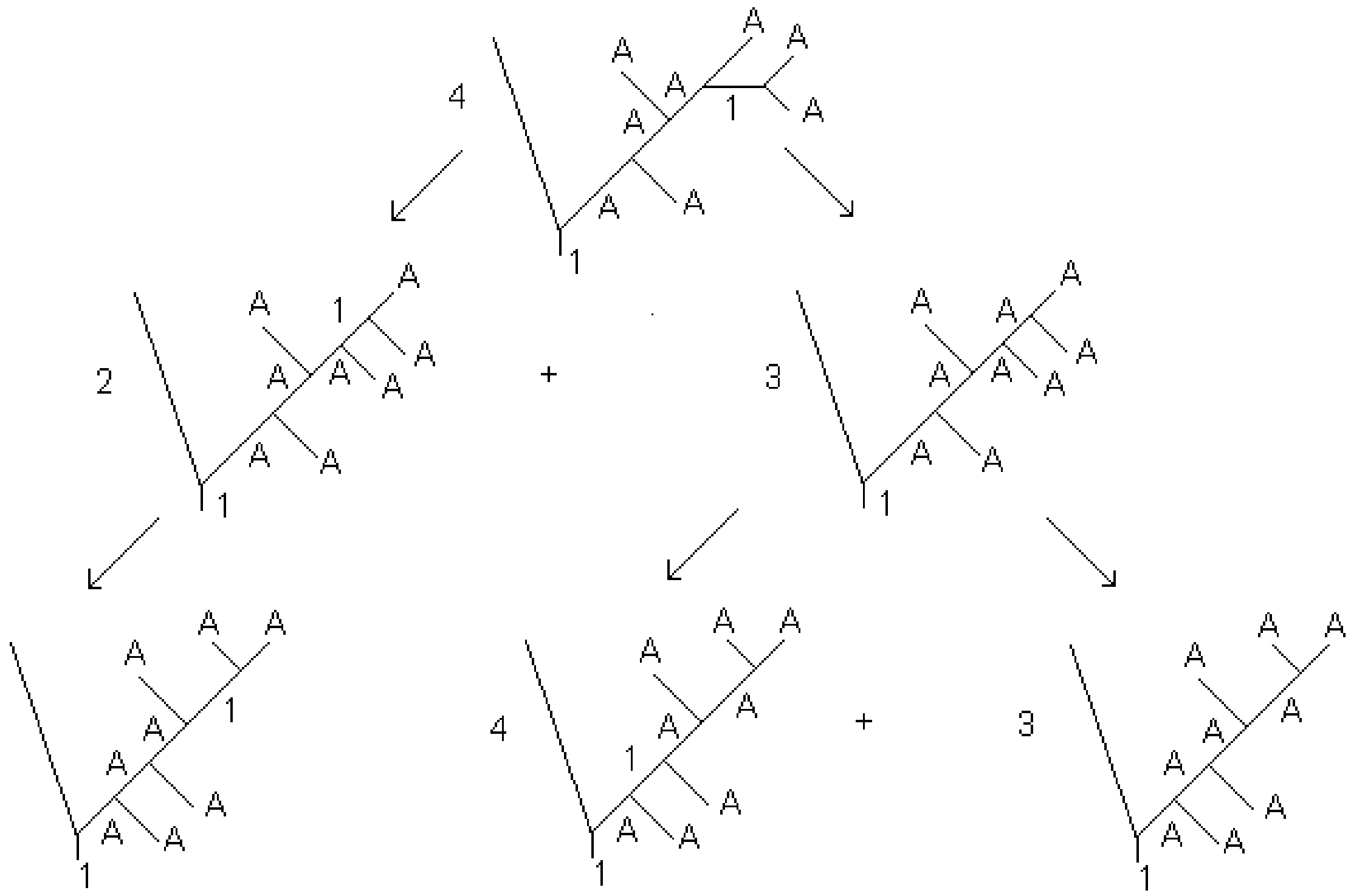

**Figure 277- Computations- example 1- overview 1**

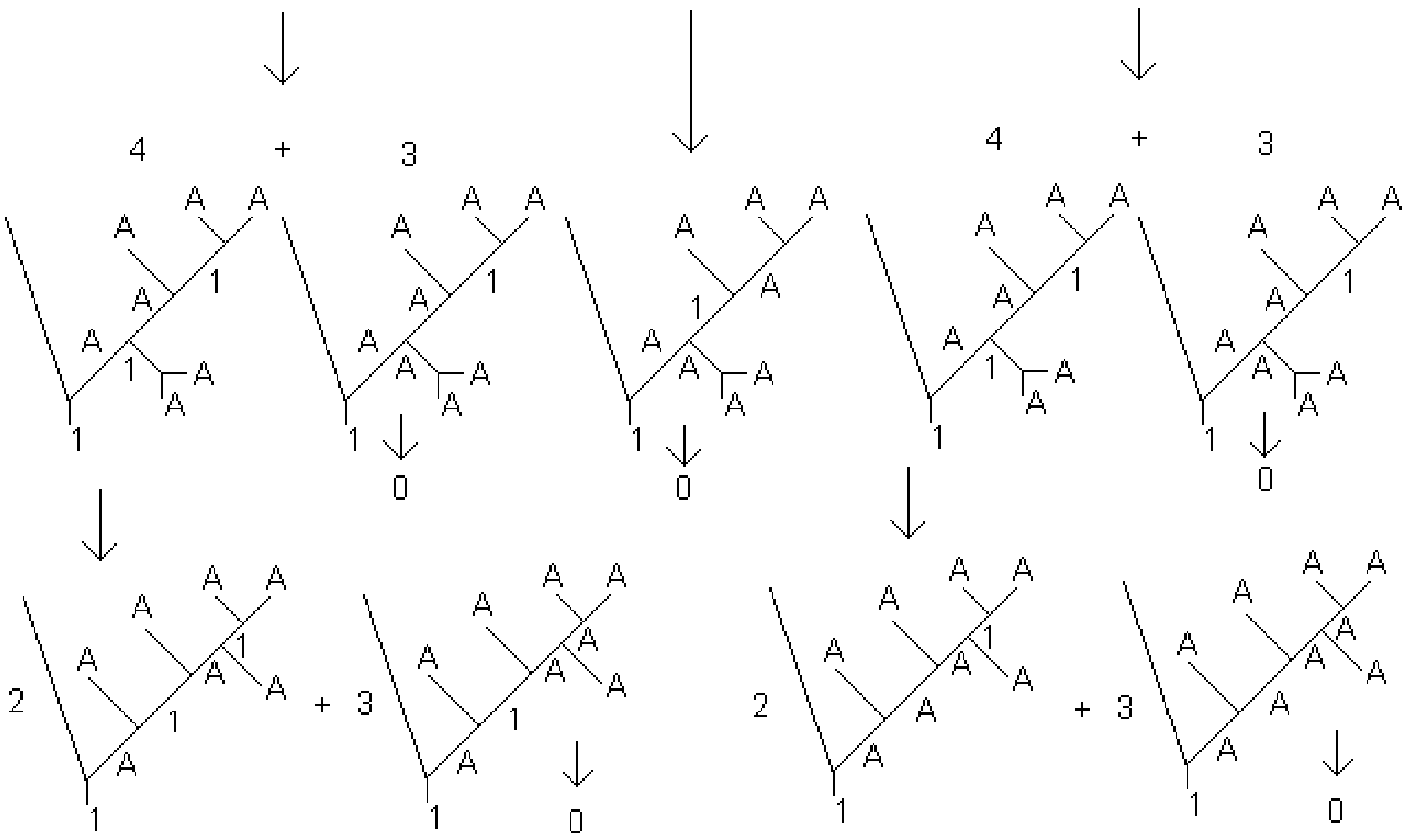

**Figure 278- Computations- example 1- overview 2**

The last row describes the step before we apply the form $\lambda$, after these we get no further coefficients and the second and last roottree will be annihilated after this step. We get:

For the first end roottree:

$4 \cdot 2 \cdot 4 \cdot 2 = 4 \mod Z_5$

for the second end roottree:
$4 \cdot 3 \cdot 3 \cdot 4 \cdot 2 = 3 \bmod Z_5$

At the start of the computation, we restrict the trace unit to the A-A fork, since if we take a 1-1 fork for the trace unit, after splitting we get 2 trees as in the next figure. However the state module of the second splitted tree is
Hom ( 1,1⊗A) = Hom (1,A) = 0, since the tensor category is semisimple:

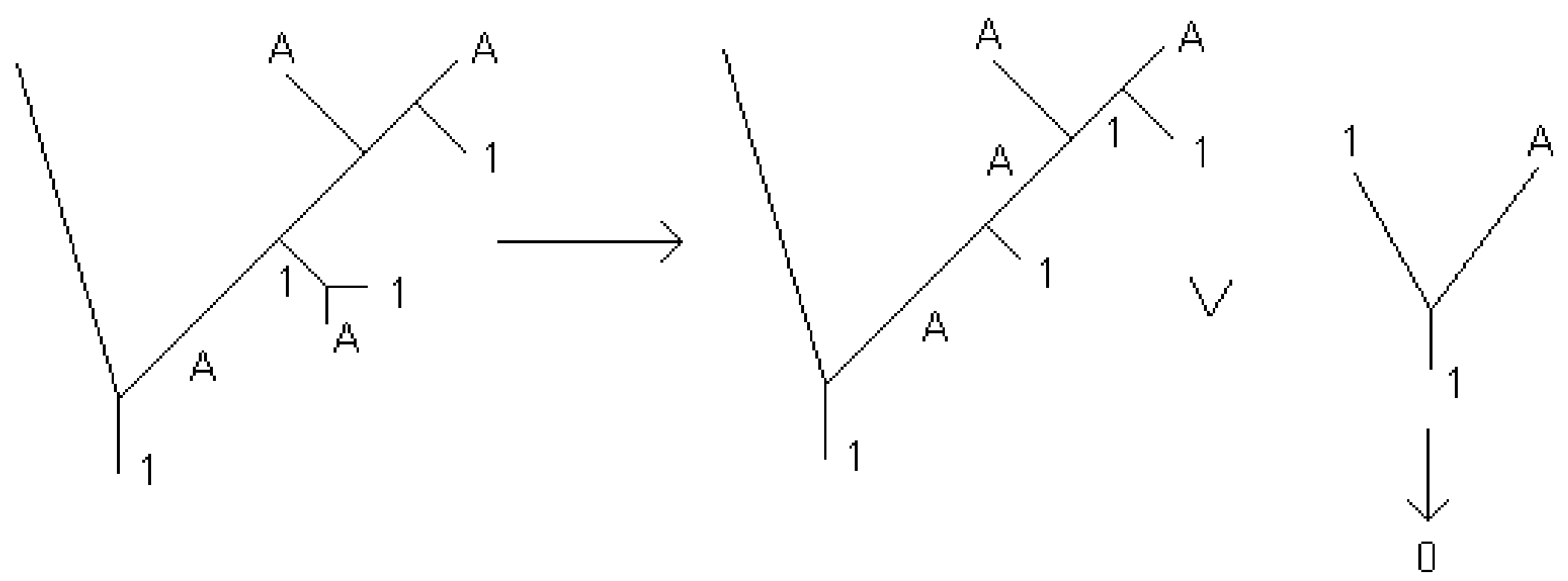

**Figure 279- Computations- example 1- 1-1 fork of trace unit not relevant**

Since disjoint roottrees lead to the tensor product of their associated state modules, the whole roottree can be omitted.

We have to determine the coefficients and roottrees corresponding to the sequence of passing a vertex, but this is only a change under associativity:

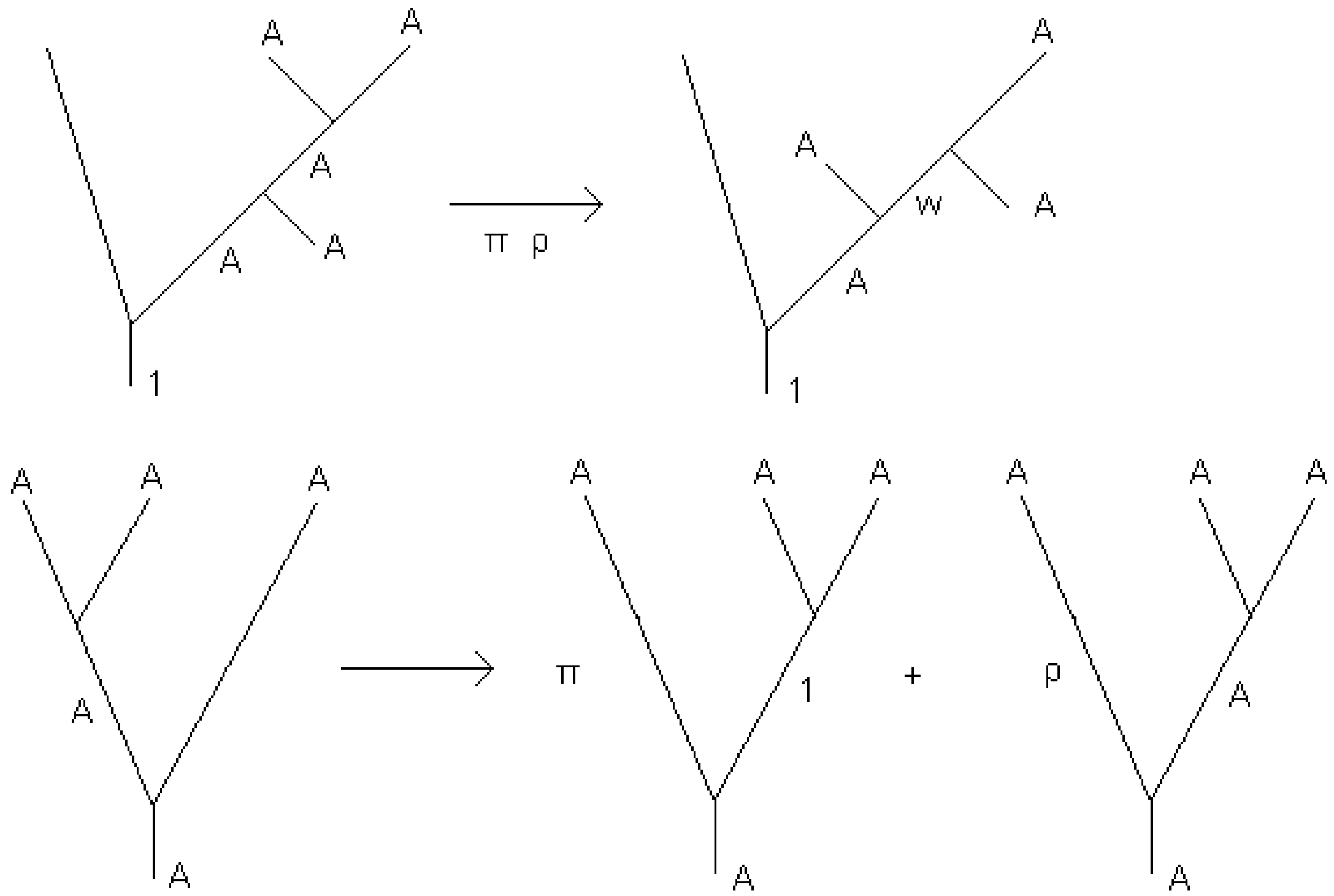

**Figure 280- Computations- example 1- passing a vertex**

From earlier computations we know that π = 4 and ρ = 3.

Hence the chosen start roottree maps to the same linear combination of end roottrees, as in the new sequence !!!

## *8.6 Computations - example 2*

We regard another start roottree, again apply both sequences on it and compare the results. Since there is nothing new in comparison to example 1, we insert only few words to clarify the single step.

We start with the new sequence:

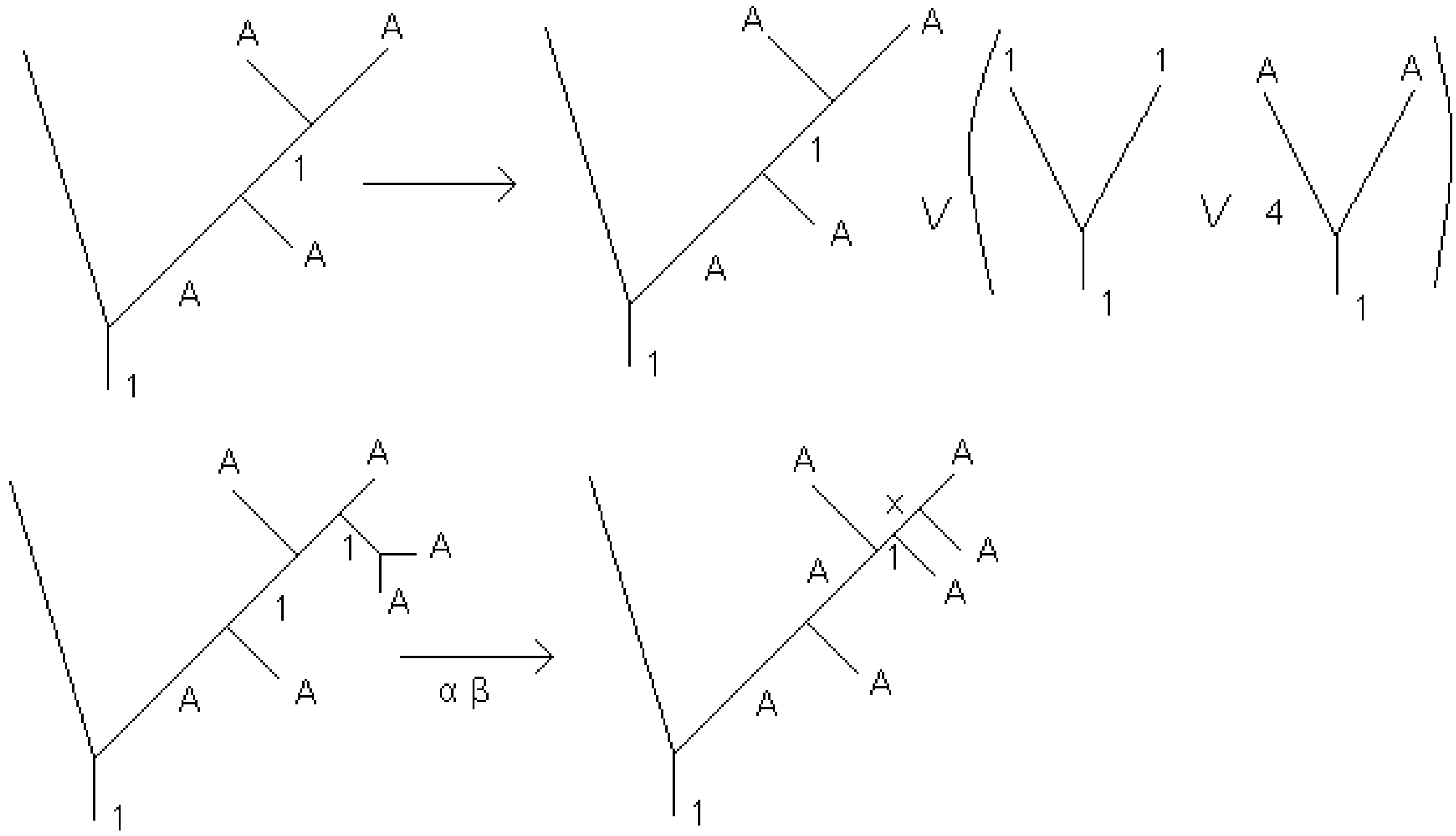

**Figure 281- Computations- example 2- trace unit**

The associativity matrix shows, that α = 2 for x = 1 and β = 3 for x = A.

Remark: As we have shown in example 1, the 1-1 fork of the trace unit can be omitted.

We get 2 cases, one for α and one for β, where we apply the next move:

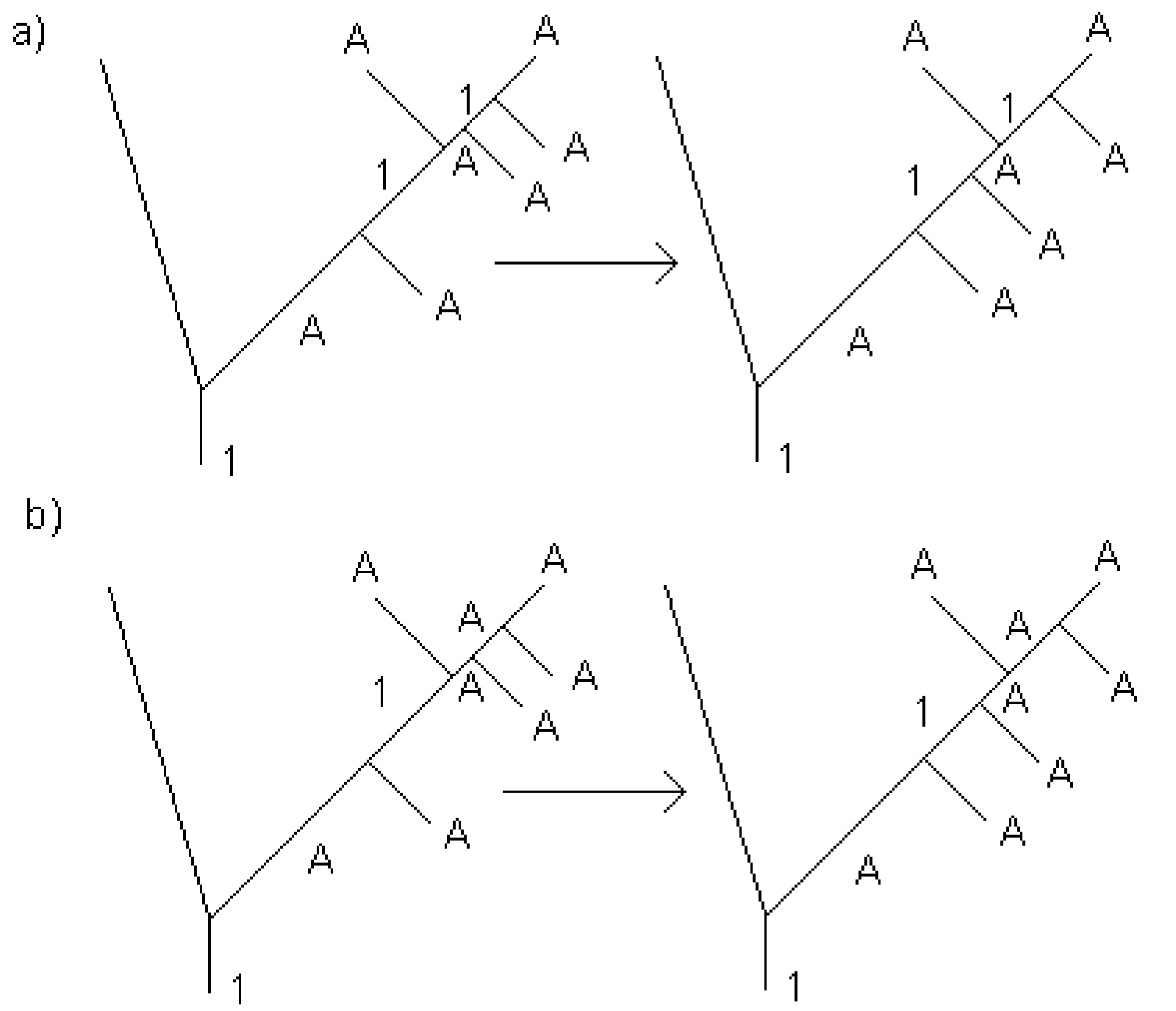

**Figure 282- Computations- example 2- associativity 2**

Both cases lead to a unique determined simple object, because of the 1 object (not the root) in the start roottree. This makes the computation easier than the first one !!! The next move prepares the split of the trees:

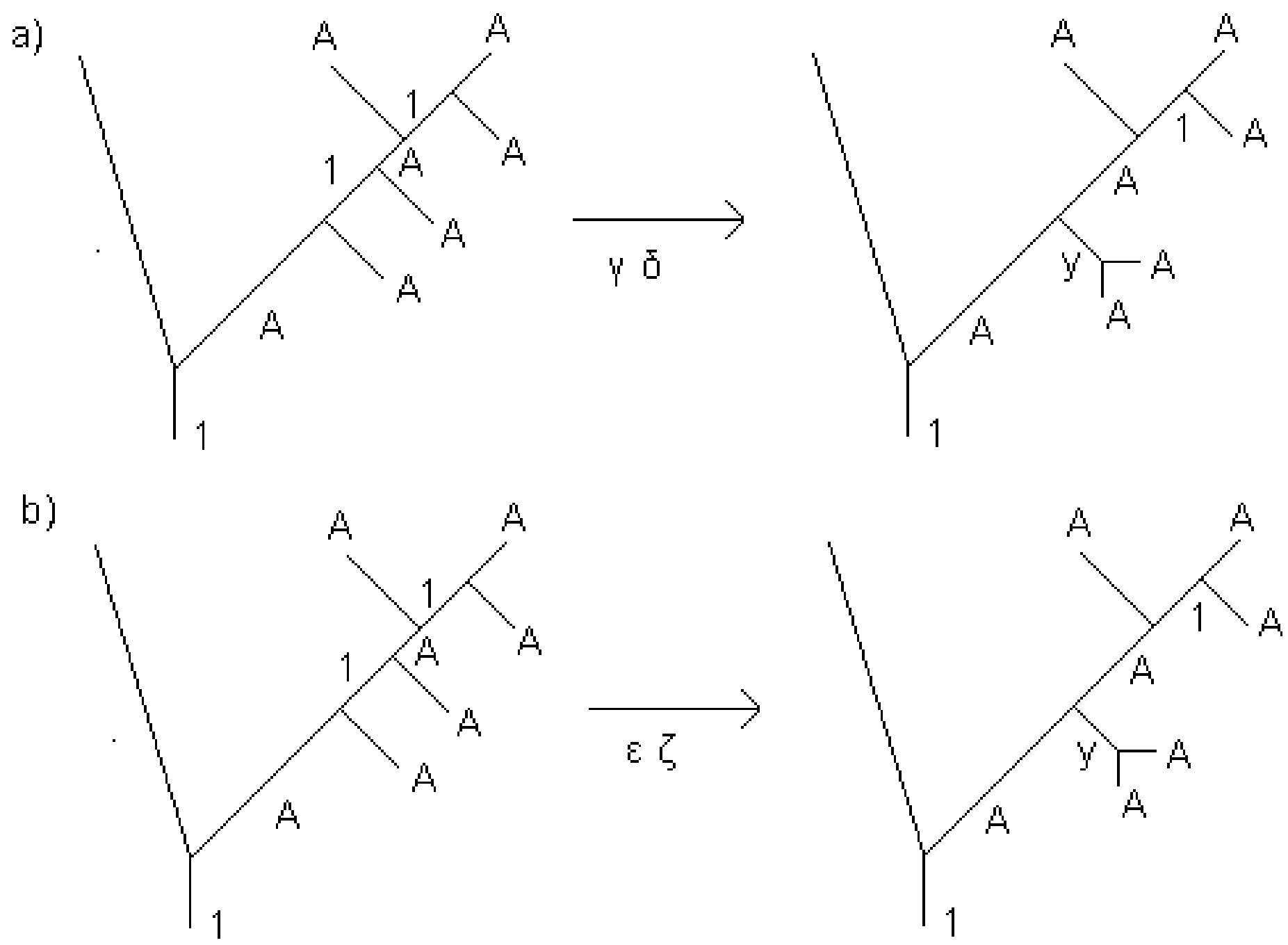

**Figure 283- Computations- example 2- before split**

Since we have the same simple objects on the changed part of the roottree, we get:

For y = 1 → γ = 2 = ε and for y = A → δ = 3 = ζ.

We get 4 cases, but since we just only split at 1, the cases 1) and 3) remain:

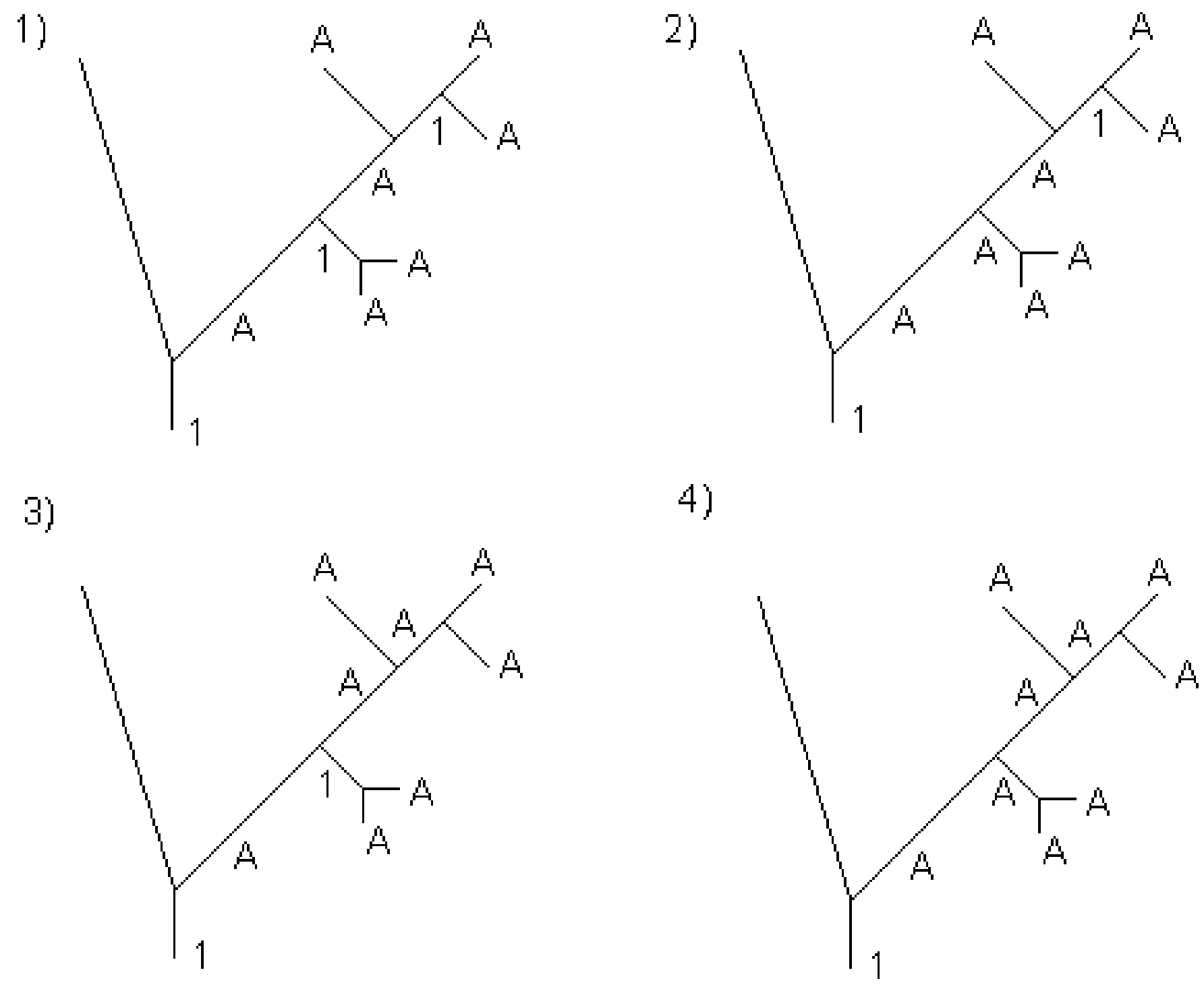

**Figure 284- Computations- example 2- collect the cases**

Continuing the sequence, we want to glue the roottrees together after splitting at the boundary points corresponding to the branches of the trace unit:

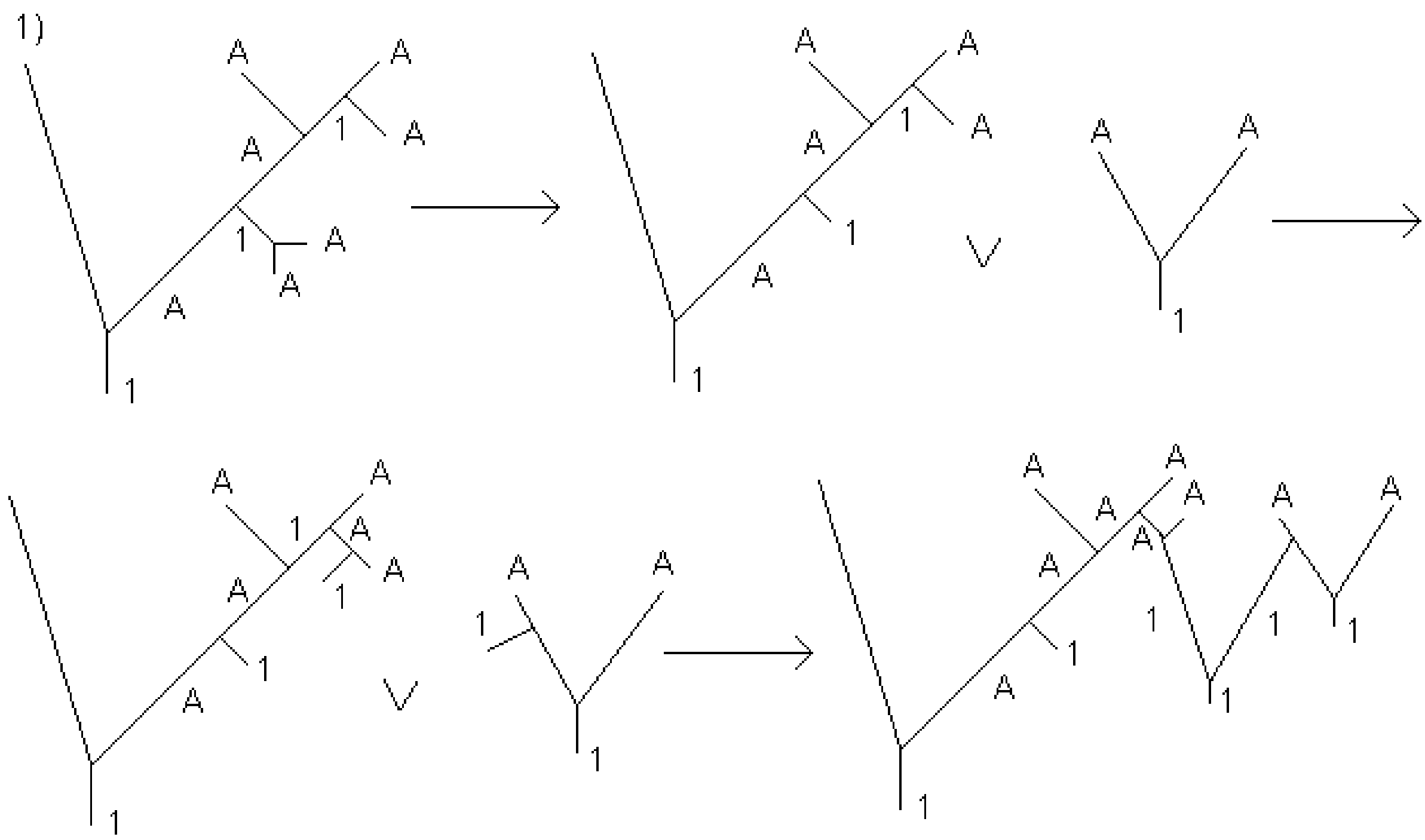

**Figure 285- Computations- example 2- glue splitted roottrees- case 1**

After glueing, we normalize the roottree, apply associativity to move the branches of the trace unit together and apply the form λ on it, which is only nonzero for u = 1:

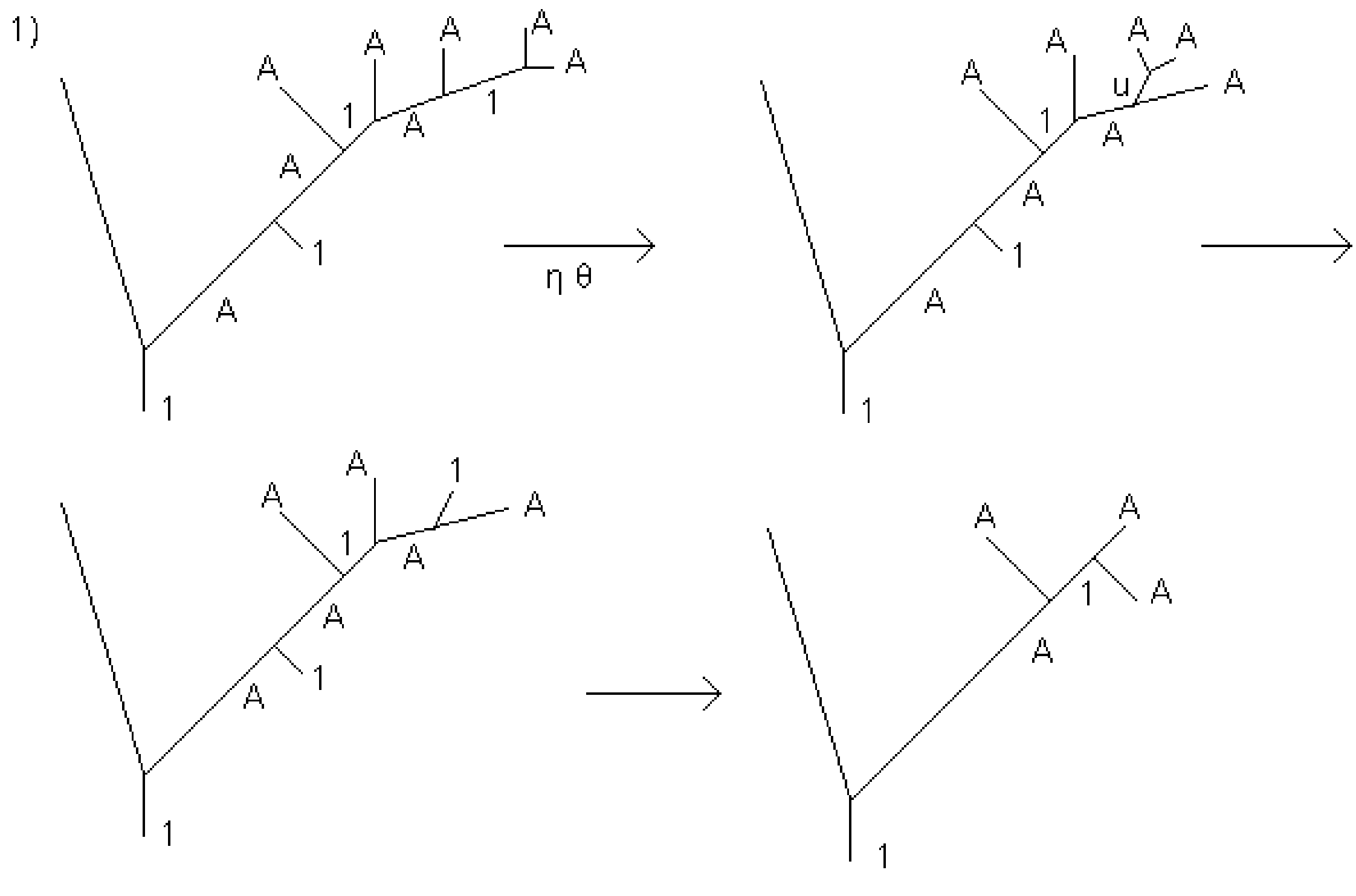

**Figure 286- Computations- example 2- glue splitted roottrees- case 1- apply the form**

The relevant coefficient η is 2.

Repeat the process for the other case:
Glue the splitted trees together:

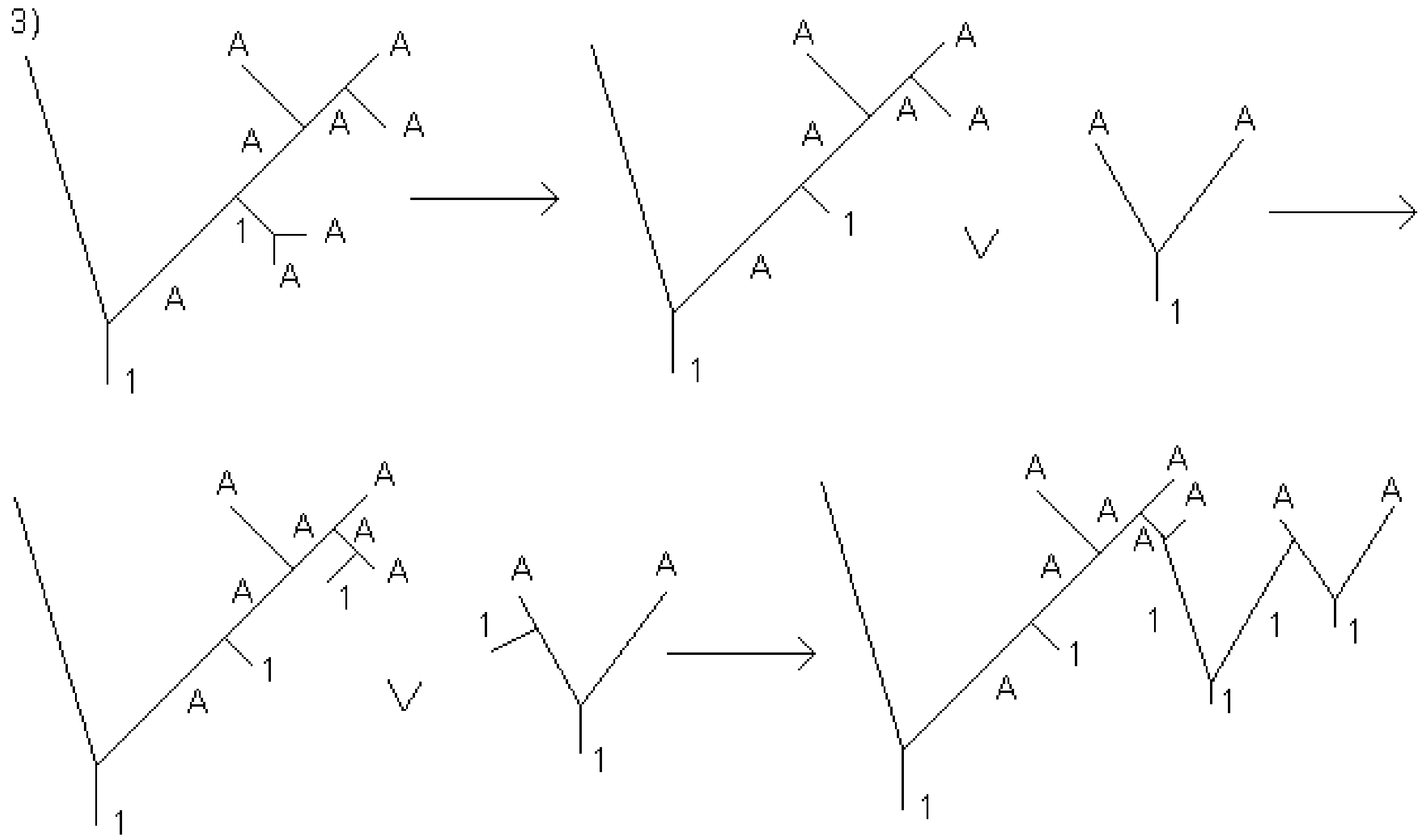

**Figure 287- Computations- example 2- glue splitted roottrees- case 3**

Normalize, apply associativity and the form λ:

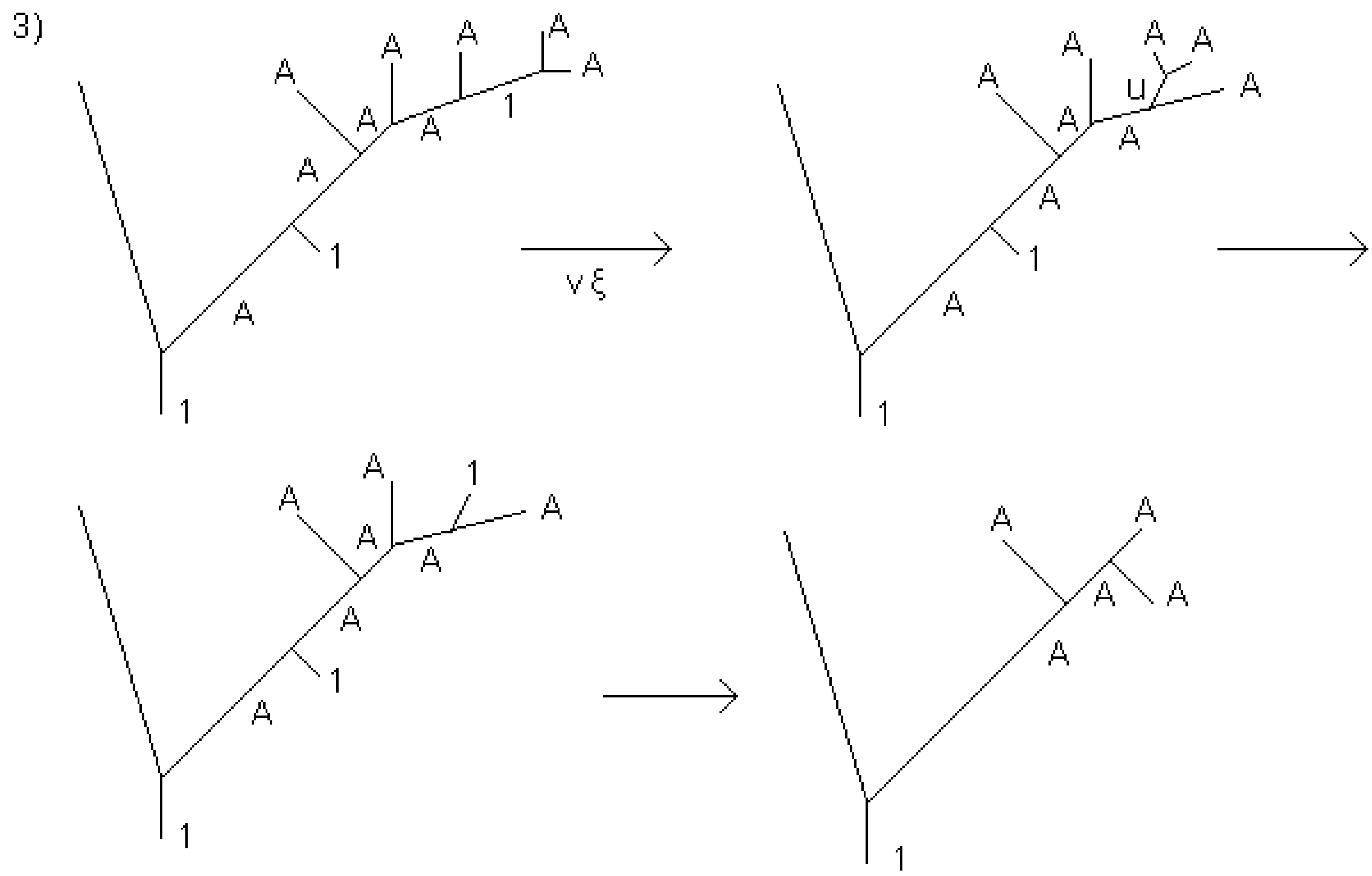

**Figure 288- Computations- example 2- glue splitted roottrees- case 3- apply the form**

The relevant coefficient for u = 1 is v = 2.

We give an overview to determine the linear combination of the end roottrees:

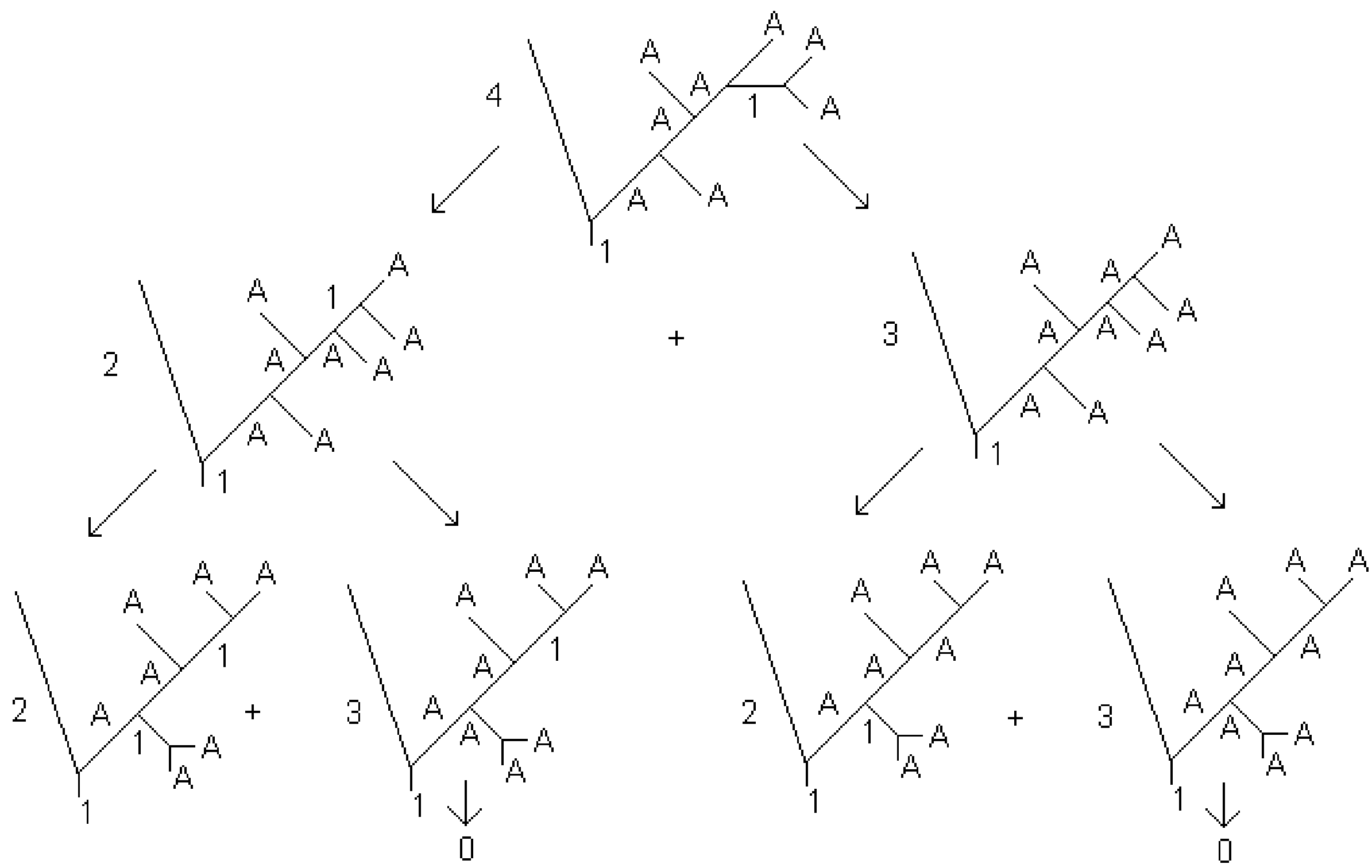

**Figure 289- Computations- example 2- overview 1**

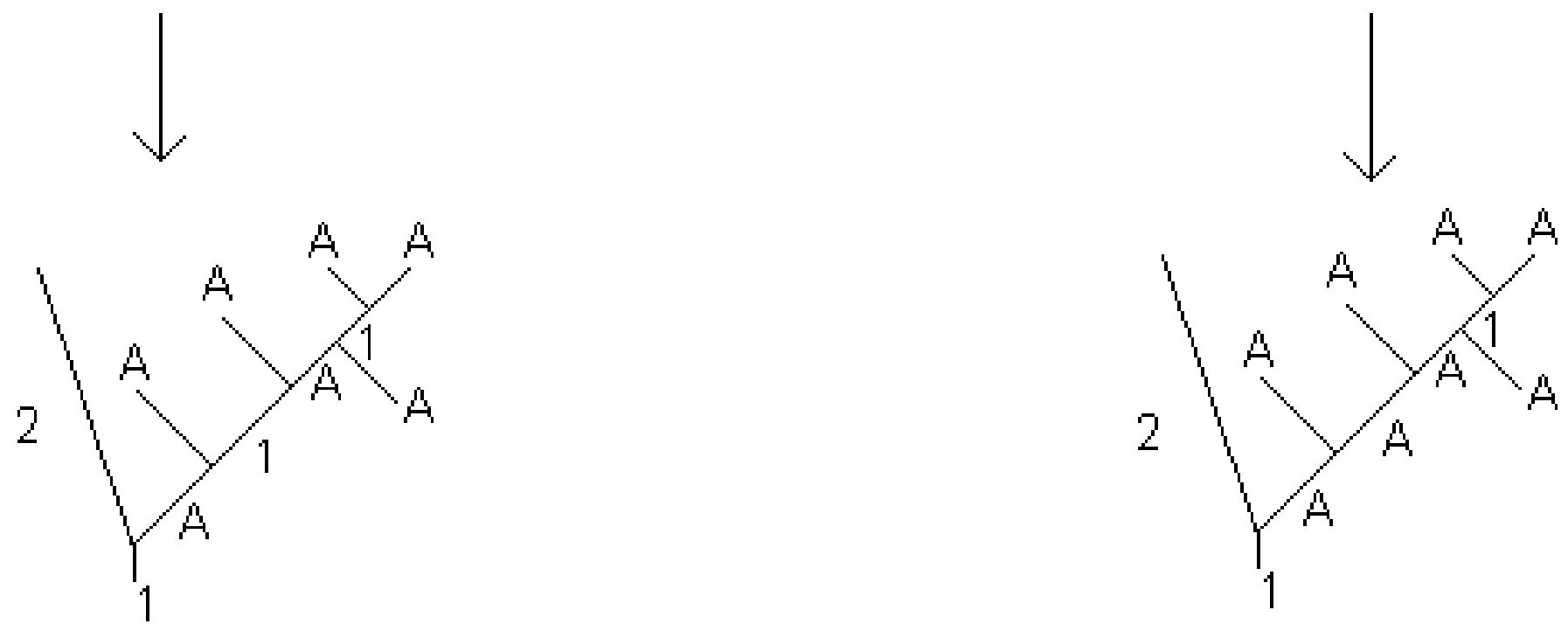

**Figure 290- Computations- example 2- overview 2**

The coefficient of the first end roottree is:
4 •  2 • 2 • 2 = 2 mod Z $_5$

The coefficient of the second end roottree is:
4 •  3 • 2 • 2 = 3 mod Z $_5$

Determine the linear combination of end roottrees repectively for the sequence of passing a vertex:

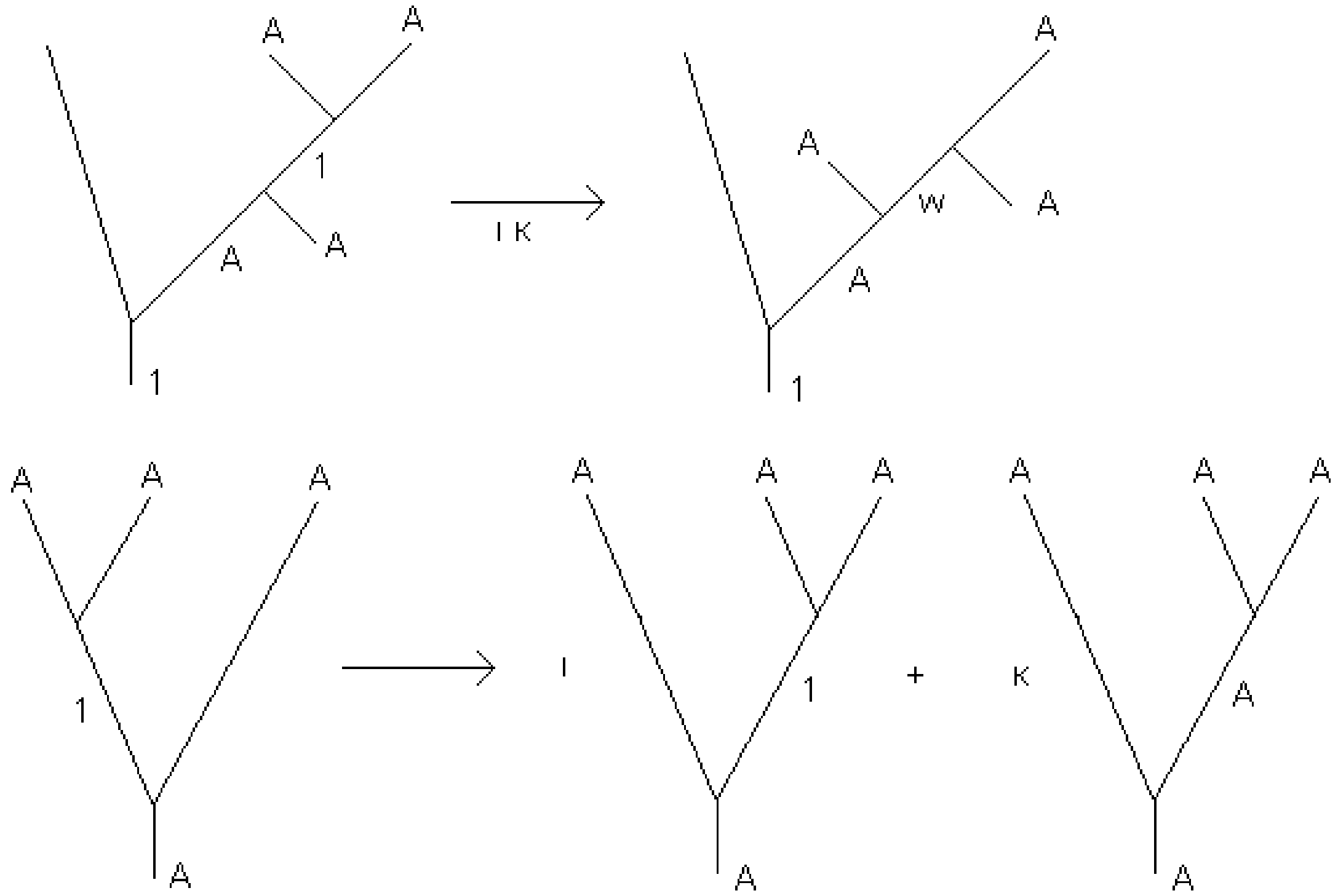

**Figure 291- Computations- example 2- passing a vertex**

We get ι = 2 and κ = 3, again both sequences result in the same linear combination of end roottrees.

## *8.7  Computations – example 3*

We chose the start roottree and apply the sequence corresponding to the new sequence:

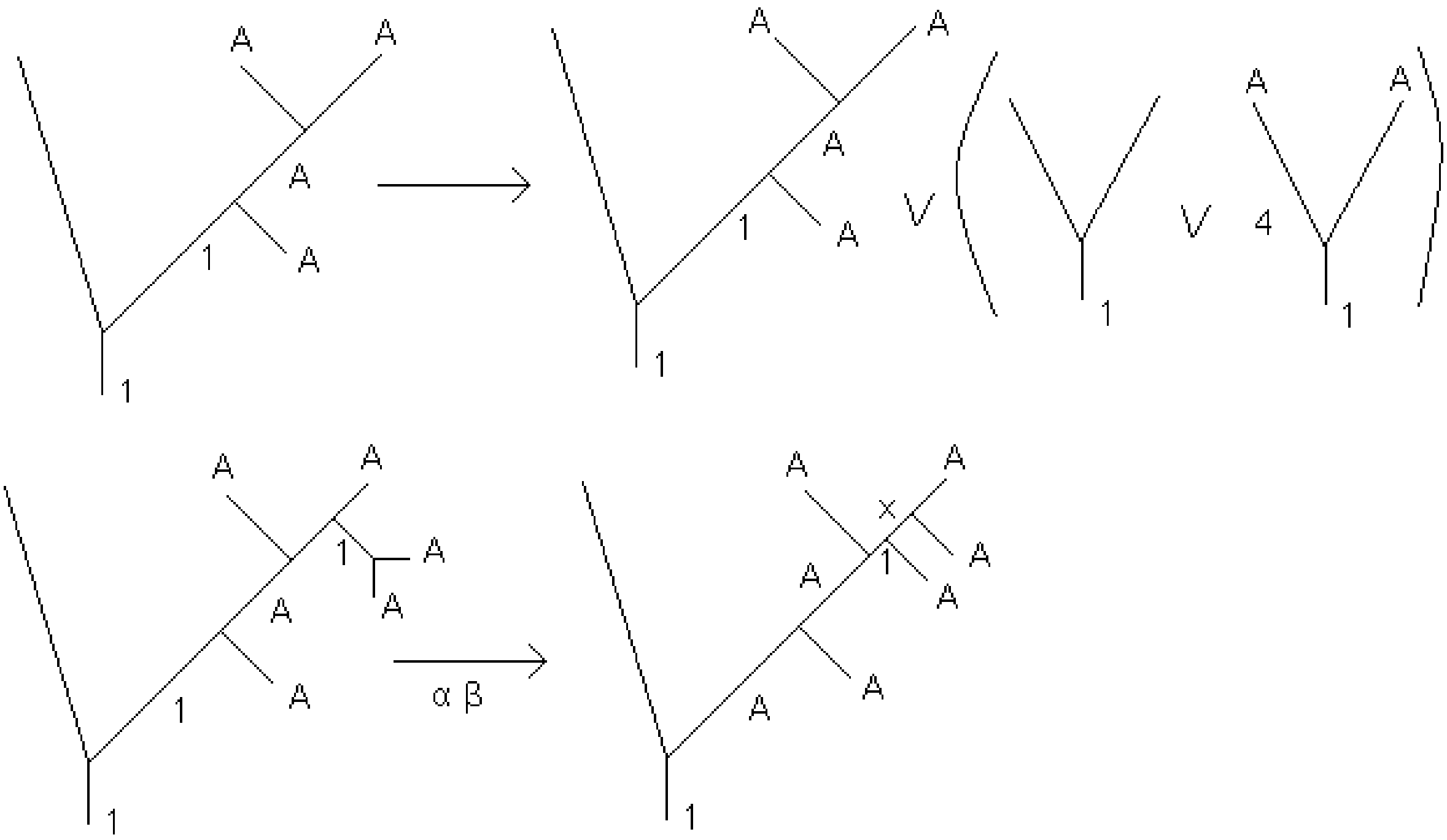

**Figure 292- Computations- example 3- the trace unit**

For $x = 1 \rightarrow \alpha = 2$ and for $x = A \rightarrow \beta = 3$.

We go on with the 2 cases:

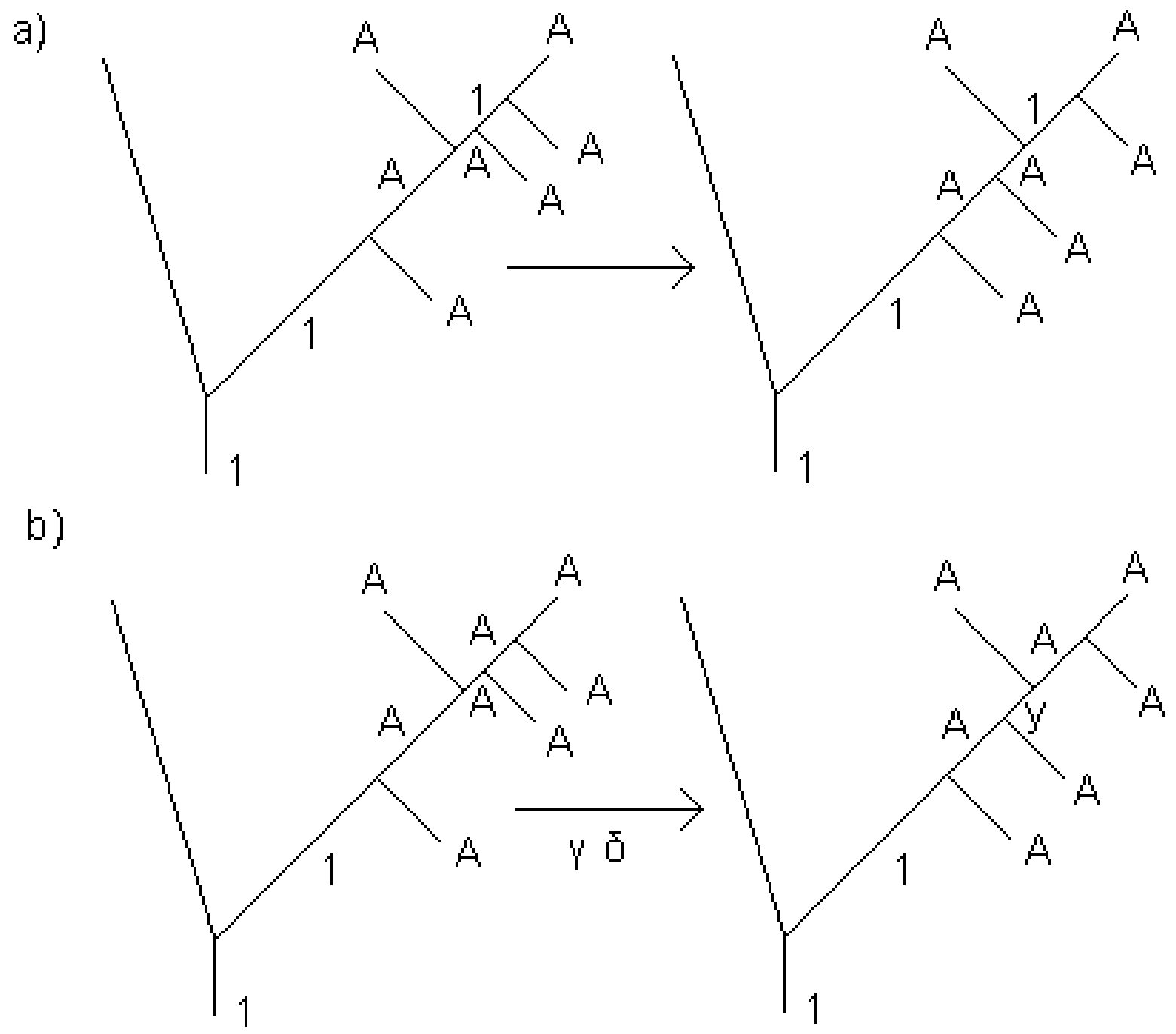

**Figure 293- Computations- example 3- associativity 2**

Case a) leads to a unique simple object and for case b) we get:

y =1 → γ = 4 and δ = 3.

This leads to 3 cases where we have to prepare the split:

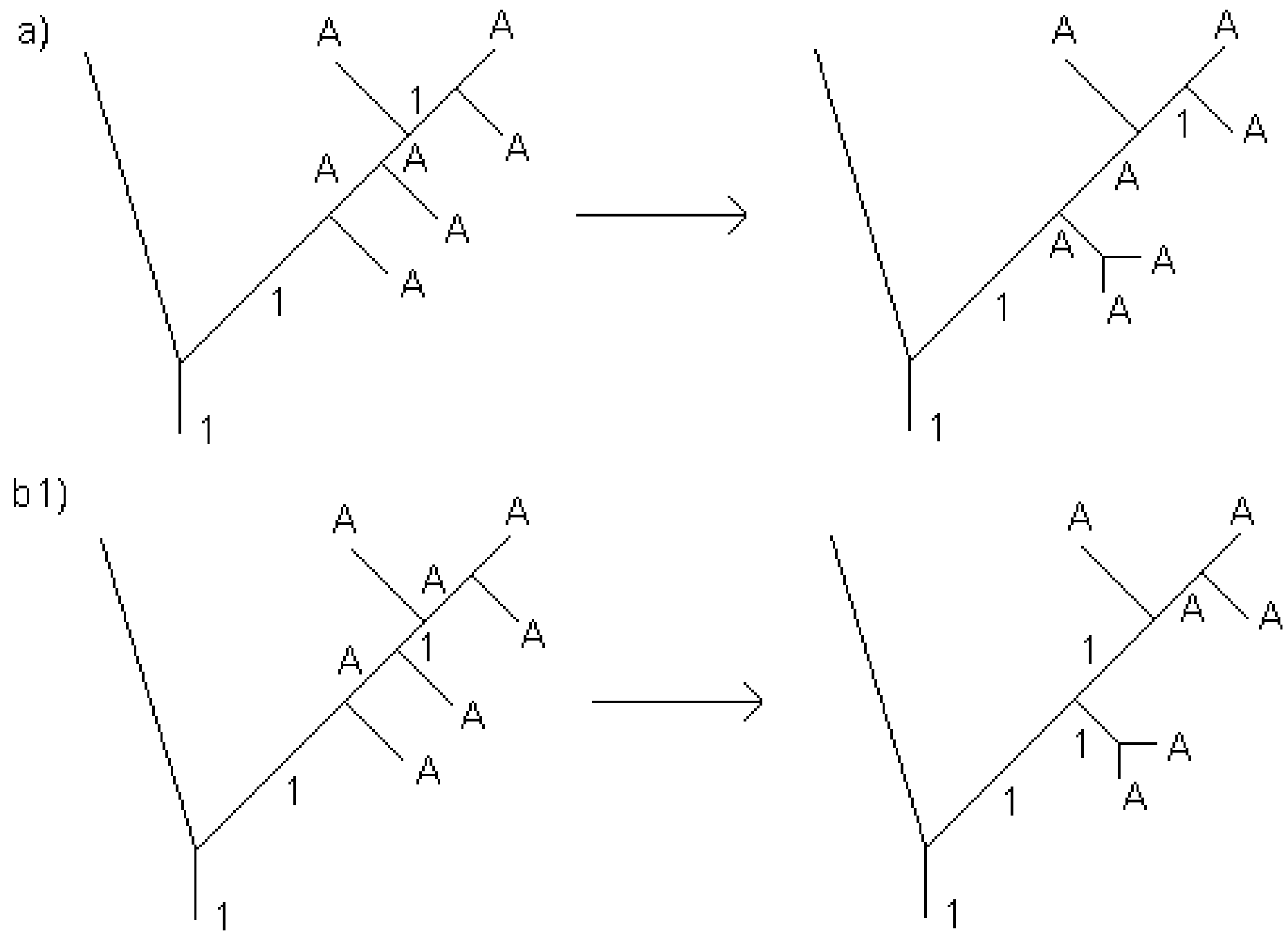

**Figure 294- Computations- example 3- before split- case a and b1**

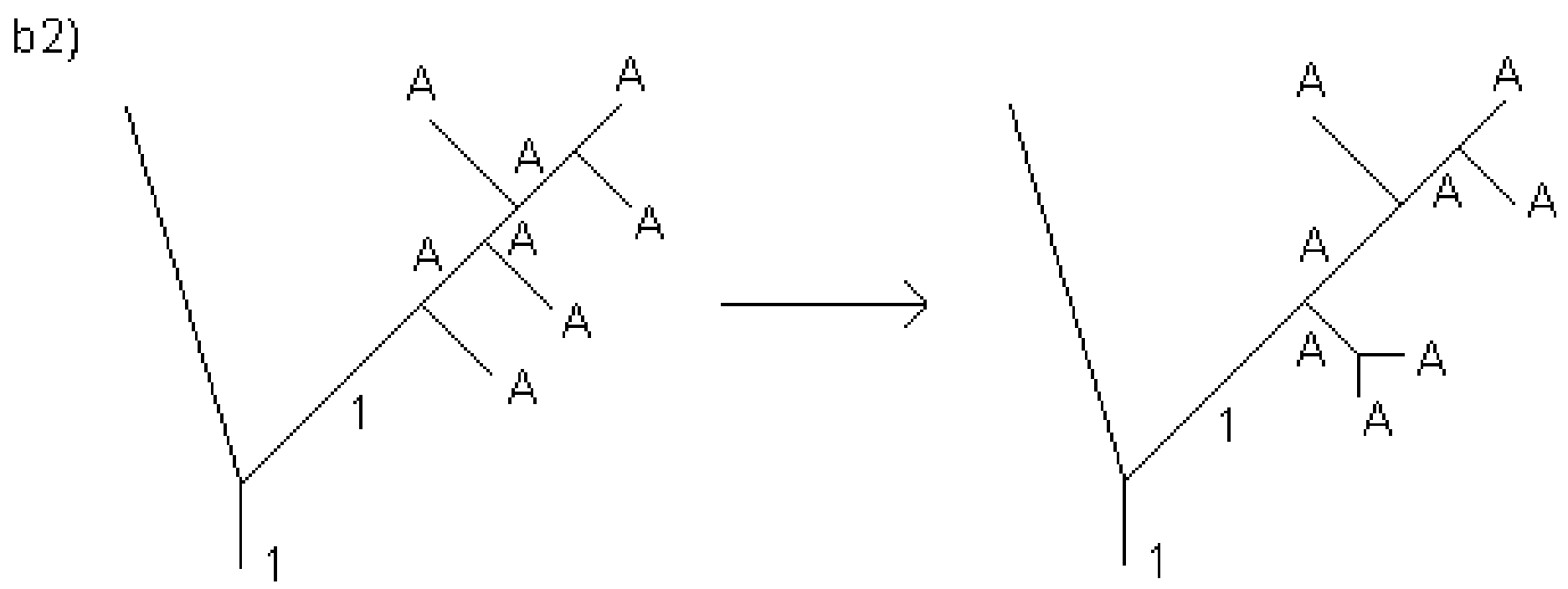

**Figure 295- Computations- example 3- before split - case b2**

Since the state modules of the other roottrees are zero, we only have to consider the case b1). Split the roottree and glue them together at the boundary points respectively to the branches of the trace unit:

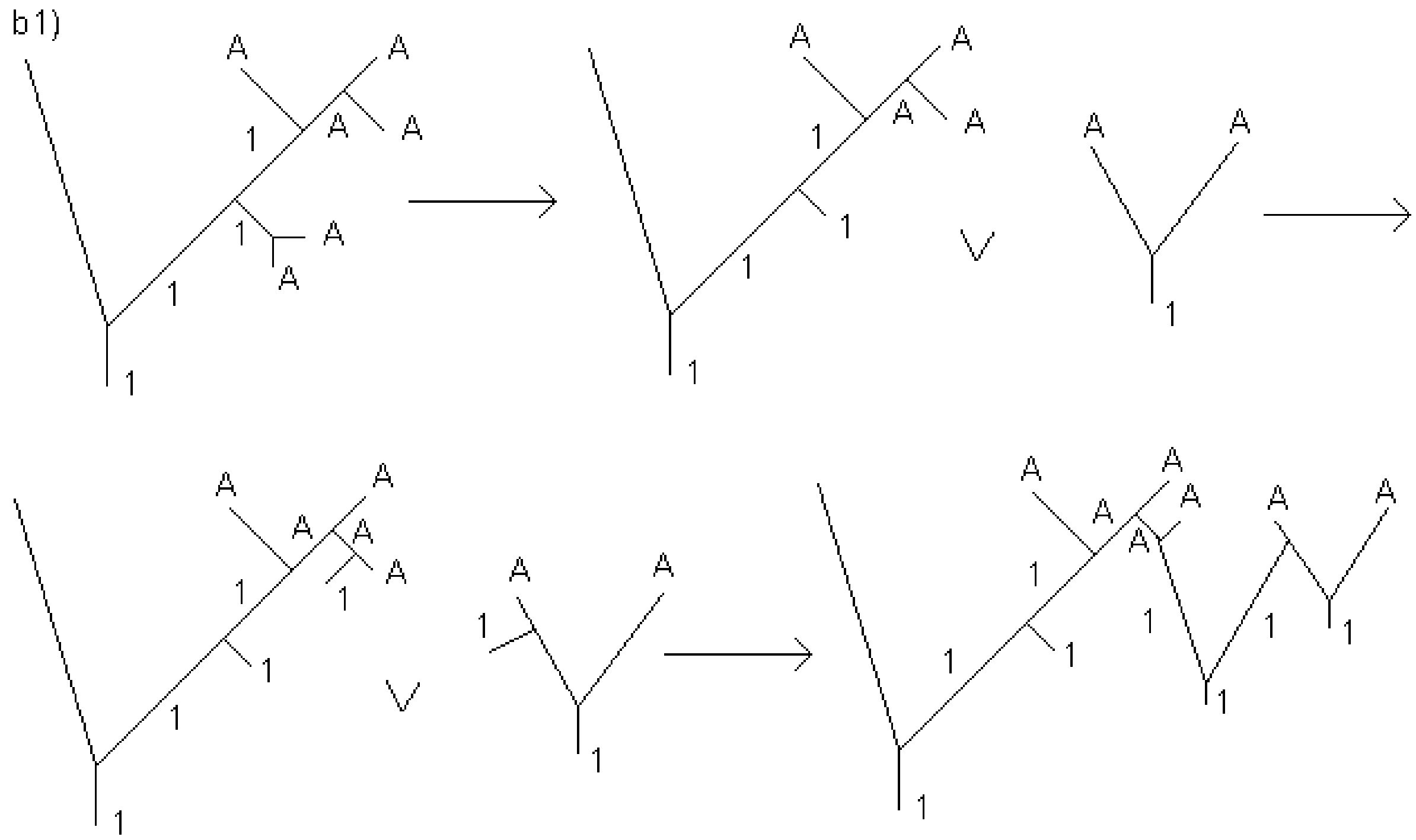

**Figure 296- Computations- example 3- glue the splitted roottree**

Normalize the roottree, apply associativity and the form λ:

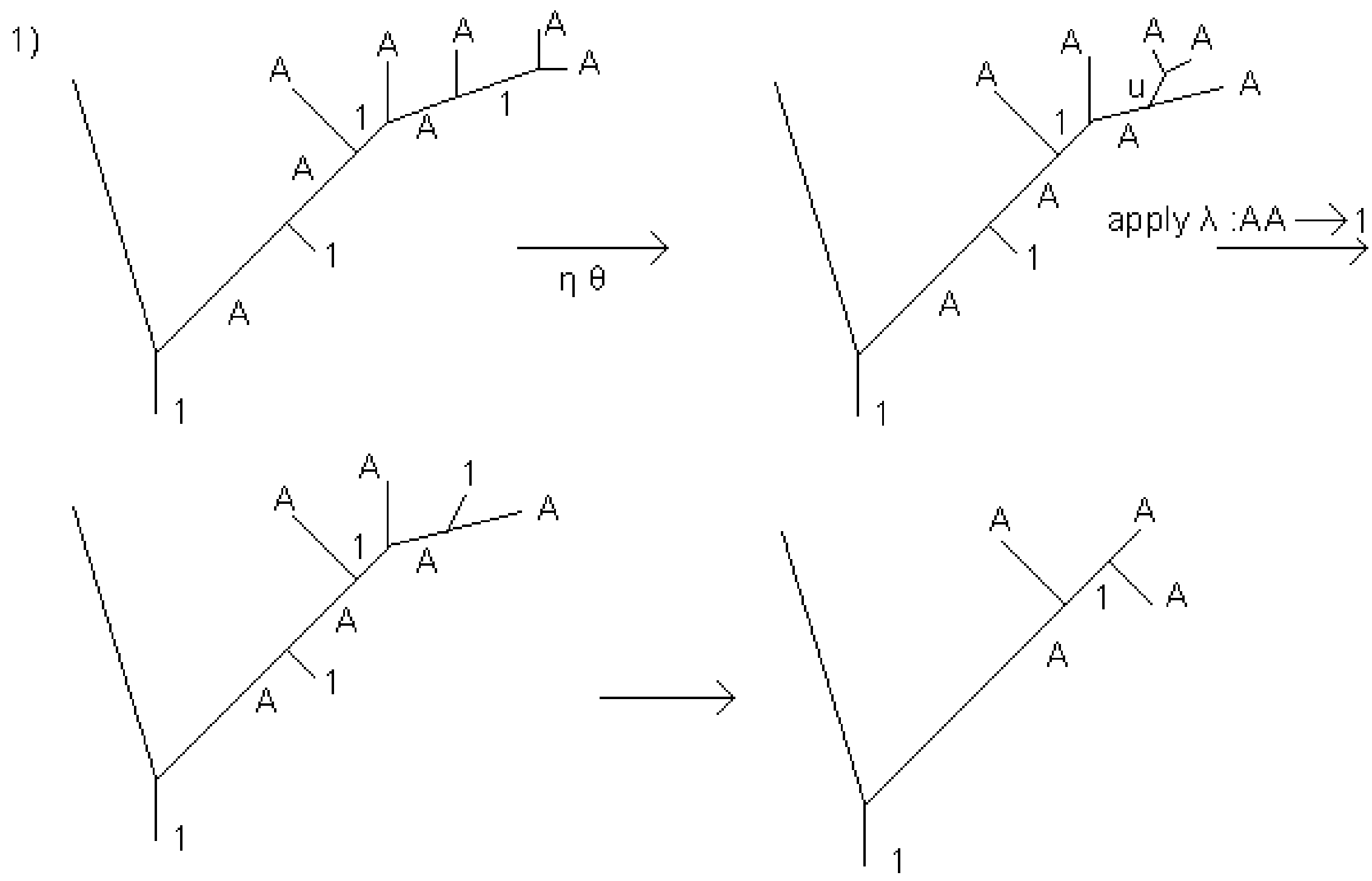

**Figure 297- Computations- example 3- apply the form**

The form λ does not vanish, only for u = 1 where η = 2.

To determine the coefficient of the end roottree, we give an overview:

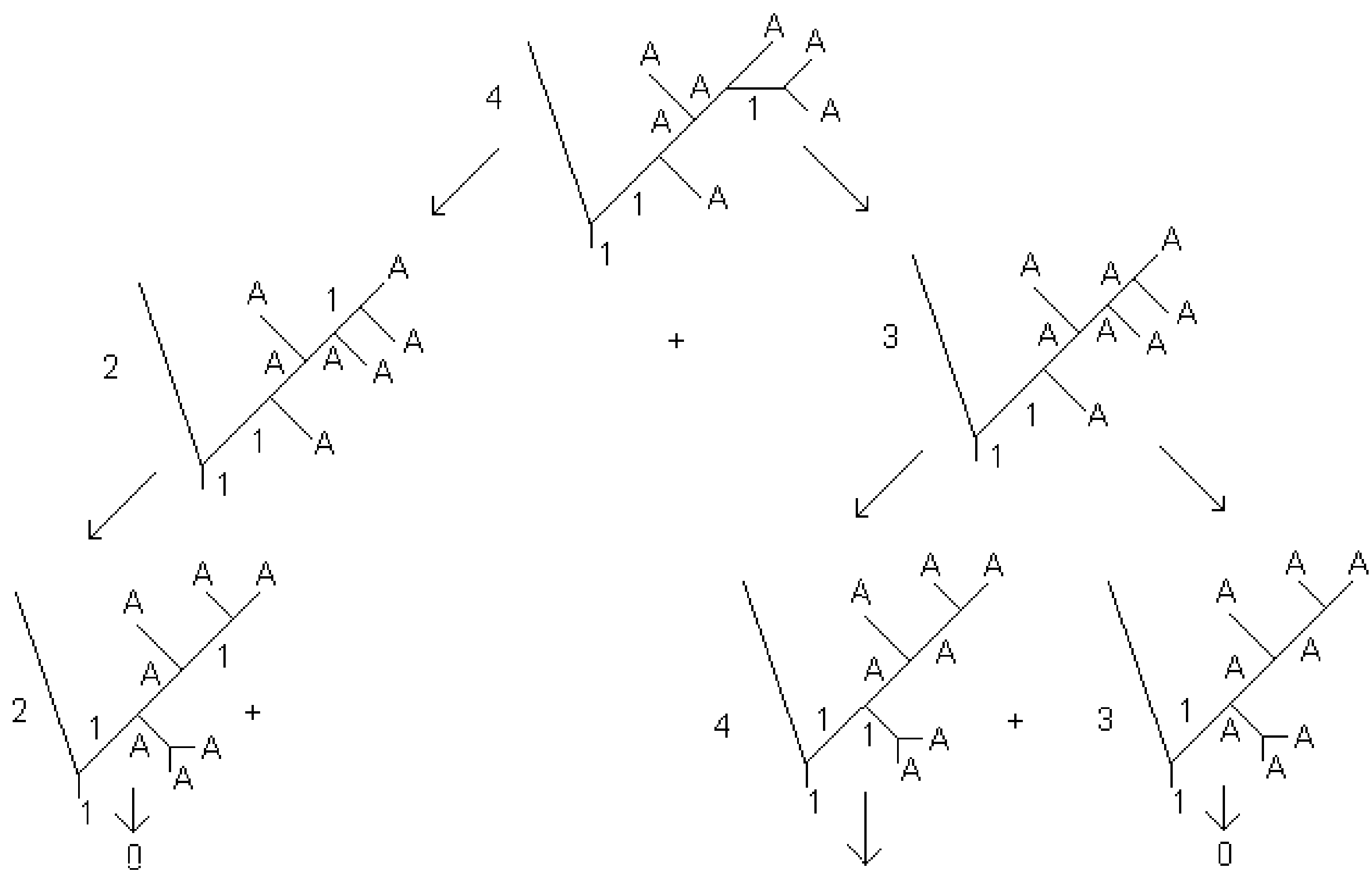

**Figure 298- Computations- example 3- overview 1**

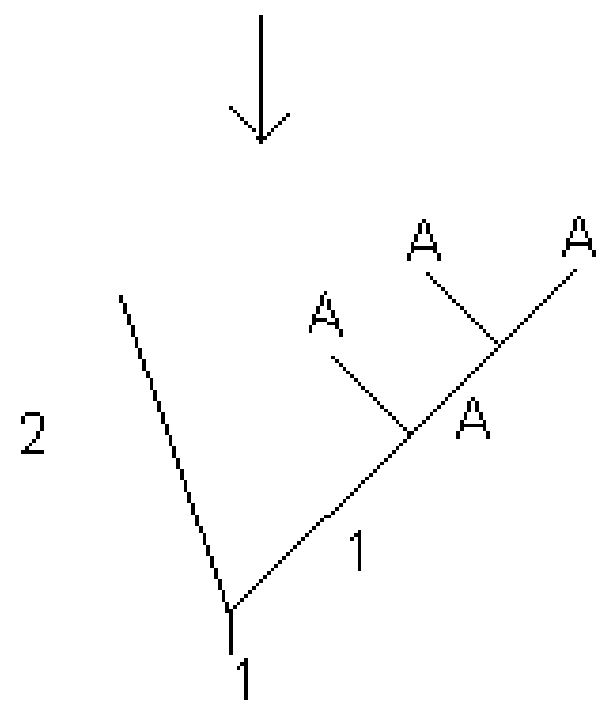

**Figure 299- Computations- example 3- overview 2**

Compute the coefficient:
$4 \cdot 3 \cdot 4 \cdot 2 = 1 \bmod Z_5$

Compare these with the sequence of passing a vertex:

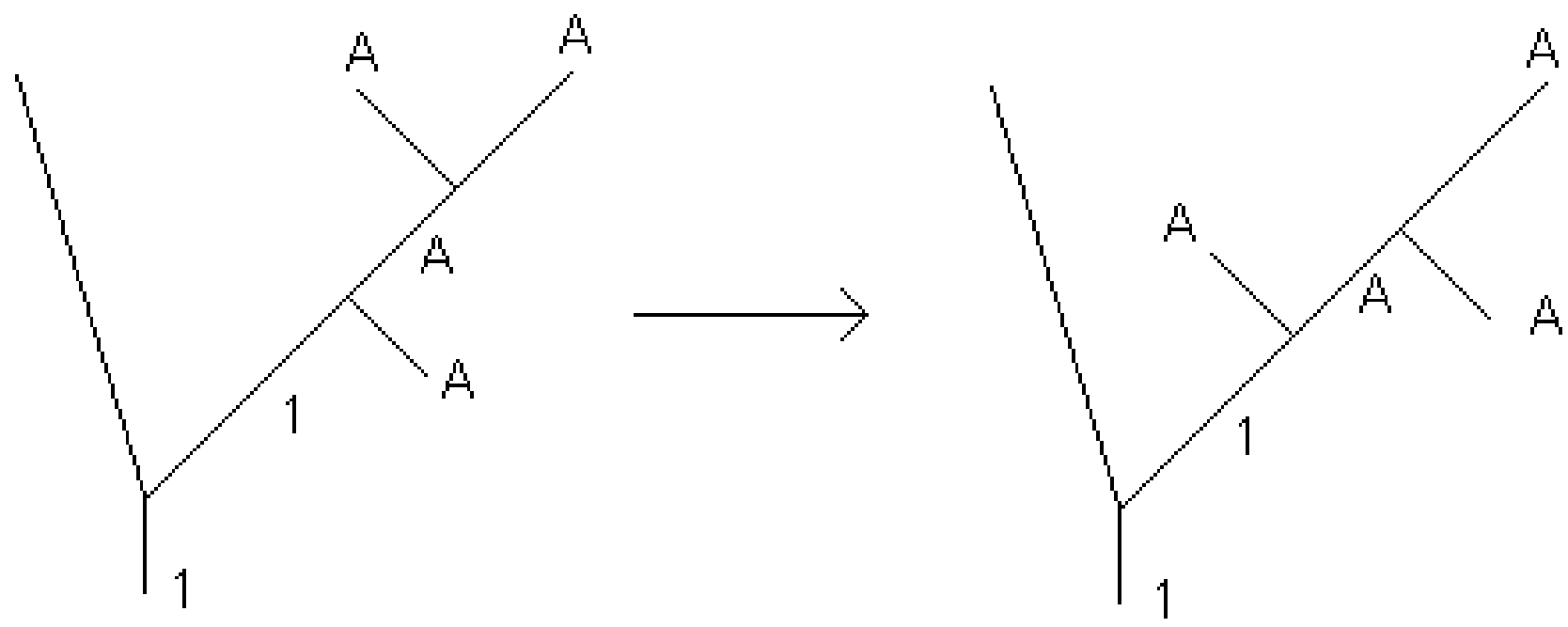

**Figure 300- Computations- example 3- passing a vertex**

We get a unique simple object under associativity and the coefficient is also 1 as in the sequence before.

We have dealt with all cases where the last branches are A-A forks, so it remains to consider the cases where the last branches are 1-1 forks. Note, that for the last branches we only need to regard forks with the same objects at the boundary, this corresponds to the circle of the generator cylinder as part of a slice.

## *8.8 Computations – example 4*

The start roottree has an 1-1 fork at the last branches; apply the trace unit:

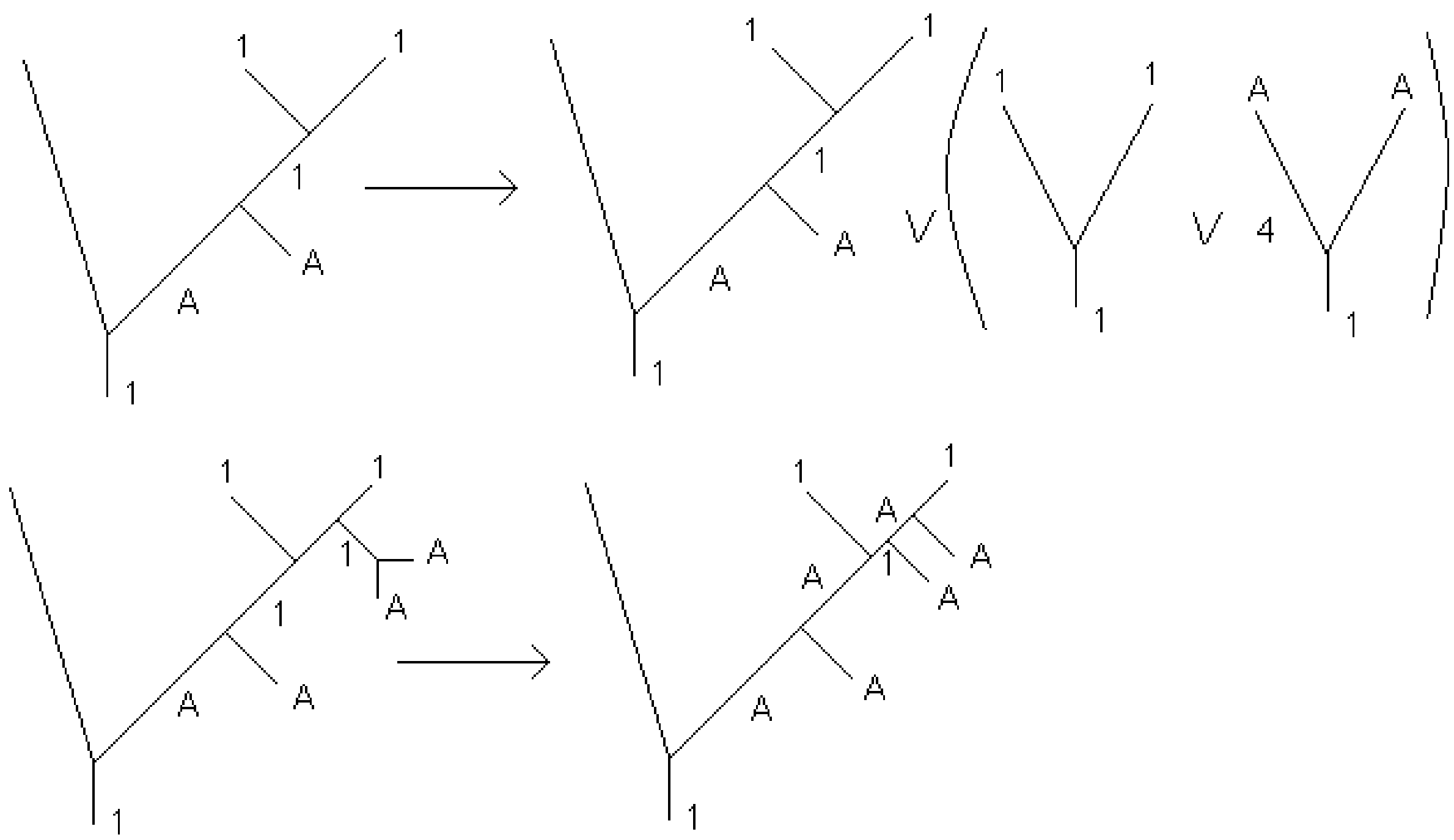

**Figure 301- Computations- example 4- trace unit**

Use associativity and prepare the split:

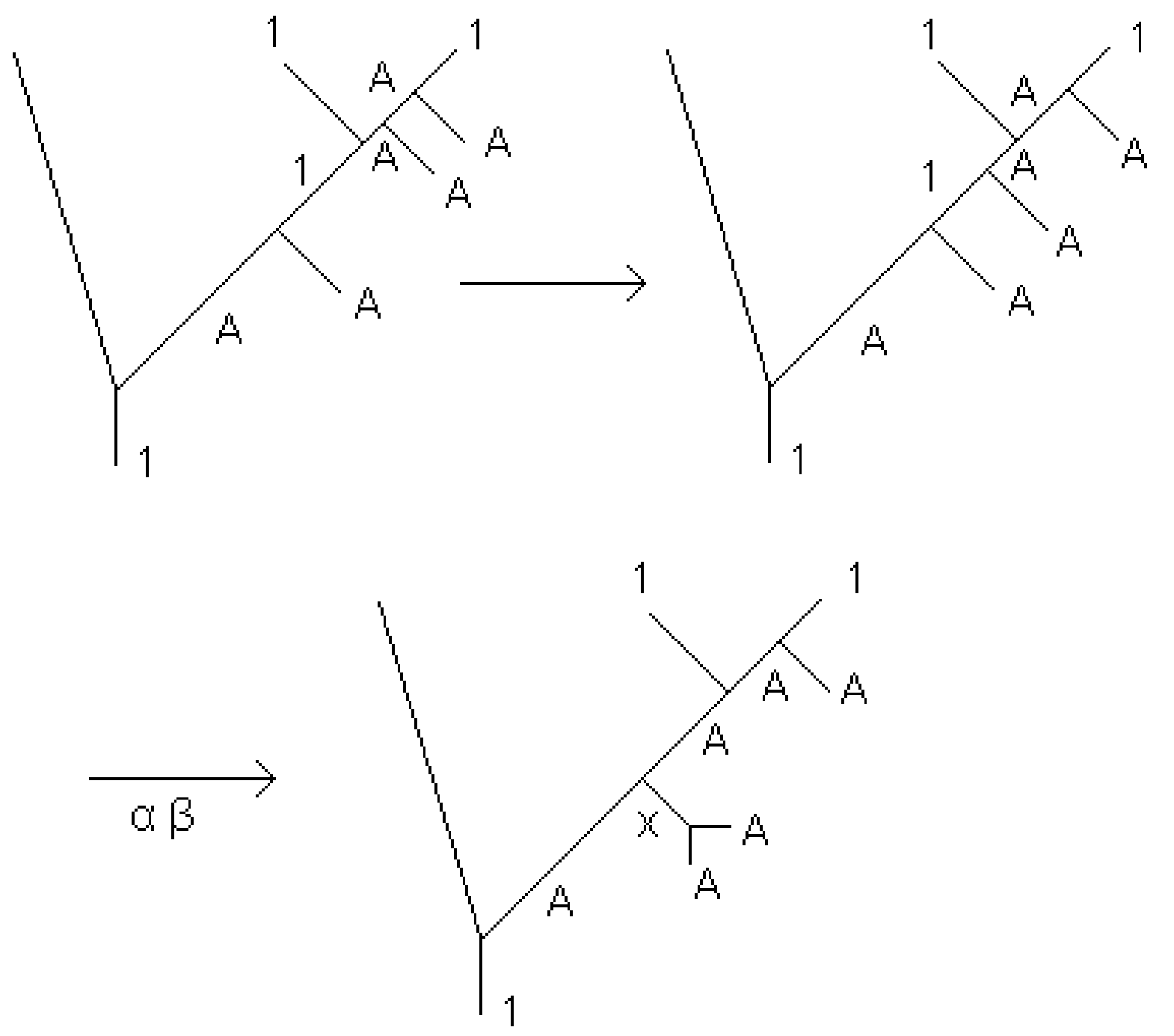

**Figure 302- Computations- example 4- before split**

For x = 1 → α = 2 and for x = A → β = 2.

Perform the split for x = 1 and glue the splitted trees together at the bounday of the branches coming from the trace unit:

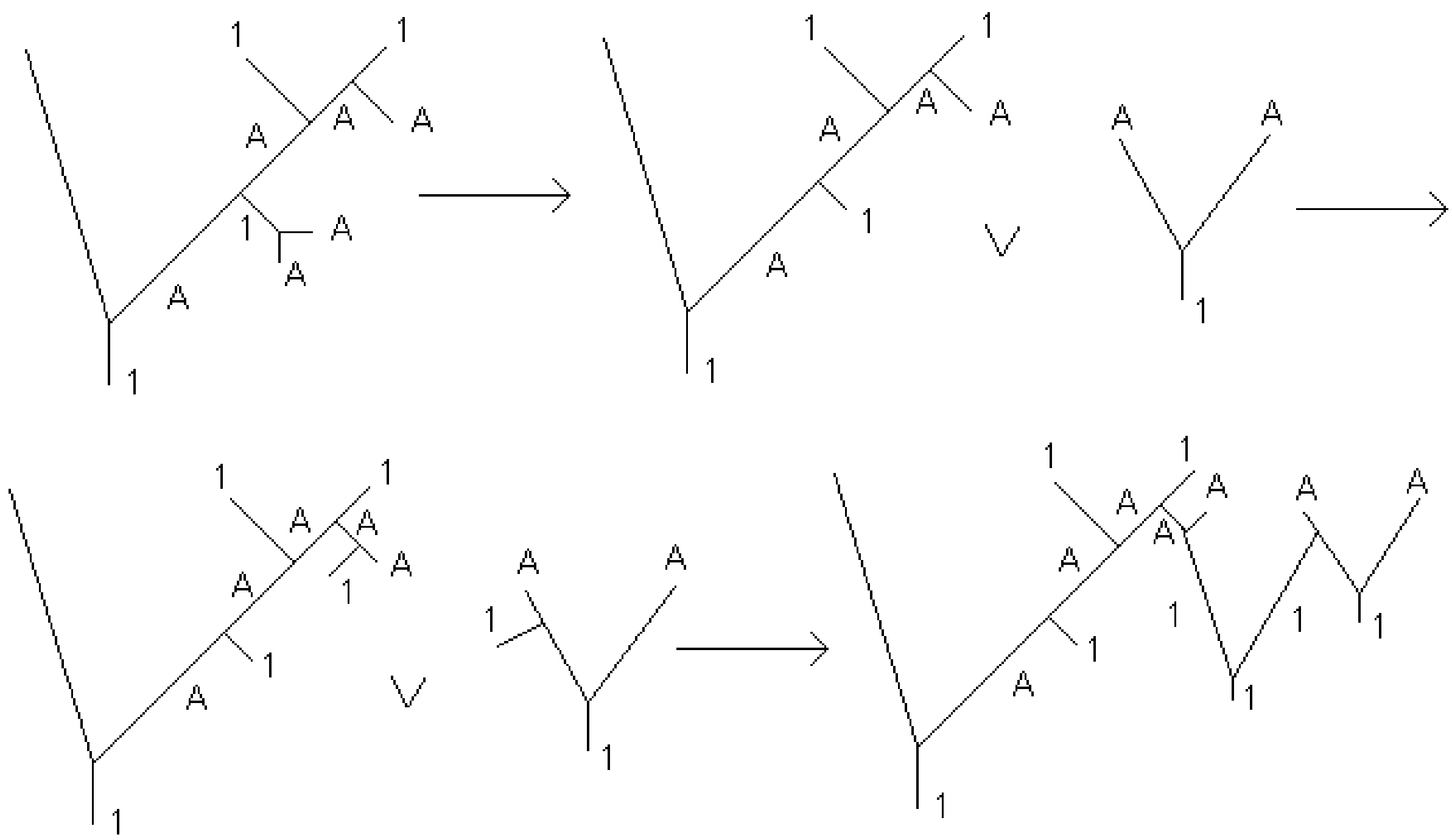

**Figure 303- Computations- example 4- glue splitted trees**

Normalize the roottree, apply associativity to move the branches of the trace unit together and then apply the form λ on that A-A fork:

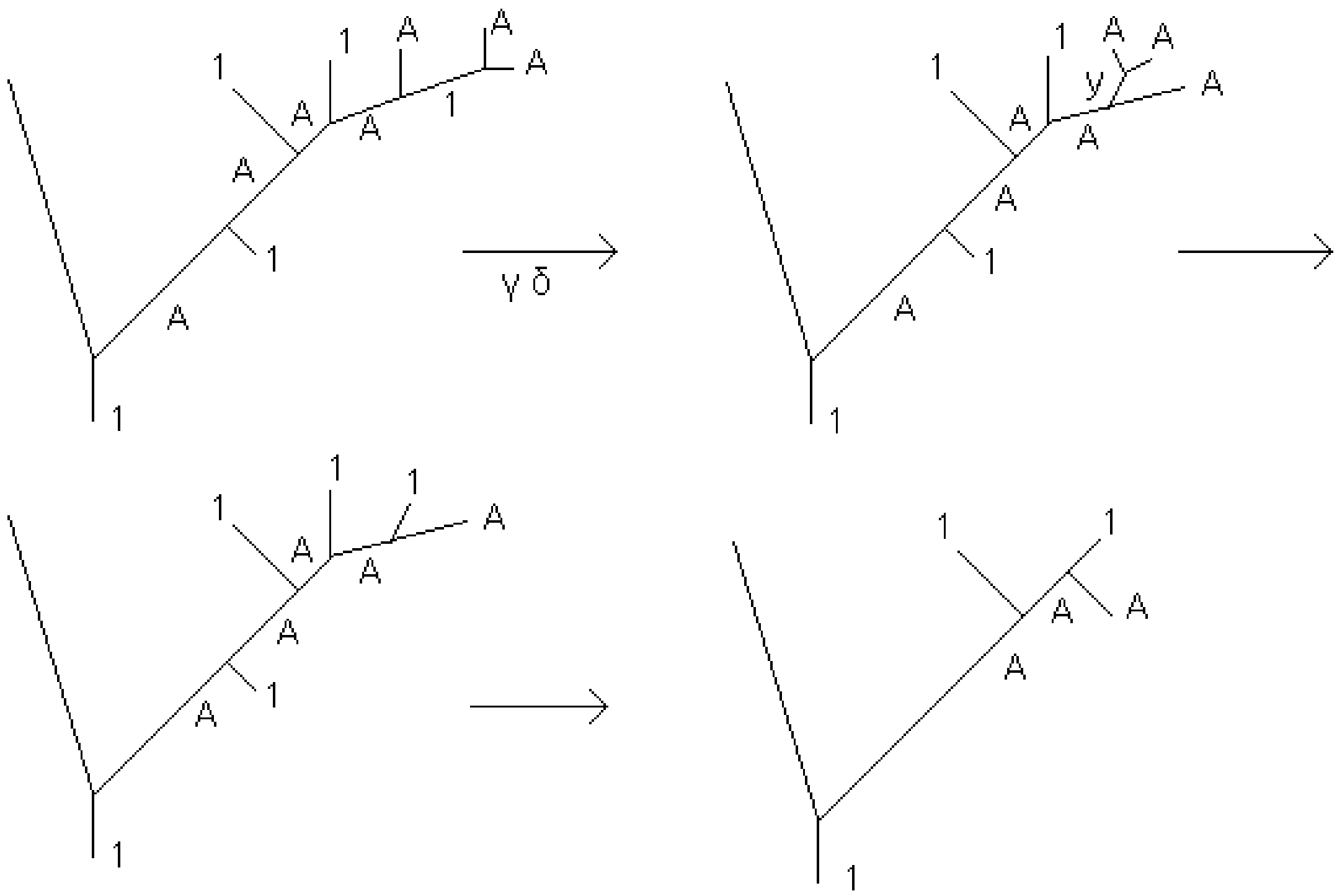

**Figure 304- Computations- example 4- apply the form**

For y = 1 → γ = 2 and for y = A → δ = 3.
But the form λ does not vanish only for the case y = 1.

The coefficient of the end roottree is:
$4 \cdot 2 \cdot 2 = 1 \bmod Z_5$

Determine the end rootree for the sequence of passing a vertex:

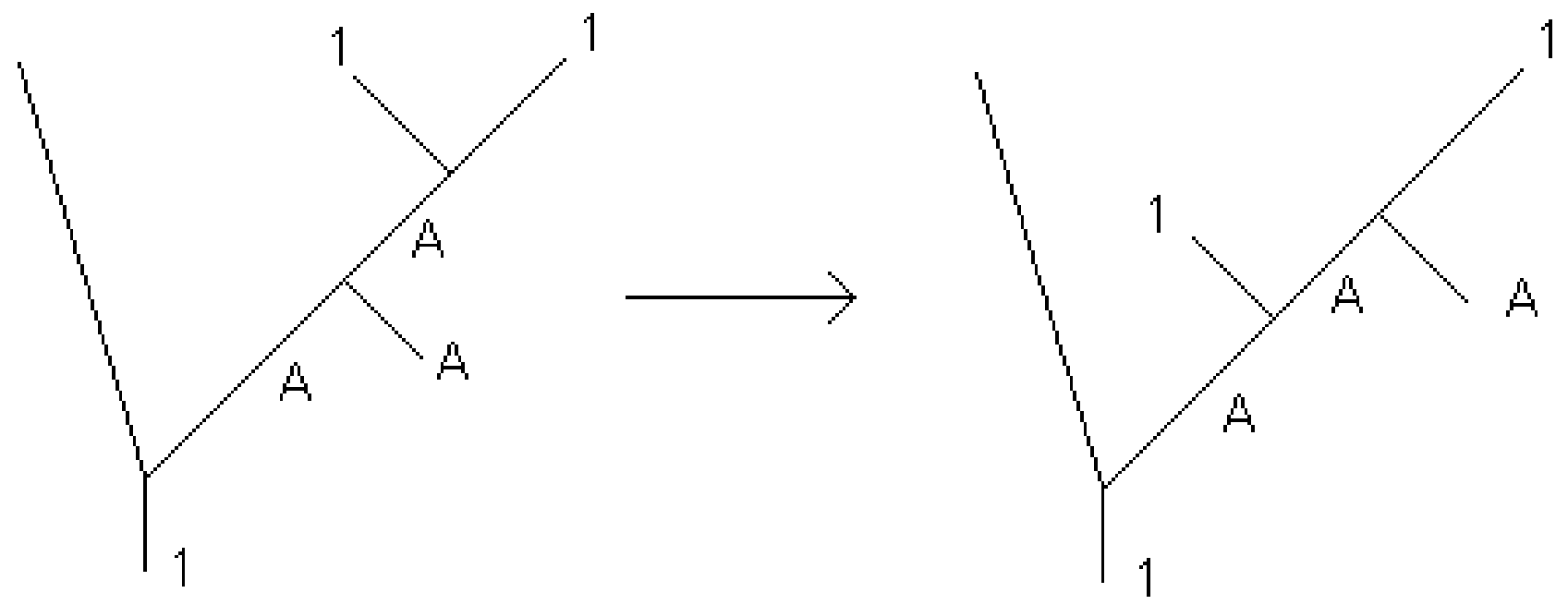

**Figure 305- Computations- example 4- pass a vertex**

Since the coefficient is also 1, both sequences result in the same linear combination of end roottrees.

We close our considerations with a few comments to the last case, where all branches and roots are 1:
If we attach the A-A fork this will vanish at splitting, because we would split at A.
The 1-1 fork of the trace unit works, but since the entry in the associativity matrix for the 1 object is 1, the coefficient is also 1 for the end roottree.

# 9 Calculus – the relation that builds a bubble

We compute for the tensor category, presented in the example above, that the relation building a bubble is an algebraic relation, that means the corresponding sequence of rootrees maps the start rootree with identity to the end rootree.

## *9.1 The Relation presented as slices and roottrees*

We remember, that the relation arises by pushing a bubble out of a rectangle and add little flanges to get the local vertex model:

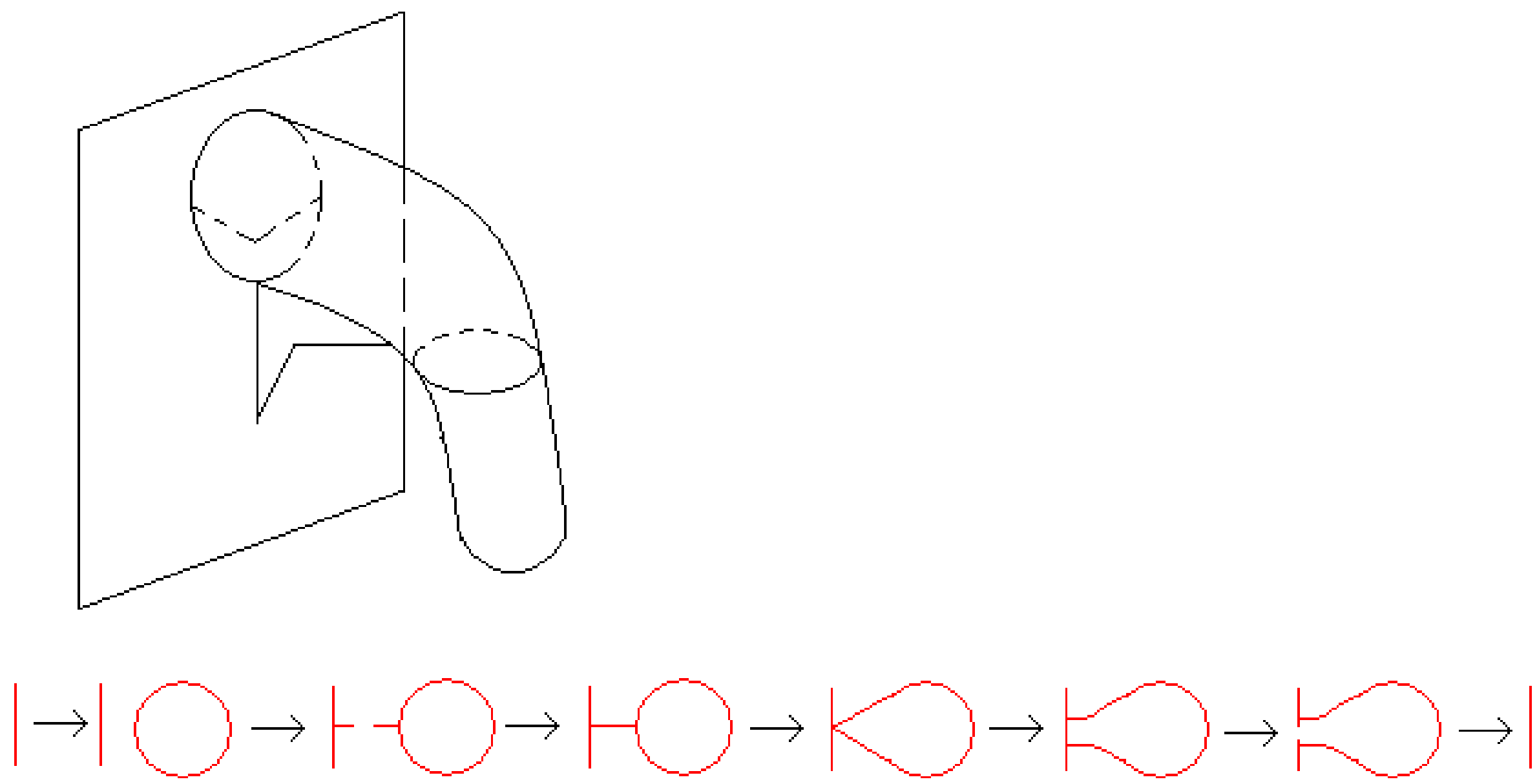

**Figure 306- The relation as slices and roottrees- the sequence of slices**

We transfer this to a sequene of rootrees:
For $S^1$ we use the trace unit and to break up the interval in the second last slice we split again the corresponding branch as justifid in chapter 8.3.
The remaining little branches come from the added flange and after splitting collaps to a point, only the object 1 at the boundary makes sense:

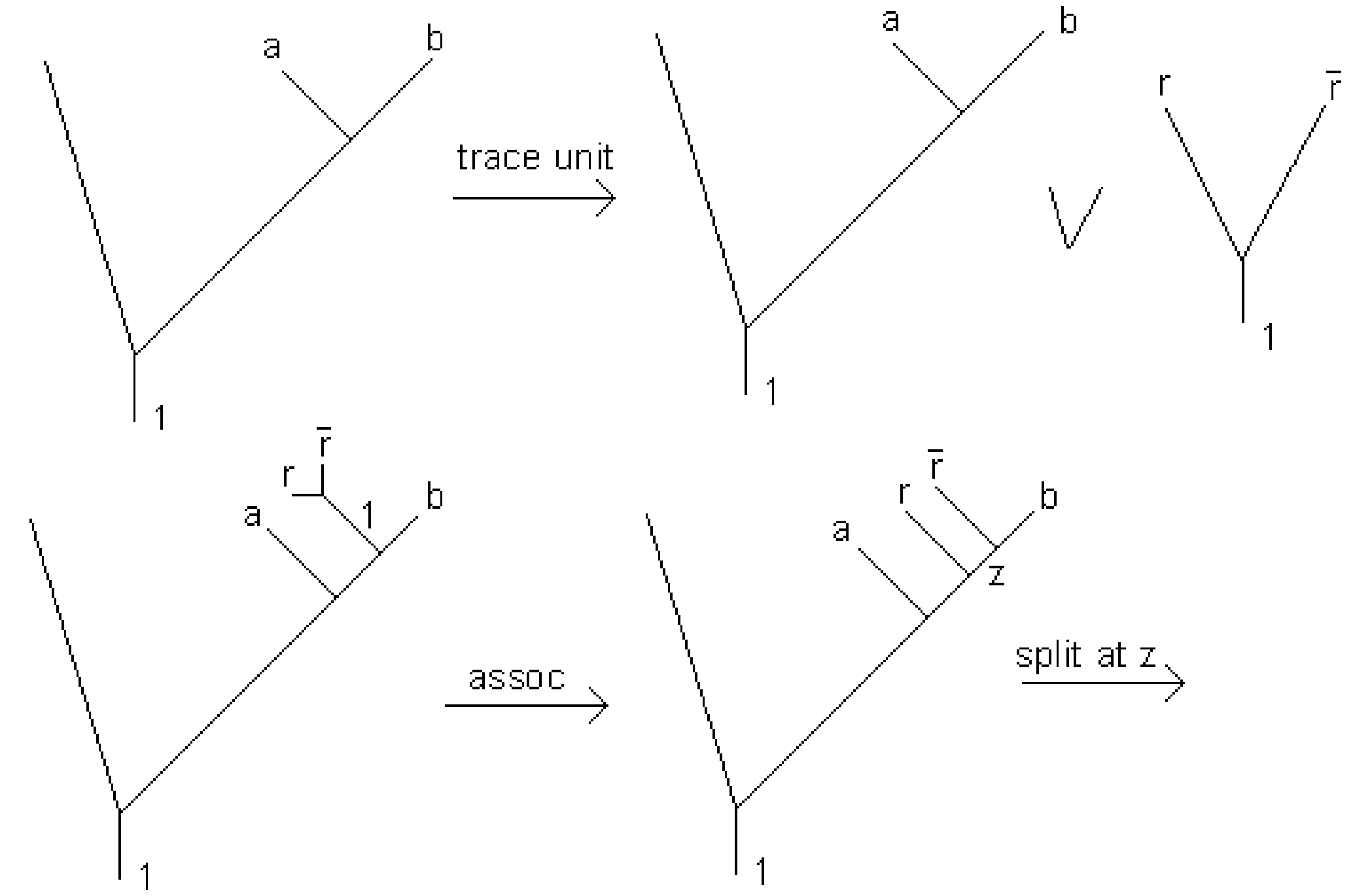

**Figure 307- The relation as slices and roottrees- the roottrees- trace unit**

Glue together the splitted trees and identify the boundary of the branches with respect to the trace unit:

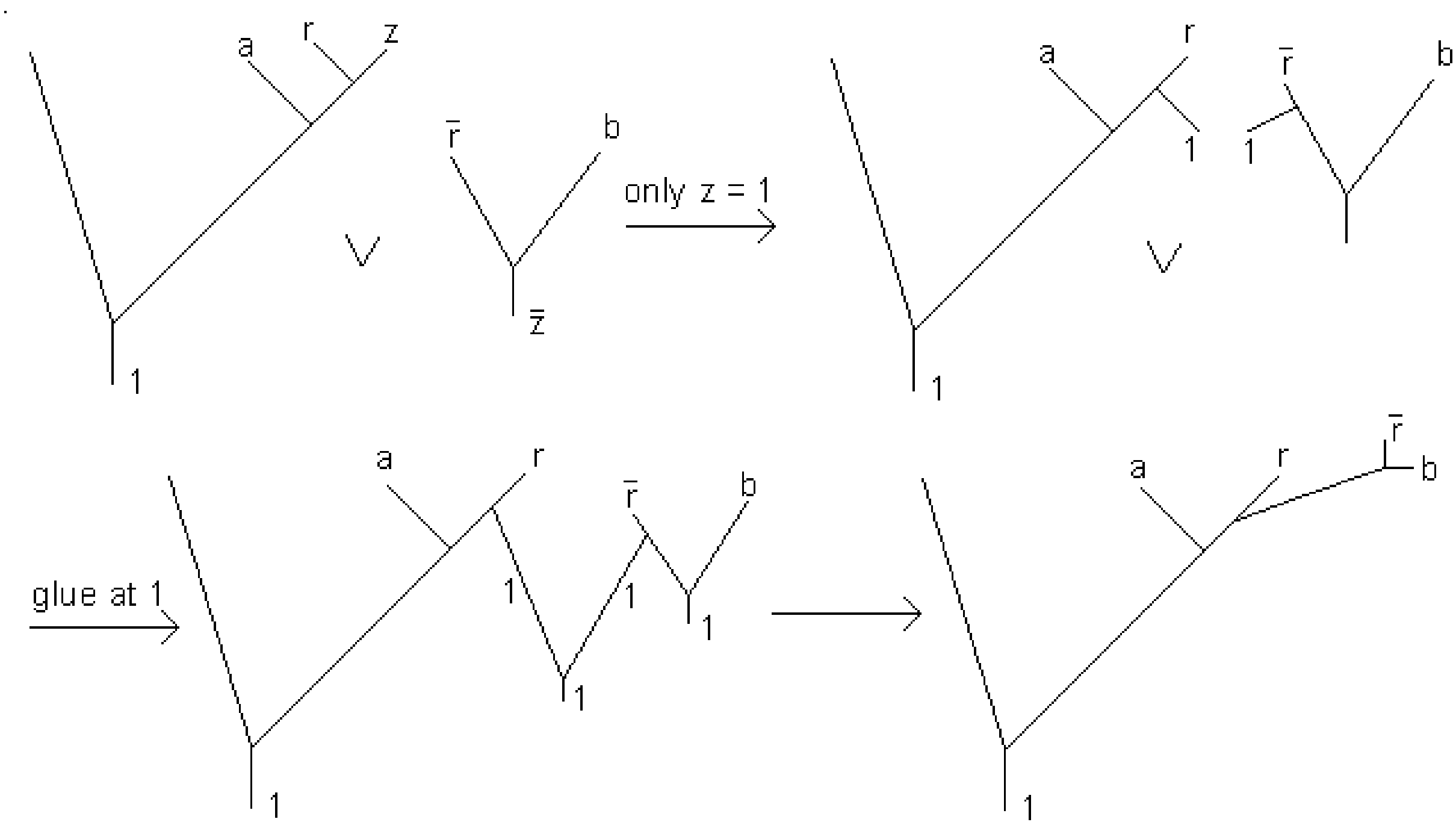

**Figure 308- The relation as slices and roottrees- the roottrees- glue splitted roottrees**

Normalize the roottree, apply associativity and the form λ:

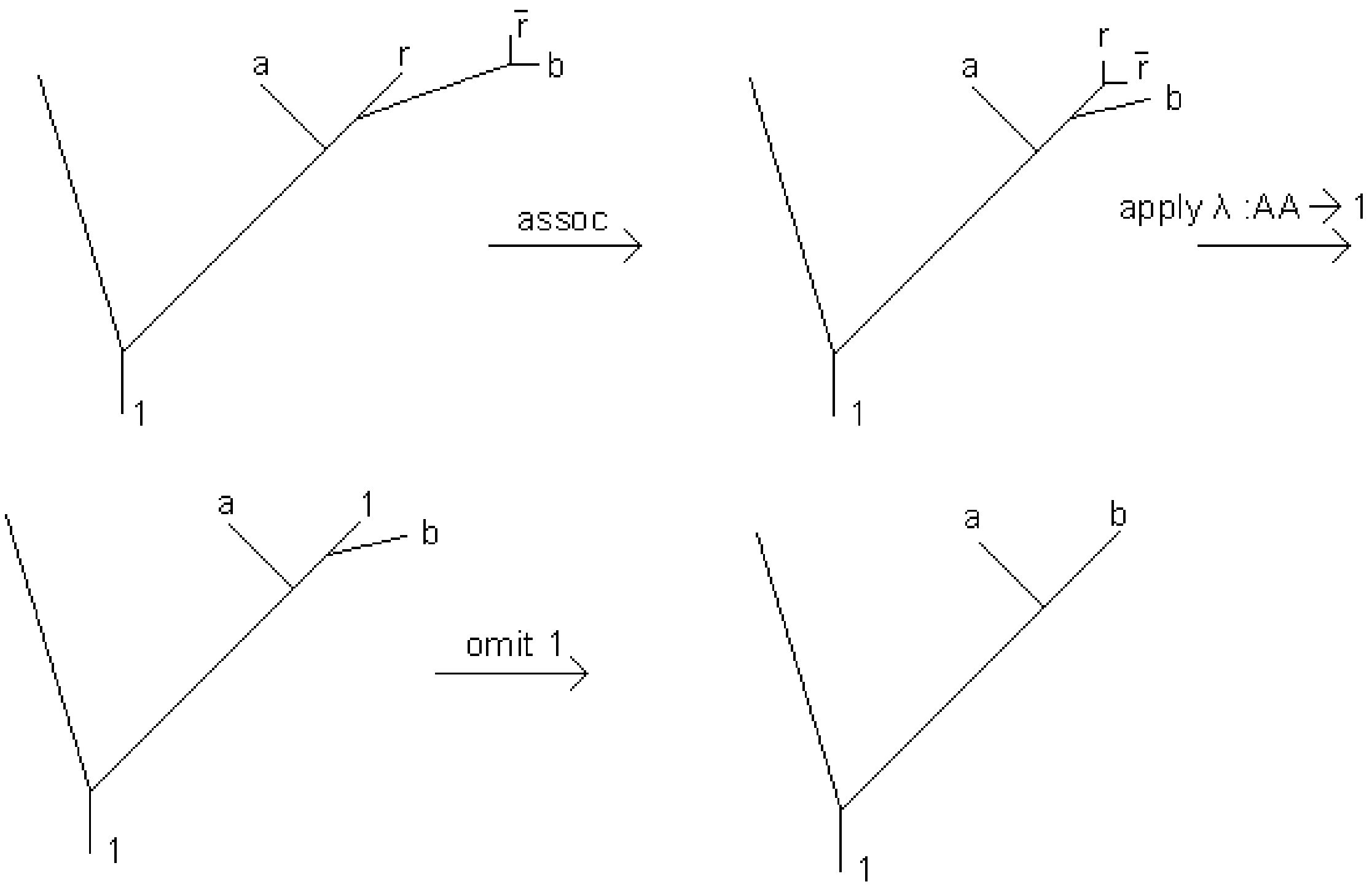

**Figure 309- The relation as slices and roottrees- the roottrees- apply the form**

The form λ: $r \otimes \bar{r} \rightarrow 1$ does not vanish only for the root 1 of the r-$\bar{r}$ fork.

## *9.2 Verify the relation*

We construct the trace unit and apply associativity to prepare the split:

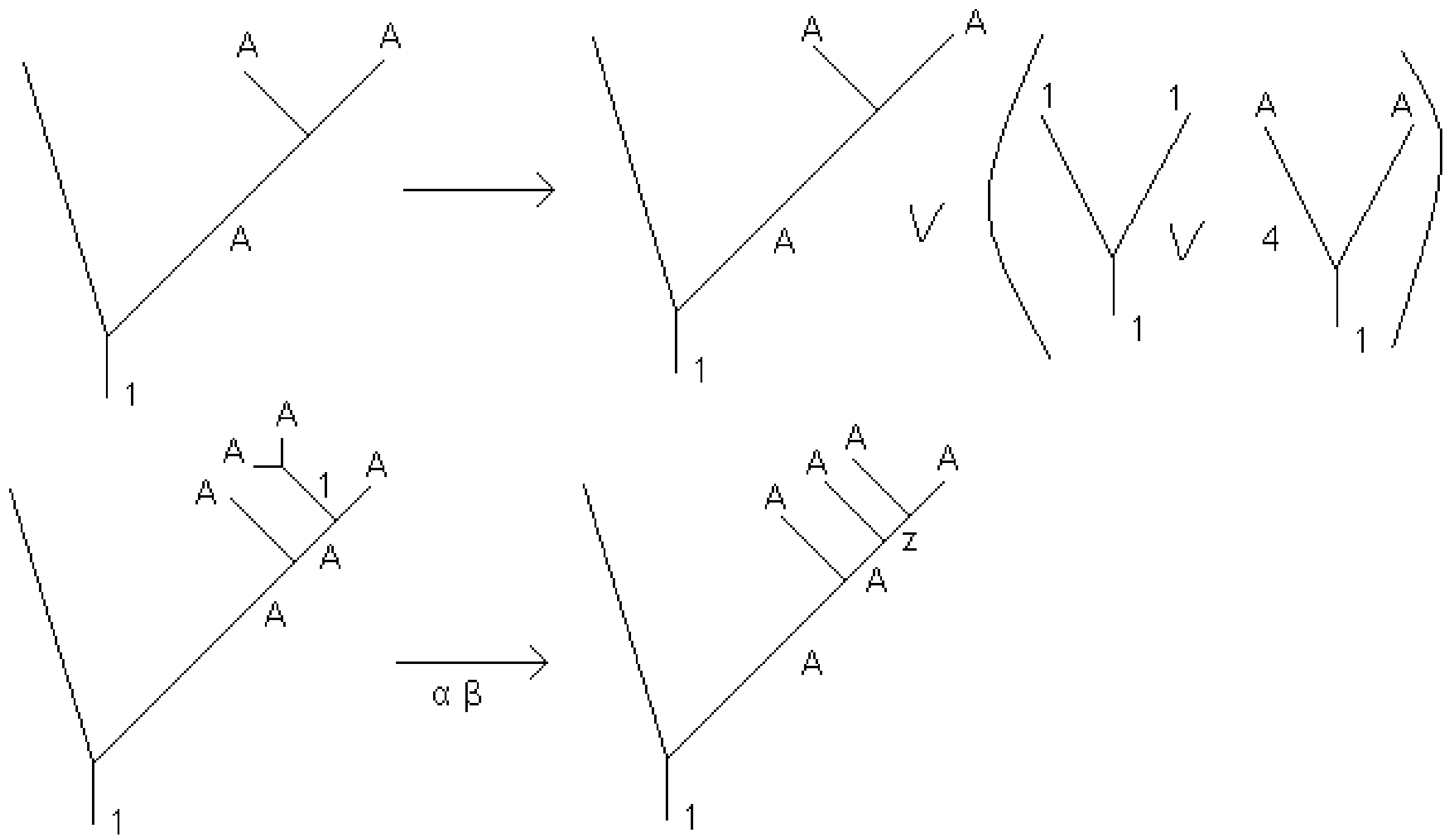

**Figure 310- Verify relation- trace unit**

For z = 1 → α = 2 and z = A → β = 3.

For the split consider only the case z = 1 and glue the splitted trees together:

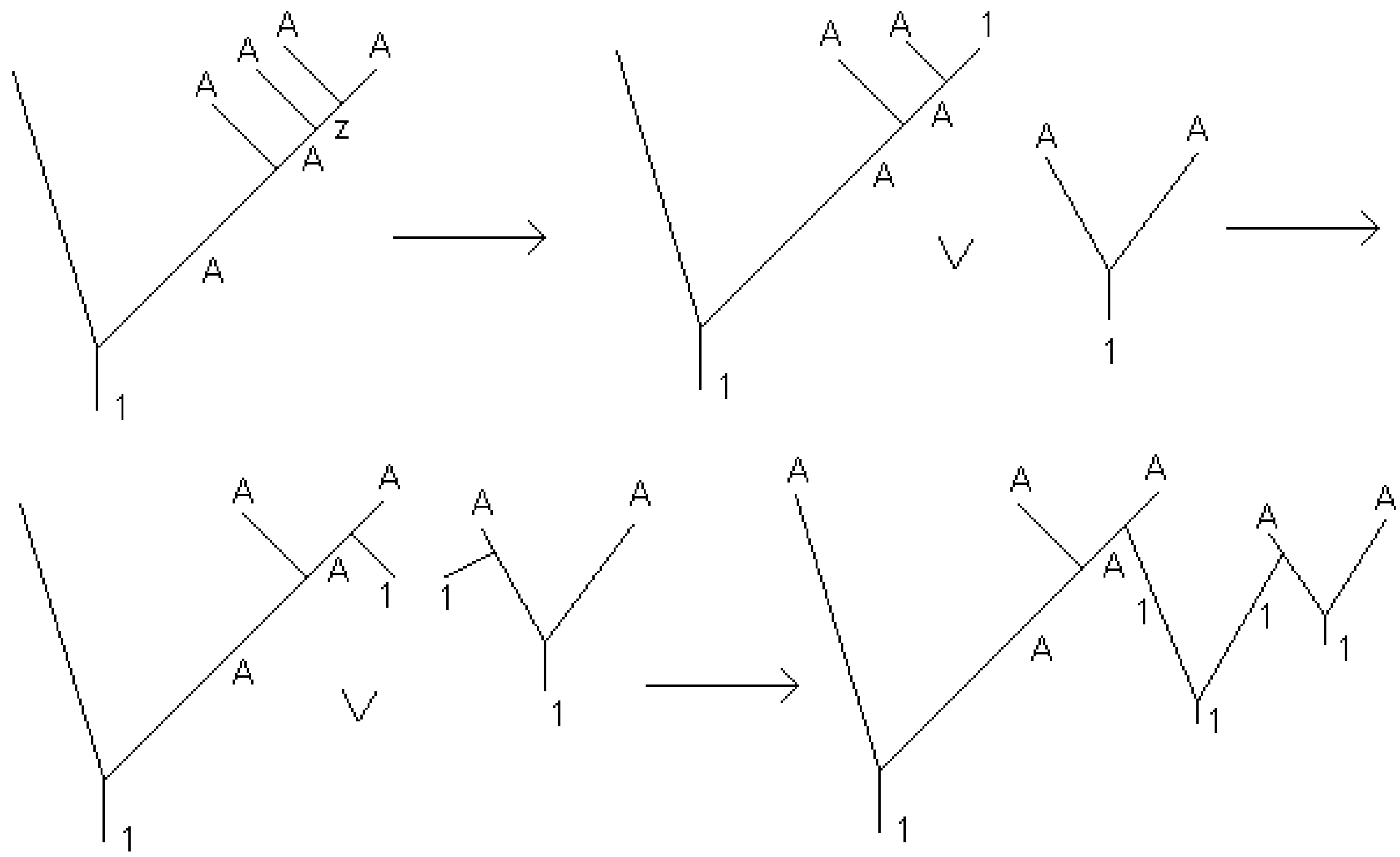

**Figure 311- Verify relation- glue splitted roottrees**

Normalize the roottree, apply associativity and the form λ:

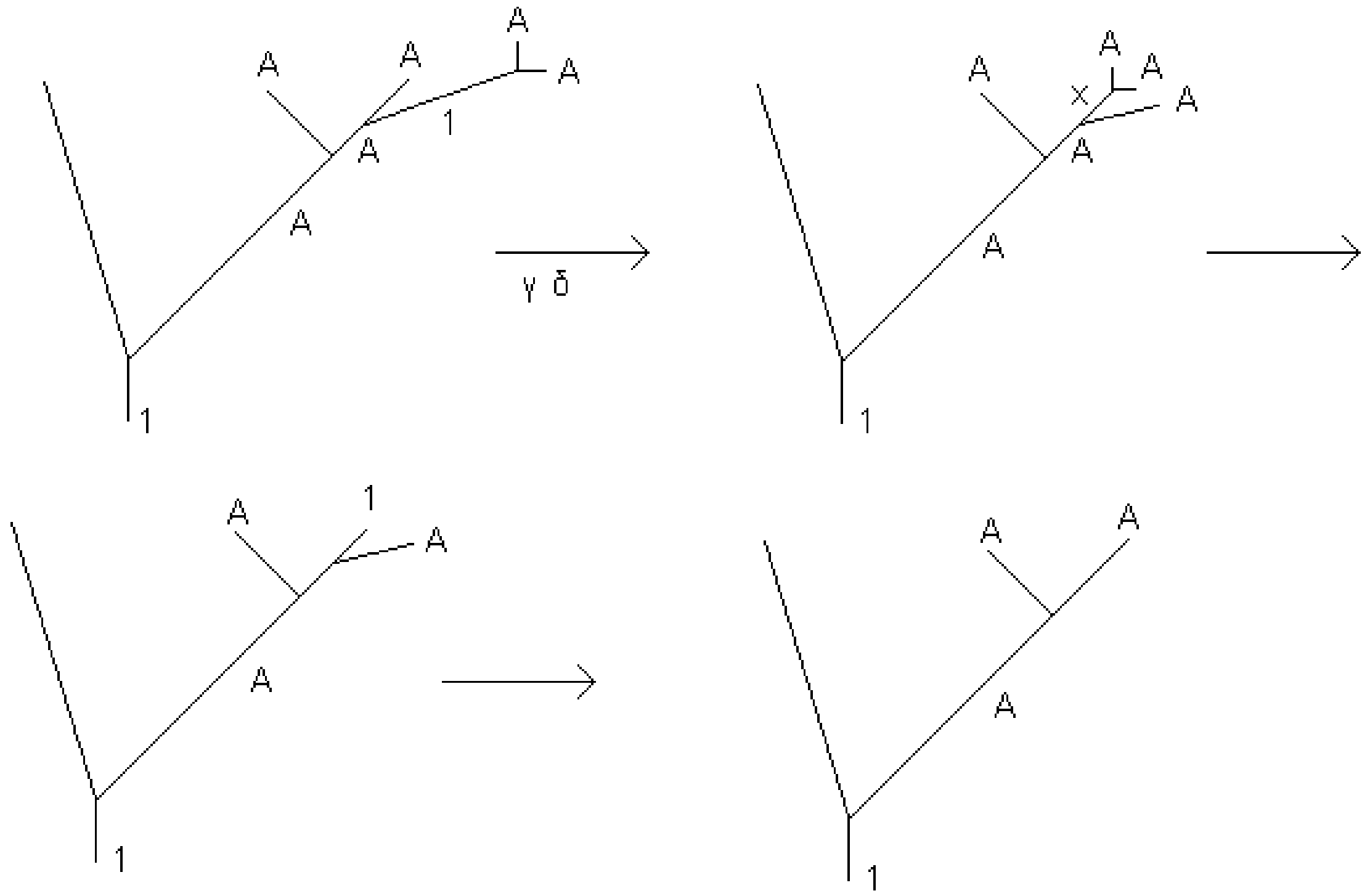

**Figure 312- Verify relation- apply the form**

We get x = 1 → γ = 2 and for x = A → δ = 3, but only x = 1 survives the form λ: x → 1 since the tensor category is semisimple.

the resulting coefficient of the end roottree is :
$4 \cdot 2 \cdot 2 = 1 \bmod Z_5$

Hence we get the identity map.

For the next case we only compare the first step, we see that the part, where we split is assigned with the same simple objects as in the previous case. The first fork with 1 as root corresponds to the start roottree.

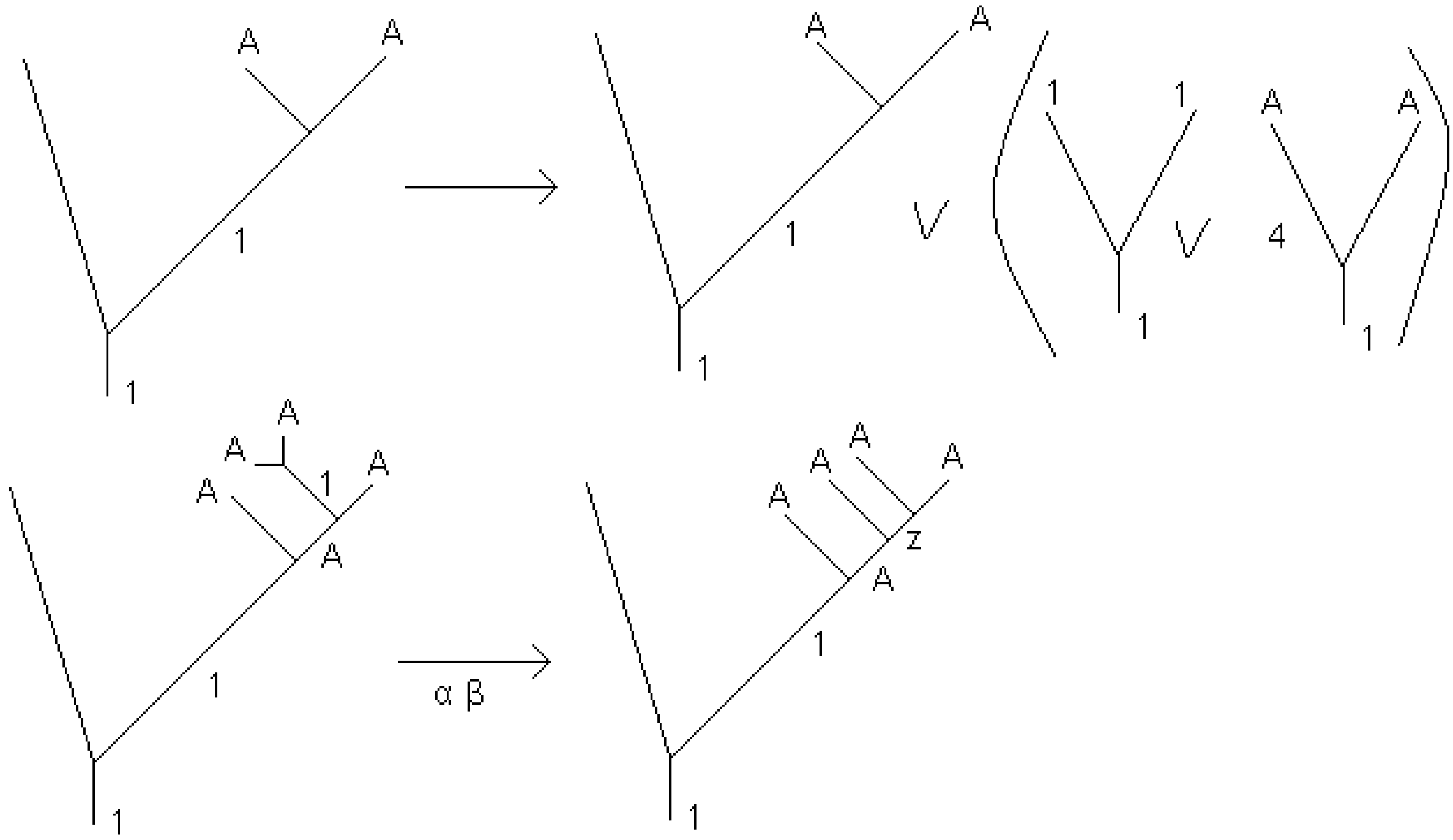

**Figure 313- Verify relation- the other case**

# 10 Stages for further research

## 10.1 A-C–invariants based on sliced 2-complexes

The list of topological relations which are required for A-C–invariance based on sliced 2-complexes in the Quinn model can serve as a checklist for:

a) consisting A-C –invariants
b) new constructed A-C –invariants.

For example in chapters 8.5 – 8.7 we have veryfied such a point from our list for an TQFT example.

## 10.2 The s-move and attached 3-cells

One problem of the invariants based on sliced 2-complexes is, that there is no expression in terms of the invariant for the simple-homotopy equivalence between the considered 2-complexes. It could be a way to construct A-C-invariants based on 3-cells, attached to the 2-complexes, or to be more precisely:
Consider the attached 3-cell as underlying space with the sliced relator discs on the boundary as a carrier for potential A-C-invariants based on sliced 3-cells.
I'm very far away to figure that out, but the basics behind that idea are described in this paper and for further details see [HoMeSier] and [Q3] .

We state a result of [HoMeSier] , where Metzler and Hog-Angeloni give an algebraic (on presentation oriented) criteria for the simple-homotopy equivalence between 2-complexes $K^2$ and $L^2$:

**Theorem**
Let $K^2$ and $L^2$ are compact connected CW complexes with presentation P respectively Q, $K^2$ is simple-homotopy equivalent to $L^2$ $\Leftrightarrow$
$K^2$ and $L^2$ can be 3–deformed, such that $P = \langle a_1,\ldots,a_n / R_1,\ldots,R_m \rangle$ and $Q = \langle a_1,\ldots,a_n / S_1,\ldots,S_m \rangle$ and for each μ:

$$R_\mu S_\mu^{\,-1} \in [\, N, N \,]$$

where $N = N(R_{\mu k}) = N(S_{\mu k})$ is the normal subgroup of the relators. Since we have $N = N(R_{\mu k}) = N(S_{\mu k})$ we can put the R relators in the first factor and the S relators in the second one:

$$R_\mu S_\mu^{\,-1} \in [\, N(R_{\mu,k}), N(S_{\mu,k}) \,]$$

Since the commutator of a product is a conjugate product of commutators, we can substitute the right side of the equation:

$$R_\mu S_\mu^{\,-1} = \prod_{k=1}^{p(k,\mu)} [\alpha_{\mu k} R^{*}_{r(\mu k)} \alpha_{\mu k}^{\,-1}, \beta_{\mu k} S^{*}_{s(\mu,k)} \beta_{\mu k}^{\,-1}] \qquad * = +-1$$

The product of commutators on the right side describes the relation in the fundamentalgroup of an orientable surface with genus g = p(k,μ), therefore the system of equations can be translated into a topologically context. We give an example [Q3] :

$R_1S_1^{-1} = [\alpha_{11}R_2\,\alpha_{11}^{-1}, \beta_{11}\,S_1\,\beta_{11}^{-1}]\,[\alpha_{12}\,R_1\,\alpha_{12}^{-1}, \beta_{12}\,S_2\,\beta_{12}^{\ -1}]$
$R_2S_2^{-1} = [\alpha_{21}R_1\,\alpha_{21}^{-1}, \beta_{21}\,S_2\,\beta_{21}^{-1}]$

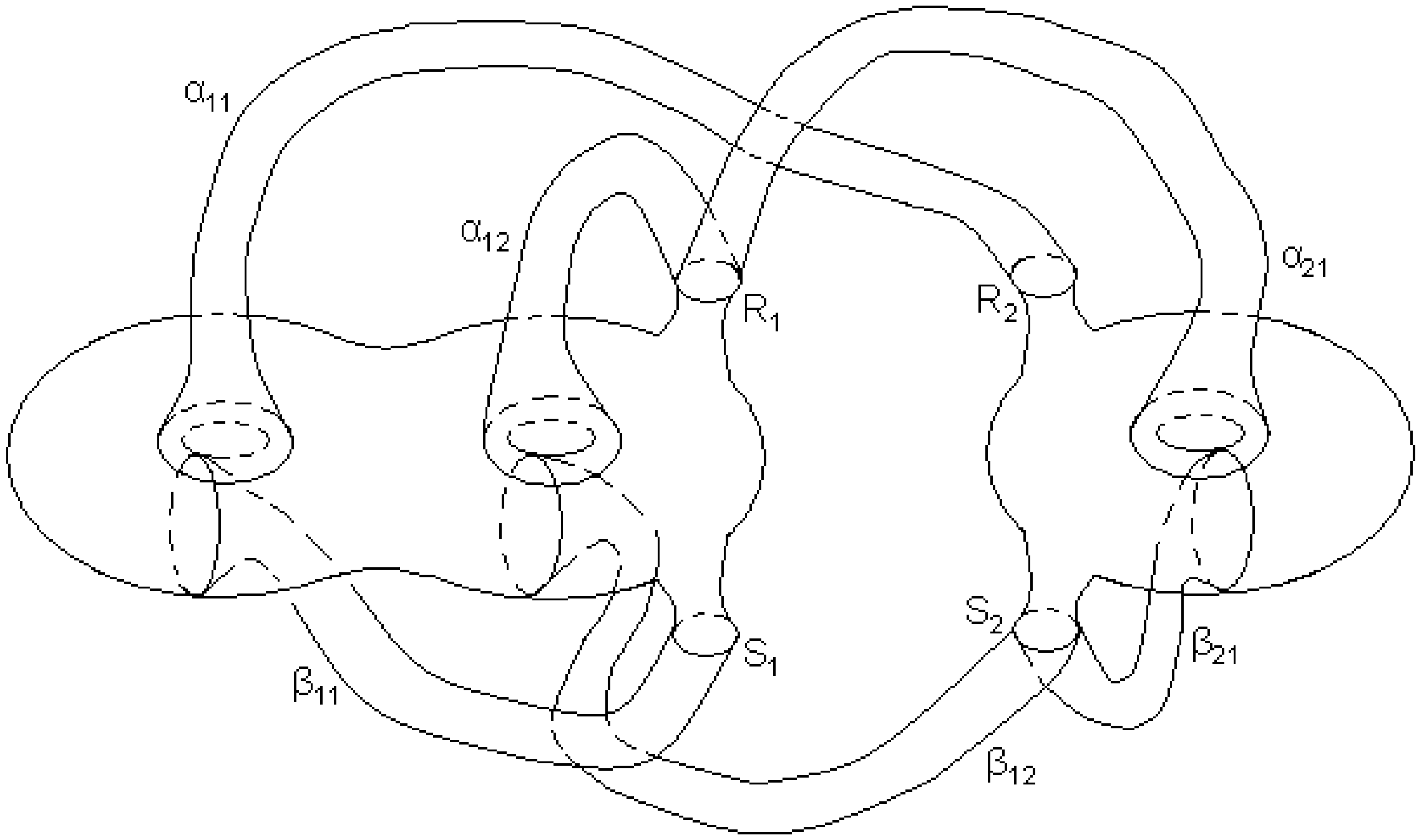

**Figure 314- s-move- difference as a commutator product**

We derive some facts from the equation (compare with the example):
Each commutator define a pair of generating curves of the corresponding orientable surface. Its genus is the number p(k,μ) of commutators in the product. Each generating curve is connected to a relator by exactly one annuli; a meridian connects to a single relator S* and a longitudinal connects to a single relator R*.
If r(μ,k) = s(μ,k), then R* and S* belong to the same surface, otherwise to different surfaces.
If r(μ,k) = μ, then R* (is equal to $R_\mu$) and its (via annuli) connected longitudinal belongs to the same surface (otherwise to different surfaces). Similar for S*.

Note, that this is an equation in the free group $F = \langle a_1,\ldots,a_n / - \rangle$, hence the surface together with the added annulis are mapped (in correspondence to their words) into the bouquet of $a_1,\ldots,a_n$ so it is a singular map to the common 1-skeleton of $K^2$ and $L^2$:

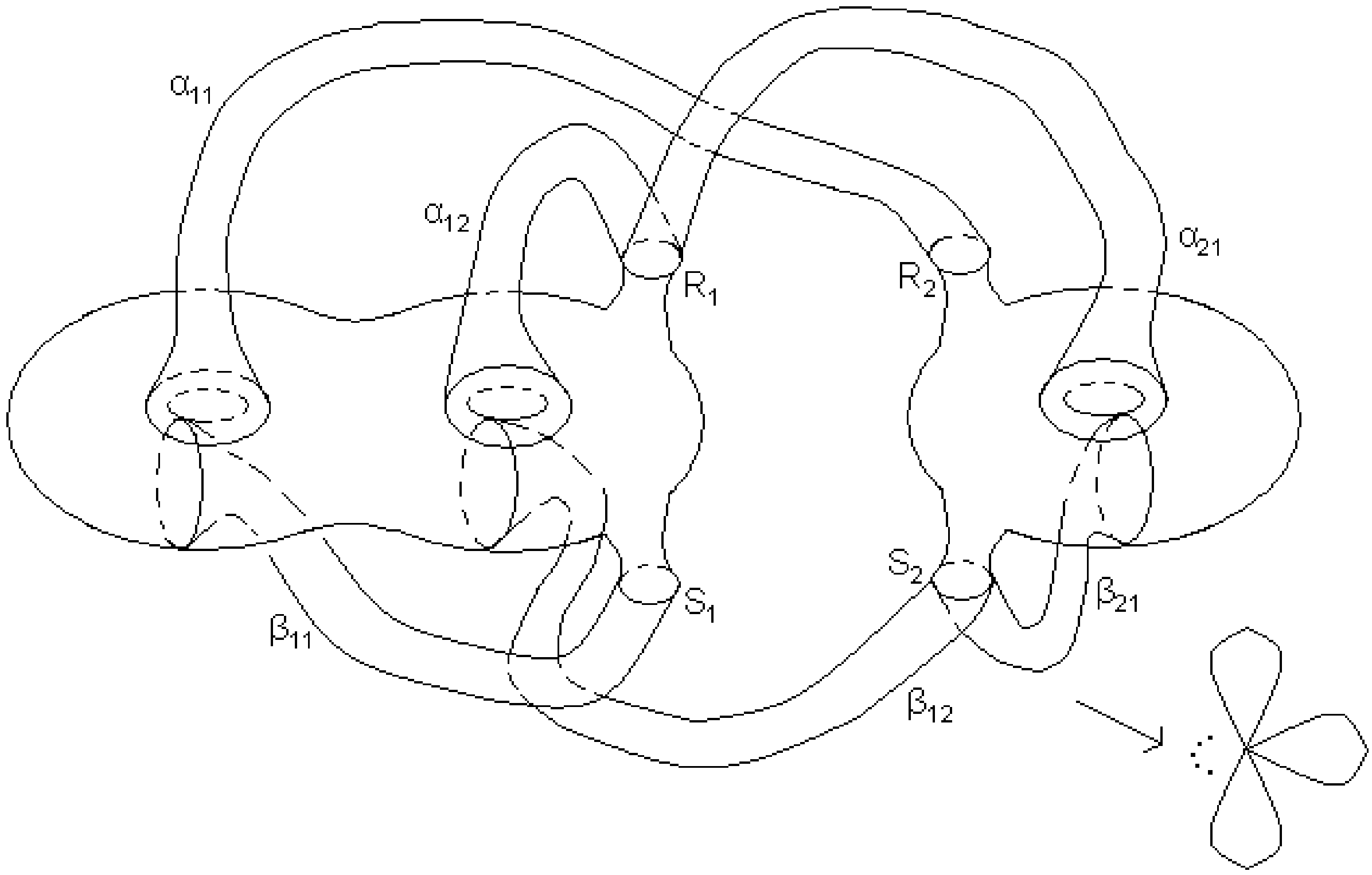

**Figure 315- s-move- singular map to 1-skeleton**

Following the notation of F. Quinn [Q3], we glue all relator discs of $K^2$ to the $R_\mu$ and all relator discs of $L^2$ to the $S_\mu$, we denote the resulting 2-complexes $\overline{K}^2$ and $\overline{L}^2$ as related by an s-move. Since the geometric difference between them (relator discs of $K^2$ respectively $L^2$ are hidden in the longitudinal curves respectively the meridian curves of the surface associated via annuli) maps singular to the 1-skeleton, the underlying space of $\overline{K}^2$ is still $K^2$ and of $\overline{L}^2$ is still $L^2$. However we can understand $\overline{K}^2$, $\overline{L}^2$ as images of 2-sphere maps, hence we have an attaching map for a 3-cell for each one.

In a simplified case we illustrate, how to attach 3-cells as an elementary expansion to $\overline{K}^2$ respectively $\overline{L}^2$.

Note, that we have to explain the attaching 2-sphere maps:

- Elementary 3 expansion from $\overline{K}^2$ to $\overline{K}^2 \vee \overline{L}^2$

Consider the 2-sphere with 3 subdics and a hole, first identify the subdics at the poles as indicated by the arrow, then we get a torus with a disc attached at the longitudinal curve, with remaining subdisc and a hole. Identify the remaining subdisc with the longitudinal disc, this connects both discs by an annuli. Identify the boundary of the hole with the meridian curve, this connects both curves by a second annuli. Compose the map of the sphere with the map of the constructed torus with added disc and annuli together to get the singular map to $\overline{K}^2$, where the longitudinal disc maps to an $R^*_{r(\mu k).}$ Fill in the 3-ball. Then the disc that fills the meridian curve is a free 2-cell and maps to an $S^*_{s(\mu k)}$ so it is a conjugate of a relator in $L^2$ via the annuli:

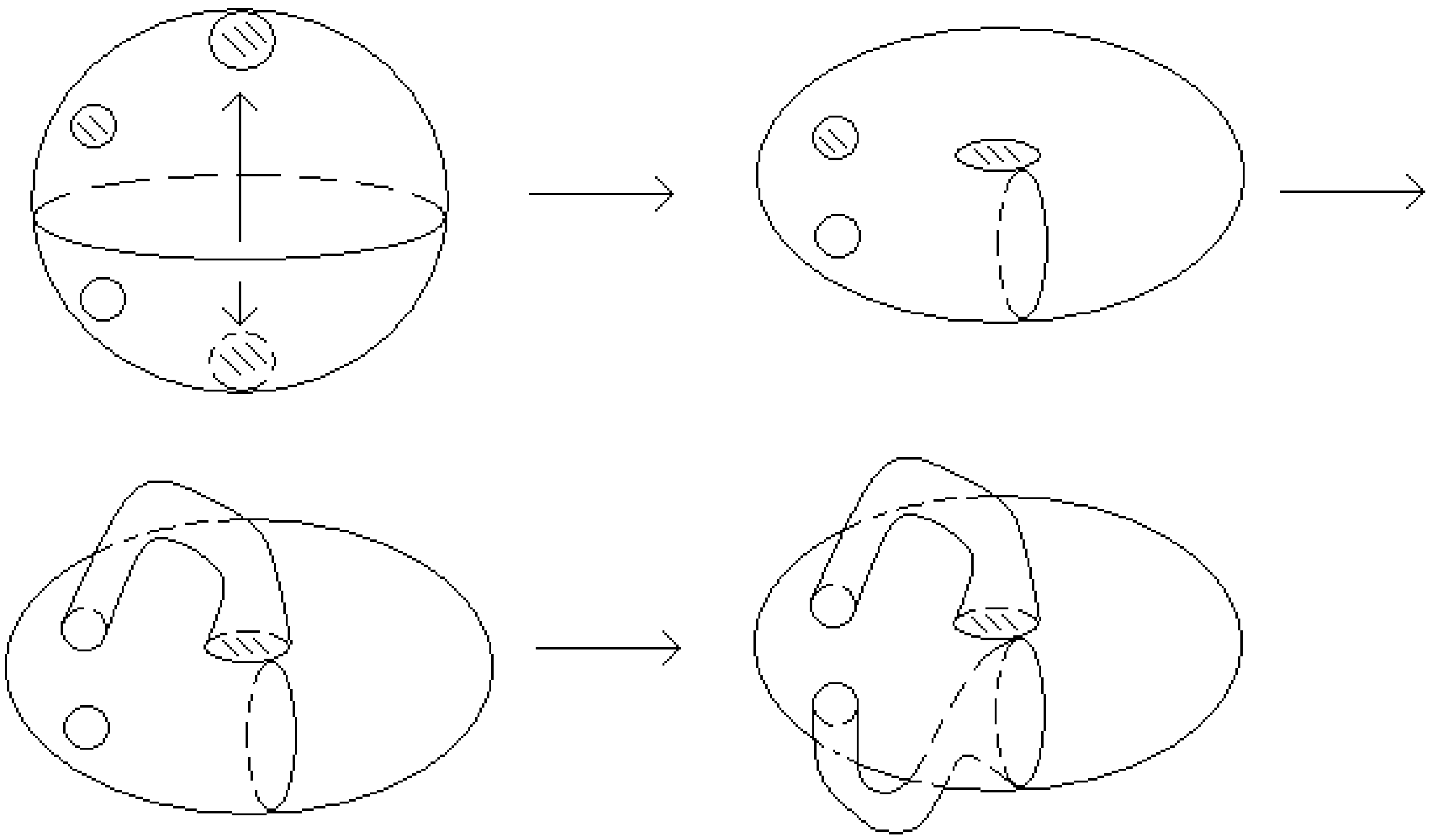

**Figure 316- s-move- identify longitudinal disc**

- Elementary 3 expansion from $\overline{L}^2$ to $\overline{L}^2 \vee \overline{K}^2$

Consider the 2-sphere with 3 subdics and a hole, first identify the subdics at the poles as indicated by the arrow, then we get a torus with a disc attached at the meridian curve, with remaining subdisc and a hole. Identify the remaining subdisc with the longitudinal disc, this connects both discs by an annuli. Identify the boundary of the hole with the longitudinal curve, this connects both curves by a second annuli. Compose the map of the sphere with the map of the constructed torus with added disc and annuli together to the singular map to $\overline{L}^2$, where the meridian disc maps to an $S^*_{s(\mu k)}$. Fill in the 3-ball. Then the disc that fills the longitudinal curve is a free 2-cell and maps to an $R^*_{r(\mu k)}$. so it is a conjugate of a relator in $K^2$ via the annuli:

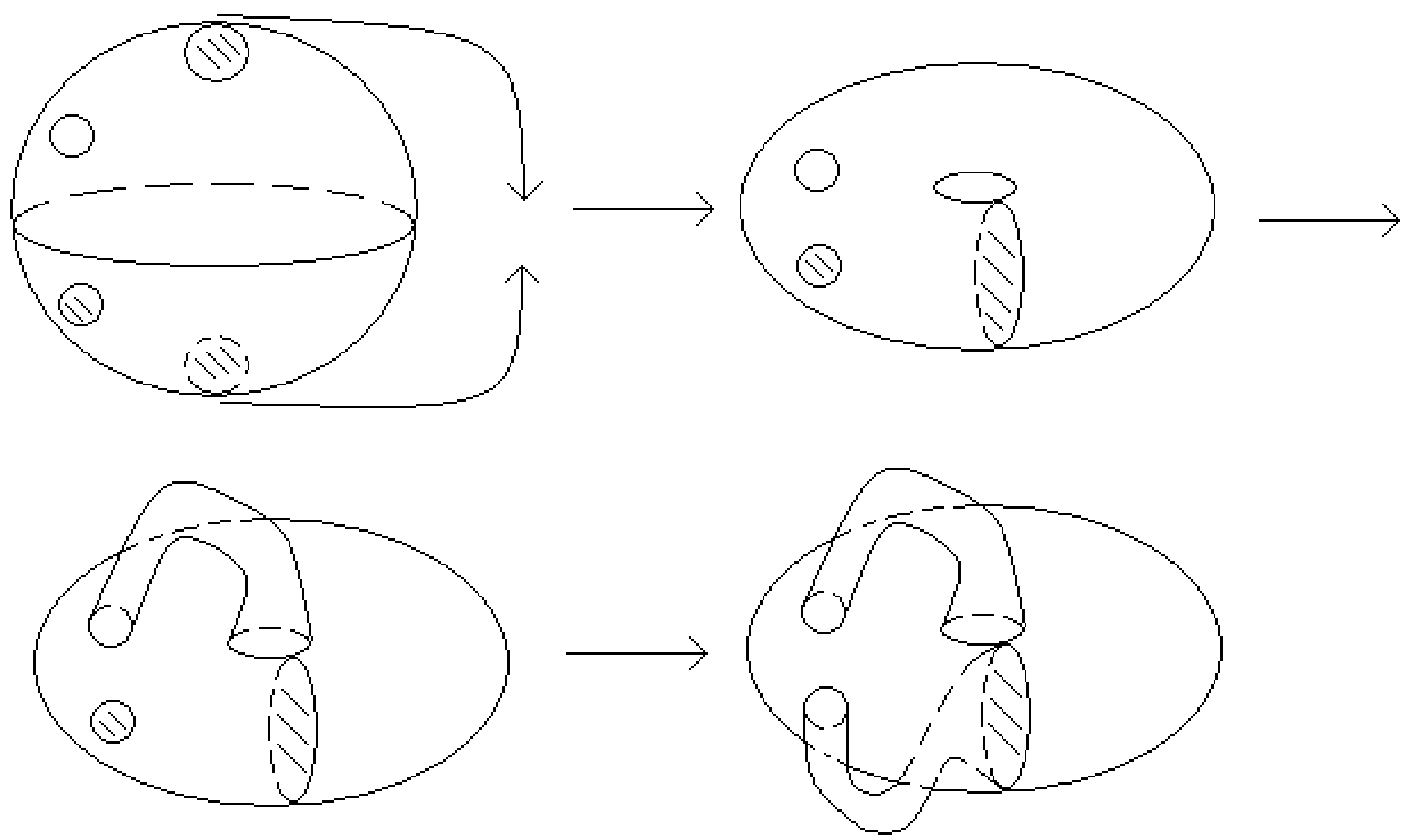

**Figure 317- s-move- identify meridian disc**

The general case:
For each orientable surface identify pairs of discs on the 2-sphere to get an orientable surface with connected annuli and discs filled in the half set of generating curves, one case for the set of meridian curves and the other case for the set of longitudinal curves. Connect the remaining set of generating curves by disjoint annulis. By filling in a 3-ball, the discs filling in the remaining set of generating curves, are free 2-cells, because each generating curve is connected to a relator by exactly one annuli.

2-complexes, which are related by an s-move are simple-homotopy equivalent, and this homotopy can be described geometrically.
It is sufficient to show, that the 2-sphere attaching maps are homotopic, but instead of $\overline{L}^2 \vee \overline{K}^2$ we consider a regular neigbourhood $N(\overline{L}^2 \vee \overline{K}^2)$. I would like thank S. Matveev, who explained that to me at a conference 1996 in Cheljabinsk.
Consider the identification map of $S^2 \times I$, the inner $S^2$ corresponds to the identification of the meridian disc, the outer $S^2$ to the identification of the longitudinal disc. The homotopy can be achieved by pushing of the outer sphere across the meridian and longitudinal discs into the inner sphere. In the next picture we illustrate that and draw the inner and outer sphere without identifications. Clearly this homotopy has to be composed with the singular map above:

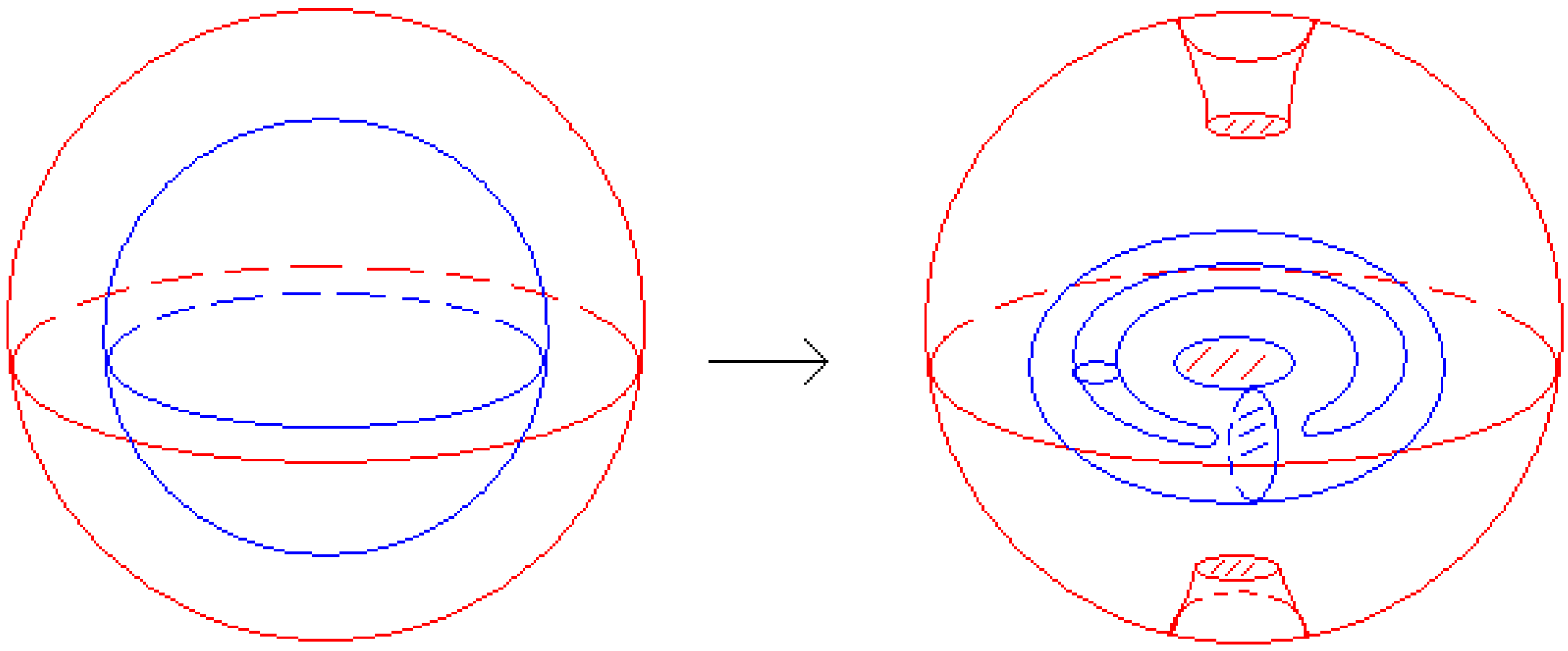

**Figure 318- s-move- homotopy for s-move**

The topological question is, how to slice the 3-cells in 2-dimensional pieces. We need a similar model (to the Quinn model of a 2-complex) for the 3-cells. Especially each step of the homotopy itself is an attaching map of a 3-cell, which then can be sliced. Hence the family of transitions of the 2-dimensional pieces during that homotopy may be translated in an expression for simple-homotopy equivalence.

## 11 List of References